\documentclass[fontsize=11pt,english,reqno]{amsart} 
\usepackage{ucs}
\usepackage{ae}
\usepackage{microtype}
\usepackage[utf8]{inputenc}
\usepackage[T1]{fontenc}
\usepackage[english]{babel}

\usepackage{amsmath}
\usepackage{amssymb}
\usepackage{amstext}
\usepackage{amsthm}

\usepackage{enumitem}

\usepackage{scrlayer-scrpage}
\usepackage{geometry}
\usepackage{xcolor}
\newcommand{\TODO}[1]{{\color{orange}*** #1 ***}}

\usepackage{mathtools}
\usepackage{slashed}
\usepackage{yhmath}
\usepackage{stmaryrd}

\usepackage{mathscinet}
\usepackage[numbers]{natbib}

\usepackage{graphicx}
\usepackage{tikz-cd}
\usepackage{caption}

\usepackage{imakeidx}
\makeindex[intoc]

\usepackage{tabularx}

\usepackage[colorlinks=true,hyperindex, linkcolor=magenta, urlcolor=black, pagebackref=false, citecolor=cyan,pdfpagelabels]{hyperref}
\usepackage[capitalize]{cleveref}

\newif\ifcomment
\newcommand{\comment}[1]{\ifcomment#1\fi}

\newlength{\currentparindent}
\makeatletter
\newcommand{\@minipagerestore}{\setlength{\parindent}{\currentparindent}}
\makeatother

\newcommand{\nospacepunct}[1]{\makebox[0pt][l]{#1}}

\usepackage{relsize}
\usepackage[bbgreekl]{mathbbol}
\usepackage{amsfonts}

\usepackage{accents}

\newcommand{\overcirc}[1]{\accentset{\circ}{#1}}

\DeclareSymbolFontAlphabet{\mathbb}{AMSb}
\DeclareSymbolFontAlphabet{\mathbbl}{bbold}

\newcommand{\prism}{\mathbbl{\Delta}}

\newcommand{\solid}{\mathsmaller{\square}}

\DeclareMathOperator{\Ker}{Ker}

\DeclareMathOperator{\id}{id}
\DeclareMathOperator{\Gal}{Gal}
\DeclareMathOperator{\Aut}{Aut}

\DeclareMathOperator{\Hom}{Hom}
\DeclareMathOperator{\Ext}{Ext}

\DeclareMathOperator{\Spec}{Spec}
\DeclareMathOperator{\Spa}{Spa}
\DeclareMathOperator{\Pic}{Pic}
\DeclareMathOperator{\Cone}{Cone}
\DeclareMathOperator{\Mod}{Mod}
\DeclareMathOperator{\Perf}{Perf}
\DeclareMathOperator{\gr}{gr}
\DeclareMathOperator{\Fil}{Fil}
\DeclareMathOperator{\Rep}{Rep}
\DeclareMathOperator{\Fun}{Fun}

\DeclareMathOperator{\Vect}{Vect}

\DeclareMathOperator{\Loc}{Loc}

\DeclareMathOperator{\fib}{fib}
\DeclareMathOperator{\cofib}{cofib}
\DeclareMathOperator{\RHom}{RHom}

\DeclareMathOperator{\Map}{Map}

\DeclareMathOperator{\Sym}{Sym}
\DeclareMathOperator{\Tot}{Tot}

\DeclareMathOperator{\Nil}{Nil}
\DeclareMathOperator{\AnSpec}{AnSpec}
\DeclareMathOperator{\Spd}{Spd}

\DeclareMathOperator{\LocConst}{LocConst}
\DeclareMathOperator{\GSpec}{GSpec}
\DeclareMathOperator{\DR}{DR}
\DeclareMathOperator{\Shv}{Shv}
\DeclareMathOperator{\Sing}{Sing}

\DeclareMathOperator{\SD}{SD}

\newcommand{\colim}{\operatornamewithlimits{colim}}
\newcommand{\Lim}{\operatornamewithlimits{lim}}
\let\lim\Lim

\newcommand{\PrL}{\mathrm{Pr}^\mathrm{L}}

\newcommand{\sHom}{\underline{\mathrm{Hom}}}

\newcommand{\Der}{\mathcal{D}\kern -.5pt er}

\newcommand\A{\mathbb{A}}
\newcommand\F{\mathbb{F}}
\newcommand\G{\mathbb{G}}
\newcommand\Q{\mathbb{Q}}
\newcommand\V{\mathbb{V}}
\newcommand\Z{\mathbb{Z}}
\newcommand\B{\mathbb{B}}
\newcommand\T{\mathbb{T}}
\newcommand\NN{\mathbb{N}}
\newcommand\DD{\mathbb{D}}
\newcommand\C{\mathbb{C}}
\newcommand\TT{\mathbb{T}}

\let\T\TT
\let\P\PP

\newcommand\D{\mathcal{D}}
\newcommand\DF{\mathcal{DF}}

\let\O\cO

\newcommand\crys{\mathrm{crys}}
\newcommand\dR{\mathrm{dR}}
\newcommand\N{\mathrm{N}}
\newcommand\Syn{\mathrm{Syn}}
\newcommand\et{\mathrm{\acute{e}t}}
\newcommand\proet{\mathrm{pro\acute{e}t}}

\newcommand\op{\mathrm{op}}

\let\inf\Ainf
\newcommand\HT{\mathrm{HT}}

\newcommand\red{\mathrm{red}}
\newcommand\Hod{\mathrm{Hod}}

\newcommand\can{\mathrm{can}}

\newcommand\an{\mathrm{an}}
\newcommand\Div{\mathrm{Div}}

\newcommand\cycl{\mathrm{cycl}}
\newcommand\FF{\mathrm{FF}}

\newcommand\la{\mathrm{la}}

\newcommand\HK{\mathrm{HK}}
\newcommand\Betti{\mathrm{Betti}}

\newcommand\st{\mathrm{st}}

\newcommand\NilPerfd{\mathrm{NilPerfd}}

\newcommand\arcStk{\mathrm{arcStk}}

\newcommand\sm{\mathrm{sm}}
\newcommand\cont{\mathrm{cont}}

\newcommand\GelfStk{\mathrm{GelfStk}}
\newcommand\un{\mathrm{un}}
\newcommand\pcrys{\mathrm{pcrys}}

\newcommand\AbGrp{\mathrm{AbGrp}}
\newcommand\arithm{\mathrm{arithm}}
\newcommand\geom{\mathrm{geom}}
\newcommand\mHK{{\HK, \mathrm{mock}}}
\newcommand\mDiv{{\Div^1, \mathrm{mock}}}
\newcommand\GelfRing{\mathrm{GelfRing}}
\newcommand\Ani{\mathrm{Ani}}

\newcommand\alg{\mathrm{alg}}

\let\epsilon\varepsilon
\let\phi\varphi
\let\ol\overline
\let\ul\underline
\let\tensor\otimes

\let\cal\mathcal
\let\diamond\diamondsuit

\let\Map\map

\newtheorem{thm}{Theorem}[section]
\newtheorem{prop}[thm]{Proposition}
\newtheorem{lem}[thm]{Lemma}
\newtheorem{cor}[thm]{Corollary}
\newtheorem{conj}[thm]{Conjecture}

\theoremstyle{definition}
\newtheorem{defi}[thm]{Definition}

\newtheorem{warn}[thm]{Warning}
\newtheorem{exx}[thm]{Example}
\newenvironment{ex}
  {\pushQED{\qed}\exx}
  {\popQED\endexx}
\newtheorem{remm}[thm]{Remark}
\newenvironment{rem}
  {\pushQED{\qed}\remm}
  {\popQED\endremm}
  
\AddToHook{env/prop/begin}{\crefalias{thm}{prop}}
\AddToHook{env/lem/begin}{\crefalias{thm}{lem}}
\AddToHook{env/cor/begin}{\crefalias{thm}{cor}}
\AddToHook{env/conj/begin}{\crefalias{thm}{conj}}
\AddToHook{env/defi/begin}{\crefalias{thm}{defi}}
\AddToHook{env/thmdefi/begin}{\crefalias{thm}{thmdefi}}
\AddToHook{env/lemdefi/begin}{\crefalias{thm}{lemdefi}}
\AddToHook{env/conv/begin}{\crefalias{thm}{conv}}
\AddToHook{env/rem/begin}{\crefalias{thm}{remm}}
\AddToHook{env/warn/begin}{\crefalias{thm}{warn}}
\AddToHook{env/ex/begin}{\crefalias{thm}{exx}}

\numberwithin{equation}{subsection}

\makeatletter
\renewcommand{\address}[1]{\gdef\@address{#1}}
\renewcommand{\email}[1]{\gdef\@email{\url{#1}}}
\newcommand{\@endstuff}{\par\vspace{\baselineskip}\noindent\small
\begin{tabular}{@{}l}\scshape\@address\\\textit{E-mail address:} \@email\end{tabular}}
\AtEndDocument{\@endstuff}
\makeatother

\title[Rational syntomic cohomology]{Rational analytic syntomic cohomology}
\author{Maximilian Hauck}
\address{Max-Planck-Institut f\"ur Mathematik, Vivatsgasse 7, 53111 Bonn, Germany}
\email{max.hauck01@gmail.com}

\keywords{$p$-adic cohomology, $p$-adic Hodge theory, étale cohomology, de Rham cohomology, prismatic cohomology, de Rham local systems}
\subjclass[2020]{14F20, 14F30, 14F40, 14G45}

\begin{document}

\begin{abstract}
We define and study the rational analytic syntomification $X^\Syn$ of a partially proper rigid-analytic variety $X$ over $\Q_p$. We establish Poincaré duality and a theory of first Chern classes for the resulting cohomology theory, identify vector bundles on $X^\Syn$ with de Rham bundles on the Fargues--Fontaine curve of $X^\diamond$ and recover several classical comparison theorems in $p$-adic Hodge theory. We also develop analogues of our results and constructions over $\C_p$.
\end{abstract}

\maketitle
{\hypersetup{hidelinks} \tableofcontents}

\newpage

\section{Introduction}

Let $p$ be a prime number.

\subsection{Motivation and background}

Let $X$ be a smooth proper variety over $\Q_p$ and let $\ol{\Q}_p$ denote a fixed algebraic closure of $\Q_p$. The central question of $p$-adic Hodge theory is the following:

\bigskip

\begin{center}
What can be said about the relation between de Rham cohomology $R\Gamma_\dR(X)$ \\
and $p$-adic étale cohomology $R\Gamma_\et(X_{\ol{\Q}_p}, \Z_p)$ of $X$?
\end{center}

\bigskip

The first concrete conjectures aimed at answering this question were formulated by Fontaine in \cite{Fontaine} and subsequently studied by many authors, most notably Faltings \cite{Faltings}, Tsuji \cite{Tsuji} and Nizio{\l} \cite{Niziol}. Following the introduction of perfectoid spaces in \cite{Perfectoid}, similar questions for rigid-analytic varieties instead of algebraic varieties have also been studied in detail, particularly by Scholze \cite{PAdicHodgeTheory}, Colmez--Nizio{\l} \cite{padicComparisons}, \cite{BasicComparison} and Bosco \cite{BoscopAdic}. Indeed, the focus of this paper will be this latter case of rigid-analytic varieties.

These developments have spawned a whole zoo of new $p$-adic cohomology theories for algebraic or, more generally, rigid-analytic varieties over $p$-adic fields. On the de Rham side, this includes the so-called log-crystalline cohomology of Hyodo and Kato from \cite{Hyodo}, \cite{HyodoKato}, see also \cite{Beilinson}, which is nowadays often simply called \emph{Hyodo--Kato cohomology} and has been adapted to the rigid-analytic setting by Colmez--Nizio{\l} in \cite{padicComparisons}, and, in the geometric setting, i.e.\ over the completion $\C_p$ of $\ol{\Q}_p$, the \emph{$B_\dR^+$-cohomology} introduced in \cite[§13]{IntegralpAdicHT}, see also \cite{GuoBdR}. On the étale side, Scholze has introduced a variant of étale cohomology of rigid-analytic varieties called \emph{proétale cohomology} in \cite{PAdicHodgeTheory}. Similarly to the guiding question above, much of $p$-adic Hodge theory is centred around proving comparison theorems between all these cohomology theories and trying to understand how they are related.

The idea behind syntomic cohomology is that there should be a single cohomology theory that mediates all of these comparison theorems. It was first introduced in the algebraic setting by Fontaine--Messing in \cite{FontaineMessing} and later defined for $p$-adic formal schemes by Bhatt--Morrow--Scholze in \cite[§7.4]{THHAndPAdicHT} as well as for rigid-analytic varieties by Colmez--Nizio{\l}, see \cite{padicComparisons}. It is the $p$-adic analogue of Deligne cohomology over the complex numbers and intimately related to $K$-theory and $p$-adic regulators, see e.g.\ \cite[§7.4]{THHAndPAdicHT} or \cite{NekovarNiziol}, and hence also of interest beyond the scope of $p$-adic comparison theorems.

The major trend in $p$-adic Hodge theory in recent years has been geometrisation. This is the idea that, instead of directly studying some cohomology theory attached to a variety $X$, one constructs a new space attached to $X$ which is a geometric avatar of the cohomology theory under consideration; for this process, Bhatt has coined the term \emph{transmutation}. In this approach, one recovers the original cohomology theory by taking cohomology of the structure sheaf on the transmuted space. In the case of de Rham cohomology, this idea goes back to work of Simpson, see \cite{Simpson} and \cite{Simpson2}, and it has entered the sphere of $p$-adic Hodge theory through recent work of Drinfeld and Bhatt--Lurie on such a ``stacky approach'' for prismatic and syntomic cohomology of $p$-adic formal schemes, see \cite{Prismatization}, \cite{APC}, \cite{PFS} and \cite{FGauges}. 

This geometrisation approach is advantageous for several reasons: Conceptually, it allows one to think about cohomology theories in a geometric way and, in particular, to reinterpret comparison theorems as certain statements about the geometry of the ``transmuted'' spaces. It also provides a theory of coefficients for a cohomology theory essentially for free: one simply replaces the structure sheaf on the transmuted space by a vector bundle or even a perfect complex. Finally, assuming that one has a six-functor formalism for quasicoherent sheaves on the transmuted spaces, one immediately obtains a full six-functor formalism for the original cohomology theory.  

While this picture has been almost fully developed for $p$-adic cohomology theories of $p$-adic formal schemes, much less is known for rigid-analytic varieties. The first steps in this direction have been taken by Rodríguez Camargo, who has introduced the analytic de Rham stack in \cite{dRStack}, see also the recent update \cite{dRFF}, and Anschütz--Le Bras--Rodríguez Camargo--Scholze, who have proposed a theory of \emph{analytic prismatisation} for rigid-analytic varieties with $\Q_p$-coefficients. However, the full picture for syntomic cohomology has been missing. In this paper, we develop the theory of \emph{analytic syntomification} for rigid-analytic varieties with $\Q_p$-coefficients, making use of the framework of \emph{Gelfand stacks} recently introduced in \cite{dRFF}. This is a variant of Clausen--Scholze's analytic stacks tailored to $p$-adic analytic geometry and admits a full six-functor formalism for quasicoherent sheaves.

More precisely, for any partially proper rigid-analytic variety $X$ over $\Q_p$, we define its analytic syntomification $X^\Syn$, which is a Gelfand stack over $\Q_p$, and taking quasicoherent cohomology on $X^\Syn$ yields a notion of rational analytic syntomic cohomology of $X$. For this cohomology theory, we prove Poincaré duality as well as the classical comparison theorems with Hyodo--Kato cohomology and proétale cohomology by studying the geometry of the stack $X^\Syn$. Using similar methods, we also recover the de Rham comparison theorem of Scholze from \cite{PAdicHodgeTheory}. Moreover, if $X$ is smooth, we identify the category of vector bundles on $X^\Syn$ with a certain full subcategory of ``de Rham'' vector bundles on the relative Fargues--Fontaine curve of the diamond of $X$. Finally, we also provide variants of these results and constructions for rigid-analytic varieties over $\C_p$.

\subsection{Rational analytic syntomification over $\Q_p$}

We now briefly describe the construction of $X^\Syn$, which works in the generality of $X$ being any Gelfand stack over $\Q_p$. Namely, similarly to the algebraic story developed by Drinfeld and Bhatt--Lurie in \cite{Prismatization} and \cite{FGauges}, we will construct an \emph{analytic Nygaardification} $X^\N$ of $X$ which receives two maps 
\begin{equation*}
j_\dR, j_\HT: X^\prism\rightarrow X^\N
\end{equation*}
from the \emph{analytic prismatisation} $X^\prism$ of $X$, the theory of which has been developed by Anschütz--Le Bras--Rodríguez Camargo--Scholze. Then $X^\Syn$ is obtained by gluing the two copies of $X^\prism$ embedded into $X^\N$ via $j_\dR$ and $j_\HT$.

Crucial for our constructions will be the notion of \emph{nilperfectoid rings} from \cite[Def.\ 4.6.1]{dRFF}. These are Gelfand rings $A$ whose $\dagger$-reduction $\ol{A}\coloneqq A^{\dagger\text{-red}}$ is perfectoid, and they form a basis of the $!$-topology on Gelfand rings. In particular, we may define Gelfand stacks by only describing their $A$-valued points for nilperfectoid rings $A$.

With this in mind, let us first recall the definition of the analytic prismatisation $\Q_p^\prism$ of $\Q_p$ in broad strokes. For any nilperfectoid $A$, one constructs a $\dagger$-nilpotent thickening $Y_A$ of the curve $Y_{\ol{A}}$ from \cite[§4.4]{dRFF}. Most importantly, this thickening $Y_A$ is equipped with a Frobenius $\phi: Y_A\rightarrow Y_A$ and a map $\iota: \GSpec A\rightarrow Y_A$ lifting the usual such structures in the perfectoid case $A=\ol{A}$. Then the Gelfand stack $\Q_p^\prism$ is obtained by sheafifying the assignment
\begin{equation*}
\Q_p^\prism(A)\coloneqq \{\text{degree $1$ Cartier divisors $D\subseteq Y_A$}\}
\end{equation*}
for nilperfectoids $A$. Here, a Cartier divisor $D\subseteq Y_A$ is called \emph{of degree $1$} if its pullback to $Y_{\ol{A}}$ has degree $1$ in the sense of \cite[Def.\ II.1.19]{FarguesScholze}.

To construct $\Q_p^\N$, the most important feature of $\Q_p^\prism$ is that it admits a map $\widetilde{\mu}: \Q_p^\prism\rightarrow (\A^1/\G_m)^\dR$ given by pulling back a degree $1$ Cartier divisor along the map $\phi^{-1}\circ\iota: \GSpec\ol{A}\rightarrow Y_A$, which is a Frobenius twist of the map $\iota$ from above. In fact, we will show that this may be refined to a map $\widetilde{\mu}: \Q_p^\prism\rightarrow (\overcirc{\DD}/\ol{\T})^\dR$, where $\overcirc{\DD}$ denotes the open unit disk and $\ol{\T}$ is the overconvergent unit torus over $\Q_p$. Mimicking an alternative definition of the algebraic filtered prismatisation recently given by Gardner--Madapusi in \cite[§6]{GardnerMadapusi}, we then define $\Q_p^\N$ as the pullback
\begin{equation*}
\begin{tikzcd}
\Q_p^\N\ar[r, "\pi"]\ar[d, "{(t, u)}", swap] & \Q_p^\prism\ar[d, "\widetilde{\mu}"] \\
\ol{\DD}/\ol{\T}\times (\ol{\DD}/\ol{\T})^\dR\ar[r, "\mathrm{mult}"] & (\ol{\DD}/\ol{\T})^\dR\nospacepunct{\;,}
\end{tikzcd}
\end{equation*}
where $(-)^\dR$ denotes the analytic de Rham stack from \cite{dRFF}, $\ol{\DD}$ is the overconvergent unit disk and the map $\mathrm{mult}$ is induced by the multiplication map on $\ol{\DD}$.

Following Bhatt's transmutation philosophy, defining $X^\N$ for any Gelfand stack $X$ over $\Q_p$ now comes down to defining a (Gelfand) ring stack $\G_a^\N$ over $\Q_p^\N$; then $X^\N$ should parametrise maps from $\GSpec\G_a^\N(A)$ to $X$ for any Gelfand ring $A$ over $\Q_p^\N$. In fact, we will proceed slightly differently and directly define $\GSpec\G_a^\N(A)$. Namely, given an $A$-valued point of $\Q_p^\N$ for $A$ nilperfectoid, we will construct a $\dagger$-nilpotent thickening $D^+$ of the divisor $D\subseteq Y_A$ classified by the induced map $\GSpec A\rightarrow\Q_p^\prism$. Then $X^\N$ is defined as the Gelfand stack over $\Q_p^\N$ obtained as the sheafification of
\begin{equation*}
X^\N(\GSpec A\rightarrow\Q_p^\N)\coloneqq \{\text{maps }D^+\rightarrow X\}\;.
\end{equation*}

Finally, recall the maps $t: \Q_p^\N\rightarrow \ol{\DD}/\ol{\T}$ and $u: \Q_p^\N\rightarrow (\ol{\DD}/\ol{\T})^\dR$ from above, which induce corresponding maps on $X^\N$. Noting that the norms of $t$ and $u$ are well-defined, it makes sense to consider the loci $\{|t|=1\}$ and $\{|u|=1\}$ inside $X^\N$. We will show that the corresponding substacks of $X^\N$ are both isomorphic to $X^\prism$, and this is what defines the maps $j_\dR$ and $j_\HT$, respectively.

To define analytic syntomic cohomology of any Gelfand stack $X$ from here, we construct a line bundle $\O\{1\}$ on $\Q_p^\N$ playing the role of a Breuil--Kisin twist and show that it descends to $\Q_p^\Syn$. Then we set
\begin{equation*}
R\Gamma_\Syn(X, \Q_p(i))\coloneqq R\Gamma(X^\Syn, \O\{i\})
\end{equation*}
for any $i\in\Z$ and call it the \emph{rational analytic syntomic cohomology} of $X$ of weight $i$. Note that we may just as well take cohomology of any other vector bundle or even perfect complex on $X^\Syn$, and this yields a theory of analytic syntomic cohomology with coefficients. In keeping with the algebraic case from \cite{FGauges}, the category $\D(X^\Syn)$ of such coefficients will be called the category of \emph{analytic (prismatic) $F$-gauges} on $X$.

\begin{rem}
We briefly sketch how to extend the definition of $X^\Syn$ to integral coefficients. In fact, the most difficult part of this is to establish a theory of analytic prismatisation with integral coefficients, which is the subject of ongoing joint work of the author with Anschütz--Le Bras--Rodríguez Camargo--Scholze. Once this is in place, the definitions of $X^\N$ and $X^\Syn$ carry over verbatim. Indeed, in our setting where $X$ is a Gelfand stack over $\Q_p$, it turns out that the added complexity of $X^\Syn$ in the integral coefficient case then purely comes from $X^\prism$. This allows one to extend many of our results to the case of integral coefficients.
\end{rem}

\begin{rem}
The fact that our construction of $X^\Syn$ takes as an input a Gelfand stack $X$ over $\Q_p$ explains our restriction to rigid-analytic spaces which are partially proper: As explained in \cite{dRFF}, partially proper rigid spaces are the ones for which the functor to Gelfand stacks is fully faithful. In general, passing to Gelfand stacks identifies a rigid-analytic space with its Huber compactification; e.g., the rigid disk $\Spa(\Q_p\langle T\rangle, \Z_p\langle T\rangle)$ is identified with $\Spa(\Q_p\langle T\rangle, \Z_p)$, which has an additional rank-$2$-point.
\end{rem}

\begin{rem}
Continuing the preceding remark, let us note that we will in fact mostly work with (derived) Berkovich spaces $X$ over $\Q_p$ in the sense of \cite[§4.3]{dRFF}, which form a full subcategory of Gelfand stacks by \cite[Prop.\ 4.3.3]{dRFF}. However, we will assume $X$ to be a partially proper rigid space whenever we want to connect to results in the $p$-adic Hodge theory of rigid-analytic varieties. We point out that, as a third option, we could also work with $\dagger$-rigid spaces in the sense of \cite{GrosseKlonne} as these also form a full subcategory of derived Berkovich spaces; however, we will mostly refrain from doing this.
\end{rem}

\subsection{Main results}

We now state the main results of the paper, which roughly fall into three categories: (1) Interaction of the syntomification with six functors, (2) comparison theorems for $p$-adic cohomology theories and (3) variants in the geometric setting, i.e.\ over $\C_p$.

\subsubsection{Poincaré duality and Chern classes}

On the foundational side, we establish basic properties of the functor $X\mapsto X^\Syn$ with respect to the six-functor formalism for quasi-coherent sheaves on Gelfand stacks. The first of these concerns the case $X=\GSpec\Q_p$ and yields a form of Poincaré duality for cohomology on the stack $\Q_p^\Syn$.

\begin{thm}[\cref{thm:sixfunctors-qpsynsmooth}]
The map $f: \Q_p^\Syn\rightarrow\GSpec\Q_p$ is cohomologically smooth with dualising sheaf $\O\{1\}[3]$, i.e. for any $E\in\Perf(\Q_p^\Syn)$, there is an isomorphism
\begin{equation*}
f_*(E^\vee\{1\})[3]\cong (f_!E)^\vee\;.
\end{equation*}
\end{thm}

\begin{rem}
Normally, one would expect the dualising sheaf to be concentrated in degree $-2$ instead of $-3$, e.g.\ by comparing with the duality obtained by Bhatt--Lurie for the algebraic syntomification of $\Z_p$ in \cite[Cor.\ 6.4.6]{FGauges} or with local Tate duality for Galois representations. The extra shift here is due to the fact that $\Q_p^\Syn$ has an extra ``topological dimension'' that roughly arises from the Betti stack of an interval. However, this ``extra interval'' does not affect cohomology of perfect complexes on $\Q_p^\Syn$ by contractibility.
\end{rem}

\begin{rem}
We do \emph{not} expect $\Q_p^\Syn$ to be cohomologically proper (or even prim). Nevertheless, cohomology on $\Q_p^\Syn$ preserves perfect complexes, which can be deduced from \cref{thm:intro-hkcomp} below.
\end{rem}

The next result concerns the relative situation, i.e.\ it is about induced maps between syntomifications.

\begin{thm}[\cref{thm:sixfunctors-smoothmaps}, \cref{thm:sixfunctors-propermaps}]
Let $f: X\rightarrow Y$ be a map of derived Berkovich spaces in the sense of \cite[§4.3]{dRFF}.
\begin{enumerate}[label=(\alph*)]
\item If $f$ is rigid smooth, then the induced morphism $f^\Syn: X^\Syn\rightarrow Y^\Syn$ is cohomologically smooth. If $X$ is moreover a smooth partially proper rigid space over $\Q_p$ and $f$ is pure of relative dimension $d$, the dualising sheaf of $f^\Syn$ is given by $\O\{d\}[2d]$.
\item If $f$ is locally of finite presentation and quasicompact, the induced morphism $f^\Syn: X^\Syn\rightarrow Y^\Syn$ is weakly cohomologically proper.
\end{enumerate}
In particular, if $f: X\rightarrow Y$ is a smooth proper map of smooth partially proper rigid spaces pure of relative dimension $d$, then $f^\Syn_*$ preserves perfect complexes and we have
\begin{equation*}
f^\Syn_*(E^\vee\{d\})[2d]\cong (f^\Syn_*E)^\vee
\end{equation*}
for any $E\in\Perf(X^\Syn)$.
\end{thm}

\begin{rem}
In part (a), we expect that the assumption that $X$ is a smooth partially proper rigid space should not be necessary for the identification of the dualising sheaf to hold.
\end{rem}

Finally, we show that syntomic cohomology of derived Berkovich spaces admits a strong theory of first Chern classes in the sense of \cite[§5]{Zavyalov}. More precisely, this means the following:

\begin{prop}[\cref{prop:sixfunctors-chern}]
For any derived Berkovich space $X$, there is a natural morphism
\begin{equation*}
c_1^\Syn: R\Gamma(X_{\mathrm{Berk}\text{-}\et}, \G_m)[1]\rightarrow R\Gamma(X^\Syn, \O\{1\}[2])\;,
\end{equation*}
where the left-hand side denotes the cohomology of $\G_m[1]$ on the Berkovich étale site of $X$. Moreover, the projective bundle formula holds: for any $d\geq 1$ and any derived Berkovich space $X$, the induced morphism
\begin{equation*}
\sum_{0\leq k\leq d} (c_1^\Syn)^k\{d-k\}[2d-2k]: \bigoplus_{0\leq k\leq d} \O\{d-k\}[2d-2k]\rightarrow f_*^\Syn\O\{d\}[2d]
\end{equation*}
is an isomorphism, where $f: \P^d_X\rightarrow X$ is the projection.
\end{prop}

\subsubsection{Comparison theorems}

By studying the geometry of the syntomification, we can recover many of the classical comparison theorems in the $p$-adic Hodge theory of rigid-analytic varieties and, in many cases, upgrade them to coefficients in arbitrary perfect complexes on $X^\Syn$. This will culminate in a description of the category $\Vect(X^\Syn)$ for a smooth partially proper rigid space $X$ over $\Q_p$ in terms of de Rham bundles on the relative Fargues--Fontaine curve of $X^\diamond$. In particular, this shows that vector bundle analytic $F$-gauges are actually very classical objects in $p$-adic Hodge theory.

We begin with Scholze's de Rham comparison theorem from \cite{PAdicHodgeTheory}. For this, we first use the geometry of $X^\Syn$ to construct a functor
\begin{equation}
\label{eq:intro-scholzefunctor}
\Perf(X^{\dR, +})\rightarrow\Perf(X_\proet, \mathbb{B}_\dR^+)
\end{equation}
for any smooth partially proper rigid space $X$ over $\Q_p$. Here, $X^{\dR, +}$ denotes the filtered analytic de Rham stack of $X$ from \cite[Def.\ 5.2.2]{dRStack}, perfect complexes on which identify with filtered perfect complexes with connection on $X$. 

\begin{thm}[\cref{thm:padicht-scholzefunctor}]
The functor
\begin{equation*}
\{\text{filtered vector bundles with connection on $X$}\}\rightarrow\{\text{$\mathbb{B}_\dR^+$-local systems on $X_\proet$}\}
\end{equation*}
obtained by restricting (\ref{eq:intro-scholzefunctor}) to $\Vect(X^{\dR, +})$ identifies with Scholze's functor
\begin{equation*}
(\Fil^\bullet E, \nabla)\mapsto \Fil^0(E\tensor_{\O_X} \O\mathbb{B}_\dR)^{\nabla=0}
\end{equation*}
from \cite[§7]{PAdicHodgeTheory}.
\end{thm}

Using this, the version of the de Rham comparison theorem we recover takes the following form:

\begin{thm}[\cref{thm:padicht-drcomp}]
Let $f: X\rightarrow Y$ be a smooth proper morphism of smooth partially proper rigid spaces over $\Q_p$ and let $\mathbb{L}$ be a de Rham local system on $X$ with associated filtered vector bundle with connection $E\in\Vect(X^{\dR, +})$. Then the pushforward $f^{\dR, +}_*E\in\D(Y^{\dR, +})$ is perfect and its image under the functor (\ref{eq:intro-scholzefunctor}) is given by $f_{\proet, *}\mathbb{L}\tensor_{\Q_p} \mathbb{B}_\dR^+$. In particular, for any smooth proper rigid space $X$ over $\Q_p$, there is an isomorphism
\begin{equation*}
R\Gamma_\proet(X_{\C_p}, \Q_p)\tensor_{\Q_p} B_\dR\cong R\Gamma_\dR(X)\tensor_{\Q_p} B_\dR
\end{equation*}
compatible with the $\Gal_{\Q_p}$-actions and the filtrations.
\end{thm}

We can also recover the classical comparison theorems of Colmez--Nizio{\l} between syntomic cohomology and Hyodo--Kato cohomology or proétale cohomology, respectively, from \cite{padicComparisons}. In fact, we are able to extend them to coefficients in an arbitrary perfect analytic $F$-gauge and reinterpret them as descriptions of the full category of perfect analytic $F$-gauges in terms of coefficients for Hyodo--Kato cohomology or proétale cohomology, respectively. 

Namely, in the Hyodo--Kato case, we construct realisation functors
\begin{equation*}
\begin{tikzcd}[column sep=huge]
& \Perf(X^\HK)\\
\Perf(X^\Syn)\ar[ur, "T_\HK"]\ar[r, "T_{\dR, +}"]\ar[dr, "T_\dR", swap] & \Perf(X^{\dR, +}) \\
& \Perf(X^\dR)
\end{tikzcd}
\end{equation*}
for any Gelfand stack $X$, where we recall that $X^\HK$ is the \emph{Hyodo--Kato stack} of $X$ from \cite[§6]{dRFF}. Then we prove:

\begin{thm}[\cref{thm:hkcomp-main}]
\label{thm:intro-hkcomp}
Let $X$ be a Berkovich smooth derived Berkovich space over $\Q_p$. Then there is an equivalence of categories
\begin{equation*}
\Perf(X^\Syn)\cong \Perf(X^\HK)\times_{\Perf(X^\dR)} \Perf(X^{\dR, +})
\end{equation*}
induced by the realisation functors $T_\HK$ and $T_{\dR, +}$. In particular, for any $E\in\Perf(X^\Syn)$, there is a pullback diagram
\begin{equation*}
\begin{tikzcd}
R\Gamma(X^\Syn, E)\ar[r]\ar[d] & R\Gamma(X^\HK, T_\HK(E))\ar[d] \\
R\Gamma(X^{\dR, +}, T_{\dR, +}(E))\ar[r] & R\Gamma(X^\dR, T_\dR(E))\nospacepunct{\;.}
\end{tikzcd}
\end{equation*}
\end{thm}

\begin{rem}
Specialising the result above to the case where $E=\O\{i\}$ is a Breuil--Kisin twist recovers the more classical statement \cite[Thm.\ 1.1.(4)]{padicComparisons}. In particular, this shows that our analytic syntomic cohomology coincides with Colmez--Nizio{\l} syntomic cohomology if we take coefficients in a Breuil--Kisin twist.
\end{rem}

We proceed similarly in the proétale case: Here, we construct realisation functors
\begin{equation*}
\begin{tikzcd}[column sep=huge]
& \Perf(X^{\Div^1})\\
\Perf(X^\Syn)\ar[ur, "T_{\Div^1}"]\ar[r, "T_{\HT, \dagger, +}"]\ar[dr, "T_{\HT, \dagger}", swap] & \Perf(X^{\HT, \dagger, +}) \\
& \Perf(X^{\HT, \dagger})
\end{tikzcd}
\end{equation*}
for any Gelfand stack $X$, where $X^{\Div^1}$ was defined in \cite{AnPrism} and computes proétale cohomology of $X$ as shown in loc.\ cit.\ while $X^{\HT, \dagger, +}$ and $X^{\HT, \dagger}$ are certain closed substacks of $X^\N$ to be defined in §\ref{subsect:loci} below. Then our result reads as follows:

\begin{thm}[\cref{thm:proet-main}]
\label{thm:intro-proet1}
Let $X$ be a Berkovich smooth derived Berkovich space over $\Q_p$. Then there is an equivalence of categories
\begin{equation*}
\Perf(X^\Syn)\cong \Perf(X^{\Div^1})\times_{\Perf(X^{\HT, \dagger})} \Perf(X^{\HT, \dagger, +})
\end{equation*}
induced by the realisation functors $T_{\Div^1}$ and $T_{\HT, \dagger, +}$. In particular, for any $E\in\Perf(X^\Syn)$, there is a pullback diagram
\begin{equation*}
\begin{tikzcd}
R\Gamma(X^\Syn, E)\ar[r]\ar[d] & R\Gamma(X^{\Div^1}, T_{\Div^1}(E))\ar[d] \\
R\Gamma(X^{\HT, \dagger, +}, T_{\HT, \dagger, +}(E))\ar[r] & R\Gamma(X^{\HT, \dagger}, T_{\HT, \dagger}(E))\nospacepunct{\;.}
\end{tikzcd}
\end{equation*}
\end{thm}

The conclusion above does not yet provide an immediate comparison between syntomic cohomology and proétale cohomology. However, analysing the map 
\begin{equation*}
R\Gamma(X^{\HT, \dagger, +}, T_{\HT, \dagger, +}(E))\rightarrow R\Gamma(X^{\HT, \dagger}, T_{\HT, \dagger}(E))
\end{equation*}
in the diagram above, we obtain the following variant (see §\ref{subsect:syn} for the definition of Hodge--Tate weights in this context):

\begin{thm}[\cref{thm:proet-main2}]
\label{thm:intro-proet2}
Let $X$ be a smooth partially proper rigid space over $\Q_p$. If $E\in\Vect(X^\Syn)$ is a vector bundle analytic $F$-gauge with Hodge--Tate weights all at most $-i$ for some $i\geq 0$, then the natural morphism
\begin{equation*}
R\Gamma(X^\Syn, E)\rightarrow R\Gamma(X^{\Div^1}, T_{\Div^1}(E))
\end{equation*}
is an isomorphism on $\tau^{\leq i}$ and induces an injection on $H^{i+1}$. In particular, for $E=\O\{i\}$, we have
\begin{equation*}
\tau^{\leq i} R\Gamma_\Syn(X, \Q_p(i))\cong \tau^{\leq i}R\Gamma_\proet(X, \Q_p(i))\;.
\end{equation*}
\end{thm}

\begin{rem}
Once a robust theory of analytic syntomification with $\Z_p$-coefficients is set up, we expect that the statements as well as the proofs of \cref{thm:intro-proet1} and \cref{thm:intro-proet2} carry over almost verbatim to the integral coefficient case. In contrast, this should \emph{not} be expected for \cref{thm:intro-hkcomp} as this is an inherently rational statement.
\end{rem}

Finally, we can use the results above to explicitly describe the category of vector bundles on $X^\Syn$ for smooth partially proper rigid spaces $X$ over $\Q_p$. In particular, this allows us to relate the analytic $F$-gauge coefficients that are allowed in the comparison theorems above to more classical objects in $p$-adic Hodge theory. We start with the case $X=\GSpec\Q_p$.

\begin{thm}[\cref{thm:bk-mainqp}]
There is an equivalence of categories
\begin{equation*}
\Vect(\Q_p^\Syn)\cong \Vect^\dR(\FF_{\Spd\Q_p})
\end{equation*}
between vector bundles on $\Q_p^\Syn$ and $\Gal_{\Q_p}$-equivariant vector bundles on the Fargues--Fontaine curve which are de Rham in the sense of \cite[Def.\ 15.12]{FarguesFontaine}. Moreover, if $V$ and $W$ are de Rham representations of $\Gal_{\Q_p}$, then
\begin{equation*}
\RHom_{\Q_p^\Syn}(V, W)\cong \RHom_{\Rep^\dR(\Gal_{\Q_p})}(V, W)\;.
\end{equation*}
In particular, we have
\begin{equation*}
H^1(\Q_p^\Syn, V)\cong H^1_g(\Gal_{\Q_p}, V)\;.
\end{equation*}
\end{thm}

In the above theorem, the group $H^1_g$ was defined by Bloch--Kato in \cite{BlochKato} and coincides with the subspace of $H^1$ spanned by extensions of $\Q_p$ by $V$ which are de Rham. This parallels the description of reflexive sheaves on the algebraic syntomification of $\Z_p$ in terms of lattices in crystalline representations of $\Gal_{\Q_p}$ given by Bhatt--Lurie in \cite[Thm.\ 6.6.13]{FGauges} and their recovery of the Bloch--Kato--Selmer group $H^1_f$ in terms of $H^1$ on the algebraic syntomification of $\Z_p$ from \cite[Prop.\ 6.7.3]{FGauges}.

Finally, we are able to extend the above result to arbitrary smooth partially proper rigid spaces $X$ over $\Q_p$. In particular, this relates the category of vector bundle analytic $F$-gauges to the category of de Rham local systems on $X_\proet$, which are a more classical candidate for a category of ``universal coefficients'' for $p$-adic cohomology theories attached to $X$.

\begin{thm}[\cref{thm:drbundles-main}]
\label{thm:intro-drbundles}
Let $X$ be a smooth partially proper rigid space over $\Q_p$. Then there is an equivalence of categories
\begin{equation*}
\Vect(X^\Syn)\cong \Vect^\dR(\FF_{X^\diamond})\;,
\end{equation*}
where the right-hand side denotes the full subcategory of vector bundles on the relative Fargues--Fontaine curve of $X^\diamond$ which are de Rham in the sense of \cref{defi:drbundles-drvb}. In particular, there is a fully faithful embedding
\begin{equation*}
\Loc^\dR(X_\proet, \Q_p)\hookrightarrow \Vect(X^\Syn)
\end{equation*}
of the category of de Rham local systems on $X_\proet$ into the category of vector bundle analytic $F$-gauges on $X$.
\end{thm}

\begin{rem}
There should also be a variant of the above result with integral coefficients. Namely, we expect that vector bundles on the integral analytic syntomification of a partially proper rigid space $X$ over $\Q_p$ are equivalent to lattices in de Rham local systems on $X_\proet$, with the proof being mostly analogous to the one we give for \cref{thm:intro-drbundles}. Indeed, one should think of the passage to integral coefficients as adding a slope zero condition on the Fargues--Fontaine curve, which explains the fact the we only see honest local systems in the integral coefficient case.
\end{rem}

\subsubsection{Syntomic cohomology over $\C_p$}

Finally, we develop analogues of the above results and constructions for rigid-analytic varieties over $\C_p$. For this, the key is to replace $X^\Syn$ by an appropriate base change, which we denote by $X^{\Syn/\C_p}$ for any Gelfand stack $X$ over $\C_p$. As a warning, let us already point out that the resulting syntomic cohomology over $\C_p$ we obtain is \emph{different} from the one defined by Colmez--Nizio{\l}, see \cref{cor:intro-geomhkclassical} below.

As the base change that yields $X^{\Syn/\C_p}$ does not change the overall geometry of $X^\Syn$, the proofs of the comparison theorems above mostly go through unchanged. The main difference is that we now see a new stack appearing, which we call the filtered \emph{$B_\dR^{+, \dagger}$-stack} of $X$ and denote by $X^{\dR, +/B_\dR^{+, \dagger}}$. The main part of the work lies in identifying the base change of the cohomology of the structure sheaf on this stack to $B_\dR^+$ with the $B_\dR^+$-cohomology $R\Gamma(X/B_\dR^+)$ of $X$ from \cite[§13]{IntegralpAdicHT} or \cite{GuoBdR}, which we recall also has a natural filtration.

\begin{thm}[\cref{thm:geom-bdr+coh}]
Let $X$ be a smooth partially proper rigid space over $\C_p$. Then there is a natural isomorphism
\begin{equation*}
\Fil^\bullet R\Gamma(X^{\dR, +/B_\dR^{+, \dagger}},\O)\tensor_{B_\dR^{+, \dagger}} B_\dR^+\cong \Fil^\bullet R\Gamma(X/B_\dR^+)\;,
\end{equation*}
where the left-hand side denotes the base change of the filtered object
\begin{equation*}
\dots\rightarrow R\Gamma(X^{\dR, +/B_\dR^{+, \dagger}},\O\langle -(i+1)\rangle)\rightarrow R\Gamma(X^{\dR, +/B_\dR^{+, \dagger}},\O\langle -i\rangle)\rightarrow R\Gamma(X^{\dR, +/B_\dR^{+, \dagger}},\O\langle -(i-1)\rangle)\rightarrow\dots
\end{equation*}
from the overconvergent de Rham period ring $B_\dR^{+, \dagger}$ to $B_\dR^+$.
\end{thm}

\begin{rem}
Using this stacky approach to $B_\dR^+$-cohomology, we are also able to give straightforward proofs of the basic theorems relating $B_\dR^+$-cohomology and de Rham cohomology in \cref{cor:geom-bdr+props}. Namely, for any smooth partially proper rigid space over $\C_p$, we have
\begin{equation*}
R\Gamma(X/B_\dR^+)\tensor_{B_\dR^+} \C_p\cong R\Gamma_\dR(X/\C_p)
\end{equation*}
and, if $X$ is in addition qcqs and admits a model $\ol{X}$ over some finite extension $K$ of $\Q_p$, then
\begin{equation*}
R\Gamma(X/B_\dR^+)\cong R\Gamma_\dR(\ol{X}/K)\tensor_K B_\dR^+\;. \qedhere
\end{equation*}
\end{rem}

Then the first comparison theorem we can recover is the de Rham comparison theorem over $\C_p$ from \cite[Thm.\ 13.1]{IntegralpAdicHT}, which takes the following shape:

\begin{thm}[\cref{thm:geom-drcomp}]
Let $X$ be a smooth proper rigid space over $\C_p$. Then there is a natural isomorphism
\begin{equation*}
R\Gamma_\proet(X, \Q_p)\tensor_{\Q_p} B_\dR\cong R\Gamma(X/B_\dR^+)\tensor_{B_\dR^+} B_\dR
\end{equation*}
compatible with the filtrations.
\end{thm}

We also obtain geometric analogues of the comparison theorems between syntomic cohomology and Hyodo--Kato cohomology or proétale cohomology, respectively:

\begin{thm}[\cref{thm:geom-hkcomp}]
Let $X$ be a Berkovich smooth derived Berkovich space over $\C_p$. Then there is an equivalence of categories
\begin{equation*}
\Perf(X^{\Syn/\C_p})\cong \Perf(X^{\HK/\C_p})\times_{\Perf(X^{\dR/B_\dR^{+, \dagger}})} \Perf(X^{\dR, +/B_\dR^{+, \dagger}})\;.
\end{equation*}
In particular, for any $E\in\Perf(X^{\Syn/\C_p})$, there is a pullback diagram
\begin{equation*}
\begin{tikzcd}
R\Gamma(X^{\Syn/\C_p}, E)\ar[r]\ar[d] & R\Gamma(X^{\HK/\C_p}, T_\HK(E))\ar[d] \\
R\Gamma(X^{\dR, +/B_\dR^{+, \dagger}}, T_{\dR, +}(E))\ar[r] & R\Gamma(X^{\dR/B_\dR^{+, \dagger}}, T_\dR(E))\nospacepunct{\;.}
\end{tikzcd}
\end{equation*}
\end{thm}

\begin{cor}[\cref{cor:geom-hkcompclassical}]
\label{cor:intro-geomhkclassical}
Let $X$ be a smooth partially proper qcqs rigid space over $\C_p$. For any $i\in\Z$, there is a cartesian diagram
\begin{equation*}
\begin{tikzcd}
R\Gamma_\Syn(X/\C_p, \Q_p(i))\ar[r]\ar[d] & (R\Gamma_\HK(X)\tensor_{\Q_p^\un} B_{\log})^{\phi=p^i, N=0}\ar[d] \\
\Fil^i R\Gamma(X/B_\dR^+)\ar[r] & R\Gamma(X/B_\dR^+)\nospacepunct{\;,}
\end{tikzcd}
\end{equation*}
where $R\Gamma_\Syn(X/\C_p, \Q_p(i))$ denotes the cohomology of $\O\{i\}$ on $X^{\Syn/\C_p}$ and $B_{\log}$ is the period ring from \cite[Def.\ 10.3.1]{FarguesFontaine}.
\end{cor}

Comparing the above corollary to \cite[Thm.\ 1.1.(4)]{padicComparisons}, we see that, over $\C_p$, our analytic syntomic cohomology does \emph{not} coincide with Colmez--Nizio{\l} syntomic cohomology due to the appearance of the period ring $B_{\log}$ in place of $B_\st^+$. Instead, by \cite[Thm.\ 6.3]{BoscopAdic}, we recover the \emph{syntomic Fargues--Fontaine cohomology} of Bosco.

\begin{thm}[\cref{thm:geom-proet1}]
Let $X$ be a Berkovich smooth derived Berkovich space over $\C_p$. Then there is an equivalence of categories
\begin{equation*}
\Perf(X^{\Syn/\C_p})\cong \Perf(X^{\Div^1/\C_p})\times_{\Perf(X^{\HT, \dagger/\C_p})} \Perf(X^{\HT, \dagger, +/\C_p})\;.
\end{equation*}
In particular, for any $E\in\Perf(X^{\Syn/\C_p})$, there is a pullback diagram
\begin{equation*}
\begin{tikzcd}
R\Gamma(X^{\Syn/\C_p}, E)\ar[r]\ar[d] & R\Gamma(X^{\Div^1/\C_p}, T_{\Div^1}(E))\ar[d] \\
R\Gamma(X^{\HT, \dagger, +/\C_p}, T_{\HT, \dagger, +}(E))\ar[r] & R\Gamma(X^{\HT, \dagger/\C_p}, T_{\HT, \dagger}(E))\nospacepunct{\;.}
\end{tikzcd}
\end{equation*}
\end{thm}

Again, an analysis of the map
\begin{equation*}
R\Gamma(X^{\HT, \dagger, +/\C_p}, T_{\HT, \dagger, +}(E))\rightarrow R\Gamma(X^{\HT, \dagger/\C_p}, T_{\HT, \dagger}(E))
\end{equation*}
enables us to give a more direct comparison between syntomic and proétale cohomology.

\begin{thm}[\cref{thm:geom-proet2}]
Let $X$ be a smooth partially proper rigid space over $\C_p$. If $E\in\Vect(X^{\Syn/\C_p})$ is a vector bundle analytic $F$-gauge on $X$ relative to $\C_p$ with Hodge--Tate weights all at most $-i$ for some $i\geq 0$, then the natural morphism
\begin{equation*}
R\Gamma(X^{\Syn/\C_p}, E)\rightarrow R\Gamma(X^{\Div^1/\C_p}, T_{\Div^1}(E))
\end{equation*}
is an isomorphism on $\tau^{\leq i}$ and induces an injection on $H^{i+1}$. In particular, for $E=\O\{i\}$, we obtain
\begin{equation*}
\tau^{\leq i} R\Gamma_\Syn(X/\C_p, \Q_p(i))\cong \tau^{\leq i}R\Gamma_\proet(X, \Q_p(i))\;.
\end{equation*}
\end{thm}

\subsection{Organisation of the paper}

In §2, we mostly review some aspects of the theory of Gelfand stacks and their quasicoherent sheaves. The main novelty here is a spreading out property for perfect complexes on a certain class of Gelfand stacks which we call ``nicely coverable''. This will be proved in §2.3 and is one of the main ingredients in our derivation of comparison theorems from the geometry of $X^\Syn$ later on. We also establish a version of the Rees equivalence in the analytic setting, see §2.4.

We then move on to reviewing the theory of rational analytic prismatisation as developed by Anschütz--Le Bras--Rodríguez Camargo--Scholze in §3. While everything we discuss there is due to them, we have still chosen to be as detailed as possible since a written account of the theory was not yet available at the time of writing.

With these preliminaries out of the way, we can define the analytic Nygaardification and the analytic syntomification in §4. In §4.2, we study the basic geometry of the analytic Nygaardification and how it relates to other stacks attached to $p$-adic cohomology theories. Finally, we establish two technical properties of the Nygaardification: First, we prove in §4.4 that the functor $X\mapsto X^\N$ on derived Berkovich spaces is compatible with étale localisation, which will be the key to almost all of our dévissage arguments in the rest of the paper. Second, in §4.5, we show that, after restricting to a certain locus, the stack $\Q_p^\N$ is nicely coverable in the sense of §2.3. This will be crucial for us to later be able to apply the spreading out property for perfect complexes proved earlier. 

In §5, we give explicit presentations for some Nygaardifications as quotient stacks. Namely, we start with the case $X=\GSpec\Q_p$ in §5.1 and then move on to the case $X=\G_m$ in §5.2. These presentations will be key to our proofs of the cohomological smoothness and properness results for the syntomification. Then we study a certain closed substack of $X^\Syn$ that we call $X^{\HT, \dagger, +}$ in more detail in §5.3 by not only giving an explicit presentation as a quotient, but also describing its category of perfect complexes and their cohomology in the case where $X$ admits a Berkovich étale map to the $n$-dimensional overconvergent unit torus $\ol{\T}^n$. This will be the crucial ingredient for both the derivation of the comparison between syntomic and proétale cohomology in low degrees and the description of vector bundles on $X^\Syn$ in terms of de Rham vector bundles on $\FF_{X^\diamond}$ later on.

Using the results of §5, we will be able to establish our results concerning the properties of the syntomification with respect to six functors in §6. Here we first prove cohomological smoothness of $\Q_p^\Syn$ in §6.1 and then move on to the relative case $X^\Syn\rightarrow Y^\Syn$ in §6.2. Finally, §6.3 is devoted to establishing the existence of a strong theory of first Chern classes for syntomic cohomology of derived Berkovich spaces in the sense of \cite[§5]{Zavyalov}.

In §7, we start proving comparison theorems. Namely, we discuss how to interpret Scholze's results from \cite{PAdicHodgeTheory} in terms of the syntomification by first giving an alternative construction of his functor from filtered vector bundles with connection to $\mathbb{B}_\dR^+$-local systems in §7.1 and then reproving the de Rham comparison theorem in §7.2.

We then move on to the comparison between syntomic cohomology and Hyodo--Kato cohomology in §8. We first discuss how to describe the category of perfect complexes on $X^\Syn$ in terms of perfect complexes on the Hyodo--Kato stack and on the filtered de Rham stack of $X$ in §8.1. Then we use this to relate vector bundles on $\Q_p^\Syn$ to $\Gal_{\Q_p}$-equivariant bundles on the Fargues--Fontaine curve which are de Rham and to recover Bloch--Kato's $H^1_g$ via the syntomification in §8.2.

Subsequently, we discuss the comparison between syntomic and proétale cohomology in §9. The method here will be closely analogous to the one from §8 and hence we content ourselves with outlining the steps of the proof in §9.1. Instead, the focus of this section will be §9.2, where we put this together with our results from §5.3 in order to prove the description of vector bundles on $X^\Syn$ in terms of de Rham bundles on $\FF_{X^\diamond}$.

Finally, we discuss syntomic cohomology over $\C_p$ in §10, for which we will first introduce the relevant base changes in §10.1. The main part of this section is the stacky approach to $B_\dR^+$-cohomology, which we develop in §10.2. Afterwards, the proofs of the comparison theorems go through just as in the arithmetic case and hence we only give outlines in §10.3.

\subsection*{Notations and conventions} 

Throughout, we use the solid analytic ring structure on $\Z$ and, in particular, all tensor products as well as all abelian groups are implicitly understood to be solid. Moreover, all functors, e.g.\ tensor products or pushforwards, will always implicitly be derived and hence, in particular, $\tensor$ always denotes a solid derived tensor product.

All Gelfand $\Q_p$-algebras will be assumed to be separable and qfd in the sense of \cite[Def.\ 4.2.6]{dRFF} and \cite[Def.\ 4.2.14]{dRFF}. In particular, all our Gelfand stacks will implicitly be qfd.

For any Gelfand $\Q_p$-algebra $A$ and $r\geq 0$, we write
\begin{equation*}
A\langle T\rangle_{\leq r}\coloneqq A\tensor_{\Q_p} \Q_p\langle T\rangle_{\leq r}\;, \hspace{0.3cm} A\{T\}^\dagger\coloneqq A\tensor_{\Q_p} \Q_p\langle T\rangle_{\leq 0}\;,
\end{equation*}
where $\Q_p\langle T\rangle_{\leq 0}$ is the ring of overconvergent functions on the closed disk of radius $r$ over $\Q_p$ from \cite[Def.\ 2.2.3]{dRFF}. Moreover, we will write $\ol{A}\coloneqq A^{\dagger\text{-red}}$ to denote $\dagger$-reductions in the sense of \cite[Def.\ 2.2.19]{dRFF}. 

The term ``rational localisation'' will always refer to a rational localisation on the Berkovich spectrum and hence rational localisations are always overconvergent; e.g., we have
\begin{equation*}
(\{|f|\leq 1\}\subseteq\GSpec A)\coloneqq \GSpec A\langle f\rangle_{\leq 1}
\end{equation*}
for any Gelfand ring $A$ and any $f\in A$. We warn the reader that $A\langle f\rangle_{\leq 1}\neq A\langle f\rangle$, i.e.\ this usage of the term ``rational localisation'' differs from the terminology commonly used in the context of adic spaces.

Throughout, there will be a number of different twists that play a role: The Tate twist is always denoted by $\Z_p(1)$ and the Breuil--Kisin twist (on $\Q_p^\prism, \Q_p^\N$ or $\Q_p^\Syn$) is denoted by $\O\{1\}$. Finally, we have chosen the slightly nonstandard notation $\O\langle -1\rangle$ for the tautological normed line bundle on the classifying stack $*/\ol{\T}$ of the overconvergent unit torus.

\subsection*{Acknowledgements} 

First and foremost I would like to thank my advisor Peter Scholze for lots of helpful discussions, suggestions and guidance throughout this project. I also thank Johannes Anschütz, Guido Bosco, Arthur--César Le Bras and Juan Esteban Rodríguez Camargo for useful conversations during the development of this paper. It will not be hard to recognise the tremendous intellectual debt this paper owes to \cite{dRFF} and \cite{AnPrism}. I also thank Arthur--César Le Bras for inviting me to present some of these results at the University of Strasbourg as well as Akhil Mathew for an invitation to the University of Chicago. Finally, I thank Phil Pützstück for suggesting the argument for \cref{lem:perf-betti} and Ko Aoki for suggesting the argument for \cref{lem:perf-finflatdim}. I also thank Saverio Caleca for useful comments on a draft. This paper was written during my time as a PhD student at the Max Planck Institute for Mathematics in Bonn and I wish to thank the institute for its hospitality.

\newpage

\section{Preliminaries on Gelfand stacks}
\label{sect:recall}

We start by recording some facts about Gelfand rings and Gelfand stacks over $\Q_p$. With the exception of §\ref{subsect:spreading}, most of the material presented should be considered standard and, in particular, all the definitions are taken from \cite{dRFF}.

\subsection{Quasicoherent sheaves and perfect complexes on Gelfand stacks}

Recall the definition of quasicoherent sheaves on Gelfand stacks.

\begin{defi}
Let $X$ be a Gelfand stack. Its (derived) category of \emph{quasicoherent sheaves} is defined as
\begin{equation*}
\D(X)\coloneqq \lim_{\GSpec A\rightarrow X} \D(A)\;,
\end{equation*}
where the transition maps are given by pullback and $\D(A)$ denotes the category of $A$-modules in $\D_\solid(\Q_p)$ for any Gelfand $\Q_p$-algebra $A$.
\end{defi}

We also recall that the functor
\begin{equation*}
\GelfStk\rightarrow \PrL\;, \hspace{0.3cm} X\mapsto \D(X)
\end{equation*}
promotes to a full six functor formalism in which all morphisms between affine Gelfand stacks
\begin{equation*}
f: \GSpec B=Y\rightarrow X=\GSpec A
\end{equation*}
are $!$-able and weakly cohomologically proper, i.e.\ there is an equivalence $f_!\cong f_*$. Nevertheless, throughout most of this paper, we will mostly be considering just perfect complexes instead of all quasicoherent sheaves. For this, we recall that this is a well-behaved notion in Gelfand stacks.

\begin{defi}
Let $X$ be a Gelfand stack. The category $\Perf(X)$ of \emph{perfect complexes} on $X$ is defined as the full subcategory of $\D(X)$ consisting of dualisable objects.
\end{defi}

The association $X\mapsto \Perf(X)$ automatically satisfies descent, but the above definition has the downside that dualisability in $\D(A)$ is a priori hard to control for arbitrary analytic rings $A$. This is not the case for Gelfand rings:

\begin{prop}
\label{prop:recall-fredholm}
For any Gelfand ring $A$, we have $\Perf(A)\cong \Perf(A(*))$, where the right-hand side denotes the usual category of perfect complexes over the discrete animated ring $A(*)$. Moreover, the association
\begin{equation*}
A\mapsto \Perf^{[a, b]}(A)\;,
\end{equation*}
where the right-hand side denotes the full subcategory of perfect complexes over $A(*)$ with Tor-amplitude in $[a, b]$ extends uniquely to a sheaf of categories on $\GelfStk$.
\end{prop}
\begin{proof}
See \cite[Prop.\ 3.3.5]{dRFF} and \cite[Lem.\ 4.2.11]{dRFF}.
\end{proof}

In particular, we also obtain a meaningful notion of vector bundles on Gelfand stacks:

\begin{defi}
For any Gelfand stack $X$ and integers $a\leq b$, we define
\begin{equation*}
\Perf^{[a, b]}(X)\coloneqq \lim_{\GSpec A\rightarrow X} \Perf^{[a, b]}(A)\;,
\end{equation*}
where the limit is taken along pullback functors. The category $\Vect(X)$ of \emph{vector bundles} on $X$ is defined by $\Vect(X)\coloneqq \Perf^{[0, 0]}(X)$.
\end{defi}

\comment{
\begin{defi}
Let $A$ be a solid animated $\Q_p$-algebra. The subring $A^{\leq 1}$ of \emph{elements of norm at most $1$} is defined by
\begin{equation*}
A^{\leq 1}(S)=\Map_{\Q_p\text{-Alg}}(\Q_p\langle \mathbb{N}[S]\rangle_{\leq 1}, A)\;,
\end{equation*}
for any light profinite set $S$, where $\Q_p\langle\mathbb{N}[S]\rangle_{\leq 1}=\colim_{\epsilon>0} \Q_p\langle \mathbb{N}[p^\epsilon S]\rangle$ is the algebra of functions on the overconvergent $S$-dimensional unit disk.
\end{defi}

Recall that, since $\Q_p\langle\mathbb{N}[S]\rangle_{\leq 1}$ is an idempotent algebra over $\Q_p[\mathbb{N}[S]]$, the above indeed always defines a full condensed subanima of $A$, see \cite[Rem.\ 2.2.12.(2)]{dRFF}.

\begin{defi}
Let $A$ be a solid animated $\Q_p$-algebra. It is called \dots
\begin{enumerate}[label=(\roman*)]
\item \dots \emph{bounded} if the natural map $A^{\leq 1}[\tfrac{1}{p}]\rightarrow A$ is an isomorphism.
\item \dots \emph{Gelfand} if it is bounded and $A/A^{\leq 1}$ is a discrete condensed set.
\end{enumerate}
The category of Gelfand rings is denoted by $\GelfRing$.
\end{defi}

\begin{rem}
Recall that \cite[Def.\ 4.2.6.(1)]{dRFF} and \cite[Def.\ 4.2.14.(1)]{dRFF} define what it means for a Gelfand $\Q_p$-algebra to be \emph{separable} and \emph{qfd}. Note that we will always make these finiteness assumptions, but they will hardly play an explicit role in the sequel. 
\end{rem}

Any Gelfand $\Q_p$-algebra $A$ can be seen as an analytic ring by equipping it with the induced analytic ring structure from $\Q_{p, \solid}$, i.e.\ by considering the pair $(A, \D(A))$, where
\begin{equation*}
\D(A)\coloneqq \Mod_A(\D_\solid(\Q_p))\;.
\end{equation*}
For any map $f: A\rightarrow B$, the restriction of scalars functor
\begin{equation*}
f_*=f_!: \D(B)\rightarrow \D(A)
\end{equation*}
then has both adjoints as it commutes with both limits and colimits and we denote the left-adjoint by $f^*$ while the right-adjoint is denoted $f^!$. We can extend the functor $A\mapsto \D(A)$ from Gelfand rings to all anima-valued presheaves on Gelfand rings by left Kan extension using $!$-pullback functoriality and will denote this extension by $\D^!(X)$.

\begin{defi}
A map $f: X\rightarrow Y$ in $\operatorname{PSh}(\GelfRing, \Ani)$ is called a \emph{$!$-equivalence} if $!$-pullback along any pullback of an iterated diagonal of $f$ induces an isomorphism on $\D^!$. The category $\GelfStk$ of (qfd) \emph{Gelfand stacks} over $\Q_p$ is obtained by localising $\operatorname{PSh}(\GelfRing, \Ani)$ at countably presented $!$-equivalences.
\end{defi}

Most important for us will be the fact that, up to size issues, $\GelfStk$ is an $\infty$-topos which sits strictly between sheaves and hypersheaves for the \emph{$!$-topology}, see \cite[Rem.\ 4.2.2]{dRFF}. This is the Grothendieck topology on $\GelfRing^\op$ generated by finite disjoint unions and maps $f: \GSpec B=Y\rightarrow X=\GSpec A$ satisfying \emph{$!$-descent}, which means that the induced functor
\begin{equation*}
\D(X)\xrightarrow{f^!} {\lim_\Delta}^!\, \D(Y^{\times_X n})
\end{equation*}
is an equivalence, where the limit is taken along $!$-pullback functors and over the simplex category $\Delta$. Moreover, in our setting, whether a map $A\rightarrow B$ of Gelfand rings satisfies $!$-descent is actually easily detectable: as we are only using induced analytic ring structures, this is the case if and only if the map $A\rightarrow B$ is \emph{descendable} in the sense of Mathew, see \cite[Def.\ 3.18]{MathewDescendable}. Overall, we will mostly only need to use that Gelfand stacks are sheaves for the descendable topology on Gelfand rings.

We remark that one can show that the functor sending a map $f: A\rightarrow B$ of Gelfand rings to the usual pullback $f^*: \D(A)\rightarrow \D(B)$ is a sheaf for the Grothendieck topology on $\GelfRing^\op$ generated by finite disjoint unions and maps which are effective epimorphisms in Gelfand stacks, see \cite[Rem.\ 4.2.3]{dRFF}. Thus, it makes sense to define the following:
}

\subsection{Betti stacks and their perfect complexes}

Recall that there is a functor 
\begin{equation*}
\operatorname{Pro}_\NN(\mathrm{Fin})\rightarrow\GelfRing^\op\;, \hspace{0.3cm} S\mapsto \GSpec C^\sm(S, \Q_p)
\end{equation*}
taking any light profinite set $S$ to the algebra of locally constant $\Q_p$-valued functions on $S$. 

\begin{prop}
\label{prop:recall-bettistack}
The above functor uniquely extends to a colimit-preserving functor
\begin{equation*}
(-)^\Betti: \widehat{\Shv}(\mathrm{Pro}_{\mathbb{N}}(\mathrm{Fin}))\rightarrow \GelfStk
\end{equation*}
from hypersheaves on light profinite sets to Gelfand stacks. If $S$ is a finite-dimensional metrisable compact Hausdorff space (seen as a light condensed set), then there is an equivalence
\begin{equation*}
\D(X\times S^\Betti)\cong \Shv(S, \D(X))
\end{equation*}
for any Gelfand stack $X$, where the right-hand side denotes the category of $\D(X)$-valued sheaves on $S$.
\end{prop}
\begin{proof}
See \cite[Lem.\ 4.2.9]{dRFF} and \cite[Cor.\ II.1.2]{RealLLC}.
\end{proof}

We immediately warn the reader the we will usually \emph{not} write the superscript $(-)^\Betti$ when it is clear from the context; e.g., we might write expressions like $X\times [0, 1]$ for a Gelfand stack $X$ and then it should be clear that $[0, 1]$ is seen as a Gelfand stack via its Betti stack. We next want to characterise perfect complexes on Gelfand stacks of the form $X\times S^\Betti$.

\begin{lem}
\label{lem:perf-betti}
Let $X$ be any Gelfand stack and $S$ a CW complex. Then
\begin{equation*}
\Perf(X\times S)\cong \LocConst(S, \Perf(X))\cong \Fun(\Sing(S), \Perf(X))\;,
\end{equation*}
where the middle term is the full subcategory of $\Shv(S, \D(X))$ consisting of locally constant sheaves valued in $\Perf(X)$ and $\Sing(S)$ denotes the singular complex of $S$, which we view as an $\infty$-groupoid.
\end{lem}
\begin{proof}
The second equivalence is a consequence of the results in \cite[App.\ A]{HA}. Since
\begin{equation*}
\Perf(X\times S)\cong \lim_{\GSpec A\rightarrow X} \Perf(\GSpec A\times S)
\end{equation*}
by descent and similarly for $\Perf(X)$, proving that the two outer categories above are equivalent reduces to the affine case $X=\GSpec A$. Moreover, the terms $\Perf(X\times S)$ and $\Fun(\Sing(S), \Perf(X))$ take colimits in $S$ to limits (in the former case this is because colimits are universal in $\infty$-topoi) and hence we may reduce to the case of a finite CW complex (even to the case of $S$ being an $n$-dimensional unit ball or an $n$-dimensional unit sphere).

Having made these reductions, we are now going to prove the first of the equivalences in the statement, which then implies the second by the previous paragraph. For this, recall from \cref{prop:recall-bettistack} above that there is an equivalence
\begin{equation*}
\D(X\times S)\cong \Shv(S, \D(X))
\end{equation*}
and hence, by definition of perfect complexes as dualisable objects, we have to show that the dualisable objects of $\Shv(S, \D(X))$ are exactly the locally constant sheaves with dualisable stalks. However, noting that dualisability of  $\cal{F}\in\Shv(S, \D(X))$ is local on $S$ because it is equivalent to the canonical map $\sHom(\cal{F}, 1)\tensor\cal{F}\rightarrow\sHom(\cal{F}, \cal{F})$ being an isomorphism, we easily see that $\cal{F}$ is dualisable whenever it is locally constant with dualisable stalks.

For the converse, assume that $\cal{F}$ is dualisable. By passing to stalks, we see that this in particular implies that all stalks of $\cal{F}$ are dualisable, and hence we only need to show that $\cal{F}$ is locally constant. As $S$ is compact Hausdorff, taking global sections commutes with filtered colimits and thus the unit $1\in\Shv(S, \D(X))$ is compact. Since $\Hom(\cal{F}, -)\cong \Hom(1, (-)\tensor\cal{F}^\vee)$, this implies that $\cal{F}$ is compact as well. However, noting that $\D(X)=\D(A)$ is compactly generated and hence dualisable, $\cal{F}$ being compact implies that it is locally constant by \cite[Prop.\ 6.15]{KTheory}, so we are done.
\end{proof}

\begin{cor}
\label{cor:perf-contractible}
Let $X$ be any Gelfand stack and $S$ a CW complex. If $S$ is contractible, then
\begin{equation*}
\Perf(X\times S)\cong \Perf(X)
\end{equation*}
via pullback along the projection $X\times S\rightarrow X$.
\end{cor}

We end our discussion of Betti stacks by recording the following gluing lemma for later use.

\begin{lem}
\label{lem:hkcomp-glueoverinterval}
Let $X\rightarrow (0, \infty)$ be a morphism of Gelfand stacks and denote the preimage of any interval $I\subseteq (0, \infty)$ by $X_I$. For any $r\in (0, \infty)$, the natural map induces an isomorphism
\begin{equation*}
X_{(0, r]}\coprod_{X_{[r, r]}} X_{[r, \infty)}\xrightarrow{\cong} X\;.
\end{equation*}
Similar statements hold for other intervals in place of $(0, \infty)$ and for other decompositions of such intervals.
\end{lem}
\begin{proof}
We start by noting that
\begin{equation*}
{(0, r]}\coprod_{[r, r]} {[r, \infty)}\xrightarrow{\cong} (0, \infty)
\end{equation*}
via the natural map. Since colimits are universal in $\infty$-topoi, i.e.\ they commute with pullbacks, pulling back this isomorphism along $X\rightarrow (0, \infty)$ yields the claim.
\end{proof}

\subsection{Spreading out perfect complexes}
\label{subsect:spreading}

Before we move on, we record a technical result concerning a certain ``spreading out'' property of perfect complexes on Gelfand stacks. We begin by establishing some terminology.

\begin{defi}
A $!$-hypercover $X_\bullet\rightarrow X$ of Gelfand stacks is called \emph{admissible} if the natural map
\begin{equation*}
\D(X)\rightarrow \lim \D(X_\bullet)
\end{equation*}
is an equivalence and the same holds for any base change.
\end{defi}

\begin{ex}
By definition of the $!$-topology, the \v{C}ech nerve of any surjection $Y\rightarrow X$ of Gelfand stacks is an admissible $!$-hypercover.
\end{ex}

\begin{defi}
\label{defi:perf-nicelycoverable}
A Gelfand stack $X$ is called \emph{nicely coverable} if there is an admissible $!$-hypercover
\begin{equation*}
\GSpec A_\bullet\rightarrow X
\end{equation*}
of $X$ by affine Gelfand stacks such that all $A_\bullet$ are $n$-truncated animated rings for some $n\geq 0$. Any such $!$-hypercover is called a \emph{nice} hypercover of $X$.
\end{defi}

Verifying that certain explicit Gelfand stacks of interest are nicely coverable will be the most technical part of our arguments. Thus, it will be useful to establish some permanence properties of the notion of nice coverability. In the following, when we speak about ``flatness'', we will always mean flatness with respect to the solid tensor product.

\begin{lem}
\label{lem:perf-bc}
Consider a morphism $X\rightarrow Y$ of Gelfand stacks and assume that $X$ is nicely coverable. If $Y'\rightarrow Y$ is an affine morphism of bounded flat dimension, then $X'\coloneqq X\times_Y Y'$ is nicely coverable as well.
\end{lem}
\begin{proof}
If $\GSpec A_\bullet\rightarrow X$ is a nice hypercover such that all $A_\bullet$ are $n$-truncated and the flat dimension of $Y'\rightarrow Y$ is bounded by $m\geq 0$, then $\GSpec A_\bullet\times_Y Y'\rightarrow X'$ will be an admissible hypercover where all $\GSpec A_\bullet\times_Y Y'$ are affine and $(n+m)$-truncated.
\end{proof}

The next lemma says that being nicely coverable is basically local for admissible hypercovers up to an additional uniformity assumption.

\begin{lem}
\label{lem:perf-covers}
Let $X_\bullet\rightarrow X$ be an admissible $!$-hypercover of Gelfand stacks. If there is some fixed $n\geq 0$ such that any $X_\bullet$ admits an admissible $!$-hypercover by $n$-truncated affine Gelfand stacks, then $X$ is nicely coverable.
\end{lem}
\begin{proof}
For each $j\geq 0$, let $\GSpec A_\bullet^j\rightarrow X_j$ be admissible hypercovers by affine $n$-truncated Gelfand stacks. Considering
\begin{equation*}
B_\bullet\coloneqq A_0^\bullet\times A_1^{\bullet-1}\times\dots\times A_\bullet^0\;,
\end{equation*}
we see that $\GSpec B_\bullet\rightarrow X$ is an admissible hypercover since $\GSpec B_\bullet=\bigsqcup_{i+j=\bullet} \GSpec A_i^j$ and hence
\begin{equation*}
\lim \D(B_\bullet)\cong\lim_{i, j} \D(A_i^j)\cong\lim_j \D(X_j)\cong \D(X)\;.
\end{equation*}
Clearly, all $B_\bullet$ are $n$-truncated since the $A_i^j$ are.
\end{proof}

As nice coverability is intimately related to flatness questions, it will be convenient to record the following two flatness lemmas before we move on.

\begin{lem}
\label{lem:perf-finflatdim}
For each $r\geq \epsilon>0$, the $\O(\ol{\DD}(r))$-algebra $\O(\ol{\DD}(\epsilon))$ has flat dimension at most $1$ with respect to the solid tensor product, where $\ol{\DD}(r)$ denotes the overconvergent closed disk of radius $r$.
\end{lem}
\begin{proof}
Recall that $\O(\ol{\DD}(r))$ is an idempotent $\Q_p[T]$-algebra and hence it suffices to check that $\O(\ol{\DD}(\epsilon))$ has flat dimension at most $1$ over $\Q_p[T]$. Furthermore, as $\O(\ol{\DD}(\epsilon))$ is a colimit of Tate algebras, we are reduced to checking that $\Q_p\langle T\rangle$ has flat dimension at most $1$ over $\Q_p[T]$. However, recall that $\prod_{\mathbb{N}} \Z_p$ is a flat $\Z_p$-module and hence
\begin{equation*}
\Z_p\langle U\rangle=\colim_{\substack{f: \mathbb{N}\rightarrow\mathbb{N} \\ f(n)\rightarrow\infty}} \prod_{\mathbb{N}} p^{f(n)}\Z_pU^n
\end{equation*}
is also flat over $\Z_p$. This implies that $\Q_p[T]\langle U\rangle$ is flat over $\Q_p[T]$ by base change and thus
\begin{equation*}
\begin{tikzcd}
0\ar[r] & \Q_p[T]\langle U\rangle\ar[r, "T-U"] & \Q_p[T]\langle U\rangle\ar[r] & \Q_p\langle T\rangle\ar[r] & 0
\end{tikzcd}
\end{equation*}
is a flat resolution of $\Q_p\langle T\rangle$ over $\Q_p[T]$, whence the claim. (Note that one could also have deduced the flatness of $\Q_p\langle U\rangle$ over $\Q_p$ from \cref{lem:perf-vspflat} below.)
\end{proof}

\begin{lem}
\label{lem:perf-vspflat}
Let $K$ be a nonarchimedean local field. Then any quasiseparated solid $K$-vector space is flat.
\end{lem}
\begin{proof}
See \cite[Lem.\ A.28]{BoscoDrinfeldSpace}.
\end{proof}

The following lemma shows that nice coverability is preserved under passing to finite étale covers and rational localisations.

\begin{lem}
\label{lem:perf-etclosedcoverable}
Let $X\rightarrow Y$ be a morphism of Gelfand stacks which becomes finite étale after pulling back to any totally disconnected nilperfectoid $\GSpec A\rightarrow Y$ over $Y$. If $Y$ is nicely coverable, then $X$ is nicely coverable as well. Moreover, the same conclusion holds if we replace ``finite étale'' by ``rational localisation''.
\end{lem}
\begin{proof}
Let $Y_\bullet\rightarrow Y$ be an admissible hypercover by $n$-truncated affine Gelfand stacks. We claim that there is a $!$-cover $\GSpec A\rightarrow Y_0$ with $A$ strictly totally disconnected nilperfectoid which is in addition flat. Indeed, by examining the proof of \cite[Prop.\ 6.4.6]{dRFF}, this comes down to checking that the map $\colim_n \Q_p\langle T_1, \dots, T_n\rangle\rightarrow \Q_p^\cycl\langle T_i^{1/p^\infty}: i\geq 1\rangle$ is flat, which we may in turn reduce to checking that the maps $\colim_n \Q_p\langle T_1, \dots, T_n\rangle\rightarrow \Q_p\langle T_i: i\geq 1\rangle$ and $\Q_p\langle T\rangle\rightarrow \Q_p\langle T^{1/p}\rangle$ are flat. For the latter map, this is clear as $\Q_p\langle T^{1/p}\rangle\cong \Q_p\langle T\rangle[T^{1/p}]$ is a free module of rank $p$ over $\Q_p\langle T\rangle$. For the former map, note that it suffices to check flatness of $\Q_p\langle T_1, \dots, T_n\rangle\rightarrow \Q_p\langle T_i: i\geq 1\rangle$ for each $n\geq 1$, but this follows by base change from flatness of $\Q_p\langle T_i: i>n\rangle$ over $\Q_p$, which is in turn a consequence of \cref{lem:perf-vspflat}.

Picking a flat $!$-cover $\GSpec A\rightarrow Y_0$ with $A$ strictly totally disconnected nilperfectoid as above, by \cref{lem:perf-bc}, the base change of the hypercover $Y_\bullet\rightarrow Y$ along the map $\GSpec A\rightarrow Y_0$ yields a nice hypercover of $Y$ whose zeroth term is totally disconnected nilperfectoid. In other words, we may assume that $Y_0$ is totally disconnected nilperfectoid and then the base change $X_0\coloneqq Y_0\times_Y X$ is finite étale over $Y_0$ by assumption. Consequently, $X_\bullet\coloneqq Y_\bullet\times_Y X$ is finite étale over $Y_\bullet$. Thus, each $X_\bullet$ is affine and $n$-truncated as well; moreover, as admissible hypercovers are stable under base change, we conclude that $X_\bullet\rightarrow X$ is a nice hypercover and this finishes the proof in the finite étale case.

In the case where $X\rightarrow Y$ becomes a rational localisation after base change to any strictly totally disconnected nilperfectoid over $Y$, one argues similarly: As the base change $X\times_Y \GSpec A$ is defined by finitely many inequalities of the form $|f|\leq 1$ or $|f|\geq 1$ inside $\GSpec A$, the map $X\times_Y \GSpec A\rightarrow\GSpec A$ has finite flat dimension by \cref{lem:perf-finflatdim} and thus nice hypercovers are preserved by base change along this map by \cref{lem:perf-bc}.
\end{proof}

Finally, the last missing definition in order for us to be able to state the desired spreading out property is the following:

\begin{defi}
Let $X$ be a Gelfand stack. An \emph{overconvergent normed divisor} $Z\subseteq X$ is a closed substack of $X$ arising from a pullback square
\begin{equation*}
\begin{tikzcd}
Z\ar[r]\ar[d] & X\ar[d, "f_Z"] \\
*/\ol{\T}^\dR\ar[r] & (\A^1/\ol{\T})^\dR\nospacepunct{\;,}
\end{tikzcd}
\end{equation*}
where $\ol{\T}=\GSpec \Q_p\langle s^{\pm 1}\rangle_{\leq 1}$ denotes the overconvergent unit torus. In this case, for any $\epsilon>0$, the \emph{$\epsilon$-neighbourhood} of $Z$ in $X$ is the closed substack $Z_\epsilon\subseteq X$ defined by
\begin{equation*}
\begin{tikzcd}
Z_\epsilon\ar[r]\ar[d] & X\ar[d, "f_Z"] \\
(\ol{\DD}(\epsilon)/\ol{\T})^\dR\ar[r] & (\A^1/\ol{\T})^\dR\nospacepunct{\;,}
\end{tikzcd}
\end{equation*}
where $\ol{\DD}(\epsilon)\subseteq \A^1$ denotes the overconvergent disk of radius $\epsilon$.
\end{defi}

\begin{lem}
\label{lem:perf-keylemma}
Let $X$ be a Gelfand stack and assume that $X$ is nicely coverable. If $Z\subseteq X$ is an overconvergent normed divisor, pullback induces an equivalence of categories
\begin{equation*}
\Perf(Z)\cong \colim_{\epsilon>0} \Perf(Z_\epsilon)\;.
\end{equation*}
\end{lem}
\begin{proof}
Write $Z_0\coloneqq Z$ for simplicity, i.e.\ $Z_\epsilon=Z$ if $\epsilon=0$. Let $\GSpec A_\bullet\rightarrow X$ be a nice hypercover and choose $n\geq 0$ such that all $A_\bullet$ are $n$-truncated. Then first observe that $X$ is quasi-compact by assumption and hence the map $X\rightarrow (\A^1/\ol{\T})^\dR$ classifying $Z$ factors through $(\ol{\DD}(r)/\ol{\T})^\dR$ for some $r>0$. Since the maps $\ol{\DD}(\epsilon)\rightarrow\ol{\DD}(r)$ are affine and have flat dimension at most $1$ for all $0\leq\epsilon\leq r$ by \cref{lem:perf-finflatdim}, base changing the given nice hypercover of $X$ to $Z_\epsilon$ for $\epsilon\geq 0$ yields corresponding nice hypercovers
\begin{equation*}
U_\bullet^{(\epsilon)}=\GSpec A_\bullet^{(\epsilon)}\rightarrow Z_\epsilon\;.
\end{equation*}

Now any perfect complex on $Z_\epsilon$ for $\epsilon\geq 0$ comes from a perfect complex over the animated ring $A_0^{(\epsilon)}$ and hence has bounded $\operatorname{Tor}$-amplitude, which reduces us to showing that
\begin{equation*}
\Perf^{[a, b]}(Z)\cong\colim_{\epsilon>0} \Perf^{[a, b]}(Z_\epsilon)
\end{equation*}
via pullback for all $a\leq b$. However, by descent, we have
\begin{equation*}
\Perf^{[a, b]}(Z_\epsilon)\cong\Tot(\Perf^{[a, b]}(A_\bullet^{(\epsilon)}))\cong \lim_{\Delta_{\leq N}} \Perf^{[a, b]}(A_\bullet^{(\epsilon)})
\end{equation*}
for some $N\geq 0$ which only depends on $a, b$ and $n$ since all the $A_\bullet^{(\epsilon)}$ are $n$-truncated; here, $\Delta_{\leq N}$ denotes the full subcategory of the simplex category $\Delta$ spanned by the objects $\{0<1<\dots<i\}$ with $i\leq N$. As filtered colimits commute with finite limits, we are then reduced to proving
\begin{equation*}
\Perf^{[a, b]}(A_i^{(0)})\cong \colim_{\epsilon>0} \Perf^{[a, b]}(A_i^{(\epsilon)})
\end{equation*}
for all $a\leq b$ and $i\geq 0$. However, note that $A_i^{(0)}=\colim_{\epsilon>0} A_i^{(\epsilon)}$ by construction and hence the claim follows from \cite[Tag 0BC7]{Stacks}.
\end{proof}

In applications, we will actually make use of the following easy consequence of \cref{lem:perf-keylemma}:

\begin{cor}
\label{cor:perf-corkeylemma}
Let $X$ be a nicely coverable Gelfand stack and $Z\subseteq X$ an overconvergent normed divisor. Then we have
\begin{equation*}
\Perf(X)\cong \Perf(Z)\times_{\colim_{\epsilon>0} \Perf(Z_\epsilon\setminus Z)} \Perf(X\setminus Z)
\end{equation*}
via pullback. Here, $X\setminus Z$ is defined by the pullback diagram
\begin{equation*}
\begin{tikzcd}
X\setminus Z\ar[r]\ar[d] & X\ar[d, "f_Z"] \\
(\G_m/\ol{\TT})^\dR\ar[r] & (\A^1/\ol{\TT})^\dR
\end{tikzcd}
\end{equation*}
and $Z_\epsilon\setminus Z$ is the base change of $X\setminus Z$ along $Z_\epsilon\subseteq X$.
\end{cor}
\begin{proof}
Since filtered colimits commute with finite limits, the previous lemma yields
\begin{equation*}
\begin{split}
\Perf(Z)\times_{\colim_{\epsilon>0} \Perf(Z_\epsilon\setminus Z)} \Perf(X\setminus Z)&\cong \colim_{\epsilon>0} \Perf(Z_\epsilon)\times_{\colim_{\epsilon>0} \Perf(Z_\epsilon\setminus Z)} \Perf(X\setminus Z) \\
&\cong \colim_{\epsilon>0} (\Perf(Z_\epsilon)\times_{\Perf(Z_\epsilon\setminus Z)} \Perf(X\setminus Z)) \\
&\cong \colim_{\epsilon>0} \Perf(X)=\Perf(X)\;,
\end{split}
\end{equation*}
where the penultimate step is due to the fact that $Z_\epsilon$ and $X\setminus Z$ form a cover of $X$ for any $\epsilon>0$ and intersect in $Z_\epsilon\setminus Z$.
\end{proof}

\subsection{The Rees equivalence}

Recall the classical Rees equivalence: Working in stacks on discrete static $A$-algebras for the fpqc topology, there is an equivalence
\begin{equation*}
\D(\A^{1, \alg}/\G_m^\alg)\cong \Fun((\Z, \geq), \D(A))=\DF(A)
\end{equation*}
between quasicoherent sheaves on $\A^{1, \alg}/\G_m^\alg$, where $\G_m^\alg$ acts on $\A^{1, \alg}$ by multiplication, and the filtered derived category of $A$; here we use the superscript $(-)^\alg$ to distinguish the algebraic affine line $\A^{1, \alg}$ from the analytic affine line $\A^1$ in Gelfand stacks, and similarly for $\G_m$. This equivalence takes a filtered object $\Fil^\bullet M\in\DF(A)$ to the graded $A[t]$-module
\begin{equation*}
\bigoplus_{i\in\Z} \Fil^i F\cdot t^{-i}\;,
\end{equation*}
where $t$ denotes the coordinate on $\A^{1, \alg}$. Even more fundamentally, there is an equivalence
\begin{equation*}
\D(*/\G_m^\alg)\cong \Fun(\Z, \D(A))=\D_{\gr}(A)
\end{equation*}
between quasicoherent sheaves on $*/\G_m^\alg$ and the category of $\Z$-graded objects in $\D(A)$.

Our next objective is to prove a variant of this result in analytic geometry, i.e.\ in the setting of Gelfand stacks. We start by considering just the classifying stack of $\G_m$.

\begin{lem}
\label{lem:recall-reesgm}
Let $A$ be any Gelfand ring. Then there is an equivalence between perfect complexes on $*/\G_m$ over $A$ and the full subcategory of $\D_{\gr}(A)$ consisting of objects $\bigoplus_{i\in\Z} M^i$ satisfying the following two conditions:
\begin{enumerate}[label=(\roman*)]
\item Each $M^i$ is a perfect complex over $A$.
\item We have $M^i=0$ for $i\gg 0$ and $i\ll 0$.
\end{enumerate}
The same is true for $\ol{\T}$ in place of $\G_m$.
\end{lem}
\begin{proof}
We first observe that the second part immediately follows from the first: Indeed, the map $*/\ol{\T}\rightarrow*/\G_m$ is a $(0, \infty)$-torsor by \cref{lem:defis-rhvariant} below and thus pullback induces an equivalence on perfect complexes by \cref{cor:perf-contractible}. Therefore, it suffices to prove the first assertion.

Now note that $\G_m$ is suave and hence \cite[Prop.\ 3.1.27.(1)]{dRStack} yields an identification
\begin{equation*}
\D(\GSpec A/\G_m)\cong \Mod_{f^!f_!1}(\D(A))
\end{equation*}
via $f^!$, where $f: \GSpec A\rightarrow \GSpec A/\G_m$ is the canonical map. Note that $f^!\cong f^*[1]$ as the dualising sheaf of $\G_m$ is $1[1]$ and so we also obtain an equivalence as above if we use $f^*$ instead of $f^!$. By a direct computation, one sees that 
\begin{equation*}
f^!f_!1\cong A\{x^{-1}\}^\dagger\oplus xA\{x\}^\dagger
\end{equation*}
with the multiplicative structure on the right-hand side being coefficient-wise multiplication.

Now given any $M\in\Perf(\GSpec A/\G_m)$, the actions of $x^i$ for $i\in\Z$ yield pairwise orthogonal idempotent endomorphisms of the underlying $A$-module of $M$, which we will also denote by $M$ by abuse of notation. Denoting the retracts of $M$ corresponding to these idempotents by $M^i$, we claim that only finitely many of the $M^i$ are nonzero. Once we know this, the fact that the sum of all $x^i$ acts by the identity on $M$ will imply that $M\cong \bigoplus_{i\in\Z} M^i$, as desired.

To prove the remaining claim, first note that each $M^i$ is a perfect complex as well (being a retract of $M$) and hence equivalent to a perfect complex over the discrete ring $A(*)$ by \cref{prop:recall-fredholm}, and similarly for $M$ itself. Then it suffices to show the statement we claim locally on $\Spec A(*)$: Indeed, the support of each $M^i$ will be closed and $\Spec A(*)$ is quasicompact. By the same argument, we may even reduce to checking the claim on stalks and then further to checking it on residue fields by Nakayama's lemma. Finally, over any residue field $K$, we are in the following situation: we are given a perfect complex $M$ over $K$ and a chain
\begin{equation*}
\dots\rightarrow\bigoplus_{|i|\leq n} M^i\rightarrow \bigoplus_{|i|\leq n+1} M^i\rightarrow\dots
\end{equation*}
of direct summands of $M$. Since the cohomology groups of $M$ must be finite-dimensional and only finitely many of them can be nonzero, this chain must become stationary and hence $M^i=0$ for $|i|\gg 0$.

The above discussion proves that we have a functor from the full subcategory of $\D_{\gr}(A)$ described in the statement to $\D(*/\G_m)$ and that it is essentially surjective. To see full faithfulness, note that the functor is symmetric monoidal and all objects in question are dualisable, hence it suffices to show that $\RHom$ from the unit is preserved, for which we have to check that, in the above notation
\begin{equation*}
\RHom_{A\{x^{-1}\}^\dagger\oplus xA\{x\}^\dagger}\left(A, \bigoplus_i M^i\right)=M^0\;,
\end{equation*}
where $A$ has degree zero. However, note that $A$ is a retract of $A\{x^{-1}\}^\dagger\oplus xA\{x\}^\dagger$ with the retraction given by $x^i\mapsto 0$ for all $i\neq 0$, i.e.\ $A$ is a finite projective module over $A\{x^{-1}\}^\dagger\oplus xA\{x\}^\dagger$, and then the claim follows.
\end{proof}

Before we can move on to the case of $\A^1/\G_m$, let us state the following preparatory lemma.

\begin{lem}
\label{lem:recall-reesreduction}
Let $A$ be any Gelfand ring. Then pullback along the canonical map
\begin{equation*}
\G_a^\dagger/\G_m\rightarrow \A^1/\G_m
\end{equation*}
induces an equivalence on perfect complexes.
\end{lem}
\begin{proof}
Throughout the proof, we will suppress the base $\GSpec A$ from the notation for clarity. To begin, note that there is a commutative diagram
\begin{equation*}
\begin{tikzcd}
\G_a^\dagger/\ol{\TT}\ar[d]\ar[r] & \A^1/\ol{\TT}\ar[d] \\
\G_a^\dagger/\G_m\ar[r] & \A^1/\G_m
\end{tikzcd}
\end{equation*}
and we will show that pullback along the bottom map induces an equivalence on perfect complexes by proving the same for all the other maps in the diagram. Indeed, for the vertical maps, this just follows from the fact that they are torsors for $\G_m/\ol{\TT}\cong (0, \infty)$ using \cref{cor:perf-contractible}. 

For the top horizontal map, we first observe that 
\begin{equation*}
\A^1/\ol{\TT}\cong \G_m/\ol{\TT}\coprod_{\ol{\DD}^\times/\ol{\TT}} \ol{\DD}/\ol{\TT}
\end{equation*}
by a variant of \cref{lem:hkcomp-glueoverinterval}, where $\ol{\DD}$ denotes the overconvergent unit disk and $\ol{\DD}^\times$ is its punctured analogue. As $\ol{\DD}^\times/\ol{\TT}\cong (0, 1]$ and $\G_m/\ol{\TT}\cong (0, \infty)$, another application of \cref{cor:perf-contractible} shows that pullback along the canonical map
\begin{equation*}
\ol{\DD}/\ol{\TT}\rightarrow \A^1/\ol{\TT}
\end{equation*}
induces an equivalence on perfect complexes. Finally, to see that pullback along $\G_a^\dagger/\ol{\TT}\rightarrow \ol{\DD}/\ol{\TT}$ induces an equivalence on perfect complexes, we use \cref{cor:perf-corkeylemma}: The stack $\ol{\DD}/\ol{\TT}$ is nicely coverable (use e.g.\ the \v{C}ech nerve of the surjection $\ol{\DD}\rightarrow\ol{\DD}/\ol{\TT}$) and applying loc.\ cit.\ to the overconvergent normed divisor $\{|t|=0\}$, where $t$ denotes the coordinate on $\ol{\DD}$, we obtain
\begin{equation*}
\Perf(\ol{\DD}/\ol{\TT})\cong \Perf(\G_a^\dagger/\ol{\TT})\times_{\colim_{\epsilon>0} \Perf((0, \epsilon])} \Perf((0, 1])\cong \Perf(\G_a^\dagger/\ol{\TT})\;,
\end{equation*}
where the last equivalence is due to another application of \cref{cor:perf-contractible}.
\end{proof}

The following statement was used in the proof:

\begin{lem}
\label{lem:defis-rhvariant}
The norm map induces an isomorphism $\G_m/\ol{\TT}\cong (0, \infty)$.
\end{lem}
\begin{proof}
The proof is along the lines of the proofs of \cite[Prop.\ II.1.4, Thm.\ II.3.1]{RealLLC}. First, recall from \cite[Prop.\ II.1.3]{RealLLC} that a map $\G_m\rightarrow (0, \infty)$ is equivalent to a collection $\{A_Z\}_Z$ of idempotent algebras in $\D(\G_m)$ which become simultaneously connective after pullback along some $!$-cover of $\G_m$, where the index $Z$ runs over closed subsets of $(0, \infty)$ and we demand that the association $Z\mapsto A_Z$ sends limits to colimits and finite unions to limits. In our case, it suffices for $Z$ to run over closed intervals $[a, b]$ and then the map in question is induced by sending $[a, b]$ to the algebra of functions on the overconvergent torus $\ol{\TT}(a, b)$ of inner radius $a$ and outer radius $b$. Note that the map $\G_m\rightarrow (0, \infty)$ indeed factors through the quotient of the source by $\ol{\TT}$.

We now first check that, for all $a\leq b$, the map $\ol{\TT}(a, b)\rightarrow [a, b]$ obtained in this way is surjective; since the closed intervals $[a, b]$ jointly cover $(0, \infty)$, this will imply surjectivity of $\G_m/\ol{\TT}\rightarrow (0, \infty)$. For this, note that there is a surjection $S\rightarrow [a, b]$ from the Cantor set $S=\lim_n \{0, 1\}^n$ onto $[a, b]$ such that the image of each $s_n\times_{S_n} S$ in $[a, b]$ for $s_n\in S_n=\{0, 1\}^n$ is given by a closed interval of length $2^{-n}(b-a)$ and it will be enough to show that the base change of $\ol{\TT}(a, b)\rightarrow [a, b]$ along $S\rightarrow [a, b]$ is surjective. To this end, writing 
\begin{equation*}
X_n\coloneqq \ol{\TT}(a, b)\times_{[a, b]} \bigsqcup_{s_n\in S_n} (s_n\times_{S_n} S)\;,
\end{equation*}
we see that $X_n$ is a disjoint union of $2^n$ overconvergent tori and 
\begin{equation*}
X\coloneqq \ol{\TT}(a, b)\times_{[a, b]} S
\end{equation*}
is affine with coordinate ring $\colim_n \O(X_n)$. Thus, we are done once we show that
\begin{equation*}
\colim_n \Q_p^{S_n}=C^\sm(S, \Q_p)\rightarrow \colim_n \O(X_n)
\end{equation*}
is descendable. However, by \cite[Prop.\ 2.7.2]{MannThesis}, this follows from the fact that each map $\Q_p^{S_n}\rightarrow \O(X_n)$ has a section as modules and is thus descendable of index $1$.

Having established surjectivity of $\G_m\rightarrow (0, \infty)$, it suffices to check that the fibres of this map are $\ol{\TT}$-torsors. However, note that the map in question is even a map of group stacks and, by definition, the preimage of $1\in (0, \infty)$ is given by $\ol{\TT}$, which yields the claim.
\end{proof}

Finally, we can prove the analytic version of the Rees equivalence we will need.

\begin{prop}
\label{prop:recall-rees}
Let $A$ be any Gelfand ring. Then there is an equivalence between perfect complexes on $\A^1/\G_m$ over $A$ and the full subcategory of $\DF(A)$ consisting of objects $\Fil^\bullet M$ satisfying the following two conditions:
\begin{enumerate}[label=(\roman*)]
\item Each $\Fil^i M$ is a perfect complex over $A$.
\item We have $\Fil^i M=0$ for $i\gg 0$ and $\Fil^i M\rightarrow \Fil^{i-1} M$ is an isomorphism for $i\ll 0$.
\end{enumerate}
\end{prop}
\begin{proof}
By \cref{lem:recall-reesreduction}, it suffices to prove the analogous characterisation for perfect complexes on $\G_a^\dagger/\G_m$ instead. As in the proof of \cref{lem:recall-reesgm}, the fact that $\G_m$ is suave with dualising sheaf $1[1]$ shows that
\begin{equation}
\label{eq:recall-rees}
\D(\A^1/\G_m\times\GSpec A)\cong \D(A\{x^{-1}, t\}^\dagger\oplus xA\{x, t\}^\dagger)
\end{equation}
via $f^*$, where the multiplication on the $A$-algebra $A\{x^{-1}, t\}^\dagger\oplus xA\{x, t\}^\dagger$ is noncommutative and governed by the relations 
\begin{equation*}
tx^i=x^{i+1}t\;, \hspace{0.3cm} x^ix^j=\delta_{ij}x^i\;, \hspace{0.3cm} t^it^j=t^{i+j}
\end{equation*}
for all $i, j\in\Z$, where $\delta_{ij}$ is the Kronecker delta. Note that the derived category on the right-hand side of (\ref{eq:recall-rees}) refers to the derived category of left modules.

Now given any $M\in\Perf(\A^1/\G_m\times\GSpec A)$, the actions of $x^{-i}$ for $i\in\Z$ again yield pairwise orthogonal idempotent endomorphisms of the underlying $A$-module of $M$, which we will also denote by $M$ by abuse of notation. Denoting the retracts corresponding to these idempotents by $\Fil^i M$, note that the action of $t$ induces maps $\Fil^i M\rightarrow \Fil^{i-1} M$. We claim the following:
\begin{enumerate}[label=(\roman*)]
\item The pushforward of $M$ to $\GSpec A$ is given by $\Fil^0 M$.
\item The maps $\Fil^i M\rightarrow \Fil^{i-1} M$ are isomorphisms for $i\ll 0$.
\item We have $\Fil^i M=0$ for $i\gg 0$.
\item There is an isomorphism $(\bigoplus_{i\in\Z} \Fil^i M\cdot t^{-i})\tensor_{A[t]} A\{t\}^\dagger\cong M$.
\end{enumerate} 
If we can prove this, we are done: Any filtered $A$-module satisfying the assumptions of the proposition will yield a perfect complex on $\G_a^\dagger/\G_m\times \GSpec A$ via the formula (iv) and then we know that this functor is essentially surjective by (iv) and fully faithful by (i).

To prove (i), note that the unit on $\G_a^\dagger/\G_m\times\GSpec A$ corresponds to $A\{t\}^\dagger$ with the left $A\{x^{-1}, t\}^\dagger\oplus xA\{x, t\}^\dagger$-module structure given by
\begin{equation*}
x^i\cdot t^j=\delta_{ij}t^j\;, \hspace{0.3cm} t^i\cdot t^j=t^{i+j}
\end{equation*}
and our task is to compute the $\RHom$ from $A\{t\}^\dagger$ to $M$. However, note that $A\{t\}^\dagger$ is actually a retract of $A\{x^{-1}, t\}^\dagger\oplus xA\{x, t\}^\dagger$ and thus finite projective: Indeed, the inclusion and retraction maps are given by 
\begin{alignat*}{2}
A\{t\}^\dagger&\rightarrow A\{x^{-1}, t\}^\dagger\oplus xA\{x, t\}^\dagger\;, \hspace{0.5cm}  A\{x^{-1}, t\}^\dagger\oplus xA\{x, t\}^\dagger &&\rightarrow A\{t\}^\dagger\\
t^i&\mapsto x^i t^i &&\mathllap{x^i t^j}\mapsto t^i\;.
\end{alignat*}
Now the claim follows upon noticing that the image of $1\in A\{t\}^\dagger$ already determines the map $A\{t\}^\dagger\rightarrow M$ by linearity and that $1$ is fixed by $x^0$, hence has to land in $\Fil^0 M$.

For (ii), we note that pulling back $M$ along $\GSpec A/\G_m\rightarrow \G_a^\dagger/\G_m\times\GSpec A$ yields a perfect complex and that the $M^i$ occurring in the proof of \cref{lem:recall-reesgm} are precisely the associated gradeds of the filtration $\Fil^i M$ up to a sign. As only finitely many $M^i$ are nonzero, we conclude that $\Fil^i M\rightarrow \Fil^{i-1} M$ is an isomorphism for $|i|\gg 0$.

Finally, we prove (iii) and (iv) at once. By the previous paragraph, we know that there is some $N\geq 0$ such that $\Fil^{i+1} M\rightarrow \Fil^i M$ is an isomorphism for all $i\geq N$. Moreover, by (i), \cref{lem:recall-reesreduction} and \cref{lem:recall-a1gmprim} below, we know that $\Fil^0 M$ is a perfect complex over $A$. As the argument from the previous paragraph also shows that all associated gradeds of the filtration are perfect complexes over $A$, this implies that all $\Fil^i M$ are perfect complexes over $A$ by induction. Now $(\bigoplus_{i\leq N} \Fil^i M\cdot t^{-i})\tensor_{A[t]} A\{t\}^\dagger$ defines a perfect complex on $\G_a^\dagger/\G_m\times\GSpec A$: Indeed, choosing $n\ll 0$ such that $\Fil^i M\rightarrow \Fil^{i-1} M$ is an isomorphism for all $i\leq n$, the complex we claim to be perfect receives a natural map from the evidently perfect complex $\Fil^n M\cdot t^{-n}\tensor_A A\{t\}^\dagger$ and the cofibre of this map is given by $\bigoplus_{n<i\leq N} \Fil^i M\cdot t^{-i}$, which is also perfect since each $\Fil^i M$ is a perfect complex over $A$. Moreover, $M$ receives a natural map
\begin{equation*}
\left(\bigoplus_{i\leq N} \Fil^i M\cdot t^{-i}\right)\tensor_{A[t]} A\{t\}^\dagger\rightarrow M
\end{equation*}
and since both source and target are perfect over $A\{t\}^\dagger$, checking if this is an isomorphism may be done as modules over $A\{t\}^\dagger(*)$, and after modding out the Jacobson radical (by Nakayama's lemma). However, since $t$ is $\dagger$-nilpotent in $A\{t\}^\dagger$, it is in the Jacobson radical of $A\{t\}^\dagger(*)$ and hence we may in particular check the above isomorphism after modding out $t$. In other words, we may check the above isomorphism after pullback to $\GSpec A/\G_m\rightarrow \A^1/\G_m\times\GSpec A$, but there it follows from the proof of \cref{lem:recall-reesgm} above by our choice of $N$. We conclude that $\Fil^i M=0$ for $i>N$, proving (iii), and then the above isomorphism is exactly (iv), so we are done.
\end{proof}

The following statement was used in the proof:

\begin{lem}
\label{lem:recall-a1gmprim}
For any Gelfand ring $A$, pushforward along $\A^1/\G_m\times\GSpec A\rightarrow\GSpec A$ preserves perfect complexes.
\end{lem}
\begin{proof}
We first note that $\A^1/\G_m$ is suave as both $\A^1$ and the map $\A^1\rightarrow \A^1/\G_m$ are suave and suaveness may be checked suave-locally on the source. We claim that it is also prim: Indeed, this would imply the claim by combining \cite[Lem.s 4.4.9.(ii), 4.5.11.(ii)]{HeyerMann}. Now note that we have 
\begin{equation*}
\P^1/\G_m=\A^1/\G_m\bigsqcup_{\G_m/\G_m} \A^1/\G_m
\end{equation*}
and hence $\A^1/\G_m$ is a retract of $\P^1/\G_m$: The retraction is just given by the projection of one of the two summands onto the point. Thus, it suffices to show that $\P^1/\G_m$ is prim by \cite[Prop.\ 6.17]{SixFunctors}. 

To prove this, observe that $\P^1$ is even weakly cohomologically proper and hence we are reduced to showing that $*/\G_m$ is prim. However, note that $f: *\rightarrow */\G_m$ is suave, hence $f_!1$ is prim over the point by \cite[Lem.\ 4.5.16]{HeyerMann}. As prim objects are stable under shifts and retracts, primness of $*/\G_m$ would thus follow if we can show that $1$ is a retract of $f_!1[1]$. In turn, using the calculations from the proof of \cref{lem:recall-reesgm}, this amounts to showing that $A$ is a retract of $A\{x^{-1}\}^\dagger\oplus xA\{x\}^\dagger$ as modules over $A\{x^{-1}\}^\dagger\oplus xA\{x\}^\dagger$, where $x^i$ acts trivially on $A$ for $i\neq 0$. However, this was already established in loc.\ cit., so we are done.
\end{proof}

In fact, we will typically use the Rees equivalence not for $\A^1/\G_m$, but for the quotient $\ol{\DD}/\ol{\TT}$, where $\ol{\DD}$ is the overconvergent unit disk and $\ol{\TT}$ is the overconvergent unit torus.

\begin{cor}
\label{cor:recall-reesdt}
Let $A$ be any Gelfand ring. Then there is an equivalence between perfect complexes on $\ol{\DD}/\ol{\TT}$ over $A$ and the full subcategory of $\DF(A)$ consisting of objects $\Fil^\bullet M$ satisfying the following two conditions:
\begin{enumerate}[label=(\roman*)]
\item Each $\Fil^i M$ is a perfect complex over $A$.
\item We have $\Fil^i M=0$ for $i\gg 0$ and $\Fil^i M\rightarrow \Fil^{i-1} M$ is an isomorphism for $i\ll 0$.
\end{enumerate}
\end{cor}
\begin{proof}
Combine the previous proposition with the proof of \cref{lem:recall-reesreduction}, where we have shown that $\Perf(\A^1/\G_m)\cong \Perf(\ol{\DD}/\ol{\TT})$.
\end{proof}

To conclude the section, let us record another result that, while not directly related to the Rees equivalence, still falls into the realm of the kind of ``GAGA'' results for perfect complexes we have proved above.

\begin{lem}
\label{lem:recall-cartierperf}
Let $A$ be any Gelfand ring. Pullback along the canonical map
\begin{equation*}
\GSpec A/\widehat{\G}_a\rightarrow \GSpec A/\G_a^\dagger
\end{equation*}
induces an equivalence on perfect complexes. In particular, a perfect complex on $\GSpec A/\G_a^\dagger$ identifies with a perfect complex over $A$ together with a single endomorphism.
\end{lem}
\begin{proof}
As in \cite[Prop.\ II.2.3]{RealLLC}, one shows that the pullback functor 
\begin{equation*}
\D(\GSpec A/\G_a^\dagger)\rightarrow \D(\GSpec A/\widehat{\G}_a)\cong \D(A[t])\;,
\end{equation*}
where the last equivalence is by Cartier duality, see \cite[Prop.\ 4.2.5.(1)]{dRStack}, is fully faithful and identifies the image with the full subcategory spanned by those modules killed after tensoring with the idempotent $A[t]$-algebra $A[t]\{t^{-1}\}^\dagger$. Thus, using \cref{prop:recall-fredholm}, our task is to show that if $M$ is a perfect complex over $A(*)$ equipped with an endomorphism $t: M\rightarrow M$, then $M\tensor_{A[t]} A[t]\{t^{-1}\}^\dagger=0$.

For this, fix once and for all a complex of finite projective $A(*)$-modules representing $M$ in $\D(A(*))$ and choose integers $a\leq b$ such that $M^i=0$ whenever $i\not\in [a, b]$. Then the composite map $\sigma^{\leq b-1} M\rightarrow M\xrightarrow{t} M$, where $\sigma^{\leq\bullet}$ denotes the stupid truncation, corresponds to a map $\sigma^{\leq b-1} M\rightarrow \tau^{\leq b-1} M$ by adjunction and hence in particular induces a map $\sigma^{\leq b-1} M\rightarrow \sigma^{\leq b-1} M$. Moreover, the map $\sigma^{\leq b-1} M\rightarrow M\xrightarrow{t} M\rightarrow M^b[-b]$ is nullhomotopic and hence $t$ induces a map $M^b[-b]\rightarrow M^b[-b]$. Overall, we obtain a commutative diagram
\begin{equation*}
\begin{tikzcd}
\sigma^{\leq b-1}M\ar[r]\ar[d] & M\ar[r]\ar[d, "t"] & M^b[-b]\ar[d] \\
\sigma^{\leq b-1}M\ar[r] & M\ar[r] & M^b[-b]
\end{tikzcd}
\end{equation*}
whose rows are fibre sequences and hence we may check the claim for $M^b[-b]$ and $\sigma^{\leq b-1} M$ in place of $M$. By induction, we are thus eventually reduced to the case where $M$ is a finite projective $A(*)$-module.

As the assertion may also be checked after localising on $\Spec A(*)$, we may furthermore assume that $M\cong A(*)^{\oplus n}$ is free. However, then the endomorphism $t$ is given by an $n\times n$-matrix over $A(*)$. As $A$ is bounded, each entry of this matrix has finite norm and hence the operator $t$ itself has finite norm, which implies $M\tensor_{A[t]} A[t]\{t^{-1}\}^\dagger=0$, as desired.
\end{proof}

\subsection{$p$-adic Lie groups as Gelfand stacks}

Before we finally discuss de Rham and Hyodo--Kato stacks, let us quickly recall the different ways we can view $p$-adic Lie groups as Gelfand stacks. 

The first such way in fact works for any group object $G$ in (light) condensed anima: it simply takes $G$ to the Gelfand stack $G^\Betti$, which will be a group stack by functoriality. Recall that, in the special case where $G$ is a profinite group, e.g.\ a compact $p$-adic Lie group, the latter is just given by
\begin{equation*}
G^\Betti=\GSpec C^\sm(G, \Q_p)\;,
\end{equation*}
where $C^\sm(G, \Q_p)$ denotes the algebra of locally constant $\Q_p$-valued functions on $G$. Owing to the fact that representations of (locally) profinite groups where each orbit map factors through a finite (or discrete) quotient are usually called \emph{smooth representations}, we introduce a slightly different notation for $G^\Betti$ in this context.

\begin{defi}
Let $G$ be a locally profinite group, e.g.\ a $p$-adic Lie group. Then we write $G^\sm$ to denote the group stack $G^\Betti$.
\end{defi}

In keeping with the notation, for any locally profinite group $G$, by \cite[Prop.\ 5.4.2]{SolidLocAnReps}, there is an equivalence
\begin{equation*}
\D(*/G^\sm)\cong \D(\Rep^\sm_{\Q_p}(G))
\end{equation*}
between the derived category of quasicoherent sheaves on the classifying stack $*/G^\sm$ and the derived category of smooth $G$-representations on solid $\Q_p$-vector spaces in the sense of \cite[§5]{SolidLocAnReps}.

What if we are instead interested in locally analytic representations? Following \cite{SolidLocAnReps}, we should make the following definition.

\begin{defi}
Let $G$ be a $p$-adic Lie group and $G_0\subseteq G$ a compact open subgroup. We define the group stack $G^\la$ as
\begin{equation*}
G^\la\coloneqq \coprod_{g\in G/G_0} \GSpec C^\la(G_0, \Q_p)\;,
\end{equation*}
where $C^\la(G_0, \Q_p)$ denotes the algebra of $\Q_p$-valued locally analytic functions on $G_0$.
\end{defi}

Let us point out that one can show that the above definition is independent of the choice of $G_0$. Moreover, as before in the case of smooth representations, by \cite[Thm.\ 4.3.3]{SolidLocAnReps}, for any $p$-adic Lie group $G$, there is an equivalence
\begin{equation*}
\D(*/G^\la)\cong \D(\Rep^\la_{\Q_p}(G))
\end{equation*}
between the derived categories of quasicoherent sheaves on $*/G^\la$ and locally analytic $G$-representations on solid $\Q_p$-vector spaces in the sense of \cite[§3]{SolidLocAnReps}.

\begin{defi}
Let $G$ be a locally profinite group and $G_0\subseteq G$ a compact open subgroup. The group stack $G^\cont$ is defined as 
\begin{equation*}
G^\cont\coloneqq \coprod_{g\in G/G_0} \GSpec C^\cont(G_0, \Q_p)\;,
\end{equation*}
where $C^\cont(G_0, \Q_p)$ denotes the algebra of $\Q_p$-valued continuous functions on $G_0$.
\end{defi}

Let us stress that, while one should think of quasicoherent sheaves on $*/G^\cont$ as continuous $G$-representations on solid $\Q_p$-vector spaces for a locally profinite group $G$, it is \emph{not} true that $\D(*/G^\cont)$ embeds fully faithfully into the derived category $\D(\Q_p[G])$ of modules over the solid group algebra $\Q_p[G]$, which is usually called the \emph{Iwasawa algebra} of $G$.

Now letting $G$ be a $p$-adic Lie group, the natural maps
\begin{equation*}
C^\sm(G_0, \Q_p)\rightarrow C^\la(G_0, \Q_p)\rightarrow C^\cont(G_0, \Q_p)
\end{equation*}
for a compact open subgroup $G_0$ induce corresponding maps of group stacks
\begin{equation*}
G^\cont\rightarrow G^\la\rightarrow G^\sm\;.
\end{equation*}
In particular, we obtain a chain of maps
\begin{equation*}
*/G^\cont\xrightarrow{f} */G^\la\xrightarrow{g} */G^\sm
\end{equation*}
between the corresponding classifying stacks. In terms of representations, pushforward along $f$ then corresponds to taking locally analytic vectors in the sense of \cite[§3]{SolidLocAnReps} while pushforward along $g$ corresponds to taking Lie algebra cohomology, see \cite[Prop.\ 6.2.1]{SolidLocAnReps}, while the respective pullback functors just correspond to forgetful functors.

\subsection{De Rham stacks}

Finally, we recall the definitions and some statements surrounding the analytic de Rham stack, which was first introduced in \cite{dRStack} and has recently been revisited in \cite{dRFF} in the setting of Gelfand stacks. For this, first recall that, in loc.\ cit., the authors have defined a category $\arcStk$ of (qfd) \emph{arc-stacks} over $\Q_p$, which are hypersheaves on (separable qfd) perfectoid $\Q_p$-algebras for the arc-topology introduced in \cite{BerkovichMotives}, and we will use the notation $\Spd\ol{A}$ to refer to the arc-stack represented by a perfectoid $\Q_p$-algebra $\ol{A}$.

For any Gelfand stack $X$, we can get an arc-stack denoted by $X^\diamond$ by restricting its functor of points to perfectoid rings. This yields a functor
\begin{equation*}
(-)^\diamond: \GelfStk\rightarrow\arcStk
\end{equation*}
called the \emph{diamond}, which has a left-adjoint given by extending the assignment $\Spd\ol{A}\mapsto \GSpec\ol{A}$ by colimits. While we will usually not make use of this left-adjoint directly, what we \emph{will} use is its precomposition with $(-)^\diamond$:

\begin{defi}
The postcomposition of $(-)^\diamond: \GelfStk\rightarrow\arcStk$ with its left-adjoint is called the \emph{perfectoidisation} and denoted by
\begin{equation*}
\widehat{(-)}: \GelfStk\rightarrow\GelfStk\;.
\end{equation*}
\end{defi}

Note that, by definition, any Gelfand stack $X$ receives a natural map from its perfectoidisation $\widehat{X}$. The diamond functor also has a right-adjoint.

\begin{defi}
The right-adjoint of the functor $(-)^\diamond: \GelfStk\rightarrow\arcStk$ is called the \emph{analytic de Rham stack} and denoted
\begin{equation*}
(-)^\dR: \arcStk\rightarrow\GelfStk\;.
\end{equation*}
We will also use the same name and notation to refer to the precomposition of this functor with $(-)^\diamond: \GelfStk\rightarrow\arcStk$.
\end{defi}

Again, note that, by definition, there is always a natural map $X\rightarrow X^\dR$ for any Gelfand stack $X$.

\begin{ex}
\label{ex:recall-drbetti}
Recall from \cite[Ex.\ 4.1.9]{dRFF} that there is a colimit-preserving functor
\begin{equation*}
\ul{(-)}: \widehat{\Shv}(\mathrm{Pro}_{\mathbb{N}}(\mathrm{Fin}))\rightarrow\arcStk
\end{equation*}
given by sending a light profinite set $S$ to the arc-stack $\ul{S}$ whose functor of points is given by
\begin{equation*}
\ul{S}(\ol{A})=\Hom_\cont(\cal{M}(\ol{A}), S)\;,
\end{equation*}
where $\cal{M}(\ol{A})$ denotes the Berkovich space of $\ol{A}$. By \cite[Ex.\ 5.6.8]{dRFF}, the composition 
\begin{equation*}
(-)^\dR\circ\ul{(-)}: \widehat{\Shv}(\mathrm{Pro}_{\mathbb{N}}(\mathrm{Fin}))\rightarrow\GelfStk
\end{equation*}
identifies with the Betti stack functor from \cref{prop:recall-bettistack}.
\end{ex}

\begin{ex}
\label{ex:recall-glaviadr}
Let $\mathbb{G}$ be a Berkovich smooth group object in the category of $\dagger$-rigid spaces in the sense of \cite[Def.\ 4.3.6]{dRFF}, see also \cite{GrosseKlonne}. Then $G\coloneqq \mathbb{G}(\Q_p)$ naturally admits the structure of a $p$-adic Lie group and we can give a very simple alternative description of $G^\la$ in terms of the analytic de Rham stack of $\mathbb{G}$. For this, first recall that there is a natural map $\ul{G}\rightarrow \mathbb{G}^\diamond$ and this induces a map $G^\sm\rightarrow \mathbb{G}^\dR$ upon taking de Rham stacks by the previous example. Now $G^\la$ may be described by the cartesian diagram
\begin{equation*}
\begin{tikzcd}
G^\la\ar[r]\ar[d] & G^\sm\ar[d] \\
\mathbb{G}\ar[r] & \mathbb{G}^\dR\nospacepunct{\;.}
\end{tikzcd}
\end{equation*}
Indeed, any analytic function on $\mathbb{G}$ yields a locally analytic function on $\mathbb{G}(\Q_p)$, which defines the map $G^\la\rightarrow\mathbb{G}$. Moreover, by \cite[Lem.\ 6.2.2]{dRStack}, the map $G^\la\rightarrow G^\sm$ is surjective with kernel $(1\subseteq \mathbb{G})^\dagger$, the overconvergent neighbourhood of $1$ inside $\mathbb{G}$, and the same is true for $\mathbb{G}\rightarrow\mathbb{G}^\dR$ by \cite[Ex.\ 6.1.8.(2)]{dRStack}, hence the diagram is indeed cartesian.
\end{ex}

To explicitly describe the functor of points of the de Rham stack of a given arc-stack, first recall from \cite[Prop.\ 4.6.4]{dRFF} that any Gelfand ring admits a $!$-cover by a ring of the following form.

\begin{defi}
A Gelfand ring $A$ is called \emph{nilperfectoid} if its $\dagger$-reduction $\ol{A}$ is perfectoid. It is moreover called \emph{totally disconnected} if the Berkovich space $\cal{M}(A)=\cal{M}(\ol{A})$ is a profinite set. 
\end{defi}

Now given any arc-stack $X$, it suffices to specify the functor of points of $X^\dR$ on nilperfectoid rings $A$, and there it is given by the formula
\begin{equation*}
X^\dR(A)=X(\ol{A})\;,
\end{equation*}
see \cite[Prop.\ 4.6.5]{dRFF}. If $X$ is a Gelfand stack, this shows that $X^\dR$ is obtained by sheafifying the assignment $A\mapsto X(\ol{A})$ for Gelfand rings $A$.

One of the main results of \cite{dRFF} is the fact that the functor
\begin{equation*}
(-)^\dR: \arcStk\rightarrow\GelfStk
\end{equation*}
commutes with colimits, which is also what enables \cref{ex:recall-drbetti} above. Moreover, it calculates de Rham cohomology of partially proper rigid analytic spaces. 

\begin{prop}
\label{prop:recall-dr}
For any smooth partially proper rigid space $X$ over $\Q_p$, there is a natural isomorphism
\begin{equation*}
R\Gamma(X^\dR, \O)\cong R\Gamma_\dR(X)\;,
\end{equation*}
where the right-hand side denotes usual de Rham cohomology of $X$. Moreover, there is an equivalence
\begin{equation*}
\Vect(X^\dR)\cong \{\text{vector bundles with connection on $X$}\}\;,
\end{equation*}
where the connection in the target category is always assumed to be flat.
\end{prop}
\begin{proof}
See \cite[Prop.\ 5.2.1, Rem.\ 5.2.3]{dRFF}.
\end{proof}

The analytic de Rham stack also has a filtered variant. While this is usually defined as living over the quotient $\A^1/\G_m$, it will be more convenient for us to use $\ol{\DD}/\ol{\T}$ as the base, where we recall that $\ol{\DD}$ denotes the overconvergent unit disk and $\ol{\T}$ is the overconvergent unit torus. Recall that we have already shown in \cref{lem:recall-reesreduction} that $\ol{\DD}/\ol{\T}$ and $\A^1/\G_m$ have equivalent categories of perfect complexes; more precisely, by \cref{prop:recall-rees}, perfect complexes on both of these stacks are equivalent to finitely filtered perfect complexes over $\Q_p$. 

Let us also shortly comment on how the moduli interpretations of $\ol{\DD}/\ol{\T}$ and $\A^1/\G_m$ compare. For this, recall that the latter parametrises \emph{generalised Cartier divisors}, i.e.\ maps $L\rightarrow A$, where $L$ is an invertible $A$-module (note that we implicitly need to use \cref{prop:recall-fredholm} here). A similar description holds for $A$-points of $\ol{\DD}/\ol{\T}$.

\begin{defi}
Let $A$ be a Gelfand ring. A \emph{normed line bundle} over $A$ is a line bundle $L$ over $A$ which is $!$-locally equipped with an invertible $A^{\leq 1}$-submodule $L^{\leq 1}\subseteq L$.
\end{defi}

\begin{rem}
It is probably not true that one can always obtain an actual $A^{\leq 1}$-module $L^{\leq 1}$ from a normed line bundle $L$ over $A$, i.e.\ without having to $!$-localise.
\end{rem}

\begin{defi}
Let $A$ be a Gelfand ring. A \emph{normed generalised Cartier divisor} of norm at most $r\geq 0$ is a normed line bundle $L$ over $A$ together with a map $L\rightarrow A$ such that, $!$-locally on $A$, the composite $L^{\leq 1}\rightarrow L\rightarrow A$ factors through $A^{\leq r}$.
\end{defi}

Now it is clear from the definition that the quotient $\ol{\DD}/\ol{\T}$ exactly parametrises normed generalised Cartier divisors of norm at most $1$. Moreover, given such a divisor corresponding to a map $\GSpec A\rightarrow\ol{\DD}/\ol{\T}$, postcomposing with the map $\ol{\DD}/\ol{\T}\rightarrow\A^1/\G_m$ exactly corresponds to forgetting the norm. We can now define the filtered version of $X^\dR$.

\begin{defi}
Let $X$ be a Gelfand stack over $\Q_p$. Its (analytic) \emph{Hodge-filtered de Rham stack} $X^{\dR, +}$ is the Gelfand stack over $\ol{\DD}/\ol{\T}$ obtained by sheafifying the assignment
\begin{equation*}
X^{\dR, +}(\GSpec A\rightarrow\ol{\DD}/\ol{\T})\coloneqq \{\text{maps}\,\GSpec\Cone(L\tensor_A \Nil^\dagger(A)\rightarrow A)\rightarrow X\}\;,
\end{equation*}
where $L\rightarrow A$ is the normed generalised Cartier divisor classified by the map $\GSpec A\rightarrow\ol{\DD}/\ol{\T}$.
\end{defi}

If $t$ denotes the coordinate on $\ol{\DD}$, let us write
\begin{equation*}
t: \O\langle -1\rangle\rightarrow \O
\end{equation*}
for the universal normed generalised Cartier divisor over $\ol{\DD}/\ol{\T}$. Given any sheaf $E$ on $X^{\dR, +}$ for some Gelfand stack $X$, we can then obtain a filtered complex of $\Q_p$-vector spaces by taking cohomology of twists of $E$, i.e.\ there is a diagram
\begin{equation*}
\dots\xrightarrow{t} R\Gamma(X^{\dR, +}, E\langle -n-1\rangle)\xrightarrow{t} R\Gamma(X^{\dR, +}, E\langle -n\rangle)\xrightarrow{t} R\Gamma(X^{\dR, +}, E\langle -n+1\rangle)\xrightarrow{t}\dots\;,
\end{equation*}
and we will call this the \emph{filtered cohomology} of $X^{\dR, +}$ and write $\Fil^\bullet R\Gamma(X^{\dR, +}, E)$. (This applies more generally to all Gelfand stacks over $\ol{\DD}/\ol{\T}$, not just filtered de Rham stacks.)

\begin{prop}
\label{prop:recall-fildr}
Let $X$ be a smooth partially proper rigid space over $\Q_p$. Then there is a natural isomorphism
\begin{equation*}
\Fil^\bullet R\Gamma(X^{\dR, +}, \O)\cong \Fil^\bullet_\Hod R\Gamma_\dR(X)
\end{equation*}
between filtered cohomology of $X^{\dR, +}$ and Hodge-filtered de Rham cohomology of $X$. Moreover, there is an equivalence of categories
\begin{equation*}
\Vect(X^{\dR, +})\cong \{\text{filtered vector bundles with connection on $X$}\}\;,
\end{equation*}
where the filtration in the target category is always assumed to be locally by direct summands and the connection is always assumed to be flat and Griffiths transversal with respect to the filtration, i.e.\ it carries $\Fil^i E$ into $\Fil^{i-1} E\tensor_{\O_X} \Omega_X^1$.
\end{prop}
\begin{proof}
For the first part, see \cite[Rem.\ 5.2.2]{dRFF}. For the second part, we first construct a natural functor
\begin{equation}
\label{eq:recall-fildrfunctor}
\Vect(X^{\dR, +})\rightarrow \{\text{filtered vector bundles with connection on $X$}\}\;
\end{equation}
obtained as follows: Pulling back any vector bundle on $X^{\dR, +}$ along the natural map
\begin{equation*}
X\times \ol{\DD}/\ol{\T}\rightarrow X^{\dR, +}
\end{equation*}
induced by the canonical map $A\rightarrow \Cone(L\tensor_A \Nil^\dagger(A)\rightarrow A)$ yields a filtered vector bundle $\Fil^\bullet E$ on $X$ by \cref{prop:recall-rees}. As the underlying unfiltered vector bundle arises via pullback along $X\mapsto X^\dR$, it is naturally equipped with a flat connection $\nabla: E\rightarrow E\tensor_{\O_X} \Omega_X^1$ by \cref{prop:recall-dr}.

To check the Griffiths transversality, we may then work locally on $X$ and, in particular, assume that it admits a Berkovich étale map $X\rightarrow \ol{\DD}^n$ for some $n\geq 0$. Then the diagram
\begin{equation*}
\begin{tikzcd}
X\times\ol{\DD}/\ol{\T}\ar[r]\ar[d] & X^{\dR, +}\ar[d] \\
\ol{\DD}^n\times\ol{\DD}/\ol{\T}\ar[r] & (\ol{\DD}^n)^{\dR, +}
\end{tikzcd}
\end{equation*}
is cartesian as the map $X^{\dR, +}\rightarrow (\ol{\DD}^n)^{\dR, +}$ is $\dagger$-formally étale. Moreover, we have
\begin{equation*}
(\ol{\DD}^n)^{\dR, +}\cong \ol{\DD}^n\times \ol{\DD}/\ol{\T}\;\big/\;\G_a^\dagger\langle -1\rangle^n
\end{equation*}
by definition, where $\G_a^\dagger\langle -1\rangle^n$ acts on $\ol{\DD}^n$ by $t$-scaled translation, where $t: \O\langle -1\rangle\rightarrow\O$ is the tautological normed generalised Cartier divisor on $\ol{\DD}/\ol{\T}$. Thus, we conclude that
\begin{equation*}
X^{\dR, +}\cong X\times \ol{\DD}/\ol{\T}\;\big/\;\G_a^\dagger\langle -1\rangle^n\;,
\end{equation*}
and now putting together the Rees equivalence from \cref{prop:recall-rees} with Cartier duality for $\G_a^\dagger$ as in \cref{lem:recall-cartierperf}, we see that any vector bundle on $X$ with a flat connection and a filtration arising via pullback from $X^{\dR, +}$ as above satisfies Griffiths transversality. This shows that the construction we have given above indeed yields a well-defined natural functor (\ref{eq:recall-fildrfunctor}). Moreover, the same argument then shows that this functor is actually an equivalence in the situation where $X$ is equipped with a Berkovich étale map $X\rightarrow\ol{\DD}^n$.

This is enough to conclude: As $X\mapsto X^{\dR, +}$ is compatible with rational localisations, which one deduces by base change from \cref{prop:defis-openloc} below, we may check whether (\ref{eq:recall-fildrfunctor}) is an equivalence by localising on $X$, and this finishes the proof.
\end{proof}

The proof of the first part of the above statement uses the fact that one can explicitly describe the base change of $X^{\dR, +}$ along $*/\ol{\T}\rightarrow\ol{\DD}/\ol{\T}$.

\begin{defi}
Let $X$ be a Gelfand stack. The (analytic) \emph{Hodge stack} $X^\Hod$ of $X$ is defined as
\begin{equation*}
X^\Hod\coloneqq X^{\dR, +}\times_{\ol{\DD}/\ol{\T}} */\ol{\T}\;.
\end{equation*}
\end{defi}

Then the key to the proof of the first part of \cref{prop:recall-fildr} is that, for a smooth partially proper rigid space $X$, we have
\begin{equation*}
X^\Hod\cong (X\times */\ol{\T})/\cal{T}_{X/\Q_p}^\dagger\langle -1\rangle
\end{equation*}
essentially by deformation theory, see \cite[Prop.\ 5.2.3.(2a)]{dRStack}; here, $\cal{T}_{X/\Q_p}^\dagger$ denotes the overconvergent neighbourhood of the zero section inside the restriction of $\AnSpec_X \Sym_X^\bullet \Omega^1_{X/\Q_p}$ to the test category of Gelfand rings and the action on $X\times */\ol{\T}$ is trivial. Then one deduces by Cartier duality, see \cite[Thm.\ 4.3.13]{dRStack} and \cref{lem:recall-cartierperf}, that cohomology of $X^\Hod$ computes Hodge cohomology of $X$, i.e.\ 
\begin{equation*}
\bigoplus_{i\in\Z} R\Gamma(X^\Hod, \O\langle i\rangle)\cong \bigoplus_{i\in\Z} R\Gamma(X, \Omega_X^i)\;.
\end{equation*}

We also note that the base change of $X^{\dR, +}$ along $*=\ol{\T}/\ol{\T}\rightarrow \ol{\DD}/\ol{\T}$ recovers $X^\dR$ as the normed generalised Cartier divisor $L\rightarrow A$ is an isomorphism over this locus and hence $\Cone(L\tensor_A \Nil^\dagger(A)\rightarrow A)$ simplifies to $\ol{A}$. Slightly more generally, the same argument works over the larger locus $\ol{\DD}^\times/\ol{\T}\subseteq \ol{\DD}/\ol{\T}$, where $\ol{\DD}^\times$ denotes the punctured overconvergent unit disk, and thus, using \cref{lem:defis-rhvariant}, shows that
\begin{equation*}
X^{\dR, +}\times_{\ol{\DD}/\ol{\T}} \ol{\DD}^\times/\ol{\T}\cong X^\dR\times (0, 1]\;.
\end{equation*}

\subsection{Hyodo--Kato stacks}

Let us end by shortly discussion \emph{Hyodo--Kato stacks}, which in a way also belong the circle of ideas around analytic de Rham stacks and were introduced in \cite{dRFF}. For this, first recall that, for any arc-stack $X$ over $\Q_p$, there is an arc-stack
\begin{equation*}
Y_X^\diamond\coloneqq X\times_{\Spd\F_p} \Spd\Q_p
\end{equation*}
over $\Q_p$, where the structure map is via projection onto the second factor and we point out that the fibre product of course has to be taken as arc-stacks \emph{over $\F_p$} (!), not over $\Q_p$. However, let us emphasise that all the arc-stacks that will occur in this paper will be over $\Q_p$ and arc-stacks over $\F_p$ will only feature implicitly via the above definition. Let us also note that, while $Y_X^\diamond$ is called the ``open punctured curve'' over $X$ in \cite{dRFF}, we will mostly refer to it as the \emph{$Y$-curve} of $X$.

As an arc-stack over $\F_p$, the arc-stack $X$ is equipped with a Frobenius automorphism and this induces a Frobenius 
\begin{equation*}
\phi: Y_X^\diamond\rightarrow Y_X^\diamond\;,
\end{equation*}
which is an isomorphism. Recall that the \emph{Fargues--Fontaine curve} $\FF_X^\diamond$ of $X$ is defined as
\begin{equation*}
\FF_X^\diamond\coloneqq Y_X^\diamond/\phi^\Z\;.
\end{equation*}
This leads to the following definition of Hyodo--Kato stacks:

\begin{defi}
Let $X$ be an arc-stack over $\Q_p$. The \emph{Hyodo--Kato stack} $X^\HK$ of $X$ is defined as the Gelfand stack
\begin{equation*}
X^\HK\coloneqq (\FF_X^\diamond)^\dR\;.
\end{equation*}
If $X$ is a Gelfand stack, then we will also write $X^\HK$ for the Hyodo--Kato stack of $X^\diamond$ and simply call it the Hyodo--Kato stack of $X$.
\end{defi}

The name of $X^\HK$ is of course supposed to hint at the fact that cohomology of $X^\HK$ agrees with Hyodo--Kato cohomology, which in the setting of analytic geometry was probably first introduced by Colmez--Nizio{\l} in \cite{padicComparisons}. Recall that this is a cohomology theory for smooth partially proper rigid spaces over $\Q_p$ valued in $(\phi, N, \Gal_{\Q_p})$-modules over the maximal unramified extension $\Q_p^\un$ of $\Q_p$. Indeed, by one of the main theorems of \cite{dRFF}, this is also the kind of linear algebraic structure one obtains on cohomology of Hyodo--Kato stacks:

\begin{thm}
\label{thm:recall-qphk}
The category $\Perf(\Q_p^\HK)$ of perfect complexes on $\Q_p^\HK$ is equivalent to the category of $(\phi, N)$-modules in $\Perf(\Q_p^\un)$ with a smooth $\Gal_{\Q_p}$-action. In particular, vector bundles on $\Vect(\Q_p^\HK)$ are $(\phi, N, \Gal_{\Q_p})$-modules over $\Q_p^\un$.
\end{thm}
\begin{proof}
See \cite[Thm.\ 7.1.1]{dRFF}.
\end{proof}

\begin{ex}
The isocrystal $(\Q_p, p^{-1}\phi)$ corresponds to a line bundle on $\Q_p^\HK$ which we will denote by $\O\{1\}$. For any Gelfand stack $X$, we will also use $\O\{1\}$ to denote the pullback of $\O\{1\}$ along the induced map $X^\HK\rightarrow\Q_p^\HK$.
\end{ex}

In order to apply the theorem above and extract $(\phi, N, \Gal_{\Q_p})$-modules from the cohomology of the stack $X^\HK$, we should show that pushforward along $X^\HK\rightarrow\Q_p^\HK$ preserves perfect complexes, at least when $X$ is a smooth partially proper qcqs rigid space over $\Q_p$.

\begin{lem}
\label{lem:geom-hkprim}
For any qcqs arc-stack $X$ over $\Q_p$, the map $X^\HK\rightarrow \Q_p^\HK$ is prim.
\end{lem}
\begin{proof}
By \cite[Lem.\ 4.5.8]{HeyerMann}, we may pass to a prim descendable cover of the source and thus, using \cite[Prop.\ 4.7.12]{dRFF}, we may pass to a finite strict closed cover of $X$ and hence assume that $X$ admits an étale map to $(\A^n)^\diamond$ for some $n\geq 0$. By quasicompactness, this map factors through a closed disk of finite radius, which we may without loss of generality assume to have radius $1$, and after changing coordinates on the target, we may thus assume that there is an étale map $X\rightarrow (\ol{\T}^n)^\diamond$ for some $n\geq 0$.

First observe that $X^\HK\rightarrow (\ol{\T}^n)^\HK$ is prim. Indeed, the map $\FF_X^\diamond\rightarrow \FF_{(\ol{\T}^n)^\diamond}^\diamond$ is qcqs and étale and any such map induces a prim map on de Rham stacks by the proof of \cite[Prop.\ 5.1.4]{dRFF}. Therefore, it remains to show that $(\ol{\T}^n)^\HK\rightarrow \Q_p^\HK$ is prim, which we may reduce to the case $n=1$ by compatibility of the Hyodo--Kato stack with limits.

Then recall from \cite[Rem.\ 6.2.6.(1)]{dRFF} that
\begin{equation*}
\ol{\T}^\HK\cong \T_{\infty, \Q_p^\cycl}^\HK/\Z_p(1)^\sm\rtimes\Z_p^{\times, \sm} \cong (\FF_{\Q_p^\cycl}^\diamond\times_{\Spd\Q_p} \T_\infty)^\dR/\Z_p(1)^\sm\rtimes\Z_p^{\times, \sm}\;,
\end{equation*}
where $\T_\infty\coloneqq\Spd\Q_p\langle T^{\pm 1/p^\infty}\rangle$ denotes the perfectoid $n$-dimensional torus, the semidirect product is defined by letting $\Z_p^\times$ act on $\Z_p(1)$ by multiplication, $\Z_p^\times$ acts on $\Q_p^\cycl$ via the usual Galois action and $\Z_p(1)$ acts on $\ol{\T}$ via $\theta.T^{1/p^m}=[\epsilon^{1/p^m}]^\theta T^{1/p^m}$, where $\epsilon=(1, \zeta_p, \zeta_{p^2}, \dots)\in\Q_p^{\cycl, \flat}$ is as usual. Moreover, we have
\begin{equation*}
\Q_p^\HK\cong (\FF_{\Q_p^\cycl}^\diamond)^\dR/\Z_p^{\times, \sm}\;.
\end{equation*}
Now note that
\begin{equation*}
\T_\infty^\dR\cong \ol{\T}^\la_\infty/\G_m^\dagger\;,
\end{equation*}
where $\ol{\T}^\la_\infty=\GSpec \colim_m \Q_p\langle T^{\pm 1/p^m}\rangle_{\leq 1}$, since the $p$-th power map on $\G_m^\dagger$ is an isomorphism. Thus, we see that $\T_\infty^\dR$ is even weakly cohomologically proper as the same is true for $\ol{\T}^\la_\infty$ and $*/\G_m^\dagger$, where in the latter case we use $\G_m^\dagger\cong\G_a^\dagger$ via the logarithm and \cite[Prop.\ 4.2.5.(2)]{dRStack}. As also $*/\Z_p^\sm$ is prim by \cite[Prop.\ 6.2.4]{dRStack}, we overall conclude that $\ol{\T}^\HK\rightarrow\Q_p^\HK$ is prim, as desired.
\end{proof}

\begin{cor}
\label{cor:recall-hkfin}
For any smooth partially proper qcqs rigid space $X$ over $\Q_p$, pushforward along $X^\HK\rightarrow \Q_p^\HK$ preserves perfect complexes. In particular, the pushforward of any perfect complex is naturally equipped with the structure of a $(\phi, N, \Gal_{\Q_p})$-module in $\Perf(\Q_p^\un)$.
\end{cor}
\begin{proof}
The map $X^\HK\rightarrow\Q_p^\HK$ is prim by \cref{lem:geom-hkprim} and also suave by \cite[Thm.\ 6.3.1.(2)]{dRFF}, hence pushforward preserves dualisable objects by combining \cite[Lem.s 4.4.9.(ii), 4.5.11.(ii)]{HeyerMann} and then the claim follows from \cref{prop:recall-fredholm}. The final assertion follows using \cref{thm:recall-qphk}.
\end{proof}

With the above result, it is now sensible to make the following definition.

\begin{defi}
Let $X$ be a smooth partially proper qcqs rigid space over $\Q_p$. We define the \emph{Hyodo--Kato cohomology} $R\Gamma_\HK(X)$ of $X$ as the $(\phi, N, \Gal_{\Q_p})$-module in $\Perf(\Q_p^\un)$ corresponding to the pushforward of the structure sheaf along $X^\HK\rightarrow\Q_p^\HK$.
\end{defi}

According to \cite[Rem.\ 6.1.3]{dRFF}, it is expected that the above definition agrees with classical definitions of Hyodo--Kato cohomology. More specifically, one should be able to deduce this from the fact that the cohomology of the Hyodo--Kato stack defines a realisation of Berkovich motives according to the recent PhD thesis of Ko Aoki, see \cite[Thm.\ 1.7]{AokiThesis}, together with the interpretation of Hyodo--Kato cohomology via motives from \cite{BindaGallauerVezzani}. While we will not use this expected relation in any of our arguments, we will still implicitly assume it when we formulate some of our theorems in more classical language or compare them to previous results.

\newpage

\section{Analytic prismatisation over $\Q_p$}
\label{sect:prism}

We now recall some facts about (rational) analytic prismatisation of Gelfand stacks over $\Q_p$, which has been developed by Anschütz--Le Bras--Rodríguez Camargo--Scholze. From the outset, let us emphasise that, because we restrict both our geometry and coefficients to characteristic $0$, i.e.\ we allow Gelfand stacks over $\Q_p$ as inputs and our cohomology theories will have $\Q_p$-coefficients, the overall situation is very simple: Basically, the rational analytic prismatisation $X^\prism$ of a Gelfand stack $X$ over $\Q_p$ will just combine the Hyodo--Kato stack $X^\HK$ and the stack $X^{\Div^1}$ from \cite{AnPrism}, the definition of which we will review further below. 

Roughly speaking, the construction of $X^\prism$ for a Gelfand stack $X$ over $\Q_p$ goes as follows. First, recall that nilperfectoids form a basis for the $!$-topology on Gelfand rings, hence it suffices to specify $X^\prism(A)$ for nilperfectoid rings $A$. For any such nilperfectoid, we will define a derived Berkovich space $Y_A$, which will be a $\dagger$-thickening of the $Y$-curve $Y_{\ol{A}}$ of the perfectoid $\ol{A}$ equipped with a Frobenius $\phi: Y_A\rightarrow Y_A$ compatible with the one on $Y_{\ol{A}}$. Then $\Q_p^\prism$ will be the moduli space of degree $1$ divisors $D\subseteq Y_A$, where the degree $1$ condition is asked after pullback to $Y_{\ol{A}}$, and $X^\prism$ will be defined as the Gelfand stack over $\Q_p^\prism$ given by
\begin{equation*}
X^\prism(\GSpec A\rightarrow \Q_p^\prism)=\{\text{maps}\,D\rightarrow X\}\;,
\end{equation*}
where $D\subseteq Y_A$ is the degree $1$ divisor classified by the map $\GSpec A\rightarrow\Q_p^\prism$.

For any Gelfand stack $X$, the Gelfand stack $X^\prism$ will have the following features:
\begin{enumerate}[label=(\roman*)]
\item There is an isomorphism $(X^\prism)^\dR\cong Y_{X^\diamond}^\dR$.
\item By virtue of (i), there is a radius map
\begin{equation*}
\kappa: X^\prism\rightarrow (0, \infty)\;.
\end{equation*}
\item There is a Frobenius $\phi: X^\prism\rightarrow X^\prism$ which is sent to multiplication by $p$ under $\kappa$.
\item Over $(1, \infty)$, the map $X^\prism\rightarrow Y_{X^\diamond}^\dR$ induced by (i) is an isomorphism and, in particular, Frobenius induces an isomorphism $X^\prism_{(1, \infty)}\xrightarrow{\cong} X^\prism_{(p, \infty)}$.
\item By (iv), there is a map $i_\dR: X^\dR\rightarrow X^\prism$ which factors through $X^\prism_{[p, p]}$ and there are further copies of $X^\dR$ sitting inside $X^\prism$ at all positive Frobenius translates.
\item Frobenius also induces an isomorphism $X^\prism_{(0, 1)}\xrightarrow{\cong} X^\prism_{(0, p)}$.
\item There is a map $i_\HT: X^\HT\rightarrow X^\prism$ which factors through $X^\prism_{[1, 1]}$ and there are further copies of $X^\HT$ sitting inside $X^\prism$ at all negative Frobenius translates; here, $X^\HT$ denotes the \emph{analytic Hodge--Tate stack} from \cite{AnPrism}.
\end{enumerate}
The situation is summarised by the schematic picture in \cref{fig:xprism}.

\begin{figure}

\begin{center}
\begin{tikzpicture}
  \def\r{9}
  \def\shortr{5} 
  \def\longr{9.5}
  \def\vlongr{10.5}

  \draw[dotted, thick] (0,0) -- (0,\shortr); 
  \draw[thick] (\r,0) arc[start angle=0, end angle=90, radius=\r];
  
  \draw[->, thick] (0.866*\vlongr, 1/2*\vlongr) arc[start angle=30, end angle=60, radius=\vlongr];
  
  \node at (1.4142/2*\vlongr + 0.3, 1.4142/2*\vlongr + 0.3) {$\varphi$};

  \def\shorttr{\r-0.2}
  \def\longgr{\r+0.2}

  \draw[thick] (\shorttr,0.2) arc (-180:0:0.2);
  \node at (\longr, 0) {$0$};
  
  \draw[thick] (7.5:\shorttr) -- (7.5:\longgr);
  \node at (7.5:\longr+0.3) {$1/p^3$};

  \draw[thick] (15:\shorttr) -- (15:\longgr);
  \node at (15:\longr+0.3) {$1/p^2$};

  \draw[thick] (25:\shorttr) -- (25:\longgr);
  \node at (25:\longr+0.2) {$1/p$};

  \draw[thick] (45:\shorttr) -- (45:\longgr);
  \node at (45:\longr) {$1$};

  \draw[thick] (65:\shorttr) -- (65:\longgr);
  \node at (65:\longr) {$p$};

  \draw[thick] (77.5:\shorttr) -- (77.5:\longgr);
  \node at (76.5:\longr+0.1) {$p^2$};
  
  \draw[thick] (85:\shorttr) -- (85:\longgr);
  \node at (84:\longr+0.1) {$p^3$};

  \draw[thick] (0.2, \longgr) arc (90:270:0.2);
  \node at (0, \longr) {$\infty$};

  \begin{scope}

  \clip (0,0) rectangle (\longr,\longr);
  
  \begin{scope}[transform canvas={rotate=7.5}]
        \shade[shading=axis, bottom color=white, top color=white, middle color=green, shading angle=0] (0.2, -0.1) rectangle (\shortr, 0.1);
  \end{scope} 

  \begin{scope}[transform canvas={rotate=15}]
        \shade[shading=axis, bottom color=white, top color=white, middle color=green, shading angle=0] (0.2, -0.1) rectangle (\shortr, 0.1);
  \end{scope} 

  \begin{scope}[transform canvas={rotate=25}]
        \shade[shading=axis, bottom color=white, top color=white, middle color=green, shading angle=0] (0.2, -0.1) rectangle (\shortr, 0.1);
  \end{scope}

  \begin{scope}[transform canvas={rotate=45}]
        \shade[shading=axis, bottom color=white, top color=white, middle color=green, shading angle=0] (0.2, -0.1) rectangle (\shortr, 0.1);
  \end{scope}

  \draw[line width=0.2mm, color=blue] (0.4266*0.2,0.9063*0.2) -- (0.4226*\shortr, 0.9063*\shortr);

  \draw[line width=0.2mm, color=blue] (0.21644*0.2,0.9763*0.2) -- (0.21644*\shortr, 0.9763*\shortr);

  \draw[line width=0.2mm, color=blue] (0.087156*0.2,0.9962*0.2) -- (0.087156*\shortr, 0.9962*\shortr);

  \end{scope}
  
  \def\shortishr{\shortr+0.5}
  
  \node at (7.5:\shortishr+0.1) {$X^\HT$};
  
  \node at (14.5:\shortishr+0.1) {$X^\HT$};
  
  \node at (24.5:\shortishr+0.1) {$X^\HT$};
  
  \node at (44:\shortishr+0.1) {$X^\HT$};
  
  \node at (64:\shortishr) {$X^\dR$};
  
  \node at (75.5:\shortishr) {$X^\dR$};
  
  \node at (83:\shortishr) {$X^\dR$};
  
  \draw[thick, dotted] (0,0) -- (\shortr,0);
  
  \draw[->, thick] (45:\shortishr+1) -- (45:\shorttr-0.5);
  \node at (48:7.5) {$\kappa$};

\end{tikzpicture}
\end{center}

\captionsetup{justification=centering}
\caption{A schematic picture of $X^\prism$}
\label{fig:xprism}
\end{figure}

\subsection{$Y$-curves of nilperfectoids}

Let us first recall the definition of $Y_{\ol{A}}$ for $\ol{A}$ perfectoid. For this, recall that $Y_{\ol{A}}^\diamond$ may be obtained as the arc-stack corresponding to the adic space
\begin{equation*}
\Spa \mathbb{A}_\inf(\ol{A})\setminus \{p[\pi^\flat]=0\}\;,
\end{equation*}
where $\pi\in\ol{A}$ is a perfectoid pseudouniformiser. Transporting this into the world of Gelfand stacks, we set
\begin{equation*}
Y_{\ol{A}}\coloneqq \colim_{r\rightarrow 0^+, s\rightarrow \infty} \GSpec \mathbb{B}_{[r, s]}(\ol{A})\;,
\end{equation*}
where we use the notation
\begin{equation*}
\mathbb{B}_{[r, s]}(\ol{A})\coloneqq \mathbb{A}_\inf(\ol{A})[\tfrac{1}{[\pi^\flat]}]\langle\tfrac{p^s}{[\pi^\flat]}, \tfrac{[\pi^\flat]}{p^r}\rangle
\end{equation*}
and we emphasise that this definition depends on the choice of $\pi$ even though this is not apparent from the notation on the left-hand side. By \cite[Cor.\ 4.4.10]{dRFF}, the functor $\ol{A}\mapsto Y_{\ol{A}}$ sends arc-hypercovers of totally disconnected perfectoids to $!$-equivalences of Gelfand stacks and hence satisfies hyperdescent for the arc-topology. In particular, we may extend it to a functor
\begin{equation*}
Y_{(-)}: \arcStk\rightarrow\GelfStk\;, \hspace{0.3cm} X\mapsto Y_X\;.
\end{equation*}

We point out that $Y_{\ol{A}}$ is always equipped with a Frobenius $\phi: Y_{\ol{A}}\rightarrow Y_{\ol{A}}$ induced by the one on $\A_\inf(\ol{A})$ and that this in turn induces a Frobenius $\phi: Y_X\rightarrow Y_X$ on $Y_X$ for any arc-stack $X$. Moreover, recall that Fontaine's map 
\begin{equation*}
\theta: \A_\inf(\ol{A})=W(\ol{A}^{\flat\circ})\rightarrow\ol{A}^\circ\;, \hspace{0.3cm} \sum_{n\geq 0} [a_n]p^n\mapsto \sum_{n\geq 0} a_n^\sharp p^n
\end{equation*}
induces a map
\begin{equation*}
\iota: \GSpec\ol{A}\rightarrow Y_{\ol{A}}
\end{equation*}
via which $\GSpec\ol{A}$ is a divisor in $Y_{\ol{A}}$. Consequently, for any arc-stack $X$, we obtain an analogous map $\iota: X\rightarrow Y_X$, which makes $X$ a divisor in $Y_X$, where we view $X$ as a Gelfand stack via the left-adjoint of the diamond functor.

\begin{ex}
For $X=\Spd\Q_p$, we obtain
\begin{equation*}
Y_{\Spd\Q_p}\cong Y_{\Q_p^\cycl}/\Z_p^{\times, \cont}\cong Y_{\C_p}/\Gal_{\Q_p}^\cont\;. \qedhere
\end{equation*}
\end{ex}

\begin{rem}
For any arc-stack $X$, it is indeed true that the diamond of $Y_X$ agrees with $Y_X^\diamond$ as defined previously, i.e.\ it is isomorphic to $X\times_{\Spd\F_p} \Spd\Q_p$. However, in line with \cite[Rem.\ 6.1.7]{dRFF}, let us emphasise that $Y_X$ is \emph{not} just obtained by applying the left-adjoint of the diamond functor to $Y_X^\diamond$. One way to explain this is to point out that $\Spa \A_\inf(\C_p)\setminus \{p[p^\flat]=0\}$ is not a perfectoid space, but only becomes perfectoid after a further base change along e.g.\ $\Spa\Q_p^\cycl\rightarrow\Spa\Q_p$ or $\Spa\C_p\rightarrow\Spa\Q_p$; however, 
\begin{equation*}
\widehat{\GSpec \Q_p}\neq \GSpec\Q_p\;,
\end{equation*}
see \cite[Ex.\ 4.5.6]{dRFF}. Said differently, while the image of 
\begin{equation*}
Y_X^\diamond=X\times_{\Spd\F_p} \Spd\Q_p
\end{equation*}
under the left-adjoint of the diamond functor is computed via arc-descent to perfectoids on both $X$ and $\Spd\Q_p$, the definition of $Y_X$ as above only involves arc-descent to perfectoids on $X$.
\end{rem}

Recall from above that the definition of $\mathbb{B}_{[r, s]}(\ol{A})$ for a perfectoid $\ol{A}$ depends on the chosen perfectoid pseudouniformiser $\pi$. However, note that, for perfectoids over $\Q_p(p^{1/p^\infty})$, there is a somewhat canonical choice: we can take $\pi=p$. For any such perfectoid $\ol{A}$, recall that the algebras $\mathbb{B}_{[r, s]}(\ol{A})$ are idempotent on $Y_{\ol{A}}$ by \cite[Lem.\ 4.4.9]{dRFF} and hence we obtain a canonical radius map
\begin{equation*}
\kappa: Y_{\ol{A}}\rightarrow (0, \infty)
\end{equation*}
with the property that the preimage of any interval $[r, s]\subseteq (0, \infty)$ is given by
\begin{equation*}
\mathbb{B}_{[r, s]}^\dagger(\ol{A})\coloneqq \colim_{[r, s]\subseteq [r', s']} \mathbb{B}_{[r', s']}(\ol{A})\;.
\end{equation*}
In fact, we can even make sense of this map when $\ol{A}$ is not a $\Q_p(p^{1/p^\infty})$-algebra: In that case, we can just define the map by descent thanks to \cite[Cor.\ 4.4.10]{dRFF}; indeed, for any two choices of $p^\flat=(p, p^{1/p}, p^{1/p^2}, \dots)$, the corresponding elements $[p^\flat]$ will differ by a power of $[\epsilon]=[(1, \zeta_p, , \zeta_{p^2}, \dots)]$, which has norm $1$. 

From this canonical radius map $\kappa: Y_{\ol{A}}\rightarrow (0, \infty)$, we also obtain a canonical definition of $\mathbb{B}^\dagger_{[r, s]}(\ol{A})$ for all perfectoids $\ol{A}$: We can just define
\begin{equation*}
\mathbb{B}^\dagger_{[r, s]}(\ol{A})\coloneqq \O(Y_{\ol{A}, [r, s]})\;,
\end{equation*}
where $Y_{\ol{A}, [r, s]}$ is the preimage of $[r, s]\subseteq (0, \infty)$ under $\kappa$. Thus, from now on, whenever we refer to $\mathbb{B}^\dagger_{[r, s]}(\ol{A})$ or the radius map $\kappa: Y_{\ol{A}}\rightarrow (0, \infty)$ for a perfectoid $\ol{A}$, we will always assume that they are defined with respect to $\pi=p$. In particular, this has the consequence that the image of $\iota: \GSpec\ol{A}\rightarrow Y_{\ol{A}}$ is supported above $\{1\}\subseteq (0, \infty)$ under the radius map.

Let us now define $Y_A$ for arbitrary nilperfectoids $A$; in fact, we only define $Y_A$ for totally disconnected nilperfectoids, but one could then define $Y_A$ for all nilperfectoids (and in fact all Gelfand stacks) using the descent result from \cref{lem:prism-descentya} below. For totally disconnected nilperfectoids, we first define the algebras $\mathbb{B}^\dagger_{[p^{-r}, p^r]}(A)$ for any $r\geq 0$ and then glue their associated Gelfand stacks.

\begin{defi}
Let $A$ be a totally disconnected nilperfectoid. For any $r\geq 0$, the Gelfand ring $\mathbb{B}^\dagger_{[p^{-r}, p^r]}(A)$ is defined as the pullback
\begin{equation}
\label{eq:defis-pullbackya}
\begin{tikzcd}
\mathbb{B}^\dagger_{[p^{-r}, p^r]}(A)\ar[r]\ar[d] & \mathbb{B}^\dagger_{[p^{-r}, p^r]}(\ol{A})\ar[d, "\prod_n \theta\circ F^{-n}"] \\
\prod_{n=0}^r A\ar[r] & \prod_{n=0}^r \ol{A}\nospacepunct{\;.}
\end{tikzcd}
\end{equation}
The \emph{$Y$-curve} of $A$ is defined as the Gelfand stack
\begin{equation*}
Y_A\coloneqq \colim_r \GSpec \mathbb{B}^\dagger_{[p^{-r}, p^r]}(A)\;.
\end{equation*}
\end{defi}

\begin{rem}
\label{rem:prism-uperfdya}
We could also make the above definition more generally for uniformly totally disconnected perfectoid rings $A$, i.e.\ Gelfand rings $A$ whose uniform completion $A^u$ is totally disconnected perfectoid, by replacing $\ol{A}$ by $A^u$. However, the total disconnectedness cannot be removed and is necessary for the descent statement in \cref{lem:prism-descentya} below.
\end{rem}

\begin{figure}

\begin{center}
\begin{tikzpicture}
  \def\r{9}
  \def\shortr{5} 
  \def\longr{9.5}
  \def\vlongr{10.5}

  \draw[thick] (\r,0) arc[start angle=0, end angle=90, radius=\r];
  
  \draw[->, thick] (0.866*\vlongr, 1/2*\vlongr) arc[start angle=30, end angle=60, radius=\vlongr];
  
  \node at (1.4142/2*\vlongr + 0.3, 1.4142/2*\vlongr + 0.3) {$\varphi$};

  \def\shorttr{\r-0.2}
  \def\longgr{\r+0.2}

  \draw[thick] (\shorttr,0.2) arc (-180:0:0.2);
  \node at (\longr, 0) {$0$};
  
  \draw[thick] (7.5:\shorttr) -- (7.5:\longgr);
  \node at (7.5:\longr+0.3) {$1/p^3$};

  \draw[thick] (15:\shorttr) -- (15:\longgr);
  \node at (15:\longr+0.3) {$1/p^2$};

  \draw[thick] (25:\shorttr) -- (25:\longgr);
  \node at (25:\longr+0.2) {$1/p$};

  \draw[thick] (45:\shorttr) -- (45:\longgr);
  \node at (45:\longr) {$1$};

  \draw[thick] (65:\shorttr) -- (65:\longgr);
  \node at (65:\longr) {$p$};

  \draw[thick] (77.5:\shorttr) -- (77.5:\longgr);
  \node at (76.5:\longr+0.1) {$p^2$};
  
  \draw[thick] (85:\shorttr) -- (85:\longgr);
  \node at (84:\longr+0.1) {$p^3$};

  \draw[thick] (0.2, \longgr) arc (90:270:0.2);
  \node at (0, \longr) {$\infty$};

  \begin{scope}

  \clip (0,0) rectangle (\longr,\longr);
  
  \begin{scope}[transform canvas={rotate=45}]
        \shade[shading=axis, bottom color=white, top color=white, middle color=green, shading angle=0] (0.2, -0.1) rectangle (\shortr, 0.1);
  \end{scope} 

  \begin{scope}[transform canvas={rotate=65}]
        \shade[shading=axis, bottom color=white, top color=white, middle color=green, shading angle=0] (0.2, -0.1) rectangle (\shortr, 0.1);
  \end{scope} 

  \begin{scope}[transform canvas={rotate=77.5}]
        \shade[shading=axis, bottom color=white, top color=white, middle color=green, shading angle=0] (0.2, -0.1) rectangle (\shortr, 0.1);
  \end{scope}

  \begin{scope}[transform canvas={rotate=85}]
        \shade[shading=axis, bottom color=white, top color=white, middle color=green, shading angle=0] (0.2, -0.1) rectangle (\shortr, 0.1);
  \end{scope}

  \draw[line width=0.2mm, color=blue] (7.5:0.2) -- (7.5:\shortr);

  \draw[line width=0.2mm, color=blue] (15:0.2) -- (15:\shortr);

  \draw[line width=0.2mm, color=blue] (25:0.2) -- (25:\shortr);

  \end{scope}
  
  \def\shortishr{\shortr+0.5}
  
  \node at (7.5:\shortishr+0.3) {$\GSpec \ol{A}$};
  
  \node at (14.5:\shortishr+0.3) {$\GSpec \ol{A}$};
  
  \node at (24.5:\shortishr+0.2) {$\GSpec\ol{A}$};
  
  \node at (44:\shortishr+0.1) {$\GSpec A$};
  
  \node at (64:\shortishr) {$\GSpec A$};
  
  \node at (75.5:\shortishr) {$\GSpec A$};
  
  \node at (83:\shortishr+0.4) {$\GSpec A$};
  
  \draw[thick, dotted] (0,0) -- (\shortr,0);
  \draw[dotted, thick] (0,0) -- (0,\shortr); 
  
  \draw[->, thick] (45:\shortishr+1) -- (45:\shorttr-0.5);
  \node at (48:7.5) {$\kappa$};

\end{tikzpicture}
\end{center}

\captionsetup{justification=centering}
\caption{A schematic picture of $Y_A$}
\label{fig:ya}
\end{figure}

We immediately record some structures on and properties of $Y_A$.

\begin{rem}
\label{rem:prism-properties}
\begin{enumerate}[label=(\alph*)]
\item If $A$ is perfectoid, then $Y_A=Y_{\ol{A}}$, so the notation is consistent. Indeed, in that case, the lower horizontal map in (\ref{eq:defis-pullbackya}) is an isomorphism and hence the upper horizontal map is an isomorphism as well.

\item Note that the map $\mathbb{B}^\dagger_{[p^{-r}, p^r]}(A)\rightarrow \mathbb{B}^\dagger_{[p^{-r}, p^r]}(\ol{A})$ is a $\dagger$-nilpotent thickening for any $r\geq 0$ since the same is true for the bottom row in (\ref{eq:defis-pullbackya}). In particular, for any nilperfectoid $A$, the ring $\mathbb{B}^\dagger_{[p^{-r}, p^r]}(A)$ is indeed Gelfand. Moreover, the Gelfand stack $Y_A$ is a derived Berkovich space in the sense of \cite[§4.3]{dRFF} and the natural map $Y_{\ol{A}}\rightarrow Y_A$ induces an isomorphism on underlying topological spaces.

\item There is a natural Frobenius morphism $\phi: Y_A\rightarrow Y_A$ induced by the Frobenius morphism of $Y_{\ol{A}}$ and the left shift
\begin{equation*}
\prod_{n=0}^r A\rightarrow \prod_{n=0}^{r-1} A\;.
\end{equation*}

\item Note that, by definition, there is a projection map $\theta: \mathbb{B}^\dagger_{[1, 1]}(A)\rightarrow A$ inducing a map
\begin{equation*}
\iota: \GSpec A\rightarrow Y_A\;.
\end{equation*}
As the diagram 
\begin{equation*}
\begin{tikzcd}
\GSpec\ol{A}\ar[r]\ar[d, "\iota", swap] & \GSpec A\ar[d, "\iota"] \\
Y_{\ol{A}}\ar[r] & Y_A
\end{tikzcd}
\end{equation*}
commutes, there is no clash of notation and the above recovers the previously defined map $\iota: \GSpec \ol{A}\rightarrow Y_{\ol{A}}$ if $A$ is already perfectoid. However, let us warn the reader that, contrary to what happens in the perfectoid case, the map $\iota: \GSpec A\rightarrow Y_A$ does \emph{not} generally cut out a divisor anymore.

\item The radius map $\kappa: Y_{\ol{A}}\rightarrow (0, \infty)$ extends to $Y_A$. Indeed, as $Y_A$ is a $\dagger$-nilpotent thickening of $Y_{\ol{A}}$, they have the same de Rham stack and hence we obtain a map
\begin{equation*}
\kappa: Y_A\rightarrow Y_A^\dR\cong Y_{\ol{A}}^\dR\xrightarrow{\kappa^\dR} (0, \infty)^\dR=(0, \infty)\;,
\end{equation*}
where the last step uses the idempotency of the functor $(-)^\dR$ and \cref{ex:recall-drbetti}. We note that $\GSpec \mathbb{B}^\dagger_{[p^{-r}, p^r]}(A)$ is exactly the preimage $Y_{A, [p^{-r}, p^r]}$ of the interval $[p^{-r}, p^r]\subseteq (0, \infty)$ under the map $\kappa$. We can then furthermore define $\mathbb{B}^\dagger_{[r, s]}(A)$ for any interval $[r, s]\subseteq (0, \infty)$ as the ring of global functions on the preimage of $[r, s]\subseteq (0, \infty)$ under $\kappa$. \qedhere
\end{enumerate}
\end{rem}

We have the following descent result for the association $A\mapsto Y_A$, which allows us to define $Y_A$ for arbitrary nilperfectoids $A$, and in fact for all Gelfand stacks.

\begin{lem}
\label{lem:prism-descentya}
For any $r\geq 0$, the functor $A\mapsto \mathbb{B}^\dagger_{[p^{-r}, p^r]}(A)$ from totally disconnected nilperfectoids to Gelfand rings preserves $!$-equivalences. In particular, the functor $A\mapsto Y_A$ from totally disconnected nilperfectoids to Gelfand rings preserves $!$-equivalences.
\end{lem}
\begin{proof}
The claim is true for the functor $A\mapsto \mathbb{B}^\dagger_{[p^{-r}, p^r]}(\ol{A})$ by \cite[Cor.\ 4.4.10.(1)]{dRFF} and for the functor $A\mapsto \ol{A}$ by \cite[Lem.\ 4.4.3]{dRFF}. In other words, for any $!$-equivalence $A\rightarrow A_\bullet$ of totally disconnected nilperfectoids, the pro-system $(\Tot_{\leq k} \mathbb{B}^\dagger_{[p^{-r}, p^r]}(\ol{A}_\bullet))_k$ is pro-constant and equal to $\mathbb{B}^\dagger_{[p^{-r}, p^r]}(\ol{A})$ in $\D(\mathbb{B}^\dagger_{[p^{-r}, p^r]}(\ol{A}))$ and, similarly, the pro-system $(\Tot_{\leq k} \Nil^\dagger(A_\bullet))_k$ is pro-constant and equal to $\Nil^\dagger(A)$ in $\D(A)$. From this one deduces the result using the exact sequence
\begin{equation*}
\begin{tikzcd}
0\ar[r] & \prod_{n=0}^r \Nil^\dagger(A)\ar[r] & \mathbb{B}^\dagger_{[p^{-r}, p^r]}(A)\ar[r] & \mathbb{B}^\dagger_{[p^{-r}, p^r]}(\ol{A})\ar[r] & 0\nospacepunct{\;,}
\end{tikzcd}
\end{equation*}
which is immediate from the definition.
\end{proof}

\begin{rem}
\label{rem:prism-descentuperfd}
Note that the above proof works just as well if we only assume $A$ to be uniformly totally disconnected perfectoid and use the extended definition of $Y_A$ from \cref{rem:prism-uperfdya}. Indeed, we just have to replace $\Nil^\dagger(A_\bullet)$ by $\fib(A_\bullet\rightarrow A_\bullet^u)$ in that case.
\end{rem}

The previous lemma guarantees that the association $A\mapsto Y_A$ from totally disconnected nilperfectoids $A$ to Gelfand stacks extends to a unique colimit-preserving functor
\begin{equation*}
\begin{tikzcd}
Y_{(-)}: \GelfStk\rightarrow\GelfStk\;, \hspace{0.3cm} X\mapsto Y_X
\end{tikzcd}
\end{equation*}
and we write $Y_A\coloneqq Y_{\GSpec A}$ for any Gelfand ring $A$. All the properties and structures from \cref{rem:prism-properties} extend to $Y_X$ for arbitrary Gelfand stacks $X$ by descent and by descent from the case of totally disconnected nilperfectoids we obtain a chain of maps
\begin{equation*}
\widehat{Y}_{X^\diamond}\rightarrow Y_X\rightarrow Y_X^\dR\cong Y_{X^\diamond}^\dR\;.
\end{equation*}

The association $A\mapsto Y_A$ also preserves finite étale maps and open embeddings.

\begin{lem}
\label{lem:prism-finetya}
Let $A\rightarrow B$ be a finite étale map of totally disconnected nilperfectoids. Then the induced map $Y_B\rightarrow Y_A$ is finite étale as well, i.e.\ $\mathbb{B}^\dagger_{[r, s]}(B)$ is a finite étale $\mathbb{B}^\dagger_{[r, s]}(A)$-algebra for all $[r, s]\subseteq (0, \infty)$.
\end{lem}
\begin{proof}
It suffices to check the claim for intervals of the form $[p^{-r}, p^r]$. Then first observe that $B\tensor_A \ol{A}$ is a finite étale $\ol{A}$-algebra and hence perfectoid. As the map $B\rightarrow B\tensor_A \ol{A}$ is surjective and its kernel is contained in $\Nil^\dagger(B)$, we conclude that $B\tensor_A \ol{A}\cong \ol{B}$. In particular, $\ol{B}$ is a finite étale $\ol{A}$-algebra and this implies that $\mathbb{B}^\dagger_{[p^{-r}, p^r]}(\ol{B})$ is a finite étale $\mathbb{B}^\dagger_{[p^{-r}, p^r]}(\ol{A})$-algebra. 

We now first prove that $\mathbb{B}^\dagger_{[p^{-r}, p^r]}(B)$ is a finite projective $\mathbb{B}^\dagger_{[p^{-r}, p^r]}(A)$-module. For this, after localising on $\cal{M}(\mathbb{B}^\dagger_{[p^{-r}, p^r]}(A))\cong \cal{M}(\mathbb{B}^\dagger_{[p^{-r}, p^r]}(\ol{A}))$, we may assume that $\mathbb{B}^\dagger_{[p^{-r}, p^r]}(\ol{B})$ is a free $\mathbb{B}^\dagger_{[p^{-r}, p^r]}(\ol{A})$-module. Via the maps $\theta\circ F^{-n}: \mathbb{B}^\dagger_{[p^{-r}, p^r]}(\ol{B})\rightarrow\ol{B}$ for $0\leq n\leq r$, this yields bases of $\ol{B}$ as an $\ol{A}$-module, any lift of which will yield bases of $B$ as an $A$-module: indeed, isomorphisms of finite projective modules may be checked after base change to the $\dagger$-reduction as the $\dagger$-nilradical is contained in the Jacobson radical by \cite[Prop.\ 2.6.16.(2)]{dRStack}. By definition, we then conclude that $\mathbb{B}^\dagger_{[p^{-r}, p^r]}(B)$ is a free $\mathbb{B}^\dagger_{[p^{-r}, p^r]}(A)$-module. Moreover, this argument also shows that the canonical map $\theta: \mathbb{B}^\dagger_{[p^{-r}, p^r]}(B)\rightarrow B$ induces an isomorphism
\begin{equation}
\label{eq:prism-basechangebalongtheta}
\mathbb{B}^\dagger_{[p^{-r}, p^r]}(B)\tensor_{\mathbb{B}^\dagger_{[p^{-r}, p^r]}(A)} A\cong B\;.
\end{equation}

Finally, given the finite projectivity of $\mathbb{B}^\dagger_{[p^{-r}, p^r]}(B)(*)$ as an $\mathbb{B}^\dagger_{[p^{-r}, p^r]}(A)(*)$-module, we may check whether it is étale after base change to any residue field of $\mathbb{B}^\dagger_{[p^{-r}, p^r]}(A)(*)$ and, in particular, after base change to the $\dagger$-reduction $\mathbb{B}^\dagger_{[p^{-r}, p^r]}(\ol{A})$ as the $\dagger$-nilradical is contained in the Jacobson radical. However, now we are done if we can show that the natural map induces an isomorphism
\begin{equation*}
\mathbb{B}^\dagger_{[p^{-r}, p^r]}(B)\tensor_{\mathbb{B}^\dagger_{[p^{-r}, p^r]}(A)} \mathbb{B}^\dagger_{[p^{-r}, p^r]}(\ol{A})\cong \mathbb{B}^\dagger_{[p^{-r}, p^r]}(\ol{B})\;.
\end{equation*}

To check this last isomorphism, note that the right-hand side is the quotient of $\mathbb{B}^\dagger_{[p^{-r}, p^r]}(B)$ by $\prod_{n=0}^r \Nil^\dagger(B)$ while the left-hand side is the quotient of $\mathbb{B}^\dagger_{[p^{-r}, p^r]}(B)$ by 
\begin{equation*}
\begin{split}
\mathbb{B}^\dagger_{[p^{-r}, p^r]}(B)\tensor_{\mathbb{B}^\dagger_{[p^{-r}, p^r]}(A)} \prod_{n=0}^r \Nil^\dagger(A)&\cong \left(\mathbb{B}^\dagger_{[p^{-r}, p^r]}(B)\tensor_{\mathbb{B}^\dagger_{[p^{-r}, p^r]}(A)} \prod_{n=0}^r A\right)\tensor_{\prod_{n=0}^r A} \prod_{n=0}^r \Nil^\dagger(A) \\
&\cong \prod_{n=0}^r B\tensor_{\prod_{n=0}^r A} \prod_{n=0}^r \Nil^\dagger(A)\cong \prod_{n=0}^r \Nil^\dagger(B)\;,
\end{split}
\end{equation*}
where the penultimate isomorphism is (\ref{eq:prism-basechangebalongtheta}) while the last one just amounts to $B\tensor_A \Nil^\dagger(A)\cong \Nil^\dagger(B)$, which in turn follows from $B\tensor_A \ol{A}\cong \ol{B}$ by passing to fibres of the canonical maps from $B$.
\end{proof}

\begin{lem}
\label{lem:prism-openlocya}
Let $U\rightarrow\GSpec A$ be an open embedding of derived Berkovich spaces with $A$ totally disconnected nilperfectoid. Then the induced map $Y_U\rightarrow Y_A$ is an open embedding as well. More precisely, there is a pullback square
\begin{equation*}
\begin{tikzcd}
Y_U\ar[r]\ar[d] & Y_A\ar[d] \\
{|U|}\ar[r] & \cal{M}(A)\nospacepunct{\;,}
\end{tikzcd}
\end{equation*}
where $Y_A\rightarrow \cal{M}(A)$ is given by $Y_A\rightarrow Y_A^\dR\cong Y_{\ol{A}}^\dR\rightarrow \cal{M}(\ol{A})\cong \cal{M}(A)$ for the map $Y_{\ol{A}}\rightarrow\cal{M}(\ol{A})$ coming from \cite[Lem.\ 4.4.5]{dRFF}.
\end{lem}
\begin{proof}
It suffices to check the claim for open embeddings of the form $U=\{|f|<1\}$ and $U=\{|f|>1\}$ for $f\in A$. Moreover, since $\cal{M}(A)\cong\cal{M}(\ol{A})\cong\cal{M}(\ol{A}^\flat)$, we may assume that $f$ is a lift of an element $f^{\flat\sharp}\in\ol{A}$ for some $f^\flat\in\ol{A}^\flat$. We only treat the case $U=\{|f|<1\}$, the other one is completely analogous. Moreover, we may assume that $A$ is a $\Q_p(p^{1/p^\infty})$-algebra by descent.

To begin the argument in this case, let $B_\epsilon\coloneqq A\langle f\rangle_{\leq p^{-\epsilon}}$ and note that $B_\epsilon^u\cong \ol{A}\langle p^{-\epsilon} f\rangle$ is perfectoid with tilt $\ol{A}\langle (p^\flat)^{-\epsilon}f^\flat\rangle$, hence $\mathbb{B}^\dagger_{[r, s]}(B_\epsilon^u)=\mathbb{B}^\dagger_{[r, s]}(\ol{A})\langle [(p^\flat)^{-\epsilon}f^\flat]\rangle$. Furthermore, letting $f_n$ denote an arbitrary lift of $(f^{\flat\sharp})^{1/p^n}$ to $A$ and $f_0\coloneqq f$, observe that $B_\epsilon=A\langle f_n\rangle_{\leq p^{-\epsilon/p^n}}$ by \cite[Cor.\ 2.2.14]{dRFF}. Hence, using \cref{rem:prism-descentuperfd}, we know that $\mathbb{B}^\dagger_{[p^{-r}, p^r]}(B_\epsilon)$ sits in a pullback square
\begin{equation*}
\begin{tikzcd}
\mathbb{B}^\dagger_{[p^{-r}, p^r]}(B_\epsilon)\ar[r]\ar[d] & \mathbb{B}^\dagger_{[p^{-r}, p^r]}(\ol{A})\langle (p^\flat)^{-\epsilon}f^\flat\rangle\ar[d] \\
\prod_{n=0}^r A\langle f_n\rangle_{\leq p^{-\epsilon/p^n}}\ar[r] & \prod_{n=0}^r \ol{A}\langle (p^\flat)^{-\epsilon/p^n}(f^{\flat\sharp})^{1/p^n}\rangle\;.
\end{tikzcd}
\end{equation*}

In particular, letting $f^\flat\in \mathbb{B}^\dagger_{[p^{-r}, p^r]}(A)$ denote the gluing of $f^\flat\in\mathbb{B}^\dagger_{[p^{-r}, p^r]}(\ol{A})$ and $(f_0, f_1, \dots, f_r)\in \prod_{n=0}^r A$, we see that there are canonical maps
\begin{equation*}
\mathbb{B}^\dagger_{[p^{-r}, p^r]}(A)\langle (p^\flat)^{-\epsilon} f^\flat\rangle\rightarrow \mathbb{B}^\dagger_{[p^{-r}, p^r]}(B_\epsilon)\rightarrow \mathbb{B}^\dagger_{[p^{-r}, p^r]}(A)\langle f^\flat\rangle_{\leq p^{-\epsilon}}\;,
\end{equation*} 
from which we conclude that
\begin{equation*}
Y_{U, [p^{-r}, p^r]}\cong\colim_{\epsilon>0} Y_{B_\epsilon, [p^{-r}, p^r]}=\colim_{\epsilon>0} \GSpec \mathbb{B}^\dagger_{[p^{-r}, p^r]}(B_\epsilon)\cong \colim_{\epsilon>0} \GSpec \mathbb{B}^\dagger_{[p^{-r}, p^r]}(A)\langle f^\flat\rangle_{\leq p^{-\epsilon}}
\end{equation*}
by cofinality, and the right-hand side is exactly $\{|f^\flat|<1\}\subseteq Y_{A, [p^{-r}, p^r]}$. Taking the colimit over $r$, the claim follows.
\end{proof}

Intuitively, one should think of $Y_A$ as obtained by replacing the image of $\iota: \GSpec\ol{A}\rightarrow Y_{\ol{A}}$ by the thickening $\GSpec A$ and then freely adding Frobenius translates of this thickening, see \cref{fig:ya}. This is made precise by the following proposition.

\begin{prop}
\label{prop:defis-pushoutberk}
For any totally disconnected nilperfectoid $A$, the diagram
\begin{equation*}
\begin{tikzcd}[column sep=large]
\bigsqcup_{n\geq 0} \GSpec\ol{A}\ar[r, "\sqcup_{n\geq 0} \phi^n\circ\iota"]\ar[d] & Y_{\ol{A}}\ar[d] \\
\bigsqcup_{n\geq 0} \GSpec A\ar[r, "\sqcup_{n\geq 0} \phi^n\circ\iota"] & Y_A
\end{tikzcd}
\end{equation*}
is a pushout in the category of derived Berkovich spaces. In particular, the map $Y_{\ol{A}}\rightarrow Y_A$ is an isomorphism away from $\{1, p, p^2, \dots\}\subseteq (0, \infty)$.
\end{prop}
\begin{proof}
The second assertion follows from the first since the image of the $n$-th Frobenius translate of $\GSpec A$ in $Y_A$ is supported over $\{p^n\}\subseteq (0, \infty)$ under the radius map. To prove the first assertion, we check the universal property, i.e.\ that the natural map
\begin{equation}
\label{eq:prism-yapushout}
\Hom(Y_A, X)\rightarrow \Hom(Y_{\ol{A}}, X)\times_{\prod_{n\geq 0} X(\ol{A})} \prod_{n\geq 0} X(A)
\end{equation}
is an equivalence for any derived Berkovich space $X$. For this, first note that the diagram in question induces a pushout on underlying topological spaces: Indeed, this is because both $Y_A$ and $Y_{\ol{A}}$ as well as $\GSpec A$ and $\GSpec\ol{A}$ have the same underlying topological space. Thus, it suffices to check that (\ref{eq:prism-yapushout}) separately induces equivalences on the fibres of the left- and right-hand side over $\Hom_\cont(|Y_A|, |X|)$.

Thus, fix a map $|Y_A|\rightarrow |X|$ of topological spaces. Taking an affine rational cover of $X$ whose pullback along this map can be refined by a strict rational cover of $|Y_A|$, we can reduce to the case where $X=\GSpec B$ is affine. In that case, however, the claim follows from the fact that each $\O(Y_{A, [p^{-r}, p^r]})=\mathbb{B}^\dagger_{[p^{-r}, p^r]}(A)$ is by definition the pullback of the rings of global functions on $Y_{\ol{A}, [p^{-r}, p^r]}$ and $\bigsqcup_{n=0}^r \GSpec A$ over the ring of global functions on $\bigsqcup_{n=0}^r \GSpec\ol{A}$.
\end{proof}

From the definition of $Y_A$ and the picture in \cref{fig:ya}, one should have the intuition that the Frobenius $\phi: Y_A\rightarrow Y_A$ is an isomorphism away from from the image of 
\begin{equation*}
\GSpec\ol{A}\xrightarrow{\iota} Y_{\ol{A}}\xrightarrow{\phi^{-1}} Y_{\ol{A}}\rightarrow Y_A\;,
\end{equation*}
which we will at times just denote by $\phi^{-1}\circ\iota: \GSpec\ol{A}\rightarrow Y_A$ in the future. This is made precise by the following lemma.

\begin{lem}
\label{lem:defis-frobpushout}
For any totally disconnected nilperfectoid $A$, the diagram
\begin{equation*}
\begin{tikzcd}
\GSpec\ol{A}\ar[r]\ar[d, "\phi^{-1}\circ\iota", swap] & \GSpec A\ar[d, "\iota"] \\
Y_A\ar[r, "\phi"] & Y_A
\end{tikzcd}
\end{equation*}
is a pushout in the category of derived Berkovich spaces.
\end{lem}
\begin{proof}
Note that Frobenius induces an isomorphism $\mathbb{B}^\dagger_{[p^{-r}, p^r]}(\ol{A})\xrightarrow{\cong} \mathbb{B}^\dagger_{[p^{-r-1}, p^{r-1}]}(\ol{A})$ for any $r\geq 0$. As the shift map $\prod_{n=0}^r A\rightarrow \prod_{n=0}^{r-1} A$ is surjective, we conclude that
\begin{equation*}
F: \mathbb{B}^\dagger_{[p^{-r}, p^r]}(A)\rightarrow\mathbb{B}^\dagger_{[p^{-r-1}, p^{r-1}]}(A)
\end{equation*}
is surjective as well and that the kernel agrees with the kernel of the induced map on $\dagger$-nilradicals. Since this is the shift map $\prod_{n=0}^r \Nil^\dagger(A)\rightarrow \prod_{n=0}^{r-1} \Nil^\dagger(A)$, we conclude that there is an exact sequence
\begin{equation*}
\begin{tikzcd}
0\ar[r] & \Nil^\dagger(A)\ar[r] & \mathbb{B}^\dagger_{[p^{-r}, p^r]}(A)\ar[r, "F"]\ar[r] & \mathbb{B}^\dagger_{[p^{-r-1}, p^{r-1}]}(A)\ar[r] & 0\nospacepunct{\;,}
\end{tikzcd}
\end{equation*}
which in turn implies that the diagram
\begin{equation*}
\begin{tikzcd}
\mathbb{B}^\dagger_{[p^{-r}, p^r]}(A)\ar[r, "F"]\ar[d, "\theta", swap] & \mathbb{B}^\dagger_{[p^{-r-1}, p^{r-1}]}(A)\ar[d] \\
A\ar[r] & \ol{A}
\end{tikzcd}
\end{equation*}
is cartesian. Now one deduces the desired pushout by the same argument as in the proof of \cref{prop:defis-pushoutberk}.
\end{proof}

\subsection{Analytic prismatisation}

Recall the notion of degree $1$ Cartier divisors on $Y_{\ol{A}}$ from \cite[Def.\ II.1.19]{FarguesScholze}: These are defined to be exactly those Cartier divisors arising as maps
\begin{equation*}
\GSpec \ol{A}^{\flat\sharp}\rightarrow Y_{\ol{A}}
\end{equation*}
for an untilt $\ol{A}^{\flat\sharp}$ of $\ol{A}^\flat$. Moreover, they are always cut out by an element of the form $p+a[\pi^\flat]$ for some $a\in W(\ol{A}^{\flat\circ})$ and some perfectoid pseudouniformiser $\pi\in\ol{A}$; any such element is called \emph{distinguished}. Also note that, for us, a Cartier divisor on an arbitrary Gelfand stack is defined as a map $X\rightarrow \A^1/\G_m$. By abuse of notation, we will often simply say that $D\subseteq X$ is a Cartier divisor, where $D$ denotes the pullback
\begin{equation*}
\begin{tikzcd}
D\ar[r]\ar[d] & X\ar[d] \\
*/\G_m\ar[r] & \A^1/\G_m\nospacepunct{\;.}
\end{tikzcd}
\end{equation*}
Let us point out that the above is also an abuse of notation in a second sense: the map $D\rightarrow X$ is \emph{not} an injection, i.e.\ its diagonal is not an isomorphism (note that this already happens for usual effective Cartier divisors when one works in the setting of derived algebraic geometry).

\begin{defi}
Let $A$ be totally disconnected nilperfectoid. A Cartier divisor $D\subseteq Y_A$ is called \emph{of degree $1$} if its pullback $D\times_{Y_A} Y_{\ol{A}}$ is a degree $1$ Cartier divisor in $Y_{\ol{A}}$.
\end{defi}

\begin{rem}
Note that any degree $1$ Cartier divisor $D\subseteq Y_A$ is affine. Indeed, any degree $1$ Cartier divisor in $Y_{\ol{A}}$ is fully contained in $Y_{\ol{A}, [r, s]}$ for some $[r, s]\subseteq (0, \infty)$ and hence $D$ is already a Cartier divisor in $Y_{A, [r, s]}=\GSpec\mathbb{B}^\dagger_{A, [r, s]}$.
\end{rem}

\begin{ex}
\label{ex:defis-drdivisor}
The map $\phi^{-1}\circ\iota: \GSpec\ol{A}\rightarrow Y_A$ is a degree $1$ Cartier divisor. This is because the image of this composition is supported over $\{1/p\}\subseteq (0, \infty)$ under the radius map and $Y_A$ is isomorphic to $Y_{\ol{A}}$ over $(0, 1)$ by \cref{prop:defis-pushoutberk}.
\end{ex}

We quickly record a rigidity property of degree $1$ Cartier divisors on $Y_A$ which will occasionally be useful in the sequel.

\begin{lem}
\label{lem:prism-rigidity}
Let $A$ be totally disconnected nilperfectoid. Then any map $D_1\rightarrow D_2$ between degree $1$ Cartier divisors on $Y_A$ is an isomorphism.
\end{lem}
\begin{proof}
First recall that the statement is true after pullback to $Y_{\ol{A}}$. Indeed, we have
\begin{equation*}
D_1\times_{Y_A} Y_{\ol{A}}\cong \GSpec\ol{A}^{\flat\sharp_1}\;, \hspace{0.3cm} D_2\times_{Y_A} Y_{\ol{A}}\cong\GSpec\ol{A}^{\flat\sharp_2}
\end{equation*}
for two characteristic zero untilts $\ol{A}^{\flat\sharp_1}$ and $\ol{A}^{\flat\sharp_2}$ of $\ol{A}^\flat$ and the map $D_1\rightarrow D_2$ corresponds to a map $\ol{A}^{\flat\sharp_2}\rightarrow \ol{A}^{\flat\sharp_1}$. However, as perfectoid spaces over $\ol{A}^{\flat\sharp_2}$ are equivalent to perfectoid spaces over $\ol{A}^\flat$ via tilting and the map $\ol{A}^{\flat\sharp_2}\rightarrow \ol{A}^{\flat\sharp_1}$ induces an isomorphism on tilts, it must be an isomorphism.

Returning to the situation on $Y_A$, first observe that the given map induces a map $\O(-D_2)\rightarrow\O(-D_1)$ between the ideal sheaves. Now recall that $\O(-D_1)$ and $\O(-D_2)$ yield invertible modules over any $\mathbb{B}^\dagger_{[r, s]}(A)(*)$ by \cref{prop:recall-fredholm}, and hence checking whether $\O(-D_2)\rightarrow \O(-D_1)$ is an isomorphism may be done in each residue field separately by Nakayama's lemma. However, as $\mathbb{B}^\dagger_{[r, s]}(A)\rightarrow \mathbb{B}^\dagger_{[r, s]}(\ol{A})$ is a $\dagger$-thickening and $\Nil^\dagger(\mathbb{B}^\dagger_{[r, s]}(A))(*)$ is contained in the Jacobson radical of $\mathbb{B}^\dagger_{[r, s]}(A)(*)$ by \cite[Prop.\ 2.6.16.(2)]{dRStack}, this reduces the claim to the situation after pullback to $Y_{\ol{A}}$, which we have discussed in the previous paragraph.
\end{proof}

Let us also note another useful lemma with a similar proof.

\begin{lem}
\label{lem:prism-odtrivial}
Let $A$ be totally disconnected nilperfectoid. For any degree $1$ Cartier divisor $D\subseteq Y_A$, there is an isomorphism $\O(-D)\cong \O$.
\end{lem}
\begin{proof}
First recall from \cite[Prop.\ II.2.3]{FarguesScholze} that such an isomorphism exists after pullback to $Y_{\ol{A}}$. Moreover, as $D$ is supported on $Y_{A, [r, s]}$ for some interval $[r, s]\subseteq (0, \infty)$, we have $\O(-D)\cong \O$ on $Y_{A, (0, r')}$ and on $Y_{A, (s', \infty)}$ for $r'<r$ and $s'>s$, so it suffices to trivialise $\O(-D)$ on $Y_{A, [r'-\epsilon, s'+\epsilon]}$ for some $\epsilon>0$. 

In other words, we are done as soon as we show that $\O(-D)$ yields a trivial invertible module over any $\mathbb{B}^\dagger_{[r, s]}(A)(*)$ via \cref{prop:recall-fredholm}. However, this may be checked separately in each residue field by Nakayama's lemma and hence reduces to the claim after pullback to $\mathbb{B}^\dagger_{[r, s]}(\ol{A})(*)$ as in the proof of the previous lemma, where we already know it.
\end{proof}

We are now ready to define the prismatisation. For this, recall that totally disconnected nilperfectoids form a basis for the $!$-topology on Gelfand rings and hence any Gelfand stack is uniquely determined by its points on totally disconnected nilperfectoid rings.

\begin{defi}
The (rational) \emph{analytic prismatisation} of $\Q_p$ is the Gelfand stack $\Q_p^\prism$ defined by
\begin{equation*}
\Q_p^\prism(A)=\{\text{degree $1$ Cartier divisors $D\subseteq Y_A$}\}
\end{equation*}
for totally disconnected nilperfectoids $A$. For an arbitrary Gelfand stack $X$ over $\Q_p$, we define its (rational) analytic prismatisation as the Gelfand stack $X^\prism$ over $\Q_p^\prism$ given by
\begin{equation*}
X^\prism(\GSpec A\xrightarrow{D} \Q_p^\prism)=\{\text{maps $D\rightarrow X$}\}
\end{equation*}
for totally disconnected nilperfectoids $A$.
\end{defi}

Note that \cref{lem:prism-descentya} guarantees that sheafification does not change the value of $X^\prism$ on totally disconnected nilperfectoids. Moreover, the Frobenius on $Y_A$ induces a Frobenius $\phi: \Q_p^\prism\rightarrow\Q_p^\prism$ via pullback of divisors and hence also a Frobenius $\phi: X^\prism\rightarrow X^\prism$ for any $X$. Finally, observe that it is immediate from the definition that the association $X\mapsto X^\prism$ commutes with limits.

\begin{ex}
\label{ex:defis-drmap}
The degree $1$ divisor $\phi^{-1}\circ\iota: \GSpec\ol{A}\rightarrow Y_A$ from \cref{ex:defis-drdivisor} yields a map $i_\dR: \GSpec\Q_p\rightarrow\Q_p^\prism$ called the \emph{de Rham map}.
\end{ex}

Before we move on, let us describe the de Rham stack of $X^\prism$ because it is particularly easy.

\begin{prop}
\label{prop:defis-xprismdr}
Let $X$ be any Gelfand stack over $\Q_p$. Then there is a natural map
\begin{equation*}
Y_{X^\diamond}\rightarrow X^\prism
\end{equation*}
which induces an isomorphism on de Rham stacks.
\end{prop}
\begin{proof}
To construct the desired map, we have to give a map $Y_{\ol{A}}\rightarrow X^\prism$ for any map $\GSpec\ol{A}\rightarrow X$ from a totally disconnected perfectoid to $X$. As $X^\prism$ is defined by sheafification from its value on totally disconnected perfectoids, we have to give compatible degree $1$ divisors $D\subseteq Y_{\ol{B}}$ equipped with maps $D\rightarrow X$ for every totally disconnected perfectoid $\GSpec\ol{B}\rightarrow Y_{\ol{A}}$. However, for any such $\ol{B}$, the map $\GSpec\ol{B}\rightarrow Y_{\ol{A}}$ uniquely extends to a map $Y_{\ol{B}}\rightarrow Y_{\ol{A}}$ and then we can take $D$ to be the pullback of $\iota: \GSpec\ol{A}\rightarrow Y_{\ol{A}}$ along this map, which will automatically be equipped with a map to $X$ since $\GSpec\ol{A}$ is.

Now let us see that the above defines an isomorphism on de Rham stacks. For this, note that, for any totally disconnected perfectoid $\ol{A}$, an $\ol{A}$-point of $X^\prism$ amounts to a characteristic zero untilt $\ol{A}^{\flat\sharp}$ of $\ol{A}^\flat$ together with a map $\GSpec \ol{A}^{\flat\sharp}\rightarrow X$, which we may equivalently write as $\Spd\ol{A}^{\flat\sharp}\rightarrow X^\diamond$. However, this is exactly the same as an $\ol{A}$-point of $Y_{X^\diamond}^\diamond=X^\diamond\times_{\Spd\F_p}\Spd\Q_p$, so we are done.
\end{proof}

Thanks to the radius map $\kappa: Y_{X^\diamond}\rightarrow (0, \infty)$, the previous proposition allows us to equip $X^\prism$ with a radius map for any Gelfand stack $X$ via the composition
\begin{equation*}
\kappa: X^\prism\rightarrow (X^\prism)^\dR\cong Y_{X^\diamond}^\dR\xrightarrow{\kappa^\dR} (0, \infty)^\dR\cong (0, \infty)\;,
\end{equation*}
where the last step again uses idempotency of $(-)^\dR$ and \cref{ex:recall-drbetti}. As usual, we will denote the preimage of any interval $[r, s]\subseteq (0, \infty)$ under this map by $X^\prism_{[r, s]}$.

Moving on, observe that pulling back a degree $1$ divisor along $\iota: \GSpec A\rightarrow Y_A$ induces a map
\begin{equation*}
\mu: \Q_p^\prism\rightarrow \A^1/\G_m\;.
\end{equation*}
However, we can actually do slightly better.

\begin{lem}
\label{lem:defis-mu}
The map $\mu: \Q_p^\prism\rightarrow \A^1/\G_m$ can be refined to a map $\mu: \Q_p^\prism\rightarrow \overcirc{\DD}/\ol{\T}$, where $\overcirc{\DD}$ denotes the open unit disk.
\end{lem}
\begin{proof}
As being a unit on $Y_A$ may be checked after pullback to $Y_{\ol{A}}$ by \cite[Prop.\ 2.6.16.(2)]{dRStack} and the pullback of $D$ to $Y_{\ol{A}}$ is cut out by a distinguished element, we may choose a global section of $\O(D)$ cutting out $D$ which pulls back to a distinguished element on $Y_{\ol{A}}$. We claim that, after pullback along $\iota: \GSpec A\rightarrow Y_A$, all such global sections differ only by units of norm $1$ and have themselves norm less than $1$. This will then define the desired map $\Q_p^\prism\rightarrow \overcirc{\DD}/\ol{\T}$.

As an element having norm in $[r, s]\subseteq [0, \infty)$ may be checked after $\dagger$-reduction, it will be enough to show the following assertions:
\begin{enumerate}[label=(\alph*)]
\item Any distinguished element $p+a[\pi^\flat]$ in $W(\ol{A}^{\flat\circ})$ maps to a topologically nilpotent element under Fontaine's map $\theta: W(\ol{A}^{\flat\circ})\rightarrow \ol{A}^\circ$.
\item If two distinguished elements of $W(\ol{A}^{\flat\circ})$ differ by a global unit $u$ on $Y_{\ol{A}}$, the image of $u$ under Fontaine's map has norm $1$.
\end{enumerate}
For (a), note that the kernel of $\theta$ is generated by a distinguished element $p+b[\varpi^\flat]$ as well and the difference $(p+a[\pi^\flat])-(p+b[\varpi^\flat])=a[\pi^\flat]-b[\varpi^\flat]$ is already topologically nilpotent in $W(\ol{A}^{\flat\circ})$, whence the claim. For (b), we observe that two distinguished elements $p+a[\varpi^\flat]$ and $p+b[\varpi^\flat]$ in $W(\ol{A}^{\flat\circ})$ differing by a global unit $u$ on $Y_{\ol{A}}$ will be nonvanishing in a rational neighbourhoods of $V([\pi^\flat])$ and of $V(p)$ in $\Spa W(\ol{A}^{\flat\circ})$ as they cut out a degree $1$ divisor of $Y_{\ol{A}}$. Thus, we can even assume that $u$ is a unit in $W(\ol{A}^{\flat\circ})$ and this will of course imply that $\theta(u)\in (\ol{A}^{\flat\circ})^\times$, i.e.\ $\theta(u)$ has norm $1$.
\end{proof}

For any Gelfand stack $X$, the above yields a composite map
\begin{equation*}
\mu: X^\prism\rightarrow\Q_p^\prism\xrightarrow{\mu} \overcirc{\DD}/\ol{\T}\;.
\end{equation*}
The divisor cut out by this map plays a distinguished role.

\begin{defi}
Let $X$ be a Gelfand stack. The (analytic) \emph{Hodge--Tate stack} $X^\HT$ of $X$ is defined as the pullback
\begin{equation*}
\begin{tikzcd}
X^\HT\ar[r]\ar[d] & X^\prism\ar[d, "\mu"] \\
*/\ol{\T}\ar[r] & \overcirc{\DD}/\ol{\T}\nospacepunct{\;.}
\end{tikzcd}
\end{equation*}
It is a Cartier divisor in $X^\prism$ and we let $\cal{I}\in\Pic(X^\prism)$ denote its ideal sheaf. 
\end{defi}

By \cref{prop:defis-xprismxdiv1} below, the above indeed agrees with the analytic Hodge--Tate stack as defined by Anschütz--Le Bras--Rodríguez Camargo--Scholze in \cite{AnPrism}.

\begin{defi}
\label{defi:prism-bktwist}
For any Gelfand stack $X$, the \emph{Breuil--Kisin line bundle} $\O\{1\}$ on $X^\prism$ is defined by
\begin{equation*}
\O\{1\}\coloneqq \bigotimes_{n\geq 0} (\phi^n)^*\cal{I}\;.
\end{equation*}
\end{defi}

Note that the infinite tensor product occurring in the definition above makes sense: Indeed, the divisor $X^\HT\subseteq X^\prism$ is supported over $\{1\}$ under the radius map $\kappa: X^\prism\rightarrow (0, \infty)$ and hence its inverse images under Frobenius are supported over $1/p, 1/p^2, \dots$. In particular, over each $[r, s]\subseteq (0, \infty)$, only finitely many of the line bundles $(\phi^n)^*\cal{I}$ are nontrivial and hence the infinite tensor product above is actually a finite tensor product locally on $X^\prism$.

\subsection{Relation to the Hyodo--Kato stack and the perfect prismatisation}

Recall from \cref{prop:defis-xprismdr} that, for any Gelfand stack $X$, there is a natural map $Y_{X^\diamond}\rightarrow X^\prism$ and that it induces an isomorphism on de Rham stacks. Thus, understanding where this map is an isomorphism amounts to understanding where the natural map $X^\prism\rightarrow (X^\prism)^\dR$ is an isomorphism.

\begin{prop}
\label{prop:defis-prismffdr}
For any Gelfand stack $X$, the natural map 
\begin{equation*}
X^\prism\rightarrow (X^\prism)^\dR\cong Y_{X^\diamond}^\dR
\end{equation*}
is an isomorphism over the locus $(1, \infty)$ under the radius map.
\end{prop}
\begin{proof}
We claim that a map $\GSpec A\rightarrow X^\prism$ classifying a degree $1$ Cartier divisor $D\subseteq Y_A$ together with a map $D\rightarrow X$ factors through $X^\prism_{(1, \infty)}$ if and only if $D$ is supported on $Y_{A, (0, 1)}$. Once we know this, the claim follows: Indeed, the divisor $\iota: \GSpec \ol{A}\rightarrow Y_{\ol{A}}$ is supported over $1\in (0, \infty)$ and hence its Frobenius translates are supported over $p, p^2, p^3, \dots$. By \cref{prop:defis-pushoutberk}, this implies that the natural map induces an isomorphism $Y_{\ol{A}, (0, 1)}\cong Y_{A, (0, 1)}$.

To establish the remaining claim, first note that $D$ being supported on $Y_{A, (0, 1)}$ may be checked after pullback to $Y_{\ol{A}, (0, 1)}$ as units may be detected after $\dagger$-reduction by \cite[Prop.\ 2.6.16.(2)]{dRStack}. Similarly, by definition of the radius map, $\GSpec A\rightarrow X^\prism$ factoring through $X^\prism_{(1, \infty)}$ may be checked on the level of de Rham stacks by definition of the radius map. By descent for $Y_{X^\diamond}$, we may thus reduce to the case where $A=\ol{A}$ and $X=\GSpec\ol{B}$ is affinoid perfectoid and we may even assume that both $\ol{A}$ and $\ol{B}$ are $\Q_p(p^{1/p^\infty})$-algebras. 

Now $D$ corresponds to a characteristic zero untilt $\ol{A}^{\flat\sharp}$ of $\ol{A}^\flat$ and $D$ being supported on $Y_{\ol{A}, (0, 1)}$ translates to $|p|<|[p^\flat]|$ on $\GSpec \ol{A}^{\flat\sharp}\subseteq Y_{\ol{A}}$. Similarly, the image of $\GSpec \ol{A}\rightarrow Y_{\ol{B}}$ being supported over $(1, \infty)$ translates to $|[p^\flat]|<|p|$ after pullback along
\begin{equation*}
\GSpec\ol{A}\rightarrow Y_{\ol{A}}\cong Y_{\ol{A}^{\sharp\flat}}\rightarrow Y_{\ol{B}}\;,
\end{equation*}
where the last map is induced by the given map $\GSpec\ol{A}^{\sharp\flat}\rightarrow\GSpec\ol{B}$. Indeed, these two conditions match as the isomorphism $Y_{\ol{A}}\cong Y_{\ol{A}^{\sharp\flat}}$ exactly interchanges $p$ and $[p^\flat]$, so we are done.
\end{proof}

\begin{ex}
Recall that there is a natural map $\iota: X^\diamond\rightarrow Y_{X^\diamond}$ for any Gelfand stack $X$, the image of which is supported over $1\in (0, \infty)$ under the radius map. As Frobenius is an isomorphism on $Y_{X^\diamond}$, we also obtain a map $\phi\circ\iota: X^\diamond\rightarrow Y_{X^\diamond}$, whose image is now supported over $p$. Thus, by the previous proposition, there is an induced map
\begin{equation*}
X^\dR\xrightarrow{(\phi\circ\iota)^\dR} Y_{X^\diamond, (1, \infty)}^\dR\cong X^\prism_{(1, \infty)}\rightarrow X^\prism
\end{equation*}
and one checks that this map identifies with the map $i_\dR$ from \cref{ex:defis-drmap} by unwinding definitions. Since the map $X^\dR\rightarrow Y_{X^\diamond}^\dR$ is a closed embedding as $X^\diamond\rightarrow Y_{X^\diamond}$ is a Zariski-closed immersion, see \cite[Prop.\ 4.7.12.(1)]{dRFF}, we conclude that $i_\dR$ is a closed embedding.
\end{ex}

From the above we conclude that Frobenius induces an isomorphism $\phi: X^\prism_{(1, \infty)}\xrightarrow{\cong} X^\prism_{(p, \infty)}$. Moreover, quotienting out by Frobenius on this locus exactly yields the Hyodo--Kato stack of $X$, i.e.\ the natural map $X^\prism\rightarrow (X^\prism)^\dR$ induces an isomorphism
\begin{equation*}
\operatorname{coeq}(\hspace{-0.15cm}
\begin{tikzcd}
X^\prism_{(1, \infty)}\ar[r,shift left=.75ex, "\phi"]
  \ar[r,shift right=.75ex,swap, "\id"] & X^\prism_{(1, \infty)}
\end{tikzcd}
\hspace{-0.15cm})\cong X^\HK
\end{equation*}
for any Gelfand stack $X$.

Let us now also say something about $X^\prism$ over the locus $(0, p)$. As already announced in the beginning, this will be related to the stack $X^{\Div^1}$ from \cite{AnPrism}, which is the quotient of the perfect prismatisation from loc.\ cit.\ by Frobenius and whose definition we quickly recall.

\begin{defi}
For any totally disconnected nilperfectoid $A$, we define the \emph{modified Fargues--Fontaine curve} $\FF_A$ over $A$ as the pushout
\begin{equation*}
\begin{tikzcd}
\GSpec \ol{A}\ar[r]\ar[d] & \FF_{\ol{A}}\ar[d] \\
\GSpec A\ar[r] & \FF_A
\end{tikzcd}
\end{equation*}
in the category of derived Berkovich spaces, where $\FF_{\ol{A}}=Y_{\ol{A}}/\phi^\Z$ is the usual Fargues--Fontaine curve of $\ol{A}$.
\end{defi}

\begin{rem}
\begin{enumerate}[label=(\alph*)]
\item Instead of defining $\FF_A$ as a pushout in derived Berkovich spaces, we could also have defined the rings of functions on affine pieces of $\FF_A$ by a pullback of rings similarly to how we defined $Y_A$. Indeed, the same argument as in the proof of \cref{prop:defis-pushoutberk} shows that this alternative definition is equivalent to the one given above.

\item By a similar argument as in \cref{lem:prism-descentya}, the association $A\mapsto \FF_A$ from totally disconnected nilperfectoids to Gelfand stacks preserves $!$-equivalences and hence extends to a colimit-preserving functor
\begin{equation*}
\FF_{(-)}: \GelfStk\rightarrow\GelfStk\;, \hspace{0.3cm} X\mapsto \FF_X\;. \qedhere
\end{equation*}
\end{enumerate}
\end{rem}

As in the case of $Y_A$, we define degree $1$ Cartier divisors on $\FF_A$ to be exactly those whose pullback to $\FF_{\ol{A}}$ has degree $1$. 

\begin{defi}
The Gelfand stack $\Q_p^{\Div^1}=\Div^1$ is defined by
\begin{equation*}
\Div^1(A)=\{\text{degree $1$ Cartier divisors $D\subseteq\FF_A$}\}
\end{equation*}
on totally disconnected nilperfectoids $A$. For an arbitrary Gelfand stack $X$ over $\Q_p$, we define $X^{\Div^1}$ as the Gelfand stack over $\Div^1$ given by
\begin{equation*}
X^{\Div^1}(\GSpec A\xrightarrow{D} \Div^1)=\{\text{maps $D\rightarrow X$}\}
\end{equation*}
on totally disconnected nilperfectoids $A$.
\end{defi}

Recall that pulling back a degree $1$ Cartier divisor along the natural map $\iota: \GSpec A\rightarrow \FF_A$ induces a morphism 
\begin{equation*}
\mu: \Q_p^{\Div^1}\rightarrow \A^1/\G_m
\end{equation*}
and hence we obtain $\mu: X^{\Div^1}\rightarrow \Q_p^{\Div^1}\xrightarrow{\mu} \A^1/\G_m$. Moreover, recall that, in \cite{AnPrism}, the analytic Hodge--Tate stack $X^\HT$ of $X$ is defined as the divisor in $X^{\Div^1}$ classified by this morphism and we let $\cal{I}$ denote its ideal sheaf. Using the following proposition, one can check that this is equivalent to our previous definition of $X^\HT$. 

\begin{prop}
\label{prop:defis-xprismxdiv1}
For any Gelfand stack $X$, Frobenius induces an isomorphism $\phi: X^\prism_{(0, 1)}\xrightarrow{\cong} X^\prism_{(0, p)}$ and we have
\begin{equation*}
\operatorname{coeq}(\hspace{-0.15cm}
\begin{tikzcd}
X^\prism_{(0, 1)}\ar[r,shift left=.75ex, "\phi"]
  \ar[r,shift right=.75ex,swap, "\id"] & X^\prism_{(0, p)}
\end{tikzcd}
\hspace{-0.15cm})\cong X^{\Div^1}\nospacepunct{\;.}
\end{equation*}
\end{prop}
\begin{proof}
As in the proof of \cref{prop:defis-prismffdr}, the claim reduces to the following two assertions: For any totally disconnected nilperfectoid $A$ \dots
\begin{enumerate}[label=(\roman*)]
\item \dots Frobenius induces an isomorphism $\phi: Y_{A, (1/p, \infty)}\xrightarrow{\cong} Y_{A, (1, \infty)}$.
\item \dots there is an isomorphism
\begin{equation*}
\operatorname{coeq}(\hspace{-0.15cm}
\begin{tikzcd}
Y_{A, (1/p, \infty)}\ar[r,shift left=.75ex, "\phi"]
  \ar[r,shift right=.75ex,swap, "\id"] & Y_{A, (1/p, \infty)}
\end{tikzcd}
\hspace{-0.15cm})\cong \FF_A\nospacepunct{\;.}
\end{equation*}
\end{enumerate}
However, (i) follows from \cref{lem:defis-frobpushout} as the image of $\iota: \GSpec\ol{A}\rightarrow Y_{\ol{A}}$ is supported over $1\in (0, \infty)$ under the radius map. Moreover, (ii) follows from \cref{prop:defis-pushoutberk}.
\end{proof}

Putting together \cref{prop:defis-prismffdr} and \cref{prop:defis-xprismxdiv1} above, we see that understanding $X^\prism$ exactly amounts to understanding $X^\HK$ and $X^{\Div^1}$. Thus, let us quickly recall some explicit descriptions of $X^{\Div^1}$ and $X^\HK$ in special cases.

\begin{prop}
\label{prop:prism-div1presentations}
There is an isomorphism
\begin{equation*}
\Div^1\cong \lim_{q\mapsto q^p} (1+\overcirc{\DD})\setminus \{1\} \,\Big/\, \Q_p^{\times, \la}\;,
\end{equation*}
where the action of $\Q_p^{\times, \la}$ is by exponentiation, i.e.\ $\gamma.(q, q^{1/p}, \dots)=(q^\gamma, q^{\gamma/p}, \dots)$, and the map from $\lim_{q\mapsto q^p} (1+\overcirc{\DD})$ to $\Div^1$ sends an $A$-point $(q, q^{1/p}, \dots)$ of the source to the degree $1$ Cartier divisor on $\FF_A$ cut out by the function obtained by gluing $\log[q^\flat]$ on $\FF_{\ol{A}}$ with $\log q$ on $\GSpec A$ for any totally disconnected nilperfectoid $A$. Furthermore, as normed multiplicative group stacks over $\Div^1$, we have
\begin{equation*}
\G_m^{\Div^1}\cong \lim_{x\mapsto x^p} \G_m\,\Big/\,\Z_p^\la\;,
\end{equation*}
where $\theta\in\Z_p^\la$ acts on a sequence $(x, x^{1/p}, \dots)$ of $p$-power roots via 
\begin{equation*}
\theta.(x, x^{1/p}, \dots)=(q^\theta x, q^{\theta/p}x^{1/p}, \dots)\;.
\end{equation*}
\end{prop}
\begin{proof}
See \cite{AnPrism}.
\end{proof}

As $X^\HK=(X^{\Div^1})^\dR$ for any Gelfand stack $X$ by the argument from the proof of \cref{prop:defis-xprismdr}, the above immediately implies the following presentations of Hyodo--Kato stacks.

\begin{cor}
\label{cor:prism-hkpresentations}
There are isomorphisms
\begin{equation*}
\begin{split}
\Q_p^\HK&\cong \lim_{q\mapsto q^p} (1+\overcirc{\DD})^\dR\setminus \{1\} \,\Big/\, \Q_p^{\times, \sm}\;, \\
\G_m^{\HK}&\cong \lim_{x\mapsto x^p} \G_m^\dR\,\Big/\,\Z_p^\sm\;,
\end{split}
\end{equation*}
where the actions are as in the previous proposition.
\end{cor}

Finally, recall that both the functors $X\mapsto X^\HK$ and $X\mapsto X^{\Div^1}$ have good cohomological properties; in particular, they preserve smooth maps in the sense of the theorem below. To state the result, recall that $\cal{I}$ denotes the ideal sheaf of the Cartier divisor $X^\HT\rightarrow X^{\Div^1}$ and, moreover, recall from \cite{AnPrism}, that there is a canonical map $X^{\Div^1}\rightarrow */\phi^\Z$. Via pullback along this map, any isocrystal over $\Q_p$ yields a vector bundle on $X^{\Div^1}$, and we let $\O\{1\}$ denote the line bundle on $X^{\Div^1}$ obtained by tensoring $\cal{I}$ with the line bundle obtained from the isocrystal $(\Q_p, p^{-1}\phi)$. Also recall the line bundle $\O\{1\}$ on $\Q_p^\HK$, which corresponds to the isocrystal $(\Q_p, p^{-1}\phi)$ via the equivalence from \cref{thm:recall-qphk}.

\begin{thm}
\label{thm:prism-sixfunctorsdiv1hk}
Let $f: X\rightarrow Y$ be a rigid smooth morphism of derived Berkovich spaces which is pure of relative dimension $d$. Then \dots
\begin{enumerate}[label=(\roman*)]
\item \dots the induced map $X^\HK\rightarrow Y^\HK$ is cohomologically smooth with dualising sheaf $\O\{d\}[2d]$.
\item \dots the induced map $X^{\Div^1}\rightarrow Y^{\Div^1}$ is cohomologically smooth with dualising sheaf $\O\{d\}[2d]$.
\end{enumerate}
\end{thm}
\begin{proof}
See \cite[Thm.\ 6.3.1]{dRFF} and \cite{AnPrism}.
\end{proof}

From this, we at once deduce the following properties of the functor $X\mapsto X^\prism$.

\begin{cor}
\label{cor:prism-sixfunctors}
Let $f: X\rightarrow Y$ be a rigid smooth morphism of derived Berkovich spaces which is pure of relative dimension $d$. Then the induced map $X^\prism\rightarrow Y^\prism$ is cohomologically smooth with dualising sheaf $\O\{d\}[2d]$.
\end{cor}
\begin{proof}
The claim is local on the target and hence we may separately consider the maps $X^\prism_{(1, \infty)}\rightarrow Y^\prism_{(1, \infty)}$ and $X^\prism_{(0, p)}\rightarrow Y^\prism_{(0, p)}$. However, then \cref{prop:defis-prismffdr} and \cref{prop:defis-xprismxdiv1} supply pullback diagrams
\begin{equation*}
\begin{tikzcd}[column sep=large]
X^\prism_{(0, p)}\ar[r]\ar[d] & Y^\prism_{(0, p)}\ar[d] \\
X^{\Div^1}\ar[r] & Y^{\Div^1}\nospacepunct{\;,}
\end{tikzcd}
\hspace{1cm}
\begin{tikzcd}[column sep=large]
X^\prism_{(1, \infty)}\ar[r]\ar[d] & Y^\prism_{(1, \infty)}\ar[d] \\
X^\HK\ar[r] & Y^\HK\nospacepunct{\;,}
\end{tikzcd}
\end{equation*}
which at once imply the result by \cref{thm:prism-sixfunctorsdiv1hk}.
\end{proof}

We also have a smoothness result for $\Q_p^{\Div^1}$ and $\Q_p^\HK$ over $\GSpec\Q_p$.

\begin{thm}
\label{thm:prism-div1qphkcohsmooth}
The maps $\Q_p^{\Div^1}\rightarrow\GSpec\Q_p$ and $\Q_p^\HK\rightarrow\GSpec\Q_p$ are cohomologically smooth with dualising sheaf $\O\{1\}[2]$.
\end{thm}
\begin{proof}
For the first assertion, see \cite{AnPrism}. The second assertion is obtained by combining \cite[Prop.\ 5.7.1]{dRFF} and \cite[Prop.\ 6.2.4]{dRStack} due to \cref{cor:prism-hkpresentations}.
\end{proof}

\begin{cor}
\label{cor:prism-qpprismsmooth}
The map $\Q_p^\prism\rightarrow \GSpec\Q_p$ is cohomologically smooth with dualising sheaf $\O\{1\}[2]$.
\end{cor}
\begin{proof}
We may check this on an open cover of $\Q_p^\prism$, but for $\Q_{p, (0, p)}^\prism$ and $\Q_{p, (1, \infty)}^\prism$, the claim follows from the previous theorem as the maps $\Q_{p, (0, p)}^\prism\rightarrow \Q_p^{\Div^1}$ and $\Q_{p, (1, \infty)}^\prism\rightarrow \Q_p^\HK$ are cohomologically étale by \cref{prop:defis-prismffdr} and \cref{prop:defis-xprismxdiv1}.
\end{proof}

\newpage

\section{Analytic syntomification over $\Q_p$}

We now move on to defining the analytic syntomification $X^\Syn$ of a Gelfand stack $X$ and studying its geometry. To set the stage, recall how one proceeds for the ``algebraic'' syntomification in the sense of Bhatt--Lurie and Drinfeld, see \cite{FGauges} or \cite{Prismatization}. There, the main part of the work lies in constructing a filtered refinement of the prismatisation of $\Z_p$ called the (Nygaard-)filtered prismatisation or the \emph{Nygaardification} of $\Z_p$; from here, it is not too hard to define the Nygaardification of any $p$-adic formal scheme $X$. Finally, to obtain the syntomification, one then identifies two disjoint copies of the prismatisation inside the Nygaardification and glues them together.

In our analytic situation, the picture will be entirely similar. However, to get to the filtered prismatisation, we will take a slightly different approach than \cite{FGauges} and \cite{Prismatization} and instead use a convenient reformulation of the definition of the filtered prismatisation due to Gardner--Madapusi, see \cite[§6]{GardnerMadapusi}. Namely, in loc.\ cit., they show that one can alternatively obtain the Nygaardification of $\Z_p$ directly from the prismatisation of $\Z_p$ by taking a suitable pullback, and this will be the definition that we will transport to the analytic setting.

Even more so than the analytic prismatisation, the analytic Nygaardification $X^\N$ of a Gelfand stack $X$ will come with a plethora of extra structures, some of which we list here:
\begin{enumerate}[label=(\roman*)]
\item There is a map $\pi: X^\N\rightarrow X^\prism$, via which $X^\N$ acquires a radius map defined as the composite $\kappa: X^\N\rightarrow X^\prism\rightarrow (0, \infty)$.
\item There is a map $t: X^\N\rightarrow \ol{\DD}/\ol{\T}$ and a map $u: X^\N\rightarrow (\ol{\DD}/\ol{\T})^\dR$.
\item There are isomorphisms $X^\prism\cong X^\N_{|t|=1}$ and $X^\prism\cong X^\N_{|u|=1}$ and we denote the corresponding closed embeddings $X^\prism\rightarrow X^\N$ by $j_\dR$ and $j_\HT$, respectively. Moreover, the composition $\pi\circ j_\dR$ is the identity while $\pi\circ j_\HT$ is the Frobenius on $X^\prism$.
\item There is an isomorphism $X^\N_{|ut|\neq 0} \cong (X^\prism\setminus X^\dR)\times [0, 1]$, where we view $X^\dR$ as a closed substack of $X^\prism$ via $i_\dR$.
\item There is an isomorphism $X^{\dR, +}\cong X^\N_{|u|=0}$ and we denote the corresponding closed embedding $X^{\dR, +}\rightarrow X^\N$ by $i_{\dR, +}$.
\end{enumerate}
The situation is summarised by the schematic pictures in \cref{fig:xnviatu} and \cref{fig:xnviakappa}, the first of which studies $X^\N$ via the maps $(t, u): X^\N\rightarrow \ol{\DD}/\ol{\T}\times (\ol{\DD}/\ol{\T})^\dR$ while the second studies $X^\N$ via the radius map $\kappa: X^\N\rightarrow (0, \infty)$.

At the end of the section, we will prove two somewhat technical results that will be invaluable several times later on: First, we prove that the assignment sending a derived Berkovich space $X$ to its Nygaardification $X^\N$ preserves Berkovich étale maps and that $X^\N\rightarrow Y^\N$ is a $!$-cover if $X\rightarrow Y$ is a Berkovich étale $!$-cover of derived Berkovich spaces. This will be very useful later on when we study Nygaardifications of Berkovich smooth derived Berkovich spaces as it basically allows us to reduce to the case where $X=\ol{\T}^n$ is the overconvergent $n$-dimensional torus. Second, we show that the stack $\Q_{p, [r, s]}^\N$ is nicely coverable in the sense of \cref{defi:perf-nicelycoverable} whenever $p^{1/2}<r\leq s<p^{3/2}$; using the results of the next section, this will allow us to conclude that the restriction of the Nygaardification of $\ol{\T}^n$ to any $[r, s]\subseteq (p^{1/2}, p^{3/2})$ is nicely coverable for all $n\geq 0$. This will later enable us to apply the spreading out results for perfect complexes from §\ref{sect:recall} to $X^\N$ for any smooth derived Berkovich space $X$.

\subsection{Nygaard-filtered prismatisation}

\begin{figure}

\begin{center}
\begin{tikzpicture}
  \def\size{7} 
  \def\sqr_size{0.08*\size} 

  \def\circradius{0.02*\size} 
  \draw[line width=0.2mm, color=blue] (0, \sqr_size/2) -- (0, \size); 
  \draw[thick] (0, \size) -- (\size-\circradius, \size); 
  \draw[thick] (\size, 0) -- (\size, \size-\circradius); 
  \shade[shading=axis, bottom color=white, top color=white, middle color=green, shading angle=0](\sqr_size/2, -0.1) rectangle (\size+0.02, 0.1); 
  
  \draw[thick, dotted] (\size, \size) circle (\circradius);
  
  \draw[line width=0.2mm] (-\sqr_size/2,-\sqr_size/2) rectangle (\sqr_size/2, \sqr_size/2); 

  \draw[line width=0.2mm] (-\sqr_size/2, -\sqr_size/2) -- (\sqr_size/2, \sqr_size/2); 
  \draw[line width=0.2mm] (\sqr_size/2, -\sqr_size/2) -- (-\sqr_size/2, \sqr_size/2); 
  
  \draw[thick,scale=\size,domain=0.75:1,smooth,variable=\t] plot (\t, 0.75/\t);
  \draw[thick,scale=\size,domain=0.55:1,smooth,variable=\t] plot (\t, 0.55/\t);
  \draw[thick,scale=\size,domain=0.35:1,smooth,variable=\t] plot (\t, 0.35/\t);
  \draw[thick,scale=\size,domain=0.2:1,smooth,variable=\t] plot (\t, 0.2/\t);
  \draw[thick,scale=\size,domain=0.1:1,smooth,variable=\t] plot (\t, 0.1/\t);
  \draw[thick,scale=\size,domain=0.05:1,smooth,variable=\t] plot (\t, 0.05/\t);
  \draw[thick,scale=\size,domain=0.02:1,smooth,variable=\t,samples=100] plot (\t, 0.02/\t);
  
  \draw[->, thick] (-\size/3, -0.15*\size) -- (-\size/3, 1.15*\size);
  \draw[thick] (-\size/3-0.2, 0) -- (-\size/3+0.2, 0);
  \node at (-\size/3-0.5, 0) {$0$};
  \draw[thick] (-\size/3-0.2, \size) -- (-\size/3+0.2, \size);
  \node at (-\size/3-0.5, \size) {$1$};
  \node at (-\size/3, 1.15*\size+0.5) {$t$};
  
  \draw[->, thick] (-0.15*\size, -\size/4) -- (1.15*\size, -\size/4);
  \draw[thick] (0, -\size/4-0.2) -- (0, -\size/4+0.2);
  \node at (0, -\size/4-0.5) {$0$};
  \draw[thick] (\size, -\size/4-0.2) -- (\size, -\size/4+0.2);
  \node at (\size, -\size/4-0.5) {$1$};
  \node at (1.15*\size+0.5, -\size/4) {$u$};
  
  \node at (\size/2, -\size/12) {$X^{\HT, \dagger, +}$};
  \node at (-\size/10, \size/2) {$X^{\dR, +}$};
  \node at (\size*8/7, \size/2) {$j_\HT(X^\prism)$};
  \node at (\size/2, \size*13/12) {$j_\dR(X^\prism)$};
  \node at (\size*13/12, -\size/30) {$X^\HT$};
  \node at (-\size/30, \size*21/20) {$X^\dR$};

\end{tikzpicture}
\end{center}

\captionsetup{justification=centering}
\caption{A schematic picture of $X^\N$ via the map $(t, u)$}
\label{fig:xnviatu}
\end{figure}

As already announced above, our strategy to define $\Q_p^\N$ will be to imitate the pullback definition of the algebraic Nygaardification of $\Z_p$ given in \cite[§6]{GardnerMadapusi}. For this, first recall from \cref{lem:defis-mu} that the map $\mu: \Q_p^\prism\rightarrow\A^1/\G_m$ obtained by pulling back a degree $1$ divisor along $\iota: \GSpec A\rightarrow Y_A$ can be refined to a map $\mu: \Q_p^\prism\rightarrow \overcirc{\DD}/\ol{\T}$. By the same arguments, the map $\widetilde{\mu}: \Q_p^\prism\rightarrow (\A^1/\G_m)^\dR$ obtained by instead pulling back a degree $1$ divisor along $\phi^{-1}\circ\iota: \GSpec\ol{A}\rightarrow Y_A$ may be refined to a map
\begin{equation*}
\widetilde{\mu}: \Q_p^\prism\rightarrow (\overcirc{\DD}/\ol{\T})^\dR\;.
\end{equation*}
With this in place, we can make the following definition.

\begin{defi}
The (rational) \emph{analytic Nygaard-filtered prismatisation} of $\Q_p$ is the Gelfand stack $\Q_p^\N$ defined as the pullback
\begin{equation*}
\begin{tikzcd}
\Q_p^\N\ar[r, "\pi"]\ar[d, "{(t, u)}", swap] & \Q_p^\prism\ar[d, "\widetilde{\mu}"] \\
\ol{\DD}_+/\ol{\T}\times (\ol{\DD}_-/\ol{\T})^\dR\ar[r, "\mathrm{mult}"] & (\ol{\DD}_+/\ol{\T})^\dR\nospacepunct{\;.}
\end{tikzcd}
\end{equation*}
Here, the subscript $(-)_+$ or $(-)_-$ indicates whether $\ol{\T}$ acts on $\ol{\DD}$ by multiplication or division, respectively, and the bottom map is induced by the multiplication map on $\ol{\DD}$ and the division map $\ol{\T}\times \ol{\T}\rightarrow \ol{\T}, (r, s)\mapsto r/s$.
\end{defi}

In the future, we will mostly suppress the subscripts $(-)_+$ or $(-)_-$ from the notation and only use them for particular emphasis. Moreover, in the setting above, we will often also use $t$ and $u$ not only to denote the maps from $\Q_p^\N$ to $\ol{\DD}_+/\ol{\T}$ and $(\ol{\DD}_-/\ol{\T})^\dR$ coming from the diagram above, but also the coordinates on $\ol{\DD}_+$ and $\ol{\DD}_-$, respectively. 

\begin{rem}
\label{rem:syn-supportqpn}
Note that, since $\widetilde{\mu}$ factors through $(\overcirc{\DD}/\ol{\T})^\dR\subseteq (\ol{\DD}/\ol{\T})^\dR$, the stack $\Q_p^\N$ is already supported on the locus
\begin{equation*}
\{|ut|<1\}\subseteq \ol{\DD}/\ol{\T}\times (\ol{\DD}/\ol{\T})^\dR
\end{equation*}
via the map $(t, u)$.
\end{rem}

Let us unpack the definition of $\Q_p^\N$ and describe its functor of points on totally disconnected nilperfectoids $A$. Namely, an $A$-point consists of a triple $(D\subseteq Y_A, t: L\rightarrow A, u: \ol{A}\rightarrow K)$, where \dots
\begin{enumerate}[label=(\roman*)]
\item \dots $D\subseteq Y_A$ is a degree $1$ Cartier divisor,
\item \dots $t: L\rightarrow A$ is a normed generalised Cartier divisor on $A$ of norm at most $1$,
\item \dots $u: \ol{A}\rightarrow K$ is the dual of a normed generalised Cartier divisor on $\ol{A}$ of norm at most $1$,
\end{enumerate}
together with an isomorphism
\begin{equation}
\label{eq:syn-compatibility}
D\times_{Y_A} \GSpec\ol{A}\cong \GSpec\Cone(L\tensor_A K^\vee\xrightarrow{ut} \ol{A})
\end{equation}
of normed generalised Cartier divisors on $\ol{A}$. Here, the map $\GSpec\ol{A}\rightarrow Y_A$ occurring in the pullback on the left-hand side is of course the map $\phi^{-1}\circ\iota$ also used to define $\widetilde{\mu}$ and the fact that we have refined $\widetilde{\mu}$ to a map with target $(\ol{\DD}/\ol{\T})^\dR$ endows the left-hand side of (\ref{eq:syn-compatibility}) with the structure of a normed generalised Cartier divisor on $\ol{A}$.

Keeping the notation from the previous paragraph, let us introduce the abbreviation $\ol{L}\coloneqq L\tensor_A \ol{A}$. Then the map $u: K^\vee\rightarrow\ol{A}$ induces a map
\begin{equation*}
\GSpec\Cone(\ol{L}\rightarrow\ol{A})\rightarrow\GSpec\Cone(L\tensor_A K^\vee\rightarrow\ol{A})
\end{equation*}
and, using the isomorphism (\ref{eq:syn-compatibility}), we obtain a composite map
\begin{equation}
\label{eq:syn-mapdefd+}
\GSpec\Cone(\ol{L}\rightarrow\ol{A})\rightarrow D\times_{Y_A} \GSpec\ol{A}\xrightarrow{\mathrm{pr}_1} D\;.
\end{equation}
Using this map, we can define a new Gelfand stack $D^+$ as the pushout 
\begin{equation}
\label{eq:syn-defid+}
\begin{tikzcd}
\GSpec\Cone(\ol{L}\rightarrow\ol{A})\ar[r]\ar[d] & \GSpec\Cone(L\rightarrow A)\ar[d] \\
D\ar[r] & D^+
\end{tikzcd}
\end{equation}
in the category of derived Berkovich spaces. 

\begin{figure}

\begin{center}
\begin{tikzpicture}
  \def\width{5.5} 
  \def\height{\width} 

  \draw[dotted, thick] (0, 0) -- (\width, 0); 

  

  
  \draw[thick] (\width-0.2, 0.2*\height) -- (\width+0.2, 0.2*\height); 
  \node at (\width+1.2, 0.2*\height) {$\phi^{-2}(X^\HT)$};
  \draw[thick] (\width-0.2, 0.4*\height) -- (\width+0.2, 0.4*\height); 
  \node at (\width+1.2, 0.4*\height) {$\phi^{-1}(X^\HT)$};
  \draw[thick] (\width-0.2, 0.6*\height) -- (\width+0.2, 0.6*\height); 
  \node at (\width+0.7, 0.6*\height) {$X^\HT$};
  \draw[thick] (\width-0.2, 0.8*\height) -- (\width+0.2, 0.8*\height); 
  \node at (\width+0.7, 0.8*\height) {$X^\dR$};
  \draw[thick] (\width-0.2, \height) -- (\width+0.2, \height); 
  \node at (\width+0.95, \height) {$\phi(X^\dR)$};

  \draw[thick] (-0.2, 0.2*\height) -- (0.2, 0.2*\height); 
  \node at (-1.2, 0.2*\height) {$\phi^{-1}(X^\HT)$};
  \draw[thick] (-0.2, 0.4*\height) -- (0.2, 0.4*\height); 
  \node at (-0.7, 0.4*\height) {$X^\HT$};
  \draw[thick] (-0.2, 0.6*\height) -- (0.2, 0.6*\height); 
  \node at (-0.7, 0.6*\height) {$X^\dR$};
  \draw[thick] (-0.2, 0.8*\height) -- (0.2, 0.8*\height); 
  \node at (-0.95, 0.8*\height) {$\phi(X^\dR)$};
  \draw[thick] (-0.2, \height) -- (0.2, \height); 
  \node at (-1.05, \height) {$\phi^2(X^\dR)$};
  
  \node at (0, 1.3*\height) {$j_\dR(X^\prism)$};
  
  \node at (\width, 1.3*\height) {$j_\HT(X^\prism)$};


  \draw[line width=0.2mm, color=blue] (\width/2-0.05*\width, 0.6*\height) -- (0, 0.6*\height);

  \def\sqr_size{0.1*\width} 

  \draw[line width=0.2mm] (\width/2 - \sqr_size/2, 0.6*\height - \sqr_size/2) -- (\width/2 + \sqr_size/2, 0.6*\height + \sqr_size/2); 
  \draw[line width=0.2mm] (\width/2 + \sqr_size/2, 0.6*\height - \sqr_size/2) -- (\width/2 - \sqr_size/2, 0.6*\height + \sqr_size/2); 

  
  \shade[shading=axis, bottom color=white, top color=white, middle color=green, shading angle=0](0.5*\width + \sqr_size/2, 0.6*\height - 0.1) rectangle (\width, 0.6*\height + 0.1);
  \shade[shading=axis, bottom color=white, top color=white, middle color=green, shading angle=0](0, 0.4*\height - 0.1) rectangle (\width, 0.4*\height + 0.1);
  \shade[shading=axis, bottom color=white, top color=white, middle color=green, shading angle=0](0, 0.2*\height - 0.1) rectangle (\width, 0.2*\height + 0.1);

  \draw[->, thick] (2*\width, 0) -- (2*\width, 1.2*\height);
  
  \draw[thick] (2*\width-0.2, 0) -- (2*\width+0.2, 0); 
  \node at (2*\width+0.6, 0) {$0$};
  \draw[thick] (2*\width-0.2, 0.2*\height) -- (2*\width+0.2, 0.2*\height); 
  \node at (2*\width+0.6, 0.2*\height) {$1/p$};
  \draw[thick] (2*\width-0.2, 0.4*\height) -- (2*\width+0.2, 0.4*\height); 
  \node at (2*\width+0.6, 0.4*\height) {$1$};
  \draw[thick] (2*\width-0.2, 0.6*\height) -- (2*\width+0.2, 0.6*\height); 
  \node at (2*\width+0.6, 0.6*\height) {$p$};
  \draw[thick] (2*\width-0.2, 0.8*\height) -- (2*\width+0.2, 0.8*\height); 
  \node at (2*\width+0.6, 0.8*\height) {$p^2$};
  \draw[thick] (2*\width-0.2, \height) -- (2*\width+0.2, \height); 
  \node at (2*\width+0.6, \height) {$p^3$};
  
  \draw[->, thick] (1.3*\width, 0.6*\height) -- (1.8*\width, 0.6*\height);
  \node at (1.55*\width, 0.6*\height+0.3) {$\kappa$};
  
  \draw[->, thick] (-0.5*\width, 0.35*\height) -- (-0.5*\width, 0.85*\height);
  \node at (-0.5*\width-0.4, 0.6*\height) {$\phi$};
  
  \draw[line width=0.2mm, color=blue] (\width, 0.8*\height) -- (0, 0.8*\height);
  \draw[line width=0.2mm, color=blue] (\width, \height) -- (0, \height);
  
  \draw[thick] (0, 0) -- (0, 1.2*\height); 
  \draw[thick] (\width, 0) -- (\width, 1.2*\height); 
  
  \draw[line width=0.2mm] (\width/2 - \sqr_size/2, 0.6*\height - \sqr_size/2) rectangle (\width/2 + \sqr_size/2, 0.6*\height + \sqr_size/2); 

\end{tikzpicture}
\end{center}

\captionsetup{justification=centering}
\caption{A schematic picture of $X^\N$ via the radius map}
\label{fig:xnviakappa}
\end{figure}

\begin{rem}
\label{rem:syn-d+pullbackrings}
Recall that $D$ is affine and observe that $\Cone(L\rightarrow A)$ is a $\dagger$-nilpotent thickening of $\Cone(\ol{L}\rightarrow\ol{A})$. Thus, by the same argument as in the proof of \cref{prop:defis-pushoutberk}, we could have equivalently defined $\O(D^+)$ as the pullback 
\begin{equation*}
\begin{tikzcd}
\O(D^+)\ar[r]\ar[d] & \O(D)\ar[d] \\
\Cone(L\rightarrow A)\ar[r] & \Cone(\ol{L}\rightarrow\ol{A})
\end{tikzcd}
\end{equation*}
and then set $D^+\coloneqq \GSpec\O(D^+)$. Using this, a similar argument as in the proof of \cref{lem:prism-descentya} shows that the assignment given by sending a totally disconnected nilperfectoid $A$ over $\Q_p^\N$ to $D^+$ satisfies descent.
\end{rem}

\begin{rem}
We note that there is a natural map $D^+\rightarrow Y_A$ induced by the maps
\begin{equation*}
\GSpec\Cone(L\rightarrow A)\rightarrow\GSpec A\xrightarrow{\iota} Y_A
\end{equation*}
and $D\subseteq Y_A\xrightarrow{\phi} Y_A$, which one easily checks to be compatible.
\end{rem}

The above construction allows us to define $X^\N$ for an arbitrary Gelfand stack $X$.

\begin{defi}
For any Gelfand stack $X$, we define its (rational) \emph{analytic Nygaard-filtered prismatisation} as the Gelfand stack $X^\N$ over $\Q_p^\N$ given by
\begin{equation*}
X^\N(\GSpec A\rightarrow \Q_p^\N)\coloneqq \{\text{maps}\,D^+\rightarrow X\}
\end{equation*}
for any totally disconnected nilperfectoid $A$.
\end{defi}

Note that, as in the case of $X^\prism$, it is immediate from the definition that the association $X\mapsto X^\N$ commutes with limits.

\comment{
\begin{rem}
We do not know to which extent the sheafification is really necessary. However, as 
\begin{equation*}
D^+\cong \GSpec \Cone(L\tensor_A \Nil^\dagger(A)\rightarrow A)
\end{equation*}
over the locus $\{|u|=0\}\subseteq \Q_p^\N$, as we will see in \cref{prop:defis-u0} below, any result in this direction for $X^\N$ would automatically entail an analogous result for the filtered de Rham stack $X^{\dR, +}$.
\end{rem}
}

Note that any map $D^+\rightarrow X$ induces a map $D\rightarrow X$ by precomposition with $D\rightarrow D^+$ and hence we obtain a map
\begin{equation*}
\pi: X^\N\rightarrow X^\prism
\end{equation*}
for any Gelfand stack $X$ living over the map $\pi: \Q_p^\N\rightarrow\Q_p^\prism$. Moreover, using this we can equip $X^\N$ with a radius map given by the composition
\begin{equation*}
\kappa: X^\N\xrightarrow{\pi} X^\prism\xrightarrow{\kappa} (0, \infty)\;,
\end{equation*}
see also \cref{fig:xnviakappa}.

Before we move on to discussing the geometry of $X^\N$ in more detail, let us give a convenient repackaging of the above definition in terms of yet another stack associated to $X$.

\begin{defi}
For any Gelfand stack $X$, we define its \emph{cone stack} to be the Gelfand stack $X^{\Cone}$ over $\ol{\DD}/\ol{\T}$ obtained by sheafifying the assignment
\begin{equation*}
X^{\Cone}(\GSpec A\rightarrow \ol{\DD}/\ol{\T})\coloneqq \{\text{maps}\,\GSpec\Cone(L\rightarrow A)\rightarrow X\}\;.
\end{equation*}
\end{defi}

Let us immediately emphasise that the assignment $X\mapsto X^{\Cone}$ has very bad properties: e.g., we do not even know whether it sends finite étale maps to surjections of Gelfand stacks. Nevertheless, using $X^{\Cone}$, we can give the following alternative definition of $X^\N$ for any derived Berkovich space $X$.

\begin{lem}
\label{lem:syn-xnviacone}
Let $X$ be a derived Berkovich space. Then there is a pullback diagram
\begin{equation*}
\begin{tikzcd}
X^\N\ar[r, "\pi"]\ar[d] & X^\prism\ar[d] \\
X^{\Cone}\ar[r] & (X^{\Cone})^\dR\;,
\end{tikzcd}
\end{equation*}
where the map $X^\prism\rightarrow (X^{\Cone})^\dR$ is induced by the map (\ref{eq:syn-mapdefd+}).
\end{lem}
\begin{proof}
Immediate from the definition of $D^+$ as a pushout in derived Berkovich spaces.
\end{proof}

Let us point out that, conceptually, the reason why the formal properties of the functor $X\mapsto X^\N$ are much better than those of $X\mapsto X^{\Cone}$, as we will see below, e.g.\ in \cref{thm:defis-berketalecover}, is that $X^{\Cone}$ only enters into the definition of $X^\N$ via its fibres over $(X^{\Cone})^\dR$, which are much better behaved than $X^{\Cone}$ itself.

\subsection{Interesting loci inside $X^\N$}
\label{subsect:loci}

We now turn to studying the geometry of the stack $X^\N$ by describing some interesting loci inside it.

\subsubsection{The map $j_\dR$}

Note that the composition
\begin{equation*}
\ol{\T}/\ol{\T}\times (\ol{\DD}/\ol{\T})^\dR\rightarrow \ol{\DD}/\ol{\T}\times (\ol{\DD}/\ol{\T})^\dR\xrightarrow{\mathrm{mult}} (\ol{\DD}/\ol{\T})^\dR
\end{equation*}
is an isomorphism and hence the pullback of $\Q_p^\N$ along the first map is isomorphic to $\Q_p^\prism$ via $\pi: \Q_p^\N\rightarrow\Q_p^\prism$. In other words, $\pi$ induces an isomorphism $\Q_{p, |t|=1}^\N\cong \Q_p^\prism$.

Furthermore, over the locus $\{|t|=1\}\subseteq \Q_p^\N$, the normed generalised Cartier divisor $t: L\rightarrow A$ identifies with the identity on $A$ and hence
\begin{equation*}
\GSpec\Cone(\ol{L}\rightarrow\ol{A})=\emptyset=\GSpec\Cone(L\rightarrow A)\;.
\end{equation*}
Consequently, over $\Q_{p, |t|=1}^\N$, we have $D^+\cong D$ and this implies that $X^\N_{|t|=1}\cong X^\prism$ via $\pi$ for any Gelfand stack $X$. The corresponding closed embedding of $X^\prism$ into $X^\N$ is denoted by
\begin{equation*}
j_\dR: X^\prism\rightarrow X^\N\;.
\end{equation*}

\subsubsection{The map $j_\HT$}

Now let us study the locus $\{|u|=1\}$ inside $X^\N$. Unwinding definitions, we see that an $A$-point of $\Q_{p, |u|=1}^\N$ for $A$ totally disconnected nilperfectoid corresponds to a degree $1$ divisor $D\subseteq Y_A$ together with a lift $t: L\rightarrow A$ of the normed generalised Cartier divisor $D\times_{Y_A} \GSpec\ol{A}$ on $\ol{A}$ to $A$. First noting that the diagram
\begin{equation*}
\begin{tikzcd}
\ol{\DD}/\ol{\T}\ar[r]\ar[d] & (\ol{\DD}/\ol{\T})^\dR\ar[d] \\
\A^1/\G_m\ar[r] & (\A^1/\G_m)^\dR
\end{tikzcd}
\end{equation*}
is cartesian (both the top and the bottom row are $\G_a^\dagger/\G_m^\dagger$-torsors), we see that it is equivalent to just lift $D\times_{Y_A} \GSpec \ol{A}$ from $\ol{A}$ to $A$ as a divisor, i.e.\ forgetting the norm. 

Moreover, by \cref{lem:prism-odtrivial}, we see that the divisor
\begin{equation*}
D\times_{Y_A} \GSpec\ol{A}\subseteq \GSpec\ol{A}
\end{equation*}
has trivial ideal sheaf, hence $L\cong A$ and we just have to lift the global section cutting out $D\times_{Y_A} \GSpec\ol{A}$ from $\ol{A}$ to $A$. As $D$ is also just cut out by an element of some $\mathbb{B}^\dagger_{[p^{-r}, p^r]}(A)$ by loc.\ cit., the pullback diagram
\begin{equation}
\label{eq:syn-jhtpullback}
\begin{tikzcd}
\mathbb{B}^\dagger_{[p^{-r+1}, p^{r+1}]}(A)\ar[r, "F"]\ar[d, "\theta", swap] & \mathbb{B}^\dagger_{[p^{-r}, p^{r}]}(A)\ar[d] \\
A\ar[r] & \ol{A}
\end{tikzcd}
\end{equation}
from the proof of \cref{lem:defis-frobpushout} shows that we may repackage the data $(D\subseteq Y_A, t: A\rightarrow A)$ as a divisor $D'\subseteq Y_A$ whose Frobenius pullback identifies with $D$. Moreover, note that $D'$ has degree $1$ if and only if $D$ has degree $1$ as this may be checked after pullback to $Y_{\ol{A}}$, where Frobenius is an isomorphism.

What about $D^+$? By \cref{rem:syn-d+pullbackrings}, the Gelfand stack $D^+$ is affine and there is a pullback diagram
\begin{equation*}
\begin{tikzcd}
\O(D^+)\ar[r]\ar[d] & \O(D)\ar[d] \\
\Cone(L\rightarrow A)\ar[r] & \Cone(\ol{L}\rightarrow\ol{A})\nospacepunct{\;,}
\end{tikzcd}
\end{equation*}
where the vertical map on the right just comes from pulling back a function on $D$ to a function on
\begin{equation*}
D\times_{Y_A} \GSpec\ol{A}\cong \GSpec\Cone(\ol{L}\rightarrow\ol{A})\;.
\end{equation*}
Letting $\xi\in\mathbb{B}^\dagger_{[p^{-r+1}, p^{r+1}]}(A)$ denote the function cutting out $D'$ as above, we can rewrite the pullback diagram above as
\begin{equation*}
\begin{tikzcd}
\O(D^+)\ar[r]\ar[d] & \mathbb{B}^\dagger_{[p^{-r}, p^{r}]}(A)/F(\xi)\ar[d] \\
A/\theta(\xi)\ar[r] & \ol{A}/\theta(\xi)
\end{tikzcd}
\end{equation*}
and hence (\ref{eq:syn-jhtpullback}) shows that $\O(D^+)\cong \mathbb{B}^\dagger_{[p^{-r+1}, p^{r+1}]}(A)/\xi$, i.e.\ $D^+\cong D'$. Overall, this discussion yields the following description of $X_{|u|=1}^\N$:

\begin{prop}
\label{lem:defis-jht}
For any Gelfand stack $X$, there is an isomorphism
\begin{equation*}
X^\prism\cong X^\N_{|u|=1}
\end{equation*}
whose postcomposition with \dots
\begin{enumerate}[label=(\roman*)]
\item \dots $\pi: X^\N\rightarrow X^\prism$ identifies with the Frobenius on $X^\prism$,
\item \dots $t: X^\N\rightarrow\Q_p^\N\xrightarrow{t} \ol{\DD}/\ol{\T}$ identifies with $\mu: X^\prism\rightarrow\Q_p^\prism\xrightarrow{\mu} \ol{\DD}/\ol{\T}$.
\end{enumerate}
\end{prop}

We denote the corresponding closed embedding of $X^\prism$ into $X^\N$ by
\begin{equation*}
j_\HT: X^\prism\rightarrow X^\N\;.
\end{equation*}

\subsubsection{The locus $\{|ut|\neq 0\}$}

We now want to describe the locus $\{|ut|\neq 0\}$ inside $X^\N$. For this, recall that $\ol{\DD}^\times/\ol{\T}\cong (0, 1]$ by \cref{lem:defis-rhvariant} and that $(0, 1]^\dR\cong (0, 1]$ by \cref{ex:recall-drbetti} and idempotency of the functor $(-)^\dR$. Furthermore recalling that $\Q_p^\N$ is actually supported over $\{|ut|<1\}$ by \cref{rem:syn-supportqpn}, we see that, over the locus $\{|ut|\neq 0\}$, the pullback diagram defining $\Q_p^\N$ simplifies drastically to yield
\begin{equation*}
\begin{tikzcd}
\Q_{p, |ut|\neq 0}^\N\ar[r, "\pi"]\ar[d, "{(t, u)}", swap] & \Q_{p, |\widetilde{\mu}|\neq 0}^\prism\ar[d, "\widetilde{\mu}"] \\
((0, 1]\times (0, 1])\setminus \{(1, 1)\}\ar[r, "\mathrm{mult}"] & (0, 1)\nospacepunct{\;.}
\end{tikzcd}
\end{equation*}
To move on, note that there is an isomorphism
\begin{equation*}
(0, 1)\times [0, 1]\xrightarrow{\cong} ((0, 1]\times (0, 1])\setminus \{(1, 1)\}
\end{equation*}
given by sending a pair $(x, y)$ to the unique intersection point $(t, u)$ of the hyperbola $\{ut=x\}$ and the line through $(1, 1)$ and $(1-y, y)$. Twisting the above pullback diagram by this isomorphism, we obtain the following result:

\begin{prop}
\label{prop:defis-utneq0}
For any Gelfand stack $X$, there is an isomorphism
\begin{equation*}
X^\N_{|ut|\neq 0}\cong (X^\prism\setminus X^\dR)\times [0, 1]\;,
\end{equation*}
where we view $X^\dR$ as a closed substack of $X^\prism$ via $i_\dR$. Moreover, this isomorphism identifies \dots
\begin{enumerate}[label=(\roman*)]
\item \dots $\pi$ with the first projection,
\item \dots the restriction of $j_\dR$ to $X^\prism\setminus X^\dR$ with $\id\times \{0\}$,
\item \dots the restriction of $j_\HT$ to $X^\prism\setminus X^\HT$ with $\phi\times \{1\}$.
\end{enumerate}
\end{prop}
\begin{proof}
The above discussion already shows that 
\begin{equation*}
\Q_{p, |ut|\neq 0}^\N\cong \Q_{p, |\widetilde{\mu}|\neq 0}^\prism\times [0, 1]\;.
\end{equation*}
Moreover, since the generalised Cartier divisor $t: L\rightarrow A$ is an isomorphism over $|t|\neq 0$, we see that $D^+\cong D$ over $\{|ut|\neq 0\}$, which implies the above isomorphism for any Gelfand stack $X$ in place of $\Q_p$. The only remaining claim is that $X^\prism_{|\widetilde{\mu}|\neq 0}\cong X^\prism\setminus X^\dR$ or, in other words, that $X^\dR\cong X^\prism_{|\widetilde{\mu}|=0}$ via $i_\dR$. By construction, this reduces to the case $X=\GSpec\Q_p$, where we observe that $\Q_{p, |\widetilde{\mu}|=0}^\prism$ is the moduli space of degree $1$ Cartier divisors $D\subseteq Y_A$ equipped with a factorisation of $\phi^{-1}\circ\iota: \GSpec\ol{A}\rightarrow Y_A$ through $D$. However, as $\phi^{-1}\circ\iota$ already defines a degree $1$ Cartier divisor itself, such a factorisation is actually an isomorphism $D\cong \GSpec\ol{A}$ by rigidity, see \cref{lem:prism-rigidity}, and this implies the claim.
\end{proof}

\subsubsection{The map $i_{\dR, +}$}

Let us now study the locus $\{|u|=0\}$ inside $X^\N$. For this, first note that $\Q_{p, |u|=0}^\N$ sits inside a pullback square
\begin{equation}
\label{eq:syn-idr+pullback}
\begin{tikzcd}
\Q_{p, |u|=0}^\N\ar[r]\ar[d] & \GSpec\Q_p\ar[d] \\
\ol{\DD}/\ol{\T}\times */\ol{\T}^\dR\ar[r, "\mathrm{mult}"] & */\ol{\T}^\dR\nospacepunct{\;,}
\end{tikzcd}
\end{equation}
where we have used that $\Q_{p, |\widetilde{\mu}|=0}^\prism\cong \Q_p^\dR\cong\GSpec\Q_p$ by the proof of \cref{prop:defis-utneq0} and the vertical map $\GSpec\Q_p\rightarrow */\ol{\T}^\dR$ on the right is given by the composition
\begin{equation*}
\GSpec\Q_p\cong \Q_p^\dR\xrightarrow{i_\dR} \Q_p^\prism\xrightarrow{\widetilde{\mu}} (\ol{\DD}/\ol{\T})^\dR\;,
\end{equation*}
which factors through $*/\ol{\T}^\dR$. We warn the reader that this map is \emph{not} the canonical map $*\rightarrow */\ol{\T}^\dR$! To understand this map more explicitly, let us first establish the following moduli description of $*/\ol{\T}^\dR$:

\begin{lem}
\label{lem:syn-btdr}
Let $X$ be a Gelfand stack. Then there is an equivalence between
\begin{equation*}
\{\text{maps } X\rightarrow */\ol{\T}^\dR\}\cong \{\text{normed line bundles on $X^\diamond$ for the arc-topology}\}\;.
\end{equation*}
In particular, a map $\GSpec\Q_p\rightarrow */\ol{\T}^\dR$ corresponds to a $\C_p$-linear continuous character of $\Gal_{\Q_p}$ of norm $1$, i.e.\ a continuous group morphism $\Gal_{\Q_p}\rightarrow\O_{\C_p}^\times$.
\end{lem}
\begin{proof}
Since $!$-covers induce arc-covers on diamonds, we are immediately reduced to the case $X=\GSpec B$ for $B$ nilperfectoid by descent. Then a $B$-valued point of $*/\ol{\T}^\dR$ amounts to a $\ol{B}$-valued point of $*/\ol{\T}^\diamond$, which is equivalent to a normed line bundle on $\GSpec\ol{B}$ for the arc-topology almost by definition. The final assertion follows from the fact that $\Spd\Q_p\cong \Spd\C_p/\ul{\Gal_{\Q_p}}$ as all line bundles on $\Spd\C_p$ are trivial.
\end{proof}

Under the identification from the previous lemma, the map $\GSpec\Q_p\rightarrow */\ol{\T}^\dR$ classifies the cyclotomic character $\Gal_{\Q_p}\rightarrow \Z_p^\times$, which in terms of line bundles is the Tate twist. Indeed, for $\ol{A}=\Q_p^\cycl$, pulling back the ideal sheaf of the divisor $\phi^{-1}\circ\iota: \GSpec\ol{A}\rightarrow Y_{\ol{A}}$ classified by $i_\dR$ along $\phi^{-1}\circ\iota$, we obtain the invertible $\ol{A}$-module 
\begin{equation*}
p\log[\epsilon]W(\ol{A}^\flat)[\tfrac{1}{p}]/(p\log[\epsilon])^2 W(\ol{A}^\flat)[\tfrac{1}{p}]\;,
\end{equation*}
where $\epsilon=(1, \zeta_p, \zeta_{p^2}, \dots)$, as usual, and $\Gal_{\Q_p}$ acts on this via the cyclotomic character.

Returning to the pullback diagram (\ref{eq:syn-idr+pullback}), this tells us that $\Q_{p, |u|=0}^\N$ classifies triples $(D\subseteq Y_A, t: L\rightarrow A, u: \ol{A}\rightarrow K)$ for a totally disconnected nilperfectoid $A$ together with an identification between $u$ and the zero map as well as an isomorphism $L\tensor_A K^\vee\cong \ol{A}(1)$, which we may rewrite as $K\cong \ol{L}(-1)$. In other words, we see that the map
\begin{equation*}
i_{\dR, +}: \ol{\DD}/\ol{\T}\rightarrow \Q_p^\N
\end{equation*}
given by sending a normed generalised Cartier divisor $t: L\rightarrow A$ on a totally disconnected nilperfectoid $A$ to the triple $(\GSpec\ol{A}\subseteq Y_A, t: L\rightarrow A, \ol{A}\xrightarrow{0} \ol{L}(-1))$ is an isomorphism onto $\Q_{p, |u|=0}^\N$.

What happens to $D^+$ after pullback along $i_{\dR, +}$? By \cref{rem:syn-d+pullbackrings}, the Gelfand stack $D^+$ is affine and there is a pullback diagram
\begin{equation*}
\begin{tikzcd}
\O(D^+)\ar[r]\ar[d] & \ol{A}\ar[d] \\
\Cone(L\rightarrow A)\ar[r] & \Cone(\ol{L}\rightarrow\ol{A})\nospacepunct{\;,}
\end{tikzcd}
\end{equation*}
where we have used that $D\cong\GSpec\ol{A}$ after pullback along $i_\dR$ and the vertical map on the right is the canonical map. Thus, we conclude that
\begin{equation*}
D^+\cong \GSpec\Cone(L\tensor_A \Nil^\dagger(A)\rightarrow A)\;,
\end{equation*}
whence the map $i_{\dR, +}: \ol{\DD}/\ol{\T}\rightarrow\Q_p^\N$ induces a map 
\begin{equation*}
i_{\dR, +}: X^{\dR, +}\rightarrow X^\N\;.
\end{equation*}
The above discussion then shows the following:

\begin{prop}
\label{prop:defis-u0}
For any Gelfand stack $X$, the map $i_{\dR, +}$ induces an isomorphism
\begin{equation*}
X^{\dR, +}\cong X^\N_{|u|=0}\;.
\end{equation*}
\end{prop}

\subsubsection{The locus $\{|t|=0\}$}

Finally, let us shortly discuss the locus $\{|t|=0\}$ inside $X^\N$. In some sense, this is the most mysterious and not directly comparable to any of the stacks we already know. Still, let us at least give it a name.

\begin{defi}
For any Gelfand stack $X$, we write 
\begin{equation*}
X^{\HT, \dagger, +}\coloneqq X^\N_{|t|=0}
\end{equation*}
and $X^{\HT, +}$ for the base change of $X^{\HT, \dagger, +}$ along $\{t=0\}\rightarrow \{|t|=0\}$.
\end{defi} 

The above notation is motivated by the fact that $X^{\HT, +}$ is a filtered refinement of the Hodge--Tate stack while $X^{\HT, \dagger, +}$ is a filtered refinement of the overconvergent neighbourhood
\begin{equation*}
X^{\HT, \dagger}\coloneqq (X^\HT\subseteq X^\prism)^\dagger=X^\prism_{|\mu|=0}
\end{equation*}
of $X^\HT$ inside $X^\prism$. Indeed, we can at least prove the following.

\begin{prop}
\label{prop:defis-t0uneq0}
For any Gelfand stack $X$, there is an isomorphism
\begin{equation*}
X^{\HT, \dagger, +}_{|u|\neq 0}\cong X^{\HT, \dagger}\times (0, 1]
\end{equation*}
intertwining \dots
\begin{enumerate}[label=(\roman*)]
\item \dots the map $\id\times \{1\}: X^{\HT, \dagger}\rightarrow X^{\HT, \dagger}\times (0, 1]$ with $j_\HT$,
\item \dots the map $\phi\circ\mathrm{pr}_1: X^{\HT, \dagger}\times (0, 1]\rightarrow X^\prism$ with $\pi$.
\end{enumerate}
In particular, we have $X^{\HT, +}_{|u|\neq 0} \cong X^\HT\times (0, 1]$.
\end{prop}
\begin{proof}
We first consider the case $X=\GSpec\Q_p$, where the definition of $\Q_p^\N$ supplies a pullback square
\begin{equation*}
\begin{tikzcd}
\Q_{p, |u|\neq 0}^{\HT, \dagger, +}\ar[r, "\pi"]\ar[d] & \GSpec\Q_p\ar[d, "\widetilde{\mu}"]\\
\G_a^\dagger/\ol{\T}\times (0, 1]\ar[r, "\mathrm{mult}"] & */\ol{\T}^\dR\nospacepunct{\;,}
\end{tikzcd}
\end{equation*}
where we have used that $\ol{\DD}^\times/\ol{\T}\cong (0, 1]$ by \cref{lem:defis-rhvariant} and that $\Q_{p, |\widetilde{\mu}|=0}^\prism\cong \GSpec\Q_p$ via $i_\dR$ as in the proof of \cref{prop:defis-utneq0}. As the horizontal map on the bottom agrees with
\begin{equation*}
\G_a^\dagger/\ol{\T}\times (0, 1]\xrightarrow{\mathrm{pr}_1} \G_a^\dagger/\ol{\T}\rightarrow */\ol{\T}^\dR\;,
\end{equation*}
the claim now follows from the fact that the diagram
\begin{equation*}
\begin{tikzcd}
\Q_p^{\HT, \dagger}\ar[r, "\pi"]\ar[d] & \GSpec\Q_p\ar[d, "\widetilde{\mu}"]\\
\G_a^\dagger/\ol{\T}\ar[r] & */\ol{\T}^\dR
\end{tikzcd}
\end{equation*}
is cartesian by \cref{lem:defis-jht}. The case of general $X$ now follows by verifying that the pullback of $D^+$ to $\Q_p^{\HT, \dagger}\cong \Q_{p, |u|=1}^{\HT, \dagger, +}$ identifies with the universal degree $1$ Cartier divisor $D$ over $\Q_p^{\HT, \dagger}$, but this is similarly deduced from the proof of \cref{lem:defis-jht}.
\end{proof}

\subsection{Syntomification}
\label{subsect:syn}

Finally, we can define the syntomification. As already announced above, for this, we will glue the two copies of $X^\prism$ embedded into $X^\N$ via $j_\dR$ and $j_\HT$, respectively.

\begin{defi}
Let $X$ be a Gelfand stack. The (rational) \emph{analytic syntomification} of $X$ is the Gelfand stack $X^\Syn$ defined by the coequaliser diagram
\begin{equation*}
\begin{tikzcd}
X^\prism\ar[r,shift left=.75ex,"j_\HT"]\ar[r,shift right=.75ex,swap,"j_\dR"] & X^\N\ar[r] & X^\Syn\;.
\end{tikzcd}
\end{equation*}
An object of the category $\D(X^\Syn)$ is called an \emph{analytic (prismatic) $F$-gauge} on $X$.
\end{defi}

The most important example of an analytic $F$-gauge is the Breuil--Kisin twist. Recall from \cref{defi:prism-bktwist} that we have already defined the Breuil--Kisin twist on $X^\prism$.

\begin{defi}
Let $X$ be a Gelfand stack. The \emph{Breuil--Kisin twist} is the line bundle $\O\{1\}$ on $X^\N$ defined by
\begin{equation*}
\O\{1\}\cong \pi^*\O\{1\}\tensor t^*\O\langle -1\rangle\;,
\end{equation*}
where we recall that $t: \O\langle -1\rangle\rightarrow \O$ is the universal normed generalised Cartier divisor on $\ol{\DD}_+/\ol{\T}$. It descends to $X^\Syn$ and we also use $\O\{1\}$ to denote this descended line bundle.
\end{defi}

Indeed, to see that the Breuil--Kisin twist descends to $X^\Syn$, one calculates
\begin{equation*}
j_\dR^*\O\{1\}\cong \O\{1\}\cong \phi^*\O\{1\}\tensor\cal{I}\cong j_\HT^*\O\{1\}
\end{equation*}
using \cref{lem:defis-jht} as well as the fact that $t\circ j_\HT=\mu$ classifies the Hodge--Tate divisor while $t\circ j_\dR$ factors through $\ol{\T}/\ol{\T}\subseteq \ol{\DD}/\ol{\T}$. With this, we can finally define analytic syntomic cohomology of Gelfand stacks in our sense.

\begin{defi}
Let $X$ be a Gelfand stack. For any $i\in\Z$, the \emph{analytic syntomic cohomology} of weight $i$ of $X$ is defined as
\begin{equation*}
R\Gamma_\Syn(X, \Q_p(i))\coloneqq R\Gamma(X^\Syn, \O\{i\})\;.
\end{equation*}
\end{defi}

As already hinted at by the terminology in the above definition, there is a general notion of \emph{Hodge--Tate weights} for perfect complexes on $X^\Syn$, which we shall now define. For this, note that precomposing the map $i_{\dR, +}$ with the canonical map $X^\Hod\rightarrow X^{\dR, +}$ yields a map
\begin{equation*}
i_\Hod: X^\Hod\rightarrow X^\N\;.
\end{equation*}
Moreover, recall from \cref{lem:recall-reesgm} that perfect complexes on $X\times */\ol{\T}$ identify with graded perfect complexes on $X$.

\begin{defi}
Let $X$ be a Gelfand stack and $E$ a perfect complex on $X^\N$. Identifying the pullback of $E$ along
\begin{equation*}
X\times */\ol{\T}\rightarrow X^\Hod\xrightarrow{i_\Hod} X^\N
\end{equation*}
with a graded perfect complex $M^\bullet$ on $X$, the \emph{Hodge--Tate weights} of $E$ are defined to be those integers $i\in\Z$ such that $M^i$ is nonzero. The Hodge--Tate weights of a perfect analytic $F$-gauge $E\in\Perf(X^\Syn)$ on $X$ are defined to be the Hodge--Tate weights of its pullback to $X^\N$.
\end{defi}

\begin{ex}
The line bundle $\O\{i\}$ only has a single Hodge--Tate weight, which is equal to $-i$. Indeed, one easily calculates $i_\Hod^*\O\{1\}\cong \O\langle -1\rangle$ and this yields the claim.
\end{ex}

\subsection{Compatibility with étale localisation}

We now want to prove compatibility of the functor $X\mapsto X^\N$ with Berkovich étale localisation. For this, we first extend the definition of Berkovich étale maps between derived Berkovich spaces from \cite[Def.\ 4.3.4.(3)]{dRFF} to maps between arbitrary Gelfand stacks as follows.

\begin{defi}
Let $X\rightarrow Y$ be a map of Gelfand stacks. It is called \emph{Berkovich étale} if $X\times_Y \GSpec A$ is a derived Berkovich space for any totally disconnected nilperfectoid $A$ with a map $\GSpec A\rightarrow Y$ and the map
\begin{equation*}
X\times_Y \GSpec A\rightarrow\GSpec A
\end{equation*}
is a Berkovich étale map of derived Berkovich spaces in the sense of \cite[Def.\ 4.3.4.(3)]{dRFF}. 
\end{defi}

Similarly, we can define finite étale maps between Gelfand stacks as well as open embeddings and rational localisations, and we note that all these notions will be stable under base change as an immediate consequence of the definition. With this notation in place, the result we want to prove reads as follows.

\begin{thm}
\label{thm:defis-berketalecover}
Let $f: X\rightarrow Y$ be a Berkovich étale map of derived Berkovich spaces. Then the induced map $f^\N: X^\N\rightarrow Y^\N$ is Berkovich étale as well. Moreover, if $f$ is a $!$-cover, then the same is true for the map $f^\N$.
\end{thm}

Note that, according to the discussion of various loci inside $X^\N$ above, the theorem above implies the corresponding statement for all other stacks associated to $p$-adic cohomology theories we have introduced so far, e.g.\ $X^{\dR, +}$ or $X^\prism$, by suitable base changes. 

\comment{
\begin{rem}
Note that the corresponding statement for the syntomification does \emph{not} seem to be a completely formal consequence of the one for the Nygaardification. Nevertheless, one can still deduce it from the arguments below with a bit more care by using the geometric properties of $X^\N$ discussed above. However, we are not going to use this result in the sequel.
\end{rem}
}

To prove \cref{thm:defis-berketalecover}, recall that a Berkovich étale map is defined as being finite étale locally in a strict closed cover of the source and target by rational localisations. Thus, the theorem we want to prove breaks up into two parts and we start by treating finite étale maps. Before we begin the proof, let us mention that the key is that the finite étale site is invariant under tilting and under $\dagger$-reduction, see \cite[Prop.\ 3.2.11]{dRFF} for the latter, together with the fact that finite étale maps are $\dagger$-formally étale.

\begin{prop}
\label{prop:defis-etmaps}
Let $f: X\rightarrow Y$ be a finite étale map of affine Gelfand stacks. Then the induced map $f^\N: X^\N\rightarrow Y^\N$ is finite étale. In particular, $f^\N$ is a $!$-cover.
\end{prop}
\begin{proof}
We have to show that, for any totally disconnected nilperfectoid $A$ with a map $\GSpec A\rightarrow Y^\N$, the Gelfand stack $X^\N\times_{Y^\N} \GSpec A$ is affine and represented by a finite étale $A$-algebra. Thus, take any $A$-valued point of $Y^\N$ lying over an $A$-valued point $(D\subseteq Y_A, L\rightarrow A, \ol{A}\rightarrow K)$ of $\Q_p^\N$. Then the given $A$-valued point of $Y^\N$ corresponds to a map $D^+\rightarrow Y$, which by \cref{rem:syn-d+pullbackrings} is equivalent to compatible maps $D\rightarrow Y$ and $\GSpec\Cone(L\rightarrow A)\rightarrow Y$. In particular, pulling back $D$ to $Y_{\ol{A}}$, we obtain an untilt $\ol{A}^{\flat\sharp}$ of $\ol{A}$ with a map to $Y$ and pulling this back along the map $X\rightarrow Y$ will yield a finite étale perfectoid $\ol{A}^{\flat\sharp}$-algebra $\ol{B}^{\flat\sharp}$, which will correspond to a unique finite étale perfectoid $\ol{A}$-algebra $\ol{B}$, as the notation suggests, since the finite étale site is invariant under tilting. Finally, by \cite[Prop.\ 3.2.11]{dRFF}, this is the reduction of a unique nilperfectoid finite étale $A$-algebra $B$, which will in fact be totally disconnected (as it is finite étale over a totally disconnected ring).

Now recall from \cref{lem:prism-finetya} that the map $Y_B\rightarrow Y_A$ is finite étale. Thus, applying loc.\ cit.\ once again shows that $D\times_{Y_A} Y_B\cong D\times_Y X$ since both of these lift the finite étale $\ol{A}^{\flat\sharp}$-algebra $\ol{B}^{\flat\sharp}$ to a finite étale Gelfand stack over $D$ and hence the degree $1$ divisor $D\times_{Y_A} Y_B\subseteq Y_B$ is equipped with a map to $X$. In other words, we have constructed a $B$-point of $X^\prism$ fitting into a commutative diagram
\begin{equation*}
\begin{tikzcd}
\GSpec B\ar[r]\ar[d] & X^\prism\ar[d] \\
\GSpec A\ar[r] & Y^\prism\nospacepunct{\;.}
\end{tikzcd}
\end{equation*}

We will now further produce a similar diagram for Nygaardifications. For this, observe that the data we have produced so far already yields a map
\begin{equation}
\label{eq:syn-liftmapxngiven}
\GSpec\Cone(L\tensor_A \ol{B}\rightarrow \ol{B})\rightarrow X
\end{equation}
by construction of $X^\N$ and the only missing piece is a lift of this to a map
\begin{equation}
\label{eq:syn-liftmapxn}
\GSpec\Cone(L\tensor_A B\rightarrow B)\rightarrow X\;.
\end{equation}
However, note that, again, $\Cone(L\tensor_A \ol{B}\rightarrow\ol{B})$ is a finite étale $\Cone(L\tensor_A \ol{A}\rightarrow\ol{A})$-algebra that can uniquely be lifted to a finite étale $\Cone(L\rightarrow A)$-algebra, from which we conclude that there is an isomorphism
\begin{equation*}
\GSpec\Cone(L\tensor_A B\rightarrow B)\cong \GSpec \Cone(L\rightarrow A)\times_Y X
\end{equation*}
as both the left- and the right-hand side are possible lifts. Via the expression on the right-hand side we obtain a map as in (\ref{eq:syn-liftmapxn}) which is indeed compatible with the map (\ref{eq:syn-liftmapxngiven}) we are already given and thus, we overall obtain a commutative diagram
\begin{equation*}
\begin{tikzcd}
\GSpec B\ar[r]\ar[d] & X^\N\ar[d] \\
\GSpec A\ar[r] & Y^\N\nospacepunct{\;.}
\end{tikzcd}
\end{equation*}

Let us check that this diagram is cartesian. For this, we take a totally disconnected nilperfectoid $A$-algebra $C$ and have to verify that any $C$-valued point of $X^\N$ compatible with the given $A$-valued point of $Y^\N$ uniquely comes from a map $B\rightarrow C$ via the diagram above. This $C$-valued point of $X^\N$ corresponds to compatible maps $D\times_{Y_A} Y_C\rightarrow X$ and $\GSpec\Cone(L\tensor_A C\rightarrow C)\rightarrow X$. By compatibility with the given $A$-point of $Y^\N$, the first of these maps induces a map
\begin{equation*}
D\times_{Y_A} Y_C\rightarrow D\times_Y X\cong D\times_{Y_A} Y_B\;,
\end{equation*}
where the map $D\times_{Y_A} Y_C\rightarrow D$ is the first projection, and passing to $\dagger$-reductions on both sides induces a map $\ol{B}^{\flat\sharp}\rightarrow\ol{C}^{\flat\sharp}$ for the unique untilt $\ol{C}^{\flat\sharp}$ of $\ol{C}^\flat$ over $\ol{A}^{\flat\sharp}$. By tilting, we obtain a unique map $\ol{B}\rightarrow\ol{C}$ and hence a composite map $B\rightarrow \ol{C}$. However, as $A\rightarrow B$ is finite étale and hence $\dagger$-formally étale, the map $B\rightarrow \ol{C}$ lifts uniquely to a map $B\rightarrow C$ of $A$-algebras.

Finally, we have to check that the maps $D\times_{Y_A} Y_C\rightarrow X$ and $\GSpec\Cone(L\tensor_A C\rightarrow C)\rightarrow X$ are the same as the ones obtained via pullback from the $B$-valued point of $X^\N$ we have constructed above. However, by construction, they agree after passing to the $\dagger$-reduction of the source and then the fact that $X\rightarrow Y$ is $\dagger$-formally étale implies that they must be the same as they both lift the given $A$-point of $Y^\N$.
\end{proof}

Let us move on to the case of open covers and rational localisations. Now the key is that the Berkovich space is invariant under tilting and $\dagger$-nilpotent thickenings.

\begin{prop}
\label{prop:defis-openloc}
Let $X\rightarrow Y$ be an open embedding of derived Berkovich spaces. Then the induced map $X^\N\rightarrow Y^\N$ is an open embedding as well. Moreover, if $\{X_i\rightarrow Y\}_{i\in I}$ is a family of open localisations which jointly cover $Y$, then the maps $\{X_i^\N\rightarrow Y^\N\}_{i\in I}$ also jointly cover $Y^\N$. The same is true for rational localisations and strict closed covers.
\end{prop}
\begin{proof}
For $A$ totally disconnected nilperfectoid, take any $A$-valued point of $Y^\N$, which corresponds to an $A$-valued point $(D\subseteq Y_A, L\rightarrow A, \ol{A}\rightarrow K)$ of $\Q_p^\N$ together with compatible maps $D\rightarrow Y$ and $\GSpec\Cone(L\rightarrow A)\rightarrow Y$. In particular, the pullback of $D$ to $Y_{\ol{A}}$ yields an untilt $\ol{A}^{\flat\sharp}$ of $\ol{A}^\flat$ and we let $U$ be the open subspace of $\GSpec A$ defined by the open subset
\begin{equation}
\label{eq:syn-opensubspace}
|D|\times_{|Y|} |X|\cong \cal{M}(\ol{A}^{\flat\sharp})\times_{|Y|} |X|\cong \cal{M}(A)\times_{|Y|} |X|\subseteq \cal{M}(A)\;,
\end{equation}
where we have used that the Berkovich space is invariant under tilting and $\dagger$-nilpotent thickenings. We claim that there is a cartesian diagram
\begin{equation}
\label{eq:syn-opendiagram}
\begin{tikzcd}
U\ar[r]\ar[d] & X^\N\ar[d] \\
\GSpec A\ar[r] & Y^\N\nospacepunct{\;.}
\end{tikzcd}
\end{equation}

For this, first note that, by \cite[Lem.\ 4.4.5]{dRFF}, there is a map $Y_{\ol{A}}\rightarrow \cal{M}(A)$ with the property that the preimage of any closed subset $Z\subseteq \cal{M}(A)$ is given by $Y_{\ol{A}_Z}\subseteq Y_{\ol{A}}$, where $\ol{A}_Z$ denotes the rational localisation of $\ol{A}$ corresponding to $Z\subseteq\cal{M}(A)$. Thus, by construction, we have
\begin{equation*}
D\times_{Y_A} Y_U\cong D\times_{\cal{M}(A)} |U|\cong D\times_Y X
\end{equation*}
and, moreover,
\begin{equation*}
\GSpec\Cone(L\rightarrow A)\times_{\GSpec A} U\cong \GSpec\Cone(L\rightarrow A)\times_{\cal{M}(A)} |U|\cong \GSpec\Cone(L\rightarrow A)\times_Y X\;,
\end{equation*}
where we use that there is a map $\GSpec A\rightarrow\cal{M}(A)$, see \cite[§3.2]{dRFF}. These two isomorphisms provide compatible maps $D\times_{Y_A} Y_U\rightarrow X$ and $\GSpec\Cone(L\rightarrow A)\times_{\GSpec A} U\rightarrow X$ which induce a map $U\rightarrow X^\N$ making the diagram (\ref{eq:syn-opendiagram}) commute.

Let us check that the diagram is cartesian. For this, we take a totally disconnected $A$-algebra $C$ together with a $C$-point of $X^\N$ arising from compatible maps $D\times_{Y_A} Y_C\rightarrow X$ and $\GSpec\Cone(L\tensor_A C\rightarrow C)\rightarrow X$ which lifts the given $A$-point of $Y^\N$. However, letting $\ol{C}^{\flat\sharp}$ denote the untilt of $\ol{C}^\flat$ obtained from the $\dagger$-reduction of $D\times_{Y_A} Y_C$, we see that $\cal{M}(\ol{C}^{\flat\sharp})\cong \cal{M}(C)\rightarrow\cal{M}(A)$ must factor through $|U|$ by definition and hence $\GSpec C\rightarrow\GSpec A$ factors through $U$. Moreover, the given maps $D\times_{Y_A} Y_C\rightarrow X$ and $\GSpec\Cone(L\tensor_A C\rightarrow C)\rightarrow X$ must be the same as those arising via pullback from $U$: Indeed, both are compatible with the given $A$-point of $Y^\N$, i.e.\ they yield the same maps after postcomposing with $X\rightarrow Y$, and the fact that $X\rightarrow Y$ is an open embedding implies that derived Berkovich spaces over $X$ embed fully faithfully into derived Berkovich spaces over $Y$. This shows that (\ref{eq:syn-opendiagram}) is indeed cartesian and hence $X^\N\rightarrow Y^\N$ is an open localisations, as desired.

Now assume that $\{X_i\rightarrow Y\}_{i\in I}$ form an open cover. Pulling back to a totally disconnected nilperfectoid $\GSpec A$ as above, we get open subspaces $U_i$ of $\GSpec A$ with the property that
\begin{equation*}
|U_i|=\cal{M}(A)\times_{|Y|} |X_i|
\end{equation*}
by the above argument. As the maps $\{|X_i|\rightarrow |Y|\}_{i\in I}$ are jointly surjective, we deduce the same for the maps $\{|U_i|\rightarrow \cal{M}(A)\}_{i\in I}$, and we can even reduce to $|I|<\infty$ by compactness of $\cal{M}(A)$. Finally, by our standing qfd assumption, we know that $\cal{M}(A)$ has finite cohomological dimension and hence the family of maps $\{|U_i|\rightarrow\cal{M}(A)\}_{i\in I}$ forms a $!$-cover by \cite[Prop.\ II.1.1]{RealLLC}. As the maps $U_i\rightarrow \GSpec A$ are pulled back from the maps $|U_i|\rightarrow\cal{M}(A)$, this implies the claim.

Finally, we discuss the case of rational localisations. For this, we may reduce to the case of rational localisations of the form $\{|f|\leq 1\}$ or $\{|f|\geq 1\}$ for some global function $f$ on $Y$, and we only discuss the first case, the second one is similar. In that case, note that
\begin{equation*}
\{|f|\leq 1\}=\lim_{\epsilon>0}\, \{|f|<1+\epsilon\}\;,
\end{equation*}
hence we can write $X=\lim_{\epsilon>0} X_\epsilon$ for $X_\epsilon\coloneqq \{|f|<1+\epsilon\}\subseteq Y$. As $(-)^\N$ commutes with limits and is also compatible with open localisations by the argument above, we conclude that any $A$-valued point $\GSpec A\rightarrow Y^\N$ as above gives rise to a cartesian diagram
\begin{equation*}
\begin{tikzcd}
\lim_{\epsilon>0} U_\epsilon\ar[r]\ar[d] & X^\N\ar[d] \\
\GSpec A\ar[r] & Y^\N
\end{tikzcd}
\end{equation*}
with open subspaces $U_\epsilon\subseteq\GSpec A$. Moreover, by the description of these open subspaces from (\ref{eq:syn-opensubspace}), we see that $\lim_{\epsilon>0} U_\epsilon$ is a rational localisation of $\GSpec A$, and this is what we wanted to show. As the preservation of strict closed covers now follows from the preservation of open covers, this concludes the proof.
\end{proof}

As already discussed above, the two previous propositions imply \cref{thm:defis-berketalecover} since a Berkovich étale map is defined as being finite étale locally in a strict closed cover of the source and target by rational localisations.

\subsection{Nice coverability of $\Q_p^\N$}
\label{subsect:nicecoverqpn}

To end the section, we prove that $\Q_{p, [r, s]}^\N$ is nicely coverable in the sense of \cref{defi:perf-nicelycoverable} whenever $[r, s]\subseteq (p^{1/2}, p^{3/2})$. This will be the key technical input powering our comparison theorems in §\ref{sect:hkcomp} and §\ref{sect:proet}.

\begin{prop}
\label{prop:perf-zpncover}
For any $[r, s]\subseteq (p^{1/2}, p^{3/2})$, the Gelfand stack $\Q_{p, [r, s]}^\N$ is nicely coverable.
\end{prop}
\begin{proof}
First note that $\Q_{p, [r, s]}^\prism\cong Y_{\Spd\Q_p, [r, s]}^\dR$ by \cref{prop:defis-prismffdr} and recall that
\begin{equation*}
Y_{\Spd\Q_p}^\dR\cong Y_{\Q_p^\cycl}^\dR\,/\,\Z_p^{\times, \sm}\cong \lim_{q\mapsto q^p} (1+\overcirc{\DD})^\dR\setminus\{1\}\,\Big/\,\Z_p^{\times, \sm}\;.
\end{equation*}
The quotient stack
\begin{equation*}
\lim_{q\mapsto q^p} (1+\overcirc{\DD})^\dR\setminus\{1\}\,\Big/\,\Z_p^{\times, \sm}
\end{equation*}
clearly has a $!$-cover by $\lim_{q\mapsto q^p} (1+\overcirc{\DD})^\dR\setminus\{1\}$ whose \v{C}ech nerve is given by
\begin{equation*}
\left(\lim_{q\mapsto q^p} (1+\overcirc{\DD})^\dR\setminus\{1\}\right)\times (\Z_p^{\times, \sm})^\bullet\;,
\end{equation*}
where we note that $\O(\Z_p^{\times, \sm})=C^\sm(\Z_p^\times, \Q_p)$ is a flat $\Q_p$-algebra by \cref{lem:perf-vspflat}. Moreover, note that the endomorphism $q\mapsto q^p$ of $1+\overcirc{\DD}$ is $\dagger$-formally étale and therefore
\begin{equation*}
\lim_{q\mapsto q^p} (1+\overcirc{\DD})^\dR\setminus\{1\}\cong \lim_{q\mapsto q^p} (1+\overcirc{\DD})\setminus\{1\}\,\big/\,\G_m^\dagger\;.
\end{equation*}
As also $\O(\G_m^\dagger)=\colim_n \Q_p\langle p^{-n}T\rangle$ is a flat $\Q_p$-algebra, the \v{C}ech cover
\begin{equation*}
\left(\lim_{q\mapsto q^p} (1+\overcirc{\DD})\setminus\{1\}\right)\times (\G_m^\dagger)^\bullet\rightarrow \lim_{q\mapsto q^p} (1+\overcirc{\DD})^\dR\setminus\{1\}
\end{equation*}
provides a nice hypercover of the base change of $\lim_{q\mapsto q^p} (1+\overcirc{\DD})^\dR\setminus\{1\}$ to $[r, s]\subseteq (0, \infty)$ due to the fact that $\left(\lim_{q\mapsto q^p} (1+\overcirc{\DD})\setminus \{1\}\right)_{[r, s]}$ is affine and flat over $\Q_p$.

Indeed, in \cref{lem:pres-normpflatviaeps} below, we will compute that the preimage of $[r, s]\subseteq (0, \infty)$ in the stack $\lim_{q\mapsto q^p} (1+\overcirc{\DD})\setminus\{1\}$ is cut out by the inequality
\begin{equation*}
p^{-s}\leq \lim_n |q^{1/p^n}-1|^{p^{n-1}(p-1)}\leq p^{-r}\;,
\end{equation*}
where we normalise the norm such that $|p|=1/p$, and we note that this is actually equivalent to the inequalities
\begin{equation*}
p^{-s}\leq |q^{1/p^n}-1|^{p^2(p-1)}\leq p^{-r}
\end{equation*}
for all $n\geq 3$ due to $|(1+t)^p-1|=|t|^p$ whenever $p^{-1/p}<|t|<1$. Thus, we have
\begin{equation*}
\left(\lim_{q\mapsto q^p} (1+\overcirc{\DD})\setminus\{1\}\right)_{[r, s]}\cong \lim_{n\geq 3, q\mapsto q^p} \left(1+\overline{\T}\left(p^{-s/(p^{n-1}(p-1))}, p^{-r/(p^{n-1}(p-1))}\right)\right)\;,
\end{equation*}
where $\ol{\T}(a, b)$ denotes the overconvergent torus of inner radius $a$ and outer radius $b$ for any $a<b$, and the right-hand side is affine, static and flat over $\Q_p$. Putting everything together and applying \cref{lem:perf-finflatdim} and \cref{lem:perf-covers}, we see that $\Q_{p, [r, s]}^\prism$ is nicely coverable and in fact we have even found a nice hypercover where all the algebras that occur are static and flat over $\Q_p$.

Now note that $(\ol{\DD}/\ol{\T})^\dR$ is nicely coverable: Indeed, the \v{C}ech nerve of the map $\ol{\DD}^\dR\rightarrow (\ol{\DD}/\ol{\T})^\dR$ is given by $\ol{\DD}^\dR\times (\ol{\T}^\dR)^\bullet$; moreover, $\ol{\DD}^\dR$ and $\ol{\T}^\dR$ are themselves covered by $\ol{\DD}$ and $\ol{\T}$, respectively, with the respective \v{C}ech nerves being given by
\begin{equation*}
\ol{\DD}\times (\G_a^\dagger)^\bullet\;,\hspace{0.3cm}\text{resp.}\hspace{0.3cm} \ol{\T}\times (\G_m^\dagger)^\bullet\;.
\end{equation*}
As $\ol{\DD}, \ol{\T}, \G_a^\dagger$ and $\G_m^\dagger$ are all affine static Gelfand stacks which are flat over $\Q_p$ by \cref{lem:perf-vspflat}, this yields that $(\ol{\DD}/\ol{\T})^\dR$ is nicely coverable by \cref{lem:perf-bc} and \cref{lem:perf-covers}. Clearly, also $\ol{\DD}/\ol{\T}$ is nicely coverable by the \v{C}ech cover $\ol{\DD}\times \ol{\T}^\bullet$. Again, note that we have actually produced nice hypercovers which only involve static flat $\Q_p$-algebras.

Denoting by $X_\bullet\rightarrow \Q_{p, [r, s]}^\prism$ and $Y_{\bullet'}\rightarrow \ol{\DD}/\ol{\T}\times (\ol{\DD}/\ol{\T})^\dR$ the nice hypercovers we have produced above, the definition of $\Q_p^\N$ and \cref{lem:perf-covers} now reduce us to showing that all $X_\bullet\times_{(\ol{\DD}/\ol{\T})^\dR} Y_{\bullet'}$ are nicely coverable by $n$-truncated affines for some fixed $n\geq 0$. However, note that the maps $X_\bullet\rightarrow (\ol{\DD}/\ol{\T})^\dR$ and $Y_{\bullet'}\rightarrow (\ol{\DD}/\ol{\T})^\dR$ all canonically factor through $\ol{\DD}^\dR\rightarrow(\ol{\DD}/\ol{\T})^\dR$ and hence we obtain
\begin{equation*}
\begin{split}
X_\bullet\times_{(\ol{\DD}/\ol{\T})^\dR} Y_{\bullet'}&\cong X_\bullet\times_{\ol{\DD}^\dR} (\ol{\DD}^\dR\times_{(\ol{\DD}/\ol{\T})^\dR} \ol{\DD}^\dR)\times_{\ol{\DD}^\dR} Y_{\bullet'} \\
&\cong X_\bullet\times_{\ol{\DD}^\dR} (\ol{\T}^\dR\times \ol{\DD}^\dR)\times_{\ol{\DD}^\dR} Y_{\bullet'} \\
&\cong X_\bullet\times_{\ol{\DD}^\dR} (\ol{\T}^\dR\times Y_{\bullet'})\;,
\end{split}
\end{equation*}
where the map $\ol{\T}^\dR\times Y_{\bullet'}\rightarrow \ol{\DD}^\dR$ is given by the composition
\begin{equation*}
\ol{\T}^\dR\times Y_{\bullet'}\rightarrow \ol{\T}^\dR\times \ol{\DD}^\dR\xrightarrow{\text{mult}} \ol{\DD}^\dR\;.
\end{equation*}
By another application of \cref{lem:perf-covers}, we may now check everything after passing to the nice hypercover $\ol{\T}\times (\G_m^\dagger)^{\bullet''}$ of $\ol{\T}^\dR$ from above. However, note that we can make each term in this hypercover act on $Y_{\bullet'}$ such that the map $Y_{\bullet'}\rightarrow\ol{\DD}^\dR$ becomes equivariant (e.g.\ by first projecting onto $\ol{\T}$ and then letting $\ol{\T}$ act by multiplication on one of the factors $\ol{\DD}$ in $Y_{\bullet'}$) and, after twisting by this action, we may compute the fibre product above by replacing the map $(\ol{\T}\times (\G_m^\dagger)^{\bullet''})\times Y_{\bullet'}\rightarrow \ol{\DD}^\dR$ by 
\begin{equation*}
(\ol{\T}\times (\G_m^\dagger)^{\bullet''})\times Y_{\bullet'}\xrightarrow{\mathrm{pr}_2} Y_{\bullet'}\rightarrow \ol{\DD}^\dR\;.
\end{equation*}
In particular, we may focus on the fibre products $X_\bullet\times_{\ol{\DD}^\dR} Y_{\bullet'}$ because the additional factors $\ol{\T}$ and $\G_m^\dagger$ do not play any further role for our argument due to $\O(\ol{\T})$ and $\O(\G_m^\dagger)$ being static and flat over $\Q_p$.

Now we can do a similar trick again: the maps $X_\bullet\rightarrow\ol{\DD}^\dR$ and $Y_{\bullet'}\rightarrow\ol{\DD}^\dR$ admit a canonical factorisation through $\ol{\DD}\rightarrow\ol{\DD}^\dR$, which tells us that
\begin{equation*}
X_\bullet\times_{\ol{\DD}^\dR} Y_{\bullet'}\cong X_\bullet\times_{\ol{\DD}} (\G_a^\dagger\times Y_{\bullet'})\;,
\end{equation*}
where the map $\G_a^\dagger\times Y_{\bullet'}\rightarrow\ol{\DD}$ is given by
\begin{equation*}
\G_a^\dagger\times Y_{\bullet'}\rightarrow \G_a^\dagger\times\ol{\DD}\xrightarrow{\mathrm{add}} \ol{\DD}\;.
\end{equation*}
Finally, note that up to base changing to $[r, s]\subseteq (0, \infty)$, the above pullback is of the form
\begin{equation*}
\left(\lim_{q\mapsto q^p} (1+\overcirc{\DD}\setminus\{1\})\right)\times_{\ol{\DD}} (\G_a^\dagger\times \ol{\DD}\times \ol{\DD})\times (\,\dots)\;,
\end{equation*}
where the term in brackets is a finite product of affine Gelfand stacks which are flat and static over $\Q_p$ and hence is irrelevant for our argument.

Summarising, we are done once we know that, for some $n\geq 0$, the single (!) pullback
\begin{equation*}
\left(\lim_{q\mapsto q^p} (1+\overcirc{\DD}\setminus\{1\})\right)\times_{\ol{\DD}} (\G_a^\dagger\times\ol{\DD}\times \ol{\DD})\;,
\end{equation*}
where the map $\G_a^\dagger\times\ol{\DD}\times \ol{\DD}\rightarrow \ol{\DD}$ is given by $(x, y, z)\mapsto x+yz$, is an affine $n$-truncated Gelfand stack after base changing to $[r, s]\subseteq (0, \infty)$. However, as we have already shown that $\left(\lim_{q\mapsto q^p} (1+\overcirc{\DD})\setminus\{1\}\right)_{[r, s]}$ is affine and static, this just follows from the fact that the map $\G_a^\dagger\times \ol{\DD}\times\ol{\DD}\rightarrow \ol{\DD}$ is flat: On rings, this is given by
\begin{equation*}
\Q_p\langle T\rangle_{\leq 1}\rightarrow \Q_p\{X\}^\dagger\langle Y, Z\rangle_{\leq 1}\;, \hspace{0.3cm} T\mapsto X+YZ
\end{equation*}
and both the map
\begin{equation*}
\Q_p\langle T\rangle_{\leq 1}\rightarrow \Q_p\{X\}^\dagger\langle W\rangle_{\leq 1}\;, \hspace{0.3cm} T\mapsto X+W
\end{equation*}
as well as the map
\begin{equation*}
\Q_p\langle W\rangle_{\leq 1}\rightarrow \Q_p\langle Y, Z\rangle_{\leq 1}\;, \hspace{0.3cm} W\mapsto YZ
\end{equation*}
are flat. Indeed, for the first map, this follows from $\Q_p\{X\}^\dagger\langle W\rangle_{\leq 1}\cong \Q_p\{X\}^\dagger\langle X+W\rangle_{\leq 1}$ and flatness of $\Q_p\{X\}^\dagger$ over $\Q_p$ while, for the second map, we observe that, as a $\Q_p\langle W\rangle_{\leq 1}$-module, we have
\begin{equation*}
\Q_p\langle Y, Z\rangle_{\leq 1}\cong \Q_p\langle W, Y\rangle_{\leq 1}\oplus Z\Q_p\langle W, Z\rangle_{\leq 1}
\end{equation*}
and then the claim follows from flatness of $\Q_p\langle Y\rangle_{\leq 1}$ over $\Q_p$ by base change.
\end{proof}

\begin{rem}
Modifying the above argument slightly even yields nice coverability of $\Q_{p, [r, s]}^\N$ for any $[r, s]\subseteq (1, \infty)$. Moreover, one can also prove the same for $[r, s]\subseteq (0, p)$ using that 
\begin{equation*}
\Q_{p, (0, p)}^\prism\cong \left(\lim_{q\mapsto q^p} (1+\overcirc{\DD})\setminus \{1\}\right)_{(0, p)}\,\Big/\,\Z_p^{\times, \la}\;,
\end{equation*}
which one may deduce by combining \cref{prop:defis-xprismxdiv1} and \cref{prop:prism-div1presentations}. However, we will only need the case $[r, s]\subseteq (p^{1/2}, p^{3/2})$ in the sequel.
\end{rem}

\newpage

\section{Nygaardifications of smooth derived Berkovich spaces}
\label{sect:pres}

Our aim in this section is to find explicit local presentations of some Nygaardifications. While the need for this might not be fully apparent at the moment, having such explicit presentations will be very useful later when we want to show that $\Q_p^\N$ is cohomologically smooth and that the functor $X\mapsto X^\N$ carries rigid smooth maps of derived Berkovich spaces to cohomologically smooth maps and quasicompact locally finitely presented maps of derived Berkovich spaces to weakly cohomologically proper maps.

Let us give a brief idea of what we would like to do. We immediately warn the reader that this will ultimately \emph{not} quite work and cause us to work with a slight variant of $\Q_p^\N$ instead, but it still seems valuable to explain why we introduce this variant and why it is easier to handle than $\Q_p^\N$ itself. Trying to find an explicit local presentation of $\Q_p^\N$, we first think about $\Q_p^\prism$: Recalling the presentations
\begin{equation*}
\Q_p^{\Div^1}\cong \lim_{q\mapsto q^p} (1+\overcirc{\DD})\setminus\{1\}\,\Big/\,\Q_p^{\times, \la}\;, \hspace{0.5cm} \Q_p^\HK\cong \lim_{q\mapsto q^p} (1+\overcirc{\DD})^\dR\setminus\{1\}\,\Big/\,\Q_p^{\times, \sm}
\end{equation*}
from \cref{prop:prism-div1presentations} and \cref{cor:prism-hkpresentations} and the relation between $\Q_p^\prism$ and the stacks $\Q_p^{\Div^1}$ and $\Q_p^\HK$, we are led to assume that there should be a surjection
\begin{equation}
\label{eq:pres-rhoqdr}
\rho_{q\dR}: \lim_{q\mapsto q^p} (1+\overcirc{\DD})\setminus\{1\}\rightarrow \Q_p^\prism\;.
\end{equation}
Indeed, such a map exists: For any totally disconnected nilperfectoid ring $A$ and any sequence $(q, q^{1/p}, q^{1/p^2}, \dots)$ of $p$-power roots in $1+A^{<1}$ not equal to $(1, 1, 1, \dots)$, i.e.\ any $A$-point of the source, we obtain a global section $[q^\flat]$ of $Y_A$ by gluing $[q^\flat]\in\mathbb{A}_\inf(\ol{A})$, where $q^\flat=(q, q^{1/p}, q^{1/p^2}, \dots)\in \ol{A}^{\flat\circ}$, and the element $(q, q^{1/p}, q^{1/p^2}, \dots)\in \prod_{n=0}^\infty A$; moreover, $[q^\flat]\in H^0(Y_A, \O)$ has canonical $p$-power roots obtained by shifting the sequence $(q, q^{1/p}, q^{1/p^2}, \dots)$ we started with. Then the section
\begin{equation*}
\frac{[q^\flat]-1}{[q^\flat]^{1/p}-1}=1+[q^\flat]^{1/p}+\dots+[q^\flat]^{(p-1)/p}\in H^0(Y_A, \O)
\end{equation*}
cuts out a degree $1$ Cartier divisor in $Y_A$ and this defines the map (\ref{eq:pres-rhoqdr}) above. We remark that one should think of $\rho_{q\dR}$ as analogous to the surjection from the $q$-de Rham prism onto the algebraic prismatisation of $\Z_p$ in the sense of Bhatt--Lurie and Drinfeld, see \cite[§3.8]{APC}.

One checks that the map (\ref{eq:pres-rhoqdr}) is compatible with the presentations of $\Q_p^{\Div^1}$ and $\Q_p^\HK$ above and hence, base changing to $(1, \infty)\subseteq (0, \infty)$, we obtain a map
\begin{equation*}
\left(\lim_{q\mapsto q^p} (1+\overcirc{\DD})^\dR\setminus\{1\}\right)_{(1, \infty)}\rightarrow \Q_{p, (1, \infty)}^\prism
\end{equation*}
which presents the target as the quotient of the source by $\Z_p^{\times, \sm}$. How can we obtain a presentation of $\Q_{p, (1, \infty)}^\N$ from here? For this, let us first observe that the postcomposition of the map above with $\widetilde{\mu}: \Q_p^\prism\rightarrow (\ol{\DD}/\ol{\T})^\dR$ can alternatively be factored as
\begin{equation*}
\left(\lim_{q\mapsto q^p} (1+\overcirc{\DD})^\dR\setminus\{1\}\right)_{(1, \infty)}\rightarrow \ol{\DD}^\dR\rightarrow (\ol{\DD}/\ol{\T})^\dR\;,
\end{equation*}
where the first map is given by
\begin{equation*}
(q, q^{1/p}, q^{1/p^2}, \dots)\mapsto \frac{q^p-1}{q-1}=1+q+\dots+q^{p-1}\;.
\end{equation*}

Thus, ignoring the quotient by $\Z_p^{\times, \sm}$ for the moment, we can compute the successive pullback
\begin{equation*}
\begin{tikzcd}
(\ol{\DD}_+\times\ol{\DD}^\dR_-)\times_{\ol{\DD}^\dR} \left(\lim_{q\mapsto q^p} (1+\overcirc{\DD})^\dR\setminus \{1\}\right)_{(1, \infty)}\,\Big/\,\ol{\T}\ar[r]\ar[d] & \left(\lim_{q\mapsto q^p} (1+\overcirc{\DD})^\dR\setminus \{1\}\right)_{(1, \infty)}\ar[d] \\
(\ol{\DD}_+\times\ol{\DD}^\dR_-)/\ol{\T}\ar[r, "\mathrm{mult}"]\ar[d] & \ol{\DD}^\dR\ar[d] \\
\ol{\DD}_+/\ol{\T}\times (\ol{\DD}_-/\ol{\T})^\dR\ar[r, "\mathrm{mult}"] & \ol{\DD}^\dR_+/\ol{\T}
\end{tikzcd}
\end{equation*}
and the top left term, in which $\ol{\T}$ only acts on the product $\ol{\DD}_+\times\ol{\DD}^\dR_-$, will be a $\Z_p^{\times, \sm}$-torsor over $\Q_{p, (1, \infty)}^\N$. Can we also incorporate the $\Z_p^{\times, \sm}$-action? In principle, we would now like to write a formula like
\begin{equation*}
\Q_{p, (1, \infty)}^\N\overset{?}{\cong} (\ol{\DD}\times\ol{\DD}^\dR)\times_{\ol{\DD}^\dR} \left(\lim_{q\mapsto q^p} (1+\overcirc{\DD})^\dR\setminus \{1\}\right)_{(1, \infty)}\,\Big/\,\Z_p^{\times, \la}\coprod_{\G_m^\dagger}^\AbGrp \ol{\T}\;,
\end{equation*}
where $\Z_p^{\times, \la}$ acts on $\lim_{q\mapsto q^p} (1+\overcirc{\DD})^\dR$ via the quotient $\Z_p^{\times, \la}\rightarrow\Z_p^{\times, \sm}$ while $\ol{\T}$ acts trivially on $\lim_{q\mapsto q^p} (1+\overcirc{\DD})^\dR$ and by multiplication or division on $\ol{\DD}$ or $\ol{\DD}^\dR$, respectively. However, for the formula above to make sense, we at least also have to prescribe an action of $\Z_p^{\times, \la}$ on $\ol{\DD}$ as $\G_m^\dagger$ already acts nontrivially on $\ol{\DD}$ through the action of $\ol{\T}$. Moreover, this action of $\Z_p^{\times, \la}$ on the product $\ol{\DD}\times \ol{\DD}^\dR$ has to be such that the induced action on $\ol{\DD}^\dR$ under the multiplication map agrees with the action induced from $\lim_{q\mapsto q^p} (1+\overcirc{\DD})^\dR$. This latter task, however, seems at least very hard, if not impossible: The $\Z_p^{\times, \sm}$-action on $\ol{\DD}$ induced via the map $\lim_{q\mapsto q^p} (1+\overcirc{\DD})^\dR\rightarrow \ol{\DD}^\dR$ is given by
\begin{equation*}
\gamma.\frac{q^p-1}{q-1}=\frac{q^{p\gamma}-1}{q^\gamma-1}
\end{equation*}
and it is not even clear how to rewrite the right-hand side purely in terms of $\frac{q^p-1}{q-1}$, i.e.\ without explicit reference to $q$. 

Summarising, the trouble with the above approach seems to come from the fact that the map
\begin{equation*}
\lim_{q\mapsto q^p} (1+\overcirc{\DD})\setminus\{1\}\rightarrow \ol{\DD}\;, \hspace{0.3cm} (q, q^{1/p}, q^{1/p^2}, \dots)\mapsto \frac{q^p-1}{q-1}=1+q+\dots+q^{p-1}
\end{equation*}
induced by (\ref{eq:pres-rhoqdr}) does not play very well with the $\Z_p^{\times, \la}$-action on the source. However, recall from \cref{prop:prism-div1presentations} that the map $\lim_{q\mapsto q^p} (1+\overcirc{\DD})\setminus\{1\}\rightarrow\Q_p^{\Div^1}$ is slightly different from the one to $\Q_p^\prism$: Instead of sending $(q, q^{1/p}, q^{1/p^2}, \dots)$ to the divisor cut out by $1+[q^\flat]^{1/p}+\dots+[q^\flat]^{(p-1)/p}$ in $Y_A$, it sends it the divisor cut out by $\log[q^\flat]$ in $\FF_A$! In particular, the composition
\begin{equation*}
\lim_{q\mapsto q^p} (1+\overcirc{\DD})\setminus\{1\}\rightarrow \Q_p^{\Div^1}\xrightarrow{\mu} \A^1/\G_m
\end{equation*}
factors through
\begin{equation*}
\lim_{q\mapsto q^p} (1+\overcirc{\DD})\setminus\{1\}\rightarrow \A^1\;, \hspace{0.3cm} (q, q^{1/p}, q^{1/p^2}, \dots)\mapsto \log q
\end{equation*}
which now interacts a lot better with the $\Z_p^{\times, \la}$-action: We have $\log q^\gamma=\gamma\log q$ and hence the map above intertwines the $\Z_p^{\times, \la}$-action on the source with the multiplication action on the target.

Overall, we see that we would be much better off if we could replace $1+q+\dots+q^p$ by $\log q$. However, note that the latter is not even always contained in $\ol{\DD}\subseteq \A^1$ and, in particular, will generally differ from $1+q+\dots+q^p$ by more than just multiplication by a unit of norm $1$; nevertheless, at least after restricting further to $(p^{1/2}, p^{3/2})\subseteq (1, \infty)$, we will be able to arrange that the two differ by a unit. Concluding the discussion, we see that our situation would overall be much improved if we replaced all occurrences of $\ol{\DD}/\ol{\T}$ in the definition of $\Q_p^\N$ by occurrences of $\A^1/\G_m$, which is precisely what we will do below. 

This yields a slightly different stack $\Q_p^{\ol{\N}}$ which we call the \emph{reduced Nygaardification} of $\Q_p$, for which we will construct a map
\begin{equation*}
\widetilde{\Q_{p\mathrlap{, (p^{1/2}, p^{3/2})}}^{\ol{\N}}}\hphantom{\scriptstyle{(p^{1/2}, p^{3/2})}}\coloneqq (\A^1\times \A^1)\times_{\A^1} \left(\lim_{q\mapsto q^p} (1+\overcirc{\DD})\setminus \{1\}\right)_{(p^{1/2}, p^{3/2})}\rightarrow \Q_{p, (p^{1/2}, p^{3/2})}^{\ol{\N}}
\end{equation*}
which presents the target as a quotient of the source by an explicit group stack, which we will describe below. Also note that the effect of passing from $\Q_p^\N$ to $\Q_p^{\ol{\N}}$ is not too grave: Basically, the map $\Q_p^\N\rightarrow\Q_p^{\ol{\N}}$ is a $(0, \infty)$-torsor and thus, in particular, understanding $\Q_p^{\ol{\N}}$ is still enough to prove cohomological smoothness of $\Q_p^\Syn$. Also observe that we really only need to care about $\Q_{p, (p^{1/2}, p^{3/2})}^\N$: Indeed, using that
\begin{equation*}
\Q_{p, (0, p)}^\N\cong \Q_{p, (0, p)}^\prism\times [0, 1]\;, \hspace{0.3cm} \Q_{p, (p, \infty)}^\N\cong \Q_{p, (p, \infty)}^\prism\times [0, 1]
\end{equation*}
via $\pi$ by \cref{prop:defis-utneq0}, understanding $\Q_p^\N$ over the locus $(p^{1/2}, p^{3/2})\subseteq (1, \infty)$ suffices as we already have a good understanding of $\Q_{p, (0, p)}^\prism$ and $\Q_{p, (p, \infty)}^\prism$ due to \cref{prop:defis-xprismxdiv1} and \cref{prop:defis-prismffdr}, respectively.

After this lengthy discussion of $\Q_p^\N$, let us briefly outline what happens in the rest of the section: One easily sees that the definition of the derived Berkovich space $D^+$ over $\Q_p^\N$ in fact descends to $\Q_p^{\ol{\N}}$ and, using this, one can define $X^{\ol{\N}}$ for any Gelfand stack $X$. Using the results of the discussion above, we will then be able to control $X^{\ol{\N}}\rightarrow Y^{\ol{\N}}$ for arbitrary Berkovich smooth maps $X\rightarrow Y$ between derived Berkovich spaces. Indeed, by compatibility of $X\mapsto X^{\ol{\N}}$ with Berkovich étale localisation, which is proved exactly the same as in \cref{thm:defis-berketalecover}, this formally reduces to understanding $(\G_m^n)^{\ol{\N}}\rightarrow\Q_p^{\ol{\N}}$; by commutation of $X\mapsto X^{\ol{\N}}$ with limits, we can furthermore reduce to $n=1$. In this latter case, we are going to construct a map
\begin{equation*}
\lim_{x\mapsto x^p} \G_m\times \widetilde{\Q_{p\mathrlap{, (p^{1/2}, p^{3/2})}}^{\ol{\N}}}\hphantom{\scriptstyle{(p^{1/2}, p^{3/2})}}\rightarrow \G_m^{\ol{\N}}\times_{\Q_p^{\ol{\N}}}\times \widetilde{\Q_{p\mathrlap{, (p^{1/2}, p^{3/2})}}^{\ol{\N}}}\hphantom{\scriptstyle{(p^{1/2}, p^{3/2})}}
\end{equation*}
which presents the target as the quotient of the source by an explicit group stack over $\widetilde{\Q_{p\mathrlap{, (p^{1/2}, p^{3/2})}}^{\ol{\N}}}\hphantom{\scriptstyle{(p^{1/2}, p^{3/2})}}$, which we shall describe.

Finally, we specialise to the locus $\{|t|=0\}$, i.e.\ to the stack $X^{\HT, \dagger, +}$. Note that this closed substack of $X^\N$ is already very interesting (and in fact captures almost all of the additional information that $X^\N$ holds compared to $X^\prism$): Namely, it provides a deformation from $X^{\dR, +}_{|t|=0}\cong X^\N_{|u|=|t|=0}$ to $X^{\HT, \dagger}\cong X^\N_{|t|=0, |u|=1}$. To understand this deformation, we shall want slightly finer information compared to what we obtain from the discussion above: Namely, we will establish an explicit presentation of the stack $X^{\HT, \dagger, +}$ itself in the case where $X$ is equipped with a Berkovich étale map $X\rightarrow\ol{\T}^n$. Using this, we will be able to give an explicit description of the category of perfect complexes on $X^{\HT, \dagger, +}$ for such $X$, which will be useful later when we want to describe vector bundles on $X^\Syn$ and their cohomology.

\subsection{The Nygaardification of $\Q_p$}

Let us now carry out the strategy to describe $\Q_p^\N$ more explicitly laid out above. We start by defining $\Q_p^{\ol{\N}}$ as mentioned previously.

\begin{defi}
The \emph{reduced Nygaardification} of $\Q_p$ is the Gelfand stack $\Q_p^{\ol{\N}}$ defined as the pullback
\begin{equation*}
\begin{tikzcd}
\Q_p^{\ol{\N}}\ar[r, "\pi"]\ar[d, "{(t, u)}", swap] & \Q_p^\prism\ar[d, "\widetilde{\mu}"] \\
\A^1_+/\G_m\times (\A^1_-/\G_m)^\dR\ar[r, "\mathrm{mult}"] & (\A^1_+/\G_m)^\dR\nospacepunct{\;.}
\end{tikzcd}
\end{equation*}
\end{defi}

We first want to clarify the relationship between $\Q_p^{\ol{\N}}$ and $\Q_p^\N$. For this, we introduce another definition.

\begin{defi}
The \emph{extended Nygaardification} of $\Q_p$ is the Gelfand stack $\Q_p^{\N, \mathrm{ext}}$ defined as the pullback
\begin{equation*}
\begin{tikzcd}
\Q_p^{\N, \mathrm{ext}}\ar[r, "\pi"]\ar[d, "{(t, u)}", swap] & \Q_p^\prism\ar[d, "\widetilde{\mu}"] \\
\A^1_+/\ol{\T}\times (\A^1_-/\ol{\T})^\dR\ar[r, "\mathrm{mult}"] & (\A^1_+/\ol{\T})^\dR\nospacepunct{\;.}
\end{tikzcd}
\end{equation*}
\end{defi}

Recalling that the map $\widetilde{\mu}: \Q_p^\prism\rightarrow (\A^1/\ol{\T})^\dR$ factors through $(\overcirc{\DD}/\ol{\T})^\dR$, we observe that the stack $\Q_p^{\N, \mathrm{ext}}$ is supported over the locus $\{|ut|<1\}$. Moreover, we immediately see that
\begin{equation}
\label{eq:pres-qpnfromqpnext}
\Q_p^\N\cong \Q_{p, |t|\leq 1, |u|\leq 1}^{\N, \mathrm{ext}}\;.
\end{equation}
The relation between $\Q_p^\N$ and $\Q_p^{\ol{\N}}$ is thus fully explained by the following result:

\begin{lem}
\label{lem:pres-qpredntorsor}
The canonical map $\Q_p^{\N, \mathrm{ext}}\rightarrow \Q_p^{\ol{\N}}$ induced by $\A^1/\ol{\T}\rightarrow\A^1/\G_m$ is a $(0, \infty)$-torsor.
\end{lem}
\begin{proof}
We first note that there is a pullback square
\begin{equation*}
\begin{tikzcd}
(\A^1_+/\ol{\T}\times(\A^1_-)^\dR)/\G_m^\dR\ar[r, "\mathrm{mult}"]\ar[d] & (\A^1_+/\ol{\T})^\dR\ar[d] \\
\A^1_+/\G_m\times (\A^1_-/\G_m)^\dR\ar[r, "\mathrm{mult}"] & (\A^1_+/\G_m)^\dR\nospacepunct{\;,}
\end{tikzcd}
\end{equation*}
where $\G_m^\dR$ acts diagonally on $\A^1_+/\ol{\T}\times (\A^1_-)^\dR$ in the top left corner. From this, we conclude that there is a cartesian diagram
\begin{equation*}
\begin{tikzcd}
\Q_p^{\ol{\N}}\ar[r, "\pi"]\ar[d, "{(t, u)}", swap] & \Q_p^\prism\ar[d, "\widetilde{\mu}"] \\
(\A^1_+/\ol{\T}\times(\A^1_-)^\dR)/\G_m^\dR\ar[r, "\mathrm{mult}"] & (\A^1_+/\ol{\T})^\dR\nospacepunct{\;,}
\end{tikzcd}
\end{equation*}
from which we in turn deduce a pullback square
\begin{equation*}
\begin{tikzcd}
\Q_p^{\N, \mathrm{ext}}\ar[r]\ar[d, "{(t, u)}", swap] & \Q_p^{\ol{\N}}\ar[d, "{(t, u)}"] \\
\A^1_+/\ol{\T}\times (\A^1_-/\ol{\T})^\dR\ar[r] & (\A^1_+/\ol{\T}\times(\A^1_-)^\dR)/\G_m^\dR\nospacepunct{\;.}
\end{tikzcd}
\end{equation*}
Observing that the lower horizontal map is a torsor for $\G_m^\dR/\ol{\T}^\dR\cong (0, \infty)$, where the isomorphism is due to \cref{lem:defis-rhvariant}, the claim follows.
\end{proof}

Note that, just as $\Q_p^\N$, the stacks $\Q_p^{\ol{\N}}$ and $\Q_p^{\N, \mathrm{ext}}$ are equipped with a radius map $\kappa$ towards $(0, \infty)$ via the respective maps $\pi$ towards $\Q_p^\prism$. We remark that the argument above even shows that $\Q_{p, [r, s]}^{\N, \mathrm{ext}}$ is a $(0, \infty)$-torsor over $\Q_{p, [r, s]}^{\ol{\N}}$ for any $[r, s]\subseteq (0, \infty)$. Similarly, the isomorphism (\ref{eq:pres-qpnfromqpnext}) is compatible with base change to any $[r, s]\subseteq (0, \infty)$ as well.

Moving on, recall that our goal is to obtain an explicit presentation of $\Q_{p, (p^{1/2}, p^{3/2})}^{\ol{\N}}$ starting from the map
\begin{equation*}
\rho_{q\dR}: \lim_{q\mapsto q^p} (1+\overcirc{\DD})\setminus\{1\}\rightarrow \Q_p^\prism
\end{equation*}
from (\ref{eq:pres-rhoqdr}), which we now study in slightly more detail. Let us point out here that the source should really be more accurately denoted by $\lim_{q\mapsto q^p} (1+\overcirc{\DD})\setminus\{(1, 1, 1, \dots)\}$, which in turn refers to the colimit
\begin{equation*}
\colim_n \lim (1+\overcirc{\DD}\xleftarrow{q^p\mapsfrom q} 1+\overcirc{\DD}\xleftarrow{q^p\mapsfrom q}\dots\xleftarrow{q^p\mapsfrom q} 1+\overcirc{\DD}\xleftarrow{q^p\mapsfrom q} 1+\overcirc{\DD}^\times\xleftarrow{q^p\mapsfrom q}1+\overcirc{\DD}^\times\xleftarrow{q^p\mapsfrom q}\dots)\;,
\end{equation*}
where the first occurrence of the pointed open unit disk $1+\overcirc{\DD}^\times$ around $1$ is in the $n$-th spot.

As we will eventually want to base change the map $\rho_{q\dR}$ to $(p^{1/2}, p^{3/2})$, our first task is to make the composite
\begin{equation*}
\lim_{q\mapsto q^p} (1+\overcirc{\DD})\setminus\{1\}\xrightarrow{\rho_{q\dR}} \Q_p^\prism\xrightarrow{\kappa} (0, \infty)
\end{equation*}
more explicit. For this, recall that the radius map on $\Q_p^\prism$ is induced by the one on $Y_{\Spd\Q_p}$, which in turn is obtained by descent from the map $Y_{\C_p}\rightarrow (0, \infty)$ defined by
\begin{equation*}
\kappa^{-1}([r, s])=\{|p|^s\leq |[p^\flat]|\leq |p|^r\}\subseteq Y_{\C_p}\;.
\end{equation*}
Moreover, unwinding definitions, we find that the map
\begin{equation*}
\lim_{q\mapsto q^p} (1+\overcirc{\DD})^\dR\setminus\{1\}\xrightarrow{\rho_{q\dR}^\dR} (\Q_p^\prism)^\dR\cong Y_{\Spd\Q_p}^\dR
\end{equation*}
identifies with the composition
\begin{equation*}
\lim_{q\mapsto q^p} (1+\overcirc{\DD})^\dR\setminus\{1\}\cong Y_{\Q_p^\cycl}^\dR\rightarrow Y_{\Spd\Q_p}^\dR\;,
\end{equation*}
where the first isomorphism is given by $q^{1/p^n}\mapsto [\epsilon]^{1/p^n}$ for $\epsilon=(1, \zeta_p, \zeta_{p^2}, \dots)$, as usual, see \cite[Lem.\ 6.2.1.(2)]{dRFF}. In other words, we have to understand $|[p^\flat]|$ on $Y_{\C_p}$ in terms of the $p$-power roots of $[\epsilon]$.

\begin{lem}
\label{lem:pres-normpflatviaeps}
For any continuous seminorm on $Y_{\C_p}$, we have
\begin{equation*}
|[p^\flat]|=\lim_n |[\epsilon]^{1/p^n}-1|^{p^{n-1}(p-1)}\;.
\end{equation*}
\end{lem}
\begin{proof}
First note that $(\epsilon-1)^{p-1}$ and $(p^\flat)^p$ differ by a unit in $\O_{\C_p^\flat}$. Indeed, the latter is a valuation ring whose valuation is given by $v(x)\coloneqq v(x^\sharp)$ and normalising the valuation on $\O_{\C_p}$ such that $v(p)=1$ yields
\begin{equation*}
v(\epsilon-1)=v\left(\lim_n (\zeta_{p^n}-1)^{p^n}\right)=\lim_n p^nv(\zeta_{p^n}-1)=\frac{p}{p-1}
\end{equation*}
while $v(p^\flat)=1$, which proves the claim. Thus, also $[p^\flat]^p$ and $[\epsilon-1]^{p-1}$ differ by a unit in $A_\inf$ and hence they have the same norm on $Y_{\C_p}$.

It remains to show that
\begin{equation*}
|[\epsilon-1]|=\lim_n |[\epsilon]^{1/p^n}-1|^{p^n}\;.
\end{equation*}
However, this actually holds as an equality in $A_\inf$ even before taking norms: Indeed, for each $n\geq 0$, the element $[\epsilon]^{1/p^n}-1\in A_\inf=W(\O_{\C_p^\flat})$ lifts $(\epsilon-1)^{1/p^n}=\epsilon^{1/p^n}-1\in \O_{\C_p^\flat}$ and this implies the claim.
\end{proof}

Normalising the norm such that $|p|=1/p$, the previous lemma shows that the preimage of $(p^{1/2}, p^{3/2})\subseteq (0, \infty)$ in $\lim_{q\mapsto q^p} (1+\overcirc{\DD})\setminus\{1\}$ under the radius map is described by
\begin{equation*}
p^{-p^{3/2}}<\lim_n |q^{1/p^n}-1|^{p^{n-1}(p-1)}<p^{-p^{1/2}}\;.
\end{equation*}
Noting that $|(1+t)^p-1|=|t|^p$ holds whenever $p^{-1/p}<|t|<1$, we see that the above is equivalent to requiring
\begin{equation*}
p^{-1/(p^{1/2}(p-1))}<|q^{1/p^3}-1|<p^{-1/(p^{3/2}(p-1))}\;.
\end{equation*}

\begin{cor}
\label{cor:pres-rhoqdrradius}
The base change of $\lim_{q\mapsto q^p} (1+\overcirc{\DD})\setminus\{1\}$ to $(p^{1/2}, p^{3/2})\subseteq (0, \infty)$ sits in a pullback diagram
\begin{equation*}
\begin{tikzcd}
\left(\lim_{q\mapsto q^p} (1+\overcirc{\DD})\setminus\{1\}\right)_{(p^{1/2}, p^{3/2})}\ar[r]\ar[d] & \lim_{q\mapsto q^p} (1+\overcirc{\DD})\setminus\{1\}\ar[d, "\mathrm{pr}_4"] \\
1+\overcirc{\T}(p^{-1/(p^{1/2}(p-1))}, p^{-1/(p^{3/2}(p-1))})\ar[r] & 1+\overcirc{\DD}\;,
\end{tikzcd}
\end{equation*}
where $\overcirc{\T}(a, b)$ denotes the open torus of inner radius $a$ and outer radius $b$ for any $a<b$.
\end{cor}

Our next task is to understand the composition
\begin{equation*}
\left(\lim_{q\mapsto q^p} (1+\overcirc{\DD})\setminus\{1\}\right)_{(p^{1/2}, p^{3/2})}\xrightarrow{\rho_{q\dR}} \Q_{p, (p^{1/2}, p^{3/2})}^\prism\xrightarrow{\widetilde{\mu}} (\ol{\DD}/\ol{\T})^\dR\rightarrow (\A^1/\G_m)^\dR\;.
\end{equation*}
For this, note that the first map factors through the projection onto the analytic de Rham stack of the source since $\Q_{p, (p^{1/2}, p^{3/2})}^\prism\cong (\Q_{p, (p^{1/2}, p^{3/2})}^\prism)^\dR$ by \cref{prop:defis-prismffdr} and observe that there is a commutative diagram
\begin{equation*}
\hspace{-3cm}
\begin{tikzcd}
\lim_{q\mapsto q^p} (1+\overcirc{\DD})^\dR\setminus\{1\}\ar[d]\ar[r, "\cong"] & Y_{\Q_p^\cycl}^\dR\ar[d]\ar[r] & Y_{\Spd\Q_p}^\dR\mathrlap{\cong (\Q_p^\prism)^\dR}\ar[d]  \\
\left(\lim_{q\mapsto q^p} (1+\overcirc{\DD})^\dR\setminus\{1\}\right)\,\Big/\,\phi^\Z\ar[r, "\cong"] &  \FF_{\Q_p^\cycl}^\dR\ar[r] & \FF_{\Spd\Q_p}^\dR\mathrlap{=\Q_p^\HK\cong (\Q_p^{\Div^1})^\dR\;,}
\end{tikzcd}
\end{equation*}
in which the top composition identifies with $\rho_{q\dR}^\dR$ while the bottom composition sends $(q, q^{1/p}, q^{1/p^2}, \dots)$ to the degree $1$ Cartier divisor in $\FF_{\ol{A}}$ cut out by $\log[q^\flat]$. Further noting that the diagram
\begin{equation*}
\begin{tikzcd}
\mathllap{\Q_{p, (p^{1/2}, p^{3/2})}^\prism\cong\,}(\Q_{p, (p^{1/2}, p^{3/2})}^\prism)^\dR\ar[d]\ar[r, "\widetilde{\mu}"] & (\ol{\DD}/\ol{\T})^\dR\ar[d] \\
\FF_{\Spd\Q_p}^\dR\ar[r, "\mu^\dR"] & (\A^1/\G_m)^\dR 
\end{tikzcd}
\end{equation*}
commutes, where the lower horizontal map is defined by pulling back a degree $1$ Cartier divisor in $\FF_{\ol{A}}$ along $\iota: \GSpec\ol{A}\rightarrow\FF_{\ol{A}}$ for any totally disconnected nilperfectoid $A$, we conclude:

\begin{lem}
\label{lem:pres-rhoqdrlogq}
The composition
\begin{equation*}
\left(\lim_{q\mapsto q^p} (1+\overcirc{\DD})\setminus\{1\}\right)_{(p^{1/2}, p^{3/2})}\xrightarrow{\rho_{q\dR}} \Q_{p, (p^{1/2}, p^{3/2})}^\prism\xrightarrow{\widetilde{\mu}} (\ol{\DD}/\ol{\T})^\dR\rightarrow (\A^1/\G_m)^\dR
\end{equation*}
factors through the map
\begin{equation*}
\lim_{q\mapsto q^p} (1+\overcirc{\DD})\setminus\{1\}\rightarrow \A^1\;, \hspace{0.3cm} (q, q^{1/p}, q^{1/p^2}, \dots)\mapsto \log q\;.
\end{equation*}
\end{lem}

Using this result, we can already establish the following presentation of $\Q_{p, (p^{1/2}, p^{3/2})}^{\ol{\N}}$:

\begin{prop}
\label{prop:pres-qpnquotbypushout}
There is an isomorphism
\begin{equation}
\label{eq:pres-qpnquotbypushout}
\Q_{p, (p^{1/2}, p^{3/2})}^{\ol{\N}}\cong (\A^1\times(\A^1)^\dR)\times_{(\A^1)^\dR} \left(\lim_{q\mapsto q^p} (1+\overcirc{\DD})^\dR\setminus\{1\}\right)_{(p^{1/2}, p^{3/2})}\;\Bigg/\;\Z_p^{\times, \la}\coprod_{\G_m^\dagger}^\AbGrp \G_m\;,
\end{equation}
where the right-hand side is defined as follows:
\begin{enumerate}[label=(\roman*)]
\item The maps $(\A^1)^\dR\times \A^1\rightarrow (\A^1)^\dR$ and $\lim_{q\mapsto q^p} (1+\overcirc{\DD})^\dR\rightarrow (\A^1)^\dR$ defining the pullback are given by multiplication and $\log\circ\operatorname{pr}_1$, respectively.
\item $\Z_p^{\times, \la}$ acts on $\A^1$ by multiplication and on $\lim_{q\mapsto q^p} (1+\overcirc{\DD})^\dR$ via the formula
\begin{equation*}
\gamma.(q, q^{1/p}, q^{1/p^2}, \dots)=(q^\gamma, q^{\gamma/p}, q^{\gamma/p^2}, \dots)\;.
\end{equation*}
\item $\G_m$ acts on $(\A^1)^\dR$ by division and on $\A^1$ by multiplication.
\end{enumerate}
\end{prop}
\begin{proof}
We first explain how the right-hand side of (\ref{eq:pres-qpnquotbypushout}) fits into the pullback square defining $\Q_{p, (p^{1/2}, p^{3/2})}^{\ol{\N}}$. For the map towards $\Q_{p, (p^{1/2}, p^{3/2})}^\prism$, note that the second projection ``in the numerator'' and the map 
\begin{equation*}
\Z_p^{\times, \la}\coprod_{\G_m^\dagger}^\AbGrp \G_m\rightarrow \Z_p^{\times, \sm}
\end{equation*}
induce a map from the right-hand side of (\ref{eq:pres-qpnquotbypushout}) to
\begin{equation*}
\left(\lim_{q\mapsto q^p} (1+\overcirc{\DD})^\dR\setminus\{1\}\right)_{(p^{1/2}, p^{3/2})}\,\Big/\,\Z_p^{\times, \sm}\cong \Q_{p, (p^{1/2}, p^{3/2})}^\prism\;,
\end{equation*}
where the isomorphism follows from \cref{prop:defis-prismffdr} and \cref{cor:prism-hkpresentations}. Similarly, the map
\begin{equation*}
\Z_p^{\times, \la}\coprod_{\G_m^\dagger}^\AbGrp \G_m\rightarrow \G_m^\dR
\end{equation*}
and the projection onto $(\A^1)^\dR$ ``in the numerator'' induce a map from (\ref{eq:pres-qpnquotbypushout}) to $(\A^1/\G_m)^\dR$. For the map to $\A^1/\G_m$, recall that there is a canonical map $\Z_p^{\times, \la}\rightarrow \G_m$ by \cref{ex:recall-glaviadr}, whose restriction to $\G_m^\dagger\subseteq\Z_p^{\times, \la}$ identifies with the inclusion map $\G_m^\dagger\subseteq\G_m$. Therefore, we obtain a map
\begin{equation*}
\Z_p^{\times, \la}\coprod_{\G_m^\dagger}^\AbGrp \G_m\rightarrow \G_m
\end{equation*}
and this together with the projection onto $\A^1$ ``in the numerator'' induces the map from the right-hand side of (\ref{eq:pres-qpnquotbypushout}) to $\A^1/\G_m$.

We note that the maps from the right-hand side of (\ref{eq:pres-qpnquotbypushout}) to $\A^1/\G_m\times (\A^1/\G_m)^\dR$ and $\Q_{p, (p^{1/2}, p^{3/2})}^\prism$ are indeed compatible upon mapping to $(\A^1/\G_m)^\dR$: This follows from \cref{lem:pres-rhoqdrlogq} together with the fact that the map
\begin{equation*}
\left(\lim_{q\mapsto q^p} (1+\overcirc{\DD})^\dR\setminus\{1\}\right)_{(p^{1/2}, p^{3/2})}\,\Big/\,\Z_p^{\times, \sm}\cong \Q_{p, (p^{1/2}, p^{3/2})}^\prism\xrightarrow{\widetilde{\mu}} (\A^1/\G_m)^\dR
\end{equation*}
is induced by the map $(q, q^{1/p}, q^{1/p^2}, \dots)\mapsto \log q$ ``in the numerator'' and the canonical map $\Z_p^{\times, \sm}\rightarrow \G_m^\dR$ ``in the denominator''. Overall, we thus obtain a map
\begin{equation}
\label{eq:pres-qpnpushoutiso}
(\A^1\times(\A^1)^\dR)\times_{(\A^1)^\dR} \left(\lim_{q\mapsto q^p} (1+\overcirc{\DD})^\dR\setminus\{1\}\right)_{(p^{1/2}, p^{3/2})}\;\Bigg/\;\Z_p^{\times, \la}\coprod_{\G_m^\dagger}^\AbGrp \G_m\rightarrow \Q_{p, (p^{1/2}, p^{3/2})}^{\ol{\N}}\;.
\end{equation}

To check that (\ref{eq:pres-qpnpushoutiso}) is an isomorphism, we may pull back along the $\Z_p^{\times, \sm}$-torsor
\begin{equation*}
\left(\lim_{q\mapsto q^p} (1+\overcirc{\DD})^\dR\setminus\{1\}\right)_{(p^{1/2}, p^{3/2})}\rightarrow \Q_{p, (p^{1/2}, p^{3/2})}^\prism\;.
\end{equation*}
Then the source of (\ref{eq:pres-qpnpushoutiso}) becomes
\begin{equation*}
(\A^1\times(\A^1)^\dR)/\G_m\times_{(\A^1)^\dR} \left(\lim_{q\mapsto q^p} (1+\overcirc{\DD})^\dR\setminus\{1\}\right)_{(p^{1/2}, p^{3/2})}
\end{equation*}
while target sits inside a pullback square
\begin{equation*}
\begin{tikzcd}
(\,\dots)\ar[r]\ar[d] & \left(\lim_{q\mapsto q^p} (1+\overcirc{\DD})^\dR\setminus\{1\}\right)_{(p^{1/2}, p^{3/2})}\ar[d, "\log\circ\operatorname{pr}_1"] \\
\A^1/\G_m\times (\A^1/\G_m)^\dR\ar[r, "\mathrm{mult}"] & (\A^1/\G_m)^\dR\nospacepunct{\;.}
\end{tikzcd}
\end{equation*}
Now we are done upon noting that the map $\log\circ\operatorname{pr}_1$ factors through the canonical map $(\A^1)^\dR\rightarrow (\A^1/\G_m)^\dR$ and that the diagram
\begin{equation*}
\begin{tikzcd}
(\A^1\times(\A^1)^\dR)/\G_m\ar[r]\ar[d] & (\A^1)^\dR\ar[d] \\
\A^1/\G_m\times (\A^1/\G_m)^\dR\ar[r, "\mathrm{mult}"] & (\A^1/\G_m)^\dR
\end{tikzcd}
\end{equation*}
is cartesian.
\end{proof}

Finally, our goal is to make the presentation of $\Q_{p, (p^{1/2}, p^{3/2})}^{\ol{\N}}$  from \cref{prop:pres-qpnquotbypushout} even more explicit by eliminating the de Rham stacks ``in the numerator''. For this, consider the group stack $\Z_p^{\times, \la}\times (\G_a^\dagger\rtimes\G_m)$, where $\G_m$ acts on $\G_a^\dagger$ by division, and observe that, after base change to $\A^1$, it receives the following two maps from $\G_m^\dagger$: The first one is given by 
\begin{equation*}
\G_m^\dagger\rightarrow \Z_p^{\times, \la}\times (\G_a^\dagger\rtimes \G_m)\;, \hspace{0.3cm} s\mapsto (1, u(1-s^{-1}), s)
\end{equation*}
while the second one is given by postcomposing the canonical map $\G_m^\dagger\rightarrow\Z_p^{\la, \times}$ with the inclusion of the first factor. We let $\cal{G}$ denote the coequaliser of these two maps in the category of group stacks over $\A^1$. Moreover, we define the Gelfand stack $\widetilde{\Q_{p\mathrlap{, (p^{1/2}, p^{3/2})}}^{\ol{\N}}}\hphantom{\scriptstyle{(p^{1/2}, p^{3/2})}}$ by the pullback diagram
\begin{equation*}
\begin{tikzcd}
\widetilde{\Q_{p\mathrlap{, (p^{1/2}, p^{3/2})}}^{\ol{\N}}}\hphantom{\scriptstyle{(p^{1/2}, p^{3/2})}}\ar[r]\ar[d] & \left(\lim_{q\mapsto q^p} (1+\overcirc{\DD})\setminus\{1\}\right)_{(p^{1/2}, p^{3/2})}\ar[d, "\log\circ\operatorname{pr}_1"] \\
\A^1_+\times \A^1_-\ar[r, "\mathrm{mult}"] & \A^1\nospacepunct{\;.}
\end{tikzcd}
\end{equation*}
Now our main result is the following:

\begin{thm}
\label{thm:pres-qpn}
There is an isomorphism
\begin{equation*}
\Q_{p, (p^{1/2}, p^{3/2})}^{\ol{\N}}\cong \widetilde{\Q_{p\mathrlap{, (p^{1/2}, p^{3/2})}}^{\ol{\N}}}\hphantom{\scriptstyle{(p^{1/2}, p^{3/2})}}\,/\,\cal{G}\;,
\end{equation*}
where the action of $\cal{G}$ on $\widetilde{\Q_{p\mathrlap{, (p^{1/2}, p^{3/2})}}^{\ol{\N}}}\hphantom{\scriptstyle{(p^{1/2}, p^{3/2})}}$ is given by the formula
\begin{equation*}
(\gamma, w, s).(t, u, (q^{1/p^m})_m)\coloneqq\left(\gamma s t, \frac{u}{s}+w, \left(q^{\gamma/p^m}\exp\left(\frac{\gamma w s t}{p^m}\right)\right)_m\right)
\end{equation*}
for $\gamma\in\Z_p^{\times, \la}, w\in\G_a^\dagger$ and $s\in\G_m$ and where $t$ and $u$ denote the coordinates on $\A^1_+$ and $\A^1_-$, respectively.
\end{thm}
\begin{proof}
Note that the exponential makes sense due to $\gamma w s t/p^m\in\G_a^\dagger$. We first check that the action is well-defined, i.e.\ that the actions of the two copies of $\G_m^\dagger$ which sit inside $\Z_p^{\times, \la}\times (\G_a^\dagger\rtimes\G_m)$ via the maps described above agree. For this, we calculate
\begin{equation*}
(1, u(1-s^{-1}), s).(t, u, (q^{1/p^m})_m)=\left(st, u, \left(q^{1/p^m}\exp\left(\frac{u(s-1)t}{p^m}\right)\right)_m\right)=(st, u, (q^{s/p^m})_m)\;,
\end{equation*}
where the last equality is due to $ut=\log q$ on $\widetilde{\Q_{p\mathrlap{, (p^{1/2}, p^{3/2})}}^\N}\hphantom{\scriptstyle{(p^{1/2}, p^{3/2})}}$. However, the right-hand side is the image of $(t, u, (q^{1/p^m})_m)$ under the action of $(s, 0, 1)$ and hence the two actions agree.

Now observe that there is a canonical map $\widetilde{\Q_{p\mathrlap{, (p^{1/2}, p^{3/2})}}^{\ol{\N}}}\hphantom{\scriptstyle{(p^{1/2}, p^{3/2})}}\,/\,\cal{G}\rightarrow \Q_{p, (p^{1/2}, p^{3/2})}^{\ol{\N}}$ using the presentation of the target from \cref{prop:pres-qpnquotbypushout} and checking whether it is an isomorphism may be done after pullback along
\begin{equation*}
(\A^1\times(\A^1)^\dR)\times_{(\A^1)^\dR} \left(\lim_{q\mapsto q^p} (1+\overcirc{\DD})^\dR\setminus\{1\}\right)_{(p^{1/2}, p^{3/2})}\rightarrow \Q_{p, (p^{1/2}, p^{3/2})}^{\ol{\N}}\;.
\end{equation*}
Then $\widetilde{\Q_{p\mathrlap{, (p^{1/2}, p^{3/2})}}^{\ol{\N}}}\hphantom{\scriptstyle{(p^{1/2}, p^{3/2})}}\,/\,\cal{G}$ becomes $\widetilde{\Q_{p\mathrlap{, (p^{1/2}, p^{3/2})}}^{\ol{\N}}}\hphantom{\scriptstyle{(p^{1/2}, p^{3/2})}}\,/\,\G_a^\dagger$, where the action of $\G_a^\dagger$ is given by
\begin{equation*}
w.(t, u, (q^{1/p^m})_m)=\left(t, u+w, \left(q^{1/p^m}\exp\left(\frac{wt}{p^m}\right)\right)_m\right)\;.
\end{equation*}
However, now we are done due to
\begin{equation*}
\widetilde{\Q_{p\mathrlap{, (p^{1/2}, p^{3/2})}}^{\ol{\N}}}\hphantom{\scriptstyle{(p^{1/2}, p^{3/2})}}\,/\,\G_a^\dagger\cong (\A^1\times(\A^1)^\dR)\times_{(\A^1)^\dR} \left(\lim_{q\mapsto q^p} (1+\overcirc{\DD})^\dR\setminus\{1\}\right)_{(p^{1/2}, p^{3/2})}\;.
\end{equation*}
Indeed, to see this last isomorphism, observe that, given an $A$-point $(t, u, (q^{1/p^m})_m)$ of the target for some nilperfectoid $A$, once we lift $u$ from $\ol{A}$ to $A$, the condition $ut=\log q$ in $A$ specifies a unique lift of $q$ from $\ol{A}$ to $A$ and then the fact that the map $q\mapsto q^p$ on $1+\overcirc{\DD}$ is $\dagger$-formally étale also yields unique lifts of the elements $q^{1/p}, q^{1/p^2}, \dots$. As the ambiguity of lifting $u$ from $\ol{A}$ to $A$ is precisely killed by modding out the $\G_a^\dagger$-action, this finishes the proof.
\end{proof}

\subsection{Nygaardifications of smooth maps between derived Berkovich spaces}

We now move on to studying $X^\N\rightarrow Y^\N$ for any Berkovich smooth map $X\rightarrow Y$ of derived Berkovich spaces. For this, we first observe that $X^\N$ descends along $\Q_p^\N\rightarrow \Q_p^{\ol{\N}}$ for any Gelfand stack $X$. Indeed, note that, given an $A$-point of $\Q_p^\N$ for $A$ totally disconnected nilperfectoid, the definition of the Gelfand stack $D^+$ from (\ref{eq:syn-defid+}) only makes use of the data of the $A$-point of $\Q_p^{\ol{\N}}$ obtained from the given $A$-point of $\Q_p^\N$. Thus, we can define $D^+$ by the diagram (\ref{eq:syn-defid+}) for any $A$-point of $\Q_p^{\ol{\N}}$ with $A$ totally disconnected nilperfectoid and then define, as before:

\begin{defi}
For any Gelfand stack $X$, we define its \emph{reduced Nygaardification} as the Gelfand stack $X^{\ol{\N}}$ over $\Q_p^{\ol{\N}}$ obtained by sheafifying the assignment
\begin{equation*}
X^{\ol{\N}}(\GSpec A\rightarrow \Q_p^{\ol{\N}})\coloneqq \{\text{maps}\,D^+\rightarrow X\}
\end{equation*}
for any totally disconnected nilperfectoid $A$.
\end{defi}

By definition, we have $X^\N\cong X^{\ol{\N}}\times_{\Q_p^{\ol{\N}}} \Q_p^\N$ for any Gelfand stack $X$ and hence studying the map $X^\N\rightarrow Y^\N$ for $X\rightarrow Y$ a Berkovich smooth map of derived Berkovich spaces reduces to studying the map $X^{\ol{\N}}\rightarrow Y^{\ol{\N}}$. Moreover, we point out that \cref{thm:defis-berketalecover} also holds for $X^{\ol{\N}}$ in place of $X^\N$ without any changes to the proof. By definition of Berkovich smooth maps, this further reduces us to understanding $(\G_m^n)^{\ol{\N}}\rightarrow \Q_p^{\ol{\N}}$ and, by compatibility of the functor $X\mapsto X^{\ol{\N}}$ with limits, we can even reduce to $n=1$. 

In this case, we are going to describe $\G_{m,(p^{1/2}, p^{3/2})}^{\ol{\N}}$ as an explicit quotient using the coordinates on the cover 
\begin{equation*}
(\A^1_+\times \A^1_-)\times_{\A^1} \left(\lim_{q\mapsto q^p} (1+\overcirc{\DD})\setminus\{1\}\right)_{(p^{1/2}, p^{3/2})}=\widetilde{\Q_{p\mathrlap{, (p^{1/2}, p^{3/2})}}^{\ol{\N}}}\hphantom{\scriptstyle{(p^{1/2}, p^{3/2})}}\rightarrow \Q_{p, (p^{1/2}, p^{3/2})}^{\ol{\N}}
\end{equation*}
from \cref{thm:pres-qpn}. Recalling that we denote the coordinates on $\A^1_+$ and $\A^1_-$ by $t$ and $u$, respectively, we can form the pushout $\Z_p^\la\coprod^\AbGrp_{\G_a^\dagger} \G_a^\dagger$ in the category of abelian group stacks over $\Q_{p, (p^{1/2}, p^{3/2})}^{\ol{\N}}$, where the map $\G_a^\dagger\rightarrow \G_a^\dagger$ is given by multiplication by $u$.

\begin{thm}
\label{thm:pres-tn}
Over $\Q_{p, (p^{1/2}, p^{3/2})}^{\ol{\N}}$, there is an isomorphism
\begin{equation}
\label{eq:pres-tn}
\G_m^{\ol{\N}}\cong \lim_{x\mapsto x^p} \G_m\,\Big/\,\Z_p^\la\coprod^\AbGrp_{\G_a^\dagger} \G_a^\dagger
\end{equation}
of normed multiplicative group stacks, where the action defining the right-hand side is given by the formula
\begin{equation*}
(\theta, v).x^{1/p^m}\coloneqq q^{\theta/p^m}\exp\left(\frac{tv}{p^m}\right)x^{1/p^m}\;,
\end{equation*}
where $(q^{1/p^m})_m$ are the coordinates on $\lim_{q\mapsto q^p} (1+\overcirc{\DD})$.
\end{thm}
\begin{proof}
Note that the exponential makes sense due to $tv/p^m\in\G_a^\dagger$. We first check that the action is well-defined, i.e.\ that the actions of the two copies of $\G_a^\dagger$ sitting inside $\Z_p^\la\times \G_a^\dagger$ via postcomposing the canonical map $\G_a^\dagger\rightarrow\Z_p^\la$ with the inclusion of the first factor and via postcomposing the multiplication-by-$u$-map on $\G_a^\dagger$ with the inclusion of the second factor, respectively, agree. For this, we calculate
\begin{equation*}
(0, u\theta).x^{1/p^m}=\exp\left(\frac{tu\theta}{p^m}\right)x^{1/p^m}=q^{\theta/p^m}x^{1/p^m}
\end{equation*}
for $\theta\in\G_a^\dagger$, where the last equality is due to $tu=\log q$ on $\Q_{p, (p^{1/2}, p^{3/2})}^{\ol{\N}}$. As the right-hand side is also the image of $x^{1/p^m}$ under the action of $(\theta, 0)$, the two actions indeed agree and the formula from the theorem statement defines an action of the pushout $\Z_p^\la\coprod^\AbGrp_{\G_a^\dagger} \G_a^\dagger$.

To move on, note that $\G_m^{\ol{\N}}$ sits inside a pullback square
\begin{equation}
\label{eq:pres-tnpullback}
\begin{tikzcd}
\G_m^{\ol{\N}}\ar[r]\ar[d] & \G_m^\prism\ar[d] \\
\G_m^{\Cone}\ar[r] & (\G_m^{\Cone})^\dR
\end{tikzcd}
\end{equation}
as in \cref{lem:syn-xnviacone}, where we note that the cone stack descends along $\ol{\DD}/\ol{\T}\rightarrow\A^1/\G_m$ and we use this descended version in the diagram above as well as the rest of the proof. We explain how the right-hand side of (\ref{eq:pres-tn}) fits into the same pullback diagram over $\Q_{p, (p^{1/2}, p^{3/2})}^{\ol{\N}}$. For this, first note that $\G_m^\prism\cong (\G_m^\prism)^\dR$ as we are working over $(p^{1/2}, p^{3/2})$, see \cref{prop:defis-prismffdr}, and hence \cref{cor:prism-hkpresentations} yields the presentation
\begin{equation}
\label{eq:pres-tprismpres}
\G_m^\prism\cong \lim_{x\mapsto x^p} \G_m^\dR\,/\,\Z_p^\sm
\end{equation}
over $\Q_{p, (p^{1/2}, p^{3/2})}^{\ol{\N}}$. Now note that the right-hand side of (\ref{eq:pres-tn}) maps to the right-hand side of (\ref{eq:pres-tprismpres}) via the projection onto the de Rham stack ``in the numerator'' and the map
\begin{equation*}
\Z_p^\la\coprod^\AbGrp_{\G_a^\dagger} \G_a^\dagger\rightarrow \Z_p^\sm
\end{equation*}
``in the denominator''.

To obtain the map to $\G_m^{\Cone}$, we recall from \cref{lem:pres-rhoqdrlogq} that the fact that we are working over $(p^{1/2}, p^{3/2})$ is equivalent to the inequality
\begin{equation*}
p^{-1/(p^{1/2}(p-1))}<|q^{1/p^3}-1|<p^{-1/(p^{3/2}(p-1))}\;,
\end{equation*}
from which we in particular deduce that $|q-1|<p^{-1}$. As the map $\log: 1+\ol{\DD}(p^{-1})\rightarrow \ol{\DD}(p^{-1})$ is an isomorphism with inverse given by the exponential, we obtain $q=\exp(tu)$ from the requirement $tu=\log q$ on $\Q_{p, (p^{1/2}, p^{3/2})}^{\ol{\N}}$. Thus, we can rewrite the action of $(\theta, v)\in\Z_p^\la\coprod^\AbGrp_{\G_a^\dagger} \G_a^\dagger$ on $x$ as
\begin{equation*}
(\theta, v).x=\exp(t(u\theta+v))x\;,
\end{equation*}
which makes sense due to $tu\theta\in\ol{\DD}(p^{-1})$, as we just established. Expanding the exponential on the right-hand side into its power series, we see that $(\theta, v).x=x+t(\,\dots)$, where the term in brackets is an explicit power series in $t, u, \theta, v$ and $x$. As the generalised Cartier divisor defining the cone stack over $\Q_{p, (p^{1/2}, p^{3/2})}^{\ol{\N}}$ is given by multiplication by $t$, this shows that the coordinate $x$ on $\G_m$ defines a map from the right-hand side of (\ref{eq:pres-tn}) to $\G_m^{\Cone}$ over $\Q_{p, (p^{1/2}, p^{3/2})}^{\ol{\N}}$.

As the map 
\begin{equation*}
\G_m^\prism\cong \lim_{x\mapsto x^p} \G_m^\dR\,/\,\Z_p^\sm\rightarrow (\G_m^{\Cone})^\dR
\end{equation*}
over $\Q_{p, (p^{1/2}, p^{3/2})}^{\ol{\N}}$ is induced by the coordinate on $\G_m$ as well, the maps from the right-hand side of (\ref{eq:pres-tn}) to $\G_m^\prism$ and $\G_m^{\Cone}$ we have constructed are compatible upon mapping to $(\G_m^{\Cone})^\dR$ and hence induce a map
\begin{equation*}
\lim_{x\mapsto x^p} \G_m\,\Big/\,\Z_p^\la\coprod^\AbGrp_{\G_a^\dagger} \G_a^\dagger\rightarrow \G_m^{\ol{\N}}
\end{equation*}
over $\Q_{p, (p^{1/2}, p^{3/2})}^{\ol{\N}}$.

To see that this is an isomorphism, first note that, since $\G_m$ is $\dagger$-formally smooth, there is a short exact sequence of group stacks
\begin{equation*}
\begin{tikzcd}
0\ar[r] & \G_m^\dagger/\G_a^\dagger\ar[r] & \G_m^{\Cone}\ar[r] & (\G_m^{\Cone})^\dR\ar[r] & 0\nospacepunct{\;,}
\end{tikzcd}
\end{equation*}
where $\G_a^\dagger$ acts on $\G_m^\dagger$ via $v.x\coloneqq \exp(tv)x$, and hence (\ref{eq:pres-tnpullback}) shows that $\G_m^{\ol{\N}}$ is a $\G_m^\dagger/\G_a^\dagger$-torsor over $\G_m^\prism$. However, note that also the map
\begin{equation*}
\lim_{x\mapsto x^p} \G_m\,\Big/\,\Z_p^\la\coprod^\AbGrp_{\G_a^\dagger} \G_a^\dagger\rightarrow \lim_{x\mapsto x^p} \G_m^\dR\,/\,\Z_p^\sm
\end{equation*}
is a $\G_m^\dagger/\G_a^\dagger$-torsor: Indeed, we have
\begin{equation*}
\lim_{x\mapsto x^p} \G_m^\dR\cong \lim_{x\mapsto x^p} \G_m\,\Big/\,\G_m^\dagger
\end{equation*}
as the map $x\mapsto x^p$ on $\G_m$ is $\dagger$-formally smooth. Thus, we are done upon noting that the map
\begin{equation*}
\lim_{x\mapsto x^p} \G_m\,\Big/\,\Z_p^\la\coprod^\AbGrp_{\G_a^\dagger} \G_a^\dagger\rightarrow \G_m^{\ol{\N}}
\end{equation*}
we have constructed is a map of $\G_m^\dagger/\G_a^\dagger$-torsors over $\G_m^\prism$ via the isomorphism (\ref{eq:pres-tprismpres}).
\end{proof}

\begin{cor}
\label{prop:perf-coversmoothrigid}
Let $X$ be a derived Berkovich space admitting a map $X\rightarrow \ol{\T}^n$ which is the composition of a finite étale map $X\rightarrow Z$ and a rational localisation $Z\rightarrow\ol{\T}^n$. Then the Gelfand stack $X_{[r, s]}^\N$ is nicely coverable for any $[r, s]\subseteq (p^{1/2}, p^{3/2})$.
\end{cor}
\begin{proof}
By compatibility of $X\mapsto X^{\ol{\N}}$ with finite étale maps and rational localisations, which one proves as in \cref{thm:defis-berketalecover}, we may reduce to the case $X=\ol{\T}^n$ using \cref{lem:perf-etclosedcoverable}. By commutation with limits, we then conclude from \cref{thm:pres-tn} that there is an isomorphism
\begin{equation*}
(\ol{\T}^n)^{\ol{\N}}\cong \lim_{x_i\mapsto x_i^p} \ol{\T}^n\,\Big/\,(\Z_p^n)^\la\coprod_{(\G_a^\dagger)^n}^\AbGrp (\G_a^\dagger)^n
\end{equation*}
over $\Q_{p, (p^{1/2}, p^{3/2})}^{\ol{\N}}$. Now note that the right-hand side has a \v{C}ech cover given by
\begin{equation*}
\lim_{x_i\mapsto x_i^p} \ol{\T}^n\,\times\, \left((\Z_p^n)^\la\coprod_{(\G_a^\dagger)^n}^\AbGrp (\G_a^\dagger)^n\right)^\bullet
\end{equation*}
and that the stack $\lim_{x_i\mapsto x_i^p} \ol{\T}^n$ is affine, static and flat over $\Q_p$. Thus, we will be able to conclude that $(\ol{\T}^n)^{\ol{\N}}_{[r, s]}$ is nicely coverable ``relative to $\Q_{p, (p^{1/2}, p^{3/2})}^{\ol{\N}}$'', i.e.\ that every base change to an affine truncated Gelfand stack over $\Q_{p, (p^{1/2}, p^{3/2})}^{\ol{\N}}$ is nicely coverable, once we know that $(\Z_p^n)^\la\coprod_{(\G_a^\dagger)^n}^\AbGrp (\G_a^\dagger)^n$ has a $!$-hypercover by Gelfand stacks which are affine, static and flat over $\Q_p$. 

However, this simply follows from the fact that
\begin{equation*}
(\Z_p^n)^\la\coprod_{(\G_a^\dagger)^n}^\AbGrp (\G_a^\dagger)^n\cong (\Z_p^n)^\la\times (\G_a^\dagger)^n\,\Big/\,(\G_a^\dagger)^n
\end{equation*}
as the right-hand side has a $!$-cover by $(\Z_p^n)^\la\times (\G_a^\dagger)^n$ whose \v{C}ech nerve is given by 
\begin{equation*}
(\Z_p^n)^\la\times (\G_a^\dagger)^n\times ((\G_a^\dagger)^n)^\bullet\;,
\end{equation*}
which is a nice hypercover as $\O(\Z_p^\la)=C^\la(\Z_p, \Q_p)$ and $\O(\G_a^\dagger)=\Q_p\{v\}^\dagger$ are both flat over $\Q_p$ by \cref{lem:perf-vspflat}.

Finally, as $(\ol{\T}^n)^\N\cong (\ol{\T}^n)^{\ol{\N}}\times_{\Q_p^{\ol{\N}}} \Q_p^\N$, the claim follows from \cref{lem:perf-covers} by base changing the nice hypercover of $(\ol{\T}^n)^{\ol{\N}}$ ``relative to $\Q_{p, (p^{1/2}, p^{3/2})}^{\ol{\N}}$'' from above to the nice hypercover of $\Q_{p, [r, s]}^\N$ from the proof of \cref{prop:perf-zpncover}, where we recall that all constituents of this latter hypercover are static, which yields the uniformity required in \cref{lem:perf-covers}.
\end{proof}

\subsection{The Hodge--Tate to de Rham degeneration}
\label{subsect:htdr}

We now turn our attention to the stack $X^{\HT, \dagger, +}=X^\N_{|t|=0}$, which we want to describe completely explicitly in the case where $X$ is equipped with a Berkovich étale map $X\rightarrow\ol{\T}^n$ for some $n\geq 0$. To set the stage, we begin with the case $X=\GSpec\Q_p$, which is already interesting. For this, consider the group stack $\Z_p^{\times, \la}\times (\G_a^\dagger\rtimes\ol{\T})$, where $\ol{\T}$ acts on $\G_a^\dagger$ by division. Similar to what happened above, this group stack receives two maps from $\G_m^\dagger$ after base change to $\ol{\DD}=\GSpec\Q_p\langle u\rangle_{\leq 1}$: The first one is given by 
\begin{equation*}
\G_m^\dagger\rightarrow \Z_p^{\times, \la}\times (\G_a^\dagger\rtimes \ol{\T})\;, \hspace{0.3cm} s\mapsto (1, u(1-s^{-1}), s)
\end{equation*}
while the second one is given by postcomposing the canonical map $\G_m^\dagger\rightarrow\Z_p^{\la, \times}$ with the inclusion of the first factor. We let $\cal{G}_0$ denote the coequaliser of these two maps in the category of group stacks over $\ol{\DD}$. Using this, we can describe $\Q_p^{\HT, \dagger, +}$ explicitly as follows:

\begin{prop}
\label{prop:htdr-qphtdagger+}
There is an isomorphism
\begin{equation}
\label{eq:htdr-qphtdagger+}
\Q_p^{\HT, \dagger, +}\cong \G_a^\dagger\times\ol{\DD}\times\GSpec\Q_p(\zeta_{p^\infty})\,\Big/\,\cal{G}_0\;,
\end{equation}
where the action of $\cal{G}_0$ on $\G_a^\dagger\times\ol{\DD}\times\GSpec\Q_p(\zeta_{p^\infty})$ is given by the formula
\begin{equation*}
(\gamma, w, s).(t, u, (\zeta_{p^m})_m)\coloneqq \left(\gamma s t, \frac{u}{s}+w, (\zeta_{p^m}^\gamma)_m\right)
\end{equation*}
for $\gamma\in\Z_p^{\times, \la}, w\in\G_a^\dagger$ and $s\in\ol{\T}$ and where $t$ and $u$ denote the coordinates on $\G_a^\dagger$ and $\ol{\DD}$, respectively. In particular, restricting to $|u|=1$, we have
\begin{equation*}
\Q_p^{\HT, \dagger}\cong \G_a^\dagger\times\GSpec\Q_p(\zeta_{p^\infty})\,\Big/\,\Z_p^{\times, \la}\;.
\end{equation*}
\end{prop}
\begin{proof}
From \cref{thm:pres-qpn}, we already know that
\begin{equation*}
\Q_{p, |t|=0}^{\ol{\N}}\cong (\G_a^\dagger\times \A^1)\times_{\A^1} \left(\lim_{q\mapsto q^p} (1+\overcirc{\DD})\setminus\{1\}\right)_{(p^{1/2}, p^{3/2})}\;\Big/\;\cal{G}
\end{equation*}
and, in particular, $\log q=tu$ is $\dagger$-nilpotent meaning that $q\in (\mu_{p^\infty}\subseteq 1+\overcirc{\DD})^\dagger$, where $\mu_{p^\infty}$ is the group stack of $p$-power roots of unity and $(\mu_{p^\infty}\subseteq 1+\overcirc{\DD})^\dagger$ denotes its overconvergent neighbourhood inside $1+\overcirc{\DD}$. As the fact that we are over $(p^{1/2}, p^{3/2})\subseteq (0, \infty)$ corresponds to the inequality
\begin{equation*}
p^{-1/(p^{1/2}(p-1))}<|q^{1/p^3}-1|<p^{-1/(p^{3/2}(p-1))}
\end{equation*}
by \cref{lem:pres-rhoqdrlogq}, we conclude that $q^{1/p^3}$ lies in the overconvergent neighbourhood of a primitive $p^2$-th root of unity and, consequently, $q^{1/p^m}$ lies in the overconvergent neighbourhood of a primitive $p^{m-1}$-th root of unity for all $m\geq 1$ while $q\in\G_m^\dagger$. As the map $q\mapsto q^p$ is $\dagger$-formally étale on $1+\overcirc{\DD}$, we conclude that $(q^{1/p^m})_m$ is uniquely determined by $q=\exp(tu)$ and hence
\begin{equation*}
(\G_a^\dagger\times \A^1)\times_{\A^1} \left(\lim_{q\mapsto q^p} (1+\overcirc{\DD})\setminus\{1\}\right)_{(p^{1/2}, p^{3/2})}\cong \G_a^\dagger\times\A^1\times\GSpec\Q_p(\zeta_{p^\infty})\;.
\end{equation*}

Next, note that the map
\begin{equation}
\label{eq:htdr-qphtdagger+tumap}
\G_a^\dagger\times\A^1\times\GSpec\Q_p(\zeta_{p^\infty})\cong \Q_{p, |t|=0}^{\ol{\N}}\rightarrow (\G_a^\dagger/\ol{\T}\times (\A^1_-)^\dR)/\G_m^\dR
\end{equation}
induced by the canonical maps $\G_a^\dagger\times\A^1\rightarrow \G_a^\dagger\times (\A^1_-)^\dR, \Z_p^{\times, \la}\rightarrow\ol{\T}$, see \cref{ex:recall-glaviadr}, and $\cal{G}/\Z_p^{\times, \la}\rightarrow\G_m^\dR$ agrees with the restriction of the map
\begin{equation*}
\Q_p^{\ol{\N}}\rightarrow (\A^1_+/\ol{\T}\times (\A^1_-)^\dR)/\G_m^\dR
\end{equation*}
from the proof of \cref{lem:pres-qpredntorsor}. Indeed, by loc.\ cit., this comes down to checking that the postcomposition of the map (\ref{eq:htdr-qphtdagger+tumap}) with $(\G_a^\dagger/\ol{\T}\times (\A^1_-)^\dR)/\G_m^\dR\rightarrow */\ol{\T}^\dR$ classifies the Tate twist under the equivalence from \cref{lem:syn-btdr}, but this is clear as the map $\GSpec\Q_p\rightarrow */\ol{\T}^\dR$ classifying the Tate twist factors as
\begin{equation*}
\GSpec\Q_p\rightarrow */\Z_p^{\times, \sm}\rightarrow */\ol{\T}^\dR\;,
\end{equation*}
where the first map classifies the $\Z_p^{\times, \sm}$-torsor $\GSpec\Q_p(\zeta_{p^\infty})\rightarrow\GSpec\Q_p$.

Now recall that the proof of \cref{lem:pres-qpredntorsor} yields a cartesian diagram
\begin{equation*}
\begin{tikzcd}
\Q_{p, |t|=0}^{\N, \mathrm{ext}}\ar[r]\ar[d] & \Q_{p, |t|=0}^{\ol{\N}}\ar[d] \\
\A^1_+/\ol{\T}\times (\A^1_-/\ol{\T})^\dR\ar[r] & (\A^1_+/\ol{\T}\times (\A^1_-)^\dR)/\G_m^\dR
\end{tikzcd}
\end{equation*}
and in particular shows that $\Q_{p, |t|=0}^{\N, \mathrm{ext}}\rightarrow\Q_{p, |t|=0}^{\ol{\N}}$ is a torsor for $\G_m/\ol{\T}\cong (0, \infty)$, where the isomorphism is due to \cref{lem:defis-rhvariant}. However, evidently, by the presentation of $\Q_{p, |t|=0}^{\ol{\N}}$ from (\ref{eq:htdr-qphtdagger+tumap}), the stack
\begin{equation*}
\G_a^\dagger\times\A^1\times\GSpec\Q_p(\zeta_{p^\infty})\,\Big/\,\cal{G}_0
\end{equation*}
fits into the same cartesian diagram and is a $(0, \infty)$-torsor over $\Q_{p, |t|=0}^{\ol{\N}}$ as well due to $\cal{G}/\cal{G}_0\cong \G_m/\ol{\T}\cong (0, \infty)$ and hence we conclude that the induced map
\begin{equation*}
\G_a^\dagger\times\A^1\times\GSpec\Q_p(\zeta_{p^\infty})\,\Big/\,\cal{G}_0\rightarrow \Q_{p, |t|=0}^{\N, \mathrm{ext}}
\end{equation*}
is an isomorphism upon noting that it is a map of $(0, \infty)$-torsors. Now the claimed isomorphism (\ref{eq:htdr-qphtdagger+}) follows from $\Q_p^\N\cong \Q_{p, |u|\leq 1}^\N$, see (\ref{eq:pres-qpnfromqpnext}).
\comment{
As $\Q_{p, |\widetilde{\mu}|=0}^\prism\cong \Q_p^\dR\cong\GSpec\Q_p$ by the proof of \cref{prop:defis-utneq0}, the definition of $\Q_p^\N$ shows that the stack $\Q_p^{\HT, \dagger, +}$ sits inside a pullback diagram
\begin{equation*}
\begin{tikzcd}
\Q_p^{\HT, \dagger, +}\ar[d]\ar[r] & \GSpec\Q_p\ar[d] \\
\G_a^\dagger/\ol{\T}\times (\ol{\DD}/\ol{\T})^\dR\ar[r, "\mathrm{mult}"] & */\ol{\T}^\dR\nospacepunct{\;,}
\end{tikzcd}
\end{equation*}
where the vertical map on the right classifies the Tate twist by the discussion following \cref{lem:syn-btdr}. Factoring the bottom map as
\begin{equation*}
\G_a^\dagger/\ol{\T}\times (\ol{\DD}/\ol{\T})^\dR\rightarrow */\ol{\T}^\dR\times (\ol{\DD}/\ol{\T})^\dR\xrightarrow{\mathrm{mult}} */\ol{\T}^\dR
\end{equation*}
shows that we can simplify this to a pullback diagram
\begin{equation}
\label{eq:htdr-qphtdagger+pullback}
\begin{tikzcd}
\Q_p^{\HT, \dagger, +}\ar[r]\ar[d] & (\ol{\DD}/\ol{\T})^\dR\ar[d] \\
\G_a^\dagger/\ol{\T}\ar[r] & */\ol{\T}^\dR\;,
\end{tikzcd}
\end{equation}
where the map on the right-hand side is the first Tate twist of the canonical map.

We now explain how the right-hand side of (\ref{eq:htdr-qphtdagger+}) fits into the same pullback diagram. For this, first note that the right-hand side of (\ref{eq:htdr-qphtdagger+}) maps to $(\ol{\DD}/\ol{\T})^\dR$ by projecting onto $\ol{\DD}$ and then postcomposing with the canonical map $\ol{\DD}\rightarrow\ol{\DD}^\dR$ ``in the numerator'' and using the canonical map $\cal{G}_0\rightarrow \ol{\T}^\dR$ ``in the denominator''. For the map to $\G_a^\dagger/\ol{\T}$, we note that the map $\Z_p^{\times, \la}\times (\G_a^\dagger\rtimes\ol{\T})\rightarrow \ol{\T}$ given by $(\gamma, w, s)\mapsto \gamma s$ descends to a map $\cal{G}_0\rightarrow\ol{\T}$. Together with the projection onto $\G_a^\dagger$ ``in the numerator'', this induces a map from the right-hand side of (\ref{eq:htdr-qphtdagger+}) to $\G_a^\dagger/\ol{\T}$. Finally, to see that the maps to $(\ol{\DD}/\ol{\T})^\dR$ and $\G_a^\dagger/\ol{\T}$ we have given are compatible upon mapping to $*/\ol{\T}^\dR$, one just has to observe that the map $\GSpec\Q_p\rightarrow */\ol{\T}^\dR$ classifying the Tate twist factors as
\begin{equation*}
\GSpec\Q_p\rightarrow */\Z_p^{\times, \sm}\rightarrow */\ol{\T}^\dR\;,
\end{equation*}
where the first map classifies the $\Z_p^{\times, \sm}$-torsor $\GSpec\Q_p(\zeta_{p^\infty})\rightarrow\GSpec\Q_p$.

Thus, from (\ref{eq:htdr-qphtdagger+pullback}), we obtain a map
\begin{equation}
\label{eq:htdr-qphtdagger+map}
\G_a^\dagger\times\ol{\DD}\times\GSpec\Q_p(\zeta_{p^\infty})\,\Big/\,\cal{G}_0\rightarrow \Q_p^{\HT, \dagger, +}\;.
\end{equation}
To see that it is an isomorphism, observe that (\ref{eq:htdr-qphtdagger+pullback}) shows that $\Q_p^{\HT, \dagger, +}$ is a $\G_a^\dagger/\G_m^\dagger$-torsor over $(\ol{\DD}/\ol{\T})^\dR$. However, the map
\begin{equation*}
\G_a^\dagger\times\ol{\DD}\times\GSpec\Q_p(\zeta_{p^\infty})\,\Big/\,\cal{G}_0\rightarrow (\ol{\DD}/\ol{\T})^\dR
\end{equation*}
we have constructed also makes the source a $\G_a^\dagger/\G_m^\dagger$-torsor over the target: Indeed, the map
\begin{equation*}
\G_a^\dagger\times\ol{\DD}\times\GSpec\Q_p(\zeta_{p^\infty})\,\Big/\,\cal{G}_0\rightarrow\ol{\DD}\times\GSpec\Q_p(\zeta_{p^\infty})\,\Big/\,\Z_p^{\times, \sm}\times (\G_a^\dagger\rtimes\ol{\T})/\G_m^\dagger\cong \ol{\DD}\,\Big/\,(\G_a^\dagger\rtimes\ol{\T})/\G_m^\dagger
\end{equation*}
is a $\G_a^\dagger/\G_m^\dagger$-torsor as well, where the map $\G_m^\dagger\rightarrow \G_a^\dagger\rtimes\ol{\T}$ is given by $s\mapsto (u(1-s^{-1}), s)$ and we have used the isomorphism $\GSpec\Q_p(\zeta_{p^\infty})\,/\,\Z_p^{\times, \sm}\cong\GSpec\Q_p$. Now the claim follows by noting that
\begin{equation*}
\ol{\DD}\,\Big/\,(\G_a^\dagger\rtimes\ol{\T})/\G_m^\dagger\cong (\ol{\DD}/\ol{\T})^\dR
\end{equation*}
via the canonical map: This may be checked after pullback along $\ol{\DD}^\dR\rightarrow (\ol{\DD}/\ol{\T})^\dR$, where it amounts to the claim that $\ol{\DD}/\G_a^\dagger\cong\ol{\DD}^\dR$ via the canonical map. Overall, we thus conclude that (\ref{eq:htdr-qphtdagger+map}) is an isomorphism upon noting that it is a map of $\G_a^\dagger/\G_m^\dagger$-torsors over $(\ol{\DD}/\ol{\T})^\dR$.
}

Finally, let us comment on how to obtain the claimed presentation of $\Q_p^{\HT, \dagger}$. For this, observe that, if $|u|=1$, we have $\Z_p^{\times, \la}\times\ol{\T}\cong\cal{G}_0$ via the natural map from left to right and hence
\begin{equation*}
\Q_p^{\HT, \dagger}\cong \G_a^\dagger\times\ol{\T}\times\GSpec\Q_p(\zeta_{p^\infty})\,\Big/\,\Z_p^{\times, \la}\times\ol{\T}\;.
\end{equation*}
Twisting the ``numerator'' of the right-hand side by the automorphism of $\G_a^\dagger\times\ol{\T}$ given by $(t, u)\mapsto (tu, u)$, the action of $\ol{\T}$ on $\G_a^\dagger$ becomes trivial and the claim follows.
\end{proof}

\begin{cor}
\label{prop:htdr-complexesqphtdagger+}
A perfect complex on $\Q_p^{\HT, \dagger, +}$ amounts to the following data:
\begin{enumerate}[label=(\roman*)]
\item A diagram
\begin{equation*}
\begin{tikzcd}
\dots \ar[r,shift left=.5ex,"u"]
  & \ar[l,shift left=.5ex, "t"] \Fil_{i-1} V \ar[r,shift left=.5ex,"u"] & \ar[l,shift left=.5ex, "t"] \Fil_i V \ar[r,shift left=.5ex,"u"] & \ar[l,shift left=.5ex, "t"] \Fil_{i+1} V\ar[r,shift left=.5ex,"u"] & \ar[l,shift left=.5ex, "t"] \dots
\end{tikzcd}
\end{equation*}
of perfect complexes on $\GSpec\Q_p(\zeta_{p^\infty})\times \GSpec\Q_p\{r\}^\dagger$, where $r=ut$, with the property that $t: \Fil_\bullet V\rightarrow \Fil_{\bullet-1} V$ becomes an isomorphism for $\bullet\ll 0$ and $u: \Fil_\bullet V\rightarrow\Fil_{\bullet+1} V$ becomes an isomorphism for $\bullet\gg 0$;
\item $\Q_p(\zeta_{p^\infty})$-linear maps $D: \Fil_i V\rightarrow \Fil_{i-1} V$ for each $i$ satisfying $Du=uD+1$ and commuting with $t$;
\item a semilinear locally analytic $\Z_p^\times$-action on $\Fil_\bullet V$ commuting with $u$ and $D$, where the action of $\Z_p^\times$ on $\Q_p(\zeta_{p^\infty})$ is the usual Galois action and the one on $t$ is given by multiplication;
\item for each $i$, an identification between the ``Sen operator'' $\Theta: \Fil_i V\rightarrow \Fil_i V$ coming from the Lie algebra action of $\Z_p^\times$ and the operator $uD-i$.
\end{enumerate}
Under this identification, restriction to $\Q_p^{\HT, \dagger}$ corresponds to forgetting the $u$-filtration and the $D$-maps while restriction to $\Q_{p, |u|=|t|=0}^\N$ corresponds to passing to the associated graded of the $u$-filtration and forgetting the $D$-maps. Finally, the isomorphism
\begin{equation*}
\Q_{p, |u|=|t|=0}^\N\cong \Q_{p, |t|=0}^{\dR, +}\cong \G_a^\dagger/\ol{\T}
\end{equation*}
from \cref{prop:defis-u0} Tate twists the $\Z_p^\times$-action on the $i$-th associated graded piece of the $u$-filtration by $i$, which trivialises the Lie algebra action and hence after descent along $\Q_p\rightarrow \Q_p(\zeta_{p^\infty})$ via the now smooth $\Z_p^\times$-action yields a filtered perfect complex of $\Q_p$-vector spaces via the $t$-maps.
\end{cor}
\begin{proof}
Item (i) is a version of the Rees equivalence in the form from \cref{cor:recall-reesdt} while the operator $D$ from (ii) comes from the $\G_a^\dagger$-factor in $\Z_p^{\times, \la}\times (\G_a^\dagger\rtimes\ol{\TT})$ via Cartier duality, see \cref{lem:recall-cartierperf}; moreover, the relation $Du=uD+1$ comes from the fact that the $\G_a^\dagger$-action on $\ol{\DD}$ corresponds to the endomorphism $\tfrac{\mathrm{d}}{\mathrm{d}u}$ of $\Q_p\langle u\rangle_{\leq 1}$ via Cartier duality and the fact that $D$ lowers the filtration degree by $1$ is due to $\ol{\T}$ acting on $\G_a^\dagger$ by division. Furthermore, (iii) is clear and the Tate twist occurring in translating between perfect complexes on $\Q_{p, |u|=|t|=0}^\N$ and $\Q_{p, |t|=0}^{\dR, +}$ is due to the fact that the isomorphism $\Q_{p, |u|=|t|=0}^\N\cong \Q_{p, |t|=0}^{\dR, +}$ itself involves a Tate twist.

The only part that requires a little more explanation is (iv). Here, we note that the Sen operator $\Theta$ is obtained by restricting the action of $\Z_p^{\times, \la}$ along the canonical map $\G_m^\dagger\rightarrow\Z_p^{\times, \la}$ and using the isomorphism $\G_m^\dagger\cong \G_a^\dagger$ via the logarithm together with Cartier duality for $\G_a^\dagger$ in the form from \cref{lem:recall-cartierperf}. Thus, we have to check that the endomorphism $\Theta$ obtained by restricting the action of $\G_a^\dagger\rtimes\ol{\T}$ along
\begin{equation*}
\G_m^\dagger\rightarrow \G_a^\dagger\rtimes\ol{\T}\;, \hspace{0.3cm} s\mapsto (u(1-s^{-1}), s)\;,
\end{equation*}
using $\G_m^\dagger\cong \G_a^\dagger$ via the logarithm and Cartier duality for $\G_a^\dagger$ identifies with $uD-i$ on $\Fil_i V$.

To see this, first note that pulling back the action of $\ol{\T}$ along $\G_m^\dagger\rightarrow\ol{\T}$ and using $\G_m^\dagger\cong\G_a^\dagger$ via the logarithm together with Cartier duality for $\G_a^\dagger$ yields an operator $\Psi$ which is given by multiplication by $-i$ on $\Fil_i V\cdot u^i$ as $\ol{\T}$ acts on $u$ by division. Thus, we can reformulate our goal as showing that $\Theta=uD+\Psi$. To check this, observe that the $\O(\G_a^\dagger\rtimes\G_a^\dagger)$-comodule structure corresponding to a $\G_a^\dagger\rtimes\G_a^\dagger$-representation, where $\G_a^\dagger$ acts on itself via $\ell.w\coloneqq w/\exp\ell$, is given by 
\begin{equation*}
\sum_{i, j\geq 0} \frac{\ell^i}{i!}\cdot\frac{w^j}{j!}\cdot \Psi^i D^j
\end{equation*}
in terms of the operators $D$ and $\Psi$. Restricting this along $\ell\mapsto (u(1-\exp(-\ell)), \ell)$ yields the $\O(\G_a^\dagger)$-comodule structure given by 
\begin{equation*}
\sum_{i, j\geq 0} \frac{\ell^i}{i!}\cdot\frac{u^j(1-\exp(-\ell))^j}{j!}\cdot \Psi^i D^j
\end{equation*}
and comparing this with the formula $\sum_{k\geq 0} \tfrac{\ell^k}{k!}\cdot\Theta^k$, which gives the $\O(\G_a^\dagger)$-comodule structure corresponding to the operator $\Theta$ under Cartier duality, we see that $\Theta=uD+\Psi$ by comparing the terms which are linear in $\ell$, as desired.
\end{proof}

Let us now move on to the higher-dimensional case and consider a derived Berkovich space $X$ equipped with a Berkovich étale map $X\rightarrow\ol{\T}^n$ for some $n\geq 0$. To explicitly describe $X^{\HT, \dagger, +}$ in this case, let us again start by introducing some notation: Let $\widetilde{\cal{G}}_n$ denote the group stack whose underlying Gelfand stack is the product $(\Z_p^\la)^n\times \Z_p^{\times, \la}\times (\G_a^\dagger)^n\times \G_a^\dagger\times\ol{\T}$ with the group operation given by the formula
\begin{equation*}
(\ul{\theta}_1, \gamma_1, \ul{v}_1, w_1, s_1)\cdot (\ul{\theta}_2, \gamma_2, \ul{v}_2, w_2, s_2)=(\ul{\theta}_1+\gamma_1\ul{\theta}_2, \gamma_1\gamma_2, \ul{v}_1+\gamma_1s_1(\ul{v}_2+w_1\ul{\theta}_2), w_1+w_2/s_1, s_1s_2)\;,
\end{equation*}
where $\ul{\theta}_i=(\theta_i^{(1)}, \dots, \theta_i^{(n)})$ and similarly for $\ul{v}_i$ and we write
\begin{equation*}
\gamma_2\ul{\theta}_1\coloneqq (\gamma_2\theta_1^{(1)}, \dots, \gamma_2\theta_1^{(n)})
\end{equation*}
and similarly in the other cases. After base change to $\ol{\DD}=\GSpec\Q_p\langle u\rangle_{\leq 1}$, the group stack $\widetilde{\cal{G}}_n$ admits two maps from $(\G_a^\dagger)^n\rtimes\G_m^\dagger$, where $\G_m^\dagger$ acts on $(\G_a^\dagger)^n$ by multiplication: The first one is given by 
\begin{equation*}
(\G_a^\dagger)^n\rtimes\G_m^\dagger\mapsto \widetilde{\cal{G}}_n\;, \hspace{0.3cm} (\ul{\theta}, s)\mapsto (0, 1, u\ul{\theta}, u(1-s^{-1}), s)
\end{equation*}
and the second one is given by $(\ul{\theta}, s)\mapsto (\ul{\theta}, s, 0, 0, 1)$. We let $\cal{G}_n$ be the coequaliser of these two maps in the category of group stacks over $\ol{\DD}$. Furthermore, for a derived Berkovich space $X$ equipped with a Berkovich étale map $X\rightarrow\ol{\T}^n$, we write
\begin{equation*}
X_\infty^\la\coloneqq X\times_{\ol{\T}^n} \lim_{x_i\mapsto x_i^p} \ol{\T}^n\;,
\end{equation*}
where the map $\lim_{x_i\mapsto x_i^p} \ol{\T}^n\rightarrow\ol{\T}^n$ is the first projection.

\begin{thm}
\label{thm:htdr-xhtdrdagger+pres}
Let $X$ be a derived Berkovich space over $\Q_p$ equipped with a Berkovich étale map $X\rightarrow \ol{\T}^n$ for some $n$. Then there is an isomorphism
\begin{equation}
\label{eq:htdr-xhtdrdagger+pres}
X^{\HT, \dagger, +}\cong X_\infty^\la\times \G_a^\dagger\times\ol{\DD}\times\GSpec\Q_p(\zeta_{p^\infty})\,\Big/\,\cal{G}_n\;,
\end{equation}
where the action of $\cal{G}_n$ on $X_\infty^\la\times \G_a^\dagger\times\ol{\DD}\times\GSpec\Q_p(\zeta_{p^\infty})$ is pulled back from the action of $\cal{G}_n$ on $\ol{\T}^{n, \la}_\infty\times\G_a^\dagger\times\ol{\DD}\times\GSpec\Q_p(\zeta_{p^\infty})$ given by the formula
\begin{equation*}
(\ul{\theta}, \gamma, \ul{v}, w, s).((\ul{x}^{1/p^m})_m, t, u, (\zeta_{p^m})_m)\coloneqq \left(\left(\zeta_{p^{m-1}}^{\ul{\theta}}\exp\left(\frac{t(u\ul{\theta}+\ul{v})}{p^m}\right)\ul{x}^{1/p^m}\right)_m, \gamma s t, \frac{u}{s}+w, (\zeta_{p^m}^\gamma)_m\right)\;,
\end{equation*}
where $(\ul{x}^{1/p^m})_m, t$ and $u$ denote the coordinates on $\ol{\T}^{n, \la}_\infty, \G_a^\dagger$ and $\ol{\DD}$, respectively, and we set $\zeta_{p^{-1}}\coloneqq 1$ for notational convenience. In particular, we have
\begin{equation*}
X^{\HT, \dagger}\cong X_\infty^\la\times \G_a^\dagger\times\GSpec\Q_p(\zeta_{p^\infty})\,\Big/\,(\Z_p^\la)^n\rtimes\Z_p^{\times, \la}\;,
\end{equation*}
where $\Z_p^{\times, \la}$ acts on $(\Z_p^\la)^n$ by multiplication.
\end{thm}
\begin{proof}
By compatibility of the functor $X\mapsto X^{\HT, \dagger, +}$ with Berkovich étale maps, which one deduces from the case $X\mapsto X^\N$ treated in \cref{thm:defis-berketalecover} by base change, we may reduce to the case $X=\ol{\T}^n$, from where we may further reduce to $n=1$ by compatibility of $X\mapsto X^{\HT, \dagger, +}$ with limits using that we have already worked out the case $X=\GSpec\Q_p$ in \cref{prop:htdr-qphtdagger+}. 

In the case $X=\ol{\T}$, recall that $\ol{\T}^\N\cong \ol{\T}^{\ol{\N}}\times_{\Q_p^{\ol{\N}}} \Q_p^\N$. Thus, putting together \cref{prop:htdr-qphtdagger+} and \cref{thm:pres-tn} at once implies the claim as soon as we check that the action
\begin{equation*}
(\theta, v).x^{1/p^m}=q^{\theta/p^m}\exp\left(\frac{tv}{p^m}\right)x^{1/p^m}
\end{equation*}
from loc.\ cit.\ transforms into the action described above after base changing along the map $\Q_p^{\HT, \dagger, +}\rightarrow \Q_{p, (p^{1/2}, p^{3/2})}^\N$. However, this simply follows by recalling that the proof of \cref{prop:htdr-qphtdagger+} shows that $q^{1/p^m}=\zeta_{p^{m-1}}\exp(tu/p^m)$ for all $m\geq 0$ on $\Q_p^{\HT, \dagger, +}$.
\comment{
In this case, recall from \cref{lem:syn-xnviacone} that the stack $\ol{\T}^{\HT, \dagger, +}$ sits in a pullback square
\begin{equation}
\label{eq:htdr-thtdagger+pullback}
\begin{tikzcd}
\ol{\T}^{\HT, \dagger, +}\ar[r]\ar[d] & \ol{\T}^\dR\ar[d] \\
\ol{\T}^{\Cone}\ar[r] & (\ol{\T}^{\Cone})^\dR
\end{tikzcd}
\end{equation}
over $\Q_p^{\HT, \dagger, +}$ and hence our task is to show how the right-hand side of (\ref{eq:htdr-xhtdrdagger+pres}) fits into the same pullback square. To start, observe that, using the presentation of $\Q_p^{\HT, \dagger, +}$ from \cref{prop:htdr-qphtdagger+}, the right-hand side of (\ref{eq:htdr-xhtdrdagger+pres}) indeed admits a map to $\Q_p^{\HT, \dagger, +}$ via projection onto the last three factors ``in the numerator'' and the map $\cal{G}_1\rightarrow\cal{G}_0$ given by dropping the $v$-coordinate ``in the denominator''.

To give a map from the right-hand side of (\ref{eq:htdr-xhtdrdagger+pres}) to $\ol{\T}^\dR$, note that the coordinate $x$ on $\ol{\T}^\la_\infty$ is well-defined up to multiplication by $\exp(t(u\theta+v))\in \G_m^\dagger$ on the stack on the right-hand side of (\ref{eq:htdr-xhtdrdagger+pres}) and hence defines a map to $\ol{\T}^\dR$. Expanding $\exp(t(u\theta+v))$ into its power series, we furthermore see that $\exp(t(u\theta+v))x=x+t(\,\dots)$, where the term in brackets is an explicit power series in $x, t, u, \theta$ and $v$, and thus $x$ also defines a map from the right-hand side of (\ref{eq:htdr-xhtdrdagger+pres}) to $\ol{\T}^{\Cone}$ as the generalised Cartier divisor defining the cone stack here is given by multiplication by $t$. As the map $\ol{\T}^\dR\rightarrow (\ol{\T}^{\Cone})^\dR$ in (\ref{eq:htdr-thtdagger+pullback}) is also induced by the coordinate on $\ol{\T}^\dR$, we conclude that the maps from the right-hand side of (\ref{eq:htdr-xhtdrdagger+pres}) to $\ol{\T}^\dR$ and $\ol{\T}^{\Cone}$ we have given are compatible upon mapping to $(\ol{\T}^{\Cone})^\dR$.

Thus, we obtain a map
\begin{equation}
\label{eq:htdr-thtdagger+map}
\ol{\T}^\la_{\infty}\times\G_a^\dagger\times\ol{\DD}\times\GSpec\Q_p(\zeta_{p^\infty})\,\Big/\,\cal{G}_1\rightarrow \ol{\T}^{\HT, \dagger, +}
\end{equation}
over $\Q_p^{\HT, \dagger, +}$. To see that this is an isomorphism, we note that the exact sequence
\begin{equation*}
\begin{tikzcd}
0\ar[r] & \G_m^\dagger/\G_a^\dagger\ar[r] & \ol{\T}^{\Cone}\ar[r] & (\ol{\T}^{\Cone})^\dR\ar[r] & 0\nospacepunct{\;,}
\end{tikzcd}
\end{equation*}
where $\G_a^\dagger$ acts on $\G_m^\dagger$ via $v.x\coloneqq \exp(tv)x$, shows that $\ol{\T}^{\HT, \dagger, +}$ is a $\G_m^\dagger/\G_a^\dagger$-torsor over $\ol{\T}^\dR\times\Q_p^{\HT, \dagger, +}$ by virtue of (\ref{eq:htdr-thtdagger+pullback}). However, also the map
\begin{equation*}
\ol{\T}^\la_{\infty}\times\G_a^\dagger\times\ol{\DD}\times\GSpec\Q_p(\zeta_{p^\infty})\,\Big/\,\cal{G}_1\rightarrow \ol{\T}^\dR\times\Q_p^{\HT, \dagger, +}
\end{equation*}
we have constructed is a $\G_m^\dagger/\G_a^\dagger$-torsor: Indeed, we have
\begin{equation*}
\ol{\T}^\la_\infty\,/\,\G_m^\dagger=\lim_{x\mapsto x^p} \ol{\T}\,\Big/\,\G_m^\dagger\cong \lim_{x\mapsto x^p} \ol{\T}^\dR
\end{equation*}
as the map $x\mapsto x^p$ on $\ol{\T}$ is $\dagger$-formally smooth and the canonical map $\cal{G}_1\rightarrow\Z_p^\sm\rtimes\cal{G}_0$ is a $\G_a^\dagger$-torsor, where $\cal{G}_0$ acts on $\Z_p^\sm$ by multiplication through the projection map $\cal{G}_0\rightarrow\Z_p^{\times, \sm}$. This now implies the claim using the presentation of $\Q_p^{\HT, \dagger, +}$ from \cref{prop:htdr-qphtdagger+} together with the fact that $\lim_{x\mapsto x^p} \ol{\T}^\dR\,/\,\Z_p^\sm\cong \ol{\T}^\dR$ over $\GSpec\Q_p(\zeta_{p^\infty})$. Overall, we conclude that (\ref{eq:htdr-thtdagger+map}) is an isomorphism upon noting that it is a map of $\G_m^\dagger/\G_a^\dagger$-torsors over $\ol{\T}^\dR\times\Q_p^{\HT, \dagger, +}$.
}
Finally, one obtains the claimed formula for $\ol{\T}^{\HT, \dagger}$ in the same way as in the proof of \cref{prop:htdr-qphtdagger+}: If $|u|=1$, we have $(\Z_p^\la\rtimes\Z_p^{\times, \la})\times\ol{\T}\cong\cal{G}_1$ via the natural map from left to right, where $\Z_p^{\times, \la}$ acts on $\Z_p^\la$ by multiplication, and hence
\begin{equation*}
\ol{\T}^{\HT, \dagger}\cong \ol{\T}^\la_\infty\times\G_a^\dagger\times\ol{\T}\times\GSpec\Q_p(\zeta_{p^\infty})\,\Big/\,(\Z_p^\la\rtimes\Z_p^{\times, \la})\times\ol{\T}\;.
\end{equation*}
Twisting the ``numerator'' of the right-hand side by the automorphism of $\G_a^\dagger\times\ol{\T}$ given by $(t, u)\mapsto (tu, u)$, the action of $\ol{\T}$ on $\G_a^\dagger$ becomes trivial and the claim follows.
\end{proof}

\begin{cor}
\label{cor:htdr-complexesxhtdrdagger+}
Let $X$ be a derived Berkovich space over $\Q_p$ equipped with a Berkovich étale map $X\rightarrow \ol{\T}^n$ for some $n\geq 0$. Then a perfect complex on $X^{\HT, \dagger, +}$ amounts to the following data:
\begin{enumerate}[label=(\roman*)]
\item A diagram
\begin{equation*}
\begin{tikzcd}
\dots \ar[r,shift left=.5ex,"u"]
  & \ar[l,shift left=.5ex, "t"] \Fil_{i-1} V \ar[r,shift left=.5ex,"u"] & \ar[l,shift left=.5ex, "t"] \Fil_i V \ar[r,shift left=.5ex,"u"] & \ar[l,shift left=.5ex, "t"] \Fil_{i+1} V\ar[r,shift left=.5ex,"u"] & \ar[l,shift left=.5ex, "t"] \dots
\end{tikzcd}
\end{equation*}
of perfect complexes on $X_\infty^\la\times\GSpec\Q_p(\zeta_{p^\infty})\times \GSpec\Q_p\{r\}^\dagger$, where $r=ut$, with the property that $t: \Fil_\bullet V\rightarrow \Fil_{\bullet-1} V$ becomes an isomorphism for $\bullet\ll 0$ and $u: \Fil_\bullet V\rightarrow\Fil_{\bullet+1} V$ becomes an isomorphism for $\bullet\gg 0$;
\item $X_\infty^\la\times\GSpec\Q_p(\zeta_{p^\infty})$-linear maps $D: \Fil_i V\rightarrow \Fil_{i-1} V$ for each $i$ satisfying $Du=uD+1$ and commuting with $t$;
\item a ``$t$-connection'' $\nabla: V\rightarrow V\tensor_{\O_{X_\infty^\la}} \Omega_{X_\infty^\la}^1(-1)$ on $V$, i.e.\ a map satisfying the $t$-deformed Leibniz rule $\nabla(fm)=f\nabla(m)+t\cdot\mathrm{d}f$, which commutes with $D$ and whose restriction to $\Fil_i V$ is equipped with a factorisation through $\Fil_{i-1}V\tensor_{\O_{X_\infty^\la}} \Omega_{X_\infty^\la}^1(-1)$; here, the twist by $-1$ indicates that $\Z_p^\times$ acts on $\Omega_{X_\infty^\la}^1$ by division, i.e.\ $\gamma.\omega=\gamma^{-1}\omega$;
\item a semilinear locally analytic $(\Z_p^n\rtimes\Z_p^\times)$-action on $\Fil_\bullet V$ commuting with all the previous data; here, $\Z_p^\times$ acts on $\Z_p^n$ by multiplication and the semilinearity is asked with respect to the action of $\Z_p^n\rtimes\Z_p^\times$ on $X^\la_\infty\times\GSpec\Q_p\{r\}^\dagger\times\GSpec\Q_p(\zeta_{p^\infty})$ pulled back from the action of $\Z_p^n\rtimes\Z_p^\times$ on $\ol{\T}^{n, \la}_\infty\times\GSpec\Q_p\{r\}^\dagger\times\GSpec\Q_p(\zeta_{p^\infty})$ given by
\begin{equation*}
(\ul{\theta}, \gamma).((\ul{x}^{1/p^m})_m, r, (\zeta_{p^m})_m)\coloneqq \left(\left(\zeta_{p^{m-1}}^{\ul{\theta}}\exp\left(\frac{r\ul{\theta}}{p^m}\right)\ul{x}^{1/p^m}\right)_m, \gamma r, (\zeta_{p^m}^\gamma)_m\right)\;.
\end{equation*}
\item an identification between the ``arithmetic Sen operator'' $\Theta^\arithm: \Fil_i V\rightarrow \Fil_i V$ coming from the Lie algebra action of $\Z_p^\times$ and the operator $uD-i$ for each $i$;
\item an identification between the ``geometric Sen operator'' $\Theta^\geom: \Fil_i V\rightarrow \Fil_i V\tensor_{\O_{X_\infty^\la}} \Omega_{X_\infty^\la}^1(-1)$ coming from the Lie algebra action of $\Z_p^n$ and the $t$-connection $u\nabla$ for each $i$.
\end{enumerate}
Under this identification, restriction to $X^{\HT, \dagger}$ corresponds to forgetting the $u$-filtration and the $D$-maps while restriction to $X^\N_{|u|=|t|=0}$ corresponds to passing to the associated graded of the $u$-filtration and forgetting the $D$-maps. Finally, the isomorphism
\begin{equation*}
X^\N_{|u|=|t|=0}\cong X^{\dR, +}_{|t|=0}
\end{equation*}
from \cref{prop:defis-u0} Tate twists the $\Z_p^\times$-action on the $i$-th associated graded piece of the $u$-filtration by $i$, which trivialises the action of the Lie algebra of $\Z_p^\times$ and hence after descent along the map
\begin{equation*}
X_\infty^\la\times \GSpec\Q_p(\zeta_{p^\infty})\rightarrow X
\end{equation*}
via the now smooth $\Z_p^n\rtimes\Z_p^\times$-action yields a filtered perfect complex $\Fil^\bullet W$ on $X$ via the $t$-maps together with a flat connection $\nabla: W\rightarrow W\tensor_{\O_X} \Omega^1_X$ satisfying Griffiths transversality, i.e.\ the restriction of $\nabla$ to $\Fil^i W$ is equipped with a factorisation through $\Fil^{i-1} W\tensor_{\O_X} \Omega^1_X$.
\end{cor}
\begin{proof}
This follows from \cref{thm:htdr-xhtdrdagger+pres} using the same kind of unraveling as in the proof of \cref{prop:htdr-complexesqphtdagger+}. Indeed, item (i) is again a version of the Rees equivalence in the form from \cref{cor:recall-reesdt} while the operators $D$ and $\nabla$ in (ii) and (iii) come from the $\G_a^\dagger$- and $(\G_a^\dagger)^n$-factors in $\widetilde{\cal{G}}_n$, respectively, via Cartier duality as in \cref{lem:recall-cartierperf}. The locally analytic $\Z_p^n\rtimes\Z_p^\times$-action in (iv) is due to the factor $(\Z_p^n)^\la\rtimes\Z_p^{\times, \la}$ in $\widetilde{\cal{G}}_n$ and the compatibilities in (v) and (vi) follow from the definition of $\cal{G}_n$ as a coequaliser by similar computations as in the proof of \cref{prop:htdr-complexesqphtdagger+}. Finally, the Tate twist occurring in translating between perfect complexes on $X_{|u|=|t|=0}^\N$ and $X_{|t|=0}^{\dR, +}$ is again due to the fact that the isomorphism $X_{|u|=|t|=0}^\N\cong X_{|t|=0}^{\dR, +}$ itself involves a Tate twist.
\end{proof}

\begin{rem}
In principle, one could compress the data contained in a perfect complex on $X^{\HT, \dagger, +}$ even further: Namely, we could replace (ii) and (v) by
\begin{itemize}
\item[(ii')] a factorisation of $\Theta^\arithm+i: \Fil_i V\rightarrow \Fil_i V$ through $\Fil_{i-1} V$, where $\Theta^\arithm$ is the ``arithmetic Sen operator'' coming from the Lie algebra action of $\Z_p^\times$
\end{itemize}
and similarly replace (iii) and (vi) by
\begin{itemize}
\item[(iii')] a factorisation of the ``geometric Sen operator'' $\Theta^\geom: \Fil_i V\rightarrow \Fil_i V\tensor_{\O_{X_\infty^\la}} \Omega_{X_\infty^\la}^1(-1)$ coming from the Lie algebra action of $\Z_p^n$ through $\Fil_{i-1} V\tensor_{\O_{X_\infty^\la}} \Omega_{X_\infty^\la}^1(-1)$.
\end{itemize}
Indeed, by (v) and (vi), these factorisations will precisely be given by $D$ and $\nabla$, respectively. However, we have chosen not to do this because we want to emphasise the crucial role that $D$ and $\nabla$ play, e.g.\ in describing the restriction functors to $X^{\HT, \dagger}$ and $X^\N_{|u|=|t|=0}$ or in describing cohomology on $X^{\HT, \dagger, +}$.
\end{rem}

\begin{rem}
We do not expect there to be an easy coordinate-free linear algebraic description of the category $\Perf(X^{\HT, \dagger, +})$ in the style of \cref{cor:htdr-complexesxhtdrdagger+}. However, one \emph{can} give a coordinate-free description of $\Perf(X^{\HT, +})$ for any Berkovich smooth derived Berkovich space $X$ over $\Q_p$. Namely, a perfect complex on $X^{\HT, +}$ is equivalent to 
\begin{enumerate}[label=(\roman*)]
\item a filtered perfect complex $\Fil_\bullet V$  of $\widehat{\O}$-modules on the proétale site $X_\proet$ of $X$ together with
\item a factorisation of the geometric Sen operator on $\Fil_i V$ through $\Fil_{i-1} V\tensor_{\O_X} \Omega_X^1(-1)$ and
\item an identification of the arithmetic Sen operator on $\gr_i V$ with multiplication by $-i$.
\end{enumerate}
Indeed, one can construct maps $\widehat{X}\times \ol{\DD}/\ol{\T}\rightarrow X^{\HT, +}\rightarrow X\times \Q_p^{\HT, +}$, the second of which is a $\cal{T}_{X/\Q_p}^\dagger\{1\}$-gerbe. This yields a functor from $\Perf(X^{\HT, +})$ to the linear algebraic category described above and then one can check locally that it is an equivalence, i.e.\ in particular after picking a toric chart. However, we will not need this description of $\Perf(X^{\HT, +})$ in the following.
\end{rem}

\begin{cor}
\label{cor:htdr-cohomologyxhtdrdagger+}
Let $X$ be a derived Berkovich space over $\Q_p$ equipped with a Berkovich étale map $X\rightarrow \ol{\T}^n$ for some $n\geq 0$. Under the identifications from \cref{cor:htdr-complexesxhtdrdagger+}, the cohomology of a perfect complex on $X^{\HT, \dagger, +}$ is computed by the (underived) $\Z_p^n\rtimes\Z_p^\times$-invariants of the cohomology on $X_\infty^\la$ of the total complex of
\begin{equation*}
\begin{tikzcd}
\Fil_0 V\ar[r, "\nabla"]\ar[d, "D", swap] & \Fil_{-1} V\tensor_{\O_{X_\infty^\la}} \Omega^1_{X_\infty^\la}(-1)\ar[d, "D", swap]\ar[r, "\nabla"] & \Fil_{-2} V\tensor_{\O_{X_\infty^\la}} \Omega^2_{X_\infty^\la}(-2)\ar[r, "\nabla"]\ar[d, "D", swap] & \dots \\
\Fil_{-1} V\ar[r, "\nabla"] & \Fil_{-2} V\tensor_{\O_{X_\infty^\la}} \Omega^1_{X_\infty^\la}(-1)\ar[r, "\nabla"] & \Fil_{-3} V\tensor_{\O_{X_\infty^\la}} \Omega^2_{X_\infty^\la}(-2)\ar[r, "\nabla"] & \dots\nospacepunct{\;.}
\end{tikzcd}
\end{equation*}
Similarly, the cohomology of the pullback of $E$ to $X^{\HT, \dagger}$ is computed by the (underived) $\Z_p^n\rtimes\Z_p^\times$-invariants of the cohomology on $X_\infty^\la$ of the total complex of
\begin{equation*}
\begin{tikzcd}
V\ar[r, "\nabla"]\ar[d, "D", swap] & V\tensor_{\O_{X_\infty^\la}} \Omega^1_{X_\infty^\la}(-1)\ar[d, "D", swap]\ar[r, "\nabla"] & V\tensor_{\O_{X_\infty^\la}} \Omega^2_{X_\infty^\la}(-2)\ar[r, "\nabla"]\ar[d, "D", swap] & \dots \\
V\ar[r, "\nabla"] & V\tensor_{\O_{X_\infty^\la}} \Omega^1_{X_\infty^\la}(-1)\ar[r, "\nabla"] & V\tensor_{\O_{X_\infty^\la}} \Omega^2_{X_\infty^\la}(-2)\ar[r, "\nabla"] & \dots\nospacepunct{\;.}
\end{tikzcd}
\end{equation*}
\end{cor}
\begin{proof}
Using the presentation of $X^{\HT, \dagger, +}$ from \cref{thm:htdr-xhtdrdagger+pres}, the canonical projection map $\cal{G}_n\rightarrow (\Z_p^n)^\sm\rtimes\Z_p^{\times, \sm}$, where $\Z_p^{\times, \sm}$ acts by multiplication on $(\Z_p^n)^\sm$, induces a map 
\begin{equation}
\label{eq:htdr-cohomologyxhtdrdagger+}
X^{\HT, \dagger, +}\rightarrow X_\infty^\la\,\big/\,((\Z_p^n)^\sm\rtimes \Z_p^{\times, \sm})\;.
\end{equation}
From the description of the category of perfect complexes on $X^{\HT, \dagger, +}$ from \cref{cor:htdr-complexesxhtdrdagger+}, it is clear that the total complex of the first double complex above equipped with its induced $(\Z_p^n)^\sm\rtimes \Z_p^{\times, \sm}$-action computes the pushforward along the map (\ref{eq:htdr-cohomologyxhtdrdagger+}). Now note that profinite groups like $\Z_p^n\rtimes\Z_p^\times$ have no higher smooth cohomology in characteristic $0$ and hence computing cohomology on the target of (\ref{eq:htdr-cohomologyxhtdrdagger+}) amounts to taking underived $\Z_p^n\rtimes\Z_p^\times$-invariants of the cohomology on $X_\infty^\la$, as claimed.

Similarly, for the second assertion, observe that, by the presentation of $X^{\HT, \dagger}$ from \cref{thm:htdr-xhtdrdagger+pres}, the pushforward along
\begin{equation*}
X^{\HT, \dagger}\rightarrow X_\infty^\la\,\big/\, ((\Z_p^n)^\sm\rtimes \Z_p^{\times, \sm})
\end{equation*}
is computed by taking Lie algebra cohomology for the locally analytic action of $\Z_p^n\rtimes\Z_p^\times$ on $V$. Recalling that the arithmetic Sen operator $\Theta^\arithm$ is given by $u\nabla$ on all $\Fil_i V$ while the geometric Sen operator $\Theta^\geom$ is given by $uD$ on $\Fil_i V(i)$, we see that this Lie algebra cohomology is exactly computed by the total complex of the second double complex above. One concludes as in the previous paragraph.
\end{proof}

Having made the data of a perfect complex $E$ on $X^{\HT, \dagger, +}$ for a Berkovich smooth derived Berkovich space $X$ completely explicit in terms of linear algebra, at least in the case where $X$ is equipped with a Berkovich étale map $X\rightarrow \ol{\T}^n$ for some $n\geq 0$, we want to end with the following somewhat surprising result, which says that all this data can in fact be recovered from the restriction of $E$ to $X^{\HT, \dagger, +}_{|u|=0}\cong X^{\dR, +}_{|t|=0}$. In §\ref{subsect:drbundles}, we will come back to this result and make it explicit in terms of the linear algebra involved.

\begin{prop}
\label{prop:htdr-perfequiv}
Let $X$ be a Berkovich smooth derived Berkovich space over $\Q_p$. Then the pullbacks along the maps $X^{\dR, +}\cong X^\N_{|u|=0}\rightarrow X^\N_{|ut|=0}$ and $X^{\HT, \dagger, +}=X^\N_{|t|=0}\rightarrow X^\N_{|ut|=0}$ induce equivalences of categories 
\begin{equation*}
\Perf(X^{\dR, +})\cong \Perf(X^\N_{|ut|=0})\cong \Perf(X^{\HT, \dagger, +})\;.
\end{equation*}
\end{prop}
\begin{proof}
We only show the first equivalence, the second one is analogous. Moreover, note that the canonical maps induce an isomorphism
\begin{equation}
\label{eq:htdr-pushoutiso}
X^\N_{|t|=0}\coprod_{X^\N_{|u|=|t|=0}} X^\N_{|u|=0}\xrightarrow{\cong} X^\N_{|ut|=0}
\end{equation}
by a variant of \cref{lem:hkcomp-glueoverinterval} and thus it suffices to show that $\Perf(X^\N_{|t|=0})\cong \Perf(X^\N_{|u|=|t|=0})$ via pullback along the canonical map $X^\N_{|u|=|t|=0}\rightarrow X^\N_{|t|=0}$. 

To this end, by compatibility of $X\mapsto X^\N$ with rational localisations, see \cref{prop:defis-openloc}, we may assume that $X$ admits a finite étale map $X\rightarrow Z$ to some rational subspace $Z\subseteq \ol{\T}^n$ for some $n\geq 0$. Then $X^\N_{[r, s]}$ is nicely coverable for any $[r, s]\subseteq (p^{1/2}, p^{3/2})$ by \cref{prop:perf-coversmoothrigid}. In particular, using \cref{lem:perf-bc} and \cref{lem:perf-finflatdim}, we conclude that $X^\N_{|t|=0}$ is nicely coverable by base change.

Now consider the overconvergent normed divisor $Z\coloneqq\{|u|=0\}\subseteq X^\N_{|t|=0}$ and observe that
\begin{equation*}
Z_\epsilon\setminus Z= X^\N_{0<|u|\leq \epsilon, |t|=0}\cong X^{\HT, \dagger}\times (0, \epsilon]
\end{equation*}
for all $0<\epsilon\leq 1$ by \cref{prop:defis-t0uneq0}, where we note that $X^\N_{|t|=0}\setminus Z=Z_1\setminus Z$. Thus, by \cref{lem:perf-betti}, we conclude that
\begin{equation*}
\colim_{\epsilon>0} \Perf(Z_\epsilon\setminus Z)\cong \Perf(X^{\HT, \dagger})\cong \Perf(X^\N_{|t|=0}\setminus Z)\;,
\end{equation*}
and now the desired equivalence $\Perf(X^\N_{|t|=0})\cong \Perf(Z)$ follows by applying \cref{cor:perf-corkeylemma}.
\end{proof}

\newpage

\section{Six functors for the syntomification}

Using the explicit presentations of Nygaardifications from the previous section, we are now in a position to prove that the functor $X\mapsto X^\Syn$ interacts nicely with the six functor formalism on Gelfand stacks. More explicitly, we want to prove the following properties:
\begin{enumerate}[label=(\arabic*)]
\item The map $\Q_p^\Syn\rightarrow\GSpec\Q_p$ is cohomologically smooth with dualising sheaf $\O\{1\}[3]$.

\item If $f: X\rightarrow Y$ is a rigid smooth map of derived Berkovich spaces over $\Q_p$ which is pure of relative dimension $d$, then the induced map $f^\Syn: X^\Syn\rightarrow Y^\Syn$ is cohomologically smooth with dualising sheaf $\O\{d\}[2d]$. 

\item If $f: X\rightarrow Y$ is a map of derived Berkovich spaces over $\Q_p$ which is locally of finite presentation and quasicompact, then $f^\Syn$ is weakly cohomologically proper in the sense of \cite[Def.\ 3.1.19]{dRStack}.
\end{enumerate}
Finally, we will show that analytic syntomic cohomology affords a \emph{strong theory of first Chern classes} in the sense of \cite[Def.\ 5.2.8]{Zavyalov}. In particular, for any derived Berkovich space $X$, there is a natural map
\begin{equation*}
c_1^\Syn: R\Gamma(X_{\mathrm{Berk}\text{-}\et}, \G_m)[1]\rightarrow R\Gamma(X^\Syn, \O\{1\}[2])\;,
\end{equation*}
where $X_{\mathrm{Berk}\text{-}\et}$ denotes the Berkovich étale site of $X$, and this induces a natural isomorphism between the pushforward of the structure sheaf along $(\P^1)^\Syn\rightarrow\Q_p^\Syn$ and the sheaf $\O\oplus \O\{-1\}[-2]$.

In particular, this will let us deduce the following form of relative Poincaré duality for analytic syntomic cohomology: For any rigid smooth morphism $f: X\rightarrow Y$ of derived Berkovich spaces which is quasicompact and any perfect analytic $F$-gauge $E\in\Perf(X^\Syn)$ on $X$, the pushforward $f^\Syn_*E$ is perfect on $Y^\Syn$ and there is an isomorphism
\begin{equation*}
f^\Syn_*(E^\vee\{d\})[2d]\cong (f^\Syn_*E)^\vee\;.
\end{equation*}
For absolute syntomic cohomology, in view of the lack of weak cohomological properness of $\Q_p^\Syn$, the only form of Poincaré duality we can obtain is an isomorphism
\begin{equation*}
R\Gamma(X^\Syn, E^\vee\{d+1\})[2d+3]\cong R\Gamma_c(X^\Syn, E)^\vee
\end{equation*}
for any rigid smooth derived Berkovich space $X$ over $\Q_p$ which is pure of relative dimension $d$ and any perfect analytic $F$-gauge $E\in\Perf(X^\Syn)$ on $X$; here, $R\Gamma_c$ denotes the $!$-pushforward to the point. Note that one can somewhat ``see'' the failure of $\Q_p^\Syn$ to be weakly cohomologically proper from the geometry: Namely, this should intuitively be explained by the fact that $\Q_p^\N$ is ``open at $\infty$'' with respect to the radius map $\kappa: \Q_p^\N\rightarrow (0, \infty)$.

Let us quickly say something about the strategy of the proofs: As we can fully describe $X^\N$ away from $p\in (0, \infty)$ in terms of $X^{\Div^1}$ and the Hyodo--Kato stack $X^\HK$ and both cohomological smoothness as well as weak cohomological properness may be checked locally on the source (in an open or finite closed cover, respectively), the main problem is to get a handle on $X^\N_{(p^{1/2}, p^{3/2})}$. However, here the previous section provides explicit presentations in the cases $X=\G_m$ and $X=\GSpec\Q_p$.

To illustrate how one finishes from here, let us shortly discuss how we prove cohomological smoothness of the stack
\begin{equation}
\label{eq:sixfunctors-presqpnintro}
\Q_{p, (p^{1/2}, p^{3/2})}^{\ol{\N}}\cong (\A^1\times \A^1)\times_{\A^1} \left(\lim_{q\mapsto q^p} (1+\overcirc{\DD})\setminus\{1\}\right)_{(p^{1/2}, p^{3/2})}\;\Big/\;\cal{G}\;,
\end{equation}
where the isomorphism is due to \cref{thm:pres-qpn}. The main issue is the infinite limit along $p$-th power maps on the right-hand side of (\ref{eq:sixfunctors-presqpnintro}): each finite stage of the limit is cohomologically smooth, but the limit is not. To deal with this, similarly to what happens in \cite[Prop.\ 5.7.1]{dRFF} or in \cite{Mikami}, the main idea is to approximate the quotient on the right-hand side of (\ref{eq:sixfunctors-presqpnintro}) by truncating the inverse limit in the ``numerator'' at finite stages and quotienting out by slightly larger groups $\cal{G}^{h\dagger}$ for $h>0$; these are defined in the same way as $\cal{G}$, but with $\Z_p^{\times, \la}$ replaced by $\Z_p^{\times, h\dagger}$ as defined in \cite[Def.\ 2.2]{Mikami}. To make this work, the crucial point will be that the algebras of functions on finite stages of the limit
\begin{equation*}
\left(\lim_{q\mapsto q^p} (1+\overcirc{\DD})\setminus\{1\}\right)_{[(p-1)/p, p^2(p-1)]}
\end{equation*}
equipped with the usual $\Z_p^\times$-action given by $\gamma.q^{1/p^m}=q^{\gamma/p^m}$ satisfy the Tate--Sen conditions from \cite[§3]{Mikami} by a classical calculation of Berger, see \cite[Prop.\ 1.1.12]{BergerBPaires} and \cite[Thm.\ 4.4]{BergerMultivariable}. Indeed, this will have the effect that the $h\dagger$-analytic vectors in the ring of functions on the limit are already contained in the ring of functions on some finite stage by \cite[Cor.\ 3.7]{Mikami}, and this is the key property that powers the approximation argument sketched above.

\subsection{Cohomological smoothness of $\Q_p^\Syn$}

We start by proving cohomological smoothness of $\Q_p^\Syn$ over $\GSpec\Q_p$. More precisely, we are going to prove:

\begin{thm}
\label{thm:sixfunctors-qpsynsmooth}
The map $\Q_p^\Syn\rightarrow\GSpec\Q_p$ is cohomologically smooth with dualising sheaf $\O\{1\}[3]$.
\end{thm}

Let us state \cref{thm:sixfunctors-qpsynsmooth} in a slightly more explicit form. For this, we first introduce some more notation.

\begin{defi}
Let $X$ be any Gelfand stack and let $f: X^\Syn\rightarrow\GSpec\Q_p$ be the structure map. For any $i\in\Z$, the \emph{analytic syntomic cohomology with compact support} of weight $i$ of $X$ is defined by
\begin{equation*}
R\Gamma_{\Syn, c}(X, \Q_p(i))\coloneqq f_!\O\{i\}\;.
\end{equation*}
\end{defi}

\begin{cor}
\label{cor:sixfunctors-qpsynsmooth}
Let $f: \Q_p^\Syn\rightarrow\GSpec\Q_p$ be the structure map. For any perfect analytic $F$-gauge $E\in\D(\Q_p^\Syn)$, there is an isomorphism $f_*(E^\vee\{1\})[3]\cong (f_!E)^\vee$. In particular, we have
\begin{equation*}
R\Gamma_\Syn(\GSpec\Q_p, \Q_p(i))\cong R\Gamma_{\Syn, c}(\GSpec\Q_p, \Q_p(1-i))^\vee[-3]
\end{equation*}
for any $i\in\Z$.
\end{cor}
\begin{proof}
Note that \cref{thm:sixfunctors-qpsynsmooth} implies
\begin{equation*}
(f_!E)^\vee=\Hom(f_!E, 1)\cong \Hom(E, f^!1)\cong f_*\sHom(E, f^!1)\cong f_*(E^\vee\{1\}[3])
\end{equation*}
by dualisability of $E$, see \cref{prop:recall-fredholm}.
\end{proof}

As already discussed above, the fact that $\Q_p^\N$ can be completely described in terms of $\Q_p^\HK$ and $\Q_p^{\Div^1}$ away from $p\in (0, \infty)$ means that proving \cref{thm:sixfunctors-qpsynsmooth} largely comes down to analysing $\Q_p^\N$ in a neighbourhood of $p\in (0, \infty)$. Namely, the key step in the proof is the following result:

\begin{prop}
\label{prop:sixfunctors-qpsynkeypart}
The map $\Q_{p, (p^{1/2}, p^{3/2})}^{\ol{\N}}\rightarrow\GSpec\Q_p$ is cohomologically smooth with dualising sheaf concentrated in degree $-2$.
\end{prop}

To prove \cref{prop:sixfunctors-qpsynkeypart}, we will use the presentation
\begin{equation*}
\Q_{p, (p^{1/2}, p^{3/2})}^{\ol{\N}}\cong (\A^1\times \A^1)\times_{\A^1} \left(\lim_{q\mapsto q^p} (1+\overcirc{\DD})\setminus\{1\}\right)_{(p^{1/2}, p^{3/2})}\;\Big/\;\cal{G}
\end{equation*}
from \cref{thm:pres-qpn}. To handle the inverse limit along $p$-th power maps, our strategy will be to invoke the very useful result \cite[Lem.\ 5.3.1]{dRFF}, which we restate here for the reader's convenience:

\begin{lem}
\label{lem:sixfunctors-drfflem}
Let $S$ be a Gelfand stack and consider a diagram
\begin{equation*}
X_\infty\rightarrow\dots\rightarrow X_{m+1}\rightarrow X_m\rightarrow\dots\rightarrow X_0
\end{equation*}
of Gelfand stacks over $S$ with $!$-able transition maps $f_{\ell m}: X_\ell\rightarrow X_m$ for $0\leq m\leq\ell\leq\infty$. Suppose that the following conditions hold:
\begin{enumerate}[label=(\alph*)]
\item All maps $f_{\ell m}$ are prim.
\item For any $m$, the natural map $\colim_{\ell\geq m} f_{\ell m, *}1\rightarrow f_{\infty m, *}1$ is an isomorphism.
\item There is a prim $!$-cover $\widetilde{X}_\infty\rightarrow X_\infty$ such that 
\begin{equation*}
\widetilde{X}_\infty\times_{X_\infty} \widetilde{X}_\infty\cong \lim_m \widetilde{X}_\infty\times_{X_m} \widetilde{X}_\infty
\end{equation*}
in the category of kernels over $S$.
\end{enumerate}
Then $X_\infty\cong\lim_m X_m$ in the category of kernels over $S$ and, in particular, 
\begin{equation*}
\D(X_\infty)\cong \mathop{\smash[b]{\sideset{}{_{!}}{\lim}}\vphantom{\lim}}_{m} \D(X_m)
\end{equation*}
via $(f_{\infty m, !})_m$, where the transition functors in the limit are lower-$!$-functors. Moreover, an object $P\in\D(X_\infty)$ is suave over $S$ if and only if $P_m\coloneqq f_{\infty m, !}P$ is suave over $S$ for all $m$ and, in that case, the suave dual is given by
\begin{equation*}
\SD(P)\cong \lim_m f_{\infty m}^\flat\SD(P_m)\;,
\end{equation*}
where $f_{\infty m}^\flat\cong \delta_{\infty m}\tensor f^*_{\infty m}$ is the left-adjoint of $f_{\infty m, !}$, $\delta_{\infty m}$ is the codualising sheaf of $f_{\infty m}$ and all suave duals $\SD(-)$ are taken over $S$.
\end{lem}
\begin{proof}
See \cite[Lem.\ 5.3.1]{dRFF}.
\end{proof}

To be able to apply \cref{lem:sixfunctors-drfflem}, let us first set up some notation: For $m\geq 3$, we let the Gelfand stack $X_m$ be defined by the pullback diagram
\begin{equation*}
\begin{tikzcd}
X_m\ar[r]\ar[d] & 1+\overcirc{\T}(p^{-1/(p^{m-5/2}(p-1))}, p^{-1/(p^{m-3/2}(p-1))})\ar[d, "q^{1/p^m}\mapsto \log q"] \\
\A^1\times\A^1\ar[r] & \A^1\nospacepunct{\;,}
\end{tikzcd}
\end{equation*}
where $q^{1/p^m}$ denotes the coordinate on $1+\overcirc{\T}(p^{-1/(p^{m-5/2}(p-1))}, p^{-1/(p^{m-3/2}(p-1))})$. Note that, by \cref{lem:pres-rhoqdrlogq} and the discussion preceding loc.\ cit., there are canonical projection maps
\begin{equation*}
X_\infty\coloneqq \widetilde{\Q_{p\mathrlap{, (p^{1/2}, p^{3/2})}}^{\ol{\N}}}\hphantom{\scriptstyle{(p^{1/2}, p^{3/2})}}\rightarrow X_m
\end{equation*}
for each $m$.

Furthermore, recall that, for any integer $h\geq 1$, we denote by $\ol{\DD}(1/p^h)$ the overconvergent disk of radius $1/p^h$, where we normalise the norm such that $|p|=1/p$; explicitly, this is the Gelfand stack $\GSpec\Q_p\langle T/p^h\rangle_{\leq 1}$. Using this and following \cite[Def.\ 2.2]{Mikami}, we define
\begin{equation*}
\Z_p^{\times, h\dagger}\coloneqq \bigsqcup_{\gamma\in\Z_p^{\times}/(1+p^h\Z_p)} \gamma(1+\ol{\DD}(1/p^h))
\end{equation*}
and observe that this canonically admits the structure of a multiplicative group stack over $\GSpec\Q_p$. Then let
\begin{equation*}
\widetilde{\cal{G}}_0^{h\dagger}\coloneqq \Z_p^{\times, h\dagger}\times (\G_a^\dagger\rtimes\ol{\T})
\end{equation*}
and note that, over $\A^1$, this receives two maps from $\G_m^\dagger$: The first one is given by $s\mapsto (1, u(1-s^{-1}), s)$, where $u$ denotes the coordinate on $\A^1$, while the second one is given by $s\mapsto (s, 0, 1)$. We define $\cal{G}_0^{h\dagger}$ to be the coequaliser of these two maps in the category of group stacks over $\A^1$. For notational convenience, we let $\cal{G}_0^{h\dagger}\coloneqq\cal{G}_0$ for $h=\infty$, which we recall was defined in the same way as $\cal{G}_0^{h\dagger}$ above, but with $\Z_p^{\times, \la}$ replacing $\Z_p^{\times, h\dagger}$. Then observe that there are canonical maps $\cal{G}_0^{h'\dagger}\rightarrow\cal{G}_0^{h\dagger}$ for all $h'>h$ induced by the canonical maps $\Z_p^{\times, h'\dagger}\rightarrow\Z_p^{\times, h\dagger}$, where $\Z_p^{\times, \infty\dagger}\coloneqq \Z_p^{\times, \la}$.

Now note that the action of $\cal{G}_0$ on $X_\infty$ from \cref{thm:pres-qpn} canonically extends to an action of $\cal{G}_0^{h\dagger}$ on $X_m$ for all $h\geq m$. Namely, the formula
\begin{equation*}
(\gamma, w, s).(t, u, q^{1/p^m})=\left(\gamma st, \frac{u}{s}+w, q^{\gamma/p^m}\exp\left(\frac{\gamma w s t}{p^m}\right)\right)_m
\end{equation*}
from loc.\ cit.\ makes sense for $\gamma\in\Z_p^{\times, h\dagger}$ whenever $h\geq m$: Indeed, first note that the exponential of $\gamma w s t/p^m$ does not play a role since $w\in\G_a^\dagger$, and then recall from the proof of \cref{thm:pres-tn} that the condition on $q^{1/p^m}$ implies that $|q-1|<p^{-1}$ and hence $|\log q|<p^{-1}$, from which we conclude that $|\log q^{1/p^m}|<p^{m-1}$. As the $p$-adic exponential has radius of convergence $p^{-1/(p-1)}$, we thus see that $\exp((\gamma-1)\log q^{1/p^m})$ defines a function on $1+\ol{\DD}(1/p^m)$, as desired.

With this notation in place, let us already say that our ultimate aim will be to apply \cref{lem:sixfunctors-drfflem} to the tower
\begin{equation*}
X_\infty/\cal{G}_0\rightarrow\dots\rightarrow X_{m+1}/\cal{G}_0^{(m+1)\dagger}\rightarrow X_m/\cal{G}_0^{m\dagger}\rightarrow\dots\rightarrow X_3/\cal{G}_0^{3\dagger}\;,
\end{equation*}
the transition maps in which we denote by $f_{\ell m}: X_\ell/\cal{G}_0^{\ell\dagger}\rightarrow X_m/\cal{G}_0^{m\dagger}$, and the $!$-cover $X_\infty\rightarrow X_\infty/\cal{G}_0$. This will almost immediately yield \cref{prop:sixfunctors-qpsynkeypart} and thus our main task is to verify the conditions of \cref{lem:sixfunctors-drfflem}. We start with the following lemma:

\begin{lem}
\label{lem:sixfunctors-limitghdagger}
Fix $m\geq 3$ and consider $\cal{G}_0^{h\dagger}\times_{\A^1} X_\infty$ as a Gelfand stack over $X_m$ via the composite map $\cal{G}_0^{h\dagger}\times_{\A^1} X_\infty\rightarrow \cal{G}_0^{h\dagger}\times_{\A^1} X_m\xrightarrow{\mathrm{act}}X_m$ for each $h\geq m$, where the map $X_m\rightarrow \A^1$ is the projection onto the $u$-coordinate. Then
\begin{equation*}
\cal{G}_0\times_{\A^1} X_\infty\cong \lim_{h\geq m}\,(\cal{G}_0^{h\dagger}\times_{\A^1} X_\infty)
\end{equation*}
in the category of kernels over $X_m$.
\end{lem}
\begin{proof}
Our goal is to apply \cref{lem:sixfunctors-drfflem} to the diagram
\begin{equation*}
\cal{G}_0\times_{\A^1} X_\infty=\cal{G}_0^{\infty\dagger}\times_{\A^1} X_\infty\rightarrow\dots\rightarrow\cal{G}_0^{(h+1)\dagger}\times_{\A^1} X_\infty\rightarrow\cal{G}_0^{h\dagger}\times_{\A^1} X_\infty\rightarrow\dots\rightarrow\cal{G}_0^{1\dagger}\times_{\A^1} X_\infty
\end{equation*}
and the prim $!$-cover $\widetilde{\cal{G}}_0\times_{\A^1} X_\infty$ of $\cal{G}_0\times_{\A^1} X_\infty$, which we note is a $\G_m^\dagger$-torsor. We let $g_{h'h}: \cal{G}_0^{h'\dagger}\times_{\A^1} X_\infty\rightarrow\cal{G}_0^{h\dagger}\times_{\A^1} X_\infty$ for $h'\geq h$ denote the transition maps in the diagram above and note that these are induced by the canonical maps $\cal{G}_0^{h'\dagger}\rightarrow \cal{G}^{h\dagger}_0$. 

Then first note that the diagram
\begin{equation}
\label{eq:sixfunctors-coverghdagger}
\begin{tikzcd}
\widetilde{\cal{G}}_0^{h'\dagger}\ar[r]\ar[d] & \widetilde{\cal{G}}_0^{h\dagger}\ar[d] \\
\cal{G}_0^{h'\dagger}\ar[r] & \cal{G}_0^{h\dagger}
\end{tikzcd}
\end{equation}
induced by the canonical map $\Z_p^{\times, h'\dagger}\rightarrow\Z_p^{\times, h\dagger}$ is cartesian as both vertical maps are $\G_m^\dagger$-torsors and the horizontal map on the top is compatible with the $\G_m^\dagger$-action; note that we have implicitly base changed $\widetilde{\cal{G}}_0^{h\dagger}$ to $\A^1$ here, and we will do this throughout. Now observe that the base change of the horizontal map on the top along any closed embedding
\begin{equation}
\label{eq:sixfunctors-closedcovertildeghdagger}
\Z_p^{\times, h\dagger}\times \G_a^\dagger\times (\{1/N\leq |s|\leq N\}\subseteq\G_m)\rightarrow \widetilde{\cal{G}}_0^{h\dagger}
\end{equation}
is a map between affine Gelfand stacks (after also base changing to $\ol{\DD}(N)\subseteq\A^1$). In other words, $\cal{G}_0^{h'\dagger}\rightarrow\cal{G}_0^{h\dagger}$ is a map between affine Gelfand stacks $!$-locally on $\cal{G}_0^{h\dagger}$ and thus it is weakly cohomologically proper by \cite[Lem.\ 3.1.21]{dRStack}. In particular, the same is true for $g_{h'h, *}$ by base change.

Next, we show that $\colim_{h'>h} g_{h'h, *}1\rightarrow g_{\infty h, *}1$ is an isomorphism for all $h\geq 1$. This may be checked $!$-locally on $\cal{G}_0^{h\dagger}\times X_\infty$ and we already know that each $g_{h'h}$ is prim, hence satisfies base change by \cite[Lem.\ 4.5.13]{HeyerMann}. Thus, letting $\widetilde{g}_{h'h}: X_\infty\times_{\A^1}\widetilde{\cal{G}}_0^{h'\dagger}\rightarrow X_\infty\times_{\A^1}\widetilde{\cal{G}}_0^{h\dagger}$ denote the canonical map, our claim reduces to checking that the natural map $\colim_{h'>h} \widetilde{g}_{h'h, *}1\rightarrow \widetilde{g}_{\infty h, *}1$ is an isomorphism by (\ref{eq:sixfunctors-coverghdagger}). Further localising along the closed embeddings (\ref{eq:sixfunctors-closedcovertildeghdagger}), which jointly cover $\widetilde{\cal{G}}_0^{h\dagger}$, as well as passing to a strict closed cover of $X_\infty$ by affine Gelfand stacks, this just reduces to the claim that
\begin{equation*}
\O(\Z_p^{\times, \la})\xrightarrow{\cong} \colim_{h>0} \O(\Z_p^{\times, h\dagger})
\end{equation*}
via the natural map, which holds by definition.

Finally, we claim that
\begin{equation}
\label{eq:sixfunctors-calgoverglimkernelcat}
(\widetilde{\cal{G}}_0\times_{\A^1} X_\infty)\times_{\cal{G}_0\times_{\A^1} X_\infty} (\widetilde{\cal{G}}_0\times_{\A^1} X_\infty)\rightarrow \lim_h\,(\widetilde{\cal{G}}_0\times_{\A^1} X_\infty)\times_{\cal{G}_0^{h\dagger}\times_{\A^1} X_\infty} (\widetilde{\cal{G}}_0\times_{\A^1} X_\infty)
\end{equation}
is an isomorphism in the category of kernels over $X_m$. For this, note that $X_\infty$ comes out of the fibre product and compute
\begin{equation*}
\widetilde{\cal{G}}_0\times_{\cal{G}_0^{h\dagger}}\widetilde{\cal{G}}_0\cong \widetilde{\cal{G}}_0\times_{\widetilde{\cal{G}}_0^{h\dagger}}(\widetilde{\cal{G}}_0^{h\dagger}\times_{\cal{G}_0^{h\dagger}}\widetilde{\cal{G}}_0^{h\dagger})\times_{\widetilde{\cal{G}}_0^{h\dagger}}\widetilde{\cal{G}}_0\cong \widetilde{\cal{G}}_0\times_{\widetilde{\cal{G}}_0^{h\dagger}}(\G_m^\dagger\times \widetilde{\cal{G}}_0^{h\dagger})\times_{\widetilde{\cal{G}}_0^{h\dagger}}\widetilde{\cal{G}}_0\cong \widetilde{\cal{G}}_0\times_{\widetilde{\cal{G}}_0^{h\dagger}}(\G_m^\dagger\times\widetilde{\cal{G}}_0)\;,
\end{equation*}
where the map $\G_m^\dagger\times\widetilde{\cal{G}}_0\rightarrow\widetilde{\cal{G}}_0^{h\dagger}$ is obtained by letting $\G_m^\dagger$ act on $\widetilde{\cal{G}}_0$ by the difference of the two morphisms $\G_m^\dagger\rightarrow\widetilde{\cal{G}}_0$ used to define $\cal{G}_0$. Twisting by this action, we conclude that
\begin{equation*}
(\widetilde{\cal{G}}_0\times_{\A^1} X_\infty)\times_{\cal{G}_0^{h\dagger}\times_{\A^1} X_\infty} (\widetilde{\cal{G}}_0\times_{\A^1} X_\infty)\cong (\widetilde{\cal{G}}_0\times_{\widetilde{\cal{G}}_0^{h\dagger}} \widetilde{\cal{G}}_0)\times_{\A^1} X_\infty\times \G_m^\dagger\;,
\end{equation*}
where the maps $\widetilde{\cal{G}}_0\rightarrow\widetilde{\cal{G}}_0^{h\dagger}$ are both the canonical maps. Observing that $\widetilde{\cal{G}}_0\cong \lim_h \widetilde{\cal{G}}_0^{h\dagger}$ as derived Berkovich spaces, we conclude that 
\begin{equation}
\label{eq:sixfunctors-calgoverglimberk}
\widetilde{\cal{G}}_0\times_{\A^1} X_\infty\times \G_m^\dagger\cong \lim_h\,(\widetilde{\cal{G}}_0\times_{\widetilde{\cal{G}}_0^{h\dagger}} \widetilde{\cal{G}}_0)\times_{\A^1} X_\infty\times \G_m^\dagger
\end{equation}
as derived Berkovich spaces via the natural map and, by the above calculation, it suffices to show that this is also an isomorphism in the kernel category over $X_m$ to be able to draw the same conclusion about (\ref{eq:sixfunctors-calgoverglimkernelcat}).

However, this last claim follows by the argument from the proof of \cite[Lem.\ 5.4.8]{dRFF}: Indeed, we may pass to strict closed covers of $X_\infty$ and $\widetilde{\cal{G}}_0$ which are compatible on the $\A^1$-factors and then all transition maps on the right-hand side of (\ref{eq:sixfunctors-calgoverglimberk}) will become maps between affine Gelfand stacks. Moreover, the structure map to $X_m$ will also factor through an affine closed subspace of $X_m$ and then the argument from the proof of \cite[Lem.\ 5.4.8]{dRFF} applies. As already announced in the beginning, applying \cref{lem:sixfunctors-drfflem} finishes the proof.
\end{proof}

\begin{cor}
\label{cor:sixfunctors-qpsmoothcond(c)}
In the category of kernels over $\GSpec\Q_p$, we have
\begin{equation*}
X_\infty\times_{X_\infty/\cal{G}_0} X_\infty\cong \lim_m X_\infty\times_{X_m/\cal{G}_0^{m\dagger}} X_\infty\;.
\end{equation*}
\end{cor}
\begin{proof}
For any $m$, we calculate 
\begin{equation*}
\begin{split}
X_\infty\times_{X_m/\cal{G}_0^{m\dagger}} X_\infty&\cong X_\infty\times_{X_m} (X_m\times_{X_m/\cal{G}_0^{m\dagger}} X_m)\times_{X_m} X_\infty \\
&\cong X_\infty\times_{X_m} (\cal{G}_0^{m\dagger}\times_{\A^1} X_m)\times_{X_m} X_\infty\cong X_\infty\times_{X_m} (\cal{G}_0^{m\dagger}\times_{\A^1} X_\infty)\;,
\end{split}
\end{equation*}
where the map $\cal{G}_0^{m\dagger}\times_{\A^1} X_\infty\rightarrow X_m$ is obtained as the composite $\cal{G}_0^{m\dagger}\times_{\A^1} X_\infty\rightarrow \cal{G}_0^{m\dagger}\times_{\A^1} X_m\xrightarrow{\mathrm{act}} X_m$. Now note that
\begin{equation*}
\lim_m\,(X_\infty\times_{X_m} (\cal{G}_0^{m\dagger}\times_{\A^1} X_\infty))\cong \lim_m\lim_{h\geq m}\,(X_\infty\times_{X_m} (\cal{G}_0^{h\dagger}\times_{\A^1} X_\infty))
\end{equation*}
in the category of kernels by cofinality and thus, it suffices to show that
\begin{equation}
\label{eq:sixfunctors-xinftyoverquotlimkernelcat}
\cal{G}_0\times_{\A^1} X_\infty\cong \lim_m\lim_{h\geq m}\,(X_\infty\times_{X_m} (\cal{G}_0^{h\dagger}\times_{\A^1} X_\infty))
\end{equation}
in the category of kernels. 

However, by \cref{lem:sixfunctors-limitghdagger}, we know that $\cal{G}_0\times_{\A^1} X_\infty\cong \lim_{h\geq m} \cal{G}_0^{h\dagger}\times_{\A^1} X_\infty$ in the category of kernels over $X_m$ and hence (\ref{eq:sixfunctors-xinftyoverquotlimkernelcat}) reduces to
\begin{equation*}
\cal{G}_0\times_{\A^1} X_\infty\cong \lim_m\,(X_\infty\times_{X_m} (\cal{G}_0\times_{\A^1} X_\infty))\;.
\end{equation*}
Now note that $\cal{G}_0$ acts on $X_\infty$ and, twisting by this action, we can achieve that the map $\cal{G}_0\times_{\A^1} X_\infty\rightarrow X_m$ becomes the precomposition of the canonical map $X_\infty\rightarrow X_m$ with the second projection. In other words, we can reduce further to showing that
\begin{equation*}
\cal{G}_0\times_{\A^1} X_\infty\cong \lim_m\,(\cal{G}_0\times_{\A^1}(X_\infty\times_{X_m} X_\infty))
\end{equation*}
and then it suffices to show that 
\begin{equation*}
X_\infty\cong \lim_m\,(X_\infty\times_{X_m} X_\infty)
\end{equation*}
in the category of kernels over $\A^1$. However, this is again deduced by the argument from \cite[Lem.\ 5.4.8]{dRFF} as at the end of the proof of \cref{lem:sixfunctors-limitghdagger}: Indeed, we have $X_\infty\cong\lim_m X_m$ and hence $X_\infty\cong \lim_m\,(X_\infty\times_{X_m} X_\infty)$ as derived Berkovich spaces.
\end{proof}

\begin{lem}
\label{lem:sixfunctors-xmgmtransitionproper}
For all $m\leq\ell\leq\infty$, the map $f_{\ell m}: X_\ell/\cal{G}_0^{\ell\dagger}\rightarrow X_m/\cal{G}_0^{m\dagger}$ is weakly cohomologically proper.
\end{lem}
\begin{proof}
Note that $X_\ell\rightarrow X_m$ becomes a map between affine Gelfand stacks after passing to a strict closed cover of $X_m$. Thus, the map $X_\ell/\cal{G}_0^{\ell\dagger}\rightarrow X_m/\cal{G}_0^{\ell\dagger}$ is a map between affine Gelfand stacks $!$-locally on the target and hence weakly cohomologically proper by \cite[Lem.\ 3.1.21]{dRStack}. It thus remains to show that $*/\cal{G}_0^{\ell\dagger}\rightarrow */\cal{G}_0^{m\dagger}$ is weakly cohomologically proper as we can then deduce the same for $X_m/\cal{G}_0^{\ell\dagger}\rightarrow X_m/\cal{G}_0^{m\dagger}$ by base change. In fact, we even claim that it is cohomologically proper.

For this, observe that $*/\cal{G}_0^{\ell\dagger}\rightarrow */\cal{G}_0^{m\dagger}$ is a $\Z_p^{\times, m\dagger}/\Z_p^{\times, \ell\dagger}$-torsor and note that $\Z_p^{\times, m\dagger}/\Z_p^{\times, \ell\dagger}$ is prim by \cite[Prop.\ 2.39.(5)]{Mikami}. Thus, it suffices to show that the diagonal of $\Z_p^{\times, m\dagger}/\Z_p^{\times, \ell\dagger}$ is cohomologically proper. However, note that there is a pullback square
\begin{equation*}
\begin{tikzcd}
\Z_p^{\times, m\dagger}\times\Z_p^{\times, \ell\dagger}\ar[r, "\mathrm{pr}_1\times \mathrm{add}"]\ar[d] & \Z_p^{\times, m\dagger}\times \Z_p^{\times, m\dagger}\ar[d] \\
\Z_p^{\times, m\dagger}/\Z_p^{\times, \ell\dagger}\ar[r, "\Delta"] & \Z_p^{\times, m\dagger}/\Z_p^{\times, \ell\dagger}\times \Z_p^{\times, m\dagger}/\Z_p^{\times, \ell\dagger}\nospacepunct{\;,}
\end{tikzcd}
\end{equation*}
in which the vertical arrow on the right is a $!$-cover and the horizontal arrow on the top is a closed embedding of derived Berkovich spaces, hence cohomologically proper: Indeed, twisting the target by the automorphism $(\gamma, \gamma')\mapsto (\gamma, \gamma'-\gamma)$ identifies $\mathrm{pr}_1\times \mathrm{add}$ with the map
\begin{equation*}
\Z_p^{\times, m\dagger}\times\Z_p^{\times, \ell\dagger}\xrightarrow{\mathrm{can}} \Z_p^{\times, m\dagger}\times \Z_p^{\times, m\dagger}\;.
\end{equation*}
Then note that $\Z_p^{\times, \ell\dagger}\rightarrow\Z_p^{\times, m\dagger}$ is a map between affine Gelfand stacks and that $\O(\Z_p^{\times, \ell\dagger})$ is an idempotent algebra over $\O(\Z_p^{\times, m\dagger})$: Indeed, if $\ell<\infty$, the map $\Z_p^{\times, \ell\dagger}\rightarrow\Z_p^{\times, m\dagger}$ is an embedding of a finite number of closed disks into a finite number of larger closed disks and if $\ell=\infty$, we have $\O(\Z_p^\la)=\colim_{\ell'<\infty} \O(\Z_p^{\ell'\dagger})$ and this is an idempotent $\O(\Z_p^{m\dagger})$-algebra as the same is true for each $\O(\Z_p^{\ell'\dagger})$. As cohomological properness is local on the target, we conclude that $\Delta$ is cohomologically proper, and this finishes the proof.
\end{proof}

The key lemma that goes into establishing \cref{prop:sixfunctors-qpsynkeypart} is the following:

\begin{lem}
\label{lem:sixfunctors-qpsynkeylem}
For any fixed $m\geq 5$, the natural map $f_{\ell m, *}1\rightarrow f_{\infty m, *}1$ is an isomorphism for $\ell\gg m$.
\end{lem}
\begin{proof}
As we already know that the $f_{\ell m}$ are prim, see \cref{lem:sixfunctors-xmgmtransitionproper}, they satisfy base change by \cite[Lem.\ 4.5.13]{HeyerMann}. Thus, we may check the claim locally on $X_m/\cal{G}_0^{m\dagger}$ and, in particular, after pullback to $X_m/\Z_p^{\times, m\dagger}$ via the canonical map $\Z_p^{\times, m\dagger}\rightarrow\cal{G}_0^{m\dagger}$. For this, note that the diagram
\begin{equation*}
\begin{tikzcd}
X_\ell/\Z_p^{\times, \ell\dagger}\ar[r, "\widetilde{f}_{\ell m}"]\ar[d] & X_m/\Z_p^{\times, m\dagger}\ar[d] \\
X_\ell/\cal{G}_0^{\ell\dagger}\ar[r, "f_{\ell m}"] & X_m/\cal{G}_0^{m\dagger}
\end{tikzcd}
\end{equation*}
is cartesian as both vertical maps are $(\G_a^\dagger\rtimes\ol{\T})/\G_m^\dagger$-torsors and the map $\widetilde{f}_{\ell m}$ is compatible with the $(\G_a^\dagger\rtimes\ol{\T})/\G_m^\dagger$-actions. 

We may also replace $X_m$ by 
\begin{equation*}
\ol{X}_m\coloneqq (\A^1\times \A^1)\times_{\A^1} (1+\T(p^{-1/p^{m-3}}, p^{-1/p^m}))\;,
\end{equation*}
where we write $\T(a, b)$ for the rigid torus of inner radius $a$ and outer radius $b$ and the coordinate on the torus is $q^{1/p^m}$. To see this, first note that $\ol{X}_m$ contains $X_m$ as an open subspace and has the feature that $\ol{X}_\ell\times_{\ol{X}_m} X_m\cong X_\ell$ for all $\ell\geq m$. We set $\ol{X}_\infty\coloneqq \lim_m \ol{X}_m$ and note that
\begin{equation*}
\lim_m \{p^{-1/p^{m-3}}\leq |q^{1/p^m}-1|\leq p^{-1/p^m}\}\cong \left(\lim_{q\mapsto q^p} (1+\overcirc{\DD})\setminus \{1\}\right)_{[(p-1)/p, p^2(p-1)]}\;.
\end{equation*}
Moreover, the formula $\gamma.q^{1/p^m}\mapsto q^{\gamma/p^m}$ defining the action of $\Z_p^{\times, m\dagger}$ on $X_m$ still makes sense on $\ol{X}_m$ due to $|q-1|\leq p^{-1}$ and hence $|\log q|\leq p^{-1}$ on $\ol{X}_m$. As $X_m\rightarrow \ol{X}_m$ is an open embedding and hence suave, by base change, see \cite[Lem.\ 4.5.3]{HeyerMann}, it then suffices to show that $\ol{f}_{\ell m, *} 1\cong \ol{f}_{\infty m, *} 1$ for $\ell\gg m$, where $\ol{f}_{\ell m}: \ol{X}_\ell/\Z_p^{\times, \ell\dagger}\rightarrow \ol{X}_m/\Z_p^{\times, m\dagger}$ is the canonical map.

Finally, we may also replace $\ol{X}_m$ by the exhaustive family of closed subspaces
\begin{equation*}
\ol{X}_{m, N}\coloneqq (\{|u|, |t|\leq N\}\subseteq \ol{X}_m)
\end{equation*}
for $N\geq 0$ and similarly for all $\ol{X}_\ell$. Then our task is to show that $\ol{f}_{\ell m, N, *}1\cong \ol{f}_{\infty m, N, *}1$ for all $N$ and all $\ell$ large enough (independent of $N$), where $\ol{f}_{\ell m, N}$ denotes the base change of $\ol{f}_{\ell m}$ to $\ol{X}_{m, N}/\Z_p^{\times, m\dagger}$.

As $\ol{X}_{\ell, N}$ is affine for all $\ell, N$, we have $\ol{f}_{\ell m, N, *}1\cong \O(\ol{X}_{\ell, N})^{m\dagger\text{-an}}$ by \cite[Cor.\ 2.50]{Mikami} and, similarly, $\ol{f}_{\infty m, N, *}1\cong \O(\ol{X}_{\infty, N})^{m\dagger\text{-an}}$. Thus, we have to prove that the natural map
\begin{equation*}
\O(\ol{X}_{\ell, N})^{m\dagger\text{-an}}\rightarrow \O(\ol{X}_{\infty, N})^{m\dagger\text{-an}}
\end{equation*}
is an isomorphism for $\ell\gg 0$. Now observe that
\begin{equation*}
\O(\ol{X}_{\ell, N})\cong \O(\ol{\DD}(N))\tensor\O(\ol{\DD}(N))\tensor_{\O(\ol{\DD}(N^2))} B_{[(p-1)/p, p^2(p-1)], \ell}
\end{equation*}
for $B_{[(p-1)/p, p^2(p-1)], \ell}$ being the ring of functions on $1+\T(p^{-1/p^{\ell-3}}, p^{-1/p^\ell})$. Recalling that $\Z_p^\times$ acts trivially on one copy of $\ol{\DD}(N)$ and by multiplication on the other, we see that the action of $\Z_p^\times$ on the $\O(\ol{\DD}(N))$-factors is even analytic globally and thus
\begin{equation*}
\O(\ol{X}_{\ell, N})^{m\dagger\text{-an}}\cong \O(\ol{\DD}(N))\tensor\O(\ol{\DD}(N))\tensor_{\O(\ol{\DD}(N^2))} B_{[(p-1)/p, p^2(p-1)], \ell}^{m\dagger\text{-an}}
\end{equation*}
by the projection formula for $m\dagger$-analytic vectors, see \cite[Prop.\ 2.39.(5)]{Mikami}.

Therefore, our claim comes down to showing that 
\begin{equation*}
B_{[(p-1)/p, p^2(p-1)], \ell}^{m\dagger\text{-an}}\rightarrow B_{[(p-1)/p, p^2(p-1)], \infty}^{m\dagger\text{-an}}
\end{equation*}
is an isomorphism for $\ell$ sufficiently large. However, this follows from \cite[Cor.\ 3.7]{Mikami} as the family of Banach algebras $\{B_{[(p-1)/p, p^2(p-1)], \ell}\}_\ell$ satisfies the Tate--Sen axioms from \cite[§3]{Mikami} by \cite[Prop.\ 1.1.12]{BergerBPaires} and \cite[Thm.\ 4.4]{BergerMultivariable}.
\end{proof}

Finally, we can prove \cref{prop:sixfunctors-qpsynkeypart}:

\begin{proof}[Proof of \cref{prop:sixfunctors-qpsynkeypart}]
We apply \cref{lem:sixfunctors-drfflem} to the tower
\begin{equation*}
X_\infty/\cal{G}_0\rightarrow\dots\rightarrow X_{m+1}/\cal{G}_0^{(m+1)\dagger}\rightarrow X_m/\cal{G}_0^{m\dagger}\rightarrow\dots\rightarrow X_3/\cal{G}_0^{3\dagger}\;.
\end{equation*}
Note that condition (a) of loc.\ cit.\ holds by \cref{lem:sixfunctors-xmgmtransitionproper} and condition (b) holds by \cref{lem:sixfunctors-qpsynkeylem}. Finally, condition (c) holds for the $!$-cover of $X_\infty/\cal{G}_0$ by $X_\infty$ due to \cref{cor:sixfunctors-qpsmoothcond(c)} using the fact that $\cal{G}_0$ is prim. Indeed, note that the canonical map $\cal{G}_0\rightarrow \Z_p^{\times, \sm}$ is a $\G_a^\dagger\rtimes\ol{\T}$-torsor and thus even weakly cohomologically proper by \cite[Lem.\ 3.1.21]{dRStack}; moreover, the stack $\Z_p^{\times, \sm}=\GSpec C^\sm(\Z_p^\times, \Q_p)$ is affine and hence weakly cohomologically proper as well.

To conclude from \cref{lem:sixfunctors-drfflem} that $X_\infty/\cal{G}_0$ is suave, we have to check that $f_{\infty m, !}1\in\D(X_m/\cal{G}_0^{m\dagger})$ is suave for all $m$. However, by \cref{lem:sixfunctors-xmgmtransitionproper} all maps $f_{\ell m}$ are weakly cohomologically proper and, in particular, using \cref{lem:sixfunctors-qpsynkeylem}, we conclude that $f_{\infty m, !}1\cong f_{\ell m, !}1$ for some finite $\ell\gg m$. As lower-$!$ along prim maps preserves suave objects by \cite[Lem.\ 4.5.16]{HeyerMann}, it thus suffices to show that $X_\ell/\cal{G}_0^{\ell\dagger}$ is suave. In fact, we are going to show that it is cohomologically smooth and that $f_{\ell m}^*\SD(1)\cong \SD(1)$ whenever $\ell\geq m$, where $\SD(-)$ denotes suave duals over the point.

For this, first observe that $X_\ell$ is a complete intersection inside the open subspace
\begin{equation*}
\A^1\times \A^1\times (1+\overcirc{\T}(p^{-1/(p^{m-5/2}(p-1))}, p^{-1/(p^{m-3/2}(p-1))}))
\end{equation*}
of $\A^3$ cut out by the equation $ut=\log q$ and hence cohomologically smooth. Moreover, the map $X_\ell\rightarrow X_m$ is finite étale for all $\ell\geq m$ and thus $X_\ell/\cal{G}_0^{\ell\dagger}\rightarrow X_m/\cal{G}_0^{\ell\dagger}$ is cohomologically étale. This together with suave base change for the cartesian diagram
\begin{equation*}
\begin{tikzcd}
X_m/\cal{G}_0^{\ell\dagger}\ar[r]\ar[d] & */\cal{G}_0^{\ell\dagger}\ar[d] \\
X_m/\cal{G}_0^{m\dagger}\ar[r] & */\cal{G}_0^{m\dagger}\nospacepunct{\;,}
\end{tikzcd}
\end{equation*}
see \cite[Lem.\ 4.5.13]{HeyerMann}, implies that the relative dualising sheaf of $X_m/\cal{G}_0^{m\dagger}\rightarrow */\cal{G}_0^{m\dagger}$ pulls back to the relative dualising sheaf of $X_\ell/\cal{G}_0^{\ell\dagger}\rightarrow */\cal{G}_0^{\ell\dagger}$.

Thus, it remains to show that $*/\cal{G}_0^{\ell\dagger}$ is cohomologically smooth and, for this, we note the factorisation
\begin{equation*}
*/\cal{G}_0^{\ell\dagger}\xrightarrow{g_\ell} */(\Z_p^{\times, \ell\dagger})^\dR\xrightarrow{h_\ell} *\;.
\end{equation*}
For $g_\ell$, note that it is a $\G_a^\dagger\rtimes\ol{\T}$-gerbe and hence cohomologically smooth: Indeed, this is because the trivial $\G_a^\dagger\rtimes\ol{\T}$-gerbe is cohomologically smooth due to the factorisation
\begin{equation*}
*/(\G_a^\dagger\rtimes\ol{\T})\rightarrow */\ol{\T}\rightarrow */\G_m\rightarrow *\;,
\end{equation*}
in which the first map is cohomologically smooth as it is a $\G_a^\dagger$-gerbe, see \cite[Thm.\ 4.3.13]{dRStack}, the second map is cohomologically smooth as it is a $(0, \infty)$-torsor by \cref{lem:defis-rhvariant} and the last map is cohomologically smooth as this may be checked smooth-locally on the source and $*\rightarrow */\G_m$ is a cohomologically smooth $!$-cover due to $\G_m$ being cohomologically smooth. Moreover, the diagram
\begin{equation*}
\begin{tikzcd}
*/\cal{G}_0^{\ell\dagger}\ar[r, "g_\ell"]\ar[d] & */(\Z_p^{\times, \ell\dagger})^\dR\ar[d] \\
*/\cal{G}_0^{m\dagger}\ar[r, "g_m"] & */(\Z_p^{\times, m\dagger})^\dR
\end{tikzcd}
\end{equation*}
is cartesian as both horizontal maps are $\G_a^\dagger\rtimes\ol{\T}$-gerbes and the vertical map on the left is compatible with the $*/\G_a^\dagger\rtimes\ol{\T}$-actions. Thus, by smooth base change, see \cite[Lem.\ 4.5.13]{HeyerMann}, we conclude that the pullback of $g_m^!1$ to $*/\cal{G}_0^{\ell\dagger}$ naturally identifies with $g_\ell^!1$.

For $h_\ell$, note that it factors as 
\begin{equation*}
*/(\Z_p^{\times, \ell\dagger})^\dR\rightarrow */(1+\ol{\DD}(1/p^\ell))^\dR\overset{\log}{\cong} */\ol{\DD}(1/p^\ell)^\dR\rightarrow */\G_a^\dR\rightarrow *\;,
\end{equation*}
where the first map is a gerbe for the finite abelian group $\Z_p^\times/(1+p^\ell\Z_p)$ and hence cohomologically étale (this may be checked étale-locally on the source and finite sets are cohomologically étale), the second map is a torsor for $\G_a^\dR/\ol{\DD}(1/p^\ell)^\dR\cong \G_a/\ol{\DD}(1/p^\ell)$ and thus cohomologically smooth by \cite[Prop.\ 4.3.7]{dRStack} and the last map is cohomologically smooth as this may be checked smooth-locally on the source and $\G_a^\dR$ is cohomologically smooth by \cite[Thm.\ 5.1.7]{dRFF}. For the compatibility with pullback, we only have to take care of the second map, where it suffices to give a $\G_a^\dR$-equivariant isomorphism $\widetilde{h}_\ell^!1\cong \widetilde{h}_{\ell m}^*\widetilde{h}_m^!1$, where we write $\widetilde{h}_\ell: \G_a^\dR/\ol{\DD}(1/p^\ell)^\dR\rightarrow *$, similarly for $m$, and $\widetilde{h}_{\ell m}: \G_a^\dR/\ol{\DD}(1/p^\ell)^\dR\rightarrow \G_a^\dR/\ol{\DD}(1/p^m)^\dR$. Using the isomorphism $\G_a^\dR/\ol{\DD}(1/p^\ell)^\dR\cong \G_a/\ol{\DD}(1/p^\ell)$, and similarly for $m$, it thus suffices to give an isomorphism $\ol{h}_\ell^!1\cong \ol{h}_{\ell m}^*\ol{h}_m^!1$, where we write $\ol{h}_\ell: */\ol{\DD}(1/p^\ell)\rightarrow *$, similarly for $m$, and $\ol{h}_{\ell m}: */\ol{\DD}(1/p^\ell)\rightarrow */\ol{\DD}(1/p^m)$. However, by \cite[Prop.\ 4.3.7]{dRStack}, we have $\ol{h}_\ell^!1\cong 1[1]$, and similarly for $m$, which clearly yields the desired isomorphism.

Overall, we conclude that $X_\infty/\cal{G}_0$ is suave. To check that the dualising sheaf is invertible, we use the formula from \cref{lem:sixfunctors-drfflem}. Namely, as all $f_{\infty m}$ are weakly cohomologically proper by \cref{lem:sixfunctors-xmgmtransitionproper}, the dualising sheaf is given by $\colim_m f_{\infty m}^*\SD(f_{\infty m, *} 1)$, where $\SD(-)$ denotes the suave dual over $\GSpec\Q_p$. Using \cref{lem:sixfunctors-qpsynkeylem} and the fact that $f_{\ell m}$ is weakly cohomologically proper by \cref{lem:sixfunctors-xmgmtransitionproper}, we obtain
\begin{equation*}
\SD(f_{\infty m, *} 1)\cong \SD(f_{\ell m, *} 1)\cong f_{\ell m, *}\SD(1)\cong f_{\ell m, *}f_{\ell m}^*\SD(1)\cong \SD(1)\tensor f_{\ell m, *}1\cong \SD(1)\tensor f_{\infty m, *}1
\end{equation*}
for each fixed $m$ and $\ell\gg m$, where we have used that the dualising sheaves on $X_\ell/\cal{G}_0^{\ell\dagger}$ and $X_m/\cal{G}_0^{m\dagger}$ are compatible under pullback by the argument above. Thus, we find that, on $X_\infty/\cal{G}_0$, we have
\begin{equation*}
\SD(1)\cong \colim_m f_{\infty m}^*(\SD(1)\tensor f_{\infty m, *}1)\cong \colim_m f_{\infty m}^*\SD(1)\tensor\colim_m f_{\infty m}^*f_{\infty m, *}1\;.
\end{equation*}
The first tensor factor is simply $f_{\infty m}^*\SD(1)$ for any $m$ by compatibility of the dualising sheaves on the $X_m/\cal{G}_0^{m\dagger}$ under pullback while the second tensor factor is the tensor unit: Indeed, as $X_\infty/\cal{G}_0\cong \lim_m X_m/\cal{G}_0^{m\dagger}$ in the category of kernels by \cref{lem:sixfunctors-drfflem} and the maps $f_{\ell m}$ are all weakly cohomologically proper, we have
\begin{equation*}
(f_{\infty m, *})_m: \D(X_\infty/\cal{G}_0)\xrightarrow{\cong} \lim\nolimits_* \D(X_m/\cal{G}_0^{m\dagger})\;,
\end{equation*}
where the transition functors in the limit are lower-$*$-functors. One readily checks that the left adjoint and hence the inverse to the functor $(f_{\infty m, *})_m$ is given by $\colim_m f_{\infty m}^*$ and this yields $\colim_m f_{\infty m}^*f_{\infty m, *}1\cong 1$, as desired. 

Overall, we conclude that $X_\infty/\cal{G}_0$ is cohomologically smooth and that the dualising sheaf is obtained by pulling back any of the dualising sheaves of the stacks $X_m/\cal{G}^{m\dagger}$. Tracing through the arguments above and counting shifts, one finds that these dualising sheaves are all concentrated in degree $-3$ and hence the dualising sheaf on $X_\infty/\cal{G}_0$ is concentrated in degree $-3$ as well.

To finish the proof, recall from \cref{thm:pres-qpn} that $\Q_{p, (p^{1/2}, p^{3/2})}^{\ol{\N}}\cong X_\infty/\cal{G}$ and observe that there is a canonical map $\cal{G}_0\rightarrow\cal{G}$ with cokernel $\G_m/\ol{\T}\cong (0, \infty)$, where we have used \cref{lem:defis-rhvariant}. Thus, the induced map $X_\infty/\cal{G}_0\rightarrow \Q_{p, (p^{1/2}, p^{3/2})}^{\ol{\N}}$ is a $(0, \infty)$-torsor and hence $\Q_{p, (p^{1/2}, p^{3/2})}^{\ol{\N}}$ is cohomologically smooth as this may be checked smooth-locally on the source. As the dualising sheaf of $(0, \infty)$ is concentrated in degree $-1$ and the dualising sheaf of $X_\infty/\cal{G}_0$ is concentrated in degree $-3$ by the previous paragraph, we conclude that the dualising sheaf of $\Q_{p, (p^{1/2}, p^{3/2})}^{\ol{\N}}$ is concentrated in degree $-2$, as desired.
\end{proof}

As already stated in the beginning, having \cref{prop:sixfunctors-qpsynkeypart} now enables us to prove \cref{thm:sixfunctors-qpsynsmooth} without too much difficulty:

\begin{proof}[Proof of \cref{thm:sixfunctors-qpsynsmooth}]
We first check that $\Q_p^\Syn$ is cohomologically smooth with dualising sheaf concentrated in degree $-3$. For this, we may work locally on an open cover of $\Q_p^\Syn$ as cohomological smoothness is smooth-local on the source. We encourage the reader to consult \cref{fig:xnviatu} and \cref{fig:xnviakappa} while following the argument.

We first show cohomological smoothness of 
\begin{equation}
\label{eq:sixfunctors-qpsynkeyopenpart}
\Q_{p, (p^{1/2}, p^{3/2}), |u|<1, |t|<1}^\N\cong\Q_{p, (p^{1/2}, p^{3/2})}^\N\setminus (j_\dR(\Q_{p, (p^{1/2}, p^{3/2})}^\prism)\cup j_\HT(\Q_{p, (p^{-1/2}, p^{1/2})}^\prism))\;,
\end{equation}
which we note is an open subspace of $\Q_p^\Syn$. For this, we observe that (\ref{eq:sixfunctors-qpsynkeyopenpart}) is open inside $\Q_{p, (p^{1/2}, p^{3/2})}^{\N, \mathrm{ext}}$, which is a $(0, \infty)$-torsor over $\Q_{p, (p^{1/2}, p^{3/2})}^{\ol{\N}}$ by \cref{lem:pres-qpredntorsor} and hence cohomologically smooth by \cref{prop:sixfunctors-qpsynkeypart}. Moreover, by loc.\ cit., the dualising sheaf is concentrated in degree $-3$ as the dualising sheaf of $(0, \infty)$ is in degree $-1$. 

Now consider $\Q_{p, |t|\neq 0, |u|<1-\epsilon}^\N$ for any $\epsilon>0$ and note that it is isomorphic to
\begin{equation*}
\{|\widetilde{\mu}|<|t|(1-\epsilon)\}\subseteq \Q_p^\prism\times (0, 1]\;,
\end{equation*}
where $|t|$ is the coordinate on $(0, 1]$. Indeed, this follows from the definition of $\Q_p^\N$ using the isomorphism
\begin{equation*}
\begin{split}
\ol{\DD}^\times/\ol{\T}\times (\overcirc{\DD}(1-\epsilon)/\ol{\T})^\dR&\cong (\{|r|<|t|(1-\epsilon)\}\subseteq (0, 1]\times (\overcirc{\DD}(1-\epsilon)/\ol{\T})^\dR) \\
(t, u)&\mapsto (t, r)\coloneqq (t, ut)\;,
\end{split}
\end{equation*}
where we have used \cref{lem:defis-rhvariant}. Similarly, using the argument from the proof of \cref{lem:defis-jht}, we have 
\begin{equation*}
\Q_{p, |u|\neq 0, |t|<1-\epsilon}^\N\cong (\{|\widetilde{\mu}|<|u|(1-\epsilon)\}\subseteq \Q_p^\prism\times (0, 1])\;,
\end{equation*}
where now $|u|$ is the coordinate on $(0, 1]$, induced by the isomorphism
\begin{equation*}
\begin{split}
\overcirc{\DD}(1-\epsilon)/\ol{\T}\times (\ol{\DD}^\times/\ol{\T})^\dR&\cong (\{|r|<|u|(1-\epsilon)\}\subseteq (\overcirc{\DD}(1-\epsilon)/\ol{\T})^\dR)\times (0, 1] \\
(t, u)&\mapsto (r, u)\coloneqq (ut, u)\;,
\end{split}
\end{equation*}
where we use that $(\ol{\DD}^\times/\ol{\T})^\dR\cong \ol{\DD}^\times/\ol{\T}$ by \cref{lem:defis-rhvariant} to make the map well-defined. In particular, we conclude that
\begin{equation*}
\begin{split}
\Q_{p, |t|>1/2, |u|<1/2}^\N&\cong (\{|\widetilde{\mu}|<|t|/2\}\subseteq \Q_p^\prism\times (1/2, 1]) \\
\Q_{p, |u|>1/2, |t|<1/2}^\N&\cong (\{|\widetilde{\mu}|<|u|/2\}\subseteq \Q_p^\prism\times (1/2, 1])
\end{split}
\end{equation*}
and gluing these two along the isomorphism
\begin{equation*}
\Q_{p, |t|=1, |u|<1/2}^\N\cong \Q_{p, |\widetilde{\mu}|<1/2}^\prism\cong \Q_{p, |u|=1, |t|<1/2}^\N
\end{equation*}
yields an open substack of $\Q_p^\Syn$ isomorphic to $\Q_{p, |\widetilde{\mu}|<1/2}^\prism\times (1/2, 3/2)$. As $\Q_p^\prism$ is cohomologically smooth with dualising sheaf concentrated in degree $-2$ by \cref{cor:prism-qpprismsmooth}, we conclude that this is cohomologically smooth as well, with dualising sheaf concentrated in degree $-3$.

The open substacks we have considered so far cover open neighbourhoods of $\{|ut|=0\}$ and $\{|u|=1\}\cup\{|t|=1\}$ inside $\Q_p^\N$, respectively. Thus, it suffices to consider $\Q_{p, 0<|u|<1, 0<|t|<1}^\N$ for some $\epsilon>0$, which we note is an open substack of $\Q_p^\Syn$ as the gluing only happens at $\{|u|=1\}$ and $\{|t|=1\}$. However, by \cref{prop:defis-utneq0}, we have
\begin{equation*}
\Q_{p, 0<|u|<1, 0<|t|<1}^\N\cong (\Q_p^\prism\setminus \Q_p^\dR)\times (0, 1)\;,
\end{equation*}
which is cohomologically smooth as $\Q_p^\prism\setminus \Q_p^\dR$ is open inside $\Q_p^\prism$. Again, the dualising sheaf is concentrated in degree $-3$ since the dualising sheaf of $\Q_p^\prism$ is concentrated in degree $-2$ by \cref{cor:prism-qpprismsmooth}.

Overall, we conclude that $\Q_p^\Syn$ is cohomologically smooth and that its dualising sheaf is a line bundle $L$ shifted in cohomological degree $-3$. By \cref{thm:drbundles-main} below, to identify the line bundle, we just have to compute the corresponding line bundle on $\Q_p^{\Div^1}$ under the realisation functor $T_{\Div^1}$ defined in §\ref{sect:proet} below. This, in turn, amounts to computing the pullback of $L$ to $\Q_{p, (0, p)}^\prism$ together with its transformation behaviour under Frobenius. For this, consider the composite map
\begin{equation*}
\Q_{p, (0, p)}^\prism\xrightarrow{\id\times \{1/2\}} \Q_{p, (0, p)}^\prism\times (0, 1)\cong \Q_{p, (0, p), |u|<1, |t|<1}^\N\xrightarrow{j} \Q_p^\Syn\;,
\end{equation*}
where the last map is an open embedding and the isomorphism comes from \cref{prop:defis-utneq0}. Letting $f: \Q_p^\Syn\rightarrow\GSpec\Q_p$ denote the structure map, we have
\begin{equation*}
\O\{1\}[2]\cong (\id\times \{1/2\})^!j^!f^!1\cong (\id\times\{1/2\})^!j^*L[3]
\end{equation*}
by \cref{cor:prism-qpprismsmooth}. To compute $(\id\times\{1/2\})^!j^*L$, let $j': (0, 1/2)\cup (1/2, 1)\rightarrow (0, 1)$ denote the canonical inclusion and recall that $\pi: \Q_{p, (0, p)}^\prism\times (0, 1)\rightarrow \Q_{p, (0, p)}^\prism$ is the first projection; moreover, observe that $j^*L\cong \pi^*L'$ for some line bundle $L'$ on $\Q_{p, (0, p)}^\prism$ by \cref{cor:perf-contractible}. Thus, we obtain
\begin{equation*}
\begin{split}
(\id\times\{1/2\})^!j^*L&\cong \pi_*(\id\times\{1/2\})_*(\id\times\{1/2\})^!\pi^*L'\cong \pi_*\fib(\pi^*L'\rightarrow j'_*{j'}^*\pi^*L') \\
&\cong \fib(L'\rightarrow L'\oplus L')\cong L'[-1]
\end{split}
\end{equation*}
by excision and \cref{cor:perf-contractible}, from which we overall conclude that $L'\cong \O\{1\}$. Thus, to conclude that $L\cong \O\{1\}$, we just have to check that the pullback of $L$ to $\Q_{p, (0, p)}^\prism$ has the correct transformation behaviour under Frobenius, but this follows from \cref{thm:prism-div1qphkcohsmooth}.
\end{proof}

\subsection{Preservation of smooth and proper maps}

Let us move on to the relative case, i.e.\ to studying $X^\Syn\rightarrow Y^\Syn$ for a map $X\rightarrow Y$ of derived Berkovich spaces. The most difficult part is to establish that $X^\Syn\rightarrow Y^\Syn$ is cohomologically smooth whenever $X\rightarrow Y$ is rigid smooth. By compatibility of $X\mapsto X^\N$ with finite étale maps and open localisations, see \cref{prop:defis-etmaps} and \cref{prop:defis-openloc}, this will quickly reduce to the case $X=\G_m$ and $Y=\GSpec\Q_p$. Moreover, similarly as in the proof of \cref{thm:sixfunctors-qpsynsmooth}, the key part will again be to understand the locus $\G_{m, (p^{1/2}, p^{/3/2})}^\N$. Thus, our first goal is to prove:

\begin{prop}
\label{prop:sixfunctors-smoothmapskeypart}
The map $\G_{m, (p^{1/2}, p^{/3/2})}^\N\rightarrow \Q_{p, (p^{1/2}, p^{3/2})}^\N$ is cohomologically smooth with dualising sheaf concentrated in degree $-2$.
\end{prop}

The method we will use is very similar to the method of the proof of \cref{prop:sixfunctors-qpsynkeypart}, and we again start by establishing some notation: Let $Y_\infty\coloneqq \lim_{x\mapsto x^p} \G_m$ and write $Y_m$ for the $m$-th stage of this limit. For any integer $h\geq 1$, we write
\begin{equation*}
\Z_p^{h\dagger}\coloneqq \bigsqcup_{\theta\in \Z_p/p^h\Z_p} (\theta+\ol{\DD}(1/p^h))
\end{equation*}
and observe that this canonically admits the structure of an additive group stack over $\GSpec\Q_p$. Put
\begin{equation*}
\widetilde{\cal{G}}_{1, \red}^{h\dagger}\coloneqq \Z_p^{h\dagger}\times\G_a^\dagger
\end{equation*}
and note that, over $\A^1$, this group receives two maps from $\G_a^\dagger$: The first one is given by $\theta\mapsto (0, u\theta)$, where $u$ denotes the coordinate on $\A^1$, while the second one is given by $\theta\mapsto (\theta, 0)$. We define $\cal{G}_{1, \red}^{h\dagger}$ to be the coequaliser of these two maps in the category of group stacks over $\A^1$. For convenience, we write $\cal{G}_{1, \red}^{\infty\dagger}\coloneqq \cal{G}_{1, \red}\coloneqq\Z_p^\la\coprod^\AbGrp_{\G_a^\dagger} \G_a^\dagger$ over $\A^1$, where the pushout is along the multiplication-by-$u$-map on $\G_a^\dagger$, as before. Observe that the canonical maps $\Z_p^{h'\dagger}\rightarrow\Z_p^{h\dagger}$ for $h'>h$ induce maps $\cal{G}_{1, \red}^{h'\dagger}\rightarrow\cal{G}_{1, \red}^{h\dagger}$. Finally, note that, over $\Q_{p, (p^{1/2}, p^{3/2})}^{\ol{\N}}$, the action of $\Z_p^\la\coprod^\AbGrp_{\G_a^\dagger} \G_a^\dagger$ on $Y_\infty$ canonically extends to an action of $\cal{G}_{1, \red}^{h\dagger}$ on $Y_m$ for all $h\geq m$ as the formula from loc.\ cit.\ still makes sense by the argument from the discussion preceding \cref{lem:sixfunctors-limitghdagger}.

Working over the cover $\widetilde{\Q_{p\mathrlap{, (p^{1/2}, p^{3/2})}}^{\ol{\N}}}\hphantom{\scriptstyle{(p^{1/2}, p^{3/2})}}$ of $\Q_{p, (p^{1/2}, p^{3/2})}^{\ol{\N}}$ from \cref{thm:pres-qpn} throughout, our aim will be to apply \cref{lem:sixfunctors-drfflem} to the tower
\begin{equation*}
\begin{split}
Y_\infty/\cal{G}_{1, \red}\rightarrow\dots\rightarrow Y_{m+1}/\cal{G}_{1, \red}^{(m+1)\dagger}\rightarrow Y_m/\cal{G}_{1, \red}^{m\dagger}\rightarrow\dots\rightarrow Y_1/\cal{G}_{1, \red}^{1\dagger}\;,
\end{split}
\end{equation*}
the transition maps in which we denote by $f_{\ell m}: Y_\ell/\cal{G}_{1, \red}^{\ell\dagger}\rightarrow Y_m/\cal{G}_{1, \red}^{m\dagger}$, and the $!$-cover of $Y_\infty/\cal{G}_{1, \red}$ by $Y_\infty$. Let us record the following features of the situation:

\begin{lem}
\label{lem:sixfunctors-smoothmapsproperties}
In the category of kernels over $\widetilde{\Q_{p\mathrlap{, (p^{1/2}, p^{3/2})}}^{\ol{\N}}}\hphantom{\scriptstyle{(p^{1/2}, p^{3/2})}}$, we have
\begin{equation*}
Y_\infty\times_{Y_\infty/\cal{G}_{1, \red}} Y_\infty\cong \lim_m Y_\infty\times_{Y_m/\cal{G}_{1, \red}^{h\dagger}} Y_\infty\;.
\end{equation*}
Moreover, for all $m\leq\ell\leq\infty$, the map $f_{\ell m}: Y_\ell/\cal{G}_{1, \red}^{\ell\dagger}\rightarrow Y_m/\cal{G}_{1, \red}^{m\dagger}$ is weakly cohomologically proper.
\end{lem}
\begin{proof}
As in \cref{cor:sixfunctors-qpsmoothcond(c)} and \cref{lem:sixfunctors-xmgmtransitionproper}.
\end{proof}

The next lemma should be compared to \cref{lem:sixfunctors-qpsynkeylem} and is again the key to proving \cref{prop:sixfunctors-smoothmapskeypart}:

\begin{lem}
\label{lem:sixfunctors-smoothmapskeylem}
For any fixed $m\geq 5$, the natural map $f_{\ell m, *}1\rightarrow f_{\infty m, *}1$ is an isomorphism for $\ell\gg m$.
\end{lem}
\begin{proof}
By \cref{lem:sixfunctors-smoothmapsproperties}, we know that the $f_{\ell m}$ are prim and hence they satisfy base change by \cite[Lem.\ 4.5.13]{HeyerMann}. Thus, as in the proof of \cref{lem:sixfunctors-qpsynkeylem}, we may replace the maps $f_{\ell m}$ by the maps
\begin{equation*}
\widetilde{f}_{\ell m}: Y_\ell/\Z_p^{\ell\dagger}\rightarrow Y_m/\Z_p^{m\dagger}\;.
\end{equation*}
Moreover, we note that everything still makes sense over the base $\ol{X}_\infty$ from loc.\ cit.\ in place of $\widetilde{\Q_{p\mathrlap{, (p^{1/2}, p^{3/2})}}^{\ol{\N}}}\hphantom{\scriptstyle{(p^{1/2}, p^{3/2})}}$ and, by base change, it thus suffices to prove the claim over $\ol{X}_\infty$. By the same argument, we may then restrict to $\ol{X}_{\infty, N}\coloneqq \lim_m \ol{X}_{m, N}$ for varying $N\gg 0$.

Finally, we may also localise on $Y_1/\Z_p^{1\dagger}$ and restrict along
\begin{equation*}
Y_{1, k}\coloneqq\GSpec\Q_p\langle p^k x^{\pm 1}\rangle\rightarrow \G_m\cong Y_1\;.
\end{equation*}
For any $\ell$, we let $Y_{\ell, k}$ denote the base change of $Y_\ell$ along this map. We write $\widetilde{f}_{\ell m, k}: Y_{\ell, k}/\Z_p^{\ell\dagger}\rightarrow Y_{m, k}/\Z_p^{m\dagger}$ for the induced transition maps.

Then we note that each $Y_{\ell, k}$ is affine over the affine base $\ol{X}_{\infty, N}$ and thus $\widetilde{f}_{\ell m, k, *}1\cong \O(Y_{\ell, k})^{m\dagger\text{-an}}$ by \cite[Cor.\ 2.50]{Mikami} and, similarly, $\widetilde{f}_{\infty m, k, *}1\cong \O(Y_{\infty, k})^{m\dagger\text{-an}}$. Thus, we have to prove that the natural map
\begin{equation*}
\O(Y_{\ell, k})^{m\dagger\text{-an}}\rightarrow \O(Y_{\infty, k})^{m\dagger\text{-an}}
\end{equation*}
is an isomorphism for $\ell\gg 0$. However, this follows from \cite[Cor.\ 3.7]{Mikami} as the family of Banach algebras $\{\O(Y_{\ell, k})\}_\ell$ satisfies the Tate--Sen axioms (C1)--(C4) from \cite[§3]{Mikami} by the following lemma.
\end{proof}

\begin{lem}
For each $k\geq 1$, the family of Banach algebras $\{\O(Y_{\ell, k})\}_\ell$ from the proof of \cref{lem:sixfunctors-smoothmapskeylem} satisfies the Tate--Sen axioms (C1)--(C4) from \cite[§3]{Mikami}.
\end{lem}
\begin{proof}
Let $A_\ell\coloneqq \O(Y_{\ell, k})$, note that
\begin{equation*}
A_\ell=\Q_p\langle p^k x^{\pm 1}, p^{k/p}x^{\pm 1/p}, \dots, p^{k/p^\ell}x^{\pm 1/p^\ell}\rangle
\end{equation*}
and let $A$ denote the uniform completion of $\colim_\ell A_\ell$. Then we have trace maps $R_\ell: A\rightarrow A_\ell$ given by $x^\lambda\mapsto x^\lambda$ if $|\lambda|_p\leq p^\ell$ and $x^\lambda\mapsto 0$ otherwise. Clearly, we have $|R_\ell(x)|\leq |x|$ and $\lim_{\ell\rightarrow\infty} R_\ell(x)=x$ for all $x\in A$ while $R_\ell(x)=x$ whenever $x\in A_\ell$, which already yields (C1). Moreover, we have already discussed above that the action of $\Z_p$ on $Y_\ell$ is $\ell\dagger$-analytic and, in particular, locally analytic.

Now note that, for all $\theta\in\Z_p$ with $v_p(\theta)<k-3$, we have
\begin{equation}
\label{eq:sixfunctors-sentateaxioms}
p^{-p^{v_p(\theta)-k+3}}\leq |q^{\theta/p^k}-1|\leq p^{-p^{v_p(\theta)-k}}
\end{equation}
on $\ol{X}_\infty$. In particular, this implies that, whenever $x\in\Ker R_\ell$, we have
\begin{equation*}
|(\theta-\id)(x)|\geq p^{-p^{v_p(\theta)-\ell+3}}|x|\;,
\end{equation*}
which yields (C2). Moreover, we observe that our work so far also shows that the algebras $\{A_\ell\}_\ell$ satisfy the Colmez--Sen--Tate axioms from \cite[Def.\ 2.2.1]{GeometricSenTheory} and hence \cite[Lem.\ 2.4.3]{GeometricSenTheory} shows that (C4) is satisfied as well. Finally, we note that the upper bound in (\ref{eq:sixfunctors-sentateaxioms}) is true for all $\theta\in\Z_p$ and $k\geq 0$, which lets us conclude that
\begin{equation*}
|(\theta-\id)(x)|\leq p^{-1/p^\ell}|x|
\end{equation*}
for all $x\in A_\ell$. In particular, we obtain (C3) and this finishes the proof.
\end{proof}

Finally, we can prove \cref{prop:sixfunctors-smoothmapskeypart}:

\begin{proof}[Proof of \cref{prop:sixfunctors-smoothmapskeypart}]
We apply \cref{lem:sixfunctors-drfflem} to the tower 
\begin{equation*}
Y_\infty/\cal{G}_{1, \red}\rightarrow\dots\rightarrow Y_{m+1}/\cal{G}_{1, \red}^{(m+1)\dagger}\rightarrow Y_m/\cal{G}_{1, \red}^{m\dagger}\rightarrow\dots\rightarrow Y_1/\cal{G}_{1, \red}^{1\dagger}
\end{equation*}
of Gelfand stacks over $S=\widetilde{\Q_{p\mathrlap{, (p^{1/2}, p^{3/2})}}^{\ol{\N}}}\hphantom{\scriptstyle{(p^{1/2}, p^{3/2})}}$ and the $!$-cover $Y_\infty\rightarrow Y_\infty/\cal{G}_{1, \red}$. Indeed, conditions (a), (b) and (c) of loc.\ cit.\ are ensured by \cref{lem:sixfunctors-smoothmapsproperties} and \cref{lem:sixfunctors-smoothmapskeylem}, where we note that $\cal{G}_{1, \red}$ is prim as it is a $\G_a^\dagger$-torsor over the affine Gelfand stack $\Z_p^\sm$.

As in the proof of \cref{prop:sixfunctors-qpsynkeypart}, the fact that the $f_{\ell m}$ are weakly cohomologically proper and $f_{\infty m, *}1\cong f_{\ell m, *} 1$ for each fixed $m$ and $\ell\gg m$ reduces the proof that $Y_\infty/\cal{G}_{1, \red}$ is cohomologically smooth with dualising sheaf concentrated in degree $-2$ to showing that each $Y_\ell/\cal{G}_{1, \red}^{\ell\dagger}$ is cohomologically smooth with dualising sheaf concentrated in degree $-2$ and that all these dualising sheaves are compatible under pullback along the transition maps $f_{\ell m}$.

For this, we first note that $X_\ell\cong \G_m$ for each $\ell$ and that the maps $X_\ell\rightarrow X_m$ are finite étale, hence cohomologically étale. Thus, we only have to worry about the classifying stacks $*/\cal{G}_{1, \red}^{\ell\dagger}$, where we note the factorisation
\begin{equation*}
*/\cal{G}_{1, \red}^{\ell\dagger}\rightarrow */(\Z_p^{\ell\dagger})^\dR\rightarrow *\;.
\end{equation*}
While the first map is a $\G_a^\dagger$-gerbe and hence cohomologically smooth by \cite[Thm.\ 4.3.13]{dRStack}, the second map is cohomologically smooth as it factors as
\begin{equation*}
*/(\Z_p^{\ell\dagger})^\dR\rightarrow */\ol{\DD}(1/p^\ell)\rightarrow */\G_a^\dR\rightarrow *\;,
\end{equation*}
where the first map is a gerbe for a finite group, hence cohomologically étale, the second map is a $\G_a^\dR/\ol{\DD}(1/p^\ell)^\dR\cong \G_a/\ol{\DD}(1/p^\ell)$-torsor, hence cohomologically smooth by \cite[Prop.\ 4.3.7]{dRStack}, and the last map is the structure map of the classifying stack of a cohomologically smooth group by \cite[Thm.\ 5.1.17]{dRFF}, hence cohomologically smooth as well. Moreover, the compatibility of the dualising sheaves under pullback follows by the same argument as in the proof of \cref{prop:sixfunctors-qpsynkeypart}. 

Counting the shifts in the argument from the previous paragraph, we conclude that $Y_\infty/\cal{G}_{1, \red}$ is cohomologically smooth with dualising sheaf concentrated in degree $-2$. To finish the proof, recall from \cref{thm:pres-tn} that $\G_{m, (p^{1/2}, p^{3/2})}^{\ol{\N}}\cong Y_\infty/\cal{G}_{1, \red}$ over $\widetilde{\Q_{p\mathrlap{, (p^{1/2}, p^{3/2})}}^{\ol{\N}}}\hphantom{\scriptstyle{(p^{1/2}, p^{3/2})}}$. As cohomological smoothness is local on the target, this shows that $\G_{m, (p^{1/2}, p^{3/2})}^{\ol{\N}}$ is cohomologically smooth over $\Q_{p, (p^{1/2}, p^{3/2})}^{\ol{\N}}$ and then cohomological smoothness of $\G_{m, (p^{1/2}, p^{3/2})}^\N$ over $\Q_{p, (p^{1/2}, p^{3/2})}^\N$ follows by base change due to 
\begin{equation*}
\G_{m, (p^{1/2}, p^{3/2})}^\N\cong \G_{m, (p^{1/2}, p^{3/2})}^{\ol{\N}}\times_{\Q_{p, (p^{1/2}, p^{3/2})}^{\ol{\N}}} \Q_{p, (p^{1/2}, p^{3/2})}^\N\;. \qedhere
\end{equation*}
\end{proof}

Having \cref{prop:sixfunctors-smoothmapskeypart}, we can now easily deduce the following general cohomological smoothness result for maps between syntomifications:

\begin{thm}
\label{thm:sixfunctors-smoothmaps}
Let $f: X\rightarrow Y$ be a rigid smooth map of derived Berkovich spaces. Then the induced map $f^\Syn: X^\Syn\rightarrow Y^\Syn$ is cohomologically smooth. If $X$ is moreover a smooth partially proper rigid space over $\Q_p$ and $f$ is pure of relative dimension $d$, then the dualising sheaf is given by $\O\{d\}[2d]$.
\end{thm}

Note that, due to $X\rightarrow Y$ being rigid smooth, the assumption that $X$ is a smooth partially proper rigid space over $\Q_p$ is in particular satisfied if $Y$ is a smooth partially proper rigid space.

\begin{proof}[Proof of \cref{thm:sixfunctors-smoothmaps}]
As cohomological smoothness is local on the target, it suffices to check the claim for $f^\N: X^\N\rightarrow Y^\N$. Then recall from \cref{prop:defis-etmaps} and \cref{prop:defis-openloc} that the functor $X\mapsto X^\N$ preserves finite étale maps and open localisations. As finite étale maps and open embeddings are cohomologically étale and cohomological smoothness may be checked locally on the target and smooth-locally on the source, we may reduce to the case $X=\A^d_Y$. By commutation of $X\mapsto X^\N$ with limits and stability of cohomological smoothness under base change and composition, we may then further reduce to $Y=\GSpec\Q_p$ and $X=\A^1$. Finally, we perform a further open localisation on the source by covering $\A^1$ by $\G_m$ and $1+\G_m$, and this allows us to reduce to $X=\G_m$.

However, now note that $\G_{m, (p^{1/2}, p^{3/2})}^\N\rightarrow \Q_{p, (p^{1/2}, p^{3/2})}^\N$ is cohomologically smooth by \cref{prop:sixfunctors-smoothmapskeypart}. Moreover, we have
\begin{equation*}
\G_{m, |ut|\neq 0}^\N\cong (\G_m^\prism\setminus \G_m^\dR)\times [0, 1]
\end{equation*}
by \cref{prop:defis-utneq0} and similarly for $\Q_{p, |ut|\neq 0}^\N$, whence the map $\G_{m, |ut|\neq 0}^\N\rightarrow \Q_{p, |ut|\neq 0}^\N$ is cohomologically smooth by \cref{cor:prism-sixfunctors}. As $\Q_{p, (p^{1/2}, p^{3/2})}^\N$ and $\Q_{p, |ut|\neq 0}^\N$ form an open cover of $\Q_p^\N$, this proves the claim by locality of cohomological smoothness on the target.

Overall, we conclude that $f^\Syn: X^\Syn\rightarrow Y^\Syn$ is cohomologically smooth and it remains to compute the dualising sheaf under the assumption that $X$ is a smooth partially proper rigid space over $\Q_p$. Going through the argument above, we already see that the dualising sheaf is a line bundle concentrated in cohomological degree $-2d$. As at the end of the proof of \cref{thm:sixfunctors-qpsynsmooth}, we can use \cref{thm:drbundles-main} below to conclude that it suffices to identify the pullback of this line bundle to $X^\prism_{(0, p)}$ together with the automorphism induced by Frobenius. However, using smooth base change, see \cite[Lem.\ 4.5.13]{HeyerMann}, this reduces to the computation of the dualising sheaf of $f^{\Div^1}: X^{\Div^1}\rightarrow Y^{\Div^1}$, and we already know that this is $\O\{d\}[2d]$ by \cref{thm:prism-sixfunctorsdiv1hk}. This finishes the proof.
\end{proof}

\begin{cor}
Let $X$ be a smooth partially proper rigid space over $\Q_p$ which is pure of dimension $d$ and let $f: X^\Syn\rightarrow\GSpec\Q_p$ be the structure map. For any perfect analytic $F$-gauge $E\in\D(X^\Syn)$, there is an isomorphism $f_*(E^\vee\{d+1\})[2d+3]\cong (f_!E)^\vee$. In particular, we have
\begin{equation*}
R\Gamma_\Syn(X, \Q_p(i))\cong R\Gamma_{\Syn, c}(X, \Q_p(2d+1-i))^\vee[-2d-3]
\end{equation*}
for any $i\in\Z$.
\end{cor}
\begin{proof}
As in \cref{cor:sixfunctors-qpsynsmooth}.
\end{proof}

\begin{rem}
Let us note that any smooth map between partially proper rigid spaces is rigid smooth and hence satisfies the assumption of \cref{thm:sixfunctors-smoothmaps}.
\end{rem}

Finally, we would also like to establish a properness statement for maps between syntomifications. For this, we first make the following definition:

\begin{defi}
Let $X\rightarrow Y$ be a map between derived Berkovich spaces. It is called \emph{locally of finite presentation} if, after passing to a strict closed cover of the source and target by rational localisations, it is of the form $\GSpec B\rightarrow\GSpec A$ for $B$ a finitely presented animated $A$-algebra.
\end{defi}

\begin{rem}
We emphasise that the finite presentation assumption above is in the derived sense. This means that a morphism of classical rings $A\rightarrow B$ is finitely presented in this sense if and only if it factors as $A\rightarrow A[x_1, \dots, x_d]\rightarrow B$, where $A[x_1, \dots, x_d]\rightarrow B$ is the quotient by a regular sequence.
\end{rem}

Using this terminology, the properness statement we prove can be formulated as follows. We note that it is likely that the finite presentation assumption below can be weakened, but we have chosen to stick to this version as it is more than enough for our purposes.

\begin{thm}
\label{thm:sixfunctors-propermaps}
Let $f: X\rightarrow Y$ be a map between derived Berkovich spaces which is quasicompact and locally of finite presentation. Then the induced map $f^\Syn: X^\Syn\rightarrow Y^\Syn$ is weakly cohomologically proper.
\end{thm}

Before we begin the proof, let us establish the following easier variant, to which we will later reduce by a dévissage:

\begin{prop}
\label{prop:sixfunctors-berksmpropermaps}
Let $f: X\rightarrow Y$ be a Berkovich smooth map between derived Berkovich spaces which is quasicompact. Then the induced map $f^\Syn: X^\Syn\rightarrow Y^\Syn$ is weakly cohomologically proper.
\end{prop}
\begin{proof}
We first check that $f^\Syn$ is prim. This may be done locally on the target and hence reduces to checking primness of $f^\N: X^\N\rightarrow Y^\N$. Now recall from \cref{prop:defis-etmaps} and \cref{prop:defis-openloc} that $X\mapsto X^\N$ preserves finite étale maps and rational localisations and note that any finite étale map or rational localisation is cohomologically proper; in particular, the functor $(-)^\N$ sends Berkovich étale maps to cohomologically proper maps. As the map $f: X\rightarrow Y$ is quasicompact and primness is local on the target as well as local for finite closed covers of the source, we may assume that $f$ factors as $X\rightarrow\A^d_Y\rightarrow Y$ with the first map being Berkovich étale, hence inducing a cohomologically proper map Nygaardifications. By quasicompactness of $f$, the map $X\rightarrow\A^d_Y$ factors through a closed polydisk, which we may assume to be $\ol{\DD}^d_Y$ after rescaling. We are thus reduced to $X=\ol{\DD}^d_Y$ and, by compatibility of $X\mapsto X^\N$ with limits and stability of primness under base change, we may further reduce to $X=\ol{\DD}$ and $Y=\GSpec\Q_p$. Finally, we may embed $\ol{\DD}$ into $\ol{\T}$ as the closed subspace $1+\ol{\DD}(1/p)$ and hence we may assume $X=\ol{\T}$.

Let us now study the case $Y=\GSpec\Q_p$ and $X=\ol{\T}$. By \cref{thm:pres-tn}, we have
\begin{equation*}
\ol{\T}^{\ol{\N}}\cong \lim_{x\mapsto x^p} \ol{\T}\;\Big/\;\Z_p^\la\coprod_{\G_a^\dagger}^\AbGrp \G_a^\dagger
\end{equation*}
after base changing to a cover of $\Q_{p, (p^{1/2}, p^{3/2})}^{\ol{\N}}$. Now note that $\lim_{x\mapsto x^p} \ol{\T}$ is affine, hence weakly cohomologically proper and the map
\begin{equation*}
*\;\Big/\;\Z_p^\la\coprod_{\G_a^\dagger}^\AbGrp \G_a^\dagger\rightarrow *
\end{equation*}
is prim as well. Indeed, this factors as
\begin{equation*}
*\;\Big/\;\Z_p^\la\coprod_{\G_a^\dagger}^\AbGrp \G_a^\dagger\rightarrow */\Z_p^\sm\rightarrow *\;,
\end{equation*}
where the first map is a $\G_a^\dagger$-gerbe, hence prim by \cite[Prop.\ 4.2.5]{dRStack}, and the second map is prim as $*\rightarrow */\Z_p^\sm$ is a weakly cohomologically proper $!$-cover such that the relative cohomology along this map is descendable in $\D(*/\Z_p^\sm)$ by \cite[Prop.\ 6.2.3]{dRStack}. Thus, we conclude that $\ol{\T}^{\ol{\N}}_{(p^{1/2}, p^{3/2})}\rightarrow \Q_{p, (p^{1/2}, p^{3/2})}^{\ol{\N}}$ is prim, whence $\ol{\T}^\N_{(p^{1/2}, p^{3/2})}\rightarrow \Q_{p, (p^{1/2}, p^{3/2})}^\N$ is prim by base change. 

As primness is local on the target, it thus suffices to show that $\ol{\T}^\N_{|ut|\neq 0}\rightarrow\Q_{p, |ut|\neq 0}^\N$ is prim as well. In this case, recall from \cref{prop:defis-utneq0} that 
\begin{equation*}
\ol{\T}^\N_{|ut|\neq 0}\cong (\ol{\T}^\prism\setminus \ol{\T}^\dR)\times [0, 1]
\end{equation*}
and similarly for $\Q_{p, |ut|\neq 0}^\N$, hence it suffices to show that $\ol{\T}^\prism\rightarrow\Q_p^\prism$ is prim. By \cref{prop:defis-prismffdr} and \cref{prop:defis-xprismxdiv1}, this in turn reduces to primness of $\ol{\T}^\HK\rightarrow\Q_p^\HK$ and $\ol{\T}^{\Div^1}\rightarrow\Q_p^{\Div^1}$. However, this can be seen from the presentations 
\begin{equation*}
\ol{\T}^\HK\cong \lim_{x\mapsto x^p} \ol{\T}^\dR\;\Big/\;\Z_p^\sm\;, \hspace{0.5cm} \ol{\T}^{\Div^1}\cong \lim_{x\mapsto x^p} \ol{\T}\;\Big/\;\Z_p^\la
\end{equation*}
from \cref{prop:prism-div1presentations} and \cref{cor:prism-hkpresentations}: Indeed, as 
\begin{equation*}
\lim_{x\mapsto x^p} \ol{\T}^\dR\cong \lim_{x\mapsto x^p} \ol{\T}\,\Big/\,\G_m^\dagger\;,
\end{equation*}
and $\lim_{x\mapsto x^p} \ol{\T}$ is affine, the claim reduces the primness of $*/\Z_p^\la, */\Z_p^\sm$ and $*/\G_m^\dagger$. However, $*/\Z_p^\la$ is a $\G_a^\dagger$-gerbe over $*/\Z_p^\sm$ and $*/\G_m^\dagger\cong */\G_a^\dagger$ via the logarithm. As $*/\G_a^\dagger$ is prim by \cite[Prop.\ 4.2.5]{dRStack} and we have already argued above that $*/\Z_p^\sm$ is prim, this proves the claim.

We conclude that $f^\Syn$ is prim and thus it remains to check that the diagonal $X^\Syn\rightarrow X^\Syn\times_{Y^\Syn} X^\Syn\cong (X\times_Y X)^\Syn$ of $f^\Syn$ is weakly cohomologically proper. We will do this by showing that, $!$-locally on the target, this factors into maps which are maps between affine Gelfand stacks $!$-locally on the target; this suffices by the proof of \cite[Lem.\ 3.1.21]{dRStack}. However, then the same argument as above again reduces us to the case where $X$ admits a Berkovich étale map to a closed polydisk $\ol{\DD}^d_Y$ over $Y$, and then by base change to $Y=\GSpec\Q_p$. To move on, note the factorisation
\begin{equation*}
X^\Syn\rightarrow X^\Syn\times_{(\ol{\DD}^d)^\Syn} X^\Syn\rightarrow X^\Syn\times_{\Q_p^\Syn} X^\Syn
\end{equation*}
of the diagonal map of $X^\Syn\rightarrow\Q_p^\Syn$. As the first map is the diagonal of a map which is Berkovich étale $!$-locally on the target, it is an isomorphism and thus we are reduced to studying the second map, which we observe is a base change of the diagonal map of $(\ol{\DD}^d)^\Syn$. By base change and composition, we reduce to $d=1$ and can then again replace $\ol{\DD}$ by $\ol{\T}$ as above. Base changing along $\Q_p^\N\rightarrow\Q_p^\Syn$, we furthermore reduce to the corresponding question for Nygaardifications.

Localising on $\Q_p^\N$ as above, we may separately check the claim for the diagonals of the maps
\begin{equation*}
\ol{\T}^\HK\rightarrow\Q_p^\HK\;, \hspace{0.3cm} \ol{\T}^{\Div^1}\rightarrow \Q_p^{\Div^1}\;, \hspace{0.3cm}\text{and}\hspace{0.3cm} \ol{\T}_{(p^{1/2}, p^{3/2})}^{\ol{\N}}\rightarrow \Q_{p, (p^{1/2}, p^{3/2})}^{\ol{\N}}\;.
\end{equation*}
Base changing the presentations of the diagonals of the maps $\ol{\T}^{\Div^1}\rightarrow\Q_p^{\Div^1}$ and $\ol{\T}^\HK\rightarrow\Q_p^\HK$ obtained from \cref{prop:prism-div1presentations} and \cref{cor:prism-hkpresentations} to $\lim_{x\mapsto x^p} \ol{\T}\times\lim_{x\mapsto x^p} \ol{\T}$, we easily see that we get maps between affine Gelfand stacks. For the map $\ol{\T}^{\ol{\N}}_{(p^{1/2}, p^{3/2})}\rightarrow\Q_{p, (p^{1/2}, p^{3/2})}^{\ol{\N}}$, \cref{thm:pres-tn} shows that the diagonal is given by
\begin{equation*}
\lim_{x\mapsto x^p} \ol{\T}\;\Big/\;\Z_p^\la\coprod^\AbGrp_{\G_a^\dagger} \G_a^\dagger\rightarrow \lim_{x\mapsto x^p} \ol{\T}\;\Big/\;\Z_p^\la\coprod^\AbGrp_{\G_a^\dagger} \G_a^\dagger\times \lim_{x\mapsto x^p} \ol{\T}\;\Big/\;\Z_p^\la\coprod^\AbGrp_{\G_a^\dagger} \G_a^\dagger
\end{equation*}
and the base change of this to $\lim_{x\mapsto x^p} \ol{\T}\times\lim_{x\mapsto x^p} \ol{\T}$ is
\begin{equation*}
\lim_{x\mapsto x^p} \ol{\T}\times \left(\Z_p^\la\coprod^\AbGrp_{\G_a^\dagger} \G_a^\dagger\right)\rightarrow \lim_{x\mapsto x^p} \ol{\T}\times \lim_{x\mapsto x^p} \ol{\T}\;.
\end{equation*}
Observing that this map admits a factorisation as
\begin{equation*}
\lim_{x\mapsto x^p} \ol{\T}\times \left(\Z_p^\la\coprod^\AbGrp_{\G_a^\dagger} \G_a^\dagger\right)\rightarrow \lim_{x\mapsto x^p} \ol{\T}\times \Z_p^\sm\rightarrow \lim_{x\mapsto x^p} \ol{\T}\times \lim_{x\mapsto x^p} \ol{\T}\;,
\end{equation*}
where the first map is a $\G_a^\dagger$-torsor and the second map is a map between affine Gelfand stacks then finally finishes the proof.
\end{proof}

We are now ready to prove \cref{thm:sixfunctors-propermaps}.

\begin{proof}[Proof of \cref{thm:sixfunctors-propermaps}]
We first show that $X^\Syn\rightarrow Y^\Syn$ is prim. As primness is local on the target, we may instead consider the map $X^\N\rightarrow Y^\N$. Since primness is furthermore local for finite closed covers of the source, the quasicompactness of $X\rightarrow Y$ and \cref{prop:defis-openloc} imply that we can reduce to the case where $Y$ is affine and $f$ factors as $X\rightarrow \A^d_Y\rightarrow Y$, where the first map is a Zariski-closed immersion given by quotienting out finitely many equations due to our finite presentation assumption. As $f$ is quasicompact, the map $X\rightarrow\A^d_Y$ factors through some closed polydisk, which we may assume to be of radius $1$. Then $(\ol{\DD}^d_Y)^\N\rightarrow Y^\N$ is prim by \cref{prop:sixfunctors-berksmpropermaps} and hence we are reduced to $X\rightarrow Y$ being a Zariski-closed immersion of affine Gelfand stacks given by finitely many equations. By composition, we may reduce to the case of a single equation and then $X\rightarrow Y$ is pulled back from $\{0\}\rightarrow \P^1$ along the map $Y\rightarrow\A^1\rightarrow\P^1$ classifying the equation cutting out $X$. Thus, it suffices to show that $\Q_p^\N\rightarrow (\P^1)^\N$ is prim, but this follows using cancellation from the fact that the diagonal of $(\P^1)^\N\rightarrow\Q_p^\N$ is prim by \cref{prop:sixfunctors-berksmpropermaps}.

Having established that $X^\Syn\rightarrow Y^\Syn$ is prim, it suffices to show that the diagonal is weakly cohomologically proper, for which, as in the proof of \cref{prop:sixfunctors-berksmpropermaps}, we show that, $!$-locally on the target, the map $X^\Syn\rightarrow X^\Syn\times_{Y^\Syn} X^\Syn$ factors into maps which are maps between affine Gelfand stacks $!$-locally on the target. The same dévissage as in the previous paragraph then reduces us to Nygaardifications and the cases $X=\GSpec\Q_p$ and $Y=\P^1$ as well as $X=\ol{\DD}^d_Y$. In the latter case, we note that the claim follows from the proof of \cref{prop:sixfunctors-berksmpropermaps}. For the map $\{0\}\rightarrow\P^1$, we first note that its diagonal agrees with the diagonal of $\{1\}\rightarrow\ol{\T}$, which, as in the proof of \cref{prop:sixfunctors-berksmpropermaps}, finally reduces us to considering the diagonals of
\begin{equation*}
\Q_p^\HK\rightarrow\ol{\T}^\HK\;, \hspace{0.3cm} \Q_p^{\Div^1}\rightarrow \ol{\T}^{\Div^1}\;, \hspace{0.3cm}\text{and}\hspace{0.3cm} \Q_{p, (p^{1/2}, p^{3/2})}^{\ol{\N}}\rightarrow \ol{\T}_{(p^{1/2}, p^{3/2})}^{\ol{\N}}\;.
\end{equation*}
As $\{1\}\rightarrow\ol{\T}$ induces a proper map on diamonds, the diagonal of the first map is even cohomologically proper by \cite[Thm.\ 6.3.1.(1)]{dRFF}. For the second map, we use the presentation from \cref{prop:prism-div1presentations} to see that the base change of its diagonal along the cover $\lim_{x\mapsto x^p} \ol{\T}\rightarrow \ol{\T}^{\Div^1}$ over $\Q_p^{\Div^1}$ identifies with
\begin{equation*}
\Z_p^\la\rightarrow \Z_p^\la\times_{\lim_{x\mapsto x^p} \ol{\T}} \Z_p^\la\;,
\end{equation*}
which is a map between affine Gelfand stacks. For the last map, \cref{thm:pres-tn} shows that the base change of its diagonal along the cover $\lim_{x\mapsto x^p} \ol{\T}\rightarrow \ol{\T}_{(p^{1/2}, p^{3/2})}^{\ol{\N}}$ over $\Q_{p, (p^{1/2}, p^{3/2})}^{\ol{\N}}$ is given by
\begin{equation*}
\left(\Z_p^\la\coprod^\AbGrp_{\G_a^\dagger} \G_a^\dagger\right)\rightarrow \left(\Z_p^\la\coprod^\AbGrp_{\G_a^\dagger} \G_a^\dagger\right)\times_{\lim_{x\mapsto x^p} \ol{\T}} \left(\Z_p^\la\coprod^\AbGrp_{\G_a^\dagger} \G_a^\dagger\right)\;.
\end{equation*}
Further base changing along the self-product of the cover $\Z_p^\la\times\G_a^\dagger\rightarrow \Z_p^\la\coprod^\AbGrp_{\G_a^\dagger} \G_a^\dagger$ over $\lim_{x\mapsto x^p} \ol{\T}$ then yields
\begin{equation*}
(\Z_p^\la\times \G_a^\dagger)\times_{\Z_p^\la\coprod^\AbGrp_{\G_a^\dagger} \G_a^\dagger} (\Z_p^\la\times\G_a^\dagger)\rightarrow (\Z_p^\la\times \G_a^\dagger)\times_{\lim_{x\mapsto x^p} \ol{\T}} (\Z_p^\la\times \G_a^\dagger)\;.
\end{equation*}
As the source is isomorphic to $(\Z_p^\la\times\G_a^\dagger)\times\G_a^\dagger$, this is a map between affine Gelfand stacks, which concludes the proof.
\end{proof}

Let us note the following consequence of \cref{thm:sixfunctors-propermaps}.

\begin{cor}
\label{cor:sixfunctors-relativepoincare}
Let $f: X\rightarrow Y$ be a rigid smooth map of derived Berkovich spaces which is quasicompact. Then the induced map $f^\Syn: X^\Syn\rightarrow Y^\Syn$ is cohomologically smooth and weakly cohomologically proper. In particular, for any perfect analytic $F$-gauge $E\in\D(X^\Syn)$ on $X$, the pushforward $f_*^\Syn E$ is perfect. If $f$ is additionally pure of relative dimension $d$, then there is an isomorphism
\begin{equation*}
f^\Syn_*(E^\vee\{d\})[2d]\cong (f^\Syn_*E)^\vee\;.
\end{equation*}
\end{cor}
\begin{proof}
For the first part, put together \cref{thm:sixfunctors-smoothmaps} and \cref{thm:sixfunctors-propermaps}. Then the second part follows from \cite[Lem.\ 4.5.16]{HeyerMann} and the third part follows by an analogous calculation as in the proof of \cref{cor:sixfunctors-qpsynsmooth} using that the dualising sheaf of $f^\Syn$ is given by $\O\{d\}[2d]$.
\end{proof}

\begin{rem}
We point out that any smooth proper map $f: X\rightarrow Y$ between partially proper rigid spaces over $\Q_p$ satisfies the assumptions of \cref{cor:sixfunctors-relativepoincare}.
\end{rem}

\subsection{Chern classes}

To conclude our discussion of the interaction of the syntomification with six functors, we want to establish that the syntomification of derived Berkovich spaces admits a strong theory of first Chern classes in the sense of \cite[§5]{Zavyalov}. Concretely, this amounts to the following statement:

\begin{prop}
\label{prop:sixfunctors-chern}
For any derived Berkovich space $X$, there is a natural morphism
\begin{equation*}
c_1^\Syn: R\Gamma(X_{\mathrm{Berk}\text{-}\et}, \G_m)[1]\rightarrow R\Gamma(X^\Syn, \O\{1\}[2])\;,
\end{equation*}
where the left-hand side denotes the cohomology of $\G_m[1]$ on the Berkovich étale site of $X$. Moreover, for any $d\geq 1$ and any derived Berkovich space $X$, the induced morphism
\begin{equation}
\label{eq:sixfunctors-cohpdviachern}
\sum_{0\leq k\leq d} (c_1^\Syn)^k\{d-k\}[2d-2k]: \bigoplus_{0\leq k\leq d} \O\{d-k\}[2d-2k]\rightarrow f_*^\Syn\O\{d\}[2d]
\end{equation}
is an isomorphism, where $f: \P^d_X\rightarrow X$ is the projection.
\end{prop}

Before we move on to the proof, let us shortly comment on how $c_1^\Syn$ induces a morphism (\ref{eq:sixfunctors-cohpdviachern}). Namely, note that $c_1^\Syn$ in particular induces a map 
\begin{equation*}
H^1((\P^d_X)_{\mathrm{Berk}\text{-}\et}, \G_m)\rightarrow H^0((\P^d_X)^\Syn, \O\{1\}[2])\;.
\end{equation*}
Thus, the tautological line bundle on $\P^d_X$ determines a class in $H^0((\P^d_X)^\Syn, \O\{1\}[2])$, which we may interpret as a map $c_1^\Syn: \O\rightarrow\O\{1\}[2]$ on $(\P^d_X)^\Syn$. By adjunction, this corresponds to a map $c_1^\Syn: \O\rightarrow f^\Syn_*\O\{1\}[2]$. Similarly, the tensor powers $(c_1^\Syn)^k: \O\rightarrow\O\{k\}[2k]$ correspond to maps $(c_1^\Syn)^k: \O\rightarrow f^\Syn_*\O\{k\}[2k]$, and these are the maps which yield (\ref{eq:sixfunctors-cohpdviachern}).

\begin{proof}[Proof of \cref{prop:sixfunctors-chern}]
Note that the functor $X\mapsto X^\Syn$ sends Berkovich étale covers to $!$-covers: Indeed, by locality on the target, this reduces to the claim for $X\mapsto X^\N$, where it is due to \cref{thm:defis-berketalecover}. Thus, the assignment $X\mapsto R\Gamma(X^\Syn, \O\{1\}[2])$ satisfies descent for Berkovich étale covers, whence it suffices to construct a natural transformation
\begin{equation*}
\tau^{\leq 0} R\Gamma(-, \G_m[1])\rightarrow \tau^{\leq 0}R\Gamma((-)^\Syn, \O\{1\}[2])\;.
\end{equation*}
Indeed, the sheafification of the source for the Berkovich étale topology will be given by the $\D(\Z)$-valued sheaf $R\Gamma((-)_{\mathrm{Berk}\text{-}\et}, \G_m[1])$, and this will yield the desired map.

To this end, by \cite[Rem.\ 5.1.10.(3)]{dRFF}, it suffices to construct a map 
\begin{equation*}
c_1^\Syn: \G_m^\Syn\rightarrow */\G_a\{1\}
\end{equation*}
of Gelfand stacks over $\Q_p^\Syn$. Then the desired natural transformation will be given by
\begin{equation*}
\begin{split}
\tau^{\leq 0} R\Gamma(X, \G_m[1])\rightarrow\Map(X, */\G_m)&\xrightarrow{\hphantom{*/c_1^\Syn}} \Map(X^\Syn, */\G_m^\Syn) \\
&\xrightarrow{*/c_1^\Syn}\Map(X^\Syn, */(*/\G_a\{1\}))\cong \tau^{\leq 0}R\Gamma(X^\Syn, \O\{1\}[2])\;.
\end{split}
\end{equation*}
Clearly, it will suffice to construct a map $c_1^\N: \G_m^\N\rightarrow */\G_a\{1\}$ of Gelfand stacks over $\Q_p^\N$ and then to check that the restrictions of this map along $j_\dR$ and $j_\HT$ induce the same map $\G_m^\prism\rightarrow */\G_a\{1\}$. Moreover, we may proceed locally on $\Q_p^\N$ as morphisms glue.

Over $\Q_{p, (0, p)}^\N$, we have $\G_{m, (0, p)}^\N\cong \G_{m, (0, p)}^\prism\times [0, 1]$ by \cref{prop:defis-utneq0}. Furthermore, using \cref{prop:defis-xprismxdiv1}, we may reduce to constructing a map $c_1^{\Div^1}: \G_m^{\Div^1}\rightarrow */\G_a\{1\}$ over $\Q_p^{\Div^1}$. However, recalling from \cref{prop:prism-div1presentations} that
\begin{equation*}
\G_m^{\Div^1}\cong \lim_{x\mapsto x^p} \G_m\,\Big/\,\Z_p^\la
\end{equation*}
and noting that the $\Z_p^\la$ that occurs here is actually implicitly Breuil--Kisin twisted by $1$, the desired map $c_1^{\Div^1}$ is induced by the canonical map $\Z_p^\la\rightarrow\G_a$ and the projection $\lim_{x\mapsto x^p} \G_m\rightarrow *$.

Similarly, over $\Q_{p, (p, \infty)}^\N$, our problem reduces to constructing a map $c_1^\HK: \G_m^\HK\rightarrow */\G_a\{1\}$ over $\Q_p^\HK$ using \cref{prop:defis-utneq0} and \cref{prop:defis-prismffdr}. However, for this we take the first Chern class for Hyodo--Kato cohomology constructed in \cite[Lem.\ 6.2.9]{dRFF}.

Finally, we have to construct the map $c_1^\N$ over $\Q_{p, (p^{1/2}, p^{3/2})}^\N$. To this end, we will actually construct a map
\begin{equation*}
c_1^{\ol{\N}}: \G_m^{\ol{\N}}\rightarrow */\G_a\{1\}
\end{equation*}
over $\Q_{p, (p^{1/2}, p^{3/2})}^{\ol{\N}}$ and then the desired map will be obtained by base changing to $\Q_{p, (p^{1/2}, p^{3/2})}^\N$. Now recall from \cref{thm:pres-tn} that
\begin{equation*}
\G_m^{\ol{\N}}\cong \lim_{x\mapsto x^p} \G_m\;\Big/\;\Z_p^\la\coprod_{\G_a^\dagger}^\AbGrp \G_a^\dagger
\end{equation*}
over this locus, where the pushout is along multiplication by $u$ on $\G_a^\dagger$. Noting that the $\Z_p^\la$  that occurs here should really be twisted by $\pi^*\O\{1\}$ while the $\G_a^\dagger$ should be twisted by $\O\{1\}\cong \pi^*\O\{1\}\tensor t^*\O(-1)$, the desired map is induced by the projection $\lim_{x\mapsto x^p} \G_m\rightarrow *$ and the maps
\begin{equation*}
\Z_p^\la\rightarrow\G_a\xrightarrow{\cdot u}\G_a\;, \hspace{0.5cm} \G_a^\dagger\rightarrow\G_a\;.
\end{equation*}

One easily checks that the maps from the previous paragraphs glue and that they induce the same map after pullback along $j_\dR$ and $j_\HT$, so that we obtain the desired map $c_1^\Syn$. Thus, it remains to verify the claim about the relative syntomic cohomology of projective space. As the map $(\P^d)^\Syn\rightarrow\Q_p^\Syn$ is weakly cohomologically proper by \cref{thm:sixfunctors-propermaps}, we are reduced to $X=\GSpec\Q_p$ by proper base change. In that case, we note that $(\P^d)^\Syn\rightarrow\Q_p^\Syn$ is also cohomologically smooth by \cref{thm:sixfunctors-smoothmaps} and thus $f^\Syn_*\O\{d\}[2d]$ is a perfect complex. Therefore, by \cref{thm:hkcomp-main} below, we may check the claimed isomorphism (\ref{eq:sixfunctors-cohpdviachern}) after applying the realisation functors $T_{\dR, +}$ and $T_\HK$ defined in §\ref{sect:hkcomp} below. However, as the syntomic first Chern class $c_1^\Syn$ we have constructed is compatible with the first Chern class $c_1^\HK$ for Hyodo--Kato cohomology from \cite[Lem.\ 6.2.9]{dRFF} and the usual first Chern class $c_1^{\dR, +}$ for filtered de Rham cohomology by \cite[Rem.\ 5.1.10]{dRFF}, the claim now follows from the fact that relative Hyodo--Kato cohomology and relative filtered de Rham cohomology of projective space may be computed in terms of Chern classes. Indeed, in the second case this is classical and in the first case it is due to \cite[Lem.\ 6.2.9]{dRFF}. This finishes the proof.
\end{proof}

\newpage

\section{$p$-adic Hodge theory for rigid-analytic varieties revisited}
\label{sect:padicht}

Having developed plenty of theory around the syntomification in the previous sections, we now want to show how to put it to use. As our first goal, we want to show how one can recover some results from Scholze's paper \cite{PAdicHodgeTheory} using the syntomification. To this end, we are going to analyse the equivalence
\begin{equation}
\label{eq:padicht-htdrequiv}
\Vect(X^{\dR, +})\cong \Vect(X^{\HT, \dagger, +})
\end{equation}
from \cref{prop:htdr-perfequiv} and relate it to the functor 
\begin{equation}
\label{eq:proet-scholzefunctorintro}
\{\text{filtered vector bundles with connection on $X$}\}\rightarrow \{\mathbb{B}_\dR^+\text{-local systems on }X_\proet\}
\end{equation}
from \cite[§7]{PAdicHodgeTheory}. 

Namely, from (\ref{eq:padicht-htdrequiv}), we obtain a natural functor
\begin{equation*}
\begin{split}
\Vect(X^{\dR, +})\cong \Vect(X^{\HT, \dagger, +})&\rightarrow\Vect(X^{\HT, \dagger}) \\
&\rightarrow \Vect((\widehat{X}\subseteq Y_{X^\diamond})^\wedge)\cong\{\mathbb{B}_\dR^+\text{-local systems on }X_\proet\}
\end{split}
\end{equation*}
via pullback and we are going to identify this functor with (\ref{eq:proet-scholzefunctorintro}). Indeed, using the results of §\ref{subsect:htdr}, we will be able to make the above functor explicit in terms of linear algebra when $X$ admits a Berkovich étale map $X\rightarrow\ol{\T}^n$ and, in this case, our claim reduces to a fairly simple computation. Then the key is to show that these local identifications glue, i.e.\ that they are independent of the chosen toric chart $X\rightarrow\ol{\T}^n$.

Finally, we will also be able to recover the classical de Rham comparison theorem
\begin{equation*}
R\Gamma_\proet(X_{\C_p}, \Q_p)\tensor_{\Q_p} B_\dR\cong R\Gamma_\dR(X)\tensor_{\Q_p} B_\dR
\end{equation*}
for smooth proper rigid spaces $X$ over $\Q_p$ using the syntomification. In fact, we will even be able to recover the relative version for smooth proper morphisms $X\rightarrow Y$ from \cite[Thm.\ 8.8]{PAdicHodgeTheory} and, as in loc.\ cit., using the results of §\ref{sect:proet}, we will also be able to allow for coefficients in any de Rham local system on $X_\proet$.

\subsection{Scholze's functor via the syntomification}

To begin, recall the \emph{structural de Rham sheaf} $\O\mathbb{B}_\dR$ on the proétale site $X_\proet$ of any smooth partially proper rigid space $X$ over $\Q_p$ defined in \cite[§6]{PAdicHodgeTheory}, see also \cite{padicHodgeTheoryErratum}. This is a $\mathbb{B}_\dR$-algebra equipped with a descending filtration $\Fil^\bullet \O\mathbb{B}_\dR$ and a $\mathbb{B}_\dR$-linear connection
\begin{equation*}
\nabla: \O\mathbb{B}_\dR\rightarrow\O\mathbb{B}_\dR\tensor_{\O_X} \Omega_X^1\;.
\end{equation*}
Using this, Scholze defines a functor
\begin{equation}
\label{eq:padicht-scholzefunctor}
\begin{split}
\{\text{filtered vector bundles with connection on $X$}\}&\rightarrow \{\mathbb{B}_\dR^+\text{-local systems on $X_\proet$}\} \\
(\Fil^\bullet E, \nabla)&\mapsto \Fil^0(E\tensor_{\O_X} \O\mathbb{B}_\dR)^{\nabla=0}\;,
\end{split}
\end{equation}
which he proves to be fully faithful, see \cite[Thm.\ 7.6]{PAdicHodgeTheory}; note that whenever we say ``filtered vector bundle with connection'', the filtration is always assumed to be locally by direct summands and the connection is always assumed to be flat and Griffiths transversal with respect to the filtration, i.e.\ it carries $\Fil^i E$ into $\Fil^{i-1} E\tensor_{\O_X} \Omega_X^1$.

Now recall that there is a natural map
\begin{equation*}
\iota: \widehat{X}\rightarrow\FF_{X^\diamond}
\end{equation*}
which identifies the source with a Cartier divisor in the target and is induced by the maps $\theta: \mathbb{A}_\inf(\ol{A})\rightarrow \ol{A}$ for any totally disconnected perfectoid $\ol{A}$ over $X$. Thus, we see that the infinitesimal neighbourhood of the image of $\iota$ in $\FF_{X^\diamond}$ is given by
\begin{equation*}
(\widehat{X}\subseteq \FF_{X^\diamond})^\wedge\cong \colim_n \colim_{\GSpec\ol{A}\rightarrow X} \GSpec \mathbb{B}_\dR^+(\ol{A})/\xi^n\;,
\end{equation*}
where the colimit runs over all totally disconnected perfectoids $\ol{A}$ over $X$ and $\xi$ denotes a generator of $\Ker\theta$. As $\mathbb{B}_\dR^+(\ol{A})$ is $\xi$-complete, we conclude that
\begin{equation}
\label{eq:padicht-bdr+locsysviaffx}
\Vect((\widehat{X}\subseteq \FF_{X^\diamond})^\wedge)\cong \{\mathbb{B}_\dR^+\text{-local systems on $X_\proet$}\}\;.
\end{equation}

Furthermore, observe that the natural map $Y_{X^\diamond}\rightarrow X^\prism$ from \cref{prop:defis-xprismdr} restricts to a map 
\begin{equation*}
(\widehat{X}\subseteq \FF_{X^\diamond})^\dagger\rightarrow X^{\HT, \dagger}
\end{equation*}
from the overconvergent neighbourhood of $\widehat{X}$ in $\FF_{X^\diamond}$ to $X^{\HT, \dagger}$. Further restricting this map to the infinitesimal neighbourhood and using the equivalence (\ref{eq:padicht-bdr+locsysviaffx}), we obtain a pullback functor
\begin{equation}
\label{eq:padicht-xhtdaggertobdr+}
\Vect(X^{\HT, \dagger})\rightarrow \{\mathbb{B}_\dR^+\text{-local systems on $X_\proet$}\}\;.
\end{equation}

Putting (\ref{eq:padicht-xhtdaggertobdr+}) together with the equivalence $\Vect(X^{\dR, +})\cong \Vect(X^{\HT, \dagger, +})$, which one deduces from \cref{prop:htdr-perfequiv} as the argument in loc.\ cit.\ just as well applies to any $\Perf^{[a, b]}(-)$ in place of $\Perf(-)$, we arrive at a composite functor
\begin{equation}
\label{eq:padicht-xdr+tobdr+}
\Vect(X^{\dR, +})\cong \Vect(X^{\HT, \dagger, +})\rightarrow \Vect(X^{\HT, \dagger})\rightarrow  \{\mathbb{B}_\dR^+\text{-local systems on $X_\proet$}\}\;,
\end{equation}
where the functor $\Vect(X^{\HT, \dagger, +})\rightarrow \Vect(X^{\HT, \dagger})$ is of course given by pullback. Our main result in this section is the following:

\begin{thm}
\label{thm:padicht-scholzefunctor}
Let $X$ be a smooth partially proper rigid space over $\Q_p$. Under the identification between vector bundles on $X^{\dR, +}$ and filtered vector bundles with connection on $X$ from \cref{prop:recall-fildr}, the functors (\ref{eq:padicht-scholzefunctor}) and (\ref{eq:padicht-xdr+tobdr+}) are naturally equivalent.
\end{thm}

\begin{rem}
\label{rem:drbundles-functortobdr+dagger}
In principle, we could extract even finer information from a vector bundle on $X^{\HT, \dagger}$. Namely, note that
\begin{equation*}
(\widehat{X}\subseteq \FF_{X^\diamond})^\dagger\cong \colim_{\GSpec\ol{A}\rightarrow X} \mathbb{B}_\dR^{+, \dagger}(\ol{A})\;,
\end{equation*}
where the colimit again runs over all totally disconnected perfectoids $\ol{A}$ over $X$ and $\mathbb{B}_\dR^{+, \dagger}$ is the \emph{overconvergent de Rham period sheaf} on $X_\proet$ from \cite[Def.\ 3.4.1]{Wiersig1}. Thus, we even obtain a pullback functor
\begin{equation*}
\D(X^{\HT, \dagger, +})\rightarrow \D(X_\proet, \mathbb{B}_\dR^{+, \dagger})
\end{equation*}
and, in particular, the functor (\ref{eq:padicht-xdr+tobdr+}) refines to a functor
\begin{equation}
\label{eq:padicht-functortobdr+dagger}
\Vect(X^{\dR, +})\cong \Vect(X^{\HT, \dagger, +})\rightarrow \Vect(X^{\HT, \dagger})\rightarrow  \{\mathbb{B}_\dR^{+, \dagger}\text{-local systems on $X_\proet$}\}\;.
\end{equation}
Note that Scholze's functor (\ref{eq:padicht-scholzefunctor}) admits such a refinement as well: Namely, this is given by
\begin{equation}
\label{eq:padicht-scholzetobdr+dagger}
\begin{split}
\{\text{filtered vector bundles with connection on $X$}\}&\rightarrow \{\mathbb{B}_\dR^+\text{-local systems on $X_\proet$}\} \\
(\Fil^\bullet W, \nabla)&\mapsto \Fil^0(W\tensor_{\O_X} \O\mathbb{B}_\dR^\dagger)^{\nabla=0}\;,
\end{split}
\end{equation}
where $\O\mathbb{B}_\dR^{+, \dagger}$ is the structural overconvergent de Rham period sheaf from \cite[§3.5]{Wiersig1}. 

Naturally, one might wonder whether \cref{thm:padicht-scholzefunctor} refines to an equivalence between the functors (\ref{eq:padicht-functortobdr+dagger}) and (\ref{eq:padicht-scholzetobdr+dagger}). Indeed, one can check that the proof below still works in this situation: One just has to replace the use of \cite[Prop.\ 6.10]{PAdicHodgeTheory} by \cite[Thm.\ 3.5.7]{Wiersig1} and, moreover, \cite[Thm.\ 7.2]{PAdicHodgeTheory} is replaced by the analogous fact that, for any $\Q_p$-algebra $R$, a finite projective $R\{X_1, \dots, X_n\}^\dagger$-module with a flat connection has enough horizontal sections, see e.g.\ the proof of \cite[Thm.\ 9.6.1]{KedlayaDiffEqs}.
\end{rem}

The key to proving \cref{thm:padicht-scholzefunctor} is the following explicit identification of the functor $\Vect(X^{\dR, +})\rightarrow \Vect(X^{\HT, \dagger})$ from (\ref{eq:padicht-xdr+tobdr+}) in the case where $X$ is equipped with a toric chart.

\begin{lem}
\label{lem:padicht-filvsplits}
Let $X$ be a derived Berkovich space over $\Q_p$ equipped with a Berkovich étale map $X\rightarrow\ol{\T}^n$. For any perfect complex $E$ on $X^{\HT, \dagger, +}$, the filtration $\Fil_\bullet V$ from \cref{cor:htdr-complexesxhtdrdagger+} splits $\Z_p^n\rtimes\Z_p^\times$-equivariantly, i.e., in the notation of loc.\ cit., there is a $\Z_p^n\rtimes\Z_p^\times$-equivariant isomorphism
\begin{equation*}
\Fil_\bullet V\cong \left(\bigoplus_{i\leq\bullet} t^{-i}(\Fil^i W\tensor_{\O_X} \O_{X_\infty^\la})\tensor_{\Q_p} \Q_p(\zeta_{p^\infty})\right)\tensor_{\Q_p[r]} \Q_p\{r\}^\dagger\;,
\end{equation*}
where $r$ acts via the transition maps $t^{-i}\Fil^i W\rightarrow t^{-i+1} \Fil^{i-1} W$, the $\Z_p^n$-action is induced by the one on $\O_{X_\infty^\la}$ and $\Z_p^\times$ acts on $t$ by multiplication and on $\Q_p(\zeta_{p^\infty})$ via the usual Galois action.
\end{lem}
\begin{proof}
Recall from loc.\ cit.\ that there are $\Z_p^n\rtimes\Z_p^\times$-equivariant isomorphisms
\begin{equation*}
\gr_i V\cong t^{-i}(\Fil^i W\tensor_{\O_X} \O_{X_\infty^\la}\tensor_{\Q_p} \Q_p(\zeta_{p^\infty}))
\end{equation*}
for all $i\in\Z$ and for the filtration to split it suffices to show that
\begin{equation*}
\RHom(\gr_j V, \gr_i V)=0
\end{equation*}
for all $i\neq j$, where the $\RHom$ is taken in $\D(X_\infty^\la\times\GSpec\Q_p(\zeta_{p^\infty})\,/\,(\Z_p^n)^\la\rtimes\Z_p^{\times, \la})$. As $\gr_j V$ is a perfect complex and hence dualisable, this amounts to showing that the cohomology of 
\begin{equation}
\label{eq:drbundles-grtensor}
\gr_i V\tensor(\gr_j V)^\vee\cong t^{-i+j}(\Fil^i W\tensor_{\O_X} (\Fil^j W)^\vee\tensor_{\O_X} \O_{X_\infty^\la} \tensor_{\Q_p} \Q_p(\zeta_{p^\infty}))
\end{equation}
on $X_\infty^\la\times\GSpec\Q_p(\zeta_{p^\infty})\,/\,(\Z_p^n)^\la\rtimes\Z_p^{\times, \la}$ vanishes. 

However, note that pushforward along
\begin{equation*}
X_\infty^\la\times\GSpec\Q_p(\zeta_{p^\infty})\,/\,(\Z_p^n)^\la\rtimes\Z_p^{\times, \la}\rightarrow X_\infty^\la\times\GSpec\Q_p(\zeta_{p^\infty})\,/\,(\Z_p^n)^\sm\rtimes\Z_p^{\times, \sm}\cong X
\end{equation*}
is given by taking Lie algebra cohomology. Moreover, by (\ref{eq:drbundles-grtensor}), the action of $\Z_p^n$ on $\gr_i V\tensor(\gr_j V)^\vee$ is smooth while the Lie algebra of $\Z_p^{\times, \la}$ acts on $\gr_i V\tensor(\gr_j V)^\vee$ by multiplication by $j-i\neq 0$. As $j-i$ is invertible, this yields
\begin{equation*}
R\Gamma(\operatorname{Lie}(\Z_p^n\rtimes\Z_p^\times), \gr_i V\tensor(\gr_j V)^\vee)=0\;,
\end{equation*}
whence the pushforward of $\gr_i V\tensor (\gr_j V)^\vee$ to $X$ and consequently also its cohomology vanishes, as desired. 
\end{proof}

The above splitting of the filtration $\Fil_\bullet V$ already allows us to prove \cref{thm:padicht-scholzefunctor} in the case where $X$ admits a toric chart $X\rightarrow\ol{\T}^n$.

\begin{prop}
\label{prop:padicht-scholzefunctorlocal}
Let $X$ be a smooth partially proper rigid space over $\Q_p$ equipped with an étale map $X\rightarrow\Spa\Q_p\langle x_1^{\pm 1}, \dots, x_n^{\pm 1}\rangle$ such that the composition 
\begin{equation*}
X\rightarrow\GSpec\Q_p\langle x_1^{\pm 1}, \dots, x_n^{\pm 1}\rangle\rightarrow\ol{\T}^n
\end{equation*}
is Berkovich étale. Then the functor
\begin{equation*}
\Vect(X^{\dR, +})\cong \Vect(X^{\HT, \dagger, +})\rightarrow \Vect(X^{\HT, \dagger})\rightarrow \{\mathbb{B}_\dR^+\text{-local systems on $X_\proet$}\}
\end{equation*}
from (\ref{eq:padicht-xdr+tobdr+}) agrees with Scholze's functor
\begin{equation*}
\{\text{filtered vector bundles with connection on $X$}\}\rightarrow \{\mathbb{B}_\dR^+\text{-local systems on $X_\proet$}\}
\end{equation*}
from (\ref{eq:padicht-scholzefunctor}) under the identification between vector bundles on $X^{\dR, +}$ and filtered vector bundles with connection on $X$ from \cref{prop:recall-fildr}.
\end{prop}
\begin{proof}
Let $X_\infty^\cycl$ denote the adic space
\begin{equation*}
X_\infty^\cycl\coloneqq X\times_{\Spa\Q_p\langle x_1^{\pm 1}, \dots, x_n^{\pm 1}\rangle} \Spa\Q_p^\cycl\langle x_1^{\pm 1/p^\infty}, \dots, x_n^{\pm 1/p^\infty}\rangle
\end{equation*}
and observe that this is perfectoid. By \cref{lem:padicht-filvsplits}, our task is to check that the base change of
\begin{equation}
\label{eq:padicht-formulav}
\left(\bigoplus_i t^{-i}(\Fil^i W\tensor_{\O_X} \O_{X_\infty^\la} \tensor_{\Q_p} \Q_p(\zeta_{p^\infty}))\right)\tensor_{\Q_p[t]} \Q_p\{t\}^\dagger
\end{equation}
along the map $\O_{X_\infty^\la}\tensor_{\Q_p} \Q_p(\zeta_{p^\infty})\{t\}^\dagger\rightarrow\mathbb{B}_{\dR, X_\infty^\cycl}^+$ identifies $\Z_p^n\rtimes\Z_p^\times$-equivariantly with
\begin{equation*}
\Fil^0(W\tensor_{\O_X} \O\mathbb{B}_{\dR, X_\infty^\cycl})^{\nabla=0}\;.
\end{equation*}
Here, we regard $\O_{X_\infty^\la}$ and $\mathbb{B}_{\dR, X_\infty^\cycl}^+$ as sheaves on the underlying topological space of $X_\infty^\cycl$; indeed, note that the underlying topological space of the derived Berkovich space $X_\infty^\la\times_{\GSpec\Q_p} \GSpec\Q_p(\zeta_{p^\infty})$ in the sense of \cite[Def.\ 4.3.1]{dRFF} is the maximal Hausdorff quotient of the underlying topological space of $X_\infty^\cycl$ and hence $\O_{X_\infty^\la}$ can be regarded as a sheaf on the latter via pullback. Moreover, we note that $\O_{X_\infty^\la}\tensor_{\Q_p} \Q_p(\zeta_{p^\infty})\{t\}^\dagger\rightarrow\mathbb{B}_{\dR, X_\infty^\cycl}^+$ is the unique map given by $x_i^{1/p^m}\mapsto [x_i^\flat]^{1/p^m}$ and $t\mapsto \log[\epsilon]\in\mathbb{B}_\dR^+(\Q_p^\cycl)$, where we use the proof of \cite[Prop.\ 6.10]{PAdicHodgeTheory} for the existence and uniqueness and write $\epsilon=(1, \zeta_p, \zeta_{p^2}, \dots)$, as usual.

With these preliminaries out of the way, let us move on to the actual proof. Examining the formula (\ref{eq:padicht-formulav}), it is clear that it will suffice to prove that there is a $\Z_p^n\rtimes\Z_p^{\times, \la}$-equivariant isomorphism
\begin{equation}
\label{eq:padicht-scholzefunctoriso}
\Fil^0(W\tensor_{\O_X} \O\mathbb{B}_\dR)^{\nabla=0}\cong \sum_i \Fil^i W\tensor_{\O_X} t^{-i}\mathbb{B}_\dR^+
\end{equation}
of sheaves on the localised site $X_\proet/X_\infty^\cycl$, where the map $\O_X\rightarrow\mathbb{B}_\dR^+$ is obtained in the same way as the map $\O_{X_\infty^\la}\rightarrow\mathbb{B}_{\dR, X_\infty^\cycl}^+$ above. However, note that, using \cite[Prop.\ 6.10]{PAdicHodgeTheory}, this map may also be obtained as the composition
\begin{equation*}
\O_X\rightarrow\O\mathbb{B}_\dR^+\cong \mathbb{B}_\dR^+[\![X_1, \dots, X_n]\!]\rightarrow\mathbb{B}_\dR^+\;,
\end{equation*}
where the last map is reduction mod $(X_1, \dots, X_n)$ for $X_i=x_i-[x_i^\flat]$. Now observe that the natural map
\begin{equation*}
\left(\sum_i \Fil^i W\tensor_{\O_X} t^{-i}\O\mathbb{B}_\dR^+\right)^{\nabla=0}\tensor_{\mathbb{B}_\dR^+} \O\mathbb{B}_\dR^+\rightarrow \sum_i \Fil^i W\tensor_{\O_X} t^{-i}\O\mathbb{B}_\dR^+
\end{equation*}
is an isomorphism by \cite[Thm.\ 7.2]{PAdicHodgeTheory}. Thus, reducing mod $(X_1, \dots, X_n)$, we see that
\begin{equation*}
\left(\sum_i \Fil^i W\tensor_{\O_X} t^{-i}\O\mathbb{B}_\dR^+\right)^{\nabla=0}\cong \sum_i \Fil^i W\tensor_{\O_X} t^{-i}\mathbb{B}_\dR^+
\end{equation*}
and this yields the desired isomorphism (\ref{eq:padicht-scholzefunctoriso}) by \cref{lem:padicht-fil0horizontal} below.
\end{proof}

The following lemma was used in the proof:

\begin{lem}
\label{lem:padicht-fil0horizontal}
Let $X$ be a smooth partially proper rigid space over $\Q_p$. For any filtered vector bundle with connection $(\Fil^\bullet W, \nabla)$ on $X$, the natural map
\begin{equation*}
\Fil^0(W\tensor_{\O_X} \O\mathbb{B}_\dR)^{\nabla=0}\xrightarrow{\cong} \left(\sum_i \Fil^i W\tensor_{\O_X} t^{-i}\O\mathbb{B}_\dR^+\right)^{\nabla=0}
\end{equation*}
is an isomorphism.
\end{lem}
\begin{proof}
The question is local on $X$ and hence we may assume that $X$ is equipped with an étale map $X\rightarrow \Spa\Q_p\langle x_1^{\pm 1}, \dots, x_n^{\pm 1}\rangle$. Letting $X^\cycl_\infty$ be the base change of $X$ to $\Spa\Q_p^\cycl\langle x_1^{\pm 1/p^\infty}, \dots, x_n^{\pm 1/p^\infty}\rangle$, we may further restrict to the localised site $X_\proet/X^\cycl_\infty$. Then recall that $\O\mathbb{B}_\dR^+\cong \mathbb{B}_\dR^+[\![X_1, \dots, X_n]\!]$, where $X_i=x_i-[x_i^\flat]$, and that the filtration on $\O\mathbb{B}_\dR=\O\mathbb{B}_\dR^+[\tfrac{1}{t}]$ is given by
\begin{equation*}
\Fil^i \O\mathbb{B}_\dR=\sum_{j\geq -i} t^{-j}(t, X_1, \dots, X_n)^{i+j}\O\mathbb{B}_\dR^+\;.
\end{equation*}
Thus, we have
\begin{equation*}
\Fil^0(W\tensor_{\O_X} \O\mathbb{B}_\dR)=\sum_i \Fil^i W\tensor_{\O_X} \left(\sum_{j\geq i} t^{-j}(t, X_1, \dots, X_n)^{j-i}\O\mathbb{B}_\dR^+\right)
\end{equation*}
and our task is to show that a horizontal section in $\Fil^i W\tensor_{\O_X} (t, X_1, \dots, X_n)^\ell\O\mathbb{B}_\dR^+$ is divisible by $t^\ell$ for all $\ell\geq 1$. Moreover, since the transition maps of the filtration $\Fil^\bullet W$ are injections, it suffices to check the above for $W$ in place of $\Fil^i W$.

Localising further on $X$ if necessary, we may assume that $W$ is trivial and let $e_1, \dots, e_m$ be a basis. Then the connection on $W$ is given by
\begin{equation*}
\nabla(e_i)=\sum_{j, k} f_{ijk}e_k\mathrm{d}x_j
\end{equation*}
for some global sections $f_{ijk}$ of $\O_X$. Now letting 
\begin{equation*}
\sum_i e_i\tensor f_i\in W\tensor_{\O_X} (t, X_1, \dots, X_n)^\ell\O\mathbb{B}_\dR^+
\end{equation*}
be a horizontal section, we obtain the differential equations
\begin{equation}
\label{eq:padicht-diffeqsbdr+}
\frac{\partial f_k}{\partial X_j}=-\sum_i f_{ijk}f_i
\end{equation}
for the elements $f_k\in\O\mathbb{B}_\dR^+\cong \mathbb{B}_\dR^+[\![X_1, \dots, X_n]\!]$, where we are using that $\nabla(X_j)=\mathrm{d} x_j$ for the connection $\nabla$ on $\O\mathbb{B}_\dR^+$.

Working modulo $t^\ell$ from now on, our assumption implies that each $f_k$ is contained in the ideal $(X_1, \dots, X_n)$. However, this means that the right-hand side of (\ref{eq:padicht-diffeqsbdr+}) has vanishing constant term, which in turn means that each $f_k$ must be contained in $(X_1, \dots, X_n)^2$. Inductively, we see that, modulo $t^\ell$, each $f_k$ is contained in all powers of the ideal $(X_1, \dots, X_n)$, which implies that $f_k$ vanishes modulo $t^\ell$. In other words, each $f_k$ is divisible by $t^\ell$ and this is what we wanted to show.
\end{proof}

Finally, we can prove \cref{thm:padicht-scholzefunctor}.

\begin{proof}[Proof of \cref{thm:padicht-scholzefunctor}]
Note that \cref{prop:padicht-scholzefunctorlocal} already establishes the claim in the situation where $X$ admits an étale map $X\rightarrow\Spa\Q_p\langle x_1^{\pm 1}, \dots, x_n^{\pm 1}\rangle$ whose postcomposition with the canonical map to $\ol{\T}^n$ is Berkovich étale. As this can always be arranged locally on $X$, our task is to check that the equivalence from loc.\ cit.\ glues across different charts. 

For this, take two étale charts $X\rightarrow\Spa \Q_p\langle x_1^{\pm 1}, \dots, x_n^{\pm 1}\rangle$ and $X\rightarrow\Spa \Q_p\langle \widetilde{x}_1^{\pm 1}, \dots, \widetilde{x}_1^{\pm 1}\rangle$ and define $X_{\infty, 0}^\cycl$ and $X_{0, \infty}^\cycl$ by
\begin{equation*}
X_{\infty, 0}^\cycl\coloneqq X\times_{\Spa\Q_p\langle x_1^{\pm 1}, \dots, x_n^{\pm 1}\rangle} \Spa\Q_p^\cycl\langle x_1^{\pm 1/p^\infty}, \dots, x_n^{\pm 1/p^\infty}\rangle\;,
\end{equation*}
and similarly for $X_{0, \infty}^\cycl$ and the coordinates $\widetilde{x}_i$; moreover, put $X_{\infty, \infty}^\cycl\coloneqq X_{\infty, 0}^\cycl\times_X X_{0, \infty}^\cycl$. Analogously, we define 
\begin{equation*}
X_{\infty, 0}^\la\coloneqq X\times_{\ol{\T}^n} \lim_{x_i\mapsto x_i^p} \ol{\T}^n\;,
\end{equation*}
where the map $X\rightarrow\ol{\T}^n$ is induced by the coordinates $x_i$, and similarly we define $X_{0, \infty}^\la$ in terms of the coordinates $\widetilde{x}_i$; again, we put $X_{\infty, \infty}^\la\coloneqq X_{\infty, 0}^\la\times_X X_{0, \infty}^\la$.

On the localised site $X_\proet/X_{\infty, \infty}^\cycl$, the proof of \cref{prop:padicht-scholzefunctorlocal} supplies isomorphisms
\begin{equation}
\label{eq:padicht-changechartcompatible}
\begin{tikzcd}[column sep=small]
& \Fil^0(W\tensor_{\O_X} \O\mathbb{B}_\dR)^{\nabla=0}\ar[dl, "\cong", swap]\ar[dr, "\cong"] \\
\displaystyle\sum_i \Fil^i W\tensor_{\O_X, x_i\mapsto [x_i^\flat]} t^{-i}\mathbb{B}_\dR^+ & & \displaystyle\sum\limits_i \Fil^i W \tensor_{\O_X, \widetilde{x}_i\mapsto [\widetilde{x}_i^\flat]} t^{-i}\mathbb{B}_\dR^+\nospacepunct{\;,}
\end{tikzcd}
\end{equation}
where we emphasise that the base change from $\O_X$ to $\mathbb{B}_\dR^+$ is along two different maps, one of them induced by $x_i\mapsto [x_i^\flat]$ and the other induced by $\widetilde{x}_i\mapsto [\widetilde{x}_i^\flat]$, as the notation suggests. Our task is to check that the isomorphism
\begin{equation}
\label{eq:padicht-changechart}
\sum_i \Fil^i W\tensor_{\O_X, x_i\mapsto [x_i^\flat]} t^{-i}\mathbb{B}_\dR^+\xrightarrow{\cong}\sum\limits_i \Fil^i W \tensor_{\O_X, \widetilde{x}_i\mapsto [\widetilde{x}_i^\flat]} t^{-i}\mathbb{B}_\dR^+
\end{equation}
induced by changing charts makes the diagram (\ref{eq:padicht-changechartcompatible}) commute. To see where this last isomorphism comes from, note that, as in the proof of \cref{prop:padicht-scholzefunctorlocal}, there is a unique map $\O_{X_{\infty, \infty}^\la}\rightarrow \mathbb{B}_{\dR, X_{\infty, \infty}^\cycl}^+$ induced by $x_i\mapsto [x_i^\flat]$ and then (\ref{eq:padicht-changechart}) is induced by the unique $\mathbb{B}_\dR^+(\Q_p^\cycl)$-linear automorphism of $\mathbb{B}_{\dR, X_{\infty, \infty}^\cycl}^+$ making the diagram
\begin{equation}
\label{eq:padicht-changechartbdr+}
\begin{tikzcd}[column sep=small]
& \O_{X_{\infty, \infty}^\la}\ar[dl, "{x_i\mapsto [x_i^\flat]}", swap]\ar[dr, "{\widetilde{x}_i\mapsto [\widetilde{x}_i^\flat]}"] & \\
\mathbb{B}_{\dR, X_{\infty, \infty}^\cycl}^+\ar[rr, "\cong"] & & \mathbb{B}_{\dR, X_{\infty, \infty}^\cycl}^+
\end{tikzcd}
\end{equation} 
commute. Indeed, such an automorphism exists uniquely: We have
\begin{equation*}
\mathbb{B}_\dR^+(\Q_p^\cycl\langle x_1^{\pm 1/p^\infty}, \dots, x_n^{\pm 1/p^\infty}\rangle)\cong \mathbb{B}_\dR^+(\Q_p^\cycl)\langle [x_1^\flat]^{1/p^\infty}, \dots, [x_n^\flat]^{1/p^\infty}\rangle^\wedge_\xi
\end{equation*}
and $\mathbb{B}_{\dR, X_{\infty, \infty}^\cycl}^+$ is a $(p, \xi)$-completed colimit of finite étale algebras over this ring as $\mathbb{B}_\dR^+(-)$ preserves finite étale maps, so there is a unique solution to the lifting problem
\begin{equation*}
\begin{tikzcd}
\mathbb{B}_\dR^+(\Q_p^\cycl\langle x_1^{\pm 1/p^\infty}, \dots, x_n^{\pm 1/p^\infty}\rangle)\ar[r]\ar[d] & \mathbb{B}_{\dR, X_{\infty, \infty}^\cycl}^+\ar[d] \\
\mathbb{B}_{\dR, X_{\infty, \infty}^\cycl}^+\ar[r]\ar[ur, dotted] & \O_{X_{\infty, \infty}^\cycl}\nospacepunct{\;.}
\end{tikzcd}
\end{equation*}

Recalling from the proof of \cref{prop:padicht-scholzefunctorlocal} that the left map in (\ref{eq:padicht-changechartcompatible}) is obtained by setting $X_1=x_1-[x_1^\flat], \dots, X_n=x_n-[x_n^\flat]$ to zero while the right map sends $\widetilde{X}_1=\widetilde{x}_1-[\widetilde{x}_1^\flat], \dots, \widetilde{X}_n=\widetilde{x}_n-[\widetilde{x}_n^\flat]$ to zero, we are reduced to checking that the diagram
\begin{equation*}
\begin{tikzcd}[column sep=small]
& \O\mathbb{B}_{\dR, X_{\infty, \infty}^\cycl}^+\ar[dl, "{x_i\mapsto [x_i^\flat]}", swap]\ar[dr, "{\widetilde{x}_i\mapsto [\widetilde{x}_i^\flat]}"] \\
\mathbb{B}_{\dR, X_{\infty, \infty}^\cycl}^+ \ar[rr, "\cong"] & & \mathbb{B}_{\dR, X_{\infty, \infty}^\cycl}^+
\end{tikzcd}
\end{equation*}
commutes, where the horizontal map on the bottom is given by (\ref{eq:padicht-changechartbdr+}). However, this is immediate from the diagram (\ref{eq:padicht-changechartbdr+}) and this concludes the proof.
\end{proof}

\subsection{The de Rham comparison theorem}

Having proved \cref{thm:padicht-scholzefunctor}, we know that the following definition is compatible with the terminology from \cite[§7]{PAdicHodgeTheory} in the case of vector bundles.

\begin{defi}
Let $X$ be a smooth partially proper rigid space over $\Q_p$. For any perfect complex $E$ on $X^{\dR, +}$, its \emph{associated perfect complex of $\mathbb{B}_\dR^+$-modules} is defined to be the image of $E$ under the functor
\begin{equation}
\label{eq:padicht-perfxdr+tobdr+}
\Perf(X^{\dR, +})\cong \Perf(X^{\HT, \dagger, +})\rightarrow \Perf(X^{\HT, \dagger})\rightarrow\Perf(X_\proet, \mathbb{B}_\dR^+)
\end{equation}
defined analogously to (\ref{eq:padicht-xdr+tobdr+}).
\end{defi}

\begin{rem}
Imitating the proof of \cref{thm:padicht-scholzefunctor}, one can show that, if $E$ is a perfect complex on $X^{\dR, +}$ corresponding to a filtered perfect complex with connection $(\Fil^\bullet W, \nabla)$ on $X$, its associated perfect complex of $\mathbb{B}_\dR^+$-modules $\mathbb{M}$ satisfies
\begin{equation*}
\Fil^\bullet W\tensor_{\O_X} \O\mathbb{B}_\dR^+\cong \mathbb{M}\tensor_{\mathbb{B}_\dR^+} \O\mathbb{B}_\dR^+\;. \qedhere
\end{equation*}
\end{rem}

\begin{rem}
We warn the reader that, contrary to the functor 
\begin{equation*}
\Vect(X^{\dR, +})\rightarrow \{\mathbb{B}_\dR^+\text{-local systems on $X_\proet$}\}
\end{equation*}
on the abelian level, the functor (\ref{eq:padicht-perfxdr+tobdr+}) is \emph{not} fully faithful. 
\end{rem}

With this terminology in place, we can prove the following version of the relative de Rham comparison theorem, which the reader should compare with \cite[Thm.\ 8.8]{PAdicHodgeTheory}.

\begin{thm}
\label{thm:padicht-drcomp}
Let $f: X\rightarrow S$ be a smooth proper morphism of smooth partially proper rigid spaces over $\Q_p$. Then the pushforward $f^{\dR, +}_*\O\langle -i\rangle\in\D(S^{\dR, +})$ is perfect and its associated perfect complex of $\mathbb{B}_\dR^+$-modules on $S_\proet$ is given by $f_{\proet, *}\Q_p(i)\tensor_{\Q_p} \mathbb{B}_\dR^+$. In particular, for any smooth proper rigid space $X$ over $\Q_p$, there is an isomorphism
\begin{equation*}
R\Gamma_\proet(X_{\C_p}, \Q_p)\tensor_{\Q_p} B_\dR\cong R\Gamma_\dR(X)\tensor_{\Q_p} B_\dR
\end{equation*}
compatible with the $\Gal_{\Q_p}$-actions and the filtrations.
\end{thm}
\begin{proof}
By \cref{cor:sixfunctors-relativepoincare}, the induced map $f^\Syn: X^\Syn\rightarrow S^\Syn$ is cohomologically smooth and weakly cohomologically proper and the pushforward $E\coloneqq f^\Syn_*\O\{i\}$ is a perfect complex on $S^\Syn$.

By proper base change, we know that $i_{\dR, +}^*E$ identifies with $f^{\dR, +}_*\O\langle -i\rangle$. Moreover, by \cite{AnPrism} and using proper base change again, we also know that the pullback of $E|_{S^\prism}$ along $Y_{S^\diamond}\rightarrow S^\prism$ identifies with the perfect complex on $Y_{S^\diamond}$ obtained from $f_{\proet, *}\Q_p(i)$ by tensoring with $\O_{Y_{S^\diamond}}$. In particular, we conclude that the functor $\Perf(S^{\dR, +})\rightarrow \Perf(S_\proet, \mathbb{B}_\dR^+)$ sends $f^{\dR, +}_*\O\langle -i\rangle$ to $f_{\proet, *}\Q_p(i)\tensor_{\Q_p}\mathbb{B}_\dR^+$ and this implies the claim. Finally, the last assertion follows by specialising to $S=\GSpec\Q_p$ using \cref{prop:recall-fildr}.
\end{proof}

\begin{rem}
In \cref{thm:drbundles-main} below, we will see that vector bundles on $X^\Syn$ identify with vector bundles on $\FF_{X^\diamond}$ which are de Rham in a certain sense, see \cref{defi:drbundles-drvb}. In particular, one concludes that there is a fully faithful embedding
\begin{equation*}
\Loc^\dR(X_\proet, \Q_p)\hookrightarrow\Vect(X^\Syn)
\end{equation*}
of the category of de Rham local systems in the sense of \cite[Def.\ 8.3]{PAdicHodgeTheory} into vector bundle analytic $F$-gauges on $X$. Then the same proof as above shows that, as in \cite[Thm.\ 8.8]{PAdicHodgeTheory}, the conclusion of \cref{thm:padicht-drcomp} stays valid if we replace the Tate twist $\Q_p(i)$ by any de Rham local system $\mathbb{L}$ on $X$.
\end{rem}

\newpage

\section{The comparison with filtered Hyodo--Kato cohomology}
\label{sect:hkcomp}

Recall that the main utility of syntomic cohomology classically is that it relates to both ``de Rham'' and ``étale data'' via comparison theorems such as the ones proved in \cite{padicComparisons} or \cite{BasicComparison}. In this spirit, we want to use the present section to show how perfect complexes on $X^\Syn$ and their cohomology can be understood through a ``de Rham lens'' for Berkovich smooth derived Berkovich spaces $X$; in particular, this includes smooth partially proper rigid spaces over $\Q_p$ and thus also relates to more classical settings of $p$-adic analytic geometry.

More specifically, we are going to show that perfect complexes on the syntomification $X^\Syn$ of a Berkovich smooth derived Berkovich space and their cohomology can be completely understood via perfect complexes on the Hyodo--Kato stack $X^\HK$ and on the filtered de Rham stack $X^{\dR, +}$ of $X$. Namely, we are going to construct a commutative diagram
\begin{equation}
\label{eq:hkcomp-square}
\begin{tikzcd}
X^\dR\ar[r]\ar[d] & X^{\dR, +}\ar[d, "i_{\dR, +}"] \\
X^\mHK\ar[r, "i_\HK"] & X^\Syn\nospacepunct{\;,}
\end{tikzcd}
\end{equation}
where the reader may mentally replace $X^\mHK$ by the usual Hyodo--Kato stack $X^\HK$ for the moment, and show that it induces an equivalence of categories
\begin{equation}
\label{eq:hkcomp-perffibreproduct}
\Perf(X^\Syn)\cong \Perf(X^\HK)\times_{\Perf(X^\dR)} \Perf(X^{\dR, +})
\end{equation}
via pullback. We note that this result should be compared to the so-called Beilinson fibre square for syntomic cohomology of $p$-adic formal schemes in its stacky formulation from \cite[§3]{ArithEtaleComp}, where the occurrence of the Hyodo--Kato stack above is replaced by the algebraic syntomification of the special fibre.

In particular, the equivalence (\ref{eq:hkcomp-perffibreproduct}) includes the statement that the syntomic cohomology of $X$ with coefficients in a perfect analytic $F$-gauge may be computed as a pullback of the Hyodo--Kato cohomology and the Hodge-filtered de Rham cohomology over the full de Rham cohomology of $X$, each with suitable coefficients obtained from the perfect analytic $F$-gauge we started with through various realisation functors, which we shall introduce. Taking coefficients in a Breuil--Kisin twist and comparing with \cite[Thm.\ 1.1]{padicComparisons}, this in particular implies that our definition of $R\Gamma_\Syn(X, \Q_p(i))$ for smooth partially proper qcqs rigid spaces $X$ recovers the one of Colmez--Nizio{\l}.

Finally, applying the above discussion to the case $X=\GSpec\Q_p$, we will be able to identify vector bundles on $\Q_p^\Syn$ with $\Gal_{\Q_p}$-equivariant vector bundles on $\FF_{\C_p}$ which are de Rham in the sense of \cite[Def.\ 15.12]{FarguesFontaine}. Moreover, in the case where the de Rham bundle arises from a $\Q_p$-representation $V$ of $\Gal_{\Q_p}$ by extending scalars, we will see that
\begin{equation*}
H^1(\Q_p^\Syn, V)\cong H^1_g(\Gal_{\Q_p}, V)\;,
\end{equation*}
where the right-hand side was defined by Bloch--Kato in \cite{BlochKato} as the subspace of $H^1(\Gal_{\Q_p}, V)$ spanned by classes corresponding to extensions of $\Q_p$ by $V$ which are de Rham. 

One should compare this to similar of result of Bhatt--Lurie about the algebraic syntomification of $\Z_p$: In \cite[§§6.6, 6.7]{FGauges}, they prove that $\Z_p$-lattices in crystalline $\Gal_{\Q_p}$-representations are equivalent to the category of so-called reflexive coherent sheaves on the algebraic syntomification of $\Z_p$ and that this equivalence identifies $H^1$ on the syntomification with the Bloch--Kato--Selmer group $H^1_f$ after rationalisation, see \cite[Prop.\ 6.7.3]{FGauges}.

\subsection{Accessing $X^\Syn$ via $X^\HK$ and $X^{\dR, +}$}

We begin by defining the stack $X^\mHK$ and the map $i_\HK: X^\mHK\rightarrow X^\Syn$ occurring in the desired commutative diagram (\ref{eq:hkcomp-square}). This works for an arbitrary Gelfand stack $X$ and we start by introducing the following auxiliary object: For any Gelfand stack $X$, the Gelfand stack $X^{\mHK, \mathrm{pre}}$ is defined as the pushout
\begin{equation*}
\begin{tikzcd}
Y_{X^\diamond, [p^{3/2}, p^{3/2}]}^\dR\ar[r, "\can"]\ar[d, "\can\times \{0\}", swap] & Y_{X^\diamond, [p^{1/2}, p^{3/2}]}^\dR\ar[d] \\
Y_{X^\diamond, [p^{3/2}, \infty)}^\dR\times [0, 1]\ar[r] & X^{\mHK, \mathrm{pre}}\nospacepunct{\;.}
\end{tikzcd}
\end{equation*}
We now obtain $X^\mHK$ from $X^{\mHK, \mathrm{pre}}$ by gluing. For this, note that the map
\begin{equation*}
Y_{X^\diamond, [p^{3/2}, \infty)}^\dR\xrightarrow{\id\times\{0\}} Y_{X^\diamond, [p^{3/2}, \infty)}^\dR\times [0, 1]
\end{equation*}
and the identity on $Y_{X^\diamond, [p^{1/2}, p^{3/2}]}^\dR$ glue to a map
\begin{equation}
\label{eq:hkcomp-idtimes0}
\id\times\{0\}: Y_{X^\diamond, [p^{1/2}, p^{3/2}]}^\dR\coprod_{Y_{X^\diamond, [p^{3/2}, p^{3/2}]}^\dR} Y_{X^\diamond, [p^{3/2}, \infty)}^\dR\rightarrow X^{\mHK, \mathrm{pre}}\;,
\end{equation}
the source of which is isomorphic to $Y_{X^\diamond, [p^{1/2}, \infty)}^\dR$ by \cref{lem:hkcomp-glueoverinterval}.

\begin{defi}
Let $X$ be any Gelfand stack. The \emph{mock Hyodo--Kato stack} $X^\mHK$ of $X$ is defined by the coequaliser diagram
\begin{equation*}
\begin{tikzcd}
Y_{X^\diamond, [p^{1/2}, \infty)}^\dR\ar[r,shift left=.75ex,"\id\times\{0\}"]\ar[r,shift right=.75ex,swap,"\phi\times\{1\}"] & X^{\mHK, \mathrm{pre}}\ar[r] & X^\mHK \nospacepunct{\;.}
\end{tikzcd}
\end{equation*}
\end{defi}

To obtain the map $i_\HK: X^\mHK\rightarrow X^\Syn$, recall that
\begin{equation*}
X^\N_{[p^{3/2}, \infty)}\cong Y_{X^\diamond, [p^{3/2}, \infty)}^\dR\times [0, 1]
\end{equation*}
via $\pi$ by \cref{prop:defis-utneq0} and \cref{prop:defis-xprismdr}. Moreover, recalling that $X^\N_{|t|=1}\cong X^\prism$ via $\pi$, another application of \cref{prop:defis-xprismdr} yields 
\begin{equation*}
X^\N_{|t|=1, [p^{1/2}, p^{3/2}]}\cong Y_{X^\diamond, [p^{1/2}, p^{3/2}]}^\dR
\end{equation*}
via $\pi$. Gluing these two isomorphisms along $Y_{X^\diamond, [p^{3/2}, p^{3/2}]}^\dR$, we obtain a map
\begin{equation*}
X^{\mHK, \mathrm{pre}}\rightarrow X^\N
\end{equation*}
which is an isomorphism onto its image, see \cref{fig:xmhkpre}. Finally, this will descend to the desired map
\begin{equation*}
i_\HK: X^\mHK\rightarrow X^\Syn\;,
\end{equation*}
which we will call the \emph{Hyodo--Kato map}. 

To get the desired diagram (\ref{eq:hkcomp-square}), recall that the map $\phi\circ\iota: \widehat{X}\rightarrow Y_{X^\diamond}$ factors over the preimage of $\{p\}\subseteq (0, \infty)$ under the radius map and hence induces a map $X^\dR\rightarrow Y_{X^\diamond, [p^{1/2}, p^{3/2}]}^\dR$, which in turn yields a composite map
\begin{equation}
\label{eq:hkcomp-mapxdrxmhk}
X^\dR\rightarrow Y_{X^\diamond, [p^{1/2}, p^{3/2}]}^\dR\rightarrow X^{\mHK, \mathrm{pre}}\rightarrow X^\mHK\;.
\end{equation}
Observing that the postcomposition of this map with $i_\HK$ identifies with the precomposition of $i_{\dR, +}$ with the canonical map $X^\dR\rightarrow X^{\dR, +}$, we obtain the desired commutative diagram
\begin{equation*}
\begin{tikzcd}
X^\dR\ar[r]\ar[d] & X^{\dR, +}\ar[d, "i_{\dR, +}"] \\
X^{\mHK}\ar[r, "i_\HK"] & X^\Syn\nospacepunct{\;.}
\end{tikzcd}
\end{equation*}

For this to be of any use to us in the sequel, we of course need to relate $X^\mHK$ and, in particular, the category of perfect complexes on $X^\mHK$, to the actual Hyodo--Kato stack $X^\HK$ and perfect complexes thereon.

\begin{figure}

\begin{center}
\begin{tikzpicture}
  \def\width{5.5} 
  \def\height{1.2*\width} 

  \draw[dotted, thick] (0, 0) -- (\width, 0); 

  

  
  \draw[thick] (\width-0.2, 0.2*\height) -- (\width+0.2, 0.2*\height); 
  \node at (\width+1.2, 0.2*\height) {$\phi^{-1}(X^\HT)$};
  \draw[thick] (\width-0.2, 0.4*\height) -- (\width+0.2, 0.4*\height); 
  \node at (\width+0.7, 0.4*\height) {$X^\HT$};
  \draw[thick] (\width-0.2, 0.6*\height) -- (\width+0.2, 0.6*\height); 
  \node at (\width+0.7, 0.6*\height) {$X^\dR$};
  \draw[thick] (\width-0.2, 0.8*\height) -- (\width+0.2, 0.8*\height); 
  \node at (\width+0.95, 0.8*\height) {$\phi(X^\dR)$};

  \draw[thick] (-0.2, 0.2*\height) -- (0.2, 0.2*\height); 
  \node at (-0.7, 0.2*\height) {$X^\HT$};
  \draw[thick] (-0.2, 0.4*\height) -- (0.2, 0.4*\height); 
  \node at (-0.7, 0.4*\height) {$X^\dR$};
  \draw[thick] (-0.2, 0.6*\height) -- (0.2, 0.6*\height); 
  \node at (-0.95, 0.6*\height) {$\phi(X^\dR)$};
  \draw[thick] (-0.2, 0.8*\height) -- (0.2, 0.8*\height); 
  \node at (-1.05, 0.8*\height) {$\phi^2(X^\dR)$};


  \draw[line width=0.2mm, color=blue] (\width/2-0.05*\width, 0.4*\height) -- (0, 0.4*\height);

  \def\sqr_size{0.1*\width} 

  \draw[line width=0.2mm] (\width/2 - \sqr_size/2, 0.4*\height - \sqr_size/2) -- (\width/2 + \sqr_size/2, 0.4*\height + \sqr_size/2); 
  \draw[line width=0.2mm] (\width/2 + \sqr_size/2, 0.4*\height - \sqr_size/2) -- (\width/2 - \sqr_size/2, 0.4*\height + \sqr_size/2); 

  
  \shade[shading=axis, bottom color=white, top color=white, middle color=green, shading angle=0](0.5*\width + \sqr_size/2, 0.4*\height - 0.1) rectangle (\width, 0.4*\height + 0.1);
  \shade[shading=axis, bottom color=white, top color=white, middle color=green, shading angle=0](0, 0.2*\height - 0.1) rectangle (\width, 0.2*\height + 0.1);

  \draw[->, thick] (2*\width, 0) -- (2*\width, \height);
  
  \draw[thick] (2*\width-0.2, 0) -- (2*\width+0.2, 0); 
  \node at (2*\width+0.6, 0) {$0$};
  \draw[thick] (2*\width-0.2, 0.2*\height) -- (2*\width+0.2, 0.2*\height); 
  \node at (2*\width+0.6, 0.2*\height) {$1$};
  \draw[thick] (2*\width-0.2, 0.4*\height) -- (2*\width+0.2, 0.4*\height); 
  \node at (2*\width+0.6, 0.4*\height) {$p$};
  \draw[thick] (2*\width-0.2, 0.6*\height) -- (2*\width+0.2, 0.6*\height); 
  \node at (2*\width+0.6, 0.6*\height) {$p^2$};
  \draw[thick] (2*\width-0.2, 0.8*\height) -- (2*\width+0.2, 0.8*\height); 
  \node at (2*\width+0.6, 0.8*\height) {$p^3$};
  
  \draw[->, thick] (1.3*\width, \height/2) -- (1.8*\width, \height/2);
  \node at (1.55*\width, \height/2+0.3) {$\kappa$};
  
  \draw[line width=1mm] (0, 0.3*\height) -- (0, \height);
  \draw[line width=1mm] (0, 0.5*\height) -- (\width+0.05, 0.5*\height);
  \draw[line width=1mm] (\width, 0.5*\height) -- (\width, \height);
  
  \fill[color=gray!50] (0.05, 0.5*\height+0.05) rectangle (\width-0.05, \height);
  
  \draw[->, thick] (-0.5*\width, 0.3*\height) -- (-0.5*\width, 0.7*\height);
  \node at (-0.5*\width-0.4, 0.5*\height) {$\phi$};
  
  \draw[line width=0.2mm, color=blue] (\width, 0.6*\height) -- (0, 0.6*\height);
  \draw[line width=0.2mm, color=blue] (\width, 0.8*\height) -- (0, 0.8*\height);
  
  \draw[thick] (0, 0) -- (0, \height); 
  \draw[thick] (\width, 0) -- (\width, \height); 
  
  \draw[line width=0.2mm] (\width/2 - \sqr_size/2, 0.4*\height - \sqr_size/2) rectangle (\width/2 + \sqr_size/2, 0.4*\height + \sqr_size/2); 

\end{tikzpicture}
\end{center}

\captionsetup{justification=centering}
\caption{A schematic picture of $X^\N$ with the image of $X^{\mHK, \mathrm{pre}}$ \\ outlined in bold (interior shaded in grey)}
\label{fig:xmhkpre}
\end{figure}

\begin{prop}
\label{prop:hkcomp-hkpshk}
Let $X$ be any Gelfand stack. Then there is an equivalence of categories
\begin{equation*}
\Perf(X^\mHK)\cong\Perf(X^\HK)
\end{equation*}
induced by the diagram
\begin{equation*}
\begin{tikzcd}
& Y_{X^\diamond, [p^{1/2}, \infty)}^\dR\ar[rd]\ar[ld, "\id\times \{0\}", swap] & \\
X^\mHK & & X^\HK\nospacepunct{\;.}
\end{tikzcd}
\end{equation*}
\end{prop}
\begin{proof}
By \cref{cor:perf-contractible}, we have
\begin{equation*}
\Perf(Y_{X^\diamond, [p^{3/2}, \infty)}^\dR\times [0, 1])\cong \Perf(Y_{X^\diamond, [p^{3/2}, \infty)}^\dR)
\end{equation*}
and hence the definition of $X^{\mHK, \mathrm{pre}}$ shows that pullback along the map
\begin{equation*}
\id\times\{0\}: Y_{X^\diamond, [p^{1/2}, p^{3/2}]}^\dR\coprod_{Y_{X^\diamond, [p^{3/2}, p^{3/2}]}^\dR} Y_{X^\diamond, [p^{3/2}, \infty)}^\dR\rightarrow X^{\mHK, \mathrm{pre}}
\end{equation*}
from (\ref{eq:hkcomp-idtimes0}) induces an equivalence on perfect complexes. As the pushout on the left is isomorphic to $Y_{X^\diamond, [p^{1/2}, \infty)}^\dR$ by \cref{lem:hkcomp-glueoverinterval}, we conclude that the categories of perfect complexes of $X^\mHK$ and
\begin{equation*}
\operatorname{coeq}(\hspace{-0.15cm}
\begin{tikzcd}
Y_{X^\diamond, [p^{1/2}, \infty)}^\dR\ar[r,shift left=.75ex, "\id"]
  \ar[r,shift right=.75ex,swap, "\phi"] & Y_{X^\diamond, [p^{1/2}, \infty)}^\dR
\end{tikzcd}
\hspace{-0.15cm})
\end{equation*}
are equivalent. However, the latter coequaliser is actually isomorphic to $X^\HK$: Indeed, since the radius map $\kappa: Y_{X^\diamond}^\dR\rightarrow (0, \infty)$ intertwines $\phi$ with multiplication by $p$, this follows by base changing the isomorphism
\begin{equation*}
\operatorname{coeq}(\hspace{-0.15cm}
\begin{tikzcd}
{[p^{1/2}, \infty)}\ar[r,shift left=.75ex, "\id"]
  \ar[r,shift right=.75ex,swap, "p"] & {[p^{1/2}, \infty)}
\end{tikzcd}
\hspace{-0.15cm})
\cong 
\operatorname{coeq}(\hspace{-0.15cm}
\begin{tikzcd}
{(0, \infty)}\ar[r,shift left=.75ex, "\id"]
  \ar[r,shift right=.75ex,swap, "p"] & {(0, \infty)}
\end{tikzcd}
\hspace{-0.15cm})
\end{equation*}
to $Y_{X^\diamond}^\dR$ along $\kappa$ using universality of colimits in $\infty$-topoi.
\end{proof}

Using \cref{prop:hkcomp-hkpshk}, the diagram (\ref{eq:hkcomp-square}) now induces various realisation functors for perfect analytic $F$-gauges. Firstly, pullback along the map $i_\HK$ induces a realisation functor
\begin{equation*}
T_\HK: \Perf(X^\Syn)\rightarrow\Perf(X^\HK)
\end{equation*}
from perfect analytic $F$-gauges to perfect complexes on $X^\HK$, which should be thought of as coefficients for Hyodo--Kato cohomology of $X$ following \cite[Rem.\ 6.1.3]{dRFF}. Secondly, via pullback along the filtered de Rham map $i_{\dR, +}$, we obtain a realisation functor
\begin{equation*}
T_{\dR, +}: \D(X^\Syn)\rightarrow \D(X^{\dR, +})
\end{equation*}
from analytic $F$-gauges to quasicoherent sheaves on $X^{\dR, +}$, i.e.\ a category of filtered analytic $D$-modules on $X$, and, lastly, postcomposing $T_{\dR, +}$ with pullback along the canonical map $X^\dR\rightarrow X^{\dR, +}$ yields a realisation functor
\begin{equation*}
T_\dR: \D(X^\Syn)\rightarrow \D(X^\dR)\;.
\end{equation*}
Note that pullback along the map $X^\dR\rightarrow X^\mHK$ from (\ref{eq:hkcomp-mapxdrxmhk}) also induces a functor $\Perf(X^\HK)\rightarrow\Perf(X^\dR)$ whose precomposition with $T_\HK$ identifies with $T_\dR$ by commutativity of the diagram (\ref{eq:hkcomp-square}). Using this terminology, we can formulate the main result of this section as follows:

\begin{thm}
\label{thm:hkcomp-main}
Let $X$ be a Berkovich smooth derived Berkovich space over $\Q_p$. Then there is an equivalence of categories
\begin{equation*}
\Perf(X^\Syn)\cong \Perf(X^\HK)\times_{\Perf(X^\dR)} \Perf(X^{\dR, +})
\end{equation*}
induced by the realisation functors $T_\HK$ and $T_{\dR, +}$. In particular, for any $E\in\Perf(X^\Syn)$, there is a pullback diagram
\begin{equation*}
\begin{tikzcd}
R\Gamma(X^\Syn, E)\ar[r]\ar[d] & R\Gamma(X^\HK, T_\HK(E))\ar[d] \\
R\Gamma(X^{\dR, +}, T_{\dR, +}(E))\ar[r] & R\Gamma(X^\dR, T_\dR(E))\nospacepunct{\;.}
\end{tikzcd}
\end{equation*}
\end{thm}

Using the expected relation between the stack $X^\HK$ and Hyodo--Kato cohomology of $X$, comparing the result above with \cite[Thm.\ 1.1]{padicComparisons} in particular shows that our definition of syntomic cohomology recovers the more classical definition of syntomic cohomology of Colmez--Nizio{\l} for smooth partially proper qcqs rigid spaces over $\Q_p$, see e.g.\ \cite[§4.4]{padicComparisons}. Indeed, applying the above result to $E=\O\{i\}$ and using \cref{prop:recall-fildr}, we obtain:

\begin{cor}
\label{cor:hkcomp-classical}
Let $X$ be a smooth partially proper qcqs rigid space over $\Q_p$. For any $i\in\Z$, there is a cartesian diagram
\begin{equation*}
\begin{tikzcd}
R\Gamma_\Syn(X, \Q_p(i))\ar[r]\ar[d] & R\Gamma_\HK(X)^{\phi=p^i, N=0, \Gal_{\Q_p}}\ar[d] \\
\Fil_\Hod^i R\Gamma_\dR(X)\ar[r] & R\Gamma_\dR(X)\nospacepunct{\;,}
\end{tikzcd}
\end{equation*}
where we note that the $\Gal_{\Q_p}$-invariants in the top right corner should be interpreted as smooth group cohomology and are hence \emph{underived}.
\end{cor}

Let us now move on to the proof of \cref{thm:hkcomp-main}. For this, we first record the following lemma, which says that, if we ``chop up'' $X^\N$ along the radius map, gluing the pieces back together recovers $X^\N$.

\begin{lem}
\label{lem:hkcomp-chopupxn}
Let $X$ be any Gelfand stack. The canonical inclusion maps induce an isomorphism between the iterated pushout
\begin{equation*}
X^\N_{[p^{1/2}, \infty)}\coprod_{X^\N_{[p^{1/2}, p^{1/2}]}} X^\N_{[p^{-1/2}, p^{1/2}]}\coprod_{X^\N_{[p^{-1/2}, p^{-1/2}]}} X^\N_{[p^{-3/2}, p^{-1/2}]} \coprod_{X^\N_{[p^{-3/2}, p^{-3/2}]}} \dots
\end{equation*}
and the stack $X^\N$.
\end{lem}
\begin{proof}
This follows by repeatedly applying \cref{lem:hkcomp-glueoverinterval}.
\end{proof}

\begin{figure}

\begin{center}
\begin{tikzpicture}
  \def\width{5.5} 
  \def\height{1.2*\width} 

  \draw[line width=1mm] (0, -0.05) -- (0, 0.6*\height); 
  \draw[line width=1mm] (\width, 0.2*\height) -- (2*\width, 0.2*\height); 
  \draw[line width=1mm] (\width, 0.2*\height-0.05) -- (\width, 0.6*\height); 
  \draw[line width=1mm] (0, 0) -- (2*\width, 0); 

  

  
  \draw[thick] (\width-0.2, 0.3*\height) -- (\width+0.2, 0.3*\height); 
  \node at (\width+0.95, 0.3*\height) {$X^\dR$};
  \draw[thick] (\width-0.2, 0.5*\height) -- (\width+0.2, 0.5*\height); 
  \node at (\width+1.05, 0.5*\height) {$\phi(X^\dR)$};

  \draw[thick] (-0.2, 0.1*\height) -- (0.2, 0.1*\height); 
  \node at (-0.95, 0.1*\height) {$X^\dR$};
  \draw[thick] (-0.2, 0.3*\height) -- (0.2, 0.3*\height); 
  \node at (-1.05, 0.3*\height) {$\phi(X^\dR)$};
  \draw[thick] (-0.2, 0.5*\height) -- (0.2, 0.5*\height); 
  \node at (-1.05, 0.5*\height) {$\phi^2(X^\dR)$};


  \draw[line width=0.2mm, color=blue] (\width/2-0.05*\width, 0.1*\height) -- (0, 0.1*\height);
  \draw[line width=0.2mm, color=blue] (\width, 0.3*\height) -- (0, 0.3*\height);
  \draw[line width=0.2mm, color=blue] (\width, 0.5*\height) -- (0, 0.5*\height);

  \def\sqr_size{0.1*\width} 

  \draw[line width=0.2mm] (\width/2 - \sqr_size/2, 0.1*\height - \sqr_size/2) -- (\width/2 + \sqr_size/2, 0.1*\height + \sqr_size/2); 
  \draw[line width=0.2mm] (\width/2 + \sqr_size/2, 0.1*\height - \sqr_size/2) -- (\width/2 - \sqr_size/2, 0.1*\height + \sqr_size/2); 

  
  \shade[shading=axis, bottom color=white, top color=white, middle color=green, shading angle=0](0.5*\width + \sqr_size/2, 0.1*\height - 0.1) rectangle (2*\width, 0.1*\height + 0.1);

  \draw[line width=0.2mm] (\width/2 - \sqr_size/2, 0.1*\height - \sqr_size/2) rectangle (\width/2 + \sqr_size/2, 0.1*\height + \sqr_size/2); 

\end{tikzpicture}
\end{center}

\captionsetup{justification=centering}
\caption{A schematic picture of $X^{\Syn, \mathrm{pre}, \HK}$ with \\ the two copies of $\partial X^{\Syn, \mathrm{pre}, \HK}$ in bold}
\label{fig:xsynprehk}
\end{figure}

To move on, we introduce the following notation: For any Gelfand stack $X$, let $X^{\Syn, \mathrm{pre}, \HK}$ be defined as the pushout
\begin{equation*}
\begin{tikzcd}
X^\prism_{[p^{-1/2}, p^{1/2}]}\ar[r, "\id\times\{0\}"]\ar[d, "j_\HT", swap] & X^\prism_{[p^{-1/2}, p^{1/2}]}\times [0, \infty)\ar[d] \\
X^\N_{[p^{1/2}, \infty)}\ar[r] & X^{\Syn, \mathrm{pre}, \HK}\nospacepunct{\;,}
\end{tikzcd}
\end{equation*}
see \cref{fig:xsynprehk}. Then $X^{\Syn, \mathrm{pre}, \HK}$ receives two maps from the pushout
\begin{equation*}
\begin{tikzcd}
X^\prism_{[p^{1/2}, p^{1/2}]}\ar[r, "\id\times\{0\}"]\ar[d, "\can", swap] & X^\prism_{[p^{1/2}, p^{1/2}]}\times [0, \infty)\ar[d] \\
X^\prism_{[p^{1/2}, \infty)}\ar[r] & \partial X^{\Syn, \mathrm{pre}, \HK}\nospacepunct{\;.}
\end{tikzcd}
\end{equation*}
Namely, one of them is induced by $j_\HT: X^\prism_{[p^{1/2}, \infty)}\rightarrow X^\N_{[p^{1/2}, \infty)}$ and the canonical inclusion $X^\prism_{[p^{1/2}, p^{1/2}]}\rightarrow X^\prism_{[p^{-1/2}, p^{1/2}]}$ and we denote this map by
\begin{equation*}
j_\HT: \partial X^{\Syn, \mathrm{pre}, \HK}\rightarrow X^{\Syn, \mathrm{pre}, \HK}
\end{equation*}
as well. The other one is induced by $j_\dR: X^\prism_{[p^{1/2}, \infty)}\rightarrow X^\N_{[p^{1/2}, \infty)}$ and the isomorphism
\begin{equation}
\label{eq:hkcomp-jdrxsynpre}
X^\prism_{[p^{1/2}, p^{1/2}]}\times [0, \infty)\xrightarrow{\cong} (X^\prism_{[p^{1/2}, p^{1/2}]}\times [0, 1])\coprod_{X^\prism_{[p^{-1/2}, p^{-1/2}]}} (X^\prism_{[p^{-1/2}, p^{-1/2}]}\times [0, \infty))\;,
\end{equation}
where the maps defining the pushout on the right-hand side are $\phi\times \{1\}$ ``to the left'' and $\id\times \{0\}$ ``to the right'', respectively, and we use the isomorphism
\begin{equation*}
X^\N_{[p^{1/2}, p^{1/2}]}\cong X^\prism_{[p^{1/2}, p^{1/2}]}\times [0, 1]
\end{equation*}
from \cref{prop:defis-utneq0} to view the target of (\ref{eq:hkcomp-jdrxsynpre}) as a closed substack of $X^{\Syn, \mathrm{pre}, \HK}$. We denote this second map by
\begin{equation*}
j_\dR: \partial X^{\Syn, \mathrm{pre}, \HK}\rightarrow X^{\Syn, \mathrm{pre}, \HK}\;.
\end{equation*}
In this notation, our next lemma reads as follows:

\begin{figure}

\begin{center}
\begin{tikzpicture}
  \def\width{5.5} 
  \def\height{1*\width} 

  \draw[line width=1mm, color=red] (0, 0) -- (0, 0.1*\height); 
  \draw[thick] (\width, 0) -- (\width, 0.1*\height); 
  \draw[thick] (0, 0.1*\height) -- (\width, 0.1*\height);
  \draw[dotted, thick] (0, 0) -- (\width, 0); 

  

  
  \draw[thick] (\width-0.2, 0.3*\height) -- (\width+0.2, 0.3*\height); 
  \node at (\width+1.2, 0.3*\height) {$\phi^{-2}(X^\HT)$};
  \draw[thick] (\width-0.2, 0.6*\height) -- (\width+0.2, 0.6*\height); 
  \node at (\width+1.2, 0.6*\height) {$\phi^{-1}(X^\HT)$};
  \draw[thick] (\width-0.2, 0.9*\height) -- (\width+0.2, 0.9*\height); 
  \node at (\width+0.7, 0.9*\height) {$X^\HT$};
  \draw[thick] (\width-0.2, 1.1*\height) -- (\width+0.2, 1.1*\height); 
  \node at (\width+0.95, 1.1*\height) {$X^\dR$};
  \draw[thick] (\width-0.2, 1.3*\height) -- (\width+0.2, 1.3*\height); 
  \node at (\width+1.05, 1.3*\height) {$\phi(X^\dR)$};

  \draw[thick] (-0.2, 0.3*\height) -- (0.2, 0.3*\height); 
  \node at (-1.2, 0.3*\height) {$\phi^{-1}(X^\HT)$};
  \draw[thick] (-0.2, 0.6*\height) -- (0.2, 0.6*\height); 
  \node at (-0.7, 0.6*\height) {$X^\HT$};
  \draw[thick] (-0.2, 0.9*\height) -- (0.2, 0.9*\height); 
  \node at (-0.95, 0.9*\height) {$X^\dR$};
  \draw[thick] (-0.2, 1.1*\height) -- (0.2, 1.1*\height); 
  \node at (-1.05, 1.1*\height) {$\phi(X^\dR)$};
  \draw[thick] (-0.2, 1.3*\height) -- (0.2, 1.3*\height); 
  \node at (-1.05, 1.3*\height) {$\phi^2(X^\dR)$};


  \draw[line width=0.2mm, color=blue] (\width/2-0.05*\width, 0.9*\height) -- (0, 0.9*\height);
  \draw[line width=0.2mm, color=blue] (\width, 1.1*\height) -- (0, 1.1*\height);
  \draw[line width=0.2mm, color=blue] (\width, 1.3*\height) -- (0, 1.3*\height);

  \def\sqr_size{0.1*\width} 

  \draw[line width=0.2mm] (\width/2 - \sqr_size/2, 0.9*\height - \sqr_size/2) -- (\width/2 + \sqr_size/2, 0.9*\height + \sqr_size/2); 
  \draw[line width=0.2mm] (\width/2 + \sqr_size/2, 0.9*\height - \sqr_size/2) -- (\width/2 - \sqr_size/2, 0.9*\height + \sqr_size/2); 

  
  \shade[shading=axis, bottom color=white, top color=white, middle color=green, shading angle=0](0.5*\width + \sqr_size/2, 0.9*\height - 0.1) rectangle (\width-0.05, 0.9*\height + 0.1);
  \shade[shading=axis, bottom color=white, top color=white, middle color=green, shading angle=0](0.05, 0.6*\height - 0.1) rectangle (\width-0.05, 0.6*\height + 0.1);
  \shade[shading=axis, bottom color=white, top color=white, middle color=green, shading angle=0](0.05, 0.3*\height - 0.1) rectangle (\width-0.05, 0.3*\height + 0.1);
  
  
  \draw[line width=1mm, color=orange] (0, 0.2*\height) -- (0, 0.4*\height);  
  \draw[line width=1mm, color=red] (\width, 0.2*\height) -- (\width, 0.4*\height);
  \draw[thick] (0, 0.2*\height) -- (\width, 0.2*\height);
  \draw[thick] (0, 0.4*\height) -- (\width, 0.4*\height);
  
  \draw[line width=1mm, color=yellow] (0, 0.5*\height) -- (0, 0.7*\height);  
  \draw[line width=1mm, color=orange] (\width, 0.5*\height) -- (\width, 0.7*\height);
  \draw[thick] (0, 0.5*\height) -- (\width, 0.5*\height);
  \draw[thick] (0, 0.7*\height) -- (\width, 0.7*\height);
  
  \draw[thick] (0, 0.8*\height) -- (0, 1.4*\height);  
  \draw[line width=1mm, color=yellow] (\width, 0.8*\height) -- (\width, \height);
  \draw[thick] (\width, 1.0*\height) -- (\width, 1.4*\height);
  \draw[thick] (0, 0.8*\height) -- (\width, 0.8*\height);

  \draw[line width=0.2mm] (\width/2 - \sqr_size/2, 0.9*\height - \sqr_size/2) rectangle (\width/2 + \sqr_size/2, 0.9*\height + \sqr_size/2); 

\end{tikzpicture}
\end{center}

\captionsetup{justification=centering}
\caption{Proof of \cref{lem:hkcomp-syndifferentgluing}: chop up $X^\N$, then glue along \\ the yellow, orange and red pieces to get $X^{\Syn, \mathrm{pre}, \HK}$}
\label{fig:xsynviaxsynprehk}
\end{figure}

\begin{lem}
\label{lem:hkcomp-syndifferentgluing}
Let $X$ be any Gelfand stack. Then $X^\Syn$ is isomorphic to the coequaliser of the two maps
\begin{equation*}
\begin{tikzcd}
\partial X^{\Syn, \mathrm{pre}, \HK}\ar[r,shift left=.75ex, "j_\HT"]
  \ar[r,shift right=.75ex,swap, "j_\dR"] & X^{\Syn, \mathrm{pre}, \HK}\nospacepunct{\;.}
\end{tikzcd}
\end{equation*}
\end{lem}
\begin{proof}
See also \cref{fig:xsynviaxsynprehk}. By definition, $X^\Syn$ is obtained as the coequaliser of the maps $j_\HT, j_\dR: X^\prism\rightarrow X^\N$. However, by \cref{lem:hkcomp-chopupxn}, we can view $X^\N$ and $X^\prism$ as glued together from the smaller pieces 
\begin{equation*}
X^\N_{[p^{1/2}, \infty)}, X^\N_{[p^{-1/2}, p^{1/2}]}, \dots\;, \hspace{0.3cm}\text{and}\hspace{0.3cm}X^\prism_{[p^{1/2}, \infty)}, X^\prism_{[p^{-1/2}, p^{1/2}]}, \dots\;,
\end{equation*}
respectively. To obtain $X^\Syn$ from $X^\N_{[p^{1/2}, \infty)}, X^\N_{[p^{-1/2}, p^{1/2}]}, \dots$, we have to perform the following gluings:
\begin{enumerate}[label=(\arabic*)]
\item For each $k\geq -1$ odd, glue the two copies of $X^\N_{[p^{-k/2}, p^{-k/2}]}$ inside $X^\N_{[p^{-k/2-1}, p^{-k/2}]}$ and $X^\N_{[p^{-k/2}, p^{-k/2+1}]}$ (or $X^\N_{[p^{1/2}, \infty)}$ if $k=-1$), respectively;
\item for each $k\geq 1$ odd, glue the two copies of $X^\prism_{[p^{-k/2-1}, p^{-k/2}]}$ embedded into $X^\N_{[p^{-k/2}, p^{-k/2+1}]}$ and $X^\N_{[p^{-k/2-1}, p^{-k/2}]}$ via $j_\HT$ and $j_\dR$, respectively;
\item glue the two copies of $X^\prism_{[p^{-1/2}, p^{1/2}]}$ which are embedded into $X^\N_{[p^{1/2}, \infty)}$ and $X^\N_{[p^{-1/2}, p^{1/2}]}$ via $j_\HT$ and $j_\dR$, respectively;
\item glue the two copies of $X^\prism_{[p^{1/2}, \infty)}$ which are embedded into $X^\N_{[p^{1/2}, \infty)}$ via $j_\HT$ and $j_\dR$, respectively.
\end{enumerate}

We now change the order of these gluings: For this, first note that
\begin{equation*}
X^\N_{[p^{-k/2}, p^{-k/2+1}]}\cong X^\prism_{[p^{-k/2}, p^{-k/2+1}]}\times [0, 1]\cong X^\prism_{[p^{-1/2}, p^{1/2}]}\times [0, 1]
\end{equation*}
for each $k\geq 1$ odd, where the first isomorphism is due to \cref{prop:defis-utneq0} and the second one is via a $(k-1)/2$-fold application of $\phi$, which we recall defines an isomorphism $X^\prism_{(0, 1)}\cong X^\prism_{(0, p)}$ by \cref{prop:defis-xprismxdiv1}. Thus, first performing the gluings (2) and (3), we obtain $X^{\Syn, \mathrm{pre}, \HK}$ as defined above.

Observing that
\begin{equation*}
X^\N_{[p^{-k/2}, p^{-k/2}]}\cong X^\prism_{[p^{-k/2}, p^{-k/2}]}\times [0, 1]\cong X^\prism_{[p^{1/2}, p^{1/2}]}\times [0, 1]
\end{equation*}
by another application of \cref{prop:defis-utneq0} and \cref{prop:defis-xprismxdiv1}, we now see that subsequently performing the gluings (1) and (4) equivalently amounts to gluing the two copies of $\partial X^{\Syn, \mathrm{pre}, \HK}$ embedded into $X^{\Syn, \mathrm{pre}, \HK}$ via $j_\HT$ and $j_\dR$. This proves the claim.
\end{proof}

\begin{lem}
\label{lem:hkcomp-xsynviarestrxn}
For any Gelfand stack $X$, pullback along the morphism
\begin{equation*}
\operatorname{coeq}(\hspace{-0.15cm}
\begin{tikzcd}
X^\prism_{[p^{1/2}, \infty)}\ar[r,shift left=.75ex, "j_\HT"]
  \ar[r,shift right=.75ex,swap, "j_\dR"] & X^\N_{[p^{1/2}, \infty)}
\end{tikzcd}
\hspace{-0.15cm})\rightarrow X^\Syn
\end{equation*}
induces an equivalence on perfect complexes.
\end{lem}
\begin{proof}
This immediately follows from the previous lemma once we know that
\begin{equation*}
\Perf(X^{\Syn, \mathrm{pre}, \HK})\cong \Perf(X^\N_{[p^{1/2}, \infty)}) \hspace{0.3cm}\text{and}\hspace{0.3cm} \Perf(\partial X^{\Syn, \mathrm{pre}, \HK})\cong \Perf(X^\prism_{[p^{1/2}, \infty)})
\end{equation*}
via pullback along the natural maps $X^\N_{[p^{1/2}, \infty)}\rightarrow X^{\Syn, \mathrm{pre}, \HK}$ and $X^\prism_{[p^{1/2}, \infty)}\rightarrow \partial X^{\Syn, \mathrm{pre}, \HK}$. However, by definition of $X^{\Syn, \mathrm{pre}, \HK}$, the first equivalence follows from the fact that pullback along $\id\times\{0\}$ induces an equivalence
\begin{equation*}
\Perf(X^\prism_{[p^{-1/2}, p^{1/2}]}\times [0, \infty))\cong \Perf(X^\prism_{[p^{-1/2}, p^{1/2}]})\;,
\end{equation*}
which is a consequence of \cref{cor:perf-contractible}. Similarly, the second one follows from the fact that pullback along $\id\times\{0\}$ induces an equivalence
\begin{equation*}
\Perf(X^\prism_{[p^{1/2}, p^{1/2}]}\times [0, \infty))\cong \Perf(X^\prism_{[p^{1/2}, p^{1/2}]})\;,
\end{equation*}
which is proved by the same argument.
\end{proof}

We now come to the key lemma. This uses the spreading out results for perfect complexes we have proved in §\ref{subsect:spreading} as well as the nice coverability results from §§\ref{subsect:nicecoverqpn}, \ref{subsect:htdr}. Note that these were already key to the proof of \cref{prop:htdr-perfequiv}, which we will make use of as well.

\begin{lem}
\label{lem:hkcomp-keylem}
Let $X$ be a Berkovich smooth derived Berkovich space over $\Q_p$. Then pullback along the map
\begin{equation*}
X^\prism_{[p^{1/2}, p^{3/2}]}\coprod_{X^\dR} X^{\dR, +}\rightarrow X^\N_{[p^{1/2}, p^{3/2}]}\;,
\end{equation*}
which is induced by $j_\dR$ and $i_{\dR, +}$, induces an equivalence on perfect complexes.
\end{lem}
\begin{proof}
First note that we may replace the interval $[p^{1/2}, p^{3/2}]$ by an interval $[r, s]\subseteq (p^{1/2}, p^{3/2})$ with $r-p^{1/2}$ and $p^{3/2}-s$ sufficiently small. Indeed, by \cref{lem:hkcomp-glueoverinterval}, we have 
\begin{equation*}
X^\N_{[p^{1/2}, p^{3/2}]}\cong X^\N_{[p^{1/2}, r]}\coprod_{X^\N_{[r, r]}} X^\N_{[r, s]}\coprod_{X^\N_{[s, s]}} X^\N_{[s, p^{3/2}]}
\end{equation*}
and the same is true for $(-)^\prism$ in place of $(-)^\N$. Then recalling from \cref{prop:defis-utneq0} that 
\begin{equation*}
X^\N_{[r', s']}\cong X^\prism_{[r', s']}\times [0, 1]
\end{equation*}
for $[r', s']$ being any of the intervals $[p^{1/2}, r], [r, r], [s, s]$ or $[s, p^{3/2}]$, we obtain
\begin{equation*}
\Perf(X^\N_{[r', s']})\cong \Perf(X^\prism_{[r', s']})
\end{equation*}
via pullback along $j_\dR$ from \cref{lem:perf-betti}, which yields the desired reduction.

Next, note that, by compatibility of $X\mapsto X^\N$ with rational localisations, see \cref{prop:defis-openloc}, we may assume that $X$ admits a finite étale map $X\rightarrow Z$ to some rational localisation $Z$ of $\ol{\T}^n$ for some $n\geq 0$. Then $X^\N_{[r, s]}$ is nicely coverable by \cref{prop:perf-coversmoothrigid} and our task is to show that the diagram
\begin{equation}
\label{eq:hkcomp-keylemcartesian}
\begin{tikzcd}
\Perf(X^\N_{[r, s]})\ar[r, "j_\dR^*"]\ar[d, "i_{\dR, +}^*", swap] & \Perf(X^\prism_{[r, s]})\ar[d, "i_\dR^*"] \\
\Perf(X^{\dR, +})\ar[r] & \Perf(X^\dR)
\end{tikzcd}
\end{equation}
is cartesian.

To see this, consider the overconvergent normed divisor $Z\coloneqq\{|ut|=0\}\subseteq X^\N_{[r, s]}$ and note that, for each $\epsilon>0$ small enough, we have $Z_\epsilon\setminus Z\cong (X^\prism_{|\widetilde{\mu}|\leq\epsilon}\setminus X^\dR)\times [0, 1]$ by \cref{prop:defis-utneq0}, whence \cref{lem:perf-betti} yields
\begin{equation*}
\Perf(Z_\epsilon\setminus Z)\cong \Perf(X^\prism_{|\widetilde{\mu}|\leq\epsilon}\setminus X^\dR)\;.
\end{equation*}
Similarly, we obtain
\begin{equation*}
\Perf(X^\N_{[r, s]}\setminus Z)\cong \Perf(X^\prism_{[r, s]}\setminus X^\dR)\;.
\end{equation*}
Finally, from \cref{prop:htdr-perfequiv}, we know that
\begin{equation*}
\Perf(Z)\cong \Perf(X^{\dR, +})\;.
\end{equation*}

Applying \cref{cor:perf-corkeylemma} to the overconvergent normed divisor $Z$, the results of the previous paragraph yield a cartesian diagram
\begin{equation*}
\begin{tikzcd}
\Perf(X^\N_{[r, s]})\ar[r, "j_\dR^*"]\ar[d, "i_{\dR, +}^*", swap] & \Perf(X^\prism_{[r, s]}\setminus X^\dR)\ar[d] \\
\Perf(X^{\dR, +})\ar[r] & \colim_{\epsilon>0} \Perf(X^\prism_{|\widetilde{\mu}|\leq\epsilon}\setminus X^\dR)\;.
\end{tikzcd}
\end{equation*}
Putting this together with the cartesian diagram
\begin{equation*}
\begin{tikzcd}
\Perf(X^\prism_{[r, s]})\ar[r]\ar[d, "i_\dR^*", swap] & \Perf(X^\prism_{[r, s]}\setminus X^\dR)\ar[d] \\
\Perf(X^\dR)\ar[r] & \colim_{\epsilon>0} \Perf(X^\prism_{|\widetilde{\mu}|\leq\epsilon}\setminus X^\dR)
\end{tikzcd}
\end{equation*}
obtained by applying \cref{cor:perf-corkeylemma} to the overconvergent normed divisor $X^\dR\cong (\{|\widetilde{\mu}|=0\}\subseteq X^\prism_{[r, s]})$, where the isomorphism is via $i_\dR$, then finally yields the desired cartesian diagram (\ref{eq:hkcomp-keylemcartesian}). We point out that $X^\prism_{[r, s]}$ is indeed nicely coverable: This is deduced from nice coverability of $X^\N_{[r, s]}$ by base change using \cref{lem:perf-bc} and \cref{lem:perf-finflatdim}.
\end{proof}

Finally, we can put everything together and prove \cref{thm:hkcomp-main}.

\begin{proof}[Proof of \cref{thm:hkcomp-main}]
First recall from \cref{prop:defis-utneq0} that $X^\N_{[p^{3/2}, p^{3/2}]}\cong X^\prism_{[p^{3/2}, p^{3/2}]}\times [0, 1]$, which yields
\begin{equation*}
\Perf(X^\N_{[p^{3/2}, p^{3/2}]})\cong \Perf(X^\prism_{[p^{3/2}, p^{3/2}]})
\end{equation*}
via pullback along $j_\dR$ by \cref{lem:perf-betti}. Using \cref{lem:hkcomp-glueoverinterval}, we hence deduce from \cref{lem:hkcomp-keylem} that pullback along the map
\begin{equation}
\label{eq:hkcomp-xnphalfviaxmhkpre}
X^\N_{[p^{3/2}, \infty)}\coprod_{X^\prism_{[p^{3/2}, p^{3/2}]}} X^\prism_{[p^{1/2}, p^{3/2}]}\coprod_{X^\dR} X^{\dR, +}\rightarrow X^\N_{[p^{1/2}, \infty)}
\end{equation}
induces an equivalence on perfect complexes, where the map ``to the left'' in the leftmost pushout is induced by $j_\dR$. Recalling that $X^\prism_{[p^{3/2}, p^{3/2}]}\cong Y_{X^\diamond, [p^{3/2}, p^{3/2}]}^\dR$ and 
\begin{equation*}
X^\N_{[p^{3/2}, \infty)}\cong Y_{X^\diamond, [p^{3/2}, \infty)}^\dR\times [0, 1]\;,
\end{equation*}
which follow from \cref{prop:defis-prismffdr} and \cref{prop:defis-utneq0}, we recognise the source of (\ref{eq:hkcomp-xnphalfviaxmhkpre}) as the pushout of $X^{\mHK, \mathrm{pre}}$ along $X^\dR\rightarrow X^{\dR, +}$. Furthermore observing that $X^\prism_{[p^{1/2}, \infty)}\cong Y_{X^\diamond, [p^{1/2}, \infty)}^\dR$, which is due to \cref{prop:defis-prismffdr} as well, we deduce from (\ref{eq:hkcomp-xnphalfviaxmhkpre}) and \cref{lem:hkcomp-xsynviarestrxn} that pullback along
\begin{equation*}
\operatorname{coeq}(\hspace{-0.15cm}
\begin{tikzcd}
Y_{X^\diamond, [p^{1/2}, \infty)}^\dR\ar[r,shift left=.75ex, "\id\times \{0\}"]
  \ar[r,shift right=.75ex,swap, "\phi\times \{1\}"] & X^{\mHK, \mathrm{pre}}
\end{tikzcd}
\hspace{-0.15cm})\coprod_{X^\dR} X^{\dR, +}\rightarrow X^\Syn
\end{equation*}
induces an equivalence on perfect complexes as well, where we note that we have switched the order of the colimits on the left-hand side. However, by definition of $X^\mHK$, the source of this last map is exactly
\begin{equation*}
X^\mHK\coprod_{X^\dR} X^{\dR, +}
\end{equation*}
and this concludes the proof.
\end{proof}

\subsection{Bloch--Kato's $H^1_g$ via the syntomification}

We now specialise to the case $X=\GSpec\Q_p$ and use \cref{thm:hkcomp-main} to explicitly describe vector bundles on $\Q_p^\Syn$ and their cohomology in terms of $\Gal_{\Q_p}$-equivariant vector bundles on the Fargues--Fontaine curve.

\begin{thm}
\label{thm:bk-mainqp}
There is an equivalence of categories
\begin{equation*}
\Vect(\Q_p^\Syn)\cong \Vect^\dR(\FF_{\Spd\Q_p})
\end{equation*}
between vector bundles on $\Q_p^\Syn$ and vector bundles on the Fargues--Fontaine curve of $\Spd\Q_p$ which are de Rham in the sense of \cite[Def.\ 15.12]{FarguesFontaine}. Moreover, if $V$ and $W$ are de Rham representations of $\Gal_{\Q_p}$, then
\begin{equation*}
\RHom_{\Q_p^\Syn}(V, W)\cong \RHom_{\Rep^\dR(\Gal_{\Q_p})}(V, W)\;.
\end{equation*}
In particular, we have
\begin{equation*}
H^1(\Q_p^\Syn, V)\cong H^1_g(\Gal_{\Q_p}, V)\;.
\end{equation*}
\end{thm}
\begin{proof}
By \cite[Thm.\ 7.1.1]{dRFF}, the category of vector bundles of $\Q_p^\HK$ is equivalent to the category of $(\phi, N, \Gal_{\Q_p})$-modules over $\Q_p^\un$. Using that vector bundles on $\Q_p^{\dR, +}\cong \ol{\DD}/\ol{\T}$ are equivalent to finitely filtered $\Q_p$-vector spaces by \cref{prop:recall-rees}, we conclude from \cref{thm:hkcomp-main} that vector bundles on $\Q_p^\Syn$ are equivalent to filtered $(\phi, N, \Gal_{\Q_p})$-modules over $\Q_p^\un$; indeed, note that the proof of loc.\ cit.\ works just as well if we replace $\Perf(-)$ by $\Perf^{[a, b]}(-)$ for all $a\leq b$. Recalling that filtered $(\phi, N, \Gal_{\Q_p})$-modules over $\Q_p^\un$ are equivalent to potentially semistable vector bundles on $\FF_{\Spd\Q_p}$ by \cite[Prop.\ 10.6.7]{FarguesFontaine}, the first part of the theorem now follows by the theorem that ``de Rham $=$ potentially semistable'', see \cite[Thm.\ 10.6.10]{FarguesFontaine}.

For the second part about Ext groups, we can reduce to the case $V=\Q_p$ by replacing $W$ by $W\tensor V^\vee$ and then \cref{thm:hkcomp-main} tells us that we have to check whether the pullback
\begin{equation*}
\begin{tikzcd}
R\Gamma(\Q_p^\Syn, W)\ar[r]\ar[d] & R\Gamma(\Q_p^\HK, T_\HK(W))\ar[d] \\
R\Gamma(\Q_p^{\dR, +}, T_{\dR, +}(W))\ar[r] & R\Gamma(\Q_p^\dR, T_\dR(W))
\end{tikzcd}
\end{equation*}
computes the Exts between $\Q_p$ and $W$ in de Rham representations. Recalling that cohomology on $\Q_p^\HK$ computes cohomology of $(\phi, N, \Gal_{\Q_p})$-modules by \cite[Thm.\ 7.1.1]{dRFF} and using the theorem ``de Rham $=$ potentially semistable'', the claim follows by Galois descent from a computation of Emerton--Kisin of the Ext groups of semistable representations in terms of filtered $(\phi, N)$-modules, see \cite[Cor.\ 2.4.4]{CrystallineExt}.
\end{proof}

\newpage

\section{The comparison with proétale cohomology}
\label{sect:proet}

We now play the game from the previous section in the other direction: Instead of approaching the syntomification of Berkovich smooth derived Berkovich spaces from the de Rham side, we will now exploit the relation between $X^\Syn$ and ``étale data'' of $X$. In fact, this will work using almost the same arguments as in the previous section, but ``in the other direction'', i.e.\ we have to swap the roles of $u$ and $t$ and reverse the direction on the $\kappa$-axis.

Namely, we are going to prove that, for a Berkovich smooth derived Berkovich space $X$, perfect complexes on $X^\Syn$ and their cohomology can be understood via perfect complexes on $X^{\Div^1}$ and on the stack $X^{\HT, \dagger, +}$ that we have introduced previously. This will again take the form of a commutative diagram
\begin{equation}
\label{eq:proet-square}
\begin{tikzcd}
X^{\HT, \dagger}\ar[r]\ar[d] & X^{\HT, \dagger, +}\ar[d, "i_{\HT, \dagger, +}"] \\
X^\mDiv\ar[r, "i_{\Div^1}"] & X^\Syn\nospacepunct{\;,}
\end{tikzcd}
\end{equation}
where $X^\mDiv$ should again be thought of as $X^{\Div^1}$ for all intents and purposes, for which we are going to show that it induces an equivalence of categories
\begin{equation}
\label{eq:proet-perffibreproduct}
\Perf(X^\Syn)\cong \Perf(X^{\Div^1})\times_{\Perf(X^{\HT, \dagger})} \Perf(X^{\HT, \dagger, +})\;.
\end{equation}
Let us point out that, given the fact that the corresponding result (\ref{eq:hkcomp-perffibreproduct}) in the ``de Rham case'' has an analogue for $p$-adic formal schemes in terms of the algebraic syntomification of Drinfeld and Bhatt--Lurie, it seems to be an interesting question to investigate to which extent the naive analogue of (\ref{eq:proet-perffibreproduct}) holds for the algebraic syntomification of smooth $p$-adic formal schemes.

As the equivalence (\ref{eq:proet-perffibreproduct}) again entails a fibre product formula expressing syntomic cohomology of $X$ with coefficients in a perfect analytic $F$-gauge in terms of cohomology on $X^{\Div^1}, X^{\HT, \dagger, +}$ and $X^{\HT, \dagger}$, we will be able to recover the comparison between syntomic cohomology and proétale cohomology of a smooth partially proper rigid space $X$ in a suitable range from \cite[Thm.\ 1.1]{padicComparisons}. Indeed, as cohomology on $X^{\Div^1}$ agrees with proétale cohomology by \cite{AnPrism}, this comes down to controlling the difference between cohomology on $X^{\HT, \dagger, +}$ and $X^{\HT, \dagger}$, which is precisely responsible for the necessary truncation in the comparison theorem.

Finally, we will make use of the equivalence (\ref{eq:proet-perffibreproduct}) above to explicitly describe vector bundles on the syntomification of any smooth partially proper rigid space $X$ over $\Q_p$. Namely, we are going to prove that there is an equivalence 
\begin{equation*}
\Vect(X^\Syn)\cong \Vect^\dR(\FF_{X^\diamond})
\end{equation*}
between vector bundles on $X^\Syn$ and vector bundles on $\FF_{X^\diamond}$ which are de Rham in a certain sense to be defined below, see \cref{defi:drbundles-drvb}. For this, we are going to show that the pullback functor $\Vect(X^{\HT, \dagger, +})\rightarrow\Vect(X^{\HT, \dagger})$ is fully faithful. To determine its essential image and make the de Rham condition appear, we are going to use the main result of §\ref{sect:padicht}: Namely, the functor
\begin{equation*}
\begin{split}
\Vect(X^{\dR, +})\cong \Vect(X^{\HT, \dagger, +})&\rightarrow\Vect(X^{\HT, \dagger}) \\
&\rightarrow \Vect((\widehat{X}\subseteq Y_{X^\diamond})^\wedge)\cong\{\mathbb{B}_\dR^+\text{-local systems on }X_\proet\}
\end{split}
\end{equation*}
obtained from \cref{prop:htdr-perfequiv} agrees with Scholze's functor
\begin{equation*}
\{\text{filtered vector bundles on $X$ with connection}\}\rightarrow \{\mathbb{B}_\dR^+\text{-local systems on }X_\proet\}
\end{equation*}
from \cite[§7]{PAdicHodgeTheory} and hence the functor $\Vect(X^{\HT, \dagger, +})\rightarrow\Vect(X^{\HT, \dagger})$ automatically factors through a certain full subcategory of ``generically flat'' vector bundles on $X^{\HT, \dagger}$ to be defined below, see \cref{defi:drbundles-genflatxhtdagger}. The hardest part will then be to show that this factorisation is already essentially surjective and this will involve an explicit computation using differential operators, which we carry out in \cref{lem:drbundles-keylem}.

\subsection{Accessing $X^\Syn$ via $X^\prism$ and $X^{\HT, \dagger, +}$}

We start by defining the stack $X^\mDiv$ and constructing the commutative diagram (\ref{eq:proet-square}). To this end, we introduce the following object: For any Gelfand stack $X$, the Gelfand stack $X^{\mDiv, \mathrm{pre}}$ is defined as the pushout
\begin{equation*}
\begin{tikzcd}
X^\prism_{[p^{-1/2}, p^{-1/2}]}\ar[r, "\can"]\ar[d, "\phi\times \{1\}", swap] & X^\prism_{[p^{-1/2}, p^{1/2}]}\ar[d] \\
X^\prism_{(0, p^{1/2}]}\times [0, 1]\ar[r] & X^{\mDiv, \mathrm{pre}}\nospacepunct{\;.}
\end{tikzcd}
\end{equation*}
From this, we again obtain $X^\mDiv$ via gluing. Namely, note that the map
\begin{equation*}
X^\prism_{(0, p^{-1/2}]}\xrightarrow{\phi\times\{1\}} X^\prism_{(0, p^{1/2}]}\times [0, 1]
\end{equation*}
and the identity on $X^\prism_{[p^{-1/2}, p^{1/2}]}$ glue to a map
\begin{equation*}
\phi\times\{1\}: X^\prism_{[p^{-1/2}, p^{1/2}]}\coprod_{X^\prism_{[p^{-1/2}, p^{-1/2}]}} X^\prism_{(0, p^{-1/2}]}\rightarrow X^{\mDiv, \mathrm{pre}}\;,
\end{equation*}
the source of which is isomorphic to $X^\prism_{(0, p^{1/2}]}$ by \cref{lem:hkcomp-glueoverinterval}.

\begin{defi}
Let $X$ be any Gelfand stack. The Gelfand stack $X^\mDiv$ is defined by the coequaliser diagram
\begin{equation*}
\begin{tikzcd}
X^\prism_{(0, p^{1/2}]}\ar[r,shift left=.75ex,"\id\times\{0\}"]\ar[r,shift right=.75ex,swap,"\phi\times\{1\}"] & X^{\mDiv, \mathrm{pre}}\ar[r] & X^\mDiv \nospacepunct{\;.}
\end{tikzcd}
\end{equation*}
\end{defi}

To obtain the map $i_{\Div^1}: X^\mDiv\rightarrow X^\Syn$, recall that
\begin{equation*}
X^\N_{(0, p^{1/2}]}\cong X^\prism_{(0, p^{1/2}]}\times [0, 1]
\end{equation*}
via $\pi$ by \cref{prop:defis-utneq0}. Moreover, recalling that $X^\N_{|u|=1}\cong X^\prism$ from \cref{lem:defis-jht}, we also obtain an isomorphism
\begin{equation*}
X^\N_{|u|=1, [p^{-1/2}, p^{1/2}]}\cong X^\prism_{[p^{-1/2}, p^{1/2}]}\;.
\end{equation*}
As the isomorphism $X^\N_{|u|=1}\cong X^\prism$ intertwines $\pi$ with the Frobenius on $X^\prism$, the two displayed isomorphisms above glue to give a map
\begin{equation*}
X^{\mDiv, \mathrm{pre}}\rightarrow X^\N
\end{equation*}
which is an isomorphism onto its image, see \cref{fig:xmdivpre}. Finally, this descends to the desired map
\begin{equation*}
i_{\Div^1}: X^\mDiv\rightarrow X^\Syn\;.
\end{equation*}

In order to get to the desired diagram (\ref{eq:proet-square}) from here, note that $X^{\HT, \dagger}\rightarrow X^\prism$ factors through $X^\prism_{[1, 1]}$ and hence we obtain a composite map
\begin{equation}
\label{eq:proet-mapxhtdaggerxmdiv}
X^{\HT, \dagger}\rightarrow X^\prism_{[p^{-1/2}, p^{1/2}]}\rightarrow X^{\mDiv, \mathrm{pre}}\rightarrow X^\mDiv\;.
\end{equation}
Observing that the postcomposition of this map with $i_{\Div^1}$ identifies with the precomposition of the natural map
\begin{equation*}
i_{\HT, \dagger, +}: X^{\HT, \dagger, +}\rightarrow X^\N\rightarrow X^\Syn
\end{equation*}
with the canonical map $X^{\HT, \dagger}\rightarrow X^{\HT, \dagger, +}$, we obtain the desired commutative diagram
\begin{equation*}
\begin{tikzcd}
X^{\HT, \dagger}\ar[r]\ar[d] & X^{\HT, \dagger, +}\ar[d, "i_{\HT, \dagger, +}"] \\
X^\mDiv\ar[r, "i_{\Div^1}"] & X^\Syn\nospacepunct{\;.}
\end{tikzcd}
\end{equation*}

\begin{figure}

\begin{center}
\begin{tikzpicture}
  \def\width{5.5} 
  \def\height{1.2*\width} 

  \draw[thick] (0, 0) -- (0, \height); 
  \draw[thick] (\width, 0) -- (\width, \height); 
  \draw[dotted, thick] (0, 0) -- (\width, 0); 

  

  
  \draw[thick] (\width-0.2, 0.2*\height) -- (\width+0.2, 0.2*\height); 
  \node at (\width+1.2, 0.2*\height) {$\phi^{-2}(X^\HT)$};
  \draw[thick] (\width-0.2, 0.4*\height) -- (\width+0.2, 0.4*\height); 
  \node at (\width+1.2, 0.4*\height) {$\phi^{-1}(X^\HT)$};
  \draw[thick] (\width-0.2, 0.6*\height) -- (\width+0.2, 0.6*\height); 
  \node at (\width+0.7, 0.6*\height) {$X^\HT$};
  \draw[thick] (\width-0.2, 0.8*\height) -- (\width+0.2, 0.8*\height); 
  \node at (\width+0.7, 0.8*\height) {$X^\dR$};

  \draw[thick] (-0.2, 0.2*\height) -- (0.2, 0.2*\height); 
  \node at (-1.2, 0.2*\height) {$\phi^{-1}(X^\HT)$};
  \draw[thick] (-0.2, 0.4*\height) -- (0.2, 0.4*\height); 
  \node at (-0.7, 0.4*\height) {$X^\HT$};
  \draw[thick] (-0.2, 0.6*\height) -- (0.2, 0.6*\height); 
  \node at (-0.7, 0.6*\height) {$X^\dR$};
  \draw[thick] (-0.2, 0.8*\height) -- (0.2, 0.8*\height); 
  \node at (-0.95, 0.8*\height) {$\phi(X^\dR)$};


  \draw[line width=0.2mm, color=blue] (\width/2-0.05*\width, 0.6*\height) -- (0, 0.6*\height);
  \draw[line width=0.2mm, color=blue] (\width, 0.8*\height) -- (0, 0.8*\height);

  \def\sqr_size{0.1*\width} 

  \draw[line width=0.2mm] (\width/2 - \sqr_size/2, 0.6*\height - \sqr_size/2) -- (\width/2 + \sqr_size/2, 0.6*\height + \sqr_size/2); 
  \draw[line width=0.2mm] (\width/2 + \sqr_size/2, 0.6*\height - \sqr_size/2) -- (\width/2 - \sqr_size/2, 0.6*\height + \sqr_size/2); 

  
  \shade[shading=axis, bottom color=white, top color=white, middle color=green, shading angle=0](0.5*\width + \sqr_size/2, 0.6*\height - 0.1) rectangle (\width, 0.6*\height + 0.1);

  \draw[->, thick] (2*\width, 0) -- (2*\width, \height);
  
  \draw[thick] (2*\width-0.2, 0) -- (2*\width+0.2, 0); 
  \node at (2*\width+0.6, 0) {$0$};
  \draw[thick] (2*\width-0.2, 0.2*\height) -- (2*\width+0.2, 0.2*\height); 
  \node at (2*\width+0.6, 0.2*\height) {$1/p$};
  \draw[thick] (2*\width-0.2, 0.4*\height) -- (2*\width+0.2, 0.4*\height); 
  \node at (2*\width+0.6, 0.4*\height) {$1$};
  \draw[thick] (2*\width-0.2, 0.6*\height) -- (2*\width+0.2, 0.6*\height); 
  \node at (2*\width+0.6, 0.6*\height) {$p$};
  \draw[thick] (2*\width-0.2, 0.8*\height) -- (2*\width+0.2, 0.8*\height); 
  \node at (2*\width+0.6, 0.8*\height) {$p^2$};
  
  \draw[->, thick] (1.3*\width, \height/2) -- (1.8*\width, \height/2);
  \node at (1.55*\width, \height/2+0.3) {$\kappa$};
  
  \fill[color=gray!50] (0.05, 0) rectangle (\width-0.05, 0.5*\height-0.05);
  
  \draw[->, thick] (-0.5*\width, 0.3*\height) -- (-0.5*\width, 0.7*\height);
  \node at (-0.5*\width-0.4, 0.5*\height) {$\phi$};
  
  \shade[shading=axis, bottom color=gray!50, top color=gray!50, middle color=green, shading angle=0](0, 0.4*\height - 0.1) rectangle (\width, 0.4*\height + 0.1);
  \shade[shading=axis, bottom color=gray!50, top color=gray!50, middle color=green, shading angle=0](0, 0.2*\height - 0.1) rectangle (\width, 0.2*\height + 0.1);
  
  \draw[line width=1mm] (\width, 0) -- (\width, 0.7*\height);
  \draw[line width=1mm] (0, 0) -- (0, 0.5*\height+0.05);
  \draw[line width=1mm] (0, 0.5*\height) -- (\width, 0.5*\height);
  
  \draw[line width=0.2mm] (\width/2 - \sqr_size/2, 0.6*\height - \sqr_size/2) rectangle (\width/2 + \sqr_size/2, 0.6*\height + \sqr_size/2); 

\end{tikzpicture}
\end{center}

\captionsetup{justification=centering}
\caption{A schematic picture of $X^\N$ with the image of $X^{\mDiv, \mathrm{pre}}$ \\ outlined in bold (interior shaded in grey)}
\label{fig:xmdivpre}
\end{figure}

As before with the stack $X^{\mHK}$, for the introduction of $X^\mDiv$ to have any point at all, we need to relate its category of perfect complexes to perfect complexes on $X^{\Div^1}$.

\begin{prop}
\label{prop:proet-psprism}
Let $X$ be any Gelfand stack over $\Q_p$. There is an equivalence of categories
\begin{equation*}
\Perf(X^{\mDiv})\cong\Perf(X^{\Div^1})
\end{equation*}
induced by the diagram
\begin{equation*}
\begin{tikzcd}
& X^\prism_{(0, p^{1/2}]}\ar[rd]\ar[ld, "\phi\times \{1\}", swap] & \\
X^{\mDiv} & & X^{\Div^1}\nospacepunct{\;.}
\end{tikzcd}
\end{equation*}
\end{prop}
\begin{proof}
As in \cref{prop:hkcomp-hkpshk}.
\end{proof}

Using \cref{prop:proet-psprism}, the diagram (\ref{eq:proet-square}) induces various realisation functors for perfect analytic $F$-gauges. Firstly, pullback along the map $i_{\Div^1}$ induces a realisation functor
\begin{equation*}
T_{\Div^1}: \Perf(X^\Syn)\rightarrow\Perf(X^{\Div^1})
\end{equation*}
from perfect analytic $F$-gauges to perfect complexes on $X^{\Div^1}$, which should be thought of as coefficients for proétale cohomology of $X$ following \cite{AnPrism}. Secondly, via pullback along the map $i_{\HT, \dagger, +}$ from above, we obtain a realisation functor
\begin{equation*}
T_{\HT, \dagger, +}: \D(X^\Syn)\rightarrow \D(X^{\HT, \dagger, +})
\end{equation*}
from analytic $F$-gauges to quasicoherent sheaves on $X^{\HT, \dagger, +}$ and, lastly, postcomposing $T_{\HT, \dagger, +}$ with pullback along the canonical map $X^{\HT, \dagger}\rightarrow X^{\HT, \dagger, +}$ yields a realisation functor
\begin{equation*}
T_{\HT, \dagger}: \D(X^\Syn)\rightarrow \D(X^{\HT, \dagger})\;.
\end{equation*}
Note that pullback along the map $X^{\HT, \dagger}\rightarrow X^\mDiv$ from (\ref{eq:proet-mapxhtdaggerxmdiv}) induces a functor $\Perf(X^{\Div^1})\rightarrow\Perf(X^{\HT, \dagger})$ whose precomposition with $T_\HK$ identifies with $T_{\HT, \dagger}$ by commutativity of the diagram (\ref{eq:proet-square}). This enables us to formulate the ``étale'' analogue of \cref{thm:hkcomp-main}:

\begin{thm}
\label{thm:proet-main}
Let $X$ be a Berkovich smooth derived Berkovich space over $\Q_p$. Then there is an equivalence of categories
\begin{equation*}
\Perf(X^\Syn)\cong \Perf(X^{\Div^1})\times_{\Perf(X^{\HT, \dagger})} \Perf(X^{\HT, \dagger, +})
\end{equation*}
induced by the realisation functors $T_{\Div^1}$ and $T_{\HT, \dagger, +}$. In particular, for any $E\in\Perf(X^\Syn)$, there is a pullback diagram
\begin{equation*}
\begin{tikzcd}
R\Gamma(X^\Syn, E)\ar[r]\ar[d] & R\Gamma(X^{\Div^1}, T_{\Div^1}(E))\ar[d] \\
R\Gamma(X^{\HT, \dagger, +}, T_{\HT, \dagger, +}(E))\ar[r] & R\Gamma(X^{\HT, \dagger}, T_{\HT, \dagger}(E))\nospacepunct{\;.}
\end{tikzcd}
\end{equation*}
\end{thm}

Using that cohomology on $X^{\Div^1}$ coincides with proétale cohomology for smooth partially proper rigid spaces $X$ over $\Q_p$ by \cite{AnPrism}, the result above lets us recover the classical comparison theorem between syntomic cohomology and proétale cohomology in a certain range, see e.g.\ \cite[Thm.\ 1.1.(4)]{padicComparisons}.

\begin{thm}
\label{thm:proet-main2}
Let $X$ be a smooth partially proper rigid space over $\Q_p$. If $E\in\Vect(X^\Syn)$ is a vector bundle analytic $F$-gauge with Hodge--Tate weights all at most $-i$ for some $i\geq 0$, then the natural morphism
\begin{equation*}
R\Gamma(X^\Syn, E)\rightarrow R\Gamma(X^{\Div^1}, T_{\Div^1}(E))
\end{equation*}
is an isomorphism on $\tau^{\leq i}$ and induces an injection on $H^{i+1}$. In particular, for $E=\O\{i\}$, we obtain
\begin{equation*}
\tau^{\leq i} R\Gamma_\Syn(X, \Q_p(i))\cong \tau^{\leq i}R\Gamma_\proet(X, \Q_p(i))\;.
\end{equation*}
\end{thm}
\begin{proof}
By \cref{thm:proet-main}, our task is to prove that the cofibre of the morphism
\begin{equation*}
R\Gamma(X^{\HT, \dagger, +}, T_{\HT, \dagger, +}(E))\rightarrow R\Gamma(X^{\HT, \dagger}, T_{\HT, \dagger}(E))
\end{equation*}
is concentrated in degrees at least $i+1$. By compatibility of the functors $X\mapsto X^{\HT, \dagger, +}$ and $X\mapsto X^{\HT, \dagger}$ with open localisation, which one deduces from \cref{prop:defis-openloc}, we are reduced to the case where $X$ admits a Berkovich étale map $X\rightarrow\ol{\T}^n$. However, in this case, \cref{cor:htdr-cohomologyxhtdrdagger+} gives precise double complexes computing the cohomologies in question. Indeed, inspecting these double complexes, we see that our claim follows once we show that the ascending filtration $\Fil_\bullet V$ associated to the restriction of $E$ to $X^{\HT, \dagger, +}$ via \cref{cor:htdr-complexesxhtdrdagger+} stabilises from $\Fil_{-i} V$ onwards, i.e.\
\begin{equation*}
\Fil_{-i} V=\Fil_{-i+1} V=\dots\;.
\end{equation*}

To see this, recall from \cref{cor:htdr-complexesxhtdrdagger+} that the (descendingly) filtered vector bundle with connection $(\Fil^\bullet W, \nabla)$ on $X$ corresponding to the restriction of $E$ to $X^{\dR, +}$ is obtained by Tate twisting and then descending the associated graded pieces of $\Fil_\bullet V$ to $X$, i.e.\ $\gr_k V$ descends to $\Fil^k W$ up to a Tate twist, and hence we equivalently have to show that
\begin{equation*}
\dots=\Fil^{-i+2} W=\Fil^{-i+1} W=0\;.
\end{equation*}
However, our assumption on the Hodge--Tate weights precisely means that $\gr^k W=0$ for $k\geq -i+1$ and now we are done by noting that the decreasing filtration $\Fil^\bullet W$ is separated, which follows e.g.\ from \cref{cor:htdr-complexesxhtdrdagger+}.
\end{proof}

\begin{rem}
Note that \cref{thm:proet-main} in some sense improves upon the ``classical'' comparison from \cref{thm:proet-main2} by also describing the failure of the map $R\Gamma(X^\Syn, \O\{i\})\rightarrow R\Gamma_\proet(X, \Q_p(i))$ to be an isomorphism in higher degrees: Namely, this is exactly measured by the failure of $R\Gamma(X^{\HT, \dagger, +}, \O\{i\})\rightarrow R\Gamma(X^{\HT, \dagger}, \O\{i\})$ to be an isomorphism.
\end{rem}

Of course, combining \cref{thm:proet-main2} with \cref{thm:hkcomp-main}, we immediately obtain the following extension of Colmez--Nizio{\l}'s ``basic comparison theorem'' from \cite[Thm.\ 1.3]{padicComparisons} to vector bundle $F$-gauge coefficients:

\begin{cor}
Let $X$ be a smooth partially proper qcqs rigid space over $\Q_p$. If $E\in\Vect(X^\Syn)$ is a vector bundle analytic $F$-gauge with Hodge--Tate weights all at most $-i$ for some $i\geq 0$, then there is a natural map
\begin{equation*}
\fib(R\Gamma(X^\HK, T_\HK(E))\rightarrow R\Gamma(X^\dR, T_\dR(E))/R\Gamma(X^{\dR, +}, T_{\dR, +}(E)))\rightarrow R\Gamma(X^{\Div^1}, T_{\Div^1}(E))
\end{equation*}
which induces an isomorphism on $\tau^{\leq i}$ and an injection on $H^{i+1}$. In particular, for $E=\O\{i\}$, we obtain
\begin{equation*}
\tau^{\leq i} R\Gamma_\proet(X, \Q_p(i))\cong \tau^{\leq i}\fib(R\Gamma_\HK(X)^{\phi=p^i, N=0, \Gal_{\Q_p}}\rightarrow R\Gamma_\dR(X)/\Fil^i_\Hod R\Gamma_\dR(X))\;.
\end{equation*} 
\end{cor}

Let us now say something about the proof of \cref{thm:proet-main}. As everything is entirely analogous to what happens in the proof of \cref{thm:hkcomp-main}, we just explain the steps of the proof in order to avoid too much redundancy. We first note the following analogue of \cref{lem:hkcomp-chopupxn}:

\begin{lem}
\label{lem:proet-chopupxn}
Let $X$ be any Gelfand stack. The canonical inclusion maps induce an isomorphism between the iterated pushout
\begin{equation*}
X^\N_{(0, p^{3/2}]}\coprod_{X^\N_{[p^{3/2}, p^{3/2}]}} X^\N_{[p^{3/2}, p^{5/2}]}\coprod_{X^\N_{[p^{5/2}, p^{5/2}]}} X^\N_{[p^{5/2}, p^{7/2}]}\coprod_{X^\N_{[p^{7/2}, p^{7/2}]}} \dots
\end{equation*}
and the stack $X^\N$.
\end{lem}
\begin{proof}
Repeatedly apply \cref{lem:hkcomp-glueoverinterval}.
\end{proof}

\begin{figure}

\begin{center}
\begin{tikzpicture}
  \def\width{5.5} 
  \def\height{1.2*\width} 
  \def\sqr_size{0.1*\width} 

  

  
  \draw[thick] (\width-0.2, 0.1*\height) -- (\width+0.2, 0.1*\height); 
  \node at (\width+0.7, 0.1*\height) {$X^\HT$};
  \draw[thick] (\width-0.2, -0.1*\height) -- (\width+0.2, -0.1*\height); 
  \node at (\width+1.2, -0.1*\height) {$\phi^{-1}(X^\HT)$};
  \draw[thick] (\width-0.2, -0.3*\height) -- (\width+0.2, -0.3*\height); 
  \node at (\width+1.2, -0.3*\height) {$\phi^{-2}(X^\HT)$};

  \draw[thick] (-0.2, -0.1*\height) -- (0.2, -0.1*\height); 
  \node at (-0.7, -0.1*\height) {$X^\HT$};
  \draw[thick] (-0.2, -0.3*\height) -- (0.2, -0.3*\height); 
  \node at (-1.2, -0.3*\height) {$\phi^{-1}(X^\HT)$};


  \draw[line width=0.2mm, color=blue] (\width/2-0.05*\width, 0.1*\height) -- (-\width, 0.1*\height);

  
  \shade[shading=axis, bottom color=white, top color=white, middle color=green, shading angle=0](0.5*\width + \sqr_size/2, 0.1*\height - 0.1) rectangle (\width, 0.1*\height + 0.1);
  \shade[shading=axis, bottom color=white, top color=white, middle color=green, shading angle=0](0, -0.1*\height - 0.1) rectangle (\width, -0.1*\height + 0.1);
  \shade[shading=axis, bottom color=white, top color=white, middle color=green, shading angle=0](0, -0.3*\height - 0.1) rectangle (\width, -0.3*\height + 0.1);
  
  \draw[line width=0.2mm] (\width/2 - \sqr_size/2, 0.1*\height - \sqr_size/2) rectangle (\width/2 + \sqr_size/2, 0.1*\height + \sqr_size/2); 

  \draw[line width=0.2mm] (\width/2 - \sqr_size/2, 0.1*\height - \sqr_size/2) -- (\width/2 + \sqr_size/2, 0.1*\height + \sqr_size/2); 
  \draw[line width=0.2mm] (\width/2 + \sqr_size/2, 0.1*\height - \sqr_size/2) -- (\width/2 - \sqr_size/2, 0.1*\height + \sqr_size/2); 
  
  \draw[line width=1mm] (0, 0.05) -- (0, -0.4*\height); 
  \draw[line width=1mm] (0, 0) -- (-\width, 0); 
  \draw[line width=1mm] (\width, 0.2*\height+0.05) -- (\width, -0.4*\height); 
  \draw[line width=1mm] (-\width, 0.2*\height) -- (\width, 0.2*\height); 

\end{tikzpicture}
\end{center}

\captionsetup{justification=centering}
\caption{A schematic picture of $X^{\Syn, \mathrm{pre}, \Div^1}$ with \\ the two copies of $\partial X^{\Syn, \mathrm{pre}, \Div^1}$ in bold}
\label{fig:xsynprediv}
\end{figure}

For any Gelfand stack $X$, we now define the stack $X^{\Syn, \mathrm{pre}, \Div^1}$ as the pushout
\begin{equation*}
\begin{tikzcd}
X^\prism_{[p^{1/2}, p^{3/2}]}\ar[d, "j_\dR", swap]\ar[r, "\id\times\{0\}"] & X^\prism_{[p^{1/2}, p^{3/2}]}\times (-\infty, 0]\ar[d] \\ 
X^\N_{(0, p^{3/2}]}\ar[r] & X^{\Syn, \mathrm{pre}, \Div^1}\nospacepunct{\;,}
\end{tikzcd}
\end{equation*}
see \cref{fig:xsynprediv}. Then $X^{\Syn, \mathrm{pre}, \Div^1}$ receives two maps from the pushout
\begin{equation*}
\begin{tikzcd}
X^\prism_{[p^{1/2}, p^{1/2}]}\ar[r, "\id\times \{0\}"]\ar[d, "\can"] & X^\prism_{[p^{1/2}, p^{1/2}]}\times (-\infty, 0]\ar[d] \\
X^\prism_{(0, p^{1/2}]}\ar[r] & \partial X^{\Syn, \mathrm{pre}, \Div^1}\nospacepunct{\;.}
\end{tikzcd}
\end{equation*}
Namely, one of them is induced by $j_\dR: X^\prism_{(0, p^{1/2}]}\rightarrow X^\N_{(0, p^{3/2}]}$ and the canonical inclusion $X^\prism_{[p^{1/2}, p^{1/2}]}\rightarrow X^\prism_{[p^{1/2}, p^{3/2}]}$ and we denote this map by
\begin{equation*}
j_\dR: \partial X^{\Syn, \mathrm{pre}, \Div^1}\rightarrow X^{\Syn, \mathrm{pre}, \Div^1}
\end{equation*}
as well. The other one is induced by $j_\HT: X^\prism_{(0, p^{1/2}]}\rightarrow X^\N_{(0, p^{3/2}]}$ and the isomorphism
\begin{equation}
\label{eq:proet-jhtxsynpre}
X^\prism_{[p^{1/2}, p^{1/2}]}\times (-\infty, 0]\xrightarrow{\cong} (X^\prism_{[p^{3/2}, p^{3/2}]}\times [0, 1])\coprod_{X^\prism_{[p^{1/2}, p^{1/2}]}} (X^\prism_{[p^{1/2}, p^{1/2}]}\times (-\infty, 0])\;,
\end{equation}
where the maps defining the pushout on the right-hand side are $\phi\times \{0\}$ ``to the left'' and $\id\times \{0\}$ ``to the right'', respectively, and we use the isomorphism
\begin{equation*}
X^\N_{[p^{3/2}, p^{3/2}]}\cong X^\prism_{[p^{3/2}, p^{3/2}]}\times [0, 1]
\end{equation*}
from \cref{prop:defis-utneq0} to view the target of (\ref{eq:proet-jhtxsynpre}) as a closed substack of $X^{\Syn, \mathrm{pre}, \Div^1}$. We denote this second map by
\begin{equation*}
j_\HT: \partial X^{\Syn, \mathrm{pre}, \Div^1}\rightarrow X^{\Syn, \mathrm{pre}, \Div^1}\;.
\end{equation*}
Now the next lemma reads as follows:

\begin{lem}
\label{lem:proet-syndifferentgluing}
Let $X$ be any Gelfand stack. Then $X^\Syn$ is isomorphic to the coequaliser of the two maps
\begin{equation*}
\begin{tikzcd}
\partial X^{\Syn, \mathrm{pre}, \Div^1}\ar[r,shift left=.75ex, "j_\HT"]
  \ar[r,shift right=.75ex,swap, "j_\dR"] & X^{\Syn, \mathrm{pre}, \Div^1}\nospacepunct{\;.}
\end{tikzcd}
\end{equation*}
\end{lem}
\begin{proof}
Analogous to \cref{lem:hkcomp-syndifferentgluing} replacing the use of \cref{lem:hkcomp-chopupxn} by \cref{lem:proet-chopupxn}.
\end{proof}

\begin{lem}
\label{lem:proet-xsynviarestrxn}
For any Gelfand stack $X$, pullback along the morphism
\begin{equation*}
\operatorname{coeq}(\hspace{-0.15cm}
\begin{tikzcd}
X^\prism_{(0, p^{1/2}]}\ar[r,shift left=.75ex, "j_\HT"]
  \ar[r,shift right=.75ex,swap, "j_\dR"] & X^\N_{(0, p^{3/2}]}
\end{tikzcd}
\hspace{-0.15cm})\rightarrow X^\Syn
\end{equation*}
induces an equivalence on perfect complexes.
\end{lem}
\begin{proof}
Analogous to \cref{lem:hkcomp-xsynviarestrxn} replacing the use of \cref{lem:hkcomp-syndifferentgluing} by \cref{lem:proet-syndifferentgluing}.
\end{proof}

The following lemma is again the key to proving \cref{thm:proet-main}. The proof is analogous to the one of \cref{lem:hkcomp-keylem}, but uses the part 
\begin{equation*}
\Perf((X^\N)_{|ut|=0})\cong \Perf(X^{\HT, \dagger, +})
\end{equation*}
of the statement of \cref{prop:htdr-perfequiv} instead of the equivalence between perfect complexes on $(X^\N)_{|ut|=0}$ and $X^{\dR, +}$.

\begin{lem}
\label{lem:proet-keylem}
Let $X$ be a Berkovich smooth derived Berkovich space over $\Q_p$. Then pullback along the map
\begin{equation*}
X^\prism_{[p^{-1/2}, p^{1/2}]}\coprod_{X^{\HT, \dagger}} X^{\HT, \dagger, +}\rightarrow X^\N_{[p^{1/2}, p^{3/2}]}\;,
\end{equation*}
which is induced by $j_\HT$ and $i_{\HT, \dagger, +}$, induces an equivalence on perfect complexes.
\end{lem}
\begin{proof}
Analogous to \cref{lem:hkcomp-keylem} with the modification mentioned above.
\end{proof}

Finally, one puts together \cref{lem:proet-xsynviarestrxn} and \cref{lem:proet-keylem} to obtain \cref{thm:proet-main}, analogously to what happens in the proof of \cref{thm:hkcomp-main}.

\subsection{De Rham bundles on $\FF_{X^\diamond}$ via the syntomification}
\label{subsect:drbundles}

We now want to use \cref{thm:proet-main} to explicitly describe vector bundles on $X^\Syn$ for a smooth partially proper rigid space $X$ over $\Q_p$. For this, recall the functor
\begin{equation*}
\begin{split}
\Vect(X^{\dR, +})\cong \Vect(X^{\HT, \dagger, +})&\rightarrow\Vect(X^{\HT, \dagger}) \\
&\rightarrow \Vect((\widehat{X}\subseteq Y_{X^\diamond})^\wedge)\cong\{\mathbb{B}_\dR^+\text{-local systems on }X_\proet\}
\end{split}
\end{equation*}
obtained from \cref{prop:htdr-perfequiv}, which we have shown to agree with Scholze's functor
\begin{equation}
\label{eq:proet-scholzefunctor}
\{\text{filtered vector bundles with connection on $X$}\}\rightarrow \{\mathbb{B}_\dR^+\text{-local systems on $X_\proet$}\}
\end{equation}
from \cite[§7]{PAdicHodgeTheory} in \cref{thm:padicht-scholzefunctor}. Also recall from \cite[Thm.\ 7.7]{PAdicHodgeTheory} that (\ref{eq:proet-scholzefunctor}) is fully faithful and that the image of a given filtered vector bundle with connection on $X$ under this functor is called its \emph{associated $\mathbb{B}_\dR^+$-local system}. We can thus make the following definition:

\begin{defi}
Let $X$ be a smooth partially proper rigid space over $\Q_p$. A $\mathbb{B}_\dR^+$-local system on $X_\proet$ is called \emph{generically flat} if it lies in the essential image of the functor (\ref{eq:proet-scholzefunctor}), i.e.\ if it is associated to a filtered vector bundle with connection on $X$.
\end{defi}

\begin{rem}
\label{rem:proet-genflatqp}
Note that the above definition recovers the notion of a generically flat $B_\dR^+$-representation of $\Gal_{\Q_p}$ of \cite[Def.\ 10.4.1]{FarguesFontaine} in the case $X=\GSpec\Q_p$. 
\end{rem}

\begin{rem}
Generic flatness of a $\mathbb{B}_\dR^+$-local system may be checked locally on $X$. Indeed, this is because (\ref{eq:proet-scholzefunctor}) is fully faithful and filtered vector bundles with connection satisfy descent for open covers.
\end{rem}

Now recall from (\ref{eq:padicht-bdr+locsysviaffx}) that there is an equivalence
\begin{equation*}
\Vect((\widehat{X}\subseteq \FF_{X^\diamond})^\wedge)\cong \{\mathbb{B}_\dR^+\text{-local systems on $X_\proet$}\}\;.
\end{equation*}
Thus, pullback along the map $(\widehat{X}\subseteq \FF_{X^\diamond})^\wedge\rightarrow\FF_{X^\diamond}$ yields a $\mathbb{B}_\dR^+$-local system associated to any vector bundle on $\FF_{X^\diamond}$.

\begin{defi}
\label{defi:drbundles-drvb}
Let $X$ be a smooth partially proper rigid space over $\Q_p$. A vector bundle on $\FF_{X^\diamond}$ is called \emph{de Rham} if its associated $\mathbb{B}_\dR^+$-local system is generically flat. We let $\Vect^\dR(\FF_{X^\diamond})$ be the full subcategory of vector bundles on $\FF_{X^\diamond}$ which are de Rham.
\end{defi}

\begin{rem}
Using \cref{rem:proet-genflatqp}, we see that the definition above recovers the notion of being de Rham for a $\Gal_{\Q_p}$-equivariant vector bundle on the Fargues--Fontaine curve from \cite[Def.\ 10.4.5]{FarguesFontaine} in the case $X=\GSpec\Q_p$.
\end{rem}

With the above notation in place, we can now give the desired explicit description of the category of vector bundles on $X^\Syn$.

\begin{thm}
\label{thm:drbundles-main}
Let $X$ be a smooth partially proper rigid space over $\Q_p$. Then there is an equivalence of categories
\begin{equation*}
\Vect(X^\Syn)\cong \Vect^\dR(\FF_{X^\diamond})\;.
\end{equation*}
\end{thm}

Before we move on to the proof of the theorem above, let us make some remarks about the relation between \cref{thm:drbundles-main} and the results from §\ref{sect:hkcomp}.

\begin{rem}
Contrary to the situation for $X=\GSpec\Q_p$ from \cref{thm:bk-mainqp}, cohomology on $X^\Syn$ does \emph{not} compute the $\RHom$ in the category of de Rham local systems on $X_\proet$ in general. This can already be seen in the case $X=\P^1$. Then \cref{cor:hkcomp-classical} yields a cartesian diagram
\begin{equation}
\label{eq:drbundles-cohp1}
\begin{tikzcd}
R\Gamma_\Syn(\P^1, \Q_p(1))\ar[r]\ar[d] & R\Gamma_\HK(\P^1)^{\phi=p, N=0, \Gal_{\Q_p}}\ar[d] \\
\Fil^1_\Hod R\Gamma_\dR(\P^1)\ar[r] & R\Gamma_\dR(\P^1)
\end{tikzcd}
\end{equation}
and we know that $R\Gamma_\dR(\P^1)\cong \Q_p\oplus \Q_p[-2]$ while $\Fil^1_\Hod R\Gamma_\dR(\P^1)\cong \Q_p[-2]$. Meanwhile, by \cite[Lem.\ 6.2.9]{dRFF}, we have $R\Gamma_\HK(\P^1)\cong \Q_p^\un\oplus \Q_p^\un[-2]$ with the monodromy being trivial, $\Gal_{\Q_p}$ acting in the canonical way and Frobenius fixing the basis element of $\Q_p^\un$ while multiplying the basis element of $\Q_p^\un[-2]$ by $p$. Thus, we have
\begin{equation*}
R\Gamma_\HK(\P^1)^{\phi=p, N=0, \Gal_{\Q_p}}\cong \Q_p\oplus \Q_p[-1]\oplus \Q_p[-2]^{\oplus 2}\oplus \Q_p[-3]
\end{equation*}
and therefore (\ref{eq:drbundles-cohp1}) yields
\begin{equation*}
R\Gamma_\Syn(\P^1, \Q_p(1))\cong \Q_p[-1]^{\oplus 2}\oplus \Q_p[-2]^{\oplus 2}\oplus \Q_p[-3]\;.
\end{equation*}
However, note that proétale de Rham local systems on $\P^1$ are just equivalent to de Rham $\Gal_{\Q_p}$-representations and hence the $\RHom$s in this category have vanishing $H^i$ for $i>2$ by the proof of \cref{thm:bk-mainqp}. In fact, using loc.\ cit., one easily calculates
\begin{equation*}
\RHom_{\Rep^\dR_{\Q_p}(\Gal_{\Q_p})}(\Q_p, \Q_p(1))\cong \Q_p[-1]^{\oplus 2}\oplus \Q_p[-2]\;. \qedhere
\end{equation*}
\end{rem}

\begin{rem}
Note that the case $X=\GSpec\Q_p$ of \cref{thm:drbundles-main} is already covered by \cref{thm:bk-mainqp}, but the proof we will give for \cref{thm:drbundles-main} is diametrically opposite to the one we have given for \cref{thm:bk-mainqp}: Indeed, here we will make use of the connection between $X^\Syn$ and $X^{\Div^1}$, i.e.\ the ``étale aspect'' of $X^\Syn$, while the proof of loc.\ cit.\ exploits the connection between $X^\Syn$ and $X^\HK$, i.e.\ the ``de Rham aspect'' of $X^\Syn$.
\end{rem}

\begin{rem}
Another view on the two proofs we have given for the above result in the case $X=\GSpec\Q_p$ is the following: While the proof via $\Q_p^{\Div^1}$ we will give below directly shows that vector bundles on $\Q_p^\Syn$ are de Rham bundles on $\FF_{\Spd\Q_p}$, the proof of \cref{thm:bk-mainqp} rather shows that vector bundles on $\Q_p^\Syn$ are the same as filtered $(\phi, N, \Gal_{\Q_p})$-modules, i.e.\ potentially semistable vector bundles on $\FF_{\Spd\Q_p}$. In other words, putting these two proofs together gives a new proof of the theorem that ``de Rham $=$ potentially semistable'' of Berger. Note, however, that we need the equivalence between vector bundles on $\Q_p^\HK$ and $(\phi, N, \Gal_{\Q_p})$-modules from \cite[Thm.\ 7.1.1]{dRFF} as an input for the proof of \cref{thm:bk-mainqp}; this should be seen as analogous to Berger's reduction of the $p$-adic monodromy theorem to Crew's conjecture about $p$-adic differential equations.
\end{rem}

Developing the previous remarks further, we might want to compare the results of \cref{thm:drbundles-main} and \cref{thm:hkcomp-main} more generally for smooth partially proper rigid spaces $X$ over $\Q_p$. This immediately yields the following result, which was already sketched in \cite[Rem.\ 7.5.7]{dRFF} and, as explained in loc.\ cit., in particular also implies Shimizu's theorem about de Rham local systems from \cite[Thm.\ 1.1]{Shimizu}:

\begin{cor}
Let $X$ be a smooth partially proper rigid space over $\Q_p$. Then there is an equivalence of categories
\begin{equation*}
\Vect^\dR(\FF_{X^\diamond})\cong \Vect(X^\HK)\times_{\Vect(X^\dR)} \Vect(X^{\dR, +})\;.
\end{equation*}
\end{cor}

Let us now come back to proving \cref{thm:drbundles-main}. For this, recall from \cref{prop:defis-xprismdr} that there is a natural map $Y_{X^\diamond}\rightarrow X^\prism$ for any Gelfand stack $X$ and, similarly, one proves that there is a natural map $\FF_{X^\diamond}\rightarrow X^{\Div^1}$. In our current situation, i.e.\ for $X$ a smooth partially proper rigid space over $\Q_p$, furthermore recall from \cite{AnPrism} that pullback along this map induces an equivalence
\begin{equation*}
\Vect(X^{\Div^1})\cong \Vect(\FF_{X^\diamond})\;.
\end{equation*}
Thus, in view of \cref{thm:proet-main}, the main part of the proof of \cref{thm:drbundles-main} is to understand the pullback functor $\Vect(X^{\HT, \dagger, +})\rightarrow \Vect(X^{\HT, \dagger})$. 

More precisely, we see that our task is to prove that this functor is fully faithful and that its essential image is exactly cut out by the de Rham condition in \cref{thm:proet-main}. More to the point, recall from (\ref{eq:padicht-xhtdaggertobdr+}) that there is a natural pullback functor
\begin{equation}
\label{eq:drbundles-xhtdaggertobdr+}
\Vect(X^{\HT, \dagger})\rightarrow \Vect((\widehat{X}\subseteq \FF_{X^\diamond})^\wedge)\cong \{\mathbb{B}_\dR^+\text{-local systems on $X_\proet$}\}\;,
\end{equation}
which associates a $\mathbb{B}_\dR^+$-local system on $X_\proet$ to any vector bundle on $X^{\HT, \dagger}$. This leads to the following definition:

\begin{defi}
\label{defi:drbundles-genflatxhtdagger}
Let $X$ be a smooth partially proper rigid space over $\Q_p$. A vector bundle on $X^{\HT, \dagger}$ is called \emph{generically flat} if its associated $\mathbb{B}_\dR^+$-local system on $X_\proet$ is generically flat.
\end{defi}

With this terminology in place, the above discussion shows that the proof of \cref{thm:drbundles-main} reduces to the following statement:

\begin{prop}
\label{prop:drbundles-xhtdagger+pullback}
Let $X$ be a smooth partially proper rigid space over $\Q_p$. Then the pullback functor 
\begin{equation*}
\Vect(X^{\HT, \dagger, +})\rightarrow \Vect(X^{\HT, \dagger})
\end{equation*}
is fully faithful and its essential image is given by the full subcategory of vector bundles on $X^{\HT, \dagger}$ which are generically flat.
\end{prop}

From \cref{thm:padicht-scholzefunctor}, we can already conclude that the pullback functor $\Vect(X^{\HT, \dagger, +})\rightarrow\Vect(X^{\HT, \dagger})$ factors through the full subcategory of vector bundles on $X^{\HT, \dagger}$ which are generically flat. As Scholze's functor (\ref{eq:proet-scholzefunctor}) is fully faithful and we already know that $\Vect(X^{\dR, +})\cong \Vect(X^{\HT, \dagger, +})$ from \cref{prop:htdr-perfequiv}, both the full faithfulness and essential surjectivity in \cref{prop:drbundles-xhtdagger+pullback} will follow if we can prove that the functor
\begin{equation*}
\{\text{generically flat vector bundles on $X^{\HT, \dagger}$}\}\rightarrow \{\mathbb{B}_\dR^+\text{-local systems on $X_\proet$}\}
\end{equation*}
obtained by restricting (\ref{eq:drbundles-xhtdaggertobdr+}) to generically flat vector bundles is fully faithful. By dualisability, this reduces to a cohomology computation, which we in turn would like to reduce to a question about vector bundles on $X^\HT$ and $\widehat{\O}_X$-local systems on $X_\proet$ using the $t$-adic filtration on the double complex from \cref{cor:htdr-cohomologyxhtdrdagger+}. However, note that this filtration is a priori not complete! This is remedied by the following key technical lemma:

\begin{lem}
\label{lem:drbundles-keylem}
Let $X$ be a smooth partially proper rigid space over $\Q_p$ equipped with an étale map $X\rightarrow \Spa\Q_p\langle x_1^{\pm 1}, \dots, x_n^{\pm 1}\rangle$ such that the induced map $X\rightarrow\ol{\T}^n$ is Berkovich étale. Let $E$ be a generically flat vector bundle on $X^{\HT, \dagger}$, which by \cref{cor:htdr-complexesxhtdrdagger+} corresponds to a vector bundle $V$ on $X_\infty^\la\times \GSpec\Q_p\{t\}^\dagger\times\GSpec\Q_p(\zeta_{p^\infty})$ equipped with a semilinear locally analytic $\Z_p^n\rtimes\Z_p^\times$-action. If $\Theta^\arithm: V\rightarrow V$ denotes the ``arithmetic Sen operator'' on $V$ obtained from the Lie algebra action of $\Z_p^\times$, then
\begin{equation*}
\Theta^\arithm: t^iV\rightarrow t^i V
\end{equation*}
is an isomorphism for $i\gg 0$.
\end{lem}
\begin{proof}
By generic flatness of $E$, there is a filtered vector bundle with connection $(\Fil^\bullet W, \nabla)$ on $X$ whose associated $\mathbb{B}_\dR^+$-local system agrees with the $\mathbb{B}_\dR^+$-local system associated to $E$. We will show that $\Theta^\arithm$ is an isomorphism on $t^i V$ for any $i$ such that $\Fil^i W=0$, which implies the claim.

To prove this, we may work locally on $X$ and, in particular, assume that all $\Fil^\bullet W$ are free $\O_X$-modules and that the transition maps of the filtration are given by the standard embeddings via the first couple of coordinates. Moreover, by \cite[Lem.\ 7.3]{PAdicHodgeTheory}, we may assume that the pullback of $V$ to $X_\infty^\la\times\GSpec\Q_p(\zeta_{p^\infty})$ is free after possibly further localising on $X$ and then we may even assume that $V$ itself is free since freeness of a finite projective module may be checked after $\dagger$-reduction (as the $\dagger$-nilradical is contained in the Jacobson radical). In the following, let us introduce the notation $\mathbb{B}_\dR^{+, \dagger, \la}$ for the pushforward of the structure sheaf along
\begin{equation*}
X_\infty^\la\times\GSpec\Q_p\{t\}^\dagger\times\GSpec\Q_p(\zeta_{p^\infty})\xrightarrow{\mathrm{pr}_1} X_\infty^\la\rightarrow X
\end{equation*}
and note that this may be viewed as a sheaf on the underlying topological space $|X|$ of $X$ in the sense of \cite[Def.\ 4.3.1]{dRFF}. We point out that the above map is affine in the sense that the source becomes affine after pulling back to any affine $\GSpec A\rightarrow X$. 

Now note that the semilinearity of the $\Z_p^\times$-action implies that $\Theta^\arithm$ will be a $t$-connection in the following sense: for all local sections $f$ of $\mathbb{B}_\dR^{+, \dagger, \la}$ and $v$ of $V$, we have
\begin{equation*}
\Theta^\arithm(fv)=f\Theta^\arithm(v)+t\frac{\partial f}{\partial t}v\;.
\end{equation*}
Thus, after choosing a basis of $V$, we have
\begin{equation*}
\Theta^\arithm=A+t\frac{\partial}{\partial t}
\end{equation*}
for some $m\times m$-matrix $A$ of global sections of $\mathbb{B}_\dR^{+, \dagger, \la}$, where $m$ is the rank of $V$, and our claim is that the operator
\begin{equation}
\label{eq:htdr-aplusiinvertible}
(A+i)+t\frac{\partial}{\partial t}: (\mathbb{B}_\dR^{+, \dagger, \la})^{\oplus m}\rightarrow (\mathbb{B}_\dR^{+, \dagger, \la})^{\oplus m}
\end{equation}
is invertible for all $i$ as in the first paragraph.

To this end, writing 
\begin{equation*}
A=A_0+A_1t+A_2t^2+\dots
\end{equation*}
for $m\times m$-matrices $A_\ell$ of global functions on $X_\infty^\la\times\GSpec\Q_p(\zeta_{p^\infty})$, our first claim is the following:

\bigskip

\textbf{Claim 1.} For all $i$ as in the first paragraph, the matrix $A_0+i$ is invertible.

\bigskip

\textit{Proof of the claim.} Recall from the proof of \cref{prop:padicht-scholzefunctorlocal} that the $\mathbb{B}_\dR^+$-local system associated to $(\Fil^\bullet W, \nabla)$ coincides with the base change of the $\O_{X_\infty^\la}(\zeta_{p^\infty})[t]$-module
\begin{equation}
\label{eq:drbundles-keylemformulav}
\bigoplus_\bullet t^{-\bullet}(\Fil^\bullet W\tensor_{\O_X} \O_{X_\infty^\la} \tensor_{\Q_p} \Q_p(\zeta_{p^\infty}))\;,
\end{equation}
where $t$ acts via the transition maps $t^{-\bullet}\Fil^\bullet W\rightarrow t^{-\bullet+1} \Fil^{\bullet-1} W$, along the map $\O_{X_\infty^\la}(\zeta_{p^\infty})[t]\rightarrow\mathbb{B}_{\dR, X_\infty^\cycl}^+$ induced by the given toric chart. Here, (\ref{eq:drbundles-keylemformulav}) is equipped with an action of $\Z_p^\times\rtimes\Z_p^n$ by letting $\Z_p^n$ act canonically on $\O_{X_\infty^\la}$ while $\Z_p^\times$ acts on $t$ by multiplication and on $\Q_p(\zeta_{p^\infty})$ via the usual Galois action. As we have assumed that all $\Fil^\bullet W$ are free and that the transition maps are the standard inclusions, this means that the image of $(\Fil^\bullet W, \nabla)$ is isomorphic to $(\mathbb{B}_{\dR, X_\infty^\cycl}^+)^{\oplus m}$ on the localised site $X_\proet/X_\infty^\cycl$ and that, under this isomorphism, the action of $\Z_p^\times$  is given by $\gamma\in\Z_p^\times$ acting on the standard basis via the diagonal matrix $\mathrm{diag}(\gamma^{-a_1}, \dots, \gamma^{-a_m})$, where the $a_k$ are the indices of the jumps in the filtration $\Fil^\bullet W$ counted with multiplicity. In particular, note that $a_k<i$.

By assumption, there is an invertible $m\times m$-matrix $S$ of sections of $\mathbb{B}_\dR^+$ over $X_\infty^\cycl$ such that the map
\begin{equation*}
V\tensor_{\mathbb{B}_\dR^{+, \dagger, \la}}\mathbb{B}^+_{\dR, X^\cycl_\infty}\cong (\mathbb{B}^+_{\dR, X^\cycl_\infty})^{\oplus m}\overset{S}{\longrightarrow} (\mathbb{B}^+_{\dR, X^\cycl_\infty})^{\oplus m}
\end{equation*}
is $\Z_p^n\rtimes\Z_p^\times$-equivariant, where the right-hand side denotes the restriction of the image of $(\Fil^\bullet W, \nabla)$ under Scholze's functor to $X_\proet/X_\infty^\cycl$ with the induced $\Z_p^n\rtimes\Z_p^\times$-action while the isomorphism on the left-hand side uses our chosen basis of $V$. Writing $S=(s_{ij})_{i, j}$ and $S^{-1}=(s'_{ij})_{i, j}$, this means that the action of $\gamma\in\Z_p^\times$ on the $k$-th basis vector $e_k$ of $V$ is given by
\begin{equation}
\label{eq:drbundles-keylemactgamma}
\begin{split}
\gamma.e_j&=S^{-1}\begin{pmatrix} \gamma^{-a_1} & & \\
& \ddots & \\
& & \gamma^{-a_m}\end{pmatrix} Se_j=S^{-1}\begin{pmatrix} \gamma^{-a_1} & & \\
& \ddots & \\
& & \gamma^{-a_m}\end{pmatrix} \begin{pmatrix} s_{1j} \\ \vdots \\ s_{mj}\end{pmatrix} \\
&=S^{-1}\begin{pmatrix} \gamma^{-a_1}\cdot (\gamma.s_{1j}) \\ \vdots \\ \gamma^{-a_m}\cdot (\gamma.s_{mj})\end{pmatrix}=\left(\sum_k s'_{ik}\gamma^{-a_k}(\gamma.s_{kj})\right)_i\;.
\end{split}
\end{equation}

Since the action of $\Z_p^\times$ on $V$ is locally analytic, all entries of the last vector must be locally analytic functions in $\gamma$ and, in particular, the functions $\gamma\mapsto \gamma.s_{kj}$ must be locally analytic: Indeed, multiplying the last vector by the matrix of constants $S$ from the left yields the vector $(\gamma^{-a_k}\cdot (\gamma.s_{kj}))_k$, all of whose entries must thus be locally analytic functions in $\gamma$. As taking the derivative of the action of $\Z_p^{\times, \la}$ at $\gamma=1$ yields the action of $\Theta^\arithm$, we obtain
\begin{equation}
\label{eq:htdr-thetaarithmej}
\Theta^\arithm(e_j)=\left(\sum_k (-a_k)s'_{ik}s_{kj}+s'_{ik}\left.\frac{\mathrm{d}(\gamma.s_{kj})}{\mathrm{d}\gamma}\right\vert_{\gamma=1}\right)_i
\end{equation}
from (\ref{eq:drbundles-keylemactgamma}), i.e.\ the right-hand side is the $j$-th column of the matrix $A$.

Now note that $\mathbb{B}_\dR^+(X_\infty^\cycl)/t\cong \O(X_\infty^\cycl)$ and recall that there is a canonical map $\O(X_\infty^\la)\rightarrow\mathbb{B}_\dR^+(X_\infty^\cycl)$. Thus, after passing to a rational localisation of $X$ so that $X$ becomes affine, we may approximate $s_{kj}$ as
\begin{equation*}
s_{kj}=r_{kj, N}+p^N(\,\dots)+t(\,\dots)
\end{equation*}
for any $N\geq 0$ and some $r_{kj, N}\in \O(X_\infty^\la)$. As the action of $\Z_p^\times$ on $X_\infty^\la$ is smooth and $\gamma.t=\gamma\cdot t$, this yields
\begin{equation*}
\left.\frac{\mathrm{d}(\gamma.s_{kj})}{\mathrm{d}\gamma}\right\vert_{\gamma=1}=p^N(\,\dots)+t(\,\dots)\;,
\end{equation*}
where we note that the implicit terms in brackets might now be different ones than in the previous equation. In other words, we may conclude from (\ref{eq:htdr-thetaarithmej}) that the $j$-th column of $A_0$ coincides with $(-\sum_k a_k s'_{ik}s_{kj})_i$ up to a multiple of $p^N$. Since the ring $\O(X_\infty^\la\times\GSpec\Q_p(\zeta_{p^\infty}))$ is $p$-adically separated, this yields
\begin{equation*}
A_0=\left(-\sum_k a_k s'_{ik}s_{kj}\right)_{i, j}=S^{-1}\begin{pmatrix} -a_1 & & \\ & \ddots & \\ & & -a_m\end{pmatrix} S\;,
\end{equation*}
hence
\begin{equation*}
A_0+i=S^{-1}\begin{pmatrix} -a_1+i & & \\ & \ddots & \\ & & -a_m+i\end{pmatrix} S\;.
\end{equation*}

Now we are done: Indeed, the above yields
\begin{equation*}
\det(A_0+i)=\prod_{k=1}^m (i-a_k)
\end{equation*}
and, by assumption, we have $a_k<i$ for all $k$. We conclude that $\det(A_0+i)$ is a unit in the $\Q$-algebra $\O(X_\infty^\la\times\GSpec\Q_p(\zeta_{p^\infty}))$ and thus the claim follows. \hfill\qed

\bigskip

Recall from (\ref{eq:htdr-aplusiinvertible}) that we need to show that the operator
\begin{equation*}
\Theta^\arithm+i=(A+i)+t\frac{\partial}{\partial t}: (\mathbb{B}_\dR^{+, \dagger, \la})^{\oplus m}\rightarrow (\mathbb{B}_\dR^{+, \dagger, \la})^{\oplus m}
\end{equation*}
is invertible. We start with injectivity and, to this end, take some
\begin{equation*}
x=\sum_\ell x_\ell t^\ell\in (\mathbb{B}_\dR^{+, \dagger, \la})^{\oplus m}\;,
\end{equation*}
where each $x_\ell$ is an $m$-tuple of local sections of the pushforward of the structure sheaf along $X_\infty^\la\times\GSpec\Q_p(\zeta_{p^\infty})\rightarrow X$ viewed as a sheaf on $|X|$. Let $\ell$ be minimal such that $x_\ell\neq 0$. Then
\begin{equation*}
(\Theta^\arithm+i)(x)=(A_0+i)x_\ell t^\ell+\ell x_\ell t^\ell+t^{\ell+1}(\,\dots)=(A_0+(i+\ell))x_\ell t^\ell+t^{\ell+1}(\,\dots)
\end{equation*}
and by applying Claim 1 to $i+\ell$ in place of $i$, we see that $(A_0+(i+\ell))x_\ell$ and consequently $(\Theta^\arithm+i)(x)$ is nonzero.

For surjectivity, let $y=\sum_\ell y_\ell t^\ell\in (\mathbb{B}_\dR^{+, \dagger, \la})^{\oplus m}$. By Claim 1, we can then recursively solve the equations
\begin{equation}
\label{eq:htdr-recursiveeqxell}
(A_0+i+\ell)x_\ell+\sum_{k<\ell} A_{\ell-k}x_k=y_\ell
\end{equation}
for $\ell=0, 1, 2, \dots$ and we claim that then
\begin{equation*}
x\coloneqq \sum_\ell x_\ell t^\ell
\end{equation*}
is in $(\mathbb{B}_\dR^{+, \dagger, \la})^{\oplus m}$ and satisfies $(\Theta^\arithm+i)(x)=y$. Indeed, once we know the former, the latter is immediate since the coefficient of $t^\ell$ in the expansion of $(\Theta^\arithm+i)(x)$ is exactly given by the left-hand side of (\ref{eq:htdr-recursiveeqxell}). To establish that $x\in (\mathbb{B}_\dR^{+, \dagger, \la})^{\oplus m}$, we need to understand the growth of the $x_\ell$, for which we begin by proving the following explicit formula:

\bigskip

\textbf{Claim 2.} For all $\ell\geq 0$, we have
\begin{equation*}
x_\ell=\sum_{k\leq \ell} \left(\left(\sum_{\lambda\vdash (\ell-k)} (-1)^{\ell(\lambda)}\prod_j A_{\lambda_j}(A_0+i+k+|\lambda|_{\leq j})^{-1}\right)\cdot (A_0+i+k)^{-1}y_k\right)\;,
\end{equation*}
where the inner sum runs over all ordered (!) partitions $\lambda=(\lambda_1, \lambda_2, \dots, \lambda_s)$ of $\ell-k$ into positive parts, $\ell(\lambda)\coloneqq s$ denotes the length of such a partition and we set
\begin{equation*}
|\lambda|_{\leq j}\coloneqq \lambda_1+\dots+\lambda_j
\end{equation*}
for any $1\leq j\leq s$.

\bigskip

\textit{Proof of the claim.} We use induction on $\ell$, in the base case $\ell=0$ the claimed formula becomes $x_0=(A_0+i)^{-1}y_0$, which is true (note that $0$ admits the empty partition!). For the induction step, note that (\ref{eq:htdr-recursiveeqxell}) implies that
\begin{equation*}
x_\ell=(A_0+i+\ell)^{-1}\left(y_\ell-\sum_{k<\ell} A_{\ell-k}x_k\right)\;.
\end{equation*}
Plugging in the formulas for the $x_k$ with $k<\ell$, we see that the right-hand side is a linear combination of the $y_k$ with $k\leq \ell$. For $k=\ell$, the coefficient of $y_\ell$ is given by $(A_0+i+\ell)^{-1}$, which is in keeping with the claimed formula (again because the empty partition of $0$ is allowed), while the coefficient of $y_k$ for $k<\ell$ is given by
\begin{equation*}
-(A_0+i+\ell)^{-1}\sum_{k\leq k'<\ell} A_{\ell-k'}\left(\left(\sum_{\lambda\vdash (k'-k)} (-1)^{\ell(\lambda)}\prod_j A_{\lambda_j}(A_0+i+k+|\lambda|_{\leq j})^{-1}\right)\cdot (A_0+i+k)^{-1}\right)\;.
\end{equation*}
Noting that there is a bijection
\begin{equation*}
\begin{split}
\{\text{pairs $(k', \lambda)$ with $k\leq k'<\ell$ and $\lambda\vdash (k'-k)$}\}&\xleftrightarrow{1:1} \{\mu\vdash (\ell-k)\} \\
(k', (\lambda_1, \dots, \lambda_s))&\;\,\mapsto\;\, \mu\coloneqq (\lambda_1, \dots, \lambda_s, \ell-k')\;,
\end{split}
\end{equation*}
we see that the above expression simplifies to 
\begin{equation*}
\left(\sum_{\mu\vdash (\ell-k)} (-1)^{\ell(\mu)}\prod_j A_{\mu_j}(A_0+i+k+|\mu|_{\leq j})^{-1}\right)(A_0+i+k)^{-1},
\end{equation*}
which is exactly the coefficient of $y_k$ in the claimed formula. \hfill\qed

\bigskip

In order to show $x\in (\mathbb{B}_\dR^{+, \dagger, \la})^{\oplus m}$, we have to argue that there is some $N>0$ such that the local sections $p^{\ell N}x_\ell$ forms a nullsequence. Since $A$ is a matrix over $\mathbb{B}_\dR^{+, \dagger, \la}$ and $y\in (\mathbb{B}_\dR^{+, \dagger, \la})^{\oplus m}$, we know that there is some $M>0$ with $p^{\ell M}A_\ell\rightarrow 0$ and $p^{\ell M}y_\ell\rightarrow 0$; as the rings we are working in are bounded, after possibly enlarging $M$, we may even assume that, for $\ell\geq 1$, the sequence $(p^{\ell M}A_\ell)_\ell$ is a nullsequence in the subring of powerbounded elements. Observing that
\begin{equation*}
p^{\ell M}x_\ell=\sum_{k\leq \ell} \left(\left(\sum_{\lambda\vdash (\ell-k)} (-1)^{\ell(\lambda)}\prod_j p^{\lambda_j M}A_{\lambda_j}(A_0+i+k+|\lambda|_{\leq j})^{-1}\right)\cdot (A_0+i+k)^{-1}p^{kM} y_k\right)
\end{equation*}
by Claim 2, it thus only remains to control the inverses occurring in the formula.

However, note that $\det(A_0+i+k+|\lambda|_{\leq j})(A_0+i+k+|\lambda|_{\leq j})^{-1}$ is a matrix whose entries are degree $m$ integer polynomials in the entries of $A_0$ by Cramer's rule. Further observing that each summand in the formula above only contains at most $\ell+1$ such inverses as factors, we are thus reduced to controlling the determinants $\det(A_0+i+k+|\lambda|_{\leq j})$. Recalling that these determinants are integers by the proof of Claim 1, we have to show that the $p$-adic valuation of the integers
\begin{equation*}
\det(A_0+i+k)\prod_j \det(A_0+i+k+|\lambda|_{\leq j})
\end{equation*}
is uniformly linearly bounded in $\ell$, where $k\leq\ell$ and $\lambda\vdash(\ell-k)$. 

Recall, though, that the proof of Claim 1 gives an explicit expression for these determinants; namely, the above product is equal to
\begin{equation*}
\prod_{r=1}^m \left((i+k-a_r)\prod_j (i+k+|\lambda|_{\leq j}-a_r)\right)\;.
\end{equation*}
Observing that $1\leq |\lambda|_{\leq 1}<|\lambda|_{\leq 2}<\dots\leq \ell-k$, we see that the above is a divisor of
\begin{equation*}
\prod_{r=1}^m (i+\ell-a_r)!
\end{equation*}
upon recalling that $i>a_r$ for all $r$. As the $p$-adic valuation of this latter product is bounded by $m\cdot (i+\ell-a)/(p-1)$, where $a$ denotes the minimum of the $a_r$, we win.
\end{proof}

Finally, we can finish the proof of \cref{prop:drbundles-xhtdagger+pullback}. In particular, recall that this also immediately proves \cref{thm:drbundles-main} by the discussion preceding \cref{prop:drbundles-xhtdagger+pullback}.

\begin{proof}[Proof of \cref{prop:drbundles-xhtdagger+pullback}]
By compatibility of $X\mapsto X^{\HT, \dagger, +}$ with open localisations, see \cref{prop:defis-openloc}, and since generic flatness may be checked locally on $X$, we may assume that $X$ admits a rigid étale map $X\rightarrow\ol{\T}^n$. In this case, recall that we already know that the functor $\Vect(X^{\HT, \dagger, +})\rightarrow \Vect(X^{\HT, \dagger})$ factors through the full subcategory spanned by generically flat vector bundles on $X^{\HT, \dagger}$ by \cref{thm:padicht-scholzefunctor}. 

If we can show that the functor
\begin{equation}
\label{eq:drbundles-genflatxhtdaggertobdr+}
\{\text{generically flat vector bundles on $X^{\HT, \dagger}$}\}\rightarrow \{\text{$\mathbb{B}_\dR^+$-local systems on $X_\proet$}\}
\end{equation}
obtained by restricting (\ref{eq:drbundles-xhtdaggertobdr+}) to generically flat vector bundles is fully faithful, we are done: Full faithfulness then follows from the fact that the composite functor
\begin{equation}
\begin{split}
\label{eq:drbundles-compositefunctorgenflat}
\Vect(X^{\dR, +})\cong \Vect(X^{\HT, \dagger, +})&\rightarrow \{\text{generically flat vector bundles on $X^{\HT, \dagger}$}\} \\
&\rightarrow \{\text{$\mathbb{B}_\dR^+$-local systems on $X_\proet$}\}
\end{split}
\end{equation}
agrees with Scholze's functor (\ref{eq:proet-scholzefunctor}) by \cref{thm:padicht-scholzefunctor}, which is known to be fully faithful, see \cite[Thm.\ 7.6]{PAdicHodgeTheory}. We also get essential surjectivity: If $E$ is a generically flat vector bundle on $X^{\HT, \dagger}$, its associated $\mathbb{B}_\dR^+$-local system $\mathbb{M}$ is associated to a filtered vector bundle with connection $(\Fil^\bullet W, \nabla)$ on $X$. Then $(\Fil^\bullet W, \nabla)$ corresponds to a vector bundle $E'$ on $X^{\dR, +}$, whose image under the composite functor (\ref{eq:drbundles-compositefunctorgenflat}) is known to be isomorphic to $\mathbb{M}$ by \cref{thm:padicht-scholzefunctor}. As the last functor in (\ref{eq:drbundles-compositefunctorgenflat}) is fully faithful, this means that the image of $E'$ under 
\begin{equation*}
\Vect(X^{\dR, +})\cong \Vect(X^{\HT, \dagger, +})\rightarrow \{\text{generically flat vector bundles on $X^{\HT, \dagger}$}\}
\end{equation*}
is already isomorphic to $E$ and then the desired preimage of $E$ is the vector bundle on $X^{\HT, \dagger, +}$ corresponding to $E'$ via the equivalence $\Vect(X^{\dR, +})\cong \Vect(X^{\HT, \dagger, +})$.

To prove that (\ref{eq:drbundles-genflatxhtdaggertobdr+}) is fully faithful, first note that it is symmetric monoidal and compatible with taking duals. Thus, checking that $\Hom(F, E)\cong \Hom(\mathbb{N}, \mathbb{M})$ for $E, F\in\Vect(X^{\HT, \dagger})$ generically flat and $\mathbb{M}, \mathbb{N}$ their associated $\mathbb{B}_\dR^+$-local systems reduces to the case $F=\O$ by replacing $E$ by $E\tensor F^\vee$. In other words, it suffices to show that, for any generically flat vector bundle $E$ on $X^{\HT, \dagger}$, we have
\begin{equation}
\label{eq:htdr-cohxhtdaggerbdr+iso}
R\Gamma(X^{\HT, \dagger}, E)\cong R\Gamma(X_\proet, \mathbb{M})\;,
\end{equation}
where $\mathbb{M}$ is the $\mathbb{B}_\dR^+$-local system associated to $E$. Using the notation from \cref{cor:htdr-cohomologyxhtdrdagger+}, the left-hand side is given by the (underived) $\Z_p^n\rtimes\Z_p^\times$-invariants of the cohomology on $X_\infty^\la$ of the total complex of
\begin{equation*}
\begin{tikzcd}
V\ar[r, "\nabla"]\ar[d, "D", swap] & V\tensor_{\O_{X_\infty^\la}} \Omega^1_{X_\infty^\la}(-1)\ar[d, "D", swap]\ar[r, "\nabla"] & V\tensor_{\O_{X_\infty^\la}} \Omega^2_{X_\infty^\la}(-2)\ar[r, "\nabla"]\ar[d, "D", swap] & \dots \\
V\ar[r, "\nabla"] & V\tensor_{\O_{X_\infty^\la}} \Omega^1_{X_\infty^\la}(-1)\ar[r, "\nabla"] & V\tensor_{\O_{X_\infty^\la}} \Omega^2_{X_\infty^\la}(-2)\ar[r, "\nabla"] & \dots\nospacepunct{\;,}
\end{tikzcd}
\end{equation*}
where $D$ is just given by $\Theta^\arithm$, and then \cref{lem:drbundles-keylem} shows that the $t$-adic filtration on this total complex is complete. Moreover, since $\mathbb{B}_\dR^+$ and hence $\mathbb{M}$ is $t$-adically complete, also the filtration
\begin{equation*}
\dots\rightarrow R\Gamma(X_\proet, t^2 \mathbb{M})\rightarrow R\Gamma(X_\proet, t \mathbb{M})\rightarrow R\Gamma(X_\proet, \mathbb{M})
\end{equation*}
is complete.

Thus, we may check the isomorphism (\ref{eq:htdr-cohxhtdaggerbdr+iso}) on associated graded pieces. Each graded piece of $t^\bullet V$ will be a vector bundle $V'$ on $X_\infty^\la\times\GSpec\Q_p(\zeta_{p^\infty})$ equipped with a semilinear locally analytic $\Z_p^n\rtimes\Z_p^\times$-action and the corresponding graded piece of $t^\bullet\mathbb{M}$ will be the pullback of $V'$ to a locally free $\widehat{\O}_X$-module $\mathbb{M}'$ on $X_\proet$. Our claim then amounts to showing that the (underived) $\Z_p^n\rtimes\Z_p^\times$-invariants of the cohomology on $X_\infty^\la$ of the total complex of 
\begin{equation*}
\begin{tikzcd}
V'\ar[r, "\nabla"]\ar[d, "D", swap] & V'\tensor_{\O_{X_\infty^\la}} \Omega^1_{X_\infty^\la}(-1)\ar[d, "D", swap]\ar[r, "\nabla"] & V'\tensor_{\O_{X_\infty^\la}} \Omega^2_{X_\infty^\la}(-2)\ar[r, "\nabla"]\ar[d, "D", swap] & \dots \\
V'\ar[r, "\nabla"] & V'\tensor_{\O_{X_\infty^\la}} \Omega^1_{X_\infty^\la}(-1)\ar[r, "\nabla"] & V'\tensor_{\O_{X_\infty^\la}} \Omega^2_{X_\infty^\la}(-2)\ar[r, "\nabla"] & \dots\nospacepunct{\;,}
\end{tikzcd}
\end{equation*}
where $\nabla$ and $D$ arise from the Lie algebra actions of $\Z_p^n$ and $\Z_p^\times$, respectively, compute $R\Gamma(X_\proet, \mathbb{M}')$. However, this is now an instance of geometric Sen theory, see \cite{GeometricSenTheory}: taking locally analytic vectors for the $\Z_p^n\rtimes\Z_p^\times$-action on $\mathbb{M}'$ recovers $V'$ by \cite{AnPrism} and $\nabla$ is the corresponding geometric Sen operator while $D$ is the arithmetic Sen operator of the Galois representation $R\Gamma(X_{\C_p, \proet}, \mathbb{M}')$.
\end{proof}

\begin{rem}
\label{rem:htdr-checkgenflat}
One can also define a notion of generic flatness for $\mathbb{B}_\dR^{+, \dagger}$-local systems on $X_\proet$ using the functor from filtered vector bundles with connection on $X$ to $\mathbb{B}_\dR^{+, \dagger}$-local systems from \cref{rem:drbundles-functortobdr+dagger}. It is then not a priori clear that generic flatness of a $\mathbb{B}_\dR^{+, \dagger}$-local system may be checked after extending scalars to $\mathbb{B}_\dR^+$. However, using that vector bundles on $X^{\HT, \dagger}$ and $\mathbb{B}_\dR^{+, \dagger}$ are equivalent via pullback by \cite{AnPrism}, one should be able to use the isomorphism (\ref{eq:htdr-cohxhtdaggerbdr+iso}) from the proof above to check that this is indeed the case.
\end{rem}

\newpage

\section{Geometric syntomic cohomology}

Finally, we want to comment on how to apply the methods we have developed to the geometric case, i.e.\ to the cohomology of rigid spaces over $\C_p$. For this, first recall the situation for the de Rham stack: By \cite[Ex.\ 4.7.6.(3)]{dRFF}, we have $\C_p^\dR\cong \GSpec\ol{\Q}_p$ and hence just taking the de Rham stack cannot compute the usual de Rham cohomology of smooth partially proper rigid spaces over $\C_p$. Indeed, the problem is that the stack $X^\dR$ always computes ``absolute de Rham cohomology'', i.e.\ de Rham cohomology with respect to the base $\Q_p$, while the usual de Rham cohomology of rigid spaces over $\C_p$ is rather relative to $\C_p$. In other words, if $X$ is a Gelfand stack over $\C_p$, we should define an appropriate version $X^{\dR/\C_p}$ of $X^\dR$ ``relative to $\C_p$''. Meditating on the situation, one is eventually led to defining $X^{\dR/\C_p}$ via the pullback square
\begin{equation*}
\begin{tikzcd}
X^{\dR/\C_p}\ar[r]\ar[d] & \GSpec\C_p\ar[d] \\
X^\dR\ar[r] & \C_p^\dR\nospacepunct{\;,}
\end{tikzcd}
\end{equation*}
where $\GSpec\C_p\rightarrow\C_p^\dR$ is the canonical map, see \cite[Def.\ 4.5.11.(2)]{dRFF}. One then proves as before that the cohomology of this stack calculates de Rham cohomology if $X$ is a smooth partially proper rigid space over $\C_p$, see \cite[Prop.\ 5.2.1]{dRFF}.

Thus, the key issue in defining syntomic cohomology in the geometric situation, i.e.\ for Berkovich smooth derived Berkovich spaces over $\C_p$, is to identify the ``correct'' base change of $X^\Syn$ or $X^\N$. Below, we will explain that there is a canonical map
\begin{equation*}
\Fil Y_{\C_p}\coloneqq (\{ut=\phi^{-1}(\xi)\}\subseteq Y_{\C_p}\times\ol{\DD}_+\times\ol{\DD}_-)\,/\,\ol{\T}\rightarrow \C_p^\N\;,
\end{equation*}
where $\xi=p-[p^\flat]$, and that base changing $X^\N$ along this map yields an appropriate definition for the filtered prismatisation $X^{\N/\C_p}$ of $X$ relative to $\C_p$. We will justify this claim by showing that basically everything we have done in the ``arithmetic case'', i.e.\ for smooth partially proper rigid space over $\Q_p$, has an analogue in the geometric case with the definition of $X^{\N/\C_p}$ proposed above.

In particular, this includes the basic description of the geometry of the filtered prismatisation that we have discussed in §\ref{subsect:loci} and the local presentations of the filtered prismatisation obtained in §\ref{sect:pres}. This will enable us to transfer the results from §§\ref{sect:hkcomp}, \ref{sect:proet} and their proofs almost verbatim to the geometric setting. In particular, we recover versions of the comparison results between syntomic cohomology over $\C_p$ and filtered Hyodo--Kato cohomology and between syntomic cohomology over $\C_p$ and proétale cohomology from \cite[Thm.\ 1.3]{BasicComparison}. However, we warn the reader that, in the geometric case, syntomic cohomology as we define it will \emph{not} agree with the more classical definition as it appears e.g.\ in the work of Colmez--Nizio{\l}, see \cite{BasicComparison}, but will rather recover a slight variant that has been introduced by Bosco in \cite{BoscopAdic}, where it is called \emph{syntomic Fargues--Fontaine cohomology}, see also \cref{rem:geom-syntomicffcoh}.

On the way, we will obtain a stacky interpretation of the $B_\dR^+$-cohomology for smooth partially proper rigid spaces over $\C_p$ as e.g.\ defined in \cite[§3]{BasicComparison}, \cite{GuoBdR} or \cite[§13]{IntegralpAdicHT} through a suitable base change of the analytic de Rham stack, and this will enable us to give a simple proof for the comparison between $B_\dR^+$-cohomology and de Rham cohomology. Moreover, similarly to what happened in §\ref{sect:padicht}, we will be able to recover the geometric de Rham comparison theorem from \cite[Thm.\ 13.1]{IntegralpAdicHT}, see \cref{thm:geom-drcomp}.

\subsection{Stacks attached to $p$-adic cohomology theories relative to $\C_p$}

We start by recalling some definitions of stacks computing $p$-adic cohomology theories relative to $\C_p$. As already discussed above, the main question throughout is along which map one should base change. To this end, recall from \cref{prop:defis-xprismdr} that there is a canonical map $Y_{X^\diamond}\rightarrow X^\prism$ for any Gelfand stack $X$ and that, similarly, we have canonical maps $\FF_{X^\diamond}\rightarrow X^{\Div^1}$ and $\FF_{X^\diamond}\rightarrow X^\HK$.

\begin{defi}
Let $X$ be a Gelfand stack over $\C_p$. The \emph{relative analytic prismatisation} $X^{\prism/\C_p}$ of $X$ over $\C_p$ is defined by the pullback diagram
\begin{equation*}
\begin{tikzcd}
X^{\prism/\C_p}\ar[r]\ar[d] & Y_{\C_p}\ar[d] \\
X^\prism\ar[r] & \C_p^\prism\nospacepunct{\;.}
\end{tikzcd}
\end{equation*}
Moreover, the \emph{relative Hyodo--Kato stack} $X^{\HK/\C_p}$ of $X$ over $\C_p$ and the stack $X^{\Div^1/\C_p}$ are defined by the pullback diagrams
\begin{equation*}
\begin{tikzcd}
X^{\HK/\C_p}\ar[r]\ar[d] & \FF_{\C_p}\ar[d] \\
X^\HK\ar[r] & \C_p^\HK\nospacepunct{\;,}
\end{tikzcd}
\hspace{1cm}
\begin{tikzcd}
X^{\Div^1/\C_p}\ar[r]\ar[d] & \FF_{\C_p}\ar[d] \\
X^{\Div^1}\ar[r] & \C_p^{\Div^1}\nospacepunct{\;.}
\end{tikzcd}
\end{equation*}
\end{defi}

We remark that $X^{\prism/\C_p}$ is still equipped with a radius map $\kappa: X^{\prism/\C_p}\rightarrow (0, \infty)$ obtained by postcomposing the canonical map $X^{\prism/\C_p}\rightarrow X^\prism$ with the usual radius map on $X^\prism$. Moreover, base changing \cref{prop:defis-prismffdr} shows that the natural map 
\begin{equation*}
X^{\prism/\C_p}\rightarrow (X^{\prism/\C_p})^{\dR/Y_{\C_p}}\cong Y_{X^\diamond}^{\dR/Y_{\C_p}}
\end{equation*}
is an isomorphism over the locus $(1, \infty)$, where $(X^{\prism/\C_p})^{\dR/Y_{\C_p}}$ is defined as the pullback of $(X^{\prism/\C_p})^\dR$ along the canonical map $Y_{\C_p}\rightarrow Y_{\C_p}^\dR$; note that Frobenius is an isomorphism on $Y_{X^\diamond}^{\dR/Y_{\C_p}}$ and that $Y_{X^\diamond}^{\dR/Y_{\C_p}}/\phi^\Z\cong X^{\HK/\C_p}$. Furthermore, from \cref{prop:defis-xprismxdiv1}, we deduce that Frobenius induces an isomorphism $X^{\prism/\C_p}_{(0, 1)}\cong X^{\prism/\C_p}_{(0, p)}$ and that
\begin{equation*}
\operatorname{coeq}(\hspace{-0.15cm}
\begin{tikzcd}
X^{\prism/\C_p}_{(0, 1)}\ar[r,shift left=.75ex, "\phi"]
  \ar[r,shift right=.75ex,swap, "\id"] & X^{\prism/\C_p}_{(0, p)}
\end{tikzcd}
\hspace{-0.15cm})\cong X^{\Div^1/\C_p}\nospacepunct{\;.}
\end{equation*}

Let us also note that, as in \cref{cor:recall-hkfin}, one proves that pushforward along $X^\HK\rightarrow \C_p^\HK$ preserves perfect complexes for any smooth partially proper qcqs rigid space $X$ over $\C_p$. In particular, the pushforward of any perfect complex along this map identifies with a $(\phi, N)$-module in $\Perf(\Q_p^\un)$ by \cite[Thm.\ 7.1.1]{dRFF} and, as in the arithmetic case, this pushforward is expected to compute Hyodo--Kato cohomology of $X$ relative to $\C_p$, see \cite[Rem.\ 6.1.3]{dRFF}.

\begin{defi}
Let $X$ be a smooth partially proper qcqs rigid space over $\C_p$. We define the \emph{Hyodo--Kato cohomology} $R\Gamma_\HK(X)$ of $X$ relative to $\C_p$ as the $(\phi, N)$-module in $\Perf(\Q_p^\un)$ corresponding to the pushforward of the structure sheaf along $X^\HK\rightarrow\C_p^\HK$.
\end{defi}

We now turn to the filtered prismatisation. As already announced above, our aim will be to base change $X^\N$ for $X$ a Gelfand stack over $\C_p$ from $\C_p^\N$ to the Gelfand stack
\begin{equation*}
\Fil Y_{\C_p}\coloneqq (\{ut=\phi^{-1}(\xi)\}\subseteq Y_{\C_p}\times\ol{\DD}_+\times\ol{\DD}_-)\,/\,\ol{\T}\;,
\end{equation*}
where $t$ denotes the coordinate on $\ol{\DD}_+$, $u$ is the coordinate on $\ol{\DD}_-$ and $\xi=p-[p^\flat]$. Intuitively, one should think of this stack as the Rees stack associated to the ``$\phi^{-1}(\xi)$-adic filtration on $Y_{\C_p}$'' and we urge the reader to compare this definition to the presentation of the algebraic Nygaardification of $\O_{\C_p}$ in the sense of Drinfeld and Bhatt--Lurie from \cite[Ex.\ 5.5.6]{FGauges}.

Our task will be to construct a map $\Fil Y_{\C_p}\rightarrow \C_p^\N$ along which to base change. For this, first observe that $\Fil Y_{\C_p}$ admits a canonical map to $\ol{\DD}_+/\ol{\T}\times (\ol{\DD}_-/\ol{\T})^\dR$ via $t$ and $u$. Moreover, there is an obvious projection map $\Fil Y_{\C_p}\rightarrow Y_{\C_p}$ and postcomposing with the canonical map $Y_{\C_p}\rightarrow\C_p^\prism$ thus yields a map $\Fil Y_{\C_p}\rightarrow \C_p^\prism$. Finally, unwinding definitions, we see that the equality $ut=\phi^{-1}(\xi)$ on $\Fil Y_{\C_p}$ ensures that the composition
\begin{equation*}
\Fil Y_{\C_p}\rightarrow \C_p^\prism\rightarrow\Q_p^\prism\xrightarrow{\widetilde{\mu}} (\ol{\DD}_+/\ol{\T})^\dR
\end{equation*}
identifies with the composition
\begin{equation*}
\Fil Y_{\C_p}\rightarrow \ol{\DD}_+/\ol{\T}\times (\ol{\DD}_-/\ol{\T})^\dR\xrightarrow{\mathrm{mult}} (\ol{\DD}_+/\ol{\T})^\dR
\end{equation*}
and hence we overall obtain a map
\begin{equation*}
\Fil Y_{\C_p}\rightarrow \Q_p^\N\times_{\Q_p^\prism} \C_p^\prism\;.
\end{equation*}

To upgrade this to a map to $\C_p^\N$, note that $\Fil Y_{\C_p}$ also maps to $\C_p^{\Cone}$ over $\ol{\DD}_+/\ol{\T}$. Indeed, using the isomorphism 
\begin{equation*}
(\{\phi^{-1}(\xi)=0\}\subseteq Y_{\C_p})\cong (\{\xi=0\}\subseteq Y_{\C_p})\cong \GSpec\C_p\;,
\end{equation*}
which is induced by Frobenius, we have
\begin{equation*}
(\{t=0\}\subseteq \Fil Y_{\C_p})\cong \GSpec\C_p\times\ol{\DD}_-\,/\,\ol{\T}
\end{equation*}
and the stack on the right-hand side admits a canonical map to $\GSpec\C_p$, which overall induces a map $\Fil Y_{\C_p}\rightarrow \C_p^{\Cone}$. Recalling from \cref{lem:syn-xnviacone} that $\C_p^\N$ may be expressed as the pullback of $\C_p^\prism$ along $\C_p^{\Cone}\rightarrow (\C_p^{\Cone})^\dR$, the preceding discussion yields an induced map
\begin{equation*}
\Fil Y_{\C_p}\rightarrow \C_p^\N
\end{equation*}
upon realising that the maps $\Fil Y_{\C_p}\rightarrow \C_p^\prism$ and $\Fil Y_{\C_p}\rightarrow \C_p^{\Cone}$ we have given are compatible upon mapping to $(\C_p^{\Cone})^\dR$. Using this, we make the following definition:

\begin{defi}
Let $X$ be a Gelfand stack over $\C_p$. The \emph{relative analytic Nygaardification} $X^{\N/\C_p}$ of $X$ over $\C_p$ is defined by the pullback diagram
\begin{equation*}
\begin{tikzcd}
X^{\N/\C_p}\ar[r]\ar[d] & \Fil Y_{\C_p}\ar[d] \\
X^\N\ar[r] & \C_p^\N\nospacepunct{\;.}
\end{tikzcd}
\end{equation*}
\end{defi}

As in §\ref{subsect:loci}, let us discuss the restriction of $X^{\N/\C_p}$ to various interesting loci. As we have already done this for $X^\N$ in loc.\ cit., this reduces to investigating various restrictions of $\Fil Y_{\C_p}$. Namely, first note that there are isomorphisms
\begin{equation*}
(\{|u|=1\}\subseteq \Fil Y_{\C_p})\cong Y_{\C_p}\cong (\{|t|=1\}\subseteq \Fil Y_{\C_p})
\end{equation*}
and thus, using \cref{lem:defis-jht}, we conclude that
\begin{equation*}
X^{\N/\C_p}_{|u|=1}\cong X^{\prism/\C_p}\cong X^{\N/\C_p}_{|t|=1}
\end{equation*}
for any Gelfand stack $X$ over $\C_p$. As in the arithmetic case, we denote the corresponding closed embeddings of $X^{\prism/\C_p}$ into $X^{\N/\C_p}$ by
\begin{equation*}
j_\dR: X^{\prism/\C_p}\cong X^{\N/\C_p}_{|t|=1}\rightarrow X^{\N/\C_p}\;, \hspace{0.5cm} j_\HT: X^{\prism/\C_p}\cong X^{\N/\C_p}_{|u|=1}\rightarrow X^{\N/\C_p}\;.
\end{equation*}

Moreover, observe that \cref{lem:defis-rhvariant} yields an isomorphism
\begin{equation*}
\begin{split}
(\ol{\DD}_+^\times\times\ol{\DD}_-^\times)\,/\,\ol{\T}&\cong (\{|x|\leq |u|\}\subseteq \ol{\DD}^\times \times \ol{\DD}_-^\times/\ol{\T})\cong (\{|x|\leq |u|\}\subseteq \ol{\DD}^\times \times (0, 1]) \\
(t, u)&\mapsto (tu, u)\;,
\end{split}
\end{equation*}
where $x=tu$ is the coordinate on $\ol{\DD}^\times$ in the target and hence we obtain
\begin{equation*}
(\{|ut|\neq 0\}\subseteq \Fil Y_{\C_p})\cong (\{0<|\phi^{-1}(\xi)|\leq |u|\}\subseteq Y_{\C_p}\times (0, 1])\;.
\end{equation*}
Finally, we note that
\begin{equation*}
(Y_{\C_p}\setminus \{\phi^{-1}(\xi)=0\})\times [0, 1]\xrightarrow{\cong} \{0<|\phi^{-1}(\xi)|\leq |u|\}\subseteq Y_{\C_p}\times (0, 1])
\end{equation*}
via the map sending $y\in [0, 1]$ to the real number $|u|\in [|\phi^{-1}(\xi)|, 1]$ which is the second coordinate of the unique intersection point $(|t|, |u|)$ of the hyperbola $\{|u|\cdot |t|=|\phi^{-1}(\xi)|\}$ and the line through $(1, 1)$ and $(1-y, y)$. Together with \cref{prop:defis-utneq0}, this induces an isomorphism
\begin{equation*}
X^{\N/\C_p}_{|ut|\neq 0} \cong X^{\prism/\C_p}_{|\widetilde{\mu}|\neq 0} \times [0, 1]
\end{equation*}
with the same properties as in loc.\ cit., where we note that the structure map $\pi: X^{\N/\C_p}\rightarrow X^{\prism/\C_p}$ is induced by the projection $\Fil Y_{\C_p}\rightarrow Y_{\C_p}$ and the structure map $\pi: X^\N\rightarrow X^\prism$ in the arithmetic case.

What is the analogue of the map $i_{\dR, +}: X^{\dR, +}\rightarrow X^\N$? Here we actually encounter a new ``relative stack'', which we shall promptly define. Before, however, we recall the ring $B_\dR^{+, \dagger}$ defined by
\begin{equation*}
B_\dR^{+, \dagger}\coloneqq A_\inf\{\xi\}^\dagger[\tfrac{1}{p}]=\colim_n A_\inf\langle p^{-n}\xi\rangle[\tfrac{1}{p}]
\end{equation*}
and note that it is a (separable qfd) Gelfand ring: Indeed, we have $\xi\in\Nil^\dagger(B_\dR^{+, \dagger})$ and hence the uniform completion of $B_\dR^{+, \dagger}$ is given by $B_\dR^{+, \dagger}/\xi\cong \C_p$. Now observe that the quotient
\begin{equation*}
\GSpec B_\dR^{+, \dagger}\langle t\rangle_{\leq 1}\{u\}^\dagger/(ut-\xi)\;\big/\;\ol{\T}\;,
\end{equation*}
where $\ol{\T}$ acts by multiplication on $t$ and by division on $u$, parametrises $B_\dR^{+, \dagger}$-algebras $A$ together with a pair $(t: L\rightarrow A, u: A\rightarrow L)$ of normed (dual) generalised Cartier divisors with $|t|\leq 1$ and $|u|=0$ and an identification between the composition
\begin{equation*}
A\xrightarrow{u} L\xrightarrow{t} A
\end{equation*}
and multiplication by $\xi$. Noting that this data in particular yields a map
\begin{equation*}
\C_p\cong B_\dR^{+, \dagger}/\xi\rightarrow A/\xi\cong A/tu\rightarrow\Cone(L\tensor_A \Nil^\dagger(A)\xrightarrow{t} A)\;,
\end{equation*}
we obtain a canonical map
\begin{equation*}
\GSpec B_\dR^{+, \dagger}\langle t\rangle_{\leq 1}\{u\}^\dagger/(ut-\xi)\;\big/\;\ol{\T}\rightarrow\C_p^{\dR, +}\;.
\end{equation*}

\begin{warn}
We immediately warn the reader that the coordinate $t$ is \emph{not} to be confused with the element $\log[\epsilon]\in B_\dR^{+, \dagger}$, where $\epsilon=(1, \zeta_p, \zeta_{p^2}, \dots)$, which usually goes by the same name. However, this element will \emph{not} play a role in the sequel and the letter $t$ will always refer to the coordinate on $\ol{\DD}_+$.
\end{warn}

\begin{defi}
Let $X$ be a Gelfand stack over $\C_p$. The \emph{filtered $B_\dR^{+, \dagger}$-stack} $X^{\dR, +/B_\dR^{+, \dagger}}$ of $X$ is defined by the pullback diagram
\begin{equation*}
\begin{tikzcd}
X^{\dR, +/B_\dR^{+, \dagger}}\ar[r]\ar[d] & \GSpec B_\dR^{+, \dagger}\langle t\rangle_{\leq 1}\{u\}^\dagger/(ut-\xi)\;\big/\;\ol{\T}\ar[d] \\
X^{\dR, +}\ar[r] & \C_p^{\dR, +}\nospacepunct{\;.}
\end{tikzcd}
\end{equation*}
Moreover, the \emph{$B_\dR^{+, \dagger}$-stack} $X^{\dR/B_\dR^{+, \dagger}}$ of $X$ is defined as the base change of $X^{\dR, +/B_\dR^{+, \dagger}}$ along the canonical map $\ol{\T}/\ol{\T}\rightarrow\ol{\DD}/\ol{\T}$ while the stack $X^{\Hod/B_\dR^{+, \dagger}}$ is defined as the base change of $X^{\dR, +/B_\dR^{+, \dagger}}$ along $*/\ol{\T}\rightarrow \ol{\DD}/\ol{\T}$.
\end{defi}

From the definition of $\Fil Y_{\C_p}$, we deduce that there is an isomorphism
\begin{equation*}
(\{|u|=0\}\subseteq\Fil Y_{\C_p})\cong \GSpec B_\dR^{+, \dagger}\langle t\rangle_{\leq 1}\{u\}^\dagger/(ut-\xi)\;\big/\;\ol{\T}\;,
\end{equation*}
which we remark involves a Frobenius on $Y_{\C_p}$ and thus we conclude from \cref{prop:defis-u0} that there is an isomorphism
\begin{equation*}
X^{\dR, +/B_\dR^{+, \dagger}}\cong X^{\N/\C_p}_{|u|=0}
\end{equation*}
for any Gelfand stack $X$ over $\C_p$. In particular, as
\begin{equation*}
(\{u=0\}\subseteq \GSpec B_\dR^{+, \dagger}\langle t\rangle_{\leq 1}\{u\}^\dagger/(ut-\xi)\;\big/\;\ol{\T})\cong \ol{\DD}_{\C_p}/\ol{\T}\;,
\end{equation*}
we obtain 
\begin{equation*}
X^{\dR, +/\C_p}\cong X^{\dR, +/B_\dR^{+, \dagger}}_{u=0}\cong X^{\N/\C_p}_{u=0}\;,
\end{equation*}
where the source is defined as the base change of $X^{\dR, +}$ along the canonical map $\ol{\DD}_{\C_p}/\ol{\T}\rightarrow\C_p^{\dR, +}$, see \cite[Def.\ 5.2.2.(2)]{dRStack}.

Let us also define
\begin{equation*}
X^{\HT, \dagger, +/\C_p}\coloneqq X^{\N/\C_p}_{|t|=0}\;, \hspace{0.5cm} X^{\HT, \dagger/\C_p}\coloneqq X^{\N/\C_p}_{|t|=0, |u|=1}
\end{equation*}
and note that
\begin{equation*}
(\{|t|=0, |u|\neq 0\}\subseteq \Fil Y_{\C_p})\cong \GSpec B_\dR^{+, \dagger}\times (0, 1]
\end{equation*}
by \cref{lem:defis-rhvariant}, which together with \cref{prop:defis-t0uneq0} yields an isomorphism
\begin{equation*}
X^{\HT, \dagger, +/\C_p}_{|u|\neq 0}\cong X^{\HT, \dagger/\C_p}\times (0, 1]
\end{equation*}
with the same properties as in loc.\ cit.

Finally, we define the syntomification relative to $\C_p$. This is done in the usual way using the closed embeddings $j_\dR, j_\HT: X^{\prism/\C_p}\rightarrow X^{\N/\C_p}$ from above.

\begin{defi}
Let $X$ be a Gelfand stack over $\C_p$. The \emph{relative analytic syntomification} $X^{\Syn/\C_p}$ of $X$ over $\C_p$ is defined as the coequaliser
\begin{equation*}
\begin{tikzcd}
X^{\prism/\C_p}\ar[r,shift left=.75ex,"j_\HT"]\ar[r,shift right=.75ex,swap,"j_\dR"] & X^{\N/\C_p}\ar[r] & X^{\Syn/\C_p}\;.
\end{tikzcd}
\end{equation*}
Moreover, for any $i\in\Z$, the \emph{relative analytic syntomic cohomology} of weight $i$ of $X$ over $\C_p$ is defined as
\begin{equation*}
R\Gamma_\Syn(X/\C_p, \Q_p(i))\coloneqq R\Gamma(X^{\Syn/\C_p}, \O\{i\})\;,
\end{equation*}
where $\O\{i\}$ is pulled back along the natural map $X^{\Syn/\C_p}\rightarrow X^\Syn$.
\end{defi}

To define Hodge--Tate weights of perfect complexes on $X^{\Syn/\C_p}$, note that the above discussion shows that the locus $\{u=t=0\}$ inside $X^{\N/\C_p}$ identifies with the stack $X^{\Hod/\C_p}$, which sits in a pullback diagram
\begin{equation*}
\begin{tikzcd}
X^{\Hod/\C_p}\ar[r]\ar[d] & \GSpec\C_p/\ol{\T}\ar[d] \\
X^\Hod\ar[r] & \C_p^\Hod
\end{tikzcd}
\end{equation*}
and hence receives a natural map
\begin{equation*}
X\times */\ol{\T}\rightarrow X^{\Hod/\C_p}\;.
\end{equation*}

\begin{defi}
Let $X$ be a Gelfand stack over $\C_p$ and $E$ a perfect complex on $X^{\N/\C_p}$. Identifying the pullback of $E$ along the map
\begin{equation*}
X\times */\ol{\T}\rightarrow X^{\Hod/\C_p}\rightarrow X^{\N/\C_p}
\end{equation*}
with a graded perfect complex $M^\bullet$ on $X$, the \emph{Hodge--Tate weights} of $E$ are defined to be those integers $i\in\Z$ such that $M^i$ is nonzero. The Hodge--Tate weights of a perfect analytic $F$-gauge $E\in\Perf(X^{\Syn/\C_p})$ on $X$ relative to $\C_p$ are defined to be the Hodge--Tate weights of its pullback to the stack $X^{\N/\C_p}$.
\end{defi}

\subsection{$B_\dR^+$-cohomology via stacks}

Before we can move on to discussing the analogues of the comparison theorems from §§\ref{sect:hkcomp}, \ref{sect:proet} in the geometric setting, we first have to take a closer look at the filtered $B_\dR^{+, \dagger}$-stack $X^{\dR, +/B_\dR^{+, \dagger}}$ for a smooth partially proper rigid space $X$ over $\C_p$. Namely, we will prove that cohomology on this stack actually computes the usual $B_\dR^+$-cohomology $R\Gamma(X/B_\dR^+)$ of $X$ equipped with its filtration after base changing along $B_\dR^{+, \dagger}\rightarrow B_\dR^+$:

\begin{thm}
\label{thm:geom-bdr+coh}
Let $X$ be a smooth partially proper rigid space over $\C_p$. Then there is a natural isomorphism
\begin{equation*}
\Fil^\bullet R\Gamma(X^{\dR, +/B_\dR^{+, \dagger}},\O)\tensor_{B_\dR^{+, \dagger}} B_\dR^+\cong \Fil^\bullet R\Gamma(X/B_\dR^+)\;.
\end{equation*}
\end{thm}

Before we begin the proof, let us prepare ourselves by studying the cohomology of $X^{\Hod/B_\dR^{+, \dagger}}$, which is defined as the locus $\{t=0\}\subseteq X^{\dR, +/B_\dR^{+, \dagger}}$ and hence sits inside a pullback square
\begin{equation*}
\begin{tikzcd}
X^{\Hod/B_\dR^{+, \dagger}}\ar[r]\ar[d] & \GSpec\C_p\{u\}^\dagger/\,\ol{\T}\ar[d] \\
X^\Hod\ar[r] & \C_p^\Hod\nospacepunct{\;.}
\end{tikzcd}
\end{equation*}

\begin{lem}
\label{lem:geom-xhod}
Let $X$ be a Berkovich smooth $\dagger$-rigid space over $\C_p$ in the sense of \cite[Def.\ 4.3.6]{dRFF}. Then
\begin{equation*}
X^{\Hod/B_\dR^{+, \dagger}}\cong \left(X\times \GSpec \C_p\{u\}^\dagger/\ol{\T}\right)\big/\,\cal{T}_{X/\C_p}^\dagger\langle -1\rangle\;,
\end{equation*}
where the action of $\cal{T}_{X/\C_p}^\dagger\langle -1\rangle$ is trivial and $\cal{T}_{X/\C_p}^\dagger$ denotes the overconvergent neighbourhood of the zero section inside the restriction of $\AnSpec_X \Sym_X^\bullet \Omega^1_{X/\C_p}$ to the test category of Gelfand rings. In particular, we have
\begin{equation*}
R\Gamma(X^{\Hod/B_\dR^{+, \dagger}}, \O\langle i\rangle)\cong \bigoplus_{j\leq i} R\Gamma(X, \Omega^j_{X/\C_p})[-j]
\end{equation*}
for all $i\in\Z$.
\end{lem}
\begin{proof}
Recall that $X^{\Hod/B_\dR^{+, \dagger}}$ is the sheafification of the assignment sending any triple $(A, u: A\rightarrow L)$ of a Gelfand $\C_p$-algebra $A$ and a dual normed generalised Cartier divisor of norm zero to the anima of maps $\GSpec(A\oplus L\tensor_A \Nil^\dagger(A)[1])\rightarrow X$. By derived deformation theory, the anima of such maps which lift a fixed map $\eta: \GSpec A\rightarrow X$ is a torsor for
\begin{equation*}
\Map_A(\eta^*\Omega^1_{X/\C_p}, \Nil^\dagger(A)\tensor_A L[1])\cong \eta^*(\Omega^1_{X/\C_p})^\vee\tensor_A \Nil^\dagger(A)\tensor_A L[1]
\end{equation*}
and hence we conclude that the natural map $X^{\Hod/B_\dR^{+, \dagger}}\rightarrow X\times\GSpec\C_p\{u\}^\dagger/\ol{\T}$ is a gerbe banded by $\cal{T}_{X/\C_p}^\dagger\langle -1\rangle$. Noting that the gerbe has a splitting induced by the canonical map $A\rightarrow A\oplus L\tensor_A\Nil^\dagger(A)[1]$ proves the first assertion. Now the second assertion follows by Cartier duality, see \cite[Thm.\ 4.3.13]{dRStack} and \cref{lem:recall-cartierperf}, together with \cref{lem:recall-reesgm}.
\end{proof}

We can now move on to the proof of \cref{thm:geom-bdr+coh}.

\begin{proof}[Proof of \cref{thm:geom-bdr+coh}]
We first note that the map $\Cone(L\tensor_A \Nil^\dagger(A)\rightarrow A)\rightarrow \ol{A}$ for any Gelfand $\B_\dR^{+, \dagger}$-algebra $A$ equipped with a normed generalised Cartier divisor $L\rightarrow A$ induces a natural map
\begin{equation*}
f: X^{\dR, +/B_\dR^{+, \dagger}}\rightarrow X^\dR\times\GSpec B_\dR^{+, \dagger}\rightarrow |X|\times\GSpec B_\dR^{+, \dagger}\;.
\end{equation*}
Note that, by compatibility of the assignment $X\mapsto X^{\dR, +/B_\dR^{+, \dagger}}$ with open localisation, which one deduces from \cref{prop:defis-openloc} by base change, the preimage of any open subset of $|X|$ under this map is given by the filtered $B_\dR^{+, \dagger}$-stack of the corresponding open subspace of $X$. By \cite[Prop.\ 2.42]{BoscopAdic}, our task is to check that there is a natural isomorphism
\begin{equation*}
\Fil^\bullet f_*\O\tensor_{B_\dR^{+, \dagger}} B_\dR^+\cong \tau^{\leq 0}_{\mathrm{Beil}} \nu_*\Fil^\bullet \mathbb{B}_\dR^+\;,
\end{equation*}
where the filtration on $\mathbb{B}_\dR^+$ is the $\xi$-adic filtration, $\nu$ is the morphism of sites $\nu: X_\proet\rightarrow |X|$ and the Beilinson $t$-structure on filtered objects is taken with respect to the standard $t$-structure on sheaves on the topological space $|X|$; moreover, we write $\Fil^\bullet f_*\O$ for the filtered object 
\begin{equation*}
\dots\xrightarrow{t} f_*\O\langle -1\rangle\xrightarrow{t} f_*\O\xrightarrow{t} f_*\O\langle 1\rangle\xrightarrow{t}\dots\;.
\end{equation*}

We first construct a natural filtered morphism
\begin{equation}
\label{eq:geom-bdr+map}
\Fil^\bullet f_*\O\tensor_{B_\dR^{+, \dagger}} B_\dR^+\rightarrow \nu_*\Fil^\bullet\mathbb{B}_\dR^+\;,
\end{equation}
which in turn comes down to constructing natural filtered morphisms
\begin{equation*}
\Fil^\bullet R\Gamma(U^{\dR, +/B_\dR^{+, \dagger}}, \O)\tensor_{B_\dR^{+, \dagger}} B_\dR^+\rightarrow R\Gamma_\proet (U, \Fil^\bullet\mathbb{B}_\dR^+)\;,
\end{equation*}
for any open subspace $U\subseteq X$. This works as follows: For any affinoid perfectoid $C$ over $U$, we have a map
\begin{equation*}
\GSpec \mathbb{B}_\dR^{+, \dagger}(C)\langle t\rangle_{\leq 1}\{u\}^\dagger/(ut-\xi)\;/\;\ol{\T}\rightarrow C^{\dR, +}
\end{equation*}
defined in the same way as for $C=\C_p$ above, and together with the natural map $B_\dR^{+, \dagger}\rightarrow\mathbb{B}_\dR^{+, \dagger}(C)$, this induces a map
\begin{equation}
\label{eq:geom-bcmapforbdr+stack}
\GSpec \mathbb{B}_\dR^{+, \dagger}(C)\langle t\rangle_{\leq 1}\{u\}^\dagger/(ut-\xi)\;/\;\ol{\T}\rightarrow C^{\dR, +/B_\dR^{+, \dagger}}\;.
\end{equation}
Now observe that the filtered object calculated by the cohomology of the source of (\ref{eq:geom-bcmapforbdr+stack}) is given by $\Fil^\bullet \mathbb{B}_\dR^{+, \dagger}(C)$ and thus, by functoriality of the filtered $B_\dR^{+, \dagger}$-stack, we obtain a filtered morphism
\begin{equation*}
\begin{split}
\Fil^\bullet R\Gamma(U^{\dR, +/B_\dR^{+, \dagger}}, \O)\tensor_{B_\dR^{+, \dagger}} B_\dR^+&\rightarrow \lim_{\GSpec C\rightarrow U} \Fil^\bullet R\Gamma(C^{\dR, +/B_\dR^{+, \dagger}}, \O) \\
&\rightarrow\lim_{\GSpec C\rightarrow U} \Fil^\bullet \mathbb{B}_\dR^+(C)\cong R\Gamma_\proet(U, \Fil^\bullet \mathbb{B}_\dR^+)\;,
\end{split}
\end{equation*}
where the second arrow is given by pullback of sections along (\ref{eq:geom-bcmapforbdr+stack}).

Next, observe that $\Fil^\bullet f_*\O\tensor_{B_\dR^{+, \dagger}} B_\dR^+$ is connective for the Beilinson $t$-structure since
\begin{equation}
\label{eq:geom-bdr+gradedpieces}
\gr^i f_*\O\cong \bigoplus_{j\leq i} \Omega_{X/\C_p}^j[-j]
\end{equation}
by applying \cref{lem:geom-xhod} to open subspaces $U\subseteq X$ and thus (\ref{eq:geom-bdr+map}) naturally factors through a map
\begin{equation}
\label{eq:geom-bdr+maptrunc}
\Fil^\bullet f_*\O\tensor_{B_\dR^{+, \dagger}} B_\dR^+\rightarrow\tau^{\leq 0}_{\mathrm{Beil}} \nu_*\Fil^\bullet\mathbb{B}_\dR^+\;.
\end{equation}
Putting this together with
\begin{equation*}
\nu_*t^i\mathbb{B}_\dR^+/t^{i+1}\mathbb{B}_\dR^+\cong \bigoplus_{n\geq 0} \Omega^n_X[-n]\;,
\end{equation*}
see \cite[Prop.\ 3.23]{ScholzeSurvey}, we conclude that (\ref{eq:geom-bdr+maptrunc}) is an isomorphism on graded pieces. As all filtration steps of $\nu_*\Fil^\bullet\mathbb{B}_\dR^+$ are cohomologically bounded above by \cite[Prop.\ 2.40]{BoscopAdic}, the natural map
\begin{equation*}
\tau_{\mathrm{Beil}}^{\leq 0} \nu_*\Fil^\bullet\mathbb{B}_\dR^+\rightarrow \nu_*\Fil^\bullet\mathbb{B}_\dR^+
\end{equation*}
is an isomorphism in high enough filtration degrees. Since the filtration $\Fil^\bullet\mathbb{B}_\dR^+$ is complete, we conclude that the target of (\ref{eq:geom-bdr+maptrunc}) is complete for the filtration as well. As we have just seen that (\ref{eq:geom-bdr+maptrunc}) induces isomorphisms on graded pieces, we are done once we know that the source is complete for the filtration as well.

For this we may localise on $|X|$ and then assume that $X$ admits a smooth lift $\widetilde{X}$ to $B_\dR^{+, \dagger}$ which comes equipped with a Berkovich étale map $\widetilde{X}\rightarrow \ol{\DD}^n\times\GSpec B_\dR^{+, \dagger}$ for some $n\geq 0$. Since Berkovich étale maps are $\dagger$-formally étale, this yields a cartesian diagram
\begin{equation*}
\begin{tikzcd}
\widetilde{X}\times_{\GSpec B_\dR^{+, \dagger}} \GSpec B_\dR^{+, \dagger}\langle t\rangle_{\leq 1}\{u\}^\dagger/(ut-\xi)\;/\;\ol{\T}\ar[r]\ar[d] & X^{\dR, +/B_\dR^{+, \dagger}}\ar[d] \\
\ol{\DD}^n\times\GSpec B_\dR^{+, \dagger}\langle t\rangle_{\leq 1}\{u\}^\dagger/(ut-\xi)\;/\;\ol{\T}\ar[r] & (\ol{\DD}^n_{\C_p})^{\dR, +/B_\dR^{+, \dagger}}\nospacepunct{\;.}
\end{tikzcd}
\end{equation*}
Now the bottom map exhibits the target as the quotient of the source by the action of $\G_a^\dagger\langle -1\rangle^n$ by addition and thus we conclude
\begin{equation*}
X^{\dR, +/B_\dR^{+, \dagger}}\cong \left(\widetilde{X}\times_{\GSpec B_\dR^{+, \dagger}} \GSpec B_\dR^{+, \dagger}\langle t\rangle_{\leq 1}\{u\}^\dagger/(ut-\xi)\;/\;\ol{\T}\right)\big/\, \G_a^\dagger\langle -1\rangle^n\;.
\end{equation*}

By Cartier duality, see \cref{lem:recall-cartierperf}, and using that the addition action of $\G_a^\dagger$ on $\ol{\DD}$ corresponds to the endomorphism of $\O_{\ol{\DD}}$ given by the derivative, we conclude from the above presentation that
\begin{equation*}
\Fil^\bullet f_*\O\cong \Fil^\bullet_\Hod (\O_{\widetilde{X}}\xrightarrow{\mathrm{d}} \Omega_{\widetilde{X}/B_\dR^{+, \dagger}}^1\xrightarrow{\mathrm{d}} \Omega_{\widetilde{X}/B_\dR^{+, \dagger}}^2\xrightarrow{\mathrm{d}}\dots)\tensor_{B_\dR^{+, \dagger}} \Fil^\bullet B_\dR^{+, \dagger}
\end{equation*}
and tensoring with $B_\dR^+$ yields
\begin{equation*}
\Fil^\bullet f_*\O\tensor_{B_\dR^{+, \dagger}} B_\dR^+\cong \Fil^\bullet_\Hod (\O_{\widetilde{X}}\xrightarrow{\mathrm{d}} \Omega_{\widetilde{X}/B_\dR^{+, \dagger}}^1\xrightarrow{\mathrm{d}} \Omega_{\widetilde{X}/B_\dR^{+, \dagger}}^2\xrightarrow{\mathrm{d}}\dots)\tensor_{B_\dR^{+, \dagger}} \Fil^\bullet B_\dR^+\;.
\end{equation*}
Finally, as the Hodge filtration vanishes in high enough degrees, completeness of the filtration on the left-hand side now follows from completeness of the $\xi$-adic filtration on $B_\dR^+$. This finishes the proof.
\end{proof}

The fact that the stack $X^{\dR, +/B_\dR^{+, \dagger}}$ computes filtered $B_\dR^+$-cohomology of smooth partially proper rigid spaces over $\C_p$ by \cref{thm:geom-bdr+coh} now lets us easily deduce some well-known properties of $B_\dR^+$-cohomology.

\begin{cor}
\label{cor:geom-bdr+props}
Let $X$ be a smooth partially proper rigid space over $\C_p$. Then there is a natural isomorphism
\begin{equation*}
R\Gamma(X/B_\dR^+)\tensor_{B_\dR^+} \C_p\cong R\Gamma_\dR(X/\C_p)
\end{equation*}
compatible with the filtrations. Moreover, if $X$ is qcqs and admits a model $\ol{X}$ over some finite extension $K$ of $\Q_p$, then
\begin{equation*}
R\Gamma(X/B_\dR^+)\cong R\Gamma_\dR(\ol{X})\tensor_K B_\dR^+\;,
\end{equation*}
again compatibly with the filtrations.
\end{cor}
\begin{proof}
For the first isomorphism, observe that the big rectangle and the lower square in the diagram
\begin{equation}
\label{eq:geom-bdr+cohdrcompsquare1}
\begin{tikzcd}
X^{\dR, +/\C_p}\ar[r]\ar[d] & \GSpec \C_p\times \ol{\DD}/\ol{\T}\ar[d] \\
X^{\dR, +/B_\dR^{+, \dagger}}\ar[r]\ar[d] & \GSpec B_\dR^{+, \dagger}\langle t\rangle_{\leq 1}\{u\}^\dagger/(ut-\xi)\;/\;\ol{\T}\ar[d] \\
X^{\dR, +}\ar[r] & \C_p^{\dR, +}
\end{tikzcd}
\end{equation}
are cartesian, hence the upper square is cartesian as well. Moreover, the map 
\begin{equation*}
\GSpec\C_p\times \ol{\DD}/\ol{\T}\rightarrow \GSpec B_\dR^{+, \dagger}\langle t\rangle_{\leq 1}\{u\}^\dagger/(ut-\xi)\;/\;\ol{\T}
\end{equation*}
is cohomologically smooth since this property is local on the target and the map
\begin{equation*}
\GSpec\C_p\langle t\rangle_{\leq 1}\rightarrow \GSpec B_\dR^{+, \dagger}\langle t\rangle_{\leq 1}\{u\}^\dagger/(ut-\xi)
\end{equation*}
is a complete intersection with kernel generated by $u$. Thus, the top square in the diagram (\ref{eq:geom-bdr+cohdrcompsquare1}) satisfies base change by \cite[Lem.\ 4.5.13.(i)]{HeyerMann} and this implies the claim upon recalling that $X^{\dR, +/\C_p}$ calculates Hodge-filtered de Rham cohomology of $X$ relative to $\C_p$, see \cite[Rem.\ 5.2.2]{dRFF}.

For the second claim, we argue similarly: Both the bottom and the top square in the diagram
\begin{equation*}
\begin{tikzcd}
X^{\dR, +/B_\dR^{+, \dagger}}\ar[r]\ar[d] & \GSpec B_\dR^{+, \dagger}\langle t\rangle_{\leq 1}\{u\}^\dagger/(ut-\xi)\;/\;\ol{\T}\ar[d] \\
X^{\dR, +}\ar[r]\ar[d] & \C_p^{\dR, +}\ar[d] \\
\ol{X}^{\dR, +}\ar[r] & K^{\dR, +}
\end{tikzcd}
\end{equation*}
are cartesian and hence the same is true for the large rectangle. By \cref{lem:geom-dr+prim} below, the map $\ol{X}^{\dR, +}\rightarrow K^{\dR, +}\cong \GSpec K\times \ol{\DD}/\ol{\T}$ is prim and thus the large rectangle satisfies base change by \cite[Lem.\ 4.5.13.(ii)]{HeyerMann}, which implies the claim.
\end{proof}

The following lemma was used in the proof:

\begin{lem}
\label{lem:geom-dr+prim}
Let $X$ be a smooth partially proper qcqs rigid space over a finite extension $K$ of $\Q_p$. Then the map $X^{\dR, +}\rightarrow K^{\dR, +}\cong \GSpec K\times\ol{\DD}/\ol{\T}$ is prim.
\end{lem}
\begin{proof}
As in the proof of \cref{lem:geom-hkprim}, we may assume that $X$ admits a Berkovich étale map $X\rightarrow\ol{\DD}_K^n$ after passing to a finite strict closed cover of $X$ using compatibility of $X\mapsto X^{\dR, +}$ with rational localisations, which one deduces from \cref{prop:defis-openloc}. After possibly passing to a further finite strict closed cover, we can assume that the map $X\rightarrow\ol{\DD}_K^n$ factors as $X\rightarrow Z\rightarrow\ol{\DD}_K^n$, where $X\rightarrow Z$ is finite étale and $Z\rightarrow\ol{\DD}_K^n$ is a rational localisation. By \cref{prop:defis-openloc}, the map $Z^{\dR, +}\rightarrow(\ol{\DD}^n_K)^{\dR, +}$ is a rational localisation as well and hence prim. Moreover, as $X\rightarrow Z$ is finite étale and hence in particular $\dagger$-formally étale, the diagram
\begin{equation*}
\begin{tikzcd}
X\times \ol{\DD}/\ol{\T}\ar[r]\ar[d] & X^{\dR, +}\ar[d] \\
Z\times \ol{\DD}/\ol{\T}\ar[r] & Z^{\dR, +}
\end{tikzcd}
\end{equation*}
is cartesian and since $Z$ is $\dagger$-formally smooth (as it is a rational subspace of $\ol{\DD}_K^n$), the bottom horizontal map is a $!$-cover. Hence, $X^{\dR, +}\rightarrow Z^{\dR, +}$ is prim if the same is true for $X\rightarrow Z$, but the latter is a finite étale map, thus in particular affine, whence prim.

Thus, we are done once we show that $(\ol{\DD}_K^n)^{\dR, +}\rightarrow K^{\dR, +}$ is prim. However, we have
\begin{equation*}
(\ol{\DD}_K^n)^{\dR, +}\cong (\ol{\DD}_K^n\times \ol{\DD}/\ol{\T})\;/\;\G_a^\dagger\langle -1\rangle^n\;,
\end{equation*}
where $\G_a^\dagger\langle -1\rangle^n$ acts on $\ol{\DD}_K^n$ by translation and thus the claim follows from the fact that both $\ol{\DD}^n_K\rightarrow\GSpec K$ and $*/\G_a^\dagger\langle -1\rangle^n\rightarrow *$ are prim: Indeed, the first map is affine and for the second one this is \cite[Thm.\ 4.3.13]{dRStack}.
\end{proof}
\comment{
We used the following lemma:

\begin{lem}
\label{lem:geom-descentbdr+stack}
Let $X\rightarrow Y$ be a rational localisation of derived Berkovich spaces over $\Q_p$. Then the induced map $X^{\dR, +}\rightarrow Y^{\dR, +}$ is a rational localisation as well. Moreover, if the collection $\{X_i\rightarrow Y\}_{i\in I}$ of rational localisations forms a strict closed cover of $Y$, then the maps $\{X_i^{\dR, +}\rightarrow Y^{\dR, +}\}_{i\in I}$ jointly cover $Y^{\dR, +}$.
\end{lem}
\begin{proof}
As in the proof of \cref{prop:defis-openloc}, we may assume that we are given an $A$-valued point of $Y^{\dR, +}$ corresponding to a map
\begin{equation*}
\GSpec \Cone(L\tensor_A \Nil^\dagger(A)\xrightarrow{t} A)\rightarrow Y
\end{equation*}
for some normed generalised Cartier divisor $t: L\rightarrow A$ of norm at most $1$. In particular, we obtain a composite map 
\begin{equation*}
\GSpec\ol{A}\rightarrow \GSpec \Cone(L\tensor_A \Nil^\dagger(A)\xrightarrow{t} A)\rightarrow Y
\end{equation*}
and we let $B$ be the rational localisation of $A$ defined by the closed subset $\cal{M}(\ol{A})\times_{|Y|} |X|$ of $\cal{M}(\ol{A})\cong \cal{M}(A)$. We claim that there is a cartesian diagram
\begin{equation}
\label{eq:geom-bdr+ratloc}
\begin{tikzcd}
\GSpec B\ar[r]\ar[d] & X^{\dR, +}\ar[d] \\
\GSpec A\ar[r] & Y^{\dR, +}\nospacepunct{\;.}
\end{tikzcd}
\end{equation}
Indeed, note that
\begin{equation*}
\begin{split}
\GSpec\Cone(L\tensor_A \Nil^\dagger(A)\rightarrow A)\times_{\GSpec A} \GSpec B&\cong \GSpec\Cone(L\tensor_A \Nil^\dagger(A)\rightarrow A)\times_{\cal{M}(A)} \cal{M}(B) \\
&\cong \GSpec\Cone(L\tensor_A \Nil^\dagger(A)\rightarrow A)\times_Y X
\end{split}
\end{equation*}
by definition and hence we obtain a composite map
\begin{equation*}
\GSpec \Cone(L\tensor_A\Nil^\dagger(B)\rightarrow B)\rightarrow \GSpec\Cone(L\tensor_A \Nil^\dagger(A)\rightarrow A)\times_{\GSpec A} \GSpec B\rightarrow X\;,
\end{equation*}
which yields a $B$-point of $X^{\dR, +}$ making the diagram (\ref{eq:geom-bdr+ratloc}) commute.

To check that the diagram is cartesian, take an $A$-algebra $C$ together with a $C$-point of $X^{\dR, +}$ arising from a map $\GSpec\Cone(L\tensor_A\Nil^\dagger(C)\rightarrow C)\rightarrow X$ which lifts the given $A$-point of $Y^{\dR, +}$. Then the commutative diagram
\begin{equation*}
\begin{tikzcd}
\GSpec\ol{C}\ar[r]\ar[d] & \GSpec\Cone(L\tensor_A\Nil^\dagger(C)\rightarrow C)\ar[r]\ar[d] & X\ar[d] \\
\GSpec\ol{A}\ar[r] & \GSpec\Cone(L\tensor_A\Nil^\dagger(A)\rightarrow A)\ar[r] & Y
\end{tikzcd}
\end{equation*}
shows that the map $\cal{M}(\ol{C})\rightarrow\cal{M}(\ol{A})$ factors through $\cal{M}(B)$ and, consequently, the $A$-algebra structure on $C$ uniquely lifts to a $B$-algebra structure. Moreover, the given $C$-point of $X^{\dR, +}$ must arise from the $B$-point constructed above by compatibility with the given $A$-point of $Y^{\dR, +}$: Indeed, this is because $X\rightarrow Y$ is a rational localisation and hence derived Berkovich spaces over $X$ embed fully faithfully into derived Berkovich spaces over $Y$. Overall, we conclude that (\ref{eq:geom-bdr+ratloc}) is cartesian, as desired. Finally, the last assertion about preservation of covers is deduced as at the end of the proof of \cref{prop:defis-openloc}.
\end{proof}
}

\subsection{Comparison theorems over $\C_p$}

Having introduced all the relevant stacks in the geometric situation, i.e.\ over $\C_p$, and clarified the relation between the filtered $B_\dR^{+, \dagger}$-stack and classical $B_\dR^+$-cohomology, we are now in a position to state and prove the analogues of the main comparison theorems from §§\ref{sect:padicht}, \ref{sect:hkcomp} and \ref{sect:proet}. We start with the geometric form of the de Rham comparison theorem from \cite[Thm.\ 13.1]{IntegralpAdicHT}.

\begin{thm}
\label{thm:geom-drcomp}
Let $X$ be a smooth proper rigid space over $\C_p$. Then there is a natural isomorphism
\begin{equation*}
R\Gamma_\proet(X, \Q_p)\tensor_{\Q_p} B_\dR\cong R\Gamma(X/B_\dR^+)\tensor_{B_\dR^+} B_\dR
\end{equation*}
compatible with the filtrations.
\end{thm}
\begin{proof}
By \cref{cor:sixfunctors-relativepoincare}, the induced map $f^{\Syn/\C_p}: X^{\Syn/\C_p}\rightarrow \C_p^{\Syn/\C_p}$ is cohomologically smooth and weakly cohomologically proper and thus the pushforward $E\coloneqq f^{\Syn/\C_p}_*\O$ is a perfect complex on $\C_p^{\Syn/\C_p}$ by \cite[Lem.\ 4.5.16]{HeyerMann}.

By proper base change and \cref{thm:geom-bdr+coh}, we know that, under the Rees equivalence from \cref{prop:recall-rees}, the pullback of $E$ to
\begin{equation*}
\C_{p, |u|=0}^{\N/\C_p}=(\{|u|=0\}\subseteq\Fil Y_{\C_p})\cong \GSpec B_\dR^{+, \dagger}\langle t\rangle_{\leq 1}\{u\}^\dagger/(ut-\xi)\;\big/\;\ol{\T}
\end{equation*}
identifies with a descendingly filtered perfect complex $\Fil^\bullet V$ over $B_\dR^{+, \dagger}$ equipped with backwards maps $u: \Fil^i V\rightarrow\Fil^{i+1} V$ such that $ut=\xi$ whose base change to $B_\dR^+$ is given by $\Fil^\bullet R\Gamma(X/B_\dR^+)$. Moreover, by \cite{AnPrism} and using proper base change again, we also know that $E|_{Y_{\C_p}}$ identifies with the perfect complex on $Y_{\C_p}$ obtained from $R\Gamma_\proet(X, \Q_p)$ by tensoring with $\O_{Y_{\C_p}}$. In particular, we see that the pullback of $E$ to $\{|u|=1, |t|=0\}$ is given by $R\Gamma_\proet(X, \Q_p)\tensor_{\Q_p} B_\dR^{+, \dagger}$. As the pullback of $E$ to
\begin{equation*}
(\{|ut|=0\}\subseteq \Fil Y_{\C_p})\cong \GSpec B_\dR^{+, \dagger}\langle u, t\rangle_{\leq 1}/(ut-\xi)\;\big/\;\ol{\T}
\end{equation*}
is perfect, we conclude that
\begin{equation*}
R\Gamma_\proet(X, \Q_p)\tensor_{\Q_p} B_\dR^+\cong \sum_i \Fil^i R\Gamma(X/B_\dR^+)\tensor_{B_\dR^+} \xi^{-i}B_\dR\;,
\end{equation*}
from which the claim follows by inverting $\xi$.
\end{proof}

Let us now move on to adapting the results from §\ref{sect:hkcomp}. For this, we make the following definition:

\begin{defi}
Let $X$ be a Gelfand stack over $\C_p$. The \emph{relative mock Hyodo--Kato stack} $X^{\mHK/\C_p}$ of $X$ over $\C_p$ is defined as the base change
\begin{equation*}
\begin{tikzcd}
X^{\mHK/\C_p}\ar[r]\ar[d] & X^{\Syn/\C_p}\ar[d] \\
X^\mHK\ar[r, "i_\HK"] & X^\Syn\nospacepunct{\;.}
\end{tikzcd}
\end{equation*}
\end{defi}

As in \cref{prop:hkcomp-hkpshk}, one proves that there is an equivalence of categories
\begin{equation*}
\Perf(X^{\mHK/\C_p})\cong \Perf(X^{\HK/\C_p})
\end{equation*}
for any Gelfand stack over $\C_p$. Moreover, we obtain a commutative diagram
\begin{equation*}
\begin{tikzcd}
X^{\dR/B_\dR^{+, \dagger}}\ar[r]\ar[d] & X^{\dR, +/B_\dR^{+, \dagger}}\ar[d] \\
X^{\mHK/\C_p}\ar[r] & X^{\Syn/\C_p}\;,
\end{tikzcd}
\end{equation*}
which hence induces realisation functors
\begin{equation*}
T_\HK: \Perf(X^{\Syn/\C_p})\rightarrow \Perf(X^{\HK/\C_p})
\end{equation*}
from perfect analytic $F$-gauges on $X$ relative to $\C_p$ to coefficients for Hyodo--Kato cohomology of $X$ relative to $\C_p$ and
\begin{equation*}
T_{\dR, +}: \D(X^{\Syn/\C_p})\rightarrow \D(X^{\dR, +/B_\dR^{+, \dagger}})
\end{equation*}
from analytic $F$-gauges on $X$ relative to $\C_p$ to coefficients for filtered $B_\dR^+$-cohomology of $X$. Postcomposing the latter functor with pullback to $X^{\dR/B_\dR^{+, \dagger}}$, we obtain a functor
\begin{equation*}
T_\dR: \D(X^{\Syn/\C_p})\rightarrow \D(X^{\dR/B_\dR^{+, \dagger}})\;.
\end{equation*}
Then the analogue of \cref{thm:hkcomp-main} in the geometric setting, i.e.\ for Berkovich smooth derived Berkovich spaces over $\C_p$, is the following:

\begin{thm}
\label{thm:geom-hkcomp}
Let $X$ be a Berkovich smooth derived Berkovich space over $\C_p$. Then there is an equivalence of categories
\begin{equation*}
\Perf(X^{\Syn/\C_p})\cong \Perf(X^{\HK/\C_p})\times_{\Perf(X^{\dR/B_\dR^{+, \dagger}})} \Perf(X^{\dR, +/B_\dR^{+, \dagger}})
\end{equation*}
induced by the realisation functors $T_\HK$ and $T_{\dR, +}$. In particular, for any $E\in\Perf(X^{\Syn/\C_p})$, there is a pullback diagram
\begin{equation*}
\begin{tikzcd}
R\Gamma(X^{\Syn/\C_p}, E)\ar[r]\ar[d] & R\Gamma(X^{\HK/\C_p}, T_\HK(E))\ar[d] \\
R\Gamma(X^{\dR, +/B_\dR^{+, \dagger}}, T_{\dR, +}(E))\ar[r] & R\Gamma(X^{\dR/B_\dR^{+, \dagger}}, T_\dR(E))\nospacepunct{\;.}
\end{tikzcd}
\end{equation*}
\end{thm}
\begin{proof}
Carries over almost verbatim from \cref{thm:hkcomp-main} with the use of \cref{prop:perf-coversmoothrigid} replaced by \cref{lem:geom-coversmooth} below.
\end{proof}

The following lemma was used in the proof:

\begin{lem}
\label{lem:geom-coversmooth}
Let $X$ be a derived Berkovich space over $\C_p$ admitting a map $X\rightarrow \ol{\T}_{\C_p}^n$ which is the composition of a finite étale map $X\rightarrow Z$ and a rational localisation $Z\rightarrow\ol{\T}_{\C_p}^n$. Then the Gelfand stack $X_{[r, s]}^{\N/\C_p}$ is nicely coverable for any $[r, s]\subseteq (p^{1/2}, p^{3/2})$.
\end{lem}
\begin{proof}
It suffices to do the case $X=\GSpec\C_p$, then the general case follows by the same argument as in \cref{prop:perf-coversmoothrigid}. However, we have
\begin{equation*}
\C_{p, [r, s]}^{\N/\C_p}=(\Fil Y_{\C_p})_{[r, s]}=(\{ut=\phi^{-1}(\xi)\}\subseteq Y_{\C_p, [r, s]}\times \ol{\DD}_+\times\ol{\DD}_-)/\ol{\T}
\end{equation*}
and this is covered by the affine truncated Gelfand stack $(\{ut=\phi^{-1}(\xi)\}\subseteq Y_{\C_p, [r, s]}\times \ol{\DD}_+\times\ol{\DD}_-)$ with the \v{C}ech nerve being
\begin{equation*}
(\{ut=\phi^{-1}(\xi)\}\subseteq Y_{\C_p, [r, s]}\times \ol{\DD}_+\times\ol{\DD}_-)\times \ol{\T}^\bullet\;.
\end{equation*}
This proves the claim as $\ol{\T}$ is affine, static and flat over $\Q_p$ by \cref{lem:perf-vspflat}.
\end{proof}

As in the arithmetic case, we want to make the conclusion of \cref{thm:geom-hkcomp} explicit in the case $E=\O\{i\}$ and expect that this compares to a more classical comparison theorem in the $p$-adic Hodge theory of rigid spaces over $\C_p$. While the first guess one might have for what we obtain might be \cite[Thm.\ 1.3.(3)]{BasicComparison}, this is actually not precisely what we get: This is due to the appearance of the period ring $B_\st^+$ in loc.\ cit.\ while we will rather see $B_{\log}$ appear. Hence, our result should more appropriately be compared with \cite[Thm.\ 7.1]{BoscopAdic}.

\begin{cor}
\label{cor:geom-hkcompclassical}
Let $X$ be a smooth partially proper qcqs rigid space over $\C_p$. For any $i\in\Z$, there is a cartesian diagram
\begin{equation*}
\begin{tikzcd}
R\Gamma_\Syn(X/\C_p, \Q_p(i))\ar[r]\ar[d] & (R\Gamma_\HK(X)\tensor_{\Q_p^\un} B_{\log})^{\phi=p^i, N=0}\ar[d] \\
\Fil^i R\Gamma(X/B_\dR^+)\ar[r] & R\Gamma(X/B_\dR^+)\nospacepunct{\;,}
\end{tikzcd}
\end{equation*}
where $B_{\log}$ is the period ring from \cite[Def.\ 10.3.1]{FarguesFontaine}.
\end{cor}
\begin{proof}
Let us first say something about the bottom row. For this, note that the diagram
\begin{equation*}
\begin{tikzcd}
R\Gamma(X^{\dR, +/B_\dR^{+, \dagger}}, \O\langle -i\rangle)\ar[r]\ar[d] & R\Gamma(X^{\dR, B_\dR^{+, \dagger}}, \O)\ar[d] \\
\Fil^i R\Gamma(X/B_\dR^+)\ar[r] & R\Gamma(X/B_\dR^+)
\end{tikzcd}
\end{equation*}
obtained from \cref{thm:geom-bdr+coh} is cartesian: Indeed, this amounts to checking that the induced map between the cofibres of the horizontal maps is an isomorphism, which in turn reduces to proving that
\begin{equation*}
R\Gamma(X^{\Hod/B_\dR^{+, \dagger}}, \O\langle -i\rangle)\cong \gr^i R\Gamma(X/B_\dR^+)
\end{equation*}
by induction. However, this follows from \cref{lem:geom-xhod} and \cite[Prop.\ 3.13]{BasicComparison}.

Thus, we are done once we show that
\begin{equation*}
R\Gamma(X^{\HK/\C_p}, \O\{i\})\cong (R\Gamma_\HK(X)\tensor_{\Q_p^\un} B_{\log})^{N=0, \phi=p^i}\;.
\end{equation*}
For this, recall that $X^{\HK/\C_p}$ by definition sits in a pullback diagram
\begin{equation*}
\begin{tikzcd}
X^{\HK/\C_p}\ar[r, "f'"]\ar[d, "g'", swap] & \FF_{\C_p}\ar[d, "g"] \\
X^\HK\ar[r, "f"] & \C_p^\HK
\end{tikzcd}
\end{equation*}
and note that $X^\HK\rightarrow\C_p^\HK$ is prim by the same proof as in \cref{lem:geom-hkprim}. Thus, the above diagram satisfies base change and we obtain 
\begin{equation*}
f'_*\O\{i\}=f'_*{g'}^*\O\{i\}\cong g^*f_*\O\{i\}\;.
\end{equation*}
As $f_*\O\{i\}$ identifies with a twist of the $(\phi, N)$-module $R\Gamma_\HK(X)$ over $\Q_p^\un$ via the equivalence on perfect complexes induced by pullback along the map
\begin{equation*}
\Psi: \C_p^\HK\rightarrow \ol{\F}_p^\HK/\,\G_a\{-1\}\cong \GSpec\Q_p^\un/\,\G_a\rtimes\phi^\Z
\end{equation*}
from \cite[Thm.\ 7.1.1]{dRFF}, understanding $g^*f_*\O(i)$, amounts to understanding pullback along the map
\begin{equation*}
\Psi\circ g: \FF_{\C_p}\rightarrow \GSpec\Q_p^\un/\,\G_a\rtimes\phi^\Z\;.
\end{equation*}
However, by the construction in loc.\ cit., the $\G_a\rtimes\phi^\Z$-bundle over $\FF_{\C_p}$ this map classifies is just the pullback of the Fargues--Fontaine surface from \cite[Def.\ 10.3.7]{FarguesFontaine} along $Y_{\C_p}\rightarrow \FF_{\C_p}$ and hence the claim follows from the fact that the global sections of this bundle are given by $B_{\log}$.
\end{proof}

\begin{rem}
\label{rem:geom-syntomicffcoh}
The above result shows that $R\Gamma_\Syn(X/\C_p, \Q_p(i))$ is \emph{not} the same as the geometric syntomic cohomology of Colmez--Nizio{\l} from \cite{BasicComparison}. Indeed, this is because, for their geometric syntomic cohomology $R\Gamma_{\Syn, \mathrm{CN}}(X/\C_p, \Q_p(i))$, we have a pullback diagram
\begin{equation*}
\begin{tikzcd}
R\Gamma_{\Syn, \mathrm{CN}}(X/\C_p, \Q_p(i))\ar[r]\ar[d] & (R\Gamma_\HK(X)\tensor_{\Q_p^\un} B_{\mathrm{st}}^+)^{\phi=p^i, N=0}\ar[d] \\
\Fil^i R\Gamma(X/B_\dR^+)\ar[r] & R\Gamma(X/B_\dR^+)\nospacepunct{\;,}
\end{tikzcd}
\end{equation*}
see \cite[Thm.\ 1.3]{BasicComparison}, and replacing $B_{\mathrm{st}}^+$ by $B_{\log}$ here \emph{does} in general change the $(\phi=p^i, N=0)$-eigenspace. Indeed, the result above rather shows that our geometric syntomic cohomology $R\Gamma_\Syn(X/\C_p, \Q_p(i))$ recovers Bosco's \emph{syntomic Fargues--Fontaine cohomology} from \cite[§6]{BoscopAdic} by \cite[Thm.s 4.1, 6.3]{BoscopAdic} and, by \cite[Ex.\ 6.36]{BoscopAdic}, this already differs from Colmez--Nizio{\l} geometric syntomic cohomology in the case of the projective line.
\end{rem}

Finally, let us discuss the analogue of the comparison with proétale cohomology from §\ref{sect:proet}.

\begin{defi}
Let $X$ be a Gelfand stack over $\C_p$. The Gelfand stack $X^{\mDiv/\C_p}$ is defined as the base change
\begin{equation*}
\begin{tikzcd}
X^{\mDiv/\C_p}\ar[r]\ar[d] & X^{\Syn/\C_p}\ar[d] \\
X^\mDiv\ar[r, "i_{\Div^1}"] & X^\Syn\nospacepunct{\;.}
\end{tikzcd}
\end{equation*}
\end{defi}

Again, there is an equivalence of categories
\begin{equation*}
\Perf(X^{\mDiv/\C_p})\cong \Perf(X^{\Div/\C_p})
\end{equation*}
for any Gelfand stack over $\C_p$ and we obtain a commutative diagram
\begin{equation*}
\begin{tikzcd}
X^{\HT, \dagger/\C_p}\ar[r]\ar[d] & X^{\HT, \dagger, +/\C_p}\ar[d] \\
X^{\mDiv/\C_p}\ar[r] & X^{\Syn/\C_p}\;,
\end{tikzcd}
\end{equation*}
which hence induces realisation functors
\begin{equation*}
T_{\Div^1}: \Perf(X^{\Syn/\C_p})\rightarrow \Perf(X^{\Div^1/\C_p})
\end{equation*}
from perfect analytic $F$-gauges on $X$ relative to $\C_p$ to perfect complexes on $X^{\Div^1/\C_p}$, which according to \cite{AnPrism} should be thought of as coefficients for proétale cohomology of $X$, and
\begin{equation*}
T_{\HT, \dagger, +}: \D(X^{\Syn/\C_p})\rightarrow \D(X^{\HT, \dagger, +/\C_p})
\end{equation*}
from analytic $F$-gauges on $X$ relative to $\C_p$ to quasicoherent sheaves on $X^{\HT, \dagger, +/\C_p}$. Postcomposing the latter functor with pullback to $X^{\HT, \dagger/\C_p}$, we also obtain a functor
\begin{equation*}
T_{\HT, \dagger}: \D(X^{\Syn/\C_p})\rightarrow \D(X^{\HT, \dagger/\C_p})\;.
\end{equation*}
This lets us formulate the analogue of \cref{thm:proet-main} in the geometric setting.

\begin{thm}
\label{thm:geom-proet1}
Let $X$ be a Berkovich smooth derived Berkovich space over $\C_p$. Then there is an equivalence of categories
\begin{equation*}
\Perf(X^{\Syn/\C_p})\cong \Perf(X^{\Div^1/\C_p})\times_{\Perf(X^{\HT, \dagger/\C_p})} \Perf(X^{\HT, \dagger, +/\C_p})
\end{equation*}
induced by the realisation functors $T_{\Div^1}$ and $T_{\HT, \dagger, +}$. In particular, for any $E\in\Perf(X^{\Syn/\C_p})$, there is a pullback diagram
\begin{equation*}
\begin{tikzcd}
R\Gamma(X^{\Syn/\C_p}, E)\ar[r]\ar[d] & R\Gamma(X^{\Div^1/\C_p}, T_{\Div^1}(E))\ar[d] \\
R\Gamma(X^{\HT, \dagger, +/\C_p}, T_{\HT, \dagger, +}(E))\ar[r] & R\Gamma(X^{\HT, \dagger/\C_p}, T_{\HT, \dagger}(E))\nospacepunct{\;.}
\end{tikzcd}
\end{equation*}
\end{thm}
\begin{proof}
As in \cref{thm:proet-main}, replacing the use of \cref{prop:perf-coversmoothrigid} by \cref{lem:geom-coversmooth}.
\end{proof}

To control the cofibre of the map
\begin{equation*}
R\Gamma(X^{\HT, \dagger, +/\C_p}, T_{\HT, \dagger, +}(E))\rightarrow R\Gamma(X^{\HT, \dagger/\C_p}, T_{\HT, \dagger}(E))
\end{equation*}
and hence deduce an honest comparison theorem between our geometric syntomic cohomology and proétale cohomology in a suitable range, we again have to understand $X^{\HT, \dagger, +/\C_p}$ more explicitly, at least locally on $X$. Namely, we need the following analogue of \cref{thm:htdr-xhtdrdagger+pres}, \cref{cor:htdr-complexesxhtdrdagger+} and \cref{cor:htdr-cohomologyxhtdrdagger+}, which is proved in the same way:

\begin{prop}
\label{prop:geom-xhtdagger+}
Let $X$ be a derived Berkovich space over $\C_p$ and assume that $X$ admits a lift $\widetilde{X}$ to $B_\dR^{+, \dagger}$ equipped with a Berkovich étale map $\widetilde{X}\rightarrow\ol{\T}^n_{B_\dR^{+, \dagger}}$ for some $n\geq 0$. 
\begin{enumerate}[label=(\arabic*)]
\item There is an isomorphism 
\begin{equation*}
X^{\HT, \dagger, +/\C_p}\cong \left(\widetilde{X}_\infty^\la\times_{\GSpec B_\dR^{+, \dagger}} \GSpec B_\dR^{+, \dagger}\langle u\rangle_{\leq 1}\{t\}^\dagger/(ut-\xi)\;/\;\ol{\T}\right)\,\big/\, (\Z_p^\la)^n\coprod^\AbGrp_{(\G_a^\dagger)^n} (\G_a^\dagger)^n\;,
\end{equation*}
where the pushout is taken along the map $(\G_a^\dagger)^n\rightarrow (\G_a^\dagger)^n$ given by multiplication by $u$ and $\widetilde{X}_\infty^\la$ denotes the base change of $\widetilde{X}$ along $\ol{\T}^{n, \la}_\infty\rightarrow\ol{\T}^n$. Here, the action of $(\Z_p^\la)^n\coprod^\AbGrp_{(\G_a^\dagger)^n} (\G_a^\dagger)^n$ on $\widetilde{X}_\infty^\la$ is obtained by pulling back the action on $\ol{\T}^{n, \la}_\infty$ given by
\begin{equation*}
(\ul{\theta}, \ul{v}).(\ul{x}^{1/p^m})_m\coloneqq \left([\epsilon]^{\ul{\theta}/p^m}\exp\left(\frac{t'\ul{v}}{p^m}\right)\ul{x}^{1/p^m}\right)_m\;,
\end{equation*}
where $(\ul{x}^{1/p^m})_m$ are the coordinates on $\ol{\T}^{n, \la}_\infty$ and $t'\coloneqq t\cdot \log[\epsilon]/\xi$.

\item A perfect complex on $X^{\HT, \dagger, +/\C_p}$ amounts to the following data:
\begin{enumerate}[label=(\roman*)]
\item A diagram
\begin{equation*}
\begin{tikzcd}
\dots \ar[r,shift left=.5ex,"u"]
  & \ar[l,shift left=.5ex, "t"] \Fil_{i-1} V \ar[r,shift left=.5ex,"u"] & \ar[l,shift left=.5ex, "t"] \Fil_i V \ar[r,shift left=.5ex,"u"] & \ar[l,shift left=.5ex, "t"] \Fil_{i+1} V\ar[r,shift left=.5ex,"u"] & \ar[l,shift left=.5ex, "t"] \dots
\end{tikzcd}
\end{equation*}
of perfect complexes over $\widetilde{X}_\infty^\la$ such that $ut=\xi$ with the property that $t: \Fil_\bullet V\rightarrow \Fil_{\bullet-1} V$ becomes an isomorphism for $\bullet\ll 0$ and $u: \Fil_\bullet V\rightarrow\Fil_{\bullet+1} V$ becomes an isomorphism for $\bullet\gg 0$;
\item a $t$-connection $\nabla: V\rightarrow V\tensor_{\O_{\widetilde{X}_\infty^\la}} \Omega_{\widetilde{X}_\infty^\la/B_\dR^{+, \dagger}}^1$ on $V$, i.e.\ a map satisfying the $t$-deformed Leibniz rule $\nabla(fm)=f\nabla(m)+t\cdot\mathrm{d}f$, whose restriction to $\Fil_i V$ is equipped with a factorisation through $\Fil_{i-1} \tensor_{\O_{\widetilde{X}_\infty^\la}} \Omega_{\widetilde{X}_\infty^\la/B_\dR^{+, \dagger}}^1$;
\item a semilinear locally analytic $\Z_p^n$-action on $\Fil_\bullet V$ commuting with all the previous data;
\item an identification between the ``geometric Sen operator'' $\Theta^\geom: \Fil_i V\rightarrow \Fil_i V\tensor_{\O_{\widetilde{X}_\infty^\la}} \Omega_{\widetilde{X}_\infty^\la/B_\dR^{+, \dagger}}^1$ coming from the Lie algebra action of $\Z_p^n$ and the $t$-connection $u\nabla$ for each $i$.
\end{enumerate}
Under this identification, restriction to $X^{\HT, \dagger/\C_p}$ corresponds to forgetting the $u$-filtration.

\item Under the above identifications, the cohomology of a perfect complex on $X^{\HT, \dagger, +/\C_p}$ is computed by the (underived) $\Z_p^n$-invariants of the cohomology on $\widetilde{X}_\infty^\la$ of the total complex of
\begin{equation*}
\begin{tikzcd}
\Fil_0 V\ar[r, "\nabla"] & \Fil_{-1} V\tensor_{\O_{\widetilde{X}_\infty^\la}} \Omega_{\widetilde{X}_\infty^\la/B_\dR^{+, \dagger}}^1\ar[r, "\nabla"] & \Fil_{-2} V\tensor_{\O_{\widetilde{X}_\infty^\la}} \Omega_{\widetilde{X}_\infty^\la/B_\dR^{+, \dagger}}^2\ar[r, "\nabla"] & \dots
\end{tikzcd}
\end{equation*}
while the cohomology of its restriction to $X^{\HT, \dagger/\C_p}$ is given by the (underived) $\Z_p^n$-invariants of the cohomology on $\widetilde{X}_\infty^\la$ of the total complex of
\begin{equation*}
\begin{tikzcd}
V\ar[r, "\nabla"] & V\tensor_{\O_{\widetilde{X}_\infty^\la}} \Omega_{\widetilde{X}_\infty^\la/B_\dR^{+, \dagger}}^1\ar[r, "\nabla"] & V\tensor_{\O_{\widetilde{X}_\infty^\la}} \Omega_{\widetilde{X}_\infty^\la/B_\dR^{+, \dagger}}^2\ar[r, "\nabla"] & \dots\nospacepunct{\;,}
\end{tikzcd}
\end{equation*}
where $V$ denotes the underlying unfiltered object of $\Fil_\bullet V$.
\end{enumerate}
\end{prop}
\begin{proof}
As in \cref{thm:htdr-xhtdrdagger+pres}, \cref{cor:htdr-complexesxhtdrdagger+} and \cref{cor:htdr-cohomologyxhtdrdagger+}.
\end{proof}

Using \cref{prop:geom-xhtdagger+} and that cohomology on $X^{\Div^1/\C_p}$ computes proétale cohomology of $X$ by \cite{AnPrism} if $X$ is a smooth partially proper rigid space over $\C_p$, we deduce the following explicit comparison theorem from \cref{thm:geom-proet1}:

\begin{thm}
\label{thm:geom-proet2}
Let $X$ be a smooth partially proper rigid space over $\C_p$. If $E\in\Vect(X^{\Syn/\C_p})$ is a vector bundle analytic $F$-gauge on $X$ relative to $\C_p$ with Hodge--Tate weights all at most $-i$ for some $i\geq 0$, then the natural morphism
\begin{equation*}
R\Gamma(X^{\Syn/\C_p}, E)\rightarrow R\Gamma(X^{\Div^1/\C_p}, T_{\Div^1}(E))
\end{equation*}
is an isomorphism on $\tau^{\leq i}$ and induces an injection on $H^{i+1}$. In particular, for $E=\O\{i\}$, we obtain
\begin{equation*}
\tau^{\leq i} R\Gamma_\Syn(X/\C_p, \Q_p(i))\cong \tau^{\leq i}R\Gamma_\proet(X, \Q_p(i))\;.
\end{equation*}
\end{thm}
\begin{proof}
As in \cref{thm:proet-main2}.
\end{proof}

Finally, we can again combine the results of \cref{thm:geom-hkcomp} and \cref{thm:geom-proet2} to obtain a version of Colmez--Nizio{\l}'s ``basic comparison theorem'' in the geometric case, see \cite[Thm.\ 1.1]{BasicComparison}, with vector bundle $F$-gauge coefficients. However, as already pointed out above, our version differs from the one in loc.\ cit.\ due to the appearance of the period ring $B_{\log}$ in place of $B_\st^+$ and should thus rather be compared to the result of Bosco from \cite[Thm.\ 7.2]{BoscopAdic}.

\begin{cor}
Let $X$ be a smooth partially proper rigid space over $\C_p$. If $E\in\Vect(X^{\Syn/\C_p})$ is a vector bundle analytic $F$-gauge on $X$ relative to $\C_p$ with Hodge--Tate weights all at most $-i$ for some $i\geq 0$, then there is a natural map
\begin{equation*}
\begin{split}
\fib(R\Gamma(X^{\HK/\C_p}, T_\HK(E))\rightarrow R\Gamma(X^{\dR/B_\dR^{+, \dagger}}, T_\dR(E))/R\Gamma(X^{\dR, +/B_\dR^{+, \dagger}}, &T_{\dR, +}(E))) \\
&\rightarrow R\Gamma(X^{\Div^1}, T_{\Div^1}(E))
\end{split}
\end{equation*}
which induces an isomorphism on $\tau^{\leq i}$ and an injection on $H^{i+1}$. In particular, if $X$ is qcqs and $E=\O\{i\}$, we obtain
\begin{equation*}
\tau^{\leq i}R\Gamma_\proet(X, \Q_p(i))\cong \tau^{\leq i}((R\Gamma_\HK(X)\tensor_{\Q_p^\un} B_{\log})^{\phi=p^i, N=0}\rightarrow R\Gamma(X/B_\dR^+)/\Fil^i R\Gamma(X/B_\dR^+))\;.
\end{equation*}
\end{cor}
\begin{proof}
The first part follows by combining \cref{thm:geom-hkcomp} and \cref{thm:geom-proet2} while the second part then follows using \cref{cor:geom-hkcompclassical}.
\end{proof}

%
%
%
%
%
%
%
%

\comment{

\section{Junk}

\subsection{The Hodge--Tate to de Rham degeneration}
\label{sect:htdr}

For a Gelfand stack $X$, the closed substack $(X^\N)_{|ut|=0}$ is in itself already interesting: It provides a deformation from $(X^\N)_{|u|=1, |t|=0}\cong X^{\HT, \dagger}$ to $(X^\N)_{|t|=1, |u|=0}\cong X^\dR$.

\begin{lem}
\label{lem:htdr-perfut0}
Let $X$ be any Gelfand stack and assume that $(X^\N)_{|ut|=0}$ is nicely coverable. Then
\begin{equation*}
\Perf((X^\N)_{|ut|=0})\cong \Perf(X^{\dR, +})
\end{equation*}
via pullback along the filtered de Rham map $i_{\dR, +}$.
\end{lem}
\begin{proof}
By virtue of the pushout diagram
\begin{equation*}
\begin{tikzcd}
(X^\N)_{|u|=|t|=0}\ar[r]\ar[d] & (X^\N)_{|t|=0}\ar[d] \\
(X^\N)_{|u|=0}\ar[r] & (X^\N)_{|ut|=0}
\end{tikzcd}
\end{equation*}
and the fact that $(X^\N)_{|u|=0}\cong X^{\dR, +}$ via $i_{\dR, +}$ by \cref{prop:defis-u0}, it suffices to show that pullback along the map $(X^\N)_{|u|=|t|=0}\rightarrow (X^\N)_{|t|=0}$ induces an equivalence on perfect complexes. To this end, consider the overconvergent normed divisor $Z=\{|u|=0\}\subseteq (X^\N)_{|t|=0}$ and observe that $Z_\epsilon=(X^\N)_{|u|\leq \epsilon, |t|=0}$ for each $\epsilon>0$. This shows
\begin{equation*}
Z_\epsilon\setminus Z\cong (X^\N)_{0<|u|\leq\epsilon, |t|=0}\cong X^{\HT, \dagger}\times (0, \epsilon]\;
\end{equation*}
by \cref{prop:defis-t0uneq0}; in particular, setting $\epsilon=1$ yields
\begin{equation*}
(X^\N)_{|t|=0}\setminus Z=Z_1\setminus Z\cong X^{\HT, \dagger}\times (0, 1]\;.
\end{equation*}
Thus, by \cref{cor:perf-contractible}, we obtain
\begin{equation*}
\Perf((X^\N)_{|t|=0}\setminus Z)\cong\Perf(X^{\HT, \dagger})\cong \Perf(Z_\epsilon\setminus Z)
\end{equation*}
for each $\epsilon>0$ and, in particular,
\begin{equation*}
\Perf((X^\N)_{|t|=0}\setminus Z)\cong\colim_{\epsilon>0} \Perf(Z_\epsilon\setminus Z)
\end{equation*}
via pullback. Hence, using \cref{cor:perf-corkeylemma}, we conclude that
\begin{equation*}
\begin{split}
\Perf((X^\N)_{|t|=0})&\cong \Perf((X^\N)_{|u|=|t|=0})\times_{\colim_{\epsilon>0} \Perf(Z_\epsilon\setminus Z)} \Perf((X^\N)_{|t|=0}\setminus Z) \\
&\cong \Perf((X^\N)_{|u|=|t|=0})
\end{split}
\end{equation*}
via pullback, as claimed.
\end{proof}

Doing the same argument ``in the other direction'', i.e.\ exchanging the roles of $u$ and $t$ and instead proving $\Perf((X^\N)_{|u|=0})\cong \Perf((X^\N)_{|u|=|t|=0})$ by using the fact that $(X^\N)_{|u|=0, |t|\neq 0}\cong X^\dR\times (0, 1]$, which can easily be deduced from \cref{prop:defis-u0}, we obtain the following variant:

\begin{lem}
\label{lem:htdr-perfut0viaht+}
Let $X$ be any Gelfand stack and assume that $(X^\N)_{|ut|=0}$ is nicely coverable. Then
\begin{equation*}
\Perf((X^\N)_{|ut|=0})\cong \Perf(X^{\HT, \dagger, +})\;,
\end{equation*}
where we define $X^{\HT, \dagger, +}\coloneqq (X^\N)_{|t|=0}$.
\end{lem}

In particular, we obtain

\begin{cor}
\label{cor:htdr-perfequiv}
Let $X$ be a Berkovich smooth derived Berkovich space over $\Q_p$ or over $\C_p$. Then there is an equivalence of categories
\begin{equation*}
\Perf(X^{\dR, +})\cong \Perf(X^{\HT, \dagger, +})\;.
\end{equation*}
\end{cor}
\begin{proof}
We have natural functors $\Perf(X^{\dR, +})\rightarrow \Perf((X^\N)_{|ut|=0})$ and $\Perf(X^{\HT, \dagger, +})\rightarrow \Perf((X^\N)_{|ut|=0})$ and our goal is to show that they are both equivalences. By compatibility of everything with Berkovich étale localisations and strict closed covers, we may assume that $X$ is finite étale over a rational localisation of $\ol{\DD}^n$ or $\ol{\DD}_{\C_p}^n$, respectively. However, in these cases, the claim follows by combining the previous two lemmas with \cref{prop:perf-coversmoothrigid} or \cref{cor:perf-rigidgeomcover}, respectively.
\end{proof}

In the case of a smooth partially proper rigid space $X$ over $\Q_p$, the above result should be related to a functor introduced by Scholze in \cite[§7]{PAdicHodgeTheory}. For this, we make the following definition, see \cite[Def.\ 3.4.1]{Wiersig1}:

\begin{defi}
\label{defi:htdr-bdr+dagger}
Let $(R, R^+)$ be a perfectoid Huber pair over $\Q_p$. The \emph{overconvergent (positive) de Rham period ring} $\mathbb{B}_\dR^{+, \dagger}(R, R^+)$ is defined by 
\begin{equation*}
\mathbb{B}_\dR^{+, \dagger}(R, R^+)\coloneqq \mathbb{A}_\inf(R, R^+)[\tfrac{1}{p}]\{\xi\}^\dagger=\colim_n \mathbb{A}_\inf(R, R^+)\langle p^{-n}\xi\rangle[\tfrac{1}{p}]\;,
\end{equation*}
where $\xi$ denotes a generator of the kernel of Fontaine's map $\mathbb{A}_\inf(R, R^+)\rightarrow R^+$. By \cite[Thm.\ 3.4.2]{Wiersig1}, the association $(R, R^+)\mapsto \mathbb{B}_\dR^{+, \dagger}(R, R^+)$ determines a sheaf on the proétale site of any arc-stack called the \emph{overconvergent (positive) de Rham period sheaf} which is denoted $\mathbb{B}_\dR^{+, \dagger}$.
\end{defi}

Now observe that, if $X$ is a smooth partially proper rigid space over $\Q_p$, the base change
\begin{equation*}
(X^\diamond\times \Spd\Q_p)\times_{X^\prism} X^{\HT, \dagger}
\end{equation*}
coincides with the overconvergent neighbourhood of $X^\diamond\subseteq X^\diamond\times\Spd\Q_p$. Thus, recalling that vector bundles on $X^{\dR, +}$ are equivalent to filtered vector bundles on $X$ with a flat connection, we see that postcomposing the analogue of the equivalence from \cref{cor:htdr-perfequiv} for vector bundles in place of perfect complexes with pullback along
\begin{equation*}
(X^\diamond\subseteq X^\diamond\times\Spd\Q_p)^\dagger\rightarrow X^{\HT, \dagger}\rightarrow X^{\HT, \dagger, +}
\end{equation*}
yields a functor
\begin{equation}
\label{eq:htdr-scholzefunctor}
\{\text{filtered vector bundles on $X$ with connection}\}\rightarrow \{\mathbb{B}_\dR^{+, \dagger}\text{-local systems on }X_\proet\}\;.
\end{equation}
Of course, we can also get a functor toward $\mathbb{B}_\dR^+$-local systems on $X_\proet$ by extending scalars in the target.

Our goals in this section are twofold: 
\begin{enumerate}[label=(\arabic*)]
\item For rigid smooth derived Berkovich spaces $X$, we want to explicitly understand the stack $X^{\HT, \dagger, +}$, its category of perfect complexes and how to compute its coherent cohomology, at least locally on $X$.
\item For smooth partially proper rigid spaces $X$, we want to check that the restriction functor 
\begin{equation*}
\Vect(X^{\dR, +})\cong \Vect(X^{\HT, \dagger, +})\rightarrow \Vect(X^{\HT, \dagger})
\end{equation*}
is fully faithful and that its essential image is given by the full subcategory of vector bundles on $X^{\HT, \dagger}$ which pull back to $\mathbb{B}_\dR^+$-local systems which are de Rham.
\end{enumerate}
We recall that a $\mathbb{B}_\dR^+$-local system on $X$ is called \emph{de Rham} if it lies in the essential image of Scholze's functor
\begin{equation*}
\{\text{filtered vector bundles on $X$ with flat connection}\}\rightarrow \{\mathbb{B}_\dR^+\text{-local systems on }X_\proet\}
\end{equation*} 
from \cite[§7]{PAdicHodgeTheory}. Of course, a direct way to prove (2) would be to identify the functor (\ref{eq:htdr-scholzefunctor}) we have constructed above with Scholze's functor up to extension of scalars. While we expect this to be true, we have chosen not to pursue this strategy since it seems to require considerably more effort than the one below and, for our purposes, having (2) is already enough and we do not require the finer knowledge that our functor is the same as the one of Scholze.

\subsubsection{The case $X=\GSpec\Q_p$}

We start with the case $X=\GSpec\Q_p$, which is already interesting. In this case, we will actually be able to identify the functor (\ref{eq:htdr-scholzefunctor}) with Scholze's functor.

\begin{prop}
\label{prop:htdr-qphtdagger+}
There is an isomorphism
\begin{equation*}
\Q_p^{\HT, \dagger, +}\cong \left(\ol{\DD}^\dR\times\G_a^\dagger\times\GSpec\Q_p(\zeta_{p^\infty})\right)\;\Big/\;\Z_p^{\times, \la}\coprod^\AbGrp_{\G_m^\dagger} \G_m(1)\;.
\end{equation*}
Here, the pushout is taken in the category of abelian group stacks and the group actions defining the quotient are given as follows:
\begin{enumerate}[label=(\roman*)]
\item $\Z_p^{\times, \la}$ acts on $\GSpec\Q_p(\zeta_{p^\infty})$ via the usual Galois action and on $\G_a^\dagger$ via $\gamma.t=\gamma\cdot t$.
\item $\G_m(1)$ acts on $\ol{\DD}^\dR$ by division and on $\G_a^\dagger$ by multiplication.
\end{enumerate}
In particular, restricting to $|u|=1$, we have
\begin{equation*}
\Q_p^{\HT, \dagger}\cong \left(\G_a^\dagger\times\GSpec\Q_p(\zeta_{p^\infty})\right)\,\big/\,\Z_p^{\times, \la}\;.
\end{equation*}
\end{prop}
\begin{proof}
Note that, by definition, the stack $\Q_p^{\HT, \dagger, +}$ sits inside a pullback diagram
\begin{equation*}
\begin{tikzcd}
\Q_p^{\HT, \dagger, +}\ar[d]\ar[r] & \GSpec\Q_p\ar[d] \\
\G_a^\dagger/\G_m(1)\times (\ol{\DD}/\G_m(1))^\dR\ar[r, "\mathrm{mult}"] & */\G_m(1)^\dR\nospacepunct{\;,}
\end{tikzcd}
\end{equation*}
where the vertical map on the right classifies the Breuil--Kisin twist and factoring the bottom map as
\begin{equation*}
\G_a^\dagger/\G_m(1)\times (\ol{\DD}/\G_m(1))^\dR\rightarrow */\G_m(1)^\dR\times (\ol{\DD}/\G_m(1))^\dR\xrightarrow{\mathrm{mult}} */\G_m(1)^\dR
\end{equation*}
shows that we can simplify this to a pullback diagram
\begin{equation*}
\begin{tikzcd}
\Q_p^{\HT, \dagger, +}\ar[r]\ar[d] & (\ol{\DD}/\G_m(1))^\dR\ar[d] \\
\G_a^\dagger/\G_m(1)\ar[r] & */\G_m(1)^\dR\;,
\end{tikzcd}
\end{equation*}
where the map on the right-hand side is the first Breuil--Kisin twist of the canonical map. 

Thus, we have to show how
\begin{equation}
\label{eq:htdr-qphtdaggerpluspresentation}
\left(\ol{\DD}^\dR\times\G_a^\dagger\times\GSpec\Q_p(\zeta_{p^\infty})\right)\;\Big/\;\Z_p^{\times, \la}\coprod^\AbGrp_{\G_m^\dagger} \G_m(1)
\end{equation}
fits into the same pullback diagram. The map to $(\ol{\DD}/\G_m(1))^\dR$ is evident and just comes from projecting onto the first factor ``in the numerator'' and then using the map
\begin{equation*}
\Z_p^{\times, \la}\coprod^\AbGrp_{\G_m^\dagger} \G_m(1)\rightarrow \G_m(1)^\dR\;.
\end{equation*}
For the map to $\G_a^\dagger/\G_m(1)$, we have to give a normed line bundle on the quotient stack above together with a $\dagger$-nilpotent section. For this, first note that, on the quotient
\begin{equation*}
\left(\ol{\DD}^\dR\times\G_a^\dagger\times\GSpec\Q_p(\zeta_{p^\infty})\right)\;\Big/\; \Z_p^{\times, \la}\times \G_m(1)\;,
\end{equation*}
we have a tautological line bundle from the $\G_m(1)$-action as well as a ``Tate twist'' coming from the $\Z_p^{\times, \la}$-action and that the antidiagonal copy of $\G_m^\dagger\subseteq \Z_p^{\times, \la}\times\G_m(1)$ acts trivially on the tensor product of these two line bundles. Thus, this tensor product descends to (\ref{eq:htdr-qphtdaggerpluspresentation}) and, by definition of the action of $\Z_p^{\times, \la}$ and $\G_m(1)$ on $\G_a^\dagger$, it is clear that the coordinate $t$ on $\G_a^\dagger$ defines a section of the descended bundle. This defines the desired map to $\G_a^\dagger/\G_m(1)$ and it is clear by construction that this is compatible with the map to $(\ol{\DD}/\G_m(1))^\dR$ in the sense of the desired pullback diagram (recall that one can also interpret the Breuil--Kisin twist as a Tate twist).

Thus, by the description of $\Q_p^{\HT, \dagger, +}$ as a pullback from above, we obtain a map
\begin{equation*}
\left(\ol{\DD}^\dR\times\G_a^\dagger\times\GSpec\Q_p(\zeta_{p^\infty})\right)\;\Big/\;\Z_p^{\times, \la}\coprod^\AbGrp_{\G_m^\dagger} \G_m(1)\rightarrow \Q_p^{\HT, \dagger, +}\;.
\end{equation*}
However, the same pullback diagram presents $\Q_p^{\HT, \dagger, +}$ as a $\ol{\DD}^\dR$-torsor over $\G_a^\dagger/\G_m(1)$ and the map
\begin{equation*}
\left(\ol{\DD}^\dR\times\G_a^\dagger\times\GSpec\Q_p(\zeta_{p^\infty})\right)\;\Big/\;\Z_p^{\times, \la}\coprod^\AbGrp_{\G_m^\dagger} \G_m(1)\rightarrow \G_a^\dagger/\G_m(1)
\end{equation*}
we have constructed is a $\ol{\DD}^\dR$-torsor as well, hence we are done upon noting that the above map is actually a map of $\ol{\DD}^\dR$-torsors.

Finally, let us explain how to obtain the formula for $\Q_p^{\HT, \dagger}$ from the one for $\Q_p^{\HT, \dagger, +}$. For this, note that there is an obvious projection map from
\begin{equation*}
\left(\G_m(1)\times\G_a^\dagger\times\GSpec\Q_p(\zeta_{p^\infty})\right)\;\Big/\;\Z_p^{\times, \la}\times\G_m(1)
\end{equation*}
to the stack
\begin{equation*}
\Q_p^{\HT, \dagger}\cong (\Q_p^{\HT, \dagger, +})_{|u|=1}\cong \left(\G_m(1)^\dR\times\G_a^\dagger\times\GSpec\Q_p(\zeta_{p^\infty})\right)\;\Big/\;\Z_p^{\times, \la}\coprod^\AbGrp_{\G_m^\dagger} \G_m(1)\;,
\end{equation*}
which is actually an isomorphism: Indeed, pulling back along the surjection 
\begin{equation*}
\G_m(1)^\dR\times \G_a^\dagger\times \GSpec\Q_p(\zeta_{p^\infty})\rightarrow \Q_p^{\HT, \dagger}\;,
\end{equation*}
the above map becomes
\begin{equation*}
\left(\G_m(1)\times \G_a^\dagger\times\GSpec\Q_p(\zeta_{p^\infty})\right)\,\big/\,\G_m^\dagger\rightarrow \G_m(1)^\dR\times \G_a^\dagger\times \GSpec\Q_p(\zeta_{p^\infty})
\end{equation*}
and this is obviously an isomorphism; to see this, note that the action of $\G_m^\dagger$ on the left-hand side comes from viewing $\G_m^\dagger$ as antidiagonally embedded into $\G_m(1)\times\Z_p^{\times, \la}$ and thus is only nontrivial on $\G_m(1)$, where it is given by multiplication. Then it remains to show that
\begin{equation*}
\left(\G_m(1)\times\G_a^\dagger\times\GSpec\Q_p(\zeta_{p^\infty})\right)\;\Big/\;\Z_p^{\times, \la}\times\G_m(1)\cong \left(\G_a^\dagger\times\GSpec\Q_p(\zeta_{p^\infty})\right)\,\big/\,\Z_p^{\times, \la}\;,
\end{equation*}
but this isomorphism is induced by the automorphism of $\G_m(1)\times \G_a^\dagger$ given by $(u, t)\mapsto (u, ut)$, after which the $\G_m(1)$-action on $\G_a^\dagger$ becomes trivial.
\end{proof}

Unraveling all the linear algebraic structures encoded by the various actions, we can now describe perfect complexes on $\Q_p^{\HT, \dagger, +}$. We start with the following lemma:

\begin{lem}
\label{lem:htdr-ddrmodgmdagger}
A perfect complex on $\ol{\DD}_-^\dR/\G_m(1)$ amounts to an increasingly filtered perfect complex $u: \Fil_\bullet V\rightarrow \Fil_{\bullet+1} V$ of $\Q_p$-vector spaces together with linear operators $D: \Fil_\bullet V\rightarrow \Fil_{\bullet-1} V$ satisfying $Du=uD+1$. Moreover, under this correspondence, pullback along
\begin{equation*}
\ol{\DD}^\dR\times */\G_m^\dagger\cong \ol{\DD}^\dR/\G_m^\dagger\rightarrow \ol{\DD}^\dR/\G_m(1)
\end{equation*}
corresponds to passing to the $D$-module $\bigoplus_{i\in \Z} \Fil_i V\cdot u^i$ on $\ol{\DD}$, where the action of $\partial_u$ is given by $D$, together with the $\Q_p[u]$-linear endomorphism $\Theta$ given by $uD-i$ on $\Fil_i V\cdot u^i$.
\end{lem}
\begin{proof}
The first part is a variant of \cite[§6.5.4]{FGauges}, hence we focus on the second part. For this, first note that pullback along $\ol{\DD}^\dR/\G_m^\dagger\rightarrow\ol{\DD}^\dR/\G_m(1)$ corresponds to passing to the action of the Lie algebra of the analytic group $\G_m(1)$ and, in the present case, this corresponds to forgetting the grading on $\bigoplus_{i\in \Z} \Fil_i V\cdot u^{-i}$ and instead remembering the operator $\Psi$ which acts by multiplication by $-i$ on $\Fil_i V$ as $\G_m(1)$ acts on $u$ by $s.u=u/s$.

Thus, the nontrivial part lies in understanding the isomorphism $\ol{\DD}^\dR\times */\G_m^\dagger\cong \ol{\DD}^\dR/\G_m^\dagger$. More specifically, perfect complexes on $\ol{\DD}^\dR/\G_m^\dagger$ can be described as $D$-modules on $\ol{\DD}$ together with an operator $\Psi$ satisfying $\Psi u=u\Psi-u$ and $\Psi D=D\Psi+D$ via descent along $\ol{\DD}/\G_m^\dagger\rightarrow \ol{\DD}^\dR/\G_m^\dagger$ while perfect complexes on $\ol{\DD}^\dR\times */\G_m^\dagger$ can be described as $D$-modules on $\ol{\DD}$ together with an operator $\Theta$ commuting with $u$ and $D$ via descent along $\ol{\DD}\times */\G_m^\dagger\rightarrow \ol{\DD}^\dR\times */\G_m^\dagger$. We claim that the translation between these two descriptions is given by
\begin{equation*}
\Theta=uD+\Psi
\end{equation*}
and this will immediately imply the statement of the lemma.

Indeed, note that we can write
\begin{equation*}
\ol{\DD}^\dR/\G_m^\dagger=\ol{\DD}/(\G_a^\dagger\rtimes\G_m^\dagger)\;,
\end{equation*}
where the group structure on the semidirect product is given by
\begin{equation*}
(t_1, s_1)(t_2, s_2)=(t_1s_2+t_2, s_1s_2)
\end{equation*}
and $(t, s)$ acts on the coordinate $u$ of $\ol{\DD}$ via
\begin{equation*}
(t, s).u=(u+t)/s\;.
\end{equation*}
Moreover, the morphism $\ol{\DD}/\G_m^\dagger\rightarrow \ol{\DD}/(\G_a^\dagger\rtimes\G_m^\dagger)$ is induced by the inclusion $\G_m^\dagger\rightarrow \G_a^\dagger\rtimes\G_m^\dagger$ given by $s\mapsto (0, s)$ while the morphism $\ol{\DD}\times */\G_m^\dagger\rightarrow \ol{\DD}/(\G_a^\dagger\rtimes\G_m^\dagger)$ is induced by the morphism
\begin{equation*}
\begin{split}
\ol{\DD}\times \G_m^\dagger&\rightarrow \ol{\DD}\times (\G_a^\dagger\rtimes\G_m^\dagger) \\
(u, s)&\mapsto (u, u(s-1), s)
\end{split}
\end{equation*}
of analytic groups over $\ol{\DD}$; indeed, note that the action of $(u(s-1), s)$ fixes $u$. However, as the operators $\Psi$ and $\Theta$ both come from Cartier duality for $\G_a^\dagger\cong \G_m^\dagger$, where the isomorphism is via the exponential, we should furthermore replace occurrences of $\G_m^\dagger$ above by $\G_a^\dagger$ and replace the coordinate $s$ on $\G_m^\dagger$ by the coordinate $\ell=\log s$ on $\G_a^\dagger$. Then the above morphism is replaced by
\begin{equation*}
\begin{split}
\ol{\DD}\times \G_a^\dagger&\rightarrow \ol{\DD}\times (\G_a^\dagger\rtimes\G_a^\dagger) \\
(u, \ell)&\mapsto (u, u(\exp\ell-1), \ell)\;.
\end{split}
\end{equation*}

To put everything together, note that, by the above, the description of perfect complexes on $\ol{\DD}^\dR/\G_m^\dagger$ via $\Psi$ can be reformulated algebraically as follows: a perfect complex on $\ol{\DD}^\dR/\G_m^\dagger$ corresponds to a perfect complex $M$ over $\O(\ol{\DD})=\Q_p\langle u\rangle_{\leq 1}$ equipped with the structure of a right comodule for the (non-cocommutative!) Hopf algebra $\O(\G_a^\dagger\rtimes\G_a^\dagger)$; in terms of the operator $\Psi$, this comodule structure is given by
\begin{equation*}
m\mapsto \sum_{i, j\geq 0} \frac{t^i}{i!}\cdot\frac{\ell^j}{j!}\cdot D^i\Psi^j(m)\;.
\end{equation*}
Then note that, by the above discussion, the pullback of such a perfect complex along the map $\ol{\DD}\times */\G_m^\dagger\rightarrow\ol{\DD}^\dR/\G_m^\dagger$ is given by the perfect complex $M$ together with the right comodule structure
\begin{equation*}
M\rightarrow M\tensor_{\Q_p\langle u\rangle_{\leq 1}} \O(\ol{\DD}\times \G_a^\dagger\rtimes\G_a^\dagger)\rightarrow M\tensor_{\Q_p\langle u\rangle_{\leq 1}} \O(\ol{\DD}\times \G_a^\dagger)
\end{equation*}
induced by the previous one via restriction along the map $t\mapsto u(\exp\ell-1)$ from above, i.e.\ it is given by
\begin{equation*}
m\mapsto \sum_{i, j\geq 0} \frac{u^i(\exp\ell-1)^i}{i!}\cdot\frac{\ell^j}{j!}\cdot D^i\Psi^j(m)\;.
\end{equation*}
Comparing this with the formula
\begin{equation*}
m\mapsto \sum_{k\geq 0} \frac{\ell^k}{k!}\cdot \Theta^k(m)
\end{equation*}
we get from Cartier duality, we see that $\Theta=uD+\Psi$ looking at the linear terms, as desired.
\end{proof}

\begin{prop}
\label{prop:htdr-complexesqphtdagger+}
A perfect complex on $\Q_p^{\HT, \dagger, +}$ amounts to the following data:
\begin{enumerate}[label=(\roman*)]
\item A diagram
\begin{equation*}
\begin{tikzcd}
\dots \ar[r,shift left=.5ex,"u"]
  & \ar[l,shift left=.5ex, "t"] \Fil_{i-1} V \ar[r,shift left=.5ex,"u"] & \ar[l,shift left=.5ex, "t"] \Fil_i V \ar[r,shift left=.5ex,"u"] & \ar[l,shift left=.5ex, "t"] \Fil_{i+1} V\ar[r,shift left=.5ex,"u"] & \ar[l,shift left=.5ex, "t"] \dots
\end{tikzcd}
\end{equation*}
of perfect complexes over $\GSpec\Q_p(\zeta_{p^\infty})\times \GSpec\Q_p\{r\}^\dagger$, where $r=ut$, with the property that $t: \Fil_\bullet V\rightarrow \Fil_{\bullet-1} V$ becomes an isomorphism for $\bullet\ll 0$ and $u: \Fil_\bullet V\rightarrow\Fil_{\bullet+1} V$ becomes an isomorphism for $\bullet\gg 0$;
\item $\Q_p(\zeta_{p^\infty})$-linear maps $D: \Fil_i V\rightarrow \Fil_{i-1} V$ for each $i$ satisfying $Du=uD+1$ and commuting with $t$;
\item a semilinear locally analytic $\Z_p^\times$-action on $\Fil_\bullet V$ commuting with $u$ and $D$, where the action of $\Z_p^\times$ on $\Q_p(\zeta_{p^\infty})$ is the usual Galois action and the one on $t$ is given by $\gamma.t=\gamma\cdot t$;
\item for each $i$, an identification between the ``Sen operator'' $\Theta: \Fil_i V\rightarrow \Fil_i V$ coming from the Lie algebra action of $\Z_p^\times$ and the operator $uD-i$.
\end{enumerate}
Under this identification, restriction to $\Q_p^{\HT, \dagger}$ corresponds to forgetting the $u$-filtration and the $D$-maps while restriction to $(\Q_p^\N)_{|u|=|t|=0}$ corresponds to passing to the associated graded of the $u$-filtration and forgetting the $D$-maps. Finally, the isomorphism
\begin{equation*}
(\Q_p^\N)_{|u|=|t|=0}\cong (\Q_p^{\dR, +})_{|t|=0}=\G_a^\dagger/\G_m(1)
\end{equation*}
Tate twists the $\Z_p^\times$-action on the $i$-th associated graded piece of the $u$-filtration by $i$, which trivialises the Lie algebra action and hence after descent along $\Q_p\rightarrow \Q_p(\zeta_{p^\infty})$ via the now smooth $\Z_p^\times$-action yields a filtered perfect complex of $\Q_p$-vector spaces via the $t$-maps.
\end{prop}
\begin{proof}
Item (i) is a version of the Rees equivalence while (ii) comes from the previous lemma and (iii) is clear. Finally, (iv) comes from the fact that the actions of $\Z_p^{\times, \la}$ and $\G_m(1)$ have to be compatible upon restriction to $\G_m^\dagger$: Indeed, on the one hand, restricting along $\G_m^\dagger\rightarrow\Z_p^{\times, \la}$ corresponds to passing to the action of the Lie algebra of $\Z_p^\times$. On the other hand, the compatibility between the $\Z_p^{\times, \la}$- and $\G_m(1)$-actions on $\ol{\DD}^\dR$ uses that the action of $\G_m^\dagger$ on $\ol{\DD}^\dR$ is trivial (recall that $\Z_p^{\times, \la}$ acts trivially on $\ol{\DD}^\dR$), i.e.\ goes through the isomorphism $\ol{\DD}^\dR/\G_m^\dagger\cong \ol{\DD}^\dR\times */\G_m^\dagger$, and thus, the previous lemma shows that restricting along $\G_m^\dagger\rightarrow \G_m(1)$ corresponds to passing to the operator $uD-i$ on $\Fil_i V$.

The rest of the statement is clear except maybe the twist occurring in translating between perfect complexes on $(\Q_p^\N)_{|u|=|t|=0}$ and $(\Q_p^{\dR, +})_{|t|=0}$ via the isomorphism $(\Q_p^{\dR, +})_{|t|=0}\cong (\Q_p^\N)_{|u|=|t|=0}$. However, this just comes from the fact that this isomorphism itself involves a Breuil--Kisin twist.
\end{proof}

In view of the above discussion, we can also give the following slightly different presentation of the stack $\Q_p^{\HT, \dagger, +}$ which makes the compatibility condition in (iv) above even more apparent. For this, let $\cal{G}_0$ be the group stack over $\ol{\DD}$ defined as the coequaliser of the map
\begin{equation*}
\begin{split}
\G_m^\dagger&\rightarrow \G_a^\dagger\rtimes\G_m(1)\xrightarrow{\mathrm{incl}_2} \Z_p^{\times, \la}\times (\G_a^\dagger\rtimes\G_m(1))\;, \\
s&\mapsto ((u-1)s, s)
\end{split}
\end{equation*}
and the inclusion map $\G_m^\dagger\rightarrow \Z_p^{\times, \la}\times (\G_a^\dagger\rtimes\G_m(1))$ into the last factor in the category of group stacks over $\ol{\DD}$, where we recall that $u$ denotes the coordinate on $\ol{\DD}$ and the semidirect product is defined using the multiplication action of $\G_m(1)$ on $\G_a^\dagger$. Then we have
\begin{equation*}
\Q_p^{\HT, \dagger, +}\cong \left(\ol{\DD}\times \G_a^\dagger\times \GSpec\Q_p(\zeta_{p^\infty})\right)\,\big/\,\cal{G}_0
\end{equation*}
for $\Z_p^{\times, \la}$ and $\G_m(1)$ acting on $\G_a^\dagger\times\GSpec\Q_p(\zeta_{p^\infty})$ as before and the action of $\G_a^\dagger\rtimes\G_m(1)$ on $\ol{\DD}$ being given by the formula
\begin{equation*}
(w, s).u=(u+w)/s
\end{equation*}
while $\Z_p^{\times, \la}$ acts trivially on $\ol{\DD}$. Indeed, note that the map $\G_m^\dagger\rightarrow \G_a^\dagger\rtimes\G_m(1)$ used to define $\cal{G}_0$ is precisely designed so that $\G_m^\dagger$ acts trivially on $\ol{\DD}$ via this map; moreover, it is also the same map that occurs in the proof of \cref{lem:htdr-ddrmodgmdagger}. Thus, one might argue that this presentation explains (iv) above much better than the presentation from \cref{prop:htdr-qphtdagger+}; in any case, it will be indispensable below when we want to describe $X^{\HT, \dagger, +}$ for more general smooth rigid spaces $X$ over $\Q_p$.

Perhaps a bit surprisingly, \cref{cor:htdr-perfequiv} tells us that all the data (i)-(iv) from \cref{prop:htdr-complexesqphtdagger+} is in fact just equivalent to the filtered perfect complex of $\Q_p$-vector spaces we obtain after pulling back to $\{|u|=0\}\subseteq \Q_p^{\HT, \dagger, +}$! To see this more explicitly, recall from above that the graded pieces $\gr_\bullet V$ of the filtration $\Fil_\bullet V$ are perfect complexes of $\Q_p(\zeta_{p^\infty})$-vector spaces equipped with a locally analytic $\Z_p^\times$-action such that $t^i \gr_i V$ descends to a perfect complex $\Fil^i W$ of $\Q_p$-vector spaces (this is because $\Theta$ acts by multiplication by $-i$ on $\gr_i V$ by the above and hence the $\Z_p^\times$-action on $t^i\gr_i V$ is smooth) -- this is just the $i$-th stage of the filtration associated to the restriction along $\G_a^\dagger/\G_m(1)\cong (\Q_p^\N)_{|u|=|t|=0}\rightarrow (\Q_p^\N)_{|t|=0}$. In other words, we have 
\begin{equation*}
\gr_i V\cong t^{-i}(\Fil^i W\tensor_{\Q_p} \Q_p(\zeta_{p^\infty}))\;,
\end{equation*}
but this implies that the filtration $\Fil_\bullet V$ splits $\Z_p^\times$-equivariantly! Indeed, we have
\begin{equation*}
\RHom_{\GSpec\Q_p(\zeta_{p^\infty})/\Z_p^{\times, \la}}(\gr_j V, \gr_i V)=0
\end{equation*}
for $i\neq j$; namely, by Tate twisting appropriately, this reduces to
\begin{equation*}
R\Gamma(\GSpec\Q_p(\zeta_{p^\infty})/\Z_p^{\times, \la}, t^{-i}\Q_p(\zeta_{p^\infty}))=0
\end{equation*}
for $i\neq 0$ and this cohomology is computed by taking $\Z_p^\times$-fixed points of Lie algebra cohomology as it coincides with pushforward along
\begin{equation*}
\GSpec\Q_p(\zeta_{p^\infty})/\Z_p^{\times, \la}\rightarrow \GSpec\Q_p(\zeta_{p^\infty})/\Z_p^{\times, \sm}\cong \GSpec\Q_p\;.
\end{equation*}
However, in generaly, for any complex of $\Q_p(\zeta_{p^\infty})$-vector spaces with a locally analytic $\Z_p^\times$-action, we have
\begin{equation*}
R\Gamma(\operatorname{Lie}\Z_p^\times, V)=\fib(V\xrightarrow{\Theta} V)
\end{equation*}
and since $\Theta=-i$ acts invertibly for $V=t^{-i}\Q_p(\zeta_{p^\infty})$ as long as $i\neq 0$, this proves the claim.

The above discussion shows that, from the data of a filtered perfect complex $\Fil^\bullet W$ of $\Q_p$-vector spaces, one can uniquely recover the data of a perfect complex on $(\Q_p^\N)_{|t|=0}$ as
\begin{equation*}
\Fil_\bullet V=\left(\bigoplus_{i\leq \bullet} t^{-i}(\Fil^i W\tensor_{\Q_p} \Q_p(\zeta_{p^\infty}))\right)\tensor_{\Q_p[r]} \Q_p\{r\}^\dagger
\end{equation*}
with $r: \Fil_\bullet V\rightarrow\Fil_\bullet V$ and $t: \Fil_\bullet V\rightarrow \Fil_{\bullet -1} V$ being given by the maps $t^{-i}\Fil^i W\rightarrow t^{-i+1} \Fil^{i-1} W$, the group $\Z_p^\times$ acting trivially on $\Fil^i W$ and by multiplication on $t$, the transition maps $u: \Fil_\bullet V\rightarrow \Fil_{\bullet+1} V$ being the obvious ones and $D: \Fil_\bullet V\rightarrow \Fil_{\bullet-1} V$ acting by multiplication by $\bullet-i$ on $t^{-i}(\Fil^i W\tensor_{\Q_p} \Q_p(\zeta_{p^\infty}))$.

\begin{defi}
A $G_{\Q_p}$-representation on a finite free $B_\dR^{+, \dagger}$-module is called \emph{generically flat} if it becomes trivial after extending scalars to $B_\dR^\dagger$.
\end{defi}

\begin{lem}
\label{lem:htdr-drrepfilteredvsp}
For any $\Q_p$-vector space $W$, there is an equivalence
\begin{equation*}
\begin{split}
\{\text{decreasing filtrations on $W$}\}&\overset{1:1}{\longleftrightarrow}\{\text{$G_{\Q_p}$-stable $B_\dR^{+, \dagger}$-lattices in $W\tensor_{\Q_p} B_\dR^\dagger$}\} \\
\Fil^\bullet W&\;\;\mapsto\;\; T=\Fil^0(W\tensor_{\Q_p} B_\dR^\dagger)\subseteq W\tensor_{\Q_p} B_\dR^\dagger \\
\Fil^\bullet W=(t^\bullet T)^{G_{\Q_p}}&\;\;\mapsfrom\;\;T\subseteq W\tensor_{\Q_p} B_\dR^\dagger\;.
\end{split}
\end{equation*}
In particular, generically flat representations on finite free $B_\dR^{+, \dagger}$-modules are equivalent to filtered $\Q_p$-vector spaces.
\end{lem}
\begin{proof}
Using the fact that $B_\dR^{+, \dagger}$ is a DVR with uniformiser $t$ and the calculation of Galois cohomology of $B_\dR^{+, \dagger}$ from \cite[Thm.\ G]{Wiersig2}, one can imitate the proof of \cite[Prop.\ 10.4.3]{FarguesFontaine}.
\end{proof}

\begin{cor}
Extension of scalars induces an equivalence between generically flat representations of $G_{\Q_p}$ over $B_\dR^{+, \dagger}$ and over $B_\dR^+$.
\end{cor}
\begin{proof}
Combine the previous lemma with \cite[Prop.\ 10.4.3]{FarguesFontaine}.
\end{proof}

\begin{rem}
Let us emphasise that the above corollary does \emph{not} imply that one can check whether a $B_\dR^{+, \dagger}$-representation is generically flat after extending scalars to $B_\dR^+$. While we expect this to be true, it is not needed for our purposes, see also \cref{rem:htdr-checkgenflat}.
\end{rem}

We can now check that the functor from filtered $\Q_p$-vector spaces to de Rham representations over $B_\dR^{\dagger, +}$ coincides with the functor
\begin{equation*}
\Vect((\Q_p^\N)_{|u|=|t|=0})\cong \Vect((\Q_p^\N)_{|t|=0})\rightarrow \Vect((\Q_p^\N)_{|t|=0, |u|=1})\rightarrow \Rep_{B_\dR^{\dagger, +}}(G_{\Q_p})
\end{equation*}
from the beginning of the section. Indeed, starting with a filtered $\Q_p$-vector space $\Fil^\bullet W$, the corresponding vector bundle on $\Vect((\Q_p^\N)_{|t|=0, |u|=1})$ is given by $\Q_p(\zeta_{p^\infty})\{t\}^\dagger$-module
\begin{equation*}
V=\left(\bigoplus_i t^{-i}(\Fil^i W\tensor_{\Q_p} \Q_p(\zeta_{p^\infty}))\right)\tensor_{\Q_p[t]} \Q_p\{t\}^\dagger
\end{equation*}
with $\Z_p^\times$ acting trivially on $\Fil^\bullet W$, by $\gamma.t=\gamma\cdot t$ on $t$ and via its usual Galois action on $\Q_p(\zeta_{p^\infty})$. Now it is easy to see that the base change of $V$ along $\Q_p(\zeta_{p^\infty})\{t\}^\dagger\rightarrow B_\dR^{\dagger, +}$ coincides $G_{\Q_p}$-equivariantly with 
\begin{equation*}
\Fil^0(W\tensor_{\Q_p} B_\dR^\dagger)=\sum_i \Fil^i W\tensor_{\Q_p} t^{-i} B_\dR^{\dagger, +}\;.
\end{equation*}
Overall, we have thus proved:

\begin{thm}
The functor
\begin{equation*}
\Vect(\Q_p^{\HT, \dagger, +})\rightarrow \Vect(\Q_p^{\HT, \dagger})\rightarrow \Rep_{G_{\Q_p}}(B_\dR^{+, \dagger})
\end{equation*}
is fully faithful and identifies the source with the full subcategory of representations of $G_{\Q_p}$ on finite free $B_\dR^{+, \dagger}$-modules which are generically flat. Moreover, under this identification, the equivalence
\begin{equation*}
\Vect(\Q_p^{\dR, +})\cong \Vect(\Q_p^{\HT, \dagger, +})
\end{equation*}
from \cref{cor:htdr-perfequiv} identifies with the equivalence between filtered $\Q_p$-vector spaces and generically flat representations on finite free $B_\dR^{+, \dagger}$-modules from \cref{lem:htdr-drrepfilteredvsp}.
\end{thm}

\begin{rem}
One can check that the identification from the above theorem matches the operator $\Theta$ with Fontaine's Sen operator for $B_\dR^+$-representations from \cite[§3]{FontaineBdRSen}. More generally, one can show that $B_\dR^{+, \dagger}$-representations are equivalent to vector bundles on $\Q_p^{\HT, \dagger}$, see \TODO{Reference!!}, and then the aforementioned statement holds more generally for vector bundles on $\Q_p^{\HT, \dagger}$.
\end{rem}

\begin{rem}
One can prove a similar result for $\Q_p^{\HT, +}$ in place of $\Q_p^{\HT, \dagger, +}$. Namely, restricting the arguments above to the locus $\{t=0\}$ show that the functor
\begin{equation*}
\Vect(\Q_p^{\HT, +})\rightarrow\Vect(\Q_p^\HT)\rightarrow\Rep_{G_{\Q_p}}(\mathbb{C}_p)
\end{equation*}
is fully faithful and identifies the source with the full subcategory of $\mathbb{C}_p$-representations of $G_{\Q_p}$ which are \emph{Hodge--Tate}, i.e.\ they become trivial after tensoring with $\bigoplus_{i\in\Z} \mathbb{C}_p(i)$. Moreover, the same techniques as in the proof of \cref{lem:htdr-perfut0} show that
\begin{equation*}
\Perf(\Q_p^\Hod)\cong \Perf(\Q_p^{\HT, +})
\end{equation*}
and the identification above matches this equivalence with the equivalence between graded $\Q_p$-vector spaces and Hodge--Tate $\mathbb{C}_p$-representations of $G_{\Q_p}$.
\end{rem}

\subsubsection{The case of general $X$}

We now move on to the case of arbitrary rigid smooth derived Berkovich spaces over $\Q_p$. As is perhaps to be expected, the key case is the one of the $n$-dimensional overconvergent torus $\ol{\T}^n=\GSpec \Q_p\langle \ul{T}^{\pm 1}\rangle_{\leq 1}$, where $\ul{T}=(T_1, \dots, T_n)$. We invite the reader to think about the case $n=1$ first as it is prototypical for the case of general $n$.

Let us first introduce some notation: Let $\widetilde{\cal{G}}_n$ denote the group stack whose underlying Gelfand stack is the product $(\Z_p^\la)^n\times\Z_p^{\times, \la}\times (\G_a^\dagger)^n\times\G_a^\dagger\times\G_m(1)$ with the group operation given by the formula
\begin{equation*}
(\ul{\theta}_1, \gamma_1, \ul{v}_1, w_1, s_1)\cdot (\ul{\theta}_2, \gamma_2, \ul{v}_2, w_2, s_2)=(\ul{\theta}_1+\gamma_1\ul{\theta}_2, \gamma_1\gamma_2, \ul{v}_1+\gamma_1(s_1\ul{v}_2+w_1\ul{\theta}_2), w_1+s_1w_2, s_1s_2)\;,
\end{equation*}
where $\ul{\theta}_i=(\theta_i^{(1)}, \dots, \theta_i^{(n)})$ and similarly for $\ul{v}_i$ and we write
\begin{equation*}
\gamma_1\ul{\theta}_2\coloneqq (\gamma_1\theta_2^{(1)}, \dots, \gamma_1\theta_2^{(n)})
\end{equation*}
and similarly in the other cases. After base change to $\ol{\DD}=\GSpec\Q_p\langle u\rangle_{\leq 1}$, this admits two maps from $(\G_a^\dagger)^n\rtimes\G_m^\dagger$, where $\G_m^\dagger$ acts on $(\G_a^\dagger)^n$ by multiplication: The first one is given by the formula
\begin{equation*}
(\G_a^\dagger)^n\rtimes\G_m^\dagger\ni (\ul{\theta}, s)\mapsto (0, 1, u\ul{\theta}, u(s-1), s)\in \widetilde{\cal{G}}_n
\end{equation*}
while the second one is just given by $(\ul{\theta}, s)\mapsto (\ul{\theta}, s, 0, 0, 1)$. We let $\cal{G}_n$ be the coequaliser of these two maps in the category of group stacks over $\ol{\DD}$.

\begin{prop}
\label{prop:htdr-thtdagger+}
There is an isomorphism
\begin{equation*}
(\ol{\T}^n)^{\HT, \dagger, +}\cong \left(\ol{\T}^{n, \la}_\infty\times \ol{\DD}\times\G_a^\dagger\times\GSpec\Q_p(\zeta_{p^\infty})\right)\,\big/\,\cal{G}_n\;,
\end{equation*}
where $\ol{\T}^{n, \la}_\infty$ denotes the Gelfand stack associated to $\colim_m \Q_p\langle \ul{T}^{\pm 1/p^m}\rangle_{\leq 1}$ and the group actions defining the quotient are given as follows:
\begin{enumerate}[label=(\roman*)]
\item $\cal{G}_n$ acts on $\ol{\T}^{n, \la}_\infty$ via
\begin{equation*}
(\ul{\theta}, \gamma, \ul{v}, w, s).\ul{T}^{1/p^\infty}=\zeta_{p^\infty}^{\ul{\theta}}\exp(t(u\ul{\theta}+\ul{v})/p^\infty)\ul{T}^{1/p^\infty}\;;
\end{equation*}
\item $\cal{G}_n$ acts on $\G_a^\dagger$ via
\begin{equation*}
(\ul{\theta}, \gamma, \ul{v}, w, s).t=\gamma st\;,
\end{equation*}
where we recall that $t$ denotes the coordinate on $\G_a^\dagger$;
\item $\Z_p^{\times, \la}$ acts on $\GSpec\Q_p(\zeta_{p^\infty})$ via the usual Galois action $\gamma.\zeta_{p^\infty}=\zeta_{p^\infty}^\gamma$;
\item $\G_a^\dagger\rtimes\G_m(1)$ acts on $\ol{\DD}$ via
\begin{equation*}
(w, s).u=(u+w)/s\;.
\end{equation*}
\end{enumerate}
In particular, we have
\begin{equation*}
(\ol{\T}^n)^{\HT, \dagger}\cong \left(\ol{\T}^{n, \la}_\infty\times \G_a^\dagger\times\GSpec\Q_p(\zeta_{p^\infty})\right)\,\big/\,((\Z_p^\la)^n\rtimes \Z_p^{\times, \la})\;,
\end{equation*}
where $\Z_p^{\times, \la}$ acts on $(\Z_p^\la)^n$ by multiplication.
\end{prop}
\begin{proof}
By definition, the stack $(\ol{\T}^n)^{\HT, \dagger, +}$ sits in a pullback square
\begin{equation*}
\begin{tikzcd}
(\ol{\T}^n)^{\HT, \dagger, +}\ar[r]\ar[d] & (\ol{\T}^n)^\dR\ar[d] \\
(\ol{\T}^n)^{\Cone}\ar[r] & ((\ol{\T}^n)^{\Cone})^\dR
\end{tikzcd}
\end{equation*}
over $\Q_p^{\HT, \dagger, +}$, which we have already computed to be isomorphic to
\begin{equation*}
\left(\ol{\DD}\times\G_a^\dagger\times\GSpec\Q_p(\zeta_{p^\infty})\right)\,\big/\,\cal{G}_0\;.
\end{equation*}
Note that the quotient we claim to be isomorphic to $(\ol{\T}^n)^{\HT, \dagger, +}$ indeed admits an obvious projection map to the stack above.

Furthermore, observe that $(\ol{\T}^n)^\dR\cong \ol{\T}^n/(\G_a^\dagger)^n$, where $(\G_a^\dagger)^n$ acts on $\ol{\T}^n$ via coordinatewise multiplication by the exponential, and that, over $\Q_p^{\HT, \dagger, +}$, the generalised Cartier divisor $L\rightarrow A$ defining $(\ol{\T}^n)^{\Cone}$ becomes the zero map after $\dagger$-reduction. Thus, the fact that $\ol{\T}^n$ is $\dagger$-formally smooth over $\Q_p$ implies that $\ol{\T}^n$ surjects onto both $(\ol{\T}^n)^{\Cone}$ and $((\ol{\T}^n)^{\Cone})^\dR$ and this lets us compute the following presentations of these stacks: We have
\begin{equation*}
((\ol{\T}^n)^{\Cone})^\dR\cong \left((\ol{\T}^n)^\dR\times */\G_m(1)^\dR\right)\,\big/\,(\G_a^n)^\dR\cong \ol{\T}^n/(\G_a^\dagger)^n\times */((\G_a^n)^\dR\rtimes \G_m(1)^\dR)\;,
\end{equation*}
where the action of $(\G_a^n)^\dR$ is trivial in the first quotient, $(\G_a^\dagger)^n$ acts on $\ol{\T}^n$ via multiplication by the exponential and $\G_m(1)^\dR$ acts on $(\G_a^n)^\dR$ via multiplication to define the semidirect product; and
\begin{equation*}
(\ol{\T}^n)^{\Cone}\cong \left(\ol{\T}^n\times \G_a^\dagger/\G_m(1)\right)\,\big/\,\G_a^n\cong (\ol{\T}^n\times \G_a^\dagger)\,\big/\,(\G_a^n\rtimes \G_m(1))\;,
\end{equation*}
where $\G_a^n$ acts on $\ol{\T}^n$ via $\ul{x}.\ul{T}=\exp(t\ul{x})\ul{T}$ for $t$ being the coordinate on $\G_a^\dagger$ and $\G_m(1)$ acts on $\G_a^n$ via multiplication to define the semidirect product; we point out that this is indeed the ``same $t$'' as in the claimed presentation of $(\ol{\T}^n)^{\HT, \dagger, +}$ and the presentation of $\Q_p^{\HT, \dagger, +}$ above -- indeed, this comes from the proof of \cref{prop:htdr-qphtdagger+}.

Let us now start by giving a map from
\begin{equation}
\label{eq:htdr-tnhtdaggerpluspresentation}
\left(\ol{\T}^{n, \la}_\infty\times \ol{\DD}\times\G_a^\dagger\times\GSpec\Q_p(\zeta_{p^\infty})\right)\,\big/\,\cal{G}_n
\end{equation}
to $(\ol{\T}^n)^\dR$. For this, observe that the coordinates $\ul{T}$ on $\ol{\T}^{n, \la}_\infty$ in the above quotient are well-defined up to multiplication by $\exp(t(u\ul{\theta}+\ul{v}))\in (\G_m^\dagger)^n$ and hence define a map to $(\ol{\T}^n)^\dR$. In other words, the map from (\ref{eq:htdr-tnhtdaggerpluspresentation}) to $(\ol{\T}^n)^\dR\cong \ol{\T}^n/(\G_a^\dagger)^n$ is induced by the projection
\begin{equation*}
(\ol{\T}^n)_{\infty}^\la\times \ol{\DD}\times\G_a^\dagger\times\GSpec\Q_p(\zeta_{p^\infty})\xrightarrow{\mathrm{pr_1}} \ol{\T}^{n, \la}_\infty\rightarrow \ol{\T}^n
\end{equation*}
and the map
\begin{equation*}
\begin{split}
\cal{G}_n&\rightarrow (\G_a^\dagger)^n \\
(\ul{\theta}, \gamma, \ul{v}, w, s)&\mapsto t(u\ul{\theta}+\ul{v})\;;
\end{split}
\end{equation*}
note that, while this might not appear to be a morphism of groups at first sight, this is salvaged by the fact that everything happens over $\Q_p^{\HT, \dagger, +}$ and that the group $\cal{G}_0$ acts nontrivially on the coordinates $t$ and $u$, which appear in the formula above. Similarly, using the above presentation of $(\ol{\T}^n)^{\Cone}$, the map from (\ref{eq:htdr-tnhtdaggerpluspresentation}) to $(\ol{\T}^n)^{\Cone}$ is induced by the maps
\begin{equation*}
\ol{\T}^{n, \la}_\infty\times \ol{\DD}\times\G_a^\dagger\times\GSpec\Q_p(\zeta_{p^\infty})\xrightarrow{\mathrm{pr_{13}}} \ol{\T}^{n, \la}_\infty\times\G_a^\dagger\rightarrow \ol{\T}^n\times\G_a^\dagger
\end{equation*}
and the map
\begin{equation*}
\begin{split}
\cal{G}_n&\rightarrow (\G_a)^n\rtimes\G_m(1) \\
(\ul{\theta}, \gamma, \ul{v}, w, s)&\mapsto (u\ul{\theta}+\ul{v}, s\gamma)\;,
\end{split}
\end{equation*}
where the same warning as previously applies.

Finally, we have to check that these two maps are indeed compatible in the sense that they induce the same map to $((\ol{\T}^n)^{\Cone})^\dR$. For this, we claim that the postcomposition of the map from (\ref{eq:htdr-tnhtdaggerpluspresentation}) to $(\ol{\T}^n)^\dR$ we have defined with $(\ol{\T}^n)^\dR\rightarrow ((\ol{\T}^n)^{\Cone})^\dR$ is induced by the projection onto $\ol{\T}^n$ ``in the numerator'' and the map
\begin{equation}
\label{eq:htdr-mapfromgn}
\begin{split}
\cal{G}_n&\rightarrow (\G_a^\dagger)^n\times ((\G_a^n)^\dR\rtimes\G_m(1)^\dR) \\
(\ul{\theta}, \gamma, \ul{v}, w, s)&\mapsto (t(u\ul{\theta}+\ul{v}), u\ul{\theta}, s\gamma)\;.
\end{split}
\end{equation}
Indeed, this would show compatibility of the two maps constructed above: For this, note that $\ul{v}\in(\G_a^\dagger)^n$ and hence $u\ul{\theta}+\ul{v}$ and $u\ul{\theta}$ become equal in $(\G_a^n)^\dR$ and recall that the generalised Cartier divisor classified by the map to $\G_a^\dagger/\G_m(1)$, over which $(\ol{\T}^n)^{\Cone}$ lives, is exactly given by multiplication by $t$.

To prove the claim from the previous paragraph, note that the composite map 
\begin{equation*}
(\ol{\T}^n)^\dR\rightarrow ((\ol{\T}^n)^{\Cone})^\dR\rightarrow */\G_m(1)^\dR
\end{equation*}
classifies the first Breuil--Kisin twist of the tautological bundle over $*/\G_m(1)^\dR$; this already shows that the last entry of the map (\ref{eq:htdr-mapfromgn}) is correct. The first entry is clear and hence it suffices to check correctness of the middle entry. To this end, recall the definition of the map $(\ol{\T}^n)^\dR\rightarrow ((\ol{\T}^n)^{\Cone})^\dR$ over $\Q_p^{\HT, \dagger, +}$: For a perfectoid $\ol{A}$ over $\Q_p^\cycl$ equipped with a map $\Q_p\langle \ul{T}^{\pm 1}\rangle\rightarrow \ol{A}$, consider the extension
\begin{equation}
\label{eq:htdr-extension}
\begin{tikzcd}
0\ar[r] & \ol{A}\{-1\}\ar[r, "\cdot\xi"] & \O_{\FF_{\ol{A}}}/\xi^2\ar[r] & \ol{A}\ar[r] & 0\nospacepunct{\;,}
\end{tikzcd}
\end{equation}
where $\xi=\log [\epsilon]$ with $\epsilon=(1, \zeta_p, \zeta_{p^2}, \dots)$, as usual and take its pushout $E$ along the multiplication map $u: \ol{A}\{-1\}\rightarrow \ol{L}\{-1\}$, where $\ol{L}$ denotes the (inverse of) the tautological bundle over $*/\G_m(1)^\dR$ and $u$ is the coordinate on $\ol{\DD}$, as always. Then the $\Q_p\langle \ul{T}^{\pm 1}\rangle$-point of $((\ol{\T}^n)^{\Cone})^\dR$ obtained via the map $(\ol{\T}^n)^\dR\rightarrow ((\ol{\T}^n)^{\Cone})^\dR$ classifies the original map $\Q_p\langle \ul{T}^{\pm 1}\rangle\rightarrow\ol{A}$ together with the pullback of $E$ along this map, which is then an extension of $\Q_p\langle \ul{T}^{\pm 1}\rangle$ by $\ol{L}\{-1\}$.

However, note that the extension (\ref{eq:htdr-extension}) splits after pullback to $\Q_p\langle \ul{T}^{\pm 1}\rangle$ once we have a fixed system of $p$-power roots of $\ul{T}$ in $\ol{A}$: Indeed, a splitting is then induced by the map
\begin{equation*}
\begin{split}
\Q_p\langle \ul{T}^{\pm 1}\rangle&\rightarrow \O_{\FF_{\ol{A}}}/\xi^2\cong B_\dR^+(\ol{A})/\xi^2 \\
\ul{T}&\mapsto [\ul{T}^\flat]=[(\ul{T}, \ul{T}^{1/p}, \ul{T}^{1/p^2}, \dots)]\;.
\end{split}
\end{equation*}
As the automorphisms of the $p$-power roots $\ul{T}^{1/p^\infty}$ on the de Rham stack of (\ref{eq:htdr-tnhtdaggerpluspresentation}) are given by $(\Z_p^\la)^n$ via $\ul{\theta}.\ul{T}^{1/p^\infty}=\zeta_{p^\infty}^{\ul{\theta}} \ul{T}^{1/p^\infty}$, this shows that the extension (\ref{eq:htdr-extension}) is classified by $\ul{\theta}\in(\Z_p^\la)^n$ over (\ref{eq:htdr-tnhtdaggerpluspresentation}). Then the correctness of the middle entry in (\ref{eq:htdr-mapfromgn}) follows from the fact that the map from (\ref{eq:htdr-tnhtdaggerpluspresentation}) to $((\ol{\T}^n)^{\Cone})^\dR$ classifies the pushout of (\ref{eq:htdr-extension}) along $u$.

Summarising, our work so far has established that (\ref{eq:htdr-tnhtdaggerpluspresentation}) admits a map to $(\ol{\T}^n)^{\HT, \dagger, +}$ compatibly with the maps to $\Q_p^{\HT, \dagger, +}$. Now working over $\Q_p^{\HT, \dagger, +}$ throughout, note that the pullback diagram defining $(\ol{\T}^n)^{\HT, \dagger, +}$ from above shows that $(\ol{\T}^n)^{\HT, \dagger, +}$ is a $(\G_a^n)^\dR$-torsor over $(\ol{\T}^n)^{\Cone}\times_{\G_a^\dagger/\G_m(1)} \Q_p^{\HT, \dagger, +}$; indeed, observe that
\begin{equation*}
(\ol{\T}^n)^\dR\times \Q_p^{\HT, \dagger, +}\rightarrow ((\ol{\T}^n)^{\Cone})^\dR\times_{*/\G_m(1)^\dR} \Q_p^{\HT, \dagger, +}
\end{equation*}
is a $(\G_a^n)^\dR$-torsor by the presentation we established for $((\ol{\T}^n)^{\Cone})^\dR$ above. However, using the presentation we have obtained for $(\ol{\T}^n)^{\Cone}$, it is also apparent that the map from (\ref{eq:htdr-tnhtdaggerpluspresentation}) to $(\ol{\T}^n)^{\Cone}\times_{\G_a^\dagger/\G_m(1)} \Q_p^{\HT, \dagger, +}$ is a $(\G_a^n)^\dR$-torsor as well. We are thus done by noting that the map from (\ref{eq:htdr-tnhtdaggerpluspresentation}) to $(\ol{\T}^n)^{\HT, \dagger, +}$ we have constructed is a map of $(\G_a^n)^\dR$-torsors.
\end{proof}

\begin{prop}
\label{prop:htdr-complexesthtdagger+}
A perfect complex on $(\ol{\T}^n)^{\HT, \dagger, +}$ amounts to the following data:
\begin{enumerate}[label=(\roman*)]
\item A diagram
\begin{equation*}
\begin{tikzcd}
\dots \ar[r,shift left=.5ex,"u"]
  & \ar[l,shift left=.5ex, "t"] \Fil_{i-1} V \ar[r,shift left=.5ex,"u"] & \ar[l,shift left=.5ex, "t"] \Fil_i V \ar[r,shift left=.5ex,"u"] & \ar[l,shift left=.5ex, "t"] \Fil_{i+1} V\ar[r,shift left=.5ex,"u"] & \ar[l,shift left=.5ex, "t"] \dots
\end{tikzcd}
\end{equation*}
of perfect complexes over $\ol{\T}^{n, \la}_\infty\times\GSpec\Q_p(\zeta_{p^\infty})\times \GSpec\Q_p\{r\}^\dagger$, where $r=ut$, with the property that $t: \Fil_\bullet V\rightarrow \Fil_{\bullet-1} V$ becomes an isomorphism for $\bullet\ll 0$ and $u: \Fil_\bullet V\rightarrow\Fil_{\bullet+1} V$ becomes an isomorphism for $\bullet\gg 0$;
\item $\ol{\T}^{n, \la}_\infty\times\GSpec\Q_p(\zeta_{p^\infty})$-linear maps $D: \Fil_i V\rightarrow \Fil_{i-1} V$ for each $i$ satisfying $Du=uD+1$ and commuting with $t$;
\item a semilinear locally analytic $\Z_p^\times$-action on $\Fil_\bullet V$ commuting with $u$ and $D$, where the action of $\Z_p^\times$ on $\Q_p(\zeta_{p^\infty})$ is the usual Galois action and the one on $t$ is given by $\gamma.t=\gamma\cdot t$;
\item $n$ commuting $\GSpec\Q_p(\zeta_{p^\infty})$-linear maps $\nabla_1, \dots, \nabla_n: \Fil_i V\rightarrow \Fil_{i-1} V$ for each $i$ commuting with $u, t$ and $D$, satisfying
\begin{equation*}
\nabla_k T_\ell^{1/p^\infty}=T_\ell^{1/p^\infty}\nabla_k+\delta_{k\ell}t\cdot\frac{1}{p^\infty}T_\ell^{1/p^\infty}\;,
\end{equation*}
where $\delta_{k\ell}$ denotes the Kronecker symbol and interacting with the $\Z_p^\times$-action according to $\gamma.\nabla_k(m)=\gamma\cdot \nabla_k(\gamma.m)$;
\item a semilinear locally analytic $\Z_p^n$-action on $\Fil_\bullet V$ commuting with $u, t, D$ and the $\nabla_k$ and interacting with the $\Z_p^\times$-action from (iii) in such a way that they assemble into a $\Z_p^n\rtimes\Z_p^\times$-action, where $\Z_p^\times$ acts via multiplication on $\Z_p^n$; here, the semilinearity is asked with respect to the action of $\Z_p^n$ on $\ol{\T}^{n, \la}_\infty\times\GSpec\Q_p\{r\}^\dagger\times\GSpec\Q_p(\zeta_{p^\infty})$ given by 
\begin{equation*}
\ul{\theta}.\ul{T}^{1/p^\infty}=\zeta_{p^\infty}^{\ul{\theta}}\exp(r\ul{\theta}/p^\infty)\ul{T}^{1/p^\infty}
\end{equation*}
and the trivial action on $r$ and $\zeta_{p^\infty}$;
\item an identification between the ``arithmetic Sen operator'' $\Theta^\arithm: \Fil_i V\rightarrow \Fil_i V$ coming from the Lie algebra action of $\Z_p^\times$ and the operator $uD-i$ for each $i$;
\item an identification between the ``geometric Sen operator'' $\Theta_k^\geom: \Fil_i V\rightarrow \Fil_i V$ coming from the Lie algebra action of the $k$-th copy of $\Z_p$ and the operator $u\nabla_k$ for each $i, k$.
\end{enumerate}
Under this identification, restriction to $(\ol{\T}^n)^{\HT, \dagger}$ corresponds to forgetting the $u$-filtration and the $D$-maps while restriction to $((\ol{\T}^n)^\N)_{|u|=|t|=0}$ corresponds to passing to the associated graded of the $u$-filtration and forgetting the $D$-maps. Finally, the isomorphism
\begin{equation*}
((\ol{\T}^n)^\N)_{|u|=|t|=0}\cong ((\ol{\T}^n)^{\dR, +})_{|t|=0}
\end{equation*}
Tate twists the $\Z_p^\times$-action on the $i$-th associated graded piece of the $u$-filtration by $i$, which trivialises the Lie algebra action $\Z_p^\times$ and hence after descent along the map
\begin{equation*}
\ol{\T}^{n, \la}_\infty\times \GSpec\Q_p(\zeta_{p^\infty})\rightarrow \ol{\T}^n
\end{equation*}
via the now smooth $\Z_p^n\rtimes\Z_p^\times$-action (note that the $\Z_p^n$-action was smooth to begin with) yields a filtered perfect complex $\Fil^\bullet W$ on $\ol{\T}^n$ via the $t$-maps together with operators $\nabla_k: \Fil^i W\rightarrow \Fil^{i-1} W$ for $k=1, \dots, n$ satisfying the relation $\nabla_k T_\ell=T_\ell\nabla_k+\delta_{k\ell}tT_\ell$.
\end{prop}
\begin{proof}
Using the explicit presentation for $(\ol{\T}^n)^{\HT, \dagger, +}$ we have obtained in \cref{prop:htdr-thtdagger+}, this is now just a beefed up version of the argument from the proof of \cref{prop:htdr-qphtdagger+}: Item (i) is a version of the Rees equivalence coming from the $\G_m(1)$-factor in $\widetilde{\mathcal{G}}_n$ while the maps $D$ and $\nabla_1, \dots, \nabla_n$ come from the factors $\G_a^\dagger$ and $(\G_a^\dagger)^n$ in $\widetilde{\mathcal{G}}_n$ via Cartier duality. The locally analytic $\Z_p^n$- and $\Z_p^\times$-actions clearly come from the corresponding factors in $\widetilde{\mathcal{G}}_n$ and, finally, the compatibilities in (vi) and (vii) come from the definition of $\mathcal{G}_n$ as a coequaliser (for (vii), one performs a similar calculation as in the proof of \cref{lem:htdr-ddrmodgmdagger}).
\end{proof}

We point out that, as expected, the operators $\nabla_1, \dots, \nabla_n$ on $\Fil^\bullet W$ we obtain after restriction to $((\ol{\T}^n)^\N)_{|u|=|t|=0}\cong ((\ol{\T}^n)^{\dR, +})_{|t|=0}$ correspond to a flat connection the underlying perfect complex $W$ of $\Fil^\bullet W$ satisfying Griffiths transversality with respect to the filtration: Given $\nabla_1, \dots, \nabla_n$ as above, the maps
\begin{equation*}
\begin{split}
W&\rightarrow W\tensor_{\O_{\ol{\T}^n}}\Omega^1_{\ol{\T}^n} \\
m&\mapsto \sum_{k=1}^n \nabla_k(m)\tensor\frac{\mathrm{d}T_k}{T_k}
\end{split}
\end{equation*}
will satisfy the Leibniz rule and flatness is ensured by the fact that the $\nabla_k$ commute with each other. 

By the compatibility of the functor $X\mapsto X^\N$ with rigid étale localisations from \cref{thm:defis-berketalecover}, we thus conclude the following:

\begin{thm}
\label{thm:htdr-xhtdrdagger+}
Let $X$ be a derived Berkovich space over $\Q_p$ equipped with a rigid étale map $X\rightarrow \ol{\T}^n$ for some $n$. Then
\begin{equation*}
X^{\HT, \dagger, +}\cong \left(X_\infty^\la\times \ol{\DD}\times\G_a^\dagger\times\GSpec\Q_p(\zeta_{p^\infty})\right)\,\big/\,\cal{G}_n\;,
\end{equation*}
where $X_\infty^\la\coloneqq X\times_{\ol{\T}^n} (\ol{\T}^n)_\infty^\la$ and, in particular,
\begin{equation*}
X^{\HT, \dagger}\cong \left(X_\infty^\la\times \G_a^\dagger\times\GSpec\Q_p(\zeta_{p^\infty})\right)\,\big/\,((\Z_p^\la)^n\rtimes\Z_p^{\times, \la})\;.
\end{equation*}
Moreover, a perfect complex on $X^{\HT, \dagger, +}$ amounts to the following data:
\begin{enumerate}[label=(\roman*)]
\item A diagram
\begin{equation*}
\begin{tikzcd}
\dots \ar[r,shift left=.5ex,"u"]
  & \ar[l,shift left=.5ex, "t"] \Fil_{i-1} V \ar[r,shift left=.5ex,"u"] & \ar[l,shift left=.5ex, "t"] \Fil_i V \ar[r,shift left=.5ex,"u"] & \ar[l,shift left=.5ex, "t"] \Fil_{i+1} V\ar[r,shift left=.5ex,"u"] & \ar[l,shift left=.5ex, "t"] \dots
\end{tikzcd}
\end{equation*}
of perfect complexes over $X_\infty^\la\times\GSpec\Q_p(\zeta_{p^\infty})\times \GSpec\Q_p\{r\}^\dagger$, where $r=ut$, with the property that $t: \Fil_\bullet V\rightarrow \Fil_{\bullet-1} V$ becomes an isomorphism for $\bullet\ll 0$ and $u: \Fil_\bullet V\rightarrow\Fil_{\bullet+1} V$ becomes an isomorphism for $\bullet\gg 0$;
\item $X_\infty^\la\times\GSpec\Q_p(\zeta_{p^\infty})$-linear maps $D: \Fil_i V\rightarrow \Fil_{i-1} V$ for each $i$ satisfying $Du=uD+1$ and commuting with $t$;
\item a ``$t$-connection'' $\nabla: V\rightarrow V\tensor_{\O_{X_\infty^\la}} \Omega_{X_\infty^\la}^1(-1)$ on $V$, i.e.\ a map satisfying the $t$-deformed Leibniz rule $\nabla(fm)=f\nabla(m)+t\cdot\mathrm{d}f$, which is Griffiths transversal with respect to the filtration $\Fil_\bullet V$ and commutes with $D$; here, the twist by $-1$ indicates that $\Z_p^\times$ acts on $\Omega_{X_\infty^\la}^1$ by division, i.e.\ $\gamma.\omega=\gamma^{-1}\omega$;
\item a semilinear locally analytic $(\Z_p^n\rtimes\Z_p^\times)$-action on $\Fil_\bullet V$ commuting with all the previous data;
\item an identification between the ``arithmetic Sen operator'' $\Theta^\arithm: \Fil_i V\rightarrow \Fil_i V$ coming from the Lie algebra action of $\Z_p^\times$ and the operator $uD-i$ for each $i$;
\item an identification between the ``geometric Sen operator'' $\Theta^\geom: \Fil_i V\rightarrow \Fil_i V\tensor_{\O_{X_\infty^\la}} \Omega_{X_\infty^\la}^1(-1)$ coming from the Lie algebra action of $\Z_p^n$ and the $t$-connection $u\nabla$ for each $i$.
\end{enumerate}
Under this identification, restriction to $X^{\HT, \dagger}$ corresponds to forgetting the $u$-filtration and the $D$-maps while restriction to $(X^\N)_{|u|=|t|=0}$ corresponds to passing to the associated graded of the $u$-filtration and forgetting the $D$-maps. Finally, the isomorphism
\begin{equation*}
(X^\N)_{|u|=|t|=0}\cong (X^{\dR, +})_{|t|=0}
\end{equation*}
Tate twists the $\Z_p^\times$-action on the $i$-th associated graded piece of the $u$-filtration by $i$, which trivialises the Lie algebra action $\Z_p^\times$ and hence after descent along the map
\begin{equation*}
X_\infty^\la\times \GSpec\Q_p(\zeta_{p^\infty})\rightarrow X
\end{equation*}
via the now smooth $\Z_p^n\rtimes\Z_p^\times$-action yields a filtered perfect complex $\Fil^\bullet W$ on $X$ via the $t$-maps together with a flat connection $\nabla: W\rightarrow W\tensor_{\O_X} \Omega^1_X$ on the underlying unfiltered complex satisfying Griffiths transversality.
\end{thm}

\begin{rem}
In principle, one could ``compress'' the data contained in a perfect complex on $X^{\HT, \dagger, +}$ from above even further: Namely, we could replace (ii) and (v) by
\begin{itemize}
\item[(ii')] a factorisation of $\Theta^\arithm+i: \Fil_i V\rightarrow \Fil_i V$ through $\Fil_{i-1} V$, where $\Theta^\arithm$ is the ``arithmetic Sen operator'' coming from the Lie algebra action of $\Z_p^\times$
\end{itemize}
and similarly replace (iii) and (vi) by
\begin{itemize}
\item[(iii')] a factorisation of the ``geometric Sen operator'' $\Theta^\geom: \Fil_i V\rightarrow \Fil_i V\tensor_{\O_{X_\infty^\la}} \Omega_{X_\infty^\la}^1(-1)$ coming from the Lie algebra action of $\Z_p^n$ through $\Fil_{i-1} V\tensor_{\O_{X_\infty^\la}} \Omega_{X_\infty^\la}^1(-1)$.
\end{itemize}
Indeed, by (v) and (vi), these factorisations will precisely be given by $D$ and $\nabla$, respectively. However, we have chosen not to do this because we want to emphasise the crucial role that $D$ and $\nabla$ play, e.g.\ in describing the restriction functors to $X^{\HT, \dagger}$ and $(X^\N)_{|u|=|t|=0}$ or in describing cohomology on $X^{\HT, \dagger, +}$.
\end{rem}

From \cref{thm:htdr-xhtdrdagger+}, we immediately obtain the following description of cohomology on the stack $X^{\HT, \dagger, +}$.

\begin{cor}
\label{cor:htdr-cohomologyxhtdrdagger+}
In the situation of \cref{thm:htdr-xhtdrdagger+} and under the identifications from loc.\ cit., the cohomology of a perfect complex on $X^{\HT, \dagger, +}$ is computed by the (underived) $\Z_p^n\rtimes\Z_p^\times$-invariants of the coherent cohomology on $X_\infty^\la$ of the total complex of
\begin{equation*}
\begin{tikzcd}
\Fil_0 V\ar[r, "\nabla"]\ar[d, "D", swap] & \Fil_{-1} V\tensor_{\O_{X_\infty^\la}} \Omega^1_{X_\infty^\la}(-1)\ar[d, "D", swap]\ar[r, "\nabla"] & \Fil_{-2} V\tensor_{\O_{X_\infty^\la}} \Omega^2_{X_\infty^\la}(-2)\ar[r, "\nabla"]\ar[d, "D", swap] & \dots \\
\Fil_{-1} V\ar[r, "\nabla"] & \Fil_{-2} V\tensor_{\O_{X_\infty^\la}} \Omega^1_{X_\infty^\la}(-1)\ar[r, "\nabla"] & \Fil_{-3} V\tensor_{\O_{X_\infty^\la}} \Omega^2_{X_\infty^\la}(-2)\ar[r, "\nabla"] & \dots\nospacepunct{\;.}
\end{tikzcd}
\end{equation*}
Similarly, the cohomology of the pullback of $E$ to $X^{\HT, \dagger}$ is computed by the (underived) $\Z_p^n\rtimes\Z_p^\times$-invariants of the coherent cohomology on $X_\infty^\la$ of the total complex of
\begin{equation*}
\begin{tikzcd}
V\ar[r, "\nabla"]\ar[d, "D", swap] & V\tensor_{\O_{X_\infty^\la}} \Omega^1_{X_\infty^\la}(-1)\ar[d, "D", swap]\ar[r, "\nabla"] & V\tensor_{\O_{X_\infty^\la}} \Omega^2_{X_\infty^\la}(-2)\ar[r, "\nabla"]\ar[d, "D", swap] & \dots \\
V\ar[r, "\nabla"] & V\tensor_{\O_{X_\infty^\la}} \Omega^1_{X_\infty^\la}(-1)\ar[r, "\nabla"] & V\tensor_{\O_{X_\infty^\la}} \Omega^2_{X_\infty^\la}(-2)\ar[r, "\nabla"] & \dots\nospacepunct{\;.}
\end{tikzcd}
\end{equation*}
\end{cor}
\begin{proof}
The projection map $\widetilde{\cal{G}}_n\rightarrow (\Z_p^n)^\la\rtimes\Z_p^{\times, \la}$ induces a map $\cal{G}_n\rightarrow (\Z_p^n)^\sm\rtimes\Z_p^{\times, \sm}$ and, using the description of $X^{\HT, \dagger, +}$ from \cref{thm:htdr-xhtdrdagger+}, this in turn induces a map 
\begin{equation}
\label{eq:htdr-cohomologyxhtdrdagger+}
X^{\HT, \dagger, +}\rightarrow X_\infty^\la\,\big/\,((\Z_p^n)^\sm\rtimes \Z_p^{\times, \sm})\;.
\end{equation}
From the description of the category of perfect complexes on $X^{\HT, \dagger, +}$ given in loc.\ cit., it is clear that the total complex of the double complex above equipped with its induced $(\Z_p^n)^\sm\rtimes \Z_p^{\times, \sm}$-action computes the derived pushforward along this map. Now note that profinite groups like $\Z_p^n\rtimes\Z_p^\times$ have no higher smooth cohomology in characteristic and hence computing cohomology on the target of (\ref{eq:htdr-cohomologyxhtdrdagger+}) amounts to taking underived $\Z_p^n\rtimes\Z_p^\times$-invariants of the cohomology on $X_\infty^\la$, as claimed.

Similarly, for the second assertion, observe that the derived pushforward along
\begin{equation*}
X^{\HT, \dagger}\rightarrow X_\infty^\la\,\big/\, ((\Z_p^n)^\sm\rtimes \Z_p^{\times, \sm})
\end{equation*}
is given by taking Lie algebra cohomology for the locally analytic action of $\Z_p^n\rtimes\Z_p^\times$ by the presentation of $X^{\HT, \dagger}$ from \cref{thm:htdr-xhtdrdagger+}. Recalling that the arithmetic Sen operator $\Theta^\arithm$ is given by $u\nabla$ on all $\Fil_i V$ while the geometric Sen operator $\Theta^\geom$ is given by $uD$ on $\Fil_i V(i)$, we see that this Lie algebra cohomology is exactly computed by the total complex of the second double complex above. One concludes as in the previous paragraph.
\end{proof}

Finally, let us identify the functor
\begin{equation*}
\Vect(X^{\dR, +})\rightarrow\Vect(X^{\HT, \dagger})\rightarrow\{\text{$\mathbb{B}_\dR^+$-local systems on $X_\proet$}\}
\end{equation*}
we obtain from \cref{cor:htdr-perfequiv} with Scholze's functor from \cite[§7]{PAdicHodgeTheory} in the case where $X$ is a partially proper rigid space equipped with a rigid étale map $X\rightarrow\ol{\T}^n$. For this, first note that, by the same argument as in the case $X=\GSpec\Q_p$, the filtration $\Fil_\bullet V$ associated to a vector bundle on $X^{\HT, \dagger, +}$ always splits $\Z_p^\times$-equivariantly, i.e.\ we have
\begin{equation*}
\Fil_\bullet V=\left(\bigoplus_{i\leq \bullet} t^{-i}(\Fil^i W\tensor_{\Q_p} \Q_p(\zeta_{p^\infty})\tensor_{\O_X} \O_{X_\infty^\la})\right)\tensor_{\Q_p[r]} \Q_p\{r\}^\dagger
\end{equation*}
in the notation of \cref{thm:htdr-xhtdrdagger+}, where $r$ acts via the map $t^{-i}\Fil^i W\rightarrow t^{-i+1}\Fil^{i-1} W$, the $\Z_p^n$-action is induced by the one on $\O_{X_\infty^\la}$ and $\Z_p^\times$ acts on $\Q_p(\zeta_{p^\infty})$ via the usual Galois action and on $t$ by multiplication. Then our claim is that the base change of 
\begin{equation*}
V=\left(\bigoplus_i t^{-i}(\Fil^i W\tensor_{\O_X} \O_{X_\infty^\la}\tensor_{\Q_p} \Q_p(\zeta_{p^\infty}))\right)\tensor_{\Q_p[t]} \Q_p\{t\}^\dagger
\end{equation*}
along $\O_{X_\infty^\la}\tensor_{\Q_p} \Q_p(\zeta_{p^\infty})\{t\}^\dagger\rightarrow \mathbb{B}^+_{\dR, X^\cycl_\infty}$ identifies $\Z_p^n\rtimes\Z_p^\times$-equivariantly with
\begin{equation*}
\Fil^0(W\tensor_{\O_X} \O\mathbb{B}_{\dR, X^\cycl_\infty})^{\nabla=0}\;,
\end{equation*}
where $X_\infty$ is the adic space obtained by base changing $X_\infty^\la$ along $\Spa \Q_p\langle T^{1/p^\infty}\rangle\rightarrow\ol{\T}^n$ and $X_\infty^\cycl$ is the base change of $X_\infty$ to $\Q_p^\cycl$. Clearly, it will suffice to prove that there is a $\Z_p^n$-equivariant isomorphism
\begin{equation}
\label{eq:htdr-isoscholzecomparison}
\Fil^0(W\tensor_{\O_X} \O\mathbb{B}_\dR)^{\nabla=0}\cong \sum_i \Fil^i W\tensor_{\O_X} t^{-i}\mathbb{B}_\dR^+
\end{equation}
of sheaves on the localised site $X_\proet/X_\infty$.

\begin{lem}
\label{lem:htdr-switchfiltrations}
In the above setup, we have
\begin{equation*}
\Fil^0(W\tensor_{\O_X} \O\mathbb{B}_\dR)^{\nabla=0}\cong \left(\sum_i \Fil^i W\tensor_{\O_X} t^{-i} \O\mathbb{B}_\dR^+\right)^{\nabla=0}
\end{equation*}
via the natural map from right to left.
\end{lem}
\begin{proof}
Recall that $\O\mathbb{B}_\dR^+\cong \mathbb{B}_\dR^+[\![X_1, \dots, X_n]\!]$ on $X_\proet/X_\infty$, where $X_i=T_i-[T_i^\flat]\in\O\mathbb{B}_\inf$, and that the filtration on $\O\mathbb{B}_\dR=\O\mathbb{B}_\dR^+[\tfrac{1}{t}]$ is given by
\begin{equation*}
\Fil^i \O\mathbb{B}_\dR=\sum_{j\geq -i} t^{-j}(t, X_1, \dots, X_n)^{i+j}\O\mathbb{B}_\dR^+\;.
\end{equation*}
Thus, we have
\begin{equation*}
\Fil^0(W\tensor_{\O_X} \O\mathbb{B}_\dR)=\sum_i \Fil^i W\tensor_{\O_X} \left(\sum_{j\geq i} t^{-j} (t, X_1, \dots, X_n)^{j-i}\O\mathbb{B}_\dR^+\right)
\end{equation*}
and our task is to show that a horizontal section in 
\begin{equation*}
W\tensor_{\O_X} (t, X_1, \dots, X_n)^\ell\O\mathbb{B}_\dR^+
\end{equation*}
is divisible by $t^\ell$ for any $\ell\geq 1$. (Note that we are indeed free to just check this for $W$ instead of each $\Fil^i W$ since the transition maps of the filtration $\Fil^\bullet W$ are injections.)

Since the assertion of the lemma is local on the site $X_\proet$, we may assume that $W$ is trivial and let $e_1, \dots, e_m$ be a basis. Then the connection on $W$ is given by 
\begin{equation*}
\nabla(e_i)=\sum_{j, k} f_{ijk}e_k\mathrm{d} T_j
\end{equation*}
for some global sections $f_{ijk}$ of $\O_X$. Now letting 
\begin{equation*}
\sum_i e_i\tensor f_i\in W\tensor_{\O_X} (t, X_1, \dots, X_n)^\ell\O\mathbb{B}_\dR^+
\end{equation*}
be a horizontal section, we obtain the differential equations
\begin{equation}
\label{eq:htdr-diffeqsbdr+}
\frac{\partial f_k}{\partial X_j}=-\sum_i f_{ijk}f_i
\end{equation}
for the elements $f_k\in\O\mathbb{B}_\dR^+\cong \mathbb{B}_\dR^+[\![X_1, \dots, X_n]\!]$; here we are using that $\nabla(X_j)=\mathrm{d} T_j$ for the connection $\nabla$ on $\O\mathbb{B}_\dR^+$ and hence 
\begin{equation*}
\nabla(f_k)=\sum_j \frac{\partial f_k}{\partial X_j}\mathrm{d}T_j\;.
\end{equation*}

Working modulo $t^\ell$ from now on, our assumption implies that each $f_k$ is contained in the ideal $(X_1, \dots, X_n)$. However, this means that the right-hand side of (\ref{eq:htdr-diffeqsbdr+}) has vanishing constant term which in turn means that each $f_k$ must be contained in $(X_1, \dots, X_n)^2$. Continuing in this manner, we see that, modulo $t^\ell$, each $f_k$ is contained in all powers of the ideal $(X_1, \dots, X_n)$, which of course implies that $f_k$ vanishes modulo $t^\ell$. In other words, each $f_k$ is divisible by $t^\ell$ and this is what we wanted to show.
\end{proof}

By the above lemma, the assertion (\ref{eq:htdr-isoscholzecomparison}) now turns into the claim that
\begin{equation*}
\left(\sum_i \Fil^i W\tensor_{\O_X} t^{-i}\O\mathbb{B}_\dR^+\right)^{\nabla=0}\cong \sum_i \Fil^i W\tensor_{\O_X} t^{-i}\mathbb{B}_\dR^+\;.
\end{equation*}
Indeed, note that there is a natural map from left to right given by first mapping into $\sum_i \Fil^i W\tensor_{\O_X} t^{-i}\O\mathbb{B}_\dR^+$ and then reducing mod $(X_1, \dots, X_n)$ using the isomorphism $\O\mathbb{B}_\dR^+\cong\B_\dR^+[\![X_1, \dots, X_n]\!]$ on $X_\proet/X_\infty$. However, since the map $\mathbb{B}_\dR^+\rightarrow\O\mathbb{B}_\dR^+$ has a section which is precisely given by the reduction mod $(X_1, \dots, X_n)$, we may alternatively check the above isomorphism after tensoring with $\O\mathbb{B}_\dR^+$, i.e.\ it suffices to show that the natural map
\begin{equation*}
\left(\sum_i \Fil^i W\tensor_{\O_X} t^{-i}\O\mathbb{B}_\dR^+\right)^{\nabla=0}\tensor_{\mathbb{B}_\dR^+}\O\mathbb{B}_\dR^+\rightarrow \sum_i \Fil^i W\tensor_{\O_X} t^{-i}\O\mathbb{B}_\dR^+
\end{equation*}
is an isomorphism. However, this is just an instance of \cite[Thm.\ 7.2]{PAdicHodgeTheory}. Overall, we have proved:

\begin{prop}
\label{prop:htdr-compatibilityscholze}
Let $X$ be a partially proper rigid space equipped with a rigid étale map $X\rightarrow\ol{\T}^n$. Then the functor
\begin{equation*}
\Vect(X^{\dR, +})\cong \Vect(X^{\HT, \dagger, +})\rightarrow\Vect(X^{\HT, \dagger})\rightarrow\{\text{$\mathbb{B}_\dR^+$-local systems on $X_\proet$}\}
\end{equation*}
obtained from \cref{cor:htdr-perfequiv} agrees with the functor
\begin{equation*}
\{\text{filtered vector bundles with connection on $X$}\}\rightarrow \{\text{$\mathbb{B}_\dR^+$-local systems on $X_\proet$}\}
\end{equation*}
from \cite[§7]{PAdicHodgeTheory}.
\end{prop}

\begin{rem}
\label{rem:htdr-functortobdr+dagger}
Let us note that the proof we have given even shows the stronger claim that the functor
\begin{equation*}
\Vect(X^{\dR, +})\cong \Vect(X^{\HT, \dagger, +})\rightarrow\Vect(X^{\HT, \dagger})\rightarrow\{\text{$\mathbb{B}_\dR^{+, \dagger}$-local systems on $X_\proet$}\}
\end{equation*}
agrees with the functor defined by
\begin{equation*}
(\Fil^\bullet W, \nabla)\mapsto \Fil^0(W\tensor_{\O_X} \O\mathbb{B}_\dR^\dagger)^{\nabla=0}\;,
\end{equation*}
where $\O\mathbb{B}_\dR^{+, \dagger}$ is the structural overconvergent de Rham period sheaf from \cite[§3.5]{Wiersig1}. Indeed, in our situation, we have $\O\mathbb{B}_\dR^{+, \dagger}\cong \mathbb{B}_\dR^{+, \dagger}\{X_1, \dots, X_n\}^\dagger$ with $X_i=T_i-[T_i^\flat]$ by \cite[Thm.\ 3.5.7]{Wiersig1} and the use of \cite[Thm.\ 7.2]{PAdicHodgeTheory} can be replaced by the analogous fact that, for any $\Q_p$-algebra $R$, a finite projective $R\{X_1, \dots, X_n\}^\dagger$-module with a flat connection has enough horizontal sections, see e.g.\ the proof of \cite[Thm.\ 9.6.1]{KedlayaDiffEqs}.
\end{rem}

\begin{rem}
In fact, we expect the above result to be true for all smooth partially proper rigid spaces $X$ over $\Q_p$. However, verifying that the identification above glues seems to be somewhat tedious and is actually not needed for our purposes.
\end{rem}

\subsubsection{Characterising the essential image of $\Vect(X^{\HT, \dagger, +})\rightarrow \Vect(X^{\HT, \dagger})$}

From the above analysis, it is not too hard to deduce that the functor 
\begin{equation*}
\Vect(X^{\HT, \dagger, +})\rightarrow \Vect(X^{\HT, \dagger})
\end{equation*}
is fully faithful whenever $X$ is a rigid smooth derived Berkovich space over $\Q_p$. Indeed, by compatibility of $X\mapsto X^\N$ with rigid étale localisation, we may reduce to the case where $X$ admits a rigid étale map $X\rightarrow \ol{\T}^n$ and then the same argument as in the case $X=\GSpec\Q_p$ shows that, in the notation of \cref{thm:htdr-xhtdrdagger+}, there are no $\Z_p^\times$-equivariant maps between $\gr_i V$ and $\gr_j V$ for $i\neq j$. Together with the fact that the filtration $\Fil_\bullet V$ always splits, which we have already observed above, this shows that any map $E_1|_{X^{\HT, \dagger}}\rightarrow E_2|_{X^{\HT, \dagger}}$ for $E_1, E_2\in\Vect(X^{\HT, \dagger, +})$ is automatically compatible with the filtrations, whence the claim.

Our remaining goal in this section is to characterise the essential image of this fully faithful embedding in the case where $X$ is a smooth partially proper rigid space over $\Q_p$. For this, observe that \cref{prop:htdr-compatibilityscholze} implies that the $\mathbb{B}_\dR^+$-local system on $X_\proet$ associated to any vector bundle on $X^{\HT, \dagger}$ which comes from a vector bundle on $X^{\HT, \dagger, +}$ will have the following property:

\begin{defi}
\label{defi:htdr-genericallyflat}
Let $X$ be a smooth rigid space over $\Q_p$. A $\mathbb{B}_\dR^+$-local system on $X_\proet$ is called \emph{generically flat} if it lies in the essential image of Scholze's functor
\begin{equation*}
\{\text{filtered vector bundles with connection on $X$}\}\rightarrow\{\text{$\mathbb{B}_\dR^+$-local systems on $X_\proet$}\}\;.
\end{equation*}
\end{defi}

\begin{rem}
Note that this recovers the notion of generically flat $B_\dR^+$-representations of $G_{\Q_p}$ of \cite[Def.\ 10.4.1]{FarguesFontaine} in the case $X=\GSpec\Q_p$.
\end{rem}

\begin{rem}
Note that generic flatness of a $\mathbb{B}_\dR^+$-local system may be checked Berkovich étale locally on $X$. Indeed, this is because Scholze's functor is fully faithful and filtered vector bundles with connection satisfy Berkovich étale descent.
\end{rem}

In other words, we already know that the functor $\Vect(X^{\HT, \dagger, +})\rightarrow\Vect(X^{\HT, \dagger})$ always factors through the following full subcategory:

\begin{defi}
Let $X$ be a smooth partially proper rigid space over $\Q_p$. A vector bundle on $X^{\HT, \dagger}$ is called \emph{generically flat} if its associated $\mathbb{B}_\dR^+$-local system on $X_\proet$ is generically flat.
\end{defi}

\begin{rem}
As $X\mapsto X^\N$ is compatible with rigid étale localisation and generic flatness of $\mathbb{B}_\dR^+$-local systems may be checked étale locally, generic flatness of a vector bundle on $X^{\HT, \dagger}$ may be checked rigid étale locally on $X$.
\end{rem}

In fact, we are going to prove that the above condition exactly singles out the essential image of the functor $\Vect(X^{\HT, \dagger, +})\rightarrow\Vect(X^{\HT, \dagger})$.

\begin{thm}
\label{thm:htdr-vectxhtdagger}
Let $X$ be a smooth partially proper rigid space over $\Q_p$. Then the functor
\begin{equation*}
\Vect(X^{\HT, \dagger, +})\rightarrow\Vect(X^{\HT, \dagger})
\end{equation*}
is fully faithful and its essential image is given by the full subcategory of vector bundles on $X^{\HT, \dagger}$ which are generically flat.
\end{thm}

Before we begin the proof, let us shortly sketch our strategy: First note that the assertion may be checked rigid étale locally on $X$ and thus we may assume that $X$ is equipped with a rigid étale map $X\rightarrow\ol{\T}^n$. In that case, by \cref{prop:htdr-compatibilityscholze}, it suffices to prove that the functor 
\begin{equation*}
\{\text{generically flat vector bundles on $X^{\HT, \dagger}$}\}\rightarrow \{\text{$\mathbb{B}_\dR^+$-local systems on $X_\proet$}\}
\end{equation*}
is fully faithful. By dualisability, this reduces to a cohomology computation, which we intend to reduce to a question between vector bundles on $X^\HT$ and $\widehat{\O}_X$-local systems on $X_\proet$ via the $t$-adic filtration. However, in that case, using the explicit description of cohomology on $X^{\HT, \dagger}$ from \cref{cor:htdr-cohomologyxhtdrdagger+}, the remaining claim will follow from geometric Sen theory. The most serious obstacle to this approach is that, in the notation of \cref{thm:htdr-xhtdrdagger+}, the module $V$ is not derived $t$-complete. However, the situation is remedied by the following lemma.

\begin{lem}
\label{lem:htdr-tiVnocoh}
Assume that $X$ is a derived Berkovich space over $\Q_p$ equipped with a rigid étale map $X\rightarrow\ol{\T}^n$. Let $E$ be a generically flat vector bundle on $X^{\HT, \dagger}$, which by \cref{thm:htdr-xhtdrdagger+} corresponds to a vector bundle $V$ on $X_\infty^\la\times\GSpec \Q_p(\zeta_{p^\infty})\times\G_a^\dagger$ equipped with a semilinear locally analytic $\Z_p^n\rtimes\Z_p^\times$-action. If $\Theta^\arithm: V\rightarrow V$ denotes the ``arithmetic Sen operator'' on $V$ obtained from the Lie algebra action of $\Z_p^\times$, then
\begin{equation*}
\Theta^\arithm: t^i V\rightarrow t^i V
\end{equation*}
is an isomorphism for $i\gg 0$.
\end{lem}
\begin{proof}
By generic flatness of $E$, there is a filtered vector bundle with connection $(\Fil^\bullet W, \nabla)$ on $X$ such that Scholze's functor sends $(\Fil^\bullet W, \nabla)$ to the associated $\mathbb{B}_\dR^+$-local system of $E$. We will show that $\Theta^\arithm$ is an isomorphism on $t^i V$ for any $i$ such that $\Fil^i W=0$, which implies the claim.

To prove this, we may work rigid étale locally on $X$ and, in particular, assume that all $\Fil^\bullet W$ are free $\O_X$-modules and that the transition maps of the filtration are given by the standard embeddings via the first couple of coordinates. Moreover, by \cite[Lem.\ 7.3]{PAdicHodgeTheory}, we may assume that the pullback of $V$ to $X_\infty^\la\times\GSpec\Q_p(\zeta_{p^\infty})$ is free after possibly further localising on $X$ and then we may even assume that $V$ itself is free (after possibly further localising) since freeness of a finite projective module may be checked after $\dagger$-reduction. In the following, let us introduce the notation $\mathbb{B}_\dR^{+, \dagger, \la}$ for the pushforward of the structure sheaf along
\begin{equation*}
X_\infty^\la\times\GSpec\Q_p(\zeta_{p^\infty})\times\G_a^\dagger\xrightarrow{\mathrm{pr}_1} X_\infty^\la\rightarrow X\;;
\end{equation*}
note that this map is affine in the sense that the source becomes affine after pulling back to any affine $\GSpec A\rightarrow X$. 

Now note that the semilinearity of the $\Z_p^\times$-action implies that $\Theta^\arithm$ will be a $t$-connection in the following sense: for all local sections $f$ of $\mathbb{B}_\dR^{+, \dagger, \la}$ and $v$ of $V$, we have
\begin{equation*}
\Theta^\arithm(fv)=f\Theta^\arithm(v)+t\frac{\partial f}{\partial t}v\;,
\end{equation*}
where $t$ denotes the coordinate on $\G_a^\dagger$. Thus, after choosing a basis of $V$, we have
\begin{equation*}
\Theta^\arithm=A+t\frac{\partial}{\partial t}
\end{equation*}
for some $m\times m$-matrix $A$ of global sections of $\mathbb{B}_\dR^{+, \dagger, \la}$, where $m$ is the rank of $V$, and our claim is that the operator
\begin{equation}
\label{eq:htdr-aplusiinvertible}
(A+i)+t\frac{\partial}{\partial t}: (\mathbb{B}_\dR^{+, \dagger, \la})^{\oplus m}\rightarrow (\mathbb{B}_\dR^{+, \dagger, \la})^{\oplus m}
\end{equation}
is invertible for all $i$ as in the first paragraph.

To this end, writing 
\begin{equation*}
A=A_0+A_1t+A_2t^2+\dots
\end{equation*}
for $m\times m$-matrices $A_\ell$ of global functions on $X_\infty^\la\times\GSpec\Q_p(\zeta_{p^\infty})$, our first claim is the following:

\bigskip

\textbf{Claim 1.} For all $i$ as in the first paragraph, the matrix $A_0+i$ is invertible.

\bigskip

\textit{Proof of the claim.} Recall that we have already computed that the image of $(\Fil^\bullet W, \nabla)$ under Scholze's functor coincides with the base change of the $\Q_p[t]$-module
\begin{equation*}
\bigoplus_\bullet t^{-\bullet}\Fil^\bullet W\;,
\end{equation*}
where $t$ acts via $t^{-\bullet}\Fil^\bullet W\rightarrow t^{-\bullet+1} \Fil^{\bullet-1} W$, along $\Q_p[t]\rightarrow\mathbb{B}_\dR^+$; here, $\Z_p^\times$ acts by multiplication on $t$ and trivially on $\Fil^\bullet W$. As we have assumed that all $\Fil^\bullet W$ are free and that the transition maps are the standard inclusions, this means that the image of $(\Fil^\bullet W, \nabla)$ is isomorphic to $(\mathbb{B}_\dR^+)^{\oplus m}$ with $\gamma\in\Z_p^\times$ acting on the standard basis via the diagonal matrix $\mathrm{diag}(\gamma^{-a_1}, \dots, \gamma^{-a_m})$, where the $a_k$ are the indices of the jumps in the filtration $\Fil^\bullet W$ counted with multiplicity; in particular, we have $a_k<i$.

By assumption, there is an invertible $m\times m$-matrix $S$ of sections of $\mathbb{B}_\dR^+$ over $X_\infty^\cycl$ such that the map
\begin{equation*}
V\tensor_{\mathbb{B}_\dR^{+, \dagger, \la}}\mathbb{B}^+_{\dR, X^\cycl_\infty}\cong (\mathbb{B}^+_{\dR, X^\cycl_\infty})^{\oplus m}\overset{S}{\longrightarrow} (\mathbb{B}^+_{\dR, X^\cycl_\infty})^{\oplus m}
\end{equation*}
is $\Z_p^n\rtimes\Z_p^\times$-equivariant, where the right-hand side denotes the image of $(\Fil^\bullet W, \nabla)$ under Scholze's functor while the isomorphism on the left-hand side uses our chosen basis of $V$. Writing $S=(s_{ij})_{i, j}$ and $S^{-1}=(s'_{ij})_{i, j}$, this means that the action of $\gamma\in\Z_p^\times$ on the $k$-th basis vector $e_k$ of $V$ is given by
\begin{equation*}
\begin{split}
\gamma.e_j&=S^{-1}\begin{pmatrix} \gamma^{-a_1} & & \\
& \ddots & \\
& & \gamma^{-a_m}\end{pmatrix} Se_j=S^{-1}\begin{pmatrix} \gamma^{-a_1} & & \\
& \ddots & \\
& & \gamma^{-a_m}\end{pmatrix} \begin{pmatrix} s_{1j} \\ \vdots \\ s_{mj}\end{pmatrix} \\
&=S^{-1}\begin{pmatrix} \gamma^{-a_1}\cdot (\gamma.s_{1j}) \\ \vdots \\ \gamma^{-a_m}\cdot (\gamma.s_{mj})\end{pmatrix}=\left(\sum_k s'_{ik}\gamma^{-a_k}(\gamma.s_{kj})\right)_i\;.
\end{split}
\end{equation*}
Since the action of $\Z_p^\times$ on $V$ is locally analytic, all entries of the last vector must be locally analytic functions in $\gamma$ and, in particular, the functions $\gamma\mapsto \gamma.s_{kj}$ must be locally analytic; taking the derivative at $\gamma=1$ then yields the action of $\Theta^\arithm$ and so
\begin{equation}
\label{eq:htdr-thetaarithmej}
\Theta^\arithm(e_j)=\left(\sum_k (-a_k)s'_{ik}s_{kj}+s'_{ik}\left.\frac{\mathrm{d}(\gamma.s_{kj})}{\mathrm{d}\gamma}\right\vert_{\gamma=1}\right)_i\;,
\end{equation}
i.e.\ the right-hand side is the $j$-th column of the matrix $A$.

Now note that $\mathbb{B}_\dR^+(X_\infty^\cycl)/t\cong \O(X_\infty^\cycl)$ and recall that there is a canonical map $\O(X_\infty^\la)\rightarrow\mathbb{B}_\dR^+(X_\infty^\cycl)$. Thus, after possibly passing to a rational localisation of $X$ so that $X$ becomes affine (note that the map $X\rightarrow\ol{\T}^n$ will then not necessarily be rigid étale anymore, but we will also not need this for the remainder of the proof of the current claim), we may approximate $s_{kj}$ as
\begin{equation*}
s_{kj}=r_{kj, N}+p^N(\,\dots)+t(\,\dots)
\end{equation*}
for any $N\geq 0$ and some $r_{kj, N}\in \O(X_\infty^\la)$. As the action of $\Z_p^\times$ on $X_\infty^\la$ is smooth and $\gamma.t=\gamma\cdot t$, this yields
\begin{equation*}
\left.\frac{\mathrm{d}(\gamma.s_{kj})}{\mathrm{d}\gamma}\right\vert_{\gamma=1}=p^N(\,\dots)+t(\,\dots)\;.
\end{equation*}
In other words, we may conclude from (\ref{eq:htdr-thetaarithmej}) that the $j$-th column of $A_0$ coincides with $(-\sum_k a_k s'_{ik}s_{kj})_i$ up to a multiple of $p^N$ and since the ring $\O(X_\infty^\la\times\GSpec\Q_p(\zeta_{p^\infty}))$ is $p$-adically separated, this yields
\begin{equation*}
A_0=\left(-\sum_k a_k s'_{ik}s_{kj}\right)_{i, j}=S^{-1}\begin{pmatrix} -a_1 & & \\ & \ddots & \\ & & -a_m\end{pmatrix} S\;,
\end{equation*}
hence
\begin{equation*}
A_0+i=S^{-1}\begin{pmatrix} -a_1+i & & \\ & \ddots & \\ & & -a_m+i\end{pmatrix} S\;.
\end{equation*}

Now we are done: Indeed, the above yields
\begin{equation*}
\det(A_0+i)=\prod_{k=1}^m (i-a_k)
\end{equation*}
and, by assumption, we have $a_k<i$ for all $k$. We conclude that $\det(A_0+i)$ is a unit in the $\Q$-algebra $\O(X_\infty^\la\times\GSpec\Q_p(\zeta_{p^\infty}))$ and thus the claim follows. \hfill\qed

\bigskip

Recall from (\ref{eq:htdr-aplusiinvertible}) that we need to show that the operator
\begin{equation*}
\Theta^\arithm+i=(A+i)+t\frac{\partial}{\partial t}: (\mathbb{B}_\dR^{+, \dagger, \la})^{\oplus m}\rightarrow (\mathbb{B}_\dR^{+, \dagger, \la})^{\oplus m}
\end{equation*}
is invertible. We start with injectivity and, to this end, take some
\begin{equation*}
x=\sum_\ell x_\ell t^\ell\in (\mathbb{B}_\dR^{+, \dagger, \la})^{\oplus m}\;,
\end{equation*}
where $x_\ell\in \O(X_\infty^\la\times\GSpec\Q_p(\zeta_{p^\infty}))^{\oplus m}$. Let $\ell$ be minimal such that $x_\ell\neq 0$. Then
\begin{equation*}
(\Theta^\arithm+i)(x)=(A_0+i)x_\ell t^\ell+\ell x_\ell t^\ell+t^{\ell+1}(\,\dots)=(A_0+(i+\ell))x_\ell t^\ell+t^{\ell+1}(\,\dots)
\end{equation*}
and by applying Claim 1 to $i+\ell$ in place of $i$, we see that $(A_0+(i+\ell))x_\ell$ and consequently $(\Theta^\arithm+i)(x)$ is nonzero.

For surjectivity, let $y=\sum_\ell y_\ell t^\ell\in (\mathbb{B}_\dR^{+, \dagger, \la})^{\oplus m}$. By Claim 1, we can then recursively solve the equations
\begin{equation}
\label{eq:htdr-recursiveeqxell}
(A_0+i+\ell)x_\ell+\sum_{k<\ell} (A_{\ell-k}+i)x_k=y_\ell
\end{equation}
for $\ell=0, 1, 2, \dots$ and we claim that then
\begin{equation*}
x\coloneqq \sum_\ell x_\ell t^\ell
\end{equation*}
is in $(\mathbb{B}_\dR^{+, \dagger, \la})^{\oplus m}$ and satisfies $(\Theta^\arithm+i)(x)=y$. Indeed, once we know the former, the latter is immediate since the coefficient of $t^\ell$ in the expansion of $(\Theta^\arithm+i)(x)$ is exactly given by the left-hand side of (\ref{eq:htdr-recursiveeqxell}). To establish that $x\in (\mathbb{B}_\dR^{+, \dagger, \la})^{\oplus m}$, we need to understand the growth of the $x_\ell$; we begin by proving the following:

\bigskip

\textbf{Claim 2.} For all $\ell\geq 0$, we have
\begin{equation*}
x_\ell\coloneqq \sum_{k\leq \ell} \left(\left(\sum_{\lambda\vdash (\ell-k)} (-1)^{\ell(\lambda)}\prod_j (A_{\lambda_j}+i)(A_0+i+k+|\lambda|_{\leq j})^{-1}\right)\cdot (A_0+i+k)^{-1}y_k\right)\;,
\end{equation*}
where the inner sum runs over all ordered (!) partitions $\lambda=(\lambda_1, \lambda_2, \dots, \lambda_s)$ of $\ell-k$ into positive parts, $\ell(\lambda)\coloneqq s$ denotes the length of such a partition and we set
\begin{equation*}
|\lambda|_{\leq j}\coloneqq \lambda_1+\dots+\lambda_j
\end{equation*}
for any $1\leq j\leq s$.

\bigskip

\textit{Proof of the claim.} We use induction on $\ell$, in the base case $\ell=0$ the claimed formula becomes $x_0=(A_0+i)^{-1}y_0$, which is true (note that $0$ admits the empty partition!). For the induction step, note that (\ref{eq:htdr-recursiveeqxell}) implies that
\begin{equation*}
x_\ell=(A_0+i+\ell)^{-1}\left(y_\ell-\sum_{k<\ell} (A_{\ell-k}+i)x_k\right)\;.
\end{equation*}
Plugging in the formulas for the $x_k$ with $k<\ell$, we see that the right-hand side is a linear combination of the $y_k$ with $k\leq \ell$. For $k=\ell$, the coefficient of $y_\ell$ is given by $(A_0+i+\ell)^{-1}$, which is in keeping with the claimed formula (again because the empty partition of $0$ is allowed), while the coefficient of $y_k$ for $k<\ell$ is given by
\begin{equation*}
-(A_0+i+\ell)^{-1}\sum_{k\leq k'<\ell} (A_{\ell-k'}+i)\left(\left(\sum_{\lambda\vdash (k'-k)} (-1)^{\ell(\lambda)}\prod_j (A_{\lambda_j}+i)(A_0+i+k+|\lambda|_{\leq j})^{-1}\right)\cdot (A_0+i+k)^{-1}\right)\;.
\end{equation*}
Noting that there is a bijection
\begin{equation*}
\{\text{pairs $(k', \lambda)$ of $k\leq k'<\ell$ and $\lambda\vdash (k'-k)$}\}\xleftrightarrow{1:1} \{\mu\vdash (\ell-k)\}
\end{equation*}
given by appending $\ell-k'$ to the ordered partition $\lambda$ of $k'-k$ into positive parts to obtain an ordered partition $\mu$ of $\ell-k$ into positive parts, we see that the above expression simplifies to 
\begin{equation*}
\left(\sum_{\mu\vdash (\ell-k)} (-1)^{\ell(\mu)}\prod_j (A_{\mu_j}+i)(A_0+i+k+|\mu|_{\leq j})^{-1}\right)(A_0+i+k)^{-1},
\end{equation*}
which is exactly the coefficient of $y_k$ in the claimed formula. \hfill\qed

\bigskip

In order to show $x\in (\mathbb{B}_\dR^{+, \dagger, \la})^{\oplus m}$, we have to argue that there is some $N>0$ such that $p^{\ell N}x_\ell$ forms a nullsequence in $\O(X_\infty^\la\times \GSpec\Q_p(\zeta_{p^\infty})$. Since $A$ is a matrix over $\mathbb{B}_\dR^{+, \dagger, \la}$ and $y\in (\mathbb{B}_\dR^{+, \dagger, \la})^{\oplus m}$, we know that there is some $M>0$ with $p^{\ell M}A_\ell\rightarrow 0$ and $p^{\ell M}y_\ell\rightarrow 0$; as $\O(X_\infty^\la\times\GSpec\Q_p(\zeta_{p^\infty}))$ is a bounded ring, we may even assume that, for $\ell\geq 1$, the sequence $(p^{\ell M}A_\ell)_\ell$ is a nullsequence in the subring of powerbounded elements. Observing that
\begin{equation*}
p^{\ell M}x_\ell=\sum_{k\leq \ell} \left(\left(\sum_{\lambda\vdash (\ell-k)} (-1)^{\ell(\lambda)}\prod_j p^{\lambda_j M}(A_{\lambda_j}+i)(A_0+i+k+|\lambda|_{\leq j})^{-1}\right)\cdot (A_0+i+k)^{-1}p^{kM} y_k\right)
\end{equation*}
by Claim 2, it thus only remains to control the inverses occurring in the formula.

However, note that $\det(A_0+i)(A_0+i)^{-1}$ is a matrix whose entries are degree $m$ integer polynomials in the entries of $A_0$ by Cramer's rule. Further observing that each summand in the formula above only contains at most $\ell+1$ such inverses as factors, we are thus reduced to controlling the determinants $\det(A_0+i)$. Recalling that these determinants are integers by the proof of Claim 1, we have to show that the $p$-adic valuation of the integers
\begin{equation*}
\det(A_0+i+k)\prod_j \det(A_0+i+k+|\lambda|_{\leq j})
\end{equation*}
is uniformly linearly bounded in $\ell$, where $k\leq\ell$ and $\lambda\vdash(\ell-k)$. Recall, though, that the proof of Claim 1 gives an explicit expression for these determinants; namely, the above product is equal to
\begin{equation*}
\prod_{r=1}^m \left((i+k-a_r)\prod_j (i+k+|\lambda|_{\leq j}-a_r)\right)\;.
\end{equation*}
Observing that $1\leq |\lambda|_{\leq 1}<|\lambda|_{\leq 2}<\dots\leq \ell-k$, we see that the above is a divisor of
\begin{equation*}
\prod_{r=1}^m (i+\ell-a_r)!
\end{equation*}
upon recalling that $i>a_r$ for all $r$. As the $p$-adic valuation of this latter product is bounded by $m\cdot (i+\ell-a)/(p-1)$, where $a$ denotes the minimum of the $a_r$, we win.
\end{proof}

We can now move on to proving \cref{thm:htdr-vectxhtdagger}.

\begin{proof}[Proof of \cref{thm:htdr-vectxhtdagger}]
By compatibility of $X\mapsto X^\N$ with rigid étale localisation and since generic flatness may be checked rigid étale locally, we may assume that $X$ admits a rigid étale map $X\rightarrow\ol{\T}^n$. If we can show that the functor
\begin{equation*}
f^*: \{\text{generically flat vector bundles on $X^{\HT, \dagger}$}\}\rightarrow \{\text{$\mathbb{B}_\dR^+$-local systems on $X_\proet$}\}
\end{equation*}
is fully faithful, then the categories of generically flat vector bundles on $X^{\HT, \dagger}$ and generically flat $\mathbb{B}_\dR^+$-local systems on $X_\proet$ are equivalent. As Scholze's functor yields an equivalence between filtered vector bundles with connection on $X$ and generically flat $\mathbb{B}_\dR^+$-local systems on $X_\proet$, this would imply the claim using \cref{prop:htdr-compatibilityscholze}.

To prove that $f^*$ is fully faithful, first note that it is compatible with taking duals. Thus, checking that $\Hom(F, E)\cong \Hom(f^*F, f^*E)$ for $E, F\in\Vect(X^{\HT, \dagger})$ generically flat reduces to the case $F=\O$ by replacing $E$ by $E\tensor F^\vee$. In other words, it suffices to show that, for any generically flat vector bundle $E$ on $X^{\HT, \dagger}$, we have
\begin{equation}
\label{eq:htdr-cohxhtdaggerbdr+iso}
R\Gamma(X^{\HT, \dagger}, E)\cong R\Gamma(X_\proet, f^*E)
\end{equation}
via $f^*$. Using the notation from \cref{thm:htdr-xhtdrdagger+}, the left-hand side is given by the (underived) $\Z_p^n\rtimes\Z_p^\times$-invariants of the coherent cohomology on $X_\infty^\la$ of the total complex of
\begin{equation*}
\begin{tikzcd}
V\ar[r, "\nabla"]\ar[d, "D", swap] & V\tensor_{\O_{X_\infty^\la}} \Omega^1_{X_\infty^\la}(-1)\ar[d, "D", swap]\ar[r, "\nabla"] & V\tensor_{\O_{X_\infty^\la}} \Omega^2_{X_\infty^\la}(-2)\ar[r, "\nabla"]\ar[d, "D", swap] & \dots \\
V\ar[r, "\nabla"] & V\tensor_{\O_{X_\infty^\la}} \Omega^1_{X_\infty^\la}(-1)\ar[r, "\nabla"] & V\tensor_{\O_{X_\infty^\la}} \Omega^2_{X_\infty^\la}(-2)\ar[r, "\nabla"] & \dots\nospacepunct{\;,}
\end{tikzcd}
\end{equation*}
where $D$ is just given by $\Theta^\arithm$, and then \cref{lem:htdr-tiVnocoh} shows that the $t$-adic filtration on this total complex is complete. Moreover, since $\mathbb{B}_\dR^+$ and hence $f^*E$ is $t$-adically complete, also the filtration
\begin{equation*}
\dots\rightarrow R\Gamma(X_\proet, t^2 f^*E)\rightarrow R\Gamma(X_\proet, t f^*E)\rightarrow R\Gamma(X_\proet, f^*E)
\end{equation*}
is complete.

Thus, we may check the isomorphism (\ref{eq:htdr-cohxhtdaggerbdr+iso}) on associated graded pieces. Each graded piece of $t^\bullet V$ will be a vector bundle $V'$ on $X_\infty^\la\times\GSpec\Q_p(\zeta_{p^\infty})$ equipped with a semilinear locally analytic $\Z_p^n\rtimes\Z_p^\times$-action and the corresponding graded piece of $t^\bullet f^*E$ will be the pullback of $V'$ to a locally free $\widehat{\O}_X$-module on $X_\proet$. Our claim then amounts to showing that the (underived) $\Z_p^n\rtimes\Z_p^\times$-invariants of the coherent cohomology on $X_\infty^\la$ of the total complex of 
\begin{equation*}
\begin{tikzcd}
V'\ar[r, "\nabla"]\ar[d, "D", swap] & V'\tensor_{\O_{X_\infty^\la}} \Omega^1_{X_\infty^\la}(-1)\ar[d, "D", swap]\ar[r, "\nabla"] & V'\tensor_{\O_{X_\infty^\la}} \Omega^2_{X_\infty^\la}(-2)\ar[r, "\nabla"]\ar[d, "D", swap] & \dots \\
V'\ar[r, "\nabla"] & V'\tensor_{\O_{X_\infty^\la}} \Omega^1_{X_\infty^\la}(-1)\ar[r, "\nabla"] & V'\tensor_{\O_{X_\infty^\la}} \Omega^2_{X_\infty^\la}(-2)\ar[r, "\nabla"] & \dots\nospacepunct{\;,}
\end{tikzcd}
\end{equation*}
where $\nabla$ and $D$ arise from the Lie algebra actions of $\Z_p^n$ and $\Z_p^\times$, respectively, compute $R\Gamma(X_\proet, f^*V')$. However, this is now an instance of geometric Sen theory, see \cite{GeometricSenTheory}: taking locally analytic vectors for the $\Z_p^n\rtimes\Z_p^\times$-action on $f^*V'$ recovers $V'$ and $\nabla$ is the corresponding geometric Sen operator while $D$ is the arithmetic Sen operator of the Galois representation $R\Gamma(X_{\C_p, \proet}, f^*V')$.
\end{proof}

\begin{rem}
\label{rem:htdr-checkgenflat}
One can also define a notion of generic flatness for $\mathbb{B}_\dR^{+, \dagger}$-local systems on $X_\proet$ using the functor from filtered vector bundles with connection on $X$ to $\mathbb{B}_\dR^{+, \dagger}$-local systems from \cref{rem:htdr-functortobdr+dagger}. It is then not a priori clear that generic flatness of a $\mathbb{B}_\dR^{+, \dagger}$-local system may be checked after extending scalars to $\mathbb{B}_\dR^+$. However, using that vector bundles on $X^{\HT, \dagger}$ and $\mathbb{B}_\dR^{+, \dagger}$ are equivalent via pullback by \TODO{Reference!!}, one should be able to use the above results to check that this is indeed the case.
\end{rem}

\newpage

\subsection{The comparison with filtered Hyodo--Kato cohomology}
\label{sect:hkcomp}

Recall the isomorphism
\begin{equation*}
(X^\N)_{|ut|\neq 0}\cong (X^\prism\setminus X^\dR)\times [0, 1]
\end{equation*}
from \cref{prop:defis-utneq0}. Using that $(X^\diamond\times\Spd\Q_p)^\dR_{[p^{3/2}, \infty)}$ is a closed substack of $X^\prism\setminus X^\dR$ by \cref{prop:defis-prismffdr}, we obtain a map
\begin{equation*}
(X^\diamond\times\Spd\Q_p)^\dR_{[p^{3/2}, \infty)}\times [0, 1]\rightarrow (X^\N)_{|ut|\neq 0}\;.
\end{equation*}
Moreover, using \cref{prop:defis-prismffdr} again, we also obtain a map
\begin{equation*}
(X^\diamond\times\Spd\Q_p)^\dR_{[p^{1/2}, p^{3/2}]}\cong (X^\prism)_{[p^{1/2}, p^{3/2}]}\xrightarrow{j_\dR} X^\N
\end{equation*}
and the two maps we have constructed induce the same map
\begin{equation*}
(X^\diamond\times\Spd\Q_p)^\dR_{[p^{3/2}, p^{3/2}]}\rightarrow X^\N\;,
\end{equation*}
where we restrict the first map above along $\{0\}\subseteq [0, 1]$. Thus, we overall obtain a map from
\begin{equation*}
X^{\mHK, \mathrm{pre}}\coloneqq (X^\diamond\times\Spd\Q_p)^\dR_{[p^{1/2}, p^{3/2}]}\coprod_{(X^\diamond\times\Spd\Q_p)^\dR_{[p^{3/2}, p^{3/2}]}} (X^\diamond\times\Spd\Q_p)^\dR_{[p^{3/2}, \infty)}\times [0, 1]
\end{equation*}
to $X^\N$, which is an isomorphism onto its image; here, the map ``towards the right'' in the pushout above is via the embedding $\{0\}\subseteq [0, 1]$. Finally, this map descends to a map
\begin{equation*}
i_\HK: X^\mHK\rightarrow X^\Syn\;,
\end{equation*}
which we will call the \emph{Hyodo--Kato map}. Here, the stack $X^\mHK$ is defined as follows:

\begin{defi}
Let $X$ be any Gelfand stack over $\Q_p$. The \emph{mock Hyodo--Kato stack} $X^{\mHK}$ of $X$ is defined by the coequaliser diagram
\begin{equation*}
\begin{tikzcd}
(X^\diamond\times\Spd\Q_p)^\dR_{[p^{1/2}, \infty)}\ar[r,shift left=.75ex,"\id\times\{0\}"]\ar[r,shift right=.75ex,swap,"\phi\times\{1\}"] & X^{\mHK, \mathrm{pre}}\ar[r] & X^\mHK \nospacepunct{\;.}
\end{tikzcd}
\end{equation*}
\end{defi}

\begin{figure}

\begin{center}
\begin{tikzpicture}
  \def\width{5.5} 
  \def\height{1.2*\width} 

  \draw[thick] (0, 0) -- (0, \height); 
  \draw[thick] (\width, 0) -- (\width, \height); 
  \draw[dotted, thick] (0, 0) -- (\width, 0); 

  

  
  \draw[thick] (\width-0.2, 0.2*\height) -- (\width+0.2, 0.2*\height); 
  \node at (\width+1.2, 0.2*\height) {$\phi^{-1}(X^\HT)$};
  \draw[thick] (\width-0.2, 0.4*\height) -- (\width+0.2, 0.4*\height); 
  \node at (\width+0.7, 0.4*\height) {$X^\HT$};
  \draw[thick] (\width-0.2, 0.6*\height) -- (\width+0.2, 0.6*\height); 
  \node at (\width+0.7, 0.6*\height) {$X^\dR$};
  \draw[thick] (\width-0.2, 0.8*\height) -- (\width+0.2, 0.8*\height); 
  \node at (\width+0.95, 0.8*\height) {$\phi(X^\dR)$};

  \draw[thick] (-0.2, 0.2*\height) -- (0.2, 0.2*\height); 
  \node at (-0.7, 0.2*\height) {$X^\HT$};
  \draw[thick] (-0.2, 0.4*\height) -- (0.2, 0.4*\height); 
  \node at (-0.7, 0.4*\height) {$X^\dR$};
  \draw[thick] (-0.2, 0.6*\height) -- (0.2, 0.6*\height); 
  \node at (-0.95, 0.6*\height) {$\phi(X^\dR)$};
  \draw[thick] (-0.2, 0.8*\height) -- (0.2, 0.8*\height); 
  \node at (-1.05, 0.8*\height) {$\phi^2(X^\dR)$};


  \draw[line width=0.2mm, color=blue] (\width/2-0.05*\width, 0.4*\height) -- (0, 0.4*\height);

  \def\sqr_size{0.1*\width} 
  \draw[line width=0.2mm] (\width/2 - \sqr_size/2, 0.4*\height - \sqr_size/2) rectangle (\width/2 + \sqr_size/2, 0.4*\height + \sqr_size/2); 

  \draw[line width=0.2mm] (\width/2 - \sqr_size/2, 0.4*\height - \sqr_size/2) -- (\width/2 + \sqr_size/2, 0.4*\height + \sqr_size/2); 
  \draw[line width=0.2mm] (\width/2 + \sqr_size/2, 0.4*\height - \sqr_size/2) -- (\width/2 - \sqr_size/2, 0.4*\height + \sqr_size/2); 

  
  \shade[shading=axis, bottom color=white, top color=white, middle color=green, shading angle=0](0.5*\width + \sqr_size/2, 0.4*\height - 0.1) rectangle (\width, 0.4*\height + 0.1);
  \shade[shading=axis, bottom color=white, top color=white, middle color=green, shading angle=0](0, 0.2*\height - 0.1) rectangle (\width, 0.2*\height + 0.1);

  \draw[->, thick] (2*\width, 0) -- (2*\width, \height);
  
  \draw[thick] (2*\width-0.2, 0) -- (2*\width+0.2, 0); 
  \node at (2*\width+0.6, 0) {$0$};
  \draw[thick] (2*\width-0.2, 0.2*\height) -- (2*\width+0.2, 0.2*\height); 
  \node at (2*\width+0.6, 0.2*\height) {$1$};
  \draw[thick] (2*\width-0.2, 0.4*\height) -- (2*\width+0.2, 0.4*\height); 
  \node at (2*\width+0.6, 0.4*\height) {$p$};
  \draw[thick] (2*\width-0.2, 0.6*\height) -- (2*\width+0.2, 0.6*\height); 
  \node at (2*\width+0.6, 0.6*\height) {$p^2$};
  \draw[thick] (2*\width-0.2, 0.8*\height) -- (2*\width+0.2, 0.8*\height); 
  \node at (2*\width+0.6, 0.8*\height) {$p^3$};
  
  \draw[->, thick] (1.3*\width, \height/2) -- (1.8*\width, \height/2);
  \node at (1.55*\width, \height/2+0.3) {$\kappa$};
  
  \draw[line width=1mm] (0, 0.3*\height) -- (0, \height);
  \draw[line width=1mm] (0, 0.5*\height) -- (\width+0.05, 0.5*\height);
  \draw[line width=1mm] (\width, 0.5*\height) -- (\width, \height);
  
  \fill[color=gray!50] (0.05, 0.5*\height+0.05) rectangle (\width-0.05, \height);
  
  \draw[->, thick] (-0.5*\width, 0.3*\height) -- (-0.5*\width, 0.7*\height);
  \node at (-0.5*\width-0.4, 0.5*\height) {$\phi$};
  
  \draw[line width=0.2mm, color=blue] (\width, 0.6*\height) -- (0, 0.6*\height);
  \draw[line width=0.2mm, color=blue] (\width, 0.8*\height) -- (0, 0.8*\height);

\end{tikzpicture}
\end{center}

\captionsetup{justification=centering}
\caption{A schematic picture of $X^\N$ with the image of $X^{\mHK, \mathrm{pre}}$ \\ outlined in bold (interior shaded in grey)}
\end{figure}

The most important feature of the stack $X^\mHK$ is that its category of perfect complexes agrees with the one of the actual Hyodo--Kato stack $X^\HK=(X^\diamond\times\Spd\Q_p)^\dR/\phi^\Z$.

\begin{prop}
\label{prop:hkcomp-hkpshk}
Let $X$ be any Gelfand stack. There is an equivalence of categories
\begin{equation*}
\Perf(X^\mHK)\cong\Perf(X^\HK)
\end{equation*}
induced by the diagram
\begin{equation*}
\begin{tikzcd}
& (X^\diamond\times\Spd\Q_p)^\dR_{[p^{1/2}, \infty)}\ar[rd]\ar[ld, "\id\times \{0\}", swap] & \\
X^\mHK & & X^\HK\nospacepunct{\;.}
\end{tikzcd}
\end{equation*}
\end{prop}
\begin{proof}
By \cref{cor:perf-contractible}, we have
\begin{equation*}
\Perf((X^\diamond\times\Spd\Q_p)^\dR_{[p^{3/2}, \infty)}\times [0, 1])\cong \Perf((X^\diamond\times\Spd\Q_p)^\dR_{[p^{3/2}, \infty)})
\end{equation*}
and hence the definition of $X^{\mHK, \mathrm{pre}}$ shows that pullback along
\begin{equation*}
(X^\diamond\times\Spd\Q_p)^\dR_{[p^{1/2}, p^{3/2}]}\coprod_{(X^\diamond\times\Spd\Q_p)^\dR_{[p^{3/2}, p^{3/2}]}} (X^\diamond\times\Spd\Q_p)^\dR_{[p^{3/2}, \infty)}\xrightarrow{\id\times\{0\}} X^{\mHK, \mathrm{pre}}
\end{equation*}
induces an equivalence on perfect complexes. As the pushout on the left is isomorphic to $(X^\diamond\times\Spd\Q_p)^\dR_{[p^{1/2}, \infty)}$ by \cref{lem:hkcomp-glueoverinterval}, we conclude that the categories of perfect complexes of $X^\mHK$ and
\begin{equation*}
\operatorname{coeq}(\hspace{-0.15cm}
\begin{tikzcd}
(X^\diamond\times\Spd\Q_p)^\dR_{[p^{1/2}, \infty)}\ar[r,shift left=.75ex, "\id"]
  \ar[r,shift right=.75ex,swap, "\phi"] & (X^\diamond\times\Spd\Q_p)^\dR_{[p^{1/2}, \infty)}
\end{tikzcd}
\hspace{-0.15cm})
\end{equation*}
are equivalent. However, the latter coequaliser is actually isomorphic to $X^\HK$: Indeed, since the radius map $(X^\diamond\times\Spd\Q_p)^\dR\rightarrow (0, \infty)$ intertwines $\phi$ with multiplication by $p$, this follows from base changing the isomorphism
\begin{equation*}
\operatorname{coeq}(\hspace{-0.15cm}
\begin{tikzcd}
{[p^{1/2}, \infty)}\ar[r,shift left=.75ex, "\id"]
  \ar[r,shift right=.75ex,swap, "p"] & {[p^{1/2}, \infty)}
\end{tikzcd}
\hspace{-0.15cm})
\cong 
\operatorname{coeq}(\hspace{-0.15cm}
\begin{tikzcd}
{(0, \infty)}\ar[r,shift left=.75ex, "\id"]
  \ar[r,shift right=.75ex,swap, "p"] & {(0, \infty)}
\end{tikzcd}
\hspace{-0.15cm})
\end{equation*}
to $(X^\diamond\times\Spd\Q_p)^\dR$ along the radius map by universality of colimits in $\infty$-topoi.
\end{proof}

By \cref{prop:hkcomp-hkpshk}, pullback along the map $i_\HK$ induces a realisation functor
\begin{equation*}
T_\HK: \Perf(X^\Syn)\rightarrow\Perf(X^\HK)
\end{equation*}
from perfect analytic $F$-gauges to perfect complexes on $X^\HK$, which should be thought of as coefficients for Hyodo--Kato cohomology of $X$ following \cite[Rem.\ 6.1.3]{dRFF}. Moreover, via pullback along the filtered de Rham map $i_{\dR, +}$, we obtain a realisation functor
\begin{equation*}
T_{\dR, +}: \D(X^\Syn)\rightarrow \D(X^{\dR, +})
\end{equation*}
from analytic $F$-gauges to filtered analytic $D$-modules and forgetting the filtration yields a realisation functor
\begin{equation*}
T_\dR: \D(X^\Syn)\rightarrow \D(X^\dR)\;.
\end{equation*}

Note that $i_{\dR, +}$ and $i_\HK$ are isomorphisms onto their image and that the intersection of their images is precisely the image of the map
\begin{equation*}
i_\dR: X^\dR\rightarrow X^\prism\rightarrow X^\Syn\;.
\end{equation*}
In other words, there is a commutative diagram
\begin{equation}
\label{eq:hkcomp-hksquare}
\begin{tikzcd}
X^\dR\ar[r]\ar[d] & X^{\dR, +}\ar[d, "i_{\dR, +}"] \\
X^{\mHK}\ar[r, "i_\HK"] & X^\Syn\nospacepunct{\;.}
\end{tikzcd}
\end{equation}
The main result of this section is the following:

\begin{thm}
\label{thm:hkcomp-main}
Let $X$ be a rigid smooth derived Berkovich space over $\Q_p$. Then the diagram (\ref{eq:hkcomp-hksquare}) induces an equivalence on categories of perfect complexes, i.e.\
\begin{equation*}
\Perf(X^\Syn)\cong \Perf(X^\mHK)\times_{\Perf(X^\dR)} \Perf(X^{\dR, +})\;.
\end{equation*}
\end{thm}

\begin{cor}
\label{cor:hkcomp-main}
In the setup of the previous theorem, there is an equivalence of categories
\begin{equation*}
\Perf(X^\Syn)\cong \Perf(X^\HK)\times_{\Perf(X^\dR)} \Perf(X^{\dR, +})
\end{equation*}
induced by the realisation functors $T_\HK$ and $T_{\dR, +}$. In particular, for any $E\in\Perf(X^\Syn)$, there is a pullback diagram
\begin{equation*}
\begin{tikzcd}
R\Gamma(X^\Syn, E)\ar[r]\ar[d] & R\Gamma(X^\HK, T_\HK(E))\ar[d] \\
R\Gamma(X^{\dR, +}, T_{\dR, +}(E))\ar[r] & R\Gamma(X^\dR, T_\dR(E))\nospacepunct{\;.}
\end{tikzcd}
\end{equation*}
\end{cor}
\begin{proof}
Combine \cref{thm:hkcomp-main} and \cref{prop:hkcomp-hkpshk}.
\end{proof}

Using the expected relation between the stack $X^\HK$ and Hyodo--Kato cohomology of $X$ for smooth partially proper rigid spaces $X$ over $\Q_p$, see \cite[Rem.\ 6.1.3]{dRFF}, the result above in particular shows that our definition of syntomic cohomology recovers the more classical definition of syntomic cohomology of rigid-analytic varieties over $\Q_p$ of Colmez--Nizio{\l} in that case, see e.g.\ \cite[§4.4]{padicComparisons}. Indeed, for them, the pullback diagram above is sometimes the definition of syntomic cohomology of a rigid space over $\Q_p$ while in our approach it is a theorem. Indeed, conditional on \cite[Rem.\ 6.1.3]{dRFF}, applying the above result to $E=\O\{i\}$ and using the fact that $X^{\dR, +}$ computes the Hodge filtration on de Rham cohomology, see \cite[Rem.\ 5.2.3]{dRFF}, we obtain:

\begin{cor}
Let $X$ be a smooth partially proper qcqs rigid space over $\Q_p$. For any $i\in\Z$, there is a cartesian diagram
\begin{equation*}
\begin{tikzcd}
R\Gamma_\Syn(X, \Q_p(i))\ar[r]\ar[d] & R\Gamma_\HK(X)^{\phi=p^i, N=0}\ar[d] \\
\Fil_\Hod^i R\Gamma_\dR(X)\ar[r] & R\Gamma_\dR(X)\nospacepunct{\;.}
\end{tikzcd}
\end{equation*}
\end{cor}

\subsubsection{Proof of \cref{thm:hkcomp-main}}

Let us now move on to the proof of \cref{thm:hkcomp-main}. We first recall that $X^\N$ inherits a radius map
\begin{equation*}
X^\N\xrightarrow{\pi} X^\prism\xrightarrow{\kappa} (0, \infty)
\end{equation*}
from the one on $X^\prism$ and, for any interval $I\subseteq (0, \infty)$, we write $(X^\N)_I$ for the preimage of $I$ under this map.

\begin{lem}
\label{lem:hkcomp-chopupxn}
Let $X$ be any Gelfand stack. The natural inclusion maps induce an isomorphism between the iterated pushout
\begin{equation*}
(X^\N)_{[p^{1/2}, \infty)}\coprod_{(X^\N)_{[p^{1/2}, p^{1/2}]}} (X^\N)_{[p^{-1/2}, p^{1/2}]}\coprod_{(X^\N)_{[p^{-1/2}, p^{-1/2}]}} (X^\N)_{[p^{-3/2}, p^{-1/2}]} \coprod_{(X^\N)_{[p^{-3/2}, p^{-3/2}]}} \dots
\end{equation*}
and $X^\N$.
\end{lem}
\begin{proof}
This follows by repeatedly applying \cref{lem:hkcomp-glueoverinterval}.
\end{proof}

\begin{figure}

\begin{center}
\begin{tikzpicture}
  \def\width{5.5} 
  \def\height{1.2*\width} 

  \draw[line width=1mm] (0, -0.05) -- (0, 0.6*\height); 
  \draw[line width=1mm] (\width, 0.2*\height) -- (2*\width, 0.2*\height); 
  \draw[line width=1mm] (\width, 0.2*\height-0.05) -- (\width, 0.6*\height); 
  \draw[line width=1mm] (0, 0) -- (2*\width, 0); 

  

  
  \draw[thick] (\width-0.2, 0.3*\height) -- (\width+0.2, 0.3*\height); 
  \node at (\width+0.95, 0.3*\height) {$X^\dR$};
  \draw[thick] (\width-0.2, 0.5*\height) -- (\width+0.2, 0.5*\height); 
  \node at (\width+1.05, 0.5*\height) {$\phi(X^\dR)$};

  \draw[thick] (-0.2, 0.1*\height) -- (0.2, 0.1*\height); 
  \node at (-0.95, 0.1*\height) {$X^\dR$};
  \draw[thick] (-0.2, 0.3*\height) -- (0.2, 0.3*\height); 
  \node at (-1.05, 0.3*\height) {$\phi(X^\dR)$};
  \draw[thick] (-0.2, 0.5*\height) -- (0.2, 0.5*\height); 
  \node at (-1.05, 0.5*\height) {$\phi^2(X^\dR)$};


  \draw[line width=0.2mm, color=blue] (\width/2-0.05*\width, 0.1*\height) -- (0, 0.1*\height);
  \draw[line width=0.2mm, color=blue] (\width, 0.3*\height) -- (0, 0.3*\height);
  \draw[line width=0.2mm, color=blue] (\width, 0.5*\height) -- (0, 0.5*\height);

  \def\sqr_size{0.1*\width} 
  \draw[line width=0.2mm] (\width/2 - \sqr_size/2, 0.1*\height - \sqr_size/2) rectangle (\width/2 + \sqr_size/2, 0.1*\height + \sqr_size/2); 

  \draw[line width=0.2mm] (\width/2 - \sqr_size/2, 0.1*\height - \sqr_size/2) -- (\width/2 + \sqr_size/2, 0.1*\height + \sqr_size/2); 
  \draw[line width=0.2mm] (\width/2 + \sqr_size/2, 0.1*\height - \sqr_size/2) -- (\width/2 - \sqr_size/2, 0.1*\height + \sqr_size/2); 

  
  \shade[shading=axis, bottom color=white, top color=white, middle color=green, shading angle=0](0.5*\width + \sqr_size/2, 0.1*\height - 0.1) rectangle (2*\width, 0.1*\height + 0.1);

\end{tikzpicture}
\end{center}

\captionsetup{justification=centering}
\caption{A schematic picture of $X^{\Syn, \mathrm{pre}}$ with \\ the two copies of $\partial X^{\Syn, \mathrm{pre}}$ in bold}
\end{figure}

To move on, we introduce the following notation: For any Gelfand stack $X$, let $X^{\Syn, \mathrm{pre}}$ be defined as the pushout
\begin{equation*}
(X^\N)_{[p^{1/2}, \infty)}\coprod_{(X^\prism)_{[p^{-1/2}, p^{1/2}]}} ((X^\prism)_{[p^{-1/2}, p^{1/2}]}\times [0, \infty))\;,
\end{equation*}
where the map ``towards the left'' is induced by $j_\HT$ while the one ``towards the right'' is given by $\id\times\{0\}$. Then $X^{\Syn, \mathrm{pre}}$ receives two maps from the Gelfand stack
\begin{equation*}
\partial X^{\Syn, \mathrm{pre}}\coloneqq (X^\prism)_{[p^{1/2}, \infty)}\coprod_{(X^\prism)_{[p^{1/2}, p^{1/2}]}} ((X^\prism)_{[p^{1/2}, p^{1/2}]}\times [0, \infty))\;.
\end{equation*}
One of them is induced by $j_\HT: (X^\prism)_{[p^{1/2}, \infty)}\rightarrow (X^\N)_{[p^{1/2}, \infty)}$ and the inclusion $(X^\prism)_{[p^{1/2}, p^{1/2}]}\rightarrow (X^\prism)_{[p^{-1/2}, p^{1/2}]}$ and we denote it
\begin{equation*}
j_\HT: \partial X^{\Syn, \mathrm{pre}}\rightarrow X^{\Syn, \mathrm{pre}}
\end{equation*}
as well. The other one is induced by $j_\dR: (X^\prism)_{[p^{1/2}, \infty)}\rightarrow (X^\N)_{[p^{1/2}, \infty)}$ and the isomorphism
\begin{equation*}
(X^\prism)_{[p^{1/2}, p^{1/2}]}\times [0, \infty)\xrightarrow{\cong} ((X^\prism)_{[p^{1/2}, p^{1/2}]}\times [0, 1])\coprod_{(X^\prism)_{[p^{-1/2}, p^{-1/2}]}} ((X^\prism)_{[p^{-1/2}, p^{-1/2}]}\times [0, \infty))\;,
\end{equation*}
where the maps defining the pushout on the right-hand side are $\phi\times \{1\}$ ``to the left'' and $\id\times \{0\}$ ``to the right'', respectively, where we use that
\begin{equation*}
(X^\N)_{[p^{1/2}, p^{1/2}]}\cong (X^\prism)_{[p^{1/2}, p^{1/2}]}\times [0, 1]
\end{equation*}
by \cref{prop:defis-utneq0} to see the target above as a substack of $X^{\Syn, \mathrm{pre}}$; we denote this second map by
\begin{equation*}
j_\dR: \partial X^{\Syn, \mathrm{pre}}\rightarrow X^{\Syn, \mathrm{pre}}\;.
\end{equation*}
In this notation, our next lemma reads as follows:

\begin{figure}

\begin{center}
\begin{tikzpicture}
  \def\width{5.5} 
  \def\height{1*\width} 

  \draw[line width=1mm, color=red] (0, 0) -- (0, 0.1*\height); 
  \draw[thick] (\width, 0) -- (\width, 0.1*\height); 
  \draw[thick] (0, 0.1*\height) -- (\width, 0.1*\height);
  \draw[dotted, thick] (0, 0) -- (\width, 0); 
  
  \draw[line width=1mm, color=orange] (0, 0.2*\height) -- (0, 0.4*\height);  
  \draw[line width=1mm, color=red] (\width, 0.2*\height) -- (\width, 0.4*\height);
  \draw[thick] (0, 0.2*\height) -- (\width, 0.2*\height);
  \draw[thick] (0, 0.4*\height) -- (\width, 0.4*\height);
  
  \draw[line width=1mm, color=yellow] (0, 0.5*\height) -- (0, 0.7*\height);  
  \draw[line width=1mm, color=orange] (\width, 0.5*\height) -- (\width, 0.7*\height);
  \draw[thick] (0, 0.5*\height) -- (\width, 0.5*\height);
  \draw[thick] (0, 0.7*\height) -- (\width, 0.7*\height);
  
  \draw[thick] (0, 0.8*\height) -- (0, 1.4*\height);  
  \draw[line width=1mm, color=yellow] (\width, 0.8*\height) -- (\width, \height);
  \draw[thick] (\width, 1.0*\height) -- (\width, 1.4*\height);
  \draw[thick] (0, 0.8*\height) -- (\width, 0.8*\height);

  

  
  \draw[thick] (\width-0.2, 0.3*\height) -- (\width+0.2, 0.3*\height); 
  \node at (\width+1.2, 0.3*\height) {$\phi^{-2}(X^\HT)$};
  \draw[thick] (\width-0.2, 0.6*\height) -- (\width+0.2, 0.6*\height); 
  \node at (\width+1.2, 0.6*\height) {$\phi^{-1}(X^\HT)$};
  \draw[thick] (\width-0.2, 0.9*\height) -- (\width+0.2, 0.9*\height); 
  \node at (\width+0.7, 0.9*\height) {$X^\HT$};
  \draw[thick] (\width-0.2, 1.1*\height) -- (\width+0.2, 1.1*\height); 
  \node at (\width+0.95, 1.1*\height) {$X^\dR$};
  \draw[thick] (\width-0.2, 1.3*\height) -- (\width+0.2, 1.3*\height); 
  \node at (\width+1.05, 1.3*\height) {$\phi(X^\dR)$};

  \draw[thick] (-0.2, 0.3*\height) -- (0.2, 0.3*\height); 
  \node at (-1.2, 0.3*\height) {$\phi^{-1}(X^\HT)$};
  \draw[thick] (-0.2, 0.6*\height) -- (0.2, 0.6*\height); 
  \node at (-0.7, 0.6*\height) {$X^\HT$};
  \draw[thick] (-0.2, 0.9*\height) -- (0.2, 0.9*\height); 
  \node at (-0.95, 0.9*\height) {$X^\dR$};
  \draw[thick] (-0.2, 1.1*\height) -- (0.2, 1.1*\height); 
  \node at (-1.05, 1.1*\height) {$\phi(X^\dR)$};
  \draw[thick] (-0.2, 1.3*\height) -- (0.2, 1.3*\height); 
  \node at (-1.05, 1.3*\height) {$\phi^2(X^\dR)$};


  \draw[line width=0.2mm, color=blue] (\width/2-0.05*\width, 0.9*\height) -- (0, 0.9*\height);
  \draw[line width=0.2mm, color=blue] (\width, 1.1*\height) -- (0, 1.1*\height);
  \draw[line width=0.2mm, color=blue] (\width, 1.3*\height) -- (0, 1.3*\height);

  \def\sqr_size{0.1*\width} 
  \draw[line width=0.2mm] (\width/2 - \sqr_size/2, 0.9*\height - \sqr_size/2) rectangle (\width/2 + \sqr_size/2, 0.9*\height + \sqr_size/2); 

  \draw[line width=0.2mm] (\width/2 - \sqr_size/2, 0.9*\height - \sqr_size/2) -- (\width/2 + \sqr_size/2, 0.9*\height + \sqr_size/2); 
  \draw[line width=0.2mm] (\width/2 + \sqr_size/2, 0.9*\height - \sqr_size/2) -- (\width/2 - \sqr_size/2, 0.9*\height + \sqr_size/2); 

  
  \shade[shading=axis, bottom color=white, top color=white, middle color=green, shading angle=0](0.5*\width + \sqr_size/2, 0.9*\height - 0.1) rectangle (\width-0.05, 0.9*\height + 0.1);
  \shade[shading=axis, bottom color=white, top color=white, middle color=green, shading angle=0](0.05, 0.6*\height - 0.1) rectangle (\width-0.05, 0.6*\height + 0.1);
  \shade[shading=axis, bottom color=white, top color=white, middle color=green, shading angle=0](0.05, 0.3*\height - 0.1) rectangle (\width-0.05, 0.3*\height + 0.1);

\end{tikzpicture}
\end{center}

\captionsetup{justification=centering}
\caption{Proof of \cref{lem:hkcomp-syndifferentgluing}: chop up $X^\N$, then glue along \\ the yellow, orange and red pieces to get $X^{\Syn, \mathrm{pre}}$}
\end{figure}

\begin{lem}
\label{lem:hkcomp-syndifferentgluing}
Let $X$ be any Gelfand stack. Then $X^\Syn$ is isomorphic to the coequaliser of the two maps
\begin{equation*}
\begin{tikzcd}
\partial X^{\Syn, \mathrm{pre}}\ar[r,shift left=.75ex, "j_\HT"]
  \ar[r,shift right=.75ex,swap, "j_\dR"] & X^{\Syn, \mathrm{pre}}\nospacepunct{\;.}
\end{tikzcd}
\end{equation*}
\end{lem}
\begin{proof}
By definition, $X^\Syn$ is obtained as the coequaliser of the maps $j_\HT, j_\dR: X^\prism\rightarrow X^\N$. However, by \cref{lem:hkcomp-chopupxn}, we can view $X^\N$ and $X^\prism$ as glued together from the smaller pieces 
\begin{equation*}
(X^\N)_{[p^{1/2}, \infty)}, (X^\N)_{[p^{-1/2}, p^{1/2}]}, \dots\;, \hspace{0.3cm}\text{and}\hspace{0.3cm}(X^\prism)_{[p^{1/2}, \infty)}, (X^\prism)_{[p^{-1/2}, p^{1/2}]}, \dots\;,
\end{equation*}
respectively. To obtain $X^\Syn$ from $(X^\N)_{[p^{1/2}, \infty)}, (X^\N)_{[p^{-1/2}, p^{1/2}]}, \dots$, we have to perform the following gluings:
\begin{enumerate}[label=(\arabic*)]
\item Glue the two copies of $(X^\N)_{[p^{-k/2}, p^{-k/2}]}$ for each $k\geq -1$ odd;
\item glue the two copies of $(X^\prism)_{[p^{-k/2-1}, p^{-k/2}]}$ for each $k\geq 1$ odd (these are embedded via $j_\HT$ or $j_\dR$ into $(X^\N)_{[p^{-k/2}, p^{-k/2+1}]}$ or $(X^\N)_{[p^{-k/2-1}, p^{-k/2}]}$, respectively);
\item glue the two copies of $(X^\prism)_{[p^{-1/2}, p^{1/2}]}$, which are embedded into the pieces $(X^\N)_{[p^{1/2}, \infty)}$ and $(X^\N)_{[p^{-1/2}, p^{1/2}]}$ via $j_\HT$ and $j_\dR$, respectively;
\item glue the two copies of $(X^\prism)_{[p^{1/2}, \infty)}$, which are embedded into $(X^\N)_{[p^{1/2}, \infty)}$ via $j_\HT$ and $j_\dR$.
\end{enumerate}

We now change the order of these gluings: For this, first note that
\begin{equation*}
(X^\N)_{[p^{-k/2}, p^{-k/2+1}]}\cong (X^\prism)_{[p^{-k/2}, p^{-k/2+1}]}\times [0, 1]\cong (X^\prism)_{[p^{-1/2}, p^{1/2}]}\times [0, 1]
\end{equation*}
for each $k\geq 1$ odd, where the first isomorphism is due to \cref{prop:defis-utneq0} and the second one is via a $(k-1)/2$-fold application of $\phi$, which we recall defines an isomorphism $(X^\prism)_{(0, 1)}\cong (X^\prism)_{(0, p)}$ by \cref{prop:defis-xprismxdiv1}. Thus, first performing the gluings (2) and (3), we obtain $X^{\Syn, \mathrm{pre}}$ as defined above.

Observing that
\begin{equation*}
(X^\N)_{[p^{-k/2}, p^{-k/2}]}\cong (X^\prism)_{[p^{-k/2}, p^{-k/2}]}\times [0, 1]\cong (X^\prism)_{[p^{1/2}, p^{1/2}]}\times [0, 1]
\end{equation*}
by another application of \cref{prop:defis-utneq0} and \cref{prop:defis-xprismxdiv1}, we now see that subsequently performing the gluings (1) and (4) equivalently amounts to gluing the two copies of $\partial X^{\Syn, \mathrm{pre}}$ embedded into $X^{\Syn, \mathrm{pre}}$ via $j_\HT$ and $j_\dR$. This proves the claim.
\end{proof}

\begin{lem}
For any Gelfand stack $X$, pullback along the morphism
\begin{equation*}
\operatorname{coeq}(\hspace{-0.15cm}
\begin{tikzcd}
(X^\prism)_{[p^{1/2}, \infty)}\ar[r,shift left=.75ex, "j_\HT"]
  \ar[r,shift right=.75ex,swap, "j_\dR"] & (X^\N)_{[p^{1/2}, \infty)}
\end{tikzcd}
\hspace{-0.15cm})\rightarrow X^\Syn
\end{equation*}
induces an equivalence on perfect complexes.
\end{lem}
\begin{proof}
This immediately follows from the previous lemma once we know that
\begin{equation*}
\Perf(X^{\Syn, \mathrm{pre}})\cong \Perf((X^\N)_{[p^{1/2}, \infty)})\;, \hspace{0.3cm}\text{and}\hspace{0.3cm} \Perf(\partial X^{\Syn, \mathrm{pre}})\cong \Perf((X^\prism)_{[p^{1/2}, \infty)})
\end{equation*}
via pullback along the natural maps $(X^\N)_{[p^{1/2}, \infty)}\rightarrow X^{\Syn, \mathrm{pre}}$ and $(X^\prism)_{[p^{1/2}, \infty)}\rightarrow \partial X^{\Syn, \mathrm{pre}}$. However, using the pushout definition of $X^{\Syn, \mathrm{pre}}$, the first equivalence follows from the fact that pullback along $\id\times\{0\}$ induces an equivalence
\begin{equation*}
\Perf((X^\prism)_{[p^{-1/2}, p^{1/2}]}\times [0, \infty))\cong \Perf((X^\prism)_{[p^{-1/2}, p^{1/2}]})\;,
\end{equation*}
which is a consequence of \cref{cor:perf-contractible}. Similarly, the second one follows from the fact that pullback along $\id\times\{0\}$ induces an equivalence
\begin{equation*}
\Perf((X^\prism)_{[p^{1/2}, p^{1/2}]}\times [0, \infty))\cong \Perf((X^\prism)_{[p^{1/2}, p^{1/2}]})\;,
\end{equation*}
which is proved by the same argument.
\end{proof}

\begin{lem}
Let $X$ be any Gelfand stack and assume that $(X^\N)_{[p^{1/2}, p^{3/2}]}$ is nicely coverable. Then pullback along the map
\begin{equation*}
(X^\prism)_{[p^{1/2}, p^{3/2}]}\coprod_{X^\dR} X^{\dR, +}\rightarrow (X^\N)_{[p^{1/2}, p^{3/2}]}\;,
\end{equation*}
which is given by $i_{\dR, +}$ and $j_\dR$, induces an equivalence on perfect complexes.
\end{lem}
\begin{proof}
First note that, by the argument from the proof of \cref{lem:perf-keylemma}, the stack $(X^\N)_{[p^{1/2}, p^{3/2}]}$ being nicely coverable implies the same for $(X^\N)_{|ut|=0}=(X^\N)_{[p, p]}$ and $(X^\prism)_{[p^{1/2}, p^{3/2}]}\cong((X^\N)_{[p^{1/2}, p^{3/2}]})_{|t|=1}$, where the last isomorphism is via $j_\dR$. Now consider the overconvergent normed divisor $Z=\{|ut|=0\}\subseteq (X^\N)_{[p^{1/2}, p^{3/2}]}$. Note that, for each $\epsilon>0$ small enough, we have
\begin{equation*}
\begin{split}
\Perf(Z_\epsilon\setminus Z)&\cong \Perf(((X^\prism)_{|\widetilde{\mu}|\leq\epsilon}\setminus X^\dR)\times [0, 1]) \\
&\cong \Perf((X^\prism)_{|\widetilde{\mu}|\leq\epsilon}\setminus X^\dR)
\end{split}
\end{equation*}
by \cref{prop:defis-utneq0} and \cref{cor:perf-contractible} and, similarly, 
\begin{equation*}
\Perf((X^\N)_{[p^{1/2}, p^{3/2}]}\setminus Z)\cong \Perf((X^\prism)_{[p^{1/2}, p^{3/2}]}\setminus X^\dR)\;.
\end{equation*}
Thus, using \cref{lem:htdr-perfut0}, an application of \cref{cor:perf-corkeylemma} yields
\begin{equation*}
\Perf((X^\N)_{[p^{1/2}, p^{3/2}]})\cong \Perf(X^{\dR, +})\times_{\colim_{\epsilon>0} \Perf((X^\prism)_{|\widetilde{\mu}|\leq\epsilon}\setminus X^\dR)} \Perf((X^\prism)_{[p^{1/2}, p^{3/2}]}\setminus X^\dR)\;.
\end{equation*}
Comparing this to the isomorphism
\begin{equation*}
\Perf((X^\prism)_{[p^{1/2}, p^{3/2}]})\cong \Perf(X^\dR)\times_{\colim_{\epsilon>0} \Perf((X^\prism)_{|\widetilde{\mu}|\leq\epsilon}\setminus X^\dR)} \Perf((X^\prism)_{[p^{1/2}, p^{3/2}]}\setminus X^\dR)
\end{equation*}
obtained by applying \cref{cor:perf-corkeylemma} to the overconvergent normed divisor 
\begin{equation*}
X^\dR\cong (\{|\widetilde{\mu}|=0\}\subseteq (X^\prism)_{[p^{1/2}, p^{3/2}]})\;,
\end{equation*}
where the isomorphism is via $i_\dR$, we see that
\begin{equation*}
\Perf((X^\N)_{[p^{1/2}, p^{3/2}]})\cong \Perf(X^{\dR, +})\times_{\Perf(X^\dR)} \Perf((X^\prism)_{[p^{1/2}, p^{3/2}]})
\end{equation*}
and this concludes the proof.
\end{proof}

\begin{prop}
Let $X$ be any Gelfand stack and assume that $(X^\N)_{[p^{1/2}, p^{3/2}]}$ is nicely coverable. Then the conclusion of \cref{thm:hkcomp-main} holds.
\end{prop}
\begin{proof}
Noting that 
\begin{equation*}
(X^\N)_{[p^{1/2}, p^{3/2}]}\coprod_{(X^\N)_{[p^{3/2}, p^{3/2}]}} (X^\N)_{[p^{3/2}, \infty)}\xrightarrow{\cong} (X^\N)_{[p^{1/2}, \infty)}
\end{equation*}
via the natural map by \cref{lem:hkcomp-glueoverinterval} and similarly for $X^\prism$ in place of $X^\N$, the claim follows by combining the previous two lemmas upon recalling that 
\begin{equation*}
(X^\N)_{[p^{3/2}, \infty)}\cong (X^\diamond\times \Spd\Q_p)^\dR_{[p^{3/2}, \infty)}\times [0, 1]
\end{equation*}
and
\begin{equation*}
(X^\prism)_{[p^{1/2}, p^{3/2}]}\cong (X^\diamond\times \Spd\Q_p)^\dR_{[p^{1/2}, p^{3/2}]}
\end{equation*}
by \cref{prop:defis-prismffdr} and \cref{prop:defis-utneq0}.
\end{proof}

\begin{proof}[Proof of \cref{thm:hkcomp-main}]
By compatibility of the statement with rigid étale localisations and strict closed covers, we may assume that $X$ is finite étale over a closed subset of $\ol{\DD}^n$ for some $n$. However, in that case, the claim follows by combining the preceding proposition with \cref{prop:perf-coversmoothrigid}.
\end{proof}

\subsubsection{Bloch--Kato's $H^1_g$ via the syntomification}

Using the results we have obtained above, we can now deduce the following theorem about vector bundles on $\Q_p^\Syn$ and their cohomology. One should compare this with a result of Bhatt--Lurie which says that reflexive sheaves on the algebraic syntomification of $\Z_p$ are equivalent to lattices in crystalline Galois representations and that, after inverting $p$, the first coherent cohomology group of such a sheaf agrees with Bloch--Kato's group $H^1_f$ for the corresponding crystalline representation.

\begin{thm}
\label{thm:bk-mainqp}
There is an equivalence of categories
\begin{equation*}
\Vect(\Q_p^\Syn)\cong \Vect^\dR(\FF_{\Q_p})
\end{equation*}
between vector bundles on $\Q_p^\Syn$ and vector bundles on the Fargues--Fontaine curve of $\Q_p$ which are de Rham in the sense of \cite[Def.\ 15.12]{FarguesFontaine}. Moreover, if $V$ and $W$ are de Rham representations of $G_{\Q_p}$, then
\begin{equation*}
\Ext_{\Q_p^\Syn}(V, W)\cong \RHom_{\Rep^\dR(G_{\Q_p})}(V, W)\;.
\end{equation*}
In particular, we have
\begin{equation*}
H^1(\Q_p^\Syn, V)\cong H^1_g(G_{\Q_p}, V)\;.
\end{equation*}
\end{thm}
\comment{
One can also prove a result of a similar flavour for $\Z_p$ in place of $\Q_p$:

\begin{thm}
\label{thm:bk-mainzp}
There is an equivalence of categories
\begin{equation*}
\Vect(\Z_p^\Syn)\cong \Vect^\pcrys(\FF_{\Q_p})
\end{equation*}
between vector bundles on $\Q_p^\Syn$ and vector bundles on the Fargues--Fontaine curve of $\Q_p$ which are potentially crystalline in the sense of \cite[Def.\ 10.2.13]{FarguesFontaine}. Moreover, if $V$ and $W$ are potentially crystalline representations of $G_{\Q_p}$, then
\begin{equation*}
\Ext_{\Q_p^\Syn}(V, W)\cong \RHom_{\Rep^\pcrys(G_{\Q_p})}(V, W)\;.
\end{equation*}
\end{thm}

\begin{rem}
In view of the result above, one may wonder if there is also an interpretation of the group $H^1_e(G_{\Q_p}, V)$ in terms of the stacks $\Q_p^\Syn$ and $\Z_p^\Syn$ for any de Rham $G_{\Q_p}$-representation. Indeed, we expect that $H^1_e(G_{\Q_p}, V)$ is somehow related to cohomology of $\Q_p^\Syn$ with compact support.
\end{rem}

\begin{rem}
It is expected that the conclusion of \cref{thm:bk-mainqp} generalises to arbitrary smooth partially proper rigid spaces $X$ in the sense that there is an equivalence
\begin{equation*}
\Vect(X^\Syn)\cong \Vect^\dR(\FF_{X^\diamond})\;;
\end{equation*}
this is in keeping with \cite[Rem.\ 7.5.7]{dRFF}. Here, a vector bundle on $\FF_{X^\diamond}$ is called de Rham if it is de Rham after pullback along each classical point of $X$. However, we do \emph{not} expect that cohomology on $X^\Syn$ computes the Ext groups between de Rham local systems on $X$: Indeed, the analogous statement for the integral syntomification of Bhatt--Lurie and Drinfeld also fails, see \cite[Ex.\ 4.8]{Pentland}, and the same counterexample works in our case.
\end{rem}
}
\begin{proof}
By \cite[Thm.\ 7.1.1]{dRFF}, the category of vector bundles of $\Q_p^\HK$ is equivalent to the category of $(\phi, N, G_{\Q_p})$-modules over $\Q_p^\un$ and thus, using \cref{cor:hkcomp-main}, we see that vector bundles on $\Q_p^\Syn$ are equivalent to filtered $(\phi, N, G_{\Q_p})$-modules over $\Q_p^\un$. As these are equivalent to potentially log-crystalline vector bundles on $\FF_{\Q_p}$ by \cite[Prop.\ 10.6.7]{FarguesFontaine}, the first part of the theorem now follows by the theorem that ``de Rham $=$ potentially log-crystalline'', see \cite[Thm.\ 10.6.10]{FarguesFontaine}.

For the second part about Ext groups, we can reduce to the case $V=\Q_p$ by replacing $W$ by $W\tensor V^\vee$ and then we have to check whether the pullback
\begin{equation*}
\begin{tikzcd}
R\Gamma(\Q_p^\Syn, W)\ar[r]\ar[d] & R\Gamma(\Q_p^\HK, T_\HK(W))\ar[d] \\
R\Gamma(\Q_p^{\dR, +}, T_{\dR, +}(W))\ar[r] & R\Gamma(\Q_p^\dR, T_\dR(W))
\end{tikzcd}
\end{equation*}
computes the Exts between $\Q_p$ and $W$ in de Rham representations. Recalling that cohomology on $\Q_p^\HK$ computes cohomology of $(\phi, N, G_{\Q_p})$-modules by \cite[Thm.\ 7.1.1]{dRFF} and using the theorem ``de Rham $=$ potentially semistable'', the claim follows by Galois descent from a computation of Emerton--Kisin of the Ext groups of semistable representations in terms of filtered $(\phi, N)$-modules, see \cite[Cor.\ 2.4.4]{CrystallineExt}.
\end{proof}

\comment{
To also prove \cref{thm:bk-mainzp}, we first have to identify the category of perfect complexes on $\Z_p^\HK$ in a similar spirit to the identification of $\Perf(\Q_p^\HK)$ with the bounded derived category of $(\phi, N, G_{\Q_p})$-modules over $\Q_p^\un$ from \cite[Thm.\ 7.1.1]{dRFF}.

\begin{prop}
Pullback along the natural morphism
\begin{equation*}
\O_{\C_p}^\HK\rightarrow \ol{\F}_p^\HK\cong \GSpec\Q_p^\un/\phi^\Z
\end{equation*}
induces an equivalence on perfect complexes. In particular, for any $K\subseteq\ol{\Q}_p$, the category of perfect complexes on $\O_K^\HK$ is equivalent to the bounded derived category of $(\phi, G_K)$-modules over $\Q_p^\un$.
\end{prop}
\begin{proof}
This is merely a variant of \cite[Thm.\ 7.1.1]{dRFF} and can be deduced by following the proof in loc.\ cit. with only minor adjustments. The key point is that, whenever one encounters a punctured perfectoid open unit disk in the proof in loc.\ cit., it will now be replaced by a non-punctured perfectoid open unit disk. In particular, the key computation in Lem.\ 7.4.7 of loc.\ cit.\ needed to establish full faithfulness now reduces to the elementary fact that de Rham cohomology of the open unit disk over $\Q_p$ is concentrated in degree $0$, where it is given by $\Q_p$. Finally, we recall that the classification of vector bundles on the Fargues--Fontaine curve used in the proof of essential surjectivity in loc.\ cit.\ extends to $\phi$-modules on $YY_{(0, \infty]}$ by virtue of \cite[Thm.\ 13.2.1]{Berkeley}.
\end{proof}

\begin{proof}[Proof of \cref{thm:bk-mainzp}]
Using the previous proposition, \cref{cor:hkcomp-main} shows that vector bundles on $\Z_p^\Syn$ are equivalent to filtered $(\phi, G_{\Q_p})$-modules over $\Q_p^\un$. Now recall that, for any extension $K$ of $\Q_p$, a $G_K$-equivariant vector bundle on $\FF_{\C_p}$ is equivalent to a $B$-pair in the sense of Berger, i.e.\ a pair $(M_\dR, M_e)$ of a $B_\dR^+$- and a $B_e$-representation of $G_K$ together with an isomorphism after extending scalars to $B_\dR$. As crystalline $B_e$-respresentations of $G_K$ are equivalent to $\phi$-modules over $K_0$, the maximal unramified subextension of $K$, by \cite[Thm.\ 10.2.14]{FarguesFontaine}, and $G_K$-stable lattices in a free $B_\dR$-representation of $G_K$ are equivalent to filtered $K$-vector spaces by \cite[Prop.\ 10.4.3]{FarguesFontaine}, our claim about vector bundles on $\Z_p^\Syn$ follows. The assertion about Ext groups now follows similarly to above from a computation of Emerton--Kisin of Ext groups in crystalline representations, see \cite[Cor.\ 2.4.4]{CrystallineExt}.
\end{proof}
}

\newpage

\subsection{The comparison with proétale cohomology}
\label{sect:proet}

Note that, we can also do the argument from the previous section ``in the other direction'': Instead of performing a retraction argument as above which makes $(X^\prism)_{(p, \infty)}$ appear, we now want to see $(X^\prism)_{(0, p)}$ instead and use \cref{lem:htdr-perfut0viaht+} in place of \cref{lem:htdr-perfut0}. Indeed, this is possible and now yields a description of perfect complexes on $X^\Syn$ in terms of ``étale'' rather than ``de Rham'' data, as in \cref{cor:hkcomp-main}.

To make things more precise, recall the isomorphism
\begin{equation*}
(X^\N)_{|ut|\neq 0}\cong (X^\prism\setminus X^\dR)\times [0, 1]
\end{equation*}
from \cref{prop:defis-utneq0}, using which we obtain a map
\begin{equation*}
(X^\prism)_{(0, p^{1/2}]}\times [0, 1]\rightarrow (X^\N)_{|ut|\neq 0}\;.
\end{equation*}
Moreover, using \cref{prop:defis-xprismxdiv1}, we also obtain a map
\begin{equation*}
(X^\prism)_{[p^{-1/2}, p^{1/2}]}\xrightarrow{j_\HT} X^\N
\end{equation*}
and the two maps we have constructed induce the same map
\begin{equation*}
(X^\prism)_{[p^{-1/2}, p^{-1/2}]}\rightarrow X^\N\;,
\end{equation*}
where we restrict the first map above along $\{1\}\subseteq [0, 1]$ and precompose by Frobenius. Thus, we overall obtain a map from
\begin{equation*}
X^{\mDiv, \mathrm{pre}}\coloneqq (X^\prism)_{[p^{-1/2}, p^{1/2}]}\coprod_{(X^\prism)_{[p^{-1/2}, p^{-1/2}]}} (X^\prism)_{(0, p^{1/2}]}\times [0, 1]
\end{equation*}
to $X^\N$, which is an isomorphism onto its image; here, the map ``towards the right'' in the pushout above is via Frobenius and the embedding $\{1\}\subseteq [0, 1]$. 

Note that this pushout has two closed substacks isomorphic to $(X^\prism)_{(0, p^{1/2}]}$: There is a closed embedding
\begin{equation*}
j_\dR: (X^\prism)_{(0, p^{1/2}]}\xrightarrow{\id\times \{0\}}(X^\prism)_{(0, p^{1/2}]}\times [0, 1]\rightarrow X^{\mDiv, \mathrm{pre}}
\end{equation*}
as well as a closed embedding
\begin{equation*}
j_\HT: (X^\prism)_{(0, p^{1/2}]}\cong (X^\prism)_{[p^{-1/2}, p^{1/2}]}\coprod_{(X^\prism)_{[p^{-1/2}, p^{-1/2}]}} (X^\prism)_{(0, p^{-1/2}]}\xrightarrow{\id\coprod \phi\times \{1\}} X^{\mDiv, \mathrm{pre}}\;,
\end{equation*}
where the first isomorphism is a consequence of \cref{lem:hkcomp-glueoverinterval}. Then, finally, the map $X^{\mDiv, \mathrm{pre}}\rightarrow X^\N$ constructed above descends to a map
\begin{equation*}
i_\HK: X^\mDiv\rightarrow X^\Syn\;,
\end{equation*}
which we will call the \emph{Fargues--Fontaine map}. Here, the stack $X^\mDiv$ is defined as follows:

\begin{defi}
Let $X$ be any Gelfand stack over $\Q_p$. The \emph{mock perfect prismatisation} $X^{\mDiv}$ of $X$ is defined by the coequaliser diagram
\begin{equation*}
\begin{tikzcd}
(X^\prism)_{(0, p^{1/2}]}\ar[r,shift left=.75ex,"j_\dR"]\ar[r,shift right=.75ex,swap,"j_\HT"] & X^{\mDiv, \mathrm{pre}}\ar[r] & X^\mDiv \nospacepunct{\;.}
\end{tikzcd}
\end{equation*}
\end{defi}

\begin{figure}

\begin{center}
\begin{tikzpicture}
  \def\width{5.5} 
  \def\height{1.2*\width} 

  \draw[thick] (0, 0) -- (0, \height); 
  \draw[thick] (\width, 0) -- (\width, \height); 
  \draw[dotted, thick] (0, 0) -- (\width, 0); 

  

  
  \draw[thick] (\width-0.2, 0.2*\height) -- (\width+0.2, 0.2*\height); 
  \node at (\width+1.2, 0.2*\height) {$\phi^{-2}(X^\HT)$};
  \draw[thick] (\width-0.2, 0.4*\height) -- (\width+0.2, 0.4*\height); 
  \node at (\width+1.2, 0.4*\height) {$\phi^{-1}(X^\HT)$};
  \draw[thick] (\width-0.2, 0.6*\height) -- (\width+0.2, 0.6*\height); 
  \node at (\width+0.7, 0.6*\height) {$X^\HT$};
  \draw[thick] (\width-0.2, 0.8*\height) -- (\width+0.2, 0.8*\height); 
  \node at (\width+0.7, 0.8*\height) {$X^\dR$};

  \draw[thick] (-0.2, 0.2*\height) -- (0.2, 0.2*\height); 
  \node at (-1.2, 0.2*\height) {$\phi^{-1}(X^\HT)$};
  \draw[thick] (-0.2, 0.4*\height) -- (0.2, 0.4*\height); 
  \node at (-0.7, 0.4*\height) {$X^\HT$};
  \draw[thick] (-0.2, 0.6*\height) -- (0.2, 0.6*\height); 
  \node at (-0.7, 0.6*\height) {$X^\dR$};
  \draw[thick] (-0.2, 0.8*\height) -- (0.2, 0.8*\height); 
  \node at (-0.95, 0.8*\height) {$\phi(X^\dR)$};


  \draw[line width=0.2mm, color=blue] (\width/2-0.05*\width, 0.6*\height) -- (0, 0.6*\height);
  \draw[line width=0.2mm, color=blue] (\width, 0.8*\height) -- (0, 0.8*\height);

  \def\sqr_size{0.1*\width} 
  \draw[line width=0.2mm] (\width/2 - \sqr_size/2, 0.6*\height - \sqr_size/2) rectangle (\width/2 + \sqr_size/2, 0.6*\height + \sqr_size/2); 

  \draw[line width=0.2mm] (\width/2 - \sqr_size/2, 0.6*\height - \sqr_size/2) -- (\width/2 + \sqr_size/2, 0.6*\height + \sqr_size/2); 
  \draw[line width=0.2mm] (\width/2 + \sqr_size/2, 0.6*\height - \sqr_size/2) -- (\width/2 - \sqr_size/2, 0.6*\height + \sqr_size/2); 

  
  \shade[shading=axis, bottom color=white, top color=white, middle color=green, shading angle=0](0.5*\width + \sqr_size/2, 0.6*\height - 0.1) rectangle (\width, 0.6*\height + 0.1);

  \draw[->, thick] (2*\width, 0) -- (2*\width, \height);
  
  \draw[thick] (2*\width-0.2, 0) -- (2*\width+0.2, 0); 
  \node at (2*\width+0.6, 0) {$0$};
  \draw[thick] (2*\width-0.2, 0.2*\height) -- (2*\width+0.2, 0.2*\height); 
  \node at (2*\width+0.6, 0.2*\height) {$1/p$};
  \draw[thick] (2*\width-0.2, 0.4*\height) -- (2*\width+0.2, 0.4*\height); 
  \node at (2*\width+0.6, 0.4*\height) {$1$};
  \draw[thick] (2*\width-0.2, 0.6*\height) -- (2*\width+0.2, 0.6*\height); 
  \node at (2*\width+0.6, 0.6*\height) {$p$};
  \draw[thick] (2*\width-0.2, 0.8*\height) -- (2*\width+0.2, 0.8*\height); 
  \node at (2*\width+0.6, 0.8*\height) {$p^2$};
  
  \draw[->, thick] (1.3*\width, \height/2) -- (1.8*\width, \height/2);
  \node at (1.55*\width, \height/2+0.3) {$\kappa$};
  
  \fill[color=gray!50] (0.05, 0) rectangle (\width-0.05, 0.5*\height-0.05);
  
  \draw[->, thick] (-0.5*\width, 0.3*\height) -- (-0.5*\width, 0.7*\height);
  \node at (-0.5*\width-0.4, 0.5*\height) {$\phi$};
  
  \shade[shading=axis, bottom color=gray!50, top color=gray!50, middle color=green, shading angle=0](0, 0.4*\height - 0.1) rectangle (\width, 0.4*\height + 0.1);
  \shade[shading=axis, bottom color=gray!50, top color=gray!50, middle color=green, shading angle=0](0, 0.2*\height - 0.1) rectangle (\width, 0.2*\height + 0.1);
  
  \draw[line width=1mm] (\width, 0) -- (\width, 0.7*\height);
  \draw[line width=1mm] (0, 0) -- (0, 0.5*\height+0.05);
  \draw[line width=1mm] (0, 0.5*\height) -- (\width, 0.5*\height);

\end{tikzpicture}
\end{center}

\captionsetup{justification=centering}
\caption{A schematic picture of $X^\N$ with the image of $X^{\mDiv, \mathrm{pre}}$ \\ outlined in bold (interior shaded in grey)}
\end{figure}

As before with the stack $X^{\mHK}$, the most important feature of the stack $X^{\mDiv}$ is that its category of perfect complexes agrees with the one of the actual perfect prismatisation $X^{\Div^1}$.

\begin{prop}
\label{prop:proet-psprism}
Let $X$ be any Gelfand stack over $\Q_p$. There is an equivalence of categories
\begin{equation*}
\Perf(X^{\mDiv})\cong\Perf(X^{\Div^1})
\end{equation*}
induced by the diagram
\begin{equation*}
\begin{tikzcd}
& (X^\prism)_{(0, p^{1/2}]}\ar[rd]\ar[ld, "\id\times \{0\}", swap] & \\
X^{\mDiv} & & X^{\Div^1}\nospacepunct{\;.}
\end{tikzcd}
\end{equation*}
\end{prop}
\begin{proof}
As in \cref{prop:hkcomp-hkpshk}.
\end{proof}

Thus, pullback along the map $i_{\Div^1}$ induces a functor
\begin{equation*}
T_\FF: \Perf(X^\Syn)\rightarrow\Perf(X^{\Div^1})
\end{equation*}
which we call the \emph{Fargues--Fontaine realisation} from perfect analytic $F$-gauges to perfect complexes on $X^{\Div^1}$. Indeed, to justify the name, recall that, for any rigid smooth derived Berkovich space $X$ over $\Q_p$, the category $\Perf(X^{\Div^1})$ is equivalent to the category of perfect complexes on the Fargues--Fontaine curve $\FF_{X^\diamond}$ of $X^\diamond$ by \TODO{Reference!}. Moreover, via pullback along the map
\begin{equation*}
i_{\HT, \dagger, +}: X^{\HT, \dagger, +}\rightarrow X^\N\rightarrow X^\Syn\;,
\end{equation*}
we obtain a realisation functor
\begin{equation*}
T_{\HT, \dagger, +}: \Perf(X^\Syn)\rightarrow \Perf(X^{\HT, \dagger, +})
\end{equation*}
from perfect analytic $F$-gauges to perfect complexes on $X^{\HT, \dagger, +}$ -- recall that we have explicitly described such complexes in \cref{thm:htdr-xhtdrdagger+} for the pieces of a sufficiently fine strict closed cover of any rigid smooth derived Berkovich space over $\Q_p$. Finally, further pulling back to $X^{\HT, \dagger}$ yields a realisation functor
\begin{equation*}
T_{\HT, \dagger}: \Perf(X^\Syn)\rightarrow \Perf(X^{\HT, \dagger})
\end{equation*}
from perfect analytic $F$-gauges to perfect complexes on $X^{\HT, \dagger}$.

Note that $i_{\HT, \dagger, +}$ and $i_{\Div^1}$ are isomorphisms onto their image and that the intersection of their images is precisely the image of the map
\begin{equation*}
i_{\HT, \dagger}: X^{\HT, \dagger}\rightarrow X^\prism\rightarrow X^\Syn\;.
\end{equation*}
In other words, there is a commutative diagram
\begin{equation}
\label{eq:proet-proetsquare}
\begin{tikzcd}
X^{\HT, \dagger}\ar[r]\ar[d] & X^{\HT, \dagger, +}\ar[d, "i_{\dR, +}"] \\
X^{\mDiv}\ar[r, "i_\HK"] & X^\Syn\nospacepunct{\;.}
\end{tikzcd}
\end{equation}
The main result of this section is the following:

\begin{thm}
\label{thm:proet-main}
Let $X$ be a rigid smooth derived Berkovich space over $\Q_p$ or over $\C_p$. Then the diagram (\ref{eq:proet-proetsquare}) induces an equivalence on categories of perfect complexes, i.e.\
\begin{equation*}
\Perf(X^\Syn)\cong \Perf(X^{\mDiv})\times_{\Perf(X^{\HT, \dagger})} \Perf(X^{\HT, \dagger, +})\;.
\end{equation*}
\end{thm}
\begin{proof}
Imitate the proof of \cref{thm:hkcomp-main} with the modifications discussed above, in particular replacing the use of \cref{lem:htdr-perfut0} by \cref{lem:htdr-perfut0viaht+}.
\end{proof}

\begin{cor}
\label{cor:proet-main}
In the setup of the previous theorem, there is an equivalence of categories
\begin{equation*}
\Perf(X^\Syn)\cong \Perf(X^{\Div^1})\times_{\Perf(X^{\HT, \dagger})} \Perf(X^{\HT, \dagger, +})
\end{equation*}
induced by the realisation functors $T_\FF$ and $T_{\HT, \dagger, +}$. In particular, for any $E\in\Perf(X^\Syn)$, there is a pullback diagram
\begin{equation*}
\begin{tikzcd}
R\Gamma(X^\Syn, E)\ar[r]\ar[d] & R\Gamma(X^{\Div^1}, T_\FF(E))\ar[d] \\
R\Gamma(X^{\HT, \dagger, +}, T_{\HT, \dagger, +}(E))\ar[r] & R\Gamma(X^{\HT, \dagger}, T_{\HT, \dagger}(E))\nospacepunct{\;.}
\end{tikzcd}
\end{equation*}
\end{cor}
\begin{proof}
Combine \cref{thm:proet-main} and \cref{prop:proet-psprism}.
\end{proof}

Using the equivalence between perfect complexes on $X^{\Div^1}$ and $\FF_{X^\diamond}$ and the relation between cohomology on $\FF_{X^\diamond}$ and proétale cohomology of $X$ for any smooth partially proper rigid space $X$ over $\Q_p$ from \TODO{Reference!!}, this lets us recover the classical comparison between syntomic cohomology of rigid-analytic varieties over $\Q_p$ and their proétale cohomology in a certain range, see e.g.\ \cite[Thm.\ 1.1.(4)]{padicComparisons}.

\begin{thm}
\label{thm:proet-main2}
Let $X$ be a smooth partially proper rigid space over $\Q_p$. If $E\in\Vect(X^\Syn)$ is a vector bundle analytic $F$-gauge with Hodge--Tate weights all at most $-i$ for some $i\geq 0$, then the natural morphism
\begin{equation*}
R\Gamma(X^\Syn, E)\rightarrow R\Gamma(\FF_{X^\diamond}, T_\FF(E))
\end{equation*}
is an isomorphism on $\tau^{\leq i}$ and induces an injection on $H^{i+1}$. In particular, for $E=\O\{i\}$, we obtain
\begin{equation*}
\tau^{\leq i} R\Gamma_\Syn(X, \Q_p(i))\cong \tau^{\leq i}R\Gamma_\proet(X, \Q_p(i))\;.
\end{equation*}
\end{thm}
\begin{proof}
By \cref{cor:proet-main}, our task is to prove that the cofibre of the morphism
\begin{equation*}
R\Gamma(X^{\HT, \dagger, +}, T_{\HT, \dagger, +}(E))\rightarrow R\Gamma(X^{\HT, \dagger}, T_{\HT, \dagger}(E))
\end{equation*}
is concentrated in degrees at least $i+1$ and the compatibility of the constructions $X\mapsto X^{\HT, \dagger, +}$ and $X\mapsto X^{\HT, \dagger}$ with rigid étale localisation then reduces us to the case where $X$ admits a rigid étale map to an overconvergent $n$-dimensional torus $\ol{\T}^n$. However, in this case, \cref{cor:htdr-cohomologyxhtdrdagger+} gives precise double complexes computing the cohomologies in question. Indeed, inspecting these double complexes, we see that our claim follows once we show that the ascending filtration $\Fil_\bullet V$ associated to the restriction of $E$ to $X^{\HT, \dagger, +}$ via \cref{thm:htdr-xhtdrdagger+} stabilises from $\Fil_{-i} V$ onwards, i.e.\
\begin{equation*}
\Fil_{-i} V=\Fil_{-i+1} V=\dots\;.
\end{equation*}

However, this is ensured by our assumption on the Hodge--Tate weights of $E$: Indeed, recall from \cref{thm:htdr-xhtdrdagger+} that the (descendingly) filtered vector bundle with connection $(\Fil^\bullet W, \nabla)$ on $X$ corresponding to the restriction of $E$ to $X^{\dR, +}$ is obtained by descending the associated graded pieces of $\Fil_\bullet V$ to $X$, i.e.\ $\gr_k V$ descends to $\Fil^k W$ up to a Tate twist, and hence we equivalently have to show that
\begin{equation*}
\dots=\Fil^{-i+2} W=\Fil^{-i+1} W=0\;.
\end{equation*}

Observe, however, that our assumption on the Hodge--Tate weights precisely means that $\gr^k W=0$ for $k\geq -i+1$. Now we are done by noting that the decreasing filtration $\Fil^\bullet W$ is separated since $E$ is a vector bundle.
\end{proof}
\comment{
\begin{proof}
By \cref{cor:proet-main}, our task is to prove that the cofibre of the morphism
\begin{equation*}
R\Gamma(X^{\HT, \dagger, +}, T_{\HT, \dagger, +}(E))\rightarrow R\Gamma(X^{\HT, \dagger}, T_{\HT, \dagger}(E))
\end{equation*}
is concentrated in degrees at least $i+1$. However, by \cref{cor:htdr-perfequiv}, we know that
\begin{equation*}
R\Gamma(X^{\HT, \dagger, +}, T_{\HT, \dagger, +}(E))\cong R\Gamma(X^{\dR, +}, T_{\dR, +}(E))
\end{equation*}
and we also have
\begin{equation*}
R\Gamma(X^{\HT, \dagger}, T_{\HT, \dagger}(E))\cong R\Gamma_\proet(X, T_{\HT, \dagger}(E))\;,
\end{equation*}
where we abuse notation by also using $T_{\HT, \dagger}(E)$ to denote the corresponding $\mathbb{B}_\dR^{+, \dagger}$-local system on the proétale site of $X^\diamond$. Now we are done once we prove the following proposition.
\end{proof}

\begin{prop}
Let $X$ be a quasicompact smooth partially proper rigid space over $\Q_p$ and let $(\cal{E}, \nabla, \Fil^\bullet)$ be a filtered vector bundle on $X$ with a flat connection satisfying Griffiths transversality with associated $\mathbb{B}_\dR^{+, \dagger}$-local system $\mathbb{M}$ on $X_\proet$. If the Hodge--Tate weights of $(\cal{E}, \nabla, \Fil^\bullet)$ are all at most $-i$ for some $i\geq 0$, then the cofibre of the natural map
\begin{equation*}
\Fil^0_\Hod R\Gamma_\dR(X, \cal{E})\rightarrow R\Gamma_\proet(X, \mathbb{M})
\end{equation*}
is concentrated in degrees at least $i+1$.
\end{prop}
\begin{proof}
Recall that $\mathbb{M}=\Fil^0(\mathcal{E}\tensor_{\O_X} \O\mathbb{B}_\dR^\dagger)^{\nabla=0}$ and hence $\mathbb{M}$ itself is equipped with a decreasing filtration $\Fil^\bullet \mathbb{M}$ which stabilises from filtration degree zero onwards. Moreover, recall that the map in question comes from the map
\begin{equation*}
\DR(\cal{E})\rightarrow \DR(\cal{E})\tensor_{\O_X} \O\mathbb{B}_\dR^\dagger\cong \mathbb{M}\tensor_{\mathbb{B}_\dR^{+, \dagger}} \mathbb{B}_\dR^\dagger\;,
\end{equation*}
where $\DR(\cal{E})$ denotes the de Rham complex associated to $\cal{E}$ equipped with its natural filtration
\begin{equation*}
\Fil^\bullet\DR(\cal{E})=(\Fil^\bullet\cal{E}\xrightarrow{\nabla} \Fil^{\bullet -1}\cal{E}\tensor_{\O_X} \Omega_X^1\xrightarrow{\nabla} \Fil^{\bullet-2}\cal{E}\tensor_{\O_X}\Omega_X^2\xrightarrow{\nabla}\dots)
\end{equation*}
and the tensor product $\DR(\cal{E})\tensor_{\O_X} \O\mathbb{B}_\dR^\dagger$ refers to the complex
\begin{equation*}
(\cal{E}\tensor_{\O_X} \O\mathbb{B}_\dR^\dagger\xrightarrow{\nabla} \cal{E}\tensor_{\O_X} \O\mathbb{B}_\dR^\dagger\tensor_{\O_X} \Omega_X^1\xrightarrow{\nabla} \dots)\;.
\end{equation*}
Also recall that the isomorphism above respects the filtrations; here, $\mathbb{B}_\dR^\dagger$ and $\O\mathbb{B}_\dR^\dagger$ are of course equipped with the $t$-adic filtrations -- however, we caution the reader that, in contrast to the situation for $\mathbb{B}_\dR$ and $\O\mathbb{B}_\dR$, these filtrations are \emph{not} complete. Passing to $\Fil^0$, the map from the proposition is obtained by taking cohomology of the map
\begin{equation*}
\Fil^0\DR(\cal{E})\rightarrow \Fil^0(\DR(\cal{E})\tensor_{\O_X} \O\mathbb{B}_\dR^\dagger)\cong \Fil^0(\mathbb{M}\tensor_{\mathbb{B}_\dR^{+, \dagger}} \mathbb{B}_\dR^\dagger)\cong \mathbb{M}\;.
\end{equation*}

Now first note that the map 
\begin{equation*}
\Fil^0(\DR(\cal{E})\tensor_{\O_X} \O\mathbb{B}_\dR^\dagger)\rightarrow \Fil^0(\DR(\cal{E})\tensor_{\O_X} \O\mathbb{B}_\dR)
\end{equation*}
induces an isomorphism in cohomology. Indeed, as the filtration on $\DR(\cal{E})$ is a finite filtration by subbundles, an easy argument similar to \cite[Tag 0944]{Stacks} reduces us to the claim that $t^n\O\mathbb{B}_\dR^{\dagger, +}\rightarrow t^n\O\mathbb{B}_\dR^+$ becomes an isomorphism after (derived) pushforward along $\lambda: X_\proet\rightarrow X_\an$. As
\begin{equation*}
t^n\O\mathbb{B}_\dR^{\dagger, +}/t^{n+1}\O\mathbb{B}_\dR^{\dagger, +}\cong \widehat{\O}(n)\cong t^n\O\mathbb{B}_\dR^+/t^{n+1}\O\mathbb{B}_\dR^+\;,
\end{equation*}
we are furthermore reduced to the case $n=0$. However, in this case, we have 
\begin{equation*}
R^i\lambda_*\O\mathbb{B}_\dR^{\dagger, +}=\begin{cases} \O_X\;, \hspace{0.3cm}\text{$i=0, 1$} \\
\mathrlap{0}\hphantom{\O_X}\;, \hspace{0.3cm}\text{$i\geq 2$}\;,
\end{cases}
\end{equation*}
which is a variant of \cite[Thm.\ H]{Wiersig2}, and the same result holds for $\O\mathbb{B}_\dR^+$ by reducing to associated gradeds of the $t$-adic filtration and then using \cite[Prop.\ 6.16]{PAdicHodgeTheory}.

Thus, we may equivalently consider the map
\begin{equation*}
\Fil^0\DR(\cal{E})\rightarrow\Fil^0(\DR(\cal{E})\tensor_{\O_X} \O\mathbb{B}_\dR)\;.
\end{equation*}
What we have gained by this maneuver is that the right-hand side is now complete with respect to its filtration and hence we may argue on the level of associated gradeds. In other words, to prove our claim, it will be enough to show that, for each $n\geq 0$, the proétale cohomology of the cofibre of the map
\begin{equation*}
\gr^n\DR(\cal{E})\rightarrow \gr^n(\DR(\cal{E})\tensor_{\O_X} \O\mathbb{B}_\dR)\cong \bigoplus_{j\in\Z} \gr^j\DR(\cal{E})\tensor_{\O_X} \gr^{n-j}\O\mathbb{B}_\dR
\end{equation*}
is concentrated in cohomological degrees at least $i+1$. For this, first observe that only the summand $j=n$ contributes to the cohomology of the right-hand side: indeed, as $\gr^j\DR(\cal{E})$ is a complex of locally free $\O_X$-modules, this follows from the fact that the cohomology of $\gr^\bullet\O\mathbb{B}_\dR$ on $X_\proet$ vanishes for $\bullet\neq 0$, see \cite[Prop.\ 6.16]{PAdicHodgeTheory}. Moreover, the assumption on the Hodge--Tate weights guarantees that $\gr^n\DR(\cal{E})$ is concentrated in cohomological degrees at least $i+n$. Thus, really the only case of interest is $n=0$ and there it is enough to show that the proétale cohomology of the cofibre of the map
\begin{equation*}
\O_X\rightarrow\widehat{\O}_X\cong \gr^0\O\mathbb{B}_\dR
\end{equation*}
is concentrated in degrees at least $1$, again since $\gr^0\DR(\cal{E})$ is a complex of locally free $\O_X$-modules.

In other words, we need to show that $H^0_\proet(X, \O_X)=H^0_\proet(X, \widehat{\O}_X)$ and that $H^1_\proet(X, \O_X)$ injects into $H^1_\proet(X, \widehat{\O}_X)$. However, recalling that $H^n_\proet(X, \O_X)=H^n_\et(X, \O_X)$ for all $n\geq 0$, see \cite[Cor.\ 3.17]{PAdicHodgeTheory}, this just follows from $\nu_*\widehat{\O}_X\cong\O_X$ for $\nu: X_\proet\rightarrow X_\et$, see \cite[Cor.\ 6.19]{PAdicHodgeTheory}, using the Grothendieck spectral sequence.
\end{proof}
}

\begin{rem}
Note that \cref{cor:proet-main} in some sense improves upon the ``classical'' comparison from \cref{thm:proet-main2} by also describing the failure of the map $R\Gamma(X^\Syn, \O\{i\})\rightarrow R\Gamma_\proet(X, \Q_p(i))$ to be an isomorphism in higher degrees: Namely, this is exactly measured by the failure of $R\Gamma(X^{\HT, \dagger, +}, \O\{i\})\rightarrow R\Gamma(X^{\HT, \dagger}, \O\{i\})$ to be an isomorphism. To our knowledge, a statement of this sort has not previously appeared in the literature.
\end{rem}

Of course, combining \cref{thm:proet-main2} and \cref{cor:hkcomp-main}, we immediately obtain the following extension of Colmez--Nizio{\l}'s ``basic comparison theorem'', see \cite[Thm.\ 1.3]{padicComparisons}, to vector bundle $F$-gauge coefficients:

\begin{cor}
Let $X$ be a smooth partially proper rigid space over $\Q_p$. If $E\in\Vect(X^\Syn)$ is a vector bundle analytic $F$-gauge with Hodge--Tate weights all at most $-i$ for some $i\geq 0$, then there is a natural map
\begin{equation*}
\fib(R\Gamma(X^\HK, T_\HK(E))\rightarrow R\Gamma(X^\dR, T_\dR(E))/R\Gamma(X^{\dR, +}, T_{\dR, +}(E)))\rightarrow R\Gamma(\FF_{X^\diamond}, T_\FF(E))
\end{equation*}
which induces an isomorphism on $\tau^{\leq i}$ and an injection on $H^{i+1}$. In particular, for $E=\O\{i\}$, we obtain
\begin{equation*}
\tau^{\leq i} R\Gamma_\proet(X, \Q_p(i))\cong \tau^{\leq i}\fib(R\Gamma_\HK(X)^{N=0, \phi=p^i}\rightarrow R\Gamma_\dR(X)/\Fil^i_\Hod R\Gamma_\dR(X))\;.
\end{equation*} 
\end{cor}

\subsubsection{Vector bundles on $X^\Syn$ are de Rham bundles on $\FF_{X^\diamond}$}

We can also use \cref{cor:proet-main} to fully describe the category of vector bundles on $X^\Syn$ for any smooth partially proper rigid space $X$ over $\Q_p$. For this, we introduce the following terminology:

\begin{defi}
Let $X$ be a smooth partially proper rigid space over $\Q_p$. A vector bundle on the Fargues--Fontaine curve $\FF_{X^\diamond}$ of $X^\diamond$ is called \emph{de Rham} if its restriction to a $\mathbb{B}_\dR^+$-local system on $X_\proet$ is generically flat in the sense of \cref{defi:htdr-genericallyflat}. We denote the full subcategory of $\Vect(\FF_{X^\diamond})$ consisting of de Rham bundles by $\Vect^\dR(\FF_{X^\diamond})$.
\end{defi}

\begin{rem}
Note that this recovers \cite[Def.\ 10.4.5]{FarguesFontaine} in the case $X=\GSpec\Q_p$.
\end{rem}

Now we can prove:

\begin{thm}
\label{thm:proet-vectxsyn}
Let $X$ be a smooth partially proper rigid space over $\Q_p$. Then there is an equivalence of categories
\begin{equation*}
\Vect(X^\Syn)\cong \Vect^\dR(\FF_{X^\diamond})\;.
\end{equation*}
\end{thm}
\begin{proof}
Recall from \cref{thm:htdr-vectxhtdagger} that the pullback functor $\Vect(X^{\HT, \dagger, +})\rightarrow\Vect(X^{\HT, \dagger})$ is fully faithful and its essential image is given by those vector bundles on $X^{\HT, \dagger}$ whose associated $\mathbb{B}_\dR^+$-local system is generically flat. Now we are done by the variant of \cref{cor:proet-main} for vector bundles using the equivalence $\Vect(X^{\Div^1})\cong \Vect(\FF_{X^\diamond})$ from \TODO{Reference!!}.
\end{proof}

\begin{rem}
Note that the case $X=\GSpec\Q_p$ of the above statement is already covered by \cref{thm:bk-mainqp}, but the proof we have given here is diametrically opposite to the one given in loc.\ cit.: Indeed, here we make use of the connection between $X^\Syn$ and $X^{\Div^1}$, i.e.\ the ``étale aspect'' of $X^\Syn$, while the proof of loc.\ cit.\ exploits the connection between $X^\Syn$ and $X^\HK$, i.e.\ the ``de Rham aspect'' of $X^\Syn$.
\end{rem}

\begin{rem}
Another view on the two proofs we have given for the above result in the case $X=\GSpec\Q_p$ is the following: While the proof via $\Q_p^{\Div^1}$ directly shows that vector bundles on $\Q_p^\Syn$ are de Rham bundles on $\FF_{\Q_p}$, the proof of \cref{thm:bk-mainqp} rather shows that vector bundles on $\Q_p^\Syn$ are the same as filtered $(\phi, N, G_{\Q_p})$-modules, i.e.\ potentially semistable bundles on $\FF_{\Q_p}$. In other words, putting these two proofs together gives a new proof of the theorem that ``de Rham $=$ potentially semistable'' of Berger. Note, however, that we need the equivalence between vector bundles on $\Q_p^\HK$ and $(\phi, N, G_{\Q_p})$-modules from \cite[Thm.\ 7.1.1]{dRFF} as an input; this should probably be seen as analogous to Berger's reduction of the theorem to Crew's conjecture about $p$-adic differential equations.
\end{rem}

Developing the previous remarks further, we might want to compare the results of \cref{thm:proet-vectxsyn} and \cref{cor:hkcomp-main} more generally for smooth partially proper rigid spaces $X$ over $\Q_p$. This immediately yields the following result, which was already sketched in \cite[Rem.\ 7.5.7]{dRFF} and, as explained in loc.\ cit., in particular also implies Shimizu's theorem about de Rham local systems from \cite[Thm.\ 1.1]{Shimizu}:

\begin{cor}
Let $X$ be a smooth partially proper rigid space over $\Q_p$. Then there is an equivalence of categories
\begin{equation*}
\Vect^\dR(\FF_{X^\diamond})\cong \Vect(X^\HK)\times_{\Vect(X^\dR)} \Vect(X^{\dR, +})\;.
\end{equation*}
\end{cor}

\newpage

\subsection{Geometric syntomic cohomology}

In this section, we want to explain how to construct an analytic syntomification of rigid-analytic spaces over $\C_p$ ``relative to $\C_p$'' and record some of the analogues of the results of the previous sections in this setting. For this, consider the Gelfand stack
\begin{equation*}
\Fil Y_{\C_p}\coloneqq (\{ut=\phi^{-1}(\xi)\}\subseteq Y_{\C_p}\times \ol{\DD}_+\times\ol{\DD}_-)/\G_m(1)\;,
\end{equation*}
where $u$ denotes the coordinate on $\ol{\DD}_+$, $t$ denotes the coordinate on $\ol{\DD}_-$ and $\xi=p-[p^\flat]$, as usual. Then note that this has a map to $\ol{\DD}/\G_m(1)\times (\ol{\DD}/\G_m(1))^\dR$ via $t$ and $u$; moreover, it also maps to $Y_{\C_p}$ and hence to $\C_p^\prism$ by \cref{prop:defis-xprismdr}. Overall, we therefore obtain a map
\begin{equation*}
\Fil Y_{\C_p}\rightarrow \Q_p^\N\times_{\Q_p^\prism} \C_p^\prism\;.
\end{equation*}
Over $\Q_p^\N$, the stack $\Fil Y_{\C_p}$ also maps to $\C_p^{\Cone}$: Indeed, we have
\begin{equation*}
(\{t=0\}\subseteq \Fil Y_{\C_p})\cong (\GSpec\C_p\times \ol{\DD}_+)/\G_m(1)
\end{equation*}
and this has a natural map to $\GSpec\C_p$; note that we have used the isomorphism 
\begin{equation*}
(\{\phi^{-1}(\xi)=0\}\subseteq Y_{\C_p})\cong (\{\xi=0\}\subseteq Y_{\C_p})
\end{equation*}
here, which is induced by Frobenius. Recalling that $\C_p^\N$ is the pullback of $\C_p^\prism$ along $\C_p^{\Cone}\rightarrow (\C_p^{\Cone})^\dR$ over $\Q_p^\N$, the preceding discussion yields a map
\begin{equation*}
\Fil Y_{\C_p}\rightarrow \C_p^\N\;.
\end{equation*}
With this, we can define analytic Nygaardification relative to $\C_p$ as follows:

\begin{defi}
Let $X$ be a Gelfand stack over $\C_p$. Its \emph{relative analytic Nygaardification} $X^{\N/\C_p}$ over $\C_p$ is defined as the pullback
\begin{equation*}
\begin{tikzcd}
X^{\N/\C_p}\ar[r]\ar[d] & \Fil Y_{\C_p}\ar[d] \\
X^\N\ar[r] & \C_p^\N\nospacepunct{\;.}
\end{tikzcd}
\end{equation*}
\end{defi}

As in the first section, let us discuss the restriction of $X^{\N/\C_p}$ to various interesting loci. First, note that both $\{|u|=1\}\subseteq \Fil Y_{\C_p}$ and $\{|t|=1\}\subseteq\Fil Y_{\C_p}$ are isomorphic to $Y_{\C_p}$. Recalling that $(X^\N)_{|u|=1}\cong (X^\N)_{|t|=1}\cong X^\prism$ for any Gelfand stack $X$, we find that there are isomorphisms
\begin{equation*}
(X^{\N/\C_p})_{|u|=1}\cong (X^{\N/\C_p})_{|t|=1}\cong X^{\prism/\C_p}
\end{equation*}
for any Gelfand stack $X$ over $\C_p$. Here, the last stack is defined as follows:

\begin{defi}
Let $X$ be a Gelfand stack over $\C_p$. Its \emph{relative analytic prismatisation} $X^{\prism/\C_p}$ over $\C_p$ is defined by the pullback diagram
\begin{equation*}
\begin{tikzcd}
X^{\prism/\C_p}\ar[r]\ar[d] & Y_{\C_p}\ar[d] \\
X^\prism\ar[r] & \C_p^\prism\nospacepunct{\;,}
\end{tikzcd}
\end{equation*}
where the map $Y_{\C_p}\rightarrow\C_p^\prism$ comes from \cref{prop:defis-xprismdr}. 
\end{defi}

Moreover, observe that, similarly as in \cref{prop:defis-utneq0}, there is an 
\begin{equation*}
(\{|ut|\neq 0\}\subseteq \Fil Y_{\C_p})\cong (Y_{\C_p}\setminus \{\phi^{-1}(\xi)=0\})\times [0, 1]
\end{equation*}
induced by the projection $\Fil Y_{\C_p}\rightarrow Y_{\C_p}$. The latter map induces a ``relative version'' $\pi: X^{\N/\C_p}\rightarrow X^{\prism/\C_p}$ of the structure map and then another application of \cref{prop:defis-utneq0} shows that
\begin{equation*}
(X^{\N/\C_p})_{|ut|\neq 0}\cong (X^{\prism/\C_p}\setminus X^{\dR/B_\dR^{+, \dagger}})\times [0, 1]
\end{equation*}
via $\pi$. Here, the stack $X^{\dR/B_\dR^{+, \dagger}}$ is the preimage of $X^\dR$ under the projection map $X^{\prism/\C_p}\rightarrow X^\prism$ and may alternatively defined as the pullback
\begin{equation*}
\begin{tikzcd}
X^{\dR/B_\dR^{+, \dagger}}\ar[r]\ar[d] & \GSpec B_\dR^{+, \dagger}\ar[d] \\
X^\dR\ar[r] & \C_p^\dR\nospacepunct{\;,}
\end{tikzcd}
\end{equation*}
where the vertical map on the right classifies the map $\C_p\xrightarrow{\cong} B_\dR^{+, \dagger}/\phi^{-1}(\xi)$; we call $X^{\dR/B_\dR^{+, \dagger}}$ the \emph{analytic $B_\dR^{+, \dagger}$-stack} of $X$ over $\C_p$.

Turning our attention to the locus $\{|u|=0\}$, note that
\begin{equation*}
(\{|u|=0\}\subseteq \Fil Y_{\C_p})\cong \GSpec(B_\dR^{+, \dagger}\langle t\rangle_{\leq 1}\{u\}^\dagger/(ut-\phi^{-1}(\xi)))/\G_m(1)\;;
\end{equation*}
we immediately warn the reader that the $t$ occurring here is \emph{not} the usual element $t\in B_\dR^{+, \dagger}$, but just an additional free variable coming from the factor $\ol{\DD}_-$ above. We conclude that restricting $X^{\N/\C_p}$ to the locus $\{|u|=0\}$ yields the following object:

\begin{defi}
Let $X$ be a Gelfand stack over $\C_p$. Its \emph{filtered analytic $B_\dR^{+, \dagger}$-stack} $X^{\dR, +/B_\dR^{+, \dagger}}$ is defined by the pullback diagram
\begin{equation*}
\begin{tikzcd}
X^{\dR, +/B_\dR^{+, \dagger}}\ar[r]\ar[d] & \GSpec(B_\dR^{+, \dagger}\langle t\rangle_{\leq 1}\{u\}^\dagger/(ut-\phi^{-1}(\xi)))/\G_m(1)\ar[d] \\
X^{\dR, +}\ar[r] & \C_p^{\dR, +}\nospacepunct{\;.}
\end{tikzcd}
\end{equation*}
\end{defi}

In the above definition, the vertical map on the right is explicitly given as follows: Given an $A$-point of the top right corner, i.e.\ a $B_\dR^{+, \dagger}$-algebra $A$ equipped with a normed generalised Cartier divisor $t: L\rightarrow A$ of norm at most $1$ and a map $u: A\rightarrow \Nil^\dagger(A)\tensor_A L$ such that the composition of $u$ and $t$ is multiplication by $\phi^{-1}(\xi)$, we obtain a map
\begin{equation*}
\C_p\cong B_\dR^{+, \dagger}/\phi^{-1}(\xi)\rightarrow A/^\mathbb{L}\phi^{-1}(\xi)\rightarrow \Cone(\Nil^\dagger(A)\tensor_A L\xrightarrow{t} A)\;,
\end{equation*}
i.e.\ an $A$-point of $\C_p^{\dR, +}$. 

\begin{rem}
Note that we do \emph{not} use the notation $X^{\dR/\C_p}$ or $X^{\dR, +/\C_p}$ for the two relative de Rham stacks introduced above. This is because this notation is already taken: There are natural maps $\GSpec\C_p\rightarrow\C_p^\dR$ and $\GSpec\C_p\times \ol{\DD}/\G_m(1)\rightarrow\C_p^{\dR, +}$ and then one defines $X^{\dR/\C_p}$ and $X^{\dR, +/\C_p}$ as the pullbacks
\begin{equation*}
\begin{tikzcd}
X^{\dR/\C_p}\ar[r]\ar[d] & \GSpec\C_p\ar[d] \\
X^\dR\ar[r] & \C_p^\dR\nospacepunct{\;,}
\end{tikzcd}
\hspace{1cm}
\begin{tikzcd}
X^{\dR, +/\C_p}\ar[r]\ar[d] & \GSpec\C_p\times \ol{\DD}/\G_m(1)\ar[d] \\
X^{\dR, +}\ar[r] & \C_p^{\dR, +}\nospacepunct{\;.}
\end{tikzcd}
\end{equation*}
By \cite[Prop.\ 5.2.1, Rem.\ 5.2.3]{dRFF}, these stacks compute (Hodge-filtered) de Rham cohomology of $X$ relative to $\C_p$ if $X$ is a smooth partially proper rigid space over $\C_p$.
\end{rem}

Let us also introduce the other ``relative stacks'' that will play a role in the sequel:

\begin{defi}
Let $X$ be a Gelfand stack over $\C_p$. Its \emph{relative Hyodo--Kato stack} $X^{\HK/\C_p}$ over $\C_p$ is defined as the pullback
\begin{equation*}
\begin{tikzcd}
X^{\HK/\C_p}\ar[r]\ar[d] & \FF_{\C_p}\ar[d] \\
X^\HK\ar[r] & \C_p^\HK\nospacepunct{\;,}
\end{tikzcd}
\end{equation*}
where the vertical map on the right is given by the projection $\FF_{\C_p}\rightarrow (\FF_{\C_p})^\dR=\C_p^\HK$.
\end{defi}

\begin{defi}
Let $X$ be a Gelfand stack over $\C_p$. The \emph{relative perfect prismatisation} $X^{\Div^1/\C_p}$ of $X$ over $\C_p$ is defined by the pullback diagram
\begin{equation*}
\begin{tikzcd}
X^{\Div^1/\C_p}\ar[r]\ar[d] & \FF_{\C_p}\ar[d] \\
X^{\Div^1}\ar[r] & \C_p^{\Div^1}\nospacepunct{\;,}
\end{tikzcd}
\end{equation*}
where the map on the right comes from combining \cref{prop:defis-xprismdr} and \cref{prop:defis-xprismxdiv1}.
\end{defi}

Finally, we have to introduce the analytic syntomification relative to $\C_p$. For this, recall that we have already seen above that, for any Gelfand stack $X$ over $\C_p$, there are two closed embeddings
\begin{equation*}
j_\dR, j_\HT: X^{\prism/\C_p}\rightarrow X^{\N/\C_p}
\end{equation*}
obtained by restricting to the locus $\{|t|=1\}$ and $\{|u|=1\}$, respectively. Unsurprisingly, we thus define:

\begin{defi}
Let $X$ be a Gelfand stack over $\C_p$. The \emph{relative analytic syntomification} $X^{\Syn/\C_p}$ of $X$ over $\C_p$ is defined as the coequaliser
\begin{equation*}
\begin{tikzcd}
X^{\prism/\C_p}\ar[r,shift left=.75ex,"j_\HT"]\ar[r,shift right=.75ex,swap,"j_\dR"] & X^{\N/\C_p}\ar[r] & X^{\Syn/\C_p}\;.
\end{tikzcd}
\end{equation*}
\end{defi}

Summarising the state of affairs, we have seen that all the geometric features of $X^\N$ discussed in Section \ref{sect:defis} are retained upon passing to the relative version $X^{\N/\C_p}$ as long as one also replaces the stacks $X^\dR, X^\prism, \dots$ by suitable relative versions as well.

\subsubsection{The $B_\dR^{+, \dagger}$-stack and $B_\dR^+$-cohomology}

We start by showing that the cohomology of the stack $X^{\dR/B_\dR^{+, \dagger}}$ agrees with the $B_\dR^+$-cohomology of any smooth partially proper rigid space $X$ over $\C_p$ as defined in \cite[§3]{BasicComparison}, \cite{GuoBdR} or \cite[§13]{IntegralpAdicHT}. Namely, our main goal is to prove the following:

\begin{thm}
\label{thm:geom-bdr+main}
Let $X$ be a smooth partially proper rigid space over $\C_p$. Then there is a natural isomorphism
\begin{equation*}
R\Gamma(X^{\dR, +/B_\dR^{+, \dagger}}, \O(-i))\tensor_{B_\dR^{+, \dagger}} B_\dR^+\cong \Fil^i R\Gamma(X/B_\dR^+)\;.
\end{equation*}
In particular, we have
\begin{equation*}
R\Gamma(X^{\dR/B_\dR^{+, \dagger}}, \O)\tensor_{B_\dR^{+, \dagger}} B_\dR^+\cong R\Gamma(X/B_\dR^+)\;.
\end{equation*}
\end{thm}

Before we begin the proof, let us fix the following notation: Contrary to our conventions in the previous sections, we denote the universal generalised Cartier divisor on $\ol{\DD}_+/\G_m(1)$ by $s: \O(-1)\rightarrow\O$; we do this to clearly distinguish the universal section $s$ from the usual element $t\in B_\dR^{+, \dagger}$, which will also play a role in the following. With this notation, the stack $X^{\dR, +/B_\dR^{+, \dagger}}$ is defined as the base change of $X^{\dR, +}$ along the map
\begin{equation*}
\GSpec(B_\dR^{+, \dagger}\langle s\rangle_{\leq 1}\{r\}^\dagger/(rs-t))/\G_m(1)\rightarrow \C_p^{\dR, +}\;,
\end{equation*}
where the variable $r$ is related to the variable $u$ from above via the formula $r=t\phi^{-1}(\xi)^{-1}u$; indeed, note that $t$ and $\phi^{-1}(\xi)$ differ by a unit in $B_\dR^{+, \dagger}$.

Further, observe that the base change of the displayed map above to $\C_p^\Hod$ is given by
\begin{equation*}
\GSpec(\C_p\{r\}^\dagger)/\G_m(1)\rightarrow \C_p^\Hod\;,
\end{equation*}
where one should think of $r$ being ``Tate-twisted'' by $1$ if one wants to keep track of the Galois actions, i.e.\ $G_{\Q_p}$ acts on $r$ via the cyclotomic character, and the $\G_m(1)$-action on $r$ is by division. For any Gelfand stack $X$ over $\C_p$, let $X^{\Hod/B_\dR^{+, \dagger}}$ denote the base change of $X^{\dR, +/B_\dR^{+, \dagger}}$ to $\C_p^\Hod$.

\begin{prop}
\label{prop:geom-xhod}
Let $X$ be a rigid smooth derived Berkovich space over $\C_p$. Then 
\begin{equation*}
X^{\Hod/B_\dR^{+, \dagger}}\cong \left(X\times \GSpec(\C_p\{r\}^\dagger)/\G_m(1)\right)\big/\,\cal{T}_{X/\C_p}^\dagger(-1)\;,
\end{equation*}
where the action of $\cal{T}_{X/\C_p}^\dagger(-1)$ is trivial and $\cal{T}_{X/\C_p}$ denotes $\GSpec_X \Sym_X^\bullet \mathbb{L}_{X/\C_p}$. In particular, if $X$ is a smooth partially proper rigid space over $\C_p$, the graded $\C_p$-vector space obtained by taking cohomology of (twists of) the structure sheaf on $X^{\Hod/B_\dR^{+, \dagger}}$ is given by
\begin{equation*}
\gr^\bullet R\Gamma(X^{\Hod/B_\dR^{+, \dagger}})=\bigoplus_{i\leq\bullet} R\Gamma(X, \Omega_{X/\C_p}^i)(\bullet-i)[-i]\;,
\end{equation*}
where the twists are now Tate twists.
\end{prop}
\begin{proof}
Our assumption implies that $X$ is $\dagger$-formally smooth over $\C_p$ and then the claim is immediate by derived deformation theory, see also \cite[Prop.\ 5.2.3.(2).(b)]{dRStack}. The second part now follows by Cartier duality between $\cal{T}_{X/\C_p}^\dagger(-1)$ and $\V(\Omega^1_{X/\C_p}(1))$.
\end{proof}

\begin{proof}[Proof of \cref{thm:geom-bdr+main}]
Since both sides of the claimed isomorphism satisfy analytic descent, we may without loss of generality assume that $X$ is affine. Then let us write $\Fil^\bullet R\Gamma(X^{\dR, +/B_\dR^{+, \dagger}})$ for the filtered object obtained by taking cohomology of the $\O(-i)$ on $X^{\dR, +/B_\dR^{+, \dagger}}$, with the transition maps being given by multiplication by $s$. 

We first construct a natural filtered morphism
\begin{equation}
\label{eq:geom-bdr+map}
\Fil^\bullet R\Gamma(X^{\dR, +/B_\dR^{+, \dagger}})\tensor_{B_\dR^{+, \dagger}} B_\dR^+\rightarrow R\Gamma_\proet (X, \Fil^\bullet\mathbb{B}_\dR^+)\;,
\end{equation}
where the filtration on $\mathbb{B}_\dR^+$ is of course the $t$-adic filtration, as follows: For any affinoid perfectoid $C$ over $X$, we have a map
\begin{equation*}
\GSpec(\mathbb{B}_\dR^{+, \dagger}(C)\langle s\rangle_{\leq 1}\{r\}^\dagger/(rs-t))/\G_m(1)\rightarrow C^{\dR, +}
\end{equation*}
defined in the same way as for $C=\C_p$ above, and since there is also a natural map $B_\dR^{+, \dagger}\rightarrow \mathbb{B}_\dR^{+, \dagger}(C)$, we even obtain a map to $C^{\dR, +/B_\dR^{+, \dagger}}$. Now observe that the filtered object obtained from cohomology of
\begin{equation*}
\GSpec(\mathbb{B}_\dR^{+, \dagger}(C)\langle s\rangle_{\leq 1}\{r\}^\dagger/(rs-t))/\G_m(1)
\end{equation*}
is exactly $\Fil^\bullet \mathbb{B}_\dR^{+, \dagger}(C)$ and thus, by functoriality of the filtered $B_\dR^{+, \dagger}$-stack, we obtain a filtered morphism
\begin{equation*}
\Fil^\bullet R\Gamma(X^{\dR, +/B_\dR^{+, \dagger}})\tensor_{B_\dR^{+, \dagger}} B_\dR^+\rightarrow \lim_{\GSpec C\rightarrow X} \Fil^\bullet \mathbb{B}_\dR^+(C)\cong R\Gamma_\proet(X, \Fil^\bullet \mathbb{B}_\dR^+)\;,
\end{equation*}
as desired. 

Recalling that we have assumed that $X$ is affine, we observe that the filtered object $\Fil^\bullet R\Gamma(X^{\dR, +/B_\dR^{+, \dagger}})$ is connective for the Beilinson $t$-structure as its graded pieces are described by \cref{prop:geom-xhod}. If we can show that the map 
\begin{equation}
\label{eq:geom-bdr+maptrunc}
\Fil^\bullet R\Gamma(X^{\dR, +/B_\dR^{+, \dagger}})\tensor_{B_\dR^{+, \dagger}} B_\dR^+\rightarrow \tau^{\leq 0}_{\mathrm{Beil}}R\Gamma_\proet (X, \Fil^\bullet \mathbb{B}_\dR^+)
\end{equation}
induced by (\ref{eq:geom-bdr+map}) is an isomorphism, the claim follows by \cite[Prop.\ 2.42]{BoscopAdic}. For this, first observe that the above is an isomorphism on graded pieces since
\begin{equation*}
\gr^i R\Gamma(X^{\dR, +/B_\dR^{+, \dagger}})\cong \tau^{\leq i} \gr^i R\Gamma_\proet (X, \Fil^\bullet \mathbb{B}_\dR^+)
\end{equation*}
for all $i\in\Z$ by combining \cref{prop:geom-xhod} with the fact that 
\begin{equation*}
R\Gamma_\proet(X, t^j\mathbb{B}_\dR^+/t^{j+1}\mathbb{B}_\dR^+)\cong \bigoplus_{n\geq 0} \Omega^n_X(j-n)[-n]\;,
\end{equation*}
see \cite[Prop.\ 3.23]{ScholzeSurvey}. Noting that all filtration steps of $R\Gamma_\proet(X, \Fil^\bullet\mathbb{B}_\dR^+)$ are cohomologically bounded above, see \cite[Prop.\ 2.40]{BoscopAdic}, we see that
\begin{equation*}
\tau^{\leq 0}_{\mathrm{Beil}}R\Gamma_\proet (X, \Fil^\bullet \mathbb{B}_\dR^+)\rightarrow R\Gamma_\proet (X, \Fil^\bullet \mathbb{B}_\dR^+)
\end{equation*}
is an isomorphism for $\bullet\gg 0$ and hence the target of (\ref{eq:geom-bdr+maptrunc}) is complete for the filtration as the same is true for $\Fil^\bullet \mathbb{B}_\dR^+$. It thus remains to show that the source is complete for the filtration as well. 

For this, after possibly further analytically localising on $X$, we may assume that $X$ admits a smooth lift $\widetilde{X}$ to $B_\dR^{+, \dagger}$ which comes equipped with a rigid étale map $\widetilde{X}\rightarrow \overcirc{\DD}^n_{B_\dR^{+, \dagger}}$ for some $n$. Since rigid étale maps are $\dagger$-formally étale, this yields a cartesian diagram
\begin{equation*}
\begin{tikzcd}
\widetilde{X}\times_{\GSpec B_\dR^{+, \dagger}} \GSpec(B_\dR^{+, \dagger}\langle s\rangle_{\leq 1}\{r\}^\dagger/(rs-t))/\G_m(1)\ar[r]\ar[d] & X^{\dR, +/B_\dR^{+, \dagger}}\ar[d] \\
\overcirc{\DD}^n_{B_\dR^{+, \dagger}}\times_{\GSpec B_\dR^{+, \dagger}} \GSpec(B_\dR^{+, \dagger}\langle s\rangle_{\leq 1}\{r\}^\dagger/(rs-t))/\G_m(1)\ar[r] & (\overcirc{\DD}^n_{\C_p})^{\dR, +/B_\dR^{+, \dagger}}\nospacepunct{\;.}
\end{tikzcd}
\end{equation*}
Now the bottom map exhibits the target as the quotient of the source by the action of $\G_a^\dagger(-1)^n$ by addition, where the twist is by the tautological line bundle on $*/\G_m(1)$, and thus we conclude
\begin{equation*}
X^{\dR, +/B_\dR^{+, \dagger}}\cong \left(\widetilde{X}\times_{\GSpec B_\dR^{+, \dagger}} \GSpec(B_\dR^{+, \dagger}\langle s\rangle_{\leq 1}\{r\}^\dagger/(rs-t))/\G_m(1)\right)\big/\, \G_a^\dagger(-1)^n\;.
\end{equation*}
By Cartier duality for $\G_a^\dagger$ and using that the addition action of $\G_a^\dagger$ on $\overcirc{\DD}$ corresponds to the endomorphism of $\O_{\overcirc{\DD}}$ given by the derivative, we now see that
\begin{equation*}
\Fil^\bullet R\Gamma(X^{\dR, +/B_\dR^{+, \dagger}})\cong \Fil^\bullet_\Hod R\Gamma_\dR(\widetilde{X}/B_\dR^{+, \dagger})\tensor_{B_\dR^{+, \dagger}} \Fil^\bullet B_\dR^{+, \dagger}
\end{equation*}
and tensoring with $B_\dR^+$ yields
\begin{equation*}
\Fil^\bullet R\Gamma(X^{\dR, +/B_\dR^{+, \dagger}})\tensor_{B_\dR^{+, \dagger}} B_\dR^+\cong \Fil^\bullet_\Hod R\Gamma_\dR(\widetilde{X}/B_\dR^{+, \dagger})\tensor_{B_\dR^{+, \dagger}} \Fil^\bullet B_\dR^+\;.
\end{equation*}
Finally, as $\Fil^\bullet_\Hod R\Gamma_\dR(\widetilde{X}/B_\dR^{+, \dagger})$ vanishes for $\bullet\gg 0$, the claim now follows from completeness of the filtration $\Fil^\bullet B_\dR^+$.
\end{proof}

From the theorem above, we can easily recover the following properties of $B_\dR^+$-cohomology:

\begin{cor}
Let $X$ be a smooth partially proper rigid space over $\C_p$. Then there is a natural isomorphism
\begin{equation*}
R\Gamma(X/B_\dR^+)\tensor_{B_\dR^+} \C_p\cong R\Gamma_\dR(X/\C_p)
\end{equation*}
compatible with the filtrations. Moreover, if $X$ is qcqs and admits a model $\ol{X}$ over some finite extension $K$ of $\Q_p$, then
\begin{equation*}
R\Gamma(X/B_\dR^+)\cong R\Gamma_\dR(\ol{X})\tensor_K B_\dR^+\;,
\end{equation*}
again compatibly with the filtrations.
\end{cor}
\begin{proof}
For the first isomorphism, observe that the big rectangle and the lower square in the diagram
\begin{equation*}
\begin{tikzcd}
X^{\dR, +/\C_p}\ar[r]\ar[d] & \GSpec \C_p\times \ol{\DD}/\G_m(1)\ar[d] \\
X^{\dR, +/B_\dR^{+, \dagger}}\ar[r]\ar[d] & \GSpec(B_\dR^{+, \dagger}\langle s\rangle_{\leq 1}\{r\}^\dagger/(rs-t))/\G_m(1)\ar[d] \\
X^{\dR, +}\ar[r] & \C_p^{\dR, +}
\end{tikzcd}
\end{equation*}
are cartesian, hence the upper square is cartesian as well. Moreover, the map 
\begin{equation*}
\GSpec\C_p\times \ol{\DD}/\G_m(1)\rightarrow \GSpec(B_\dR^{+, \dagger}\langle s\rangle_{\leq 1}\{r\}^\dagger/(rs-t))/\G_m(1)
\end{equation*}
is cohomologically smooth since this property is local on the target and
\begin{equation*}
\GSpec\C_p\langle s\rangle_{\leq 1}\rightarrow \GSpec(B_\dR^{+, \dagger}\langle s\rangle_{\leq 1}\{r\}^\dagger/(rs-t))
\end{equation*}
is a complete intersection (obtained by modding out $r$). Thus, the top square in the diagram above satisfies base change by \cite[Lem.\ 4.5.13.(i)]{HeyerMann} and this implies the claim upon recalling that $X^{\dR, +/\C_p}$ calculates Hodge-filtered de Rham cohomology of $X$ relative to $\C_p$.

For the second claim, we argue similarly: Both the bottom and the top square in the diagram
\begin{equation*}
\begin{tikzcd}
X^{\dR, +/B_\dR^{+, \dagger}}\ar[r]\ar[d] & \GSpec(B_\dR^{+, \dagger}\langle s\rangle_{\leq 1}\{r\}^\dagger/(rs-t))/\G_m(1)\ar[d] \\
X^{\dR, +}\ar[r]\ar[d] & \C_p^{\dR, +}\ar[d] \\
\ol{X}^{\dR, +}\ar[r] & K^{\dR, +}
\end{tikzcd}
\end{equation*}
are cartesian and hence the same is true for the large rectangle. By a similar argument as in \cref{lem:geom-hkprim}, the map $\ol{X}^{\dR, +}\rightarrow K^{\dR, +}\cong \GSpec K\times \ol{\DD}/\G_m(1)$ is prim and thus the large rectangle satisfies base change by \cite[Lem.\ 4.5.13.(ii)]{HeyerMann}, which implies the claim.
\end{proof}

\subsubsection{The filtered Hyodo--Kato and the proétale comparison}

We now want to prove the geometric analogues of the comparison theorems from Sections \ref{sect:hkcomp} and \ref{sect:proet}. For this, we first note that the relative stacks introduced above are again nicely coverable.

\begin{lem}
Let $X$ be a Gelfand stack over $\C_p$ equipped with a map $X\rightarrow \ol{\DD}_{\C_p}^n$ which is the composition of a finite étale map with a closed localisation. Then $(X^{\N/\C_p})_{[r, s]}$ is nicely coverable for any $1<r<s<\infty$.
\end{lem}
\begin{proof}
By the argument from \cref{prop:perf-zpncover}, it suffices to prove that 
\begin{equation*}
\C_p^{\N/\C_p}=\Fil Y_{\C_p}=(\{ut=\phi^{-1}(\xi)\}\subseteq Y_{\C_p}\times \ol{\DD}_+\times\ol{\DD}_-)/\G_m(1)
\end{equation*}
is nicely coverable after restricting to some closed interval $[r, s]$. However, this is clear: $\Fil Y_{\C_p}$ is covered by $(\{ut=\phi^{-1}(\xi)\}\subseteq Y_{\C_p}\times \ol{\DD}_+\times\ol{\DD}_-)$ with the \v{C}ech nerve being
\begin{equation*}
(\{ut=\phi^{-1}(\xi)\}\subseteq Y_{\C_p}\times \ol{\DD}_+\times\ol{\DD}_-)\times \G_m(1)^\bullet
\end{equation*}
and $\G_m(1)$ is affine, static and flat over $\Q_p$. Thus, the claim reduces to the fact that 
\begin{equation*}
(\{ut=\phi^{-1}(\xi)\}\subseteq Y_{\C_p}\times \ol{\DD}_+\times\ol{\DD}_-)
\end{equation*}
is affine after restricting to some closed interval $[r, s]$ and $n$-truncated for some $n$. (In fact, it is static and flat over $\Q_p$.)
\end{proof}

With the same argument as in Section \ref{sect:hkcomp}, we obtain realisation functors
\begin{equation*}
T_{\HK/\C_p}: \Perf(X^{\Syn/\C_p})\rightarrow \Perf(X^{\HK/\C_p})
\end{equation*}
as well as
\begin{equation*}
T_{\dR, +/B_\dR^{+, \dagger}}: \D(X^{\Syn/\C_p})\rightarrow \D(X^{\dR, +/B_\dR^{+, \dagger}})
\end{equation*}
and
\begin{equation*}
T_{\dR/B_\dR^{+, \dagger}}: \D(X^{\Syn/\C_p})\rightarrow \D(X^{\dR/B_\dR^{+, \dagger}})\;.
\end{equation*}
Then the argument from loc.\ cit.\ shows the following:

\begin{thm}
\label{thm:geom-hkcomp}
Let $X$ be a rigid smooth derived Berkovich space over $\C_p$. Then the functors $T_{\HK/\C_p}$ and $T_{\dR, +/B_\dR^{+, \dagger}}$ induce an equivalence
\begin{equation*}
\Perf(X^{\Syn/\C_p})\cong \Perf(X^{\HK/\C_p})\times_{\Perf(X^{\dR/B_\dR^{+, \dagger}})} \Perf(X^{\dR, +/B_\dR^{+, \dagger}})\;.
\end{equation*}
In particular, for any $E\in\Perf(X^{\Syn/\C_p})$, there is a pullback diagram
\begin{equation*}
\begin{tikzcd}
R\Gamma(X^{\Syn/\C_p}, E)\ar[r]\ar[d] & R\Gamma(X^{\HK/\C_p}, T_{\HK/\C_p}(E))\ar[d] \\
R\Gamma(X^{\dR, +/B_\dR^{+, \dagger}}, T_{\dR, +/B_\dR^{+, \dagger}}(E))\ar[r] & R\Gamma(X^{\dR/B_\dR^{+, \dagger}}, T_{\dR/B_\dR^{+, \dagger}}(E))\nospacepunct{\;.}
\end{tikzcd}
\end{equation*}
\end{thm}

We want to make the above theorem more explicit in the case $E=\O\{i\}$. For this, let us introduce the following notation.

\begin{defi}
Let $X$ be any Gelfand stack over $\C_p$. For any $i\in\Z$, we write
\begin{equation*}
R\Gamma_\Syn(X/\C_p, \Q_p(i))\coloneqq R\Gamma(X^{\Syn/\C_p}, \O\{i\})\;.
\end{equation*}
\end{defi}

\begin{lem}
\label{lem:geom-hkprim}
Let $X$ be a smooth partially proper qcqs rigid space over $\C_p$. Then the map $X^\HK\rightarrow \C_p^\HK$ is prim.
\end{lem}
\begin{proof}
By \cite[Lem.\ 4.5.8]{HeyerMann}, we may pass to a prim descendable cover of the source to check primness and thus, after passing to a finite strict closed cover of $X$, we may assume that $X$ admits an étale map to $\A^n_{\C_p}$ for some $n\geq 0$. However, the map $(\A^n_{\C_p})^\HK\rightarrow \C_p^\HK$ is even cohomologically proper by the arguments from the proofs of \cite[Lem.s 6.2.2, 6.2.3]{dRFF} and thus it remains to show that $X^\HK\rightarrow (\A^n_{\C_p})^\HK$ is prim. However, this follows from the fact that any proétale $\omega_1$-cocompact qcqs map induces a prim map on de Rham stacks, see the proof of \cite[Prop.\ 5.1.4]{dRFF}.
\end{proof}

\begin{cor}
\label{cor:geom-hkcomptrivcoeffs}
Let $X$ be a smooth partially proper qcqs rigid space over $\C_p$. Then there is a pullback diagram
\begin{equation*}
\begin{tikzcd}
R\Gamma_\Syn(X/\C_p, \Q_p(i))\ar[r]\ar[d] & (R\Gamma_\HK(X)\tensor_{\Q_p^\un} B_{\log})^{N=0, \phi=p^i}\ar[d] \\
\Fil^i R\Gamma(X/B_\dR^+)\ar[r] & R\Gamma(X/B_\dR^+)\nospacepunct{\;,}
\end{tikzcd}
\end{equation*}
where $B_{\log}$ is the period ring from \cite[Def.\ 10.3.1]{FarguesFontaine}.
\end{cor}
\begin{proof}
As cohomology on $X^{\dR, +/B_\dR^{+, \dagger}}$ recovers filtered $B_\dR^+$-cohomology by \cref{thm:geom-bdr+main}, it only remains to show that
\begin{equation*}
R\Gamma(X^{\HK/\C_p}, \O(i))\cong (R\Gamma_\HK(X)\tensor_{\Q_p^\un} B_{\log})^{N=0, \phi=p^i}\;.
\end{equation*}
For this, recall that $X^{\HK/\C_p}$ by definition sits in a pullback diagram
\begin{equation*}
\begin{tikzcd}
X^{\HK/\C_p}\ar[r, "f'"]\ar[d, "g'", swap] & \FF_{\C_p}\ar[d, "g"] \\
X^\HK\ar[r, "f"] & \C_p^\HK
\end{tikzcd}
\end{equation*}
and note that $X^\HK\rightarrow\C_p^\HK$ is prim by \cref{lem:geom-hkprim}. Thus, the above diagram satisfies base change and we obtain 
\begin{equation*}
f'_*\O(i)=f'_*{g'}^*\O(i)\cong g^*f_*\O(i)\;.
\end{equation*}
As $X$ is qcqs, the Hyodo--Kato cohomology of $X$ is finite and thus, using \cite[Rem.\ 6.1.3]{dRFF}, the pushforward $f_*\O(i)$ identifies with the $(\phi, N)$-module $R\Gamma_\HK(X)(i)$ over $\Q_p^\un$ via the equivalence on perfect complexes induced by pullback along the map
\begin{equation*}
\Psi: \C_p^\HK\rightarrow \ol{\F}_p^\HK/\,\V(-1)\cong \GSpec\Q_p^\un/\,\V(-1)\rtimes\phi^\Z
\end{equation*}
from \cite[Thm.\ 7.1.1]{dRFF}; here, the twist by $i$ indicates that we divide the Frobenius on $R\Gamma_\HK(X)$ by $p^i$. Thus, to understand $g^*f_*\O(i)$, we equivalently have to understand pullback along the map
\begin{equation*}
\Psi\circ g: \FF_{\C_p}\rightarrow \GSpec\Q_p^\un/\,\V(-1)\rtimes\phi^\Z\;.
\end{equation*}
However, by the construction in loc.\ cit., the $\V(-1)\rtimes\phi^\Z$-bundle over $\FF_{\C_p}$ this map classifies is just the pullback of the Fargues--Fontaine surface from \cite[Def.\ 10.3.7]{FarguesFontaine} along $Y_{\C_p}\rightarrow \FF_{\C_p}$. As the global sections of this bundle are given by $B_{\log}$, the claim follows.
\end{proof}

\begin{rem}
The above result shows that $R\Gamma_\Syn(X/\C_p, \Q_p(i))$ is \emph{not} the same as the geometric syntomic cohomology of Colmez--Nizio{\l} from \cite{BasicComparison}. Indeed, this is because, for their geometric syntomic cohmology $R\Gamma_{\Syn, \mathrm{CN}}(X/\C_p, \Q_p(i))$, we have a pullback diagram
\begin{equation*}
\begin{tikzcd}
R\Gamma_{\Syn, \mathrm{CN}}(X/\C_p, \Q_p(i))\ar[r]\ar[d] & (R\Gamma_\HK(X)\tensor_{\Q_p^\un} B_{\mathrm{st}}^+)^{N=0, \phi=p^i}\ar[d] \\
\Fil^i R\Gamma(X/B_\dR^+)\ar[r] & R\Gamma(X/B_\dR^+)\nospacepunct{\;,}
\end{tikzcd}
\end{equation*}
see Thm.\ 1.3 of loc.\ cit., and replacing $B_{\mathrm{st}}^+$ by $B_{\log}$ here \emph{does} in general change the $(N=0, \phi=p^i)$-eigenspace. The result above rather shows that our geometric syntomic cohomology $R\Gamma_\Syn(X/\C_p, \Q_p(i))$ recovers Bosco's \emph{syntomic Fargues--Fontaine cohomology} from \cite[§6]{BoscopAdic} by Thm.s 4.1 and 6.3 of loc.\ cit. and, by Ex.\ 6.36 of loc.\ cit., this already differs from Colmez--Nizio{\l} geometric syntomic cohomology in the case of the projective line.
\end{rem}

To get the comparison with proétale cohomology, we introduce the following relative versions of $X^{\HT, \dagger}$ and $X^{\HT, \dagger, +}$: For any Gelfand stack $X$ over $\C_p$, we write
\begin{equation*}
X^{\HT, \dagger, +/\C_p}\coloneqq (X^{\N/\C_p})_{|t|=0}\;, \hspace{0.5cm} X^{\HT, \dagger/\C_p}\coloneqq (X^{\N/\C_p})_{|t|=0, |u|=1}\;.
\end{equation*}
Then the same arguments as in Section \ref{sect:proet} first yields a realisation functor
\begin{equation*}
T_\FF: \Perf(X^{\Syn/\C_p})\rightarrow \Perf(X^{\Div^1/\C_p})\;.
\end{equation*}
Note that we do \emph{not} call this functor $T_{\FF/\C_p}$; this is because, by \TODO{Reference!!}, pullback along the natural map $\FF_{X^\diamond}\rightarrow X^{\Div^1/\C_p}$ induced by the map from \cref{prop:defis-xprismdr} and $\FF_{X^\diamond}\rightarrow\FF_{\C_p}$ induces an equivalence of categories
\begin{equation*}
\Perf(X^{\Div^1/\C_p})\cong \Perf(\FF_{X^\diamond})
\end{equation*}
and hence we can alternatively think of the functor $T_\FF$ as taking values in perfect complexes on $\FF_{X^\diamond}$. Finally, restriction furthermore induces functors
\begin{equation*}
T_{\HT, \dagger, +/\C_p}: \D(X^{\Syn/\C_p})\rightarrow \D(X^{\HT, \dagger, +/\C_p})
\end{equation*}
and 
\begin{equation*}
T_{\HT, \dagger/\C_p}: \D(X^{\Syn/\C_p})\rightarrow \D(X^{\HT, \dagger/\C_p})\;.
\end{equation*}
Then the same argument as in loc.\ cit.\ shows the following result.

\begin{thm}
\label{thm:geom-proet}
Let $X$ be a rigid smooth derived Berkovich space over $\C_p$. Then there is an equivalence of categories
\begin{equation*}
\Perf(X^{\Syn/\C_p})\cong \Perf(X^{\Div^1/\C_p})\times_{\Perf(X^{\HT, \dagger/\C_p})} \Perf(X^{\HT, \dagger, +/\C_p})
\end{equation*}
induced by the realisation functors $T_{\FF/\C_p}$ and $T_{\HT, \dagger, +/\C_p}$. In particular, for any $E\in\Perf(X^{\Syn/\C_p})$, there is a pullback diagram
\begin{equation*}
\begin{tikzcd}
R\Gamma(X^{\Syn/\C_p}, E)\ar[r]\ar[d] & R\Gamma(X^{\Div^1/\C_p}, T_{\FF/\C_p}(E))\ar[d] \\
R\Gamma(X^{\HT, \dagger, +/\C_p}, T_{\HT, \dagger, +/\C_p}(E))\ar[r] & R\Gamma(X^{\HT, \dagger/\C_p}, T_{\HT, \dagger/\C_p}(E))\nospacepunct{\;.}
\end{tikzcd}
\end{equation*}
\end{thm}

To extract the comparison with proétale cohomology from this, we again need a little more work. Namely, we need the following analogue of \cref{thm:htdr-xhtdrdagger+} and \cref{cor:htdr-cohomologyxhtdrdagger+}, which is proved the same way:

\begin{prop}
Let $X$ be a derived Berkovich space over $\C_p$ and assume that $X$ admits a lift $\widetilde{X}$ to $B_\dR^{+, \dagger}$ equipped with a rigid étale map $\widetilde{X}\rightarrow\ol{\T}^n_{B_\dR^{+, \dagger}}$ for some $n\geq 0$. 
\begin{enumerate}[label=(\arabic*)]
\item There is an isomorphism 
\begin{equation*}
X^{\HT, \dagger, +/\C_p}\cong \left(\widetilde{X}_\infty^\la\times_{\GSpec B_\dR^{+, \dagger}} \GSpec B_\dR^{+, \dagger}\langle u\rangle_{\leq 1}\{t\}^\dagger/(ut-\phi^{-1}(\xi))\right)\,\big/\, (\Z_p^\la)^n\coprod^\AbGrp_{(\G_a^\dagger)^n} (\G_a^\dagger)^n\;,
\end{equation*}
where the pushout is taken along the map $(\G_a^\dagger)^n\rightarrow (\G_a^\dagger)^n$ given by multiplication by $u$ and $\widetilde{X}_\infty^\la$ denotes the base change of $X$ along $\ol{\T}^{n, \la}_\infty\rightarrow\ol{\T}^n$. Here, the action of $(\Z_p^\la)^n\coprod^\AbGrp_{(\G_a^\dagger)^n} (\G_a^\dagger)^n$ on $\widetilde{X}_\infty^\la$ is obtained by pulling back the action on $\ol{\T}^{n, \la}_\infty$ given by
\begin{equation*}
(\ul{\theta}, \ul{v}).\ul{T}^{1/p^\infty}=\zeta_{p^\infty}^{\ul{\theta}}\exp(t(u\ul{\theta}+\ul{v})/p^\infty)\ul{T}^{1/p^\infty}\;,
\end{equation*}
where $\ul{T}=(T_1, \dots, T_n)$ are the coordinates on $\ol{\T}^n$.

\item A perfect complex on $X^{\HT, \dagger, +/\C_p}$ amounts to the following data:
\begin{enumerate}[label=(\roman*)]
\item A diagram
\begin{equation*}
\begin{tikzcd}
\dots \ar[r,shift left=.5ex,"u"]
  & \ar[l,shift left=.5ex, "t"] \Fil_{i-1} V \ar[r,shift left=.5ex,"u"] & \ar[l,shift left=.5ex, "t"] \Fil_i V \ar[r,shift left=.5ex,"u"] & \ar[l,shift left=.5ex, "t"] \Fil_{i+1} V\ar[r,shift left=.5ex,"u"] & \ar[l,shift left=.5ex, "t"] \dots
\end{tikzcd}
\end{equation*}
of perfect complexes over $\widetilde{X}_\infty^\la$ such that $ut=\phi^{-1}(\xi)$ and with the property that $t: \Fil_\bullet V\rightarrow \Fil_{\bullet-1} V$ becomes an isomorphism for $\bullet\ll 0$ and $u: \Fil_\bullet V\rightarrow\Fil_{\bullet+1} V$ becomes an isomorphism for $\bullet\gg 0$;
\item a $t$-connection 
\begin{equation*}
\nabla: V\rightarrow V\tensor_{\O_{\widetilde{X}_\infty^\la}} \Omega_{\widetilde{X}_\infty^\la/B_\dR^{+, \dagger}}^1
\end{equation*}
on $V$, i.e.\ a map satisfying the $t$-deformed Leibniz rule $\nabla(fm)=f\nabla(m)+t\cdot\mathrm{d}f$, which is Griffiths transversal with respect to the filtration $\Fil_\bullet V$;
\item a semilinear locally analytic $\Z_p^n$-action on $\Fil_\bullet V$ commuting with all the previous data;
\item an identification between the ``geometric Sen operator'' 
\begin{equation*}
\Theta^\geom: \Fil_i V\rightarrow \Fil_i V\tensor_{\O_{\widetilde{X}_\infty^\la}} \Omega_{\widetilde{X}_\infty^\la/B_\dR^{+, \dagger}}^1
\end{equation*}
coming from the Lie algebra action of $\Z_p^n$ and the $t$-connection $u\nabla$ for each $i$.
\end{enumerate}
Under this identification, restriction to $X^{\HT, \dagger/\C_p}$ corresponds to forgetting the $u$-filtration.

\item Under the above identifications, the cohomology of a perfect complex on $X^{\HT, \dagger, +/\C_p}$ is computed by the (underived) $\Z_p^n$-invariants of the coherent cohomology on $\widetilde{X}_\infty^\la$ of the complex
\begin{equation*}
\begin{tikzcd}
\Fil_0 V\ar[r, "\nabla"] & \Fil_{-1} V\tensor_{\O_{\widetilde{X}_\infty^\la}} \Omega_{\widetilde{X}_\infty^\la/B_\dR^{+, \dagger}}^1\ar[r, "\nabla"] & \Fil_{-2} V\tensor_{\O_{\widetilde{X}_\infty^\la}} \Omega_{\widetilde{X}_\infty^\la/B_\dR^{+, \dagger}}^2\ar[r, "\nabla"] & \dots
\end{tikzcd}
\end{equation*}
while the cohomology of its restriction to $X^{\HT, \dagger/\C_p}$ is given by the (underived) $\Z_p^n$-invariants of the coherent cohomology on $\widetilde{X}_\infty^\la$ of the complex
\begin{equation*}
\begin{tikzcd}
V\ar[r, "\nabla"] & V\tensor_{\O_{\widetilde{X}_\infty^\la}} \Omega_{\widetilde{X}_\infty^\la/B_\dR^{+, \dagger}}^1\ar[r, "\nabla"] & V\tensor_{\O_{\widetilde{X}_\infty^\la}} \Omega_{\widetilde{X}_\infty^\la/B_\dR^{+, \dagger}}^2\ar[r, "\nabla"] & \dots\nospacepunct{\;,}
\end{tikzcd}
\end{equation*}
where $V$ of course denotes the underlying unfiltered object of $\Fil_\bullet V$.
\end{enumerate}
\end{prop}

Using this, we can now deduce the comparison between syntomic and proétale cohomology from \cref{thm:geom-proet} by the same argument as in the arithmetic case.

\begin{cor}
Let $X$ be a smooth partially proper rigid space over $\C_p$. If $E\in\Vect(X^{\Syn/\C_p})$ is a vector bundle analytic $F$-gauge with Hodge--Tate weights all at most $-i$ for some $i\geq 0$, then the natural morphism
\begin{equation*}
R\Gamma(X^{\Syn/\C_p}, E)\rightarrow R\Gamma(\FF_{X^\diamond}, T_\FF(E))
\end{equation*}
is an isomorphism on $\tau^{\leq i}$ and induces an injection on $H^{i+1}$. In particular, for $E=\O\{i\}$, we obtain
\begin{equation*}
\tau^{\leq i} R\Gamma_\Syn(X/\C_p, \Q_p(i))\cong \tau^{\leq i}R\Gamma_\proet(X, \Q_p(i))\;.
\end{equation*}
\end{cor}

Again, combining the previous corollary with \cref{thm:geom-hkcomp}, we obtain a variant of Colmez--Nizio{\l}'s ``basic comparison theorem'' in the geometric case, see \cite[Thm.\ 1.1]{BasicComparison}, with vector bundle $F$-gauge coefficients. However, as already pointed out above, our version differs from the one in loc.\ cit.\ due to the appearance of the period ring $B_{\log}$ in place of $B_\st^+$ and should rather be compared to the result of Bosco from \cite[Thm.\ 7.2]{BoscopAdic}.

\begin{cor}
Let $X$ be a smooth partially proper rigid space over $\C_p$. If $E\in\Vect(X^{\Syn/\C_p})$ is a vector bundle analytic $F$-gauge with Hodge--Tate weights all at most $-i$ for some $i\geq 0$, then there is a natural map
\begin{equation*}
\begin{tikzcd}
\fib(R\Gamma(X^{\HK/\C_p}, T_{\HK/\C_p}(E))\rightarrow R\Gamma(X^{\dR/B_\dR^{+, \dagger}}, T_{\dR/B_\dR^{+, \dagger}}(E))/R\Gamma(X^{\dR, +/B_\dR^{+, \dagger}}, T_{\dR, +/B_\dR^{+, \dagger}}(E)))\ar[d] \\
R\Gamma(\FF_{X^\diamond}, T_\FF(E))
\end{tikzcd}
\end{equation*}
which induces an isomorphism on $\tau^{\leq i}$ and an injection on $H^{i+1}$. In particular, if $X$ is qcqs and $E=\O\{i\}$, we obtain
\begin{equation*}
\tau^{\leq i} R\Gamma_\proet(X, \Q_p(i))\cong \tau^{\leq i}\fib((R\Gamma_\HK(X)\tensor_{\Q_p^\un} B_{\log})^{N=0, \phi=p^i}\rightarrow R\Gamma(X/B_\dR^+)/\Fil^i R\Gamma(X/B_\dR^+))\;.
\end{equation*}
\end{cor}

\comment{
\newpage

\subsection{A conjecture of Clausen about syntomic cohomology of $\C_p$}
\label{sect:clausen}

Motivated by a calculation of the $K$-theory of the category of nuclear solid modules over $\C_p$ and using an expected Atiyah--Hirzebruch spectral sequence linking this ``analytic'' $K$-theory with analytic syntomic cohomology, Clausen has made the following conjecture about analytic syntomic cohomology of $\C_p$:

\begin{conj}[Clausen]
\label{conj:clausen}
For $n\geq 2$, we have
\begin{equation*}
R\Gamma(\C_p^\Syn, \O\{n\})\cong (B_\crys^+)^{\phi=p^n}/\Z_p(n)[-1]\;.
\end{equation*}
\end{conj}

Note that this is a conjecture about analytic syntomic cohomology with integral coefficients! As we have only discussed the case of rational coefficients, we can thus only hope to prove the above conjecture after inverting $p$. Indeed, using what we have established so far, we will prove:

\begin{thm}
\label{thm:clausen-main}
After inverting $p$, \cref{conj:clausen} is true. More specifically, we have
\begin{equation*}
R\Gamma(\C_p^\Syn, \O\{n\})\cong B_\dR^+/t^n[-1]
\end{equation*}
for any $n\geq 2$.
\end{thm}

Note that the above result indeed recovers the rational version of \cref{conj:clausen} due to $(B_\crys^+)^{\phi=p^n}/\Q_p(n)\cong B_\dR^+/t^n$ by one of the fundamental exact sequences of $p$-adic Hodge theory.

\begin{proof}[Proof of \cref{thm:clausen-main}]
By \cref{cor:hkcomp-main}, we know that there is a pullback square
\begin{equation*}
\begin{tikzcd}
R\Gamma(\C_p^\Syn, \O\{n\})\ar[r]\ar[d] & R\Gamma(\C_p^\HK, \O(n))\ar[d] \\
R\Gamma(\C_p^{\dR, +}, \O(-n))\ar[r] & R\Gamma(\C_p^\dR, \O)\nospacepunct{\;.}
\end{tikzcd}
\end{equation*}
Now recall from \cite[Rem.\ 1.2.6]{dRFF} that $\C_p^\dR\cong\GSpec\ol{\Q}_p$ and, moreover, that cohomology on $\C_p^\HK$ is the same as cohomology of $(\phi, N)$-modules over $\Q_p^\un$ by \cite[Thm.\ 7.1.1]{dRFF} with $\O(n)$ corresponding to the $(\phi, N)$-module $(\Q_p^\un, \phi=p^{-n}, N=0)$. As cohomology of a $(\phi, N)$-module $V$ is computed by the (total complex associated to the) double complex
\begin{equation*}
\begin{tikzcd}
V\ar[r, "N"] & V \\
V\ar[u, "\phi-1"]\ar[r, "N"] & V\ar[u, "1-p\phi", swap]\nospacepunct{\;,}
\end{tikzcd}
\end{equation*}
cf.\ \cite[§2.4]{CrystallineExt}, we see that the assumption $n\geq 2$ precisely guarantees that the cohomology of $\O(n)$ on $\C_p^\HK$ vanishes. Thus, we are done once we prove the following theorem.
\end{proof}

\begin{thm}
\label{thm:clausen-cpdr+}
We have
\begin{equation*}
R\Gamma(\C_p^{\dR, +}, \O(-n))\cong\begin{cases}
\mathrlap{\ol{\Q}_p}\hphantom{(B_\dR^+/t^n)/\ol{\Q}_p[-1]} \hspace{0.15cm},\hspace{0.15cm} n\leq 0 \\
(B_\dR^+/t^n)/\ol{\Q}_p[-1] \hspace{0.15cm},\hspace{0.15cm} n\geq 1\;,
\end{cases}
\end{equation*}
where $s: \O(-1)\rightarrow\O$ is the tautological generalised Cartier divisor on $\ol{\DD}/\G_m(1)$.
\end{thm}

\begin{rem}
Let us point out the following curious fact, for which we could not come up with a conceptual explanation: By the results of \cite[§8.2.1.1]{FarguesFontaine}, the cohomology groups in the above theorem agree with the ones of
\begin{equation*}
\colim_{E/\Q_p\text{ finite}} R\Gamma(\FF_{\C_p, E}, \O(-n))\;,
\end{equation*}
where $\FF_{\C_p, E}$ denotes the Fargues--Fontaine of $\C_p$ relative to the ground field $E$.
\end{rem}

The proof of this theorem will occupy the remainder of this section. Our argument will be quite straightforward in that we are going to give an explicit quotient stack surjecting onto $\C_p^{\dR, +}$ and then determine the full \v{C}ech nerve of this surjection to obtain explicit complexes computing the sought after cohomology. For this last step, the following property of $B_\dR^+$ will be crucial, see also \cite{ColmezBdR+}:

\begin{prop}
\label{prop:clausen-solidbdr+}
The map $B_\dR^+\rightarrow \C_p$ is the universal pro-infinitesimal thickening of $\C_p$ in solid rings, i.e.\ for any infinitesimal thickening $R\rightarrow\C_p$ of solid rings, the lifting problem
\begin{equation*}
\begin{tikzcd}
R\ar[r] & \C_p \\
& B_\dR^+\ar[u]\ar[lu, dotted, "\exists !"]
\end{tikzcd}
\end{equation*}
admits a unique solution. 
\end{prop}
\begin{proof}
This is a consequence of the fact that the solidification $L^{\solid}_{B_\dR^+/\Z}$ of the relative cotangent complex of $B_\dR^+$ over $\Z$ vanishes. To see the latter, following \cite{SolidBdR+}, first note that $L^{\solid}_{B_\dR^+/\Z}$ is derived $t$-complete: By the construction of free solid modules and algebras, the modules of differentials of the terms of the standard simplicial resolution of $B_\dR^+$ by solid polynomial algebras are derived $t$-complete since $B_\dR^+$ is. Then the claim follows from the fact that solid tensor products of derived $t$-complete objects are again derived $t$-complete. Overall, we thus conclude that $L^{\solid}_{B_\dR^+/\Z}=0$ may be checked after extending scalars to $\C_p$.

By the cotangent sequence for the composite $\Z\rightarrow B_\dR^+\rightarrow\C_p$, we are then reduced to showing that $L^{\solid}_{\C_p/\Z}\rightarrow L^{\solid}_{\C_p/B_\dR^+}\cong tB_\dR^+/t^2B_\dR^+[1]=\C_p\{1\}[1]$ is an isomorphism. By compatibility of cotangent complexes with localisation, we it suffices to show that $L^{\solid}_{\O_{\C_p}/\Z}[1]\cong \O_{\C_p}\{1\}[1]$ via the natural map and since $L^{\solid}_{\O_{\C_p}/\Z}$ is derived $p$-complete by a similar argument as above, it in fact suffices to check this after mod $p$ reduction, i.e.\ we are reduced to checking that $L^{\solid}_{(\O_{\C_p}/p)/\F_p}\cong \O_{\C_p}/p\{1\}[1]$. However, noting that $\O_{\C_p}/p\cong \O_{\C_p}^\flat/p^\flat$, this now follows from the fact that $L^{\solid}_{\O_{\C_p}^\flat/\F_p}=0$ using the cotangent sequence: Indeed, we obtain 
\begin{equation*}
L^{\solid}_{(\O_{\C_p}^\flat/p^\flat)/\F_p}\cong L^{\solid}_{(\O_{\C_p}^\flat/p^\flat)/\O_{\C_p}^\flat}\cong p^\flat\O_{\C_p}^\flat/(p^\flat)^2\O_{\C_p}^\flat[1]=\O_{\C_p}^\flat/p^\flat\{1\}[1]\;,
\end{equation*}
as desired.
\end{proof}

In principle, we would now like to cover $\C_p^{\dR, +}$ by the stack
\begin{equation*}
\GSpec(B_\dR^+\langle s\rangle_{\leq 1}\{r\}^\dagger/(rs-t))/\G_m(1)
\end{equation*}
over $\ol{\DD}_+/\G_m(1)$, where $s$ denotes the coordinate on $\ol{\DD}_+$ and $\G_m(1)$ acts on $r$ by division. However, as $B_\dR^+$ is not a Gelfand $\Q_p$-algebra -- note that the quotient $(B_\dR^+)^{\leq 1}/p$ is not discrete, but still has the $t$-adic topology! --, this does not quite make sense in the world of Gelfand stacks. We will thus work with the following replacement, which is the case $(R, R^+)=(\C_p, \O_{\C_p})$ of \cref{defi:htdr-bdr+dagger}:

\begin{defi}
The \emph{overconvergent (positive) de Rham period ring} is the ring
\begin{equation*}
B_\dR^{+, \dagger}\coloneqq A_\inf\{\xi\}^\dagger[\tfrac{1}{p}]=\colim_n A_\inf\langle p^{-n}\xi\rangle[\tfrac{1}{p}]\;,
\end{equation*}
where $\xi$ denotes a generator of the kernel of Fontaine's map $\theta: A_\inf\rightarrow\O_{\C_p}$.
\end{defi}

The key for our purposes is that $B_\dR^{+, \dagger}$ is now Gelfand and still has the unique lifting property from \cref{prop:clausen-solidbdr+}. Indeed, for any $R\rightarrow\C_p$ as in the proposition, the unique lift $B_\dR^+\rightarrow R$ yields a composite map $B_\dR^{+, \dagger}\rightarrow B_\dR^+\rightarrow R$. However, any such lift is also unique: As $R\rightarrow\C_p$ is a nilpotent thickening, any lift $B_\dR^{+, \dagger}\rightarrow R$ necessarily factors through $B_\dR^{+, \dagger}/t^n$ for some $n$ and the same is true for $B_\dR^+$ in place of $B_\dR^{+, \dagger}$. As $B_\dR^{+, \dagger}/t^n\cong B_\dR^+/t^n$, the uniqueness of the lift for $B_\dR^{+, \dagger}$ now follows from the uniqueness of the lift for $B_\dR^+$.

Thus, the cover of $\C_p^{\dR, +}$ that we are actually going to use is the following: consider the quotient stack
\begin{equation*}
\GSpec(B_\dR^{+, \dagger}\langle s\rangle_{\leq 1}\{r\}^\dagger/(rs-t))/\G_m(1)\;,
\end{equation*}
which evidently maps to $\ol{\DD}_+/\G_m(1)$; here, $s$ denotes the coordinate on $\ol{\DD}_+$ and $\G_m(1)$ acts on $r$ by division. Given an $A$-point of this quotient, we obtain a map $B_\dR^{+, \dagger}\rightarrow A$ as well as a generalised Cartier divisor $s: L\rightarrow A$ together with a factorisation
\begin{equation*}
A\xrightarrow{r} L\tensor_A \Nil^\dagger(A)\xrightarrow{s} A
\end{equation*}
of the multiplication-by-$t$-map. In particular, we have a commutative diagram
\begin{equation*}
\begin{tikzcd}
B_\dR^{+, \dagger}\ar[r]\ar[d, "t", swap] & A\ar[r, "r"] & L\tensor_A \Nil^\dagger(A)\ar[d, "s"] \\
B_\dR^{+, \dagger}\ar[rr] && A
\end{tikzcd}
\end{equation*}
and passing to cofibres under the vertical maps yields a map 
\begin{equation*}
\C_p\cong B_\dR^{+, \dagger}/t\rightarrow \cofib(\Nil^\dagger(A)\tensor_A L\rightarrow A)\;,
\end{equation*}
i.e.\ an $A$-point of $\C_p^{\dR, +}$.

Before we show that the map we have just constructed indeed provides a cover of $\C_p^{\dR, +}$, let us introduce the following notation:

\begin{defi}
Let $A$ be a Gelfand ring and $M\rightarrow A$ a map of anima. Then we define
\begin{equation*}
A\{M\}^\dagger\coloneqq \colim_{x_1, \dots, x_n\in M} A\tensor_{A[x_1, \dots, x_n]} A\{x_1, \dots, x_n\}^\dagger\;.
\end{equation*}
Similarly, for any $s\in A$, we write
\begin{equation*}
A\{Ms^{-1}\}^\dagger\coloneqq \colim_{x_1, \dots, x_n\in M} A\tensor_{A[x_1, \dots, x_n]} A\{\tfrac{x_1}{s}, \dots, \tfrac{x_n}{s}\}^\dagger\;.
\end{equation*}
\end{defi}

\begin{rem}
We immediately want to point out the following cautionary example: For $A=\Q_p[X_1, X_2, \dots]$ and $M=\{X_1, X_2, \dots\}$, we have
\begin{equation*}
A\{M\}^\dagger\cong \colim_n \Q_p\{X_1, \dots, X_n\}^\dagger\neq \Q_p\{X_1, X_2, \dots\}^\dagger\;.
\end{equation*}
We hope that this will not cause too much confusion.
\end{rem}

\begin{rem}
\label{rem:clausen-ms1daggercolims}
Note that $A\{M\}^\dagger$ is the universal Gelfand $A$-algebra $B$ such that the map $M\rightarrow A\rightarrow B$ factors through $\Nil^\dagger(B)$. Similarly, the ring $A\{Ms^{-1}\}^\dagger$ is the universal Gelfand $A$-algebra $B$ equipped with a factorisation of the map $M\rightarrow A\rightarrow B$ through the multiplication-by-$s$-map $s: \Nil^\dagger(B)\rightarrow B$. Using this universal property, it is clear that formation of $A\{M\}^\dagger$ and $A\{Ms^{-1}\}^\dagger$ commutes with colimits in $M$.
\end{rem}

Before we begin the proof that the map constructed above indeed covers $\C_p^{\dR, +}$, let us note the following lemmas:

\begin{lem}
\label{lem:clausen-nildaggercp}
The kernel of the multiplication map $\C_p\tensor_{\ol{\Q}_p}\C_p\rightarrow\C_p$ identifies with the ideal $\Nil^\dagger(\C_p\tensor_{\ol{\Q}_p}\C_p)$. In particular, $\Nil^\dagger(\C_p\tensor_{\ol{\Q}_p} \C_p)$ is generated by elements of the form $x\tensor 1-1\tensor x$.
\end{lem}
\begin{proof}
Since the elements $x\tensor 1-1\tensor x$ generate the kernel of $\C_p\tensor_{\ol{\Q}_p}\C_p\rightarrow\C_p$ and $\C_p$ is $\dagger$-reduced, we are done once we show that $x\tensor 1-1\tensor x$ is $\dagger$-nilpotent. However, this is due to the fact that $\ol{\Q}_p$ is dense in $\C_p$: for any $\epsilon>0$, we can find $y\in\ol{\Q}_p$ such that $|x-y|<\epsilon$ and then 
\begin{equation*}
x\tensor 1-1\tensor x=(x-y)\tensor 1-1\tensor (x-y)
\end{equation*}
has norm at most $\epsilon$ as well.
\end{proof}

\begin{lem}
\label{lem:clausen-makenullseq}
Let $A$ be a Gelfand ring and $s\in A$. Assume we have $a_1, a_2, \dots\in A$ whose images $\ol{a}_1, \ol{a}_2, \dots$ in $\Cone(\Nil^\dagger(A)\xrightarrow{\cdot s} A)$ are part of a convergent sequence with limit point $\ol{a}_\infty$. Then we can always lift the convergent sequence $\ol{a}_1, \ol{a}_2, \dots\rightarrow \ol{a}_\infty$ to a convergent sequence $a_1, a_2, \dots\rightarrow a_\infty$ in $A$ locally in the descendable topology. In fact, given any collection of such lifting problems, we can actually solve all of them at once up to a single descendable cover of $A$.
\end{lem}
\begin{proof}
As convergent sequences lift against surjections of solid abelian groups, we can find a convergent sequence $a'_1, a'_2, \dots\rightarrow a_\infty$ lifting the sequence $\ol{a}_1, \ol{a}_2, \dots\rightarrow \ol{a}_\infty$. If we can show that we can locally make $t_1=a_1-a'_1, t_2=a_2-a'_2, \dots$ a nullsequence whose image in $\Cone(\Nil^\dagger(A)\xrightarrow{\cdot s} A)$ vanishes, we are done.

To do this, it suffices to see that we can locally obtain a map
\begin{equation*}
A\{p^{-n}X_n s^{-1}: n\geq 1\}^\dagger\rightarrow A\;
\end{equation*}
sending $X_n$ to $t_n$; indeed, the sequence $(X_n)_n$ is a nullsequence in the source which is divisible by $s$ in $\Nil^\dagger(A)$. Since all the $t_n$ are divisible by $s$ in $\Nil^\dagger(A)$ by construction, however, we already know that $A$ admits a map from $\colim_N A\{p^{-n}X_n s^{-1}: n\leq N\}^\dagger$ which sends $X_n$ to $t_n$ and so we are done by the fact that the map 
\begin{equation*}
\colim_N A\{p^{-n}X_n s^{-1}: n\leq N\}^\dagger\rightarrow A\{p^{-n}X_n s^{-1}: n\geq 1\}^\dagger
\end{equation*}
is descendable. Indeed, by \cite[Prop.\ 2.7.2]{MannThesis}, it suffices to see that the maps $A\{p^{-n} X_n s^{-1}: n\leq N\}^\dagger\rightarrow A\{p^{-n} X_n s^{-1}: n\geq 1\}^\dagger$ are all descendable of bounded index, but they all admit sections and are hence descendable of index $1$.

For the last assertion, assume we are given any collection $\cal{S}$ of such lifting problems. As the above argument even works for finitely many lifting problems at once, we obtain descendable $A$-algebras $A_S$ whose index of descendability is at most $2$ for any finite subset $S\subseteq\cal{S}$ such that the lifting problems in $S$ admit solutions in $A_S$. Then $\colim_{S\subseteq\cal{S}\text{ finite}} A_S$ is descendable of index at most $4$ over $A$ and admits solutions to all lifting problems in $\cal{S}$, as desired.
\end{proof}

\begin{prop}
\label{prop:clausen-covercpdr+}
The map
\begin{equation*}
\GSpec(B_\dR^{+, \dagger}\langle s\rangle_{\leq 1}\{r\}^\dagger/(rs-t))/\G_m(1)\rightarrow \C_p^{\dR, +}
\end{equation*}
constructed above is surjective.
\end{prop}
\begin{proof}
We have to show that any $A$-point of $\C_p^{\dR, +}$ arises via the construction above up to a descendable cover of $A$. Thus, take any generalised Cartier divisor $s: L\rightarrow A$ on $A$ equipped with a map $\C_p\rightarrow \cofib(\Nil^\dagger(A)\tensor_A L\rightarrow A)$. As $\cofib(\Nil^\dagger(A)\tensor_A L\rightarrow A)$ is a $\dagger$-nilpotent thickening of $A$, the lifting problem
\begin{equation*}
\begin{tikzcd}
\Q_p\ar[rr]\ar[d] && A\ar[d] \\
E\ar[r]\ar[rru, dotted] & \C_p\ar[r] & \cofib(\Nil^\dagger(A)\tensor_A L\rightarrow A)
\end{tikzcd}
\end{equation*}
admits a unique solution for any finite extension $E/\Q_p$ since solid étale maps are $\dagger$-formally étale by \cite[Prop.\ 3.7.5]{dRStack}. We conclude that $A$ uniquely admits the structure of a $\ol{\Q}_p$-algebra. Moreover, after passing to a descendable cover of $A$, we may assume that $L$ is principal and that $s$ corresponds to an element of $A$. Finally, note that we are given a map $\C_p\rightarrow\Cone(\Nil^\dagger(A)\xrightarrow{\cdot s} A)$ and hence every Cauchy sequence in $\ol{\Q}_p$ has a limit in $\Cone(\Nil^\dagger(A)\xrightarrow{\cdot s} A)$. By \cref{lem:clausen-makenullseq}, up to passing to a descendable cover, we can then arrange that any such sequence also has a limit in $A$, compatibly with the given limit in the quotient by $s\Nil^\dagger(A)$.

Now let $B\coloneqq (A\tensor_{\ol{\Q}_p} B_\dR^{+, \dagger})\{r\}^\dagger/(rs-t)$. Since both $\Nil^\dagger(A)\tensor_A B$ and $r$ are contained in the $\dagger$-nilradical of $B$, the map $B\rightarrow \Cone(\Nil^\dagger(B)\xrightarrow{\cdot s} B)$ factors as
\begin{equation*}
B\rightarrow \Cone(\Nil^\dagger(A)\xrightarrow{\cdot s} A)\tensor_{\ol{\Q}_p} \C_p\rightarrow \Cone(\Nil^\dagger(B)\xrightarrow{\cdot s} B)\;.
\end{equation*}
Moreover, the given map $\C_p\rightarrow\Cone(\Nil^\dagger(A)\xrightarrow{\cdot s} A)$ induces
\begin{equation*}
\C_p\tensor_{\ol{\Q}_p} \C_p\rightarrow \Cone(\Nil^\dagger(A)\xrightarrow{\cdot s} A)\tensor_{\ol{\Q}_p} \C_p\;.
\end{equation*}
Now let $M$ be the $\Q_p$-module defined by the pullback diagram
\begin{equation*}
\begin{tikzcd}
M\ar[r]\ar[d] & B\ar[d] \\
\Nil^\dagger(\C_p\tensor_{\ol{\Q}_p} \C_p)\ar[r] & \Cone(\Nil^\dagger(A)\xrightarrow{\cdot s} A)\tensor_{\ol{\Q}_p} \C_p
\end{tikzcd}
\end{equation*}
and denote by $F$ the fibre of the map $M\rightarrow \Nil^\dagger(\C_p\tensor_{\ol{\Q}_p} \C_p)$. Noting that the map $F\rightarrow M\rightarrow B$ is naturally equipped with a factorisation through $\Nil^\dagger(B)\xrightarrow{\cdot s} B$, we obtain a map 
\begin{equation*}
B\{Fs^{-1}\}^\dagger\rightarrow B
\end{equation*}
by the universal property of the $B$-algebra $B\{Fs^{-1}\}^\dagger$. Then we define
\begin{equation*}
C\coloneqq B\{Ms^{-1}\}^\dagger\tensor_{B\{Fs^{-1}\}^\dagger} B\;.
\end{equation*}

If we can prove that the map $A\rightarrow C$ is descendable, we are done: Indeed, note that $C$ clearly maps to $\GSpec(B_\dR^{+, \dagger}\langle s\rangle_{\leq 1}\{r\}^\dagger/(rs-t))/\G_m(1)$ and, by construction of $C$, the map
\begin{equation*}
\C_p\tensor_{\ol{\Q}_p} \C_p\rightarrow \Cone(\Nil^\dagger(A)\xrightarrow{\cdot s} A)\tensor_{\ol{\Q}_p} \C_p\rightarrow \Cone(\Nil^\dagger(B)\xrightarrow{\cdot s} B)\rightarrow\Cone(\Nil^\dagger(C)\xrightarrow{\cdot s} C)
\end{equation*}
factors through the quotient of the source by $\Nil^\dagger(\C_p\tensor_{\ol{\Q}_p} \C_p)$. By \cref{lem:clausen-nildaggercp}, this precisely means that the $C$-point of $\C_p^{\dR, +}$ induced by the given $A$-point can be identified with the $C$-point of $\C_p^{\dR, +}$ coming from the map
\begin{equation*}
\GSpec C\rightarrow \GSpec(B_\dR^{+, \dagger}\langle s\rangle_{\leq 1}\{r\}^\dagger/(rs-t))/\G_m(1)
\end{equation*}
and the desired surjectivity follows.

To prove descendability, first note that we have
\begin{equation*}
B\{Ms^{-1}\}^\dagger\cong \colim_{x_1, \dots, x_n\in \Nil^\dagger(\C_p\tensor_{\ol{\Q}_p}\C_p)} B\{(F+(x_1+F)+\dots+(x_n+F))s^{-1}\}^\dagger
\end{equation*}
by compatibility of $B\{Ms^{-1}\}^\dagger$ with colimits, see \cref{rem:clausen-ms1daggercolims}. Then further observe that we may just as well take the colimit only over $x_i$ of the form $a_i\tensor 1-1\tensor a_i$ for $a_i\in\C_p$. Indeed, such elements generate $\Nil^\dagger(\C_p\tensor_{\ol{\Q}_p}\C_p)$ by \cref{lem:clausen-nildaggercp} and the pullback 
\begin{equation*}
\begin{tikzcd}
R\ar[r]\ar[d] & B\ar[d] \\
\C_p\tensor_{\ol{\Q}_p}\C_p\ar[r] & \Cone(\Nil^\dagger(A)\xrightarrow{\cdot s} A)\tensor_{\ol{\Q}_p} \C_p
\end{tikzcd}
\end{equation*}
has a natural ring structure making $R\rightarrow B$ a ring map and inducing an action of $R$ on $M$ by multiplication. Thus, if $x_1, \dots, x_n\in \Nil^\dagger(\C_p\tensor_{\ol{\Q}_p}\C_p)$ are in the $\C_p\tensor_{\ol{\Q}_p} \C_p$-linear span of $y_1, \dots, y_m$ with $y_j=a_j\tensor 1-1\tensor a_j$, then
\begin{equation*}
\begin{split}
B\{(F+(x_1+F)+\dots+(x_n+F)&+(y_1+F)+\dots+(y_m+F))s^{-1}\}^\dagger \\
&\cong B\{(F+(y_1+F)+\dots+(y_m+F))s^{-1}\}^\dagger
\end{split}
\end{equation*}
via the natural map and this proves our claim. We conclude that
\begin{equation*}
C\cong \colim_{a_1, \dots, a_m\in \C_p} B\{(F+(y_1+F)+\dots+(y_m+F))s^{-1}\}^\dagger\tensor_{B\{Fs^{-1}\}^\dagger} B\;,
\end{equation*}
where we have written $y_j=a_j\tensor 1-1\tensor a_j$, as above, and, by \cite[Prop.\ 2.7.2]{MannThesis}, it thus suffices to show that the terms in the colimit is descendable of bounded index over $A$.

Now take finitely many $a_1, \dots, a_m\in\C_p$. By Noether normalisation, the $\ol{\Q}_p$-algebra generated by the $a_i$ is finite over a polynomial algebra $\ol{\Q}_p[t_1, \dots, t_n]$ and since descendability satisfies cancellation, it suffices to show that
\begin{equation*}
B\{(F+(y_1+F)+\dots+(y_{m+n}+F))s^{-1}\}^\dagger\tensor_{B\{Fs^{-1}\}^\dagger} B
\end{equation*}
is descendable over $A$ of index bounded independently of the $a_j$; here, we have set $y_{m+i}=t_i\tensor 1-1\tensor t_i$. For this, consider the composition
\begin{equation*}
\ol{\Q}_p[t_1,\dots, t_n]\rightarrow\C_p\rightarrow\Cone(\Nil^\dagger(A)\xrightarrow{\cdot s} A)\;.
\end{equation*}
Since we have arranged that Cauchy sequences coming from $\ol{\Q}_p$ lift from $\Cone(\Nil^\dagger(A)\xrightarrow{\cdot s} A)$ to $A$, we may find a lift of the above map to $\ol{\Q}_p[t_1, \dots, t_n]\rightarrow A$ by picking any Cauchy sequence in $\ol{\Q}_p$ converging to $t_i\in\C_p$ and sending $t_i$ to the corresponding limit in $A$. (Note that $\ol{\Q}_p[t_1, \dots, t_n]$ has the subspace topology from $\C_p$ here, i.e., as a condensed ring, it is not the same as a polynomial ring.) As being invertible may be checked after $\dagger$-reduction, we in fact obtain a map
\begin{equation*}
\ol{\Q}_p(t_1, \dots, t_n)\rightarrow A\;.
\end{equation*}
As $\ol{\Q}_p[a_1, \dots, a_m]$ is finite over $\ol{\Q}_p[t_1, \dots, t_n]$, the field extension
\begin{equation*}
\ol{\Q}_p(t_1, \dots, t_n)\rightarrow \ol{\Q}_p(t_1, \dots, t_n)(a_1, \dots, a_m)
\end{equation*}
is finite. Moreover, via $\C_p$, the target maps to $\Cone(\Nil^\dagger(A)\xrightarrow{\cdot s} A)$ and since finite étale maps are $\dagger$-formally étale by \cite[Prop.\ 3.7.5]{dRStack}, we obtain a unique lift of this map to
\begin{equation*}
\ol{\Q}_p(t_1, \dots, t_n)(a_1, \dots, a_m)\rightarrow A\;.
\end{equation*}

In particular, we can now make sense of $t_i\tensor 1-1\tensor t_i$ and $a_i\tensor 1-1\tensor a_i$ in $A\tensor_{\ol{\Q}_p} \C_p$ and let $\tilde{y}_j$ denote arbitrary lifts of these elements to $B$. Then observe that
\begin{equation*}
B\{(F+(y_1+F)+\dots+(y_{m+n}+F))s^{-1}\}^\dagger\tensor_{B\{Fs^{-1}\}^\dagger} B\cong B\{\tilde{y}_1s^{-1}, \dots, \tilde{y}_{m+n}s^{-1}\}^\dagger
\end{equation*}
and note that the right-hand side maps to
\begin{equation*}
\begin{split}
B\{\tilde{y}_1s^{-1}, &\dots, \tilde{y}_{m+n}s^{-1}\}^\dagger/(r, y_1s^{-1}, \dots, y_{m+n}s^{-1}) \\
&\cong (A\tensor_{\ol{\Q}_p} \C_p)\{\tilde{y}_1s^{-1}, \dots, \tilde{y}_{m+n}s^{-1}\}^\dagger/(y_1s^{-1}, \dots, y_{m+n}s^{-1}) \\
&\cong A\tensor_{\ol{\Q}_p(t_1, \dots, t_n)(a_1, \dots, a_m)} \C_p\;,
\end{split}
\end{equation*}
where the last isomorphism uses that the images of the $\tilde{y}_j$ in $A\tensor_{\ol{\Q}_p}\C_p$ are given by $t_i\tensor 1-1\tensor t_i$ or $a_i\tensor 1-1\tensor a_i$, respectively. Since descendability satisfies cancellation, this reduces us to showing that the map
\begin{equation*}
A\rightarrow A\tensor_{\ol{\Q}_p(t_1, \dots, t_n)(a_1, \dots, a_m)} \C_p
\end{equation*}
is descendable of bounded index -- but this is clear as $\ol{\Q}_p(t_1, \dots, t_n)(a_1, \dots, a_m)\rightarrow \C_p$ is descendable of index bounded independently of the $t_i$ and the $a_i$.

To see this last claim and thereby finish the proof, write $K=\ol{\Q}_p(t_1, \dots, t_n)(a_1, \dots, a_m)$ and choose a transcendence basis $\{b_i\}_{i\in I}$ of $\C_p$ over $K$. As $K(b_i: i\in I)\rightarrow\C_p$ is a filtered colimit of finite étale maps and these are always descendable and since these are always descendable of bounded index, we may reduce to showing that $K\rightarrow K(b_i: i\in I)$ is descendable by \cite[Prop.\ 2.7.2]{MannThesis}. However, by \cite[Prop.\ 4.6.4.(1)]{dRFF}, we can find a descendable $K$-algebra $L$ such that $L$ surjects onto its uniform completion. If we can produce a factorisation
\begin{equation*}
K\rightarrow K(b_i: i\in I)\rightarrow L
\end{equation*}
of the map $K\rightarrow L$ we will be done since descendability satisfies cancellation. For this, we pick a Cauchy sequence in $\ol{\Q}_p$ converging to $b_i\in\C_p$ for each $i\in I$ and note that this will map to a convergent sequence in $L^u=\ol{L}$ since the uniform completion of $K$ is $\C_p$. By \cref{lem:clausen-makenullseq}, we can then arrange that the corresponding sequences in $L$ still converge after possibly passing to a further descendable cover and then the map sending each $b_i$ to the corresponding limit will produce the desired factorisation.
\end{proof}

\begin{rem}
In fact, using \cref{prop:clausen-solidbdr+}, one can even show that the above map is an isomorphism after base changing along $\widehat{\G}_a/\G_m(1)\rightarrow\ol{\DD}/\G_m(1)$.
\end{rem}

\begin{lem}
\label{lem:clausen-cechnerve}
The $n$-th term of the \v{C}ech nerve of the surjection from \cref{prop:clausen-covercpdr+} is given by
\begin{equation*}
\GSpec((B_\dR^{+, \dagger}\tensor_{\ol{\Q}_p}\dots\tensor_{\ol{\Q}_p}B_\dR^{+, \dagger})\langle s\rangle_{\leq 1}\{\Nil^\dagger(B_\dR^{+, \dagger}\tensor_{\ol{\Q}_p}\dots\tensor_{\ol{\Q}_p}B_\dR^{+, \dagger})s^{-1}\}^\dagger)/\G_m(1)\;,
\end{equation*}
where the occurring tensor products are $n$-fold each. Moreover, $\Nil^\dagger(B_\dR^{+, \dagger}\tensor_{\ol{\Q}_p}\dots\tensor_{\ol{\Q}_p}B_\dR^{+, \dagger})$ agrees with the kernel of the multiplication map $B_\dR^{+, \dagger}\tensor_{\ol{\Q}_p}\dots\tensor_{\ol{\Q}_p}B_\dR^{+, \dagger}\rightarrow\C_p$.
\end{lem}
\begin{proof}
We first address the last sentence: By design, we have $t\in\Nil^\dagger(B_\dR^{+, \dagger})$ and hence the kernel of the map
\begin{equation*}
B_\dR^{+, \dagger}\tensor_{\ol{\Q}_p}\dots\tensor_{\ol{\Q}_p}B_\dR^{+, \dagger}\rightarrow \C_p\tensor_{\ol{\Q}_p}\dots\tensor_{\ol{\Q}_p}\C_p
\end{equation*}
lies in $\Nil^\dagger(B_\dR^{+, \dagger}\tensor_{\ol{\Q}_p}\dots\tensor_{\ol{\Q}_p}B_\dR^{+, \dagger})$. As the multiplication map $\C_p\tensor_{\ol{\Q}_p}\dots\tensor_{\ol{\Q}_p}\C_p\rightarrow\C_p$ coincides with the $\dagger$-reduction of the source by \cref{lem:clausen-nildaggercp}, this proves the claim.

For the main part of the proof, consider a map $\GSpec A\rightarrow\ol{\DD}/\G_m(1)$ with associated generalised Cartier divisor $s: L\rightarrow A$ and assume that $L$ is principal for simplicity (we can always ensure this after passing to a cover). An $A$-point of the $n$-th term of the \v{C}ech nerve consists of the following data:
\begin{enumerate}[label=(\arabic*)]
\item $n$ maps $B_\dR^{+, \dagger}\rightarrow A$,
\item for each of these maps a factorisation
\begin{equation*}
\begin{tikzcd}
tB_\dR^{+, \dagger}\ar[r, dotted]\ar[d] & \Nil^\dagger(A)\ar[d, "s"] \\
B_\dR^{+, \dagger}\ar[r] & A\nospacepunct{\;,}
\end{tikzcd}
\end{equation*}
\item an identification of all the $n$ induced maps $\C_p\cong B_\dR^{+, \dagger}/t\rightarrow A/^\mathbb{L}s\Nil^\dagger(A)$.
\end{enumerate}
Now observe that the common induced map $\C_p\rightarrow A/^\mathbb{L} s\Nil^\dagger(A)$ induces a unique $\ol{\Q}_p$-algebra structure on $A$ as in the proof of \cref{prop:clausen-covercpdr+} and then the piece of data (1) is equivalent to a map $B_\dR^{+, \dagger}\tensor_{\ol{\Q}_p}\dots\tensor_{\ol{\Q}_p}B_\dR^{+, \dagger}\rightarrow A$. In this language, (2) says that each element of $B_\dR^{+, \dagger}\tensor_{\ol{\Q}_p}\dots\tensor_{\ol{\Q}_p} tB_\dR^{+, \dagger}\tensor_{\ol{\Q}_p}\dots\tensor_{\ol{\Q}_p}B_\dR^{+, \dagger}$ should be mapped to an element of $A$ which is the product of $s$ with an element in $\Nil^\dagger(A)$. Once this is ensured, we obtain an induced map
\begin{equation*}
\C_p\tensor_{\ol{\Q}_p}\dots\tensor_{\ol{\Q}_p}\C_p\rightarrow A/^\mathbb{L} s\Nil^\dagger(A)
\end{equation*}
and then (3) requires a factorisation of this map through the multiplication map $\C_p\tensor_{\ol{\Q}_p}\dots\tensor_{\ol{\Q}_p}\C_p\rightarrow \C_p$. Summarising and using that $\Nil^\dagger(B_\dR^{+, \dagger}\tensor_{\ol{\Q}_p}\dots\tensor_{\ol{\Q}_p}B_\dR^{+, \dagger})$ is identified with the kernel of the multiplication map $B_\dR^{+, \dagger}\tensor_{\ol{\Q}_p}\dots\tensor_{\ol{\Q}_p}B_\dR^{+, \dagger}\rightarrow \C_p$, we see that we can reformulate (2) and (3) as requiring a factorisation
\begin{equation*}
\label{eq:clausen-factorisation}
\begin{tikzcd}
\Nil^\dagger(B_\dR^{+, \dagger}\tensor_{\ol{\Q}_p}\dots\tensor_{\ol{\Q}_p}B_\dR^{+, \dagger})\ar[r]\ar[d] & \Nil^\dagger(A)\ar[d, "s"] \\
B_\dR^{+, \dagger}\tensor_{\ol{\Q}_p}\dots\tensor_{\ol{\Q}_p}B_\dR^{+, \dagger}\ar[r] & A
\end{tikzcd}
\end{equation*}
and it is obvious that an $A$-valued point
\begin{equation*}
\GSpec A\rightarrow \GSpec((B_\dR^{+, \dagger}\tensor_{\ol{\Q}_p}\dots\tensor_{\ol{\Q}_p}B_\dR^{+, \dagger})\langle s\rangle_{\leq 1}\{\Nil^\dagger(B_\dR^{+, \dagger}\tensor_{\ol{\Q}_p}\dots\tensor_{\ol{\Q}_p}B_\dR^{+, \dagger})s^{-1}\}^\dagger)/\G_m(1)
\end{equation*}
lifting the given $A$-valued point of $\ol{\DD}/\G_m(1)$ exactly corresponds to the data of a map $B_\dR^{+, \dagger}\tensor_{\ol{\Q}_p}\dots\tensor_{\ol{\Q}_p}B_\dR^{+, \dagger}\rightarrow A$ together with a factorisation as in (\ref{eq:clausen-factorisation}).
\end{proof}

\begin{cor}
\label{cor:clausen-cpdr+complexes}
For $n\geq 1$, the cohomology of $\O(-n)$ on $\C_p^{\dR, +}$ is calculated by the following complex:
\begin{equation*}
t^nB_\dR^{+, \dagger}\rightarrow \Nil^\dagger(B_\dR^{+, \dagger}\tensor_{\ol{\Q}_p} B_\dR^{+, \dagger})^n\rightarrow \Nil^\dagger(B_\dR^{+, \dagger}\tensor_{\ol{\Q}_p} B_\dR^{+, \dagger}\tensor_{\ol{\Q}_p} B_\dR^{+, \dagger})^n\rightarrow\dots
\end{equation*}
Meanwhile, for $n\leq 0$, all $\O(-n)$ have the same cohomology on $\C_p^{\dR, +}$ and this cohomology is computed by the complex
\begin{equation*}
B_\dR^{+, \dagger}\rightarrow B_\dR^{+, \dagger}\tensor_{\ol{\Q}_p} B_\dR^{+, \dagger}\rightarrow B_\dR^{+, \dagger}\tensor_{\ol{\Q}_p} B_\dR^{+, \dagger}\tensor_{\ol{\Q}_p} B_\dR^{+, \dagger}\rightarrow\dots\nospacepunct{\;.}
\end{equation*}
\end{cor}
\begin{proof}
This follows from the previous lemma using the fact that the cohomology of $\O(-n)$ on $\GSpec R/\G_m(1)$ is given by the degree $-n$ piece of $R$, where the grading comes from the $\G_m(1)$-action. Indeed, the \v{C}ech nerve of the cover $\GSpec R\rightarrow\GSpec R/\G_m(1)$ is given by
\begin{equation*}
\begin{tikzcd}
\dots\ar[r, shift left=2]\ar[r]\ar[r, shift right=2] & \GSpec R\times \G_m(1)\ar[r, shift left=1, "\mathrm{act}"]\ar[r, shift right=1, "\mathrm{pr}_1", swap] & \GSpec R\ar[r] & \GSpec R/\G_m(1)
\end{tikzcd}
\end{equation*}
and hence we have
\begin{equation*}
\begin{tikzcd}
R\Gamma(\GSpec R/\G_m(1), \O(-n))\ar[r, "\cong"] &\operatorname{eq}(R\ar[r, shift left=1, "s^n\cdot \mathrm{act}"]\ar[r, shift right=1, "\mathrm{incl}", swap] & R\tensor_{\Q_p} \Q_p\langle s, s^{-1}\rangle_{\leq 1})\;,
\end{tikzcd}
\end{equation*}
which is the degree $-n$ piece of $R$ by definition. (Recall that $s$ has degree $1$.)
\end{proof}

Note that $\ol{\Q}_p\rightarrow B_\dR^{+, \dagger}$ is descendable by cancellation since the composition $\ol{\Q}_p\rightarrow \C_p$ with the projection to $\C_p$ is. Thus, by descendable descent for $\G_a$, the previous corollary already proves \cref{thm:clausen-cpdr+} for all $n\leq 0$. For $n\geq 1$, the key to computing the cohomology of the complex above lies in the following statement:

\begin{lem}
\label{lem:clausen-nilpowerker}
For each $n\geq 1$, the ideal $\Nil^\dagger(B_\dR^{+, \dagger}\tensor_{\ol{\Q}_p}\dots\tensor_{\ol{\Q}_p}B_\dR^{+, \dagger})^n$ identifies with the kernel of the multiplication map
\begin{equation*}
B_\dR^{+, \dagger}\tensor_{\ol{\Q}_p}\dots\tensor_{\ol{\Q}_p}B_\dR^{+, \dagger}\rightarrow B_\dR^{+, \dagger}/t^n\;.
\end{equation*}
\end{lem}
\begin{proof}
The heart of the argument lies in the following

\bigskip

\textbf{Claim.} The map $B_\dR^{+, \dagger}\rightarrow B_\dR^{+, \dagger}\tensor_{\ol{\Q}_p} B_\dR^{+, \dagger}$ given by $x\mapsto x\tensor 1-1\tensor x$ factors through $\Nil^\dagger(B_\dR^{+, \dagger}\tensor_{\ol{\Q}_p} B_\dR^{+, \dagger})^n$ for all $n\geq 1$.

\bigskip

\textit{Proof of the claim.} Observe that the map $B_\dR^{+, \dagger}\tensor_{\ol{\Q}_p} B_\dR^{+, \dagger}/\Nil^\dagger(B_\dR^{+, \dagger}\tensor_{\ol{\Q}_p} B_\dR^{+, \dagger})^n\rightarrow \C_p$ induced by the multiplication map is an infinitesimal thickening and hence, by \cref{prop:clausen-solidbdr+}, there is a unique map filling the diagram
\begin{equation*}
\begin{tikzcd}
B_\dR^{+, \dagger}\tensor_{\ol{\Q}_p} B_\dR^{+, \dagger}/\Nil^\dagger(B_\dR^{+, \dagger}\tensor_{\ol{\Q}_p} B_\dR^{+, \dagger})^n\ar[r] & \C_p \\
& B_\dR^{+, \dagger}\ar[u]\ar[lu, dotted, "\exists !"]\nospacepunct{\;.}
\end{tikzcd}
\end{equation*}
However, note that both the map $x\mapsto 1\tensor x$ and the map $x\mapsto x\tensor 1$ can be put in place of the dotted arrow to make the diagram commute, hence we must have $x\tensor 1-1\tensor x\in\Nil^\dagger(B_\dR^{+, \dagger}\tensor_{\ol{\Q}_p} B_\dR^{+, \dagger})^n$, as desired. \hfill\qed

\bigskip

We now move on to the proof of the lemma and first argue that the multiplication map
\begin{equation*}
B_\dR^{+, \dagger}\tensor_{\ol{\Q}_p}\dots\tensor_{\ol{\Q}_p}B_\dR^{+, \dagger}\rightarrow B_\dR^{+, \dagger}/t^n
\end{equation*}
factors through the quotient of the source by $\Nil^\dagger(B_\dR^{+, \dagger}\tensor_{\ol{\Q}_p}\dots\tensor_{\ol{\Q}_p}B_\dR^{+, \dagger})^n$. For this, it suffices to show that the image of any element in $\Nil^\dagger(B_\dR^{+, \dagger}\tensor_{\ol{\Q}_p}\dots\tensor_{\ol{\Q}_p}B_\dR^{+, \dagger})$ under the multiplication map 
\begin{equation*}
B_\dR^{+, \dagger}\tensor_{\ol{\Q}_p}\dots\tensor_{\ol{\Q}_p}B_\dR^{+, \dagger}\rightarrow B_\dR^{+, \dagger}
\end{equation*}
is divisible by $t$, but this is clear: any such element will be sent to zero under the composition of this map with the projection $B_\dR^{+, \dagger}\rightarrow B_\dR^{+, \dagger}/t\cong \C_p$ by the second part of \cref{prop:clausen-covercpdr+}.

Thus, we have a well-defined map
\begin{equation*}
B_\dR^{+, \dagger}\tensor_{\ol{\Q}_p}\dots\tensor_{\ol{\Q}_p}B_\dR^{+, \dagger}/\Nil^\dagger(B_\dR^{+, \dagger}\tensor_{\ol{\Q}_p}\dots\tensor_{\ol{\Q}_p}B_\dR^{+, \dagger})^n\rightarrow B_\dR^{+, \dagger}/t^n
\end{equation*}
given by $x_1\tensor\dots\tensor x_k\mapsto x_1\cdots x_k$ and we will actually produce an inverse of this map. For this, take the map
\begin{equation*}
B_\dR^{+, \dagger}\rightarrow B_\dR^{+, \dagger}\tensor_{\ol{\Q}_p}\dots\tensor_{\ol{\Q}_p}B_\dR^{+, \dagger}/\Nil^\dagger(B_\dR^{+, \dagger}\tensor_{\ol{\Q}_p}\dots\tensor_{\ol{\Q}_p}B_\dR^{+, \dagger})^n
\end{equation*}
given by $x\mapsto x\tensor 1\tensor\dots\tensor 1$ and note that this factors through $B_\dR^{+, \dagger}/t^n$ since $t\tensor 1\tensor\dots\tensor 1\in \Nil^\dagger(B_\dR^{+, \dagger}\tensor_{\ol{\Q}_p}\dots\tensor_{\ol{\Q}_p}B_\dR^{+, \dagger})$. As it is clear that the composition
\begin{equation*}
B_\dR^{+, \dagger}/t^n\rightarrow B_\dR^{+, \dagger}\tensor_{\ol{\Q}_p}\dots\tensor_{\ol{\Q}_p}B_\dR^{+, \dagger}/\Nil^\dagger(B_\dR^{+, \dagger}\tensor_{\ol{\Q}_p}\dots\tensor_{\ol{\Q}_p}B_\dR^{+, \dagger})^n\rightarrow B_\dR^{+, \dagger}/t^n
\end{equation*}
is the identity, it remains to check the other composition and this amounts to checking that
\begin{equation*}
x_1\tensor\dots\tensor x_k=(x_1\cdots x_k)\tensor 1\dots\tensor 1\in B_\dR^{+, \dagger}\tensor_{\ol{\Q}_p}\dots\tensor_{\ol{\Q}_p}B_\dR^{+, \dagger}/\Nil^\dagger(B_\dR^{+, \dagger}\tensor_{\ol{\Q}_p}\dots\tensor_{\ol{\Q}_p}B_\dR^{+, \dagger})^n\;.
\end{equation*}
However, the claim implies that 
\begin{equation*}
x_i\tensor\dots \tensor 1\tensor\dots=1\tensor\dots \tensor x_i\tensor\dots\in B_\dR^{+, \dagger}\tensor_{\ol{\Q}_p}\dots\tensor_{\ol{\Q}_p}B_\dR^{+, \dagger}/\Nil^\dagger(B_\dR^{+, \dagger}\tensor_{\ol{\Q}_p}\dots\tensor_{\ol{\Q}_p}B_\dR^{+, \dagger})^n
\end{equation*}
for any $i=1, \dots, k$ and applying this repeatedly yields
\begin{equation*}
\begin{split}
x_1\tensor\dots\tensor x_k&=(x_1x_2)\tensor 1\tensor x_3\tensor \dots\tensor x_k \\
&=(x_1x_2x_3)\tensor 1\tensor 1\tensor x_4\tensor \dots\tensor x_k \\
&=\dots=(x_1\cdots x_k)\tensor 1\tensor\dots\tensor 1\;,
\end{split}
\end{equation*}
as desired.
\end{proof}

Now we can finish as follows:

\begin{proof}[Proof of \cref{thm:clausen-cpdr+}]
Note that, for any $n\geq 1$, there is a commutative diagram of complexes
\begin{equation*}
\begin{tikzcd}
B_\dR^{+, \dagger}\ar[r]\ar[d] & B_\dR^{+, \dagger}\tensor_{\ol{\Q}_p} B_\dR^{+, \dagger}\ar[r]\ar[d] & B_\dR^{+, \dagger}\tensor_{\ol{\Q}_p} B_\dR^{+, \dagger}\tensor_{\ol{\Q}_p} B_\dR^{+, \dagger}\ar[r]\ar[d] & \dots \\
B_\dR^{+, \dagger}/t^n\ar[r, "0"] & B_\dR^{+, \dagger}/t^n\ar[r, equals] & B_\dR^{+, \dagger}/t^n\ar[r, "0"] & \dots\nospacepunct{\;,}
\end{tikzcd}
\end{equation*}
where the vertical maps are the multiplication maps. On the one hand, \cref{lem:clausen-nilpowerker} shows taking fibres of the vertical maps yields the complex computing $R\Gamma(\C_p^{\dR, +}, \O(-n))$ by \cref{cor:clausen-cpdr+complexes}. On the other hand, the top complex is quasi-isomorphic to $\ol{\Q}_p$ by descendability of $\ol{\Q}_p\rightarrow B_\dR^{+, \dagger}$ and the fact that $\G_a$ satisfies descendable descent (see above) while the bottom complex is quasi-isomorphic to $B_\dR^{+, \dagger}/t^n$. Thus, the fibre of the map of complexes above is also given by $(B_\dR^{+, \dagger}/t^n)/\ol{\Q}_p[-1]$ and this proves the theorem as $B_\dR^{+, \dagger}/t^n\cong B_\dR^+/t^n$.
\end{proof}
}

\subsection{Nice hypercovers for some Nygaardifications}

To actually apply the above results, we have to check that some Gelfand stacks of interest are indeed nicely coverable. For this, note that $X^\N$ inherits a radius map from $X^\prism$ via
\begin{equation*}
X^\N\xrightarrow{\pi} X^\prism\rightarrow (0, \infty)
\end{equation*}
for any Gelfand stack $X$.

\begin{prop}
\label{prop:perf-zpncover}
For any $1<r<s<\infty$, the Gelfand stack $(\Q_p^\N)_{[r, s]}$ is nicely coverable.
\end{prop}
\begin{proof}
First note that $(\Q_p^\prism)_{[r, s]}=(\Spd\Q_p\times\Spd\Q_p)^\dR_{[r, s]}$ by \cref{prop:defis-prismffdr} and recall that
\begin{equation*}
(\Spd\Q_p\times\Spd\Q_p)^\dR\cong \lim_{x\mapsto x^p} (1+\overcirc{\DD})^\dR\setminus\{1\}\,\big/\,\Z_p^{\times, \sm}\;,
\end{equation*}
where $\overcirc{\DD}$ denotes the open unit disk. The quotient stack
\begin{equation*}
\lim_{x\mapsto x^p} (1+\overcirc{\DD})^\dR\,\big/\,\Z_p^{\times, \sm}
\end{equation*}
clearly has a descendable cover by $\lim_{x\mapsto x^p} (1+\overcirc{\DD})^\dR$ whose \v{C}ech nerve is given by
\begin{equation*}
\left(\lim_{x\mapsto x^p} (1+\overcirc{\DD})^\dR\right)\times (\Z_p^{\times, \sm})^\bullet\;,
\end{equation*}
where we note that $\O(\Z_p^{\times, \sm})=C^\sm(\Z_p^\times, \Q_p)$ is a flat $\Q_p$-algebra by \cref{lem:perf-vspflat}. Moreover, note that the endomorphism of $1+\overcirc{\DD}$ given by $x\mapsto x^p$ is standard finite étale, hence $\dagger$-formally étale by \cite[Prop.\ 3.7.5]{dRStack}, and therefore
\begin{equation*}
\lim_{x\mapsto x^p} (1+\overcirc{\DD})^\dR\cong \lim_{x\mapsto x^p} (1+\overcirc{\DD})\,\big/\,\G_m^\dagger\;.
\end{equation*}
As also $\O(\G_m^\dagger)=\colim_n \Q_p\langle p^{-n}T\rangle$ is a flat $\Q_p$-algebra, the \v{C}ech cover
\begin{equation*}
\lim_{x\mapsto x^p} (1+\overcirc{\DD})\,\times\, (\G_m^\dagger)^\bullet\rightarrow \lim_{x\mapsto x^p} (1+\overcirc{\DD})^\dR
\end{equation*}
provides a nice hypercover of any base change of $\lim_{x\mapsto x^p} (1+\overcirc{\DD})^\dR$ to the de Rham stack of a closed annulus centered around $1$. In particular, putting everything together and applying \cref{lem:perf-finflatdim} and \cref{lem:perf-covers}, we see that $(\Q_p^\prism)_{[r, s]}$ is nicely coverable and in fact we have even found a nice hypercover where all the algebras occurring are static and flat over $\Q_p$.

Now note that $(\ol{\DD}/\G_m(1))^\dR$ is nicely coverable: Indeed, the \v{C}ech nerve of the map $\ol{\DD}^\dR\rightarrow (\ol{\DD}/\G_m(1))^\dR$ is given by $\ol{\DD}^\dR\times (\G_m(1)^\dR)^\bullet$; moreover, $\ol{\DD}^\dR$ and $\G_m(1)^\dR$ are themselves covered by $\ol{\DD}$ and $\G_m(1)$, respectively, with the respective \v{C}ech nerves being given by
\begin{equation*}
\ol{\DD}\times (\G_a^\dagger)^\bullet\;,\hspace{0.3cm}\text{resp.}\hspace{0.3cm} \G_m(1)\times (\G_m^\dagger)^\bullet\;.
\end{equation*}
As $\ol{\DD}, \G_m(1), \G_a^\dagger$ and $\G_m^\dagger$ are all affine static Gelfand stacks which are flat over $\Q_p$ by \cref{lem:perf-vspflat}, this yields that $(\ol{\DD}/\G_m(1))^\dR$ is nicely coverable by \cref{lem:perf-covers} and \cref{lem:perf-bc}. Clearly, also $\ol{\DD}/\G_m(1)$ is nicely coverable by the \v{C}ech cover $\ol{\DD}\times \G_m(1)^\bullet$. Again, note that we have actually produced nice hypercovers which only involve static flat $\Q_p$-algebras.

Denoting by $X_\bullet\rightarrow (\Q_p^\prism)_{[r, s]}$ and $Y_{\bullet'}\rightarrow \ol{\DD}/\G_m(1)\times (\ol{\DD}/\G_m(1))^\dR$ the nice hypercovers we have produced above, the pullback definition of $\Q_p^\N$ and \cref{lem:perf-covers} now reduce us to showing that all $X_\bullet\times_{(\ol{\DD}/\G_m(1))^\dR} Y_{\bullet'}$ are nicely coverable by $n$-truncated affines for some fixed $n\geq 0$. However, note that the maps $X_\bullet\rightarrow (\ol{\DD}/\G_m(1))^\dR$ and $Y_{\bullet'}\rightarrow (\ol{\DD}/\G_m(1))^\dR$ all naturally factor through $\ol{\DD}^\dR\rightarrow(\ol{\DD}/\G_m(1))^\dR$ and hence we obtain
\begin{equation*}
\begin{split}
X_\bullet\times_{(\ol{\DD}/\G_m(1))^\dR} Y_{\bullet'}&\cong X_\bullet\times_{\ol{\DD}^\dR} (\ol{\DD}^\dR\times_{(\ol{\DD}/\G_m(1))^\dR} \ol{\DD}^\dR)\times_{\ol{\DD}} Y_{\bullet'} \\
&\cong X_\bullet\times_{\ol{\DD}^\dR} (\G_m(1)^\dR\times \ol{\DD}^\dR)\times_{\ol{\DD}^\dR} Y_{\bullet'}) \\
&\cong X_\bullet\times_{\ol{\DD}^\dR} (\G_m(1)^\dR\times Y_{\bullet'})\;,
\end{split}
\end{equation*}
where the map $\G_m(1)^\dR\times Y_{\bullet'}\rightarrow \ol{\DD}^\dR$ is given by the composition
\begin{equation*}
\G_m(1)^\dR\times Y_{\bullet'}\rightarrow \G_m(1)^\dR\times \ol{\DD}^\dR\xrightarrow{\text{mult}} \ol{\DD}^\dR\;.
\end{equation*}
By another application of \cref{lem:perf-covers}, we may now check everything after passing to the nice hypercover $\G_m(1)\times (\G_m^\dagger)^{\bullet''}$ of $\G_m(1)^\dR$ from above. However, note that we can make each term in this hypercover act on $Y_{\bullet'}$ such that the map $Y_{\bullet'}\rightarrow\ol{\DD}^\dR$ becomes equivariant (e.g.\ by first projecting onto $\G_m(1)$ and then letting $\G_m(1)$ act by multiplication on one of the factors $\ol{\DD}$ in $Y_{\bullet'}$) and, after twisting by this action, we may compute the fibre product above by replacing the map $(\G_m(1)\times (\G_m^\dagger)^{\bullet''})\times Y_{\bullet'}\rightarrow \ol{\DD}^\dR$ by 
\begin{equation*}
(\G_m(1)\times (\G_m^\dagger)^{\bullet''})\times Y_{\bullet'}\xrightarrow{\mathrm{pr}_2} Y_{\bullet'}\rightarrow \ol{\DD}^\dR\;.
\end{equation*}
In particular, we may focus on the fibre products $X_\bullet\times_{\ol{\DD}^\dR} Y_{\bullet'}$ because the additional factors $\G_m(1)$ and $\G_m^\dagger$ do not play any further role for our argument due to $\O(\G_m(1))$ and $\O(\G_m^\dagger)$ being static and flat over $\Q_p$.

Now we can do a similar trick again: the maps $X_\bullet\rightarrow\ol{\DD}^\dR$ and $Y_{\bullet'}\rightarrow\ol{\DD}^\dR$ actually factor through $\ol{\DD}\rightarrow\ol{\DD}^\dR$, which tells us that
\begin{equation*}
X_\bullet\times_{\ol{\DD}^\dR} Y_{\bullet'}\cong X_\bullet\times_{\ol{\DD}} (\G_a^\dagger\times Y_{\bullet'})\;,
\end{equation*}
where the map $\G_a^\dagger\times Y_{\bullet'}\rightarrow\ol{\DD}$ is given by
\begin{equation*}
\G_a^\dagger\times Y_{\bullet'}\rightarrow \G_a^\dagger\times\ol{\DD}\xrightarrow{\mathrm{add}} \ol{\DD}\;.
\end{equation*}
Finally, note that (up to base changing to a closed annulus around $1$ in $1+\overcirc{\DD}$), the above pullback is of the form
\begin{equation*}
\left(\lim_{x\mapsto x^p} (1+\overcirc{\DD})\right)\times_{\ol{\DD}} (\G_a^\dagger\times \ol{\DD}\times \ol{\DD})\times (\,\dots)\;,
\end{equation*}
where the term in brackets is a finite product of affine Gelfand stacks which are flat and static over $\Q_p$ and hence is irrelevant for our argument.

Summarising, we are done once we know that, for some $n\geq 0$, the single (!) pullback
\begin{equation*}
\left(\lim_{x\mapsto x^p} (1+\overcirc{\DD})\right)\times_{\ol{\DD}} (\G_a^\dagger\times\ol{\DD}\times \ol{\DD})\;,
\end{equation*}
where the map $\G_a^\dagger\times\ol{\DD}\times \ol{\DD}\rightarrow \ol{\DD}$ is given by $(x, y, z)\mapsto x+yz$, is an affine $n$-truncated Gelfand stack after base changing to a closed annulus centered around $1$ inside $1+\overcirc{\DD}$. However, this just follows from the fact that the map $\G_a^\dagger\times \ol{\DD}\times\ol{\DD}\rightarrow \ol{\DD}$ is flat: On rings, this is given by
\begin{equation*}
\Q_p\langle T\rangle_{\leq 1}\rightarrow \Q_p\{X\}^\dagger\langle Y, Z\rangle_{\leq 1}\;, \hspace{0.3cm} T\mapsto X+YZ
\end{equation*}
and both the map
\begin{equation*}
\Q_p\langle T\rangle_{\leq 1}\rightarrow \Q_p\{X\}^\dagger\langle W\rangle_{\leq 1}\;, \hspace{0.3cm} T\mapsto X+W
\end{equation*}
as well as the map
\begin{equation*}
\Q_p\langle W\rangle_{\leq 1}\rightarrow \Q_p\langle Y, Z\rangle_{\leq 1}\;, \hspace{0.3cm} W\mapsto YZ
\end{equation*}
is flat. Indeed, for the first map, this follows from $\Q_p\{X\}^\dagger\langle W\rangle_{\leq 1}\cong \Q_p\{X\}^\dagger\langle X+W\rangle_{\leq 1}$ and flatness of $\Q_p\{X\}^\dagger$ over $\Q_p$ while, for the second map, we observe that, as a $\Q_p\langle W\rangle_{\leq 1}$-module, we have
\begin{equation*}
\Q_p\langle Y, Z\rangle_{\leq 1}\cong \Q_p\langle W, Y\rangle_{\leq 1}\oplus Z\Q_p\langle W, Z\rangle_{\leq 1}
\end{equation*}
and then the claim follows from flatness of $\Q_p\langle Y\rangle_{\leq 1}$ over $\Q_p$.
\end{proof}

\begin{rem}
A similar argument applies to $(\Q_p^\N)_{[r, s]}$ for any $0<r<s\leq\infty$. However, one has to account for the fact that, over $(0, p)$, the prismatisation of $\Q_p$ is isomorphic to the perfect prismatisation (up to modding out Frobenius) instead of the de Rham stack of $\Spd\Q_p\times\Spd\Q_p$ and use the presentation
\begin{equation*}
\Div^1\cong \lim_{x\mapsto x^p} (1+\overcirc{\mathbb{D}})\setminus \{1\} \,\big/\, \Q_p^{\times, \la}\;,
\end{equation*}
which we have not proved.
\end{rem}

To simplify the following proof, let us first introduce some more terminology.

\begin{defi}
Let $X$ be any Gelfand stack over $\Q_p$. The \emph{Cone stack} of $X$ is the Gelfand stack $X^{\Cone}$ over $\ol{\DD}/\G_m(1)$ defined by
\begin{equation*}
X^{\Cone}(\GSpec A\rightarrow \ol{\DD}/\G_m(1))=\{\text{maps $\GSpec\Cone(L\rightarrow A)\rightarrow X$}\}\;,
\end{equation*}
where $L\rightarrow A$ is the normed generalised Cartier divisor classified by the map $\GSpec A\rightarrow\ol{\DD}/\G_m(1)$.
\end{defi}

In this language, the definition of $\widetilde{D}$ over $\Q_p^\N$ as a pushout supplies a pullback diagram
\begin{equation*}
\begin{tikzcd}
X^\N\ar[r]\ar[d] & X^\prism\ar[d] \\
X^{\Cone}\ar[r] & (X^{\Cone})^\dR
\end{tikzcd}
\end{equation*}
for any derived Berkovich space $X$.

\begin{prop}
\label{prop:perf-coversmoothrigid}
Let $X$ be a Gelfand stack over $\Q_p$ equipped with a map $X\rightarrow \ol{\DD}^n$ which is the composition of a finite étale map with a closed localisation. Then $(X^\N)_{[r, s]}$ is nicely coverable for any $1<r<s<\infty$. 
\end{prop}
\begin{proof}
Since closed localisations have finite flat dimension by \cref{lem:perf-finflatdim} and the functor $X\mapsto X^\N$ is compatible with finite étale maps and closed localisations by \cref{prop:defis-etmaps} and \cref{prop:defis-openloc}, we may reduce to $X=\ol{\DD}^n$, from where we can further reduce to $X=\ol{\DD}$ by compatibility of $X\mapsto X^\N$ with limits.

Pulling back to the nice hypercover of $(\Q_p^\N)_{[r, s]}$ from the previous propsition, as we may do by \cref{lem:perf-covers}, the pushout definition of $\widetilde{D}$ yields a pullback square
\begin{equation*}
\begin{tikzcd}
\ol{\DD}^\N\ar[r]\ar[d] & \ol{\DD}^\prism\ar[d] \\
\ol{\DD}^{\Cone}\ar[r] & (\ol{\DD}^{\Cone})^\dR\nospacepunct{\;,}
\end{tikzcd}
\end{equation*}
where the cone stacks live over the divisor cut out by the universal section of $\ol{\DD}$, which we recall occurred as a factor in each of the terms of the nice hypercover of $(\Q_p^\N)_{[r, s]}$ we have given. To simplify the above, note that the diagram
\begin{equation*}
\begin{tikzcd}
\ol{\DD}^{\Cone}\ar[r]\ar[d] & (\ol{\DD}^{\Cone})^\dR\ar[d] \\
(\A^1)^{\Cone}\ar[r] & ((\A^1)^{\Cone})^\dR
\end{tikzcd}
\end{equation*}
is cartesian and that $(\A^1)^{\Cone}=\A^1/\A^1$, where the quotient is taken along multiplication by the universal section of $\ol{\DD}$. Finally, as we are working over $\lim_{x\mapsto x^p} (1+\overcirc{\DD})\rightarrow \Q_p^\prism$, we can write $\ol{\DD}^\prism=\ol{\DD}^{\dR, \flat\sharp}$, where the latter refers to (the elements of norm at most $1$ of) the universal untilt of the tilt of any perfectoid over $\lim_{x\mapsto x^p} (1+\overcirc{\DD})$.

Let us first explain how to control $\ol{\DD}^{\dR, \flat\sharp}$. For this, we use that
\begin{equation*}
\ol{\DD}^{\dR, \flat\sharp}=\Cone(\lim_F W^\dR\xrightarrow{1+[x^\flat]+\dots+[x^\flat]^{p-1}} \lim_F W^\dR)\;,
\end{equation*}
where $W=\GSpec \colim_n \Q_p\langle X_1, \dots, X_n\rangle_{\leq 1}$ is a version of the ring scheme of Witt vectors; in particular, note that this is again affine and flat as a $\Q_p$-algebra by \cref{lem:perf-vspflat}. 
Then first note that $\lim_F W^\dR\rightarrow \ol{\DD}^{\dR, \flat\sharp}$ is a cover with \v{C}ech nerve $\lim_F W^\dR\times (\lim_F W^\dR)^{\bullet}$. To get rid of the de Rham stacks, observe that both $W$ itself as well as the map $F: W\rightarrow W$ are $\dagger$-formally smooth, hence $\lim_F W\rightarrow\lim_F W^\dR$ is a surjection with \v{C}ech nerve
\begin{equation*}
\lim_F W\times (\lim_F \Ker(W\rightarrow W^\dR))^\bullet\;,
\end{equation*}
which is in fact a nice hypercover since $\Ker(W\rightarrow W^\dR)\cong \GSpec\colim_n \Q_p\{X_1, \dots, X_n\}^\dagger$ is affine and flat over $\Q_p$. Overall, by \cref{lem:perf-covers}, this provides a nice hypercover $X_\bullet$ of $\rightarrow\ol{\DD}^{\dR, \flat\sharp}$ by affine Gelfand stacks which are static and flat over $\Q_p$ and we may now pull back to this cover to proceed.

Now note that the maps $X_\bullet\rightarrow ((\A^1)^{\Cone})^\dR\cong (\A^1)^\dR/(\A^1)^\dR$ actually factor through $\A^1$ and, in particular, through $\A^1/\A^1$. Thus, we can write
\begin{equation*}
\begin{split}
\A^1/\A^1\times_{(\A^1)^\dR/(\A^1)^\dR} X_\bullet&\cong (\A^1/\A^1\times_{(\A^1)^\dR/(\A^1)^\dR} \A^1/\A^1)\times_{\A^1/\A^1} X_\bullet \\
&\cong (\G_a^\dagger/\G_a^\dagger\times \A^1/\A^1)\times_{\A^1/\A^1} X_\bullet\cong \G_a^\dagger/\G_a^\dagger\times X_\bullet\;,
\end{split}
\end{equation*}
where we have used that $(\A^1)^\dR/(\A^1)^\dR\cong (\A^1/\A^1)\,/\,(\G_a^\dagger/\G_a^\dagger)$. Finally, noting that $\G_a^\dagger/\G_a^\dagger$ is covered by $\G_a^\dagger$ with \v{C}ech nerve given by $\G_a^\dagger\times (\G_a^\dagger)^\bullet$ finishes the proof by \cref{lem:perf-covers} since $\O(\G_a^\dagger)$ is static and flat over $\Q_p$.
\end{proof}

Finally, we will later also be interested in the case $X=\GSpec\C_p$. 

\begin{prop}
\label{prop:perf-cpncover}
For any $1<r<s<\infty$, the stack $(\C_p^\N)_{[r, s]}$ is nicely coverable. 
\end{prop}
\begin{proof}
First note that $\C_p^\N$ lives over the stack 
\begin{equation*}
\Q_p^\N\times_{\Q_p^\prism} \C_p^\prism\cong (\ol{\DD}/\G_m(1)\times (\ol{\DD}/\G_m(1))^\dR)\times_{(\ol{\DD}/\G_m(1))^\dR} \C_p^\prism
\end{equation*}
and that $(\C_p^\prism)_{[r, s]}=(\Spd\C_p\times\Spd\Q_p)^\dR_{[r, s]}$. Now recall from \cite[§7.5]{dRFF} that $(\Spd\C_p\times\Spd\Q_p)^\dR$ is covered by $Y_{\C_p}$. More precisely, we have
\begin{equation*}
(\Spd\C_p\times\Spd\Q_p)^\dR=Y_{\C_p}/\Delta^\dagger\;,
\end{equation*}
where $\Delta^\dagger$ is the overconvergent neighbourhood of the diagonal in $Y_{\C_p}\times Y_{\C_p}$ and this in turn can be written as a limit of $Y_U$ over open neighbourhoods $U$ of the diagonal in $\Spd\C_p\times_{\Spd\F_q} \Spd\C_p$ for $q\rightarrow\infty$ a power of $p$. Thus, the \v{C}ech nerve of this cover will produce a nice hypercover $X_\bullet$ of $(\C_p^\prism)_{[r, s]}$ for us and we may check everything after base change along this cover. Furthermore, the maps $X_\bullet\rightarrow (\ol{\DD}/\G_m(1))^\dR$ all factor through $\ol{\DD}$ already and then the same argument as in the proof of \cref{prop:perf-zpncover} yields a nice hypercover $W_\bullet\rightarrow (\Q_p^\N)_{[r, s]}\times_{\Q_p^\prism} \C_p^\prism$ such that all $W_\bullet$ map to
\begin{equation*}
Y_{\C_p}\times_{\ol{\DD}} (\G_a^\dagger\times\ol{\DD}\times\ol{\DD})\;.
\end{equation*}

By definition, after pulling back to the hypercover $W_\bullet$, the stack $(\C_p^\N)_{[r, s]}$ is given by the pullback
\begin{equation*}
\begin{tikzcd}
(\C_p^\N)_{[r, s]}\ar[r]\ar[d] & *\ar[d] \\
\C_p^{\Cone}\ar[r] & (\C_p^{\Cone})^\dR\nospacepunct{\;,}
\end{tikzcd}
\end{equation*}
where the prismatisation one would normally expect in the top right corner has become a point since it is already present in the base $W_\bullet$. Indeed, the vertical map on the right is induced by the map
\begin{equation*}
Y_{\C_p}\times_{\ol{\DD}} (\G_a^\dagger\times\ol{\DD}\times\ol{\DD})\rightarrow (\C_p^{\Cone})^\dR
\end{equation*}
which takes a nilperfectoid $A$ over $Y_{\C_p}$ together with elements $\delta\in \Nil^\dagger(A)$ and $t, u\in A$ such that $\xi=ut+\delta$, where $\xi=p-[p^\flat]$ as usual, to the map
\begin{equation*}
\C_p=\O(Y_{\C_p})/^\mathbb{L}\xi\rightarrow A/^\mathbb{L} (ut+\delta)\rightarrow \ol{A}/^\mathbb{L} ut\rightarrow \ol{A}/^\mathbb{L} t\;.
\end{equation*}
Our task is to show that the fibre of $\C_p^{\Cone}$ over this fixed point of $(\C_p^{\Cone})^\dR$ is nicely coverable.

To this end, we fix the following choices for the rest of the proof: Let $(x_i)_{i\in I}$ be a transcendence basis of $\C_p$ over $\ol{\Q}_p$ and, for each $i\in I$, fix a Cauchy sequence $(x_{in})_n$ in $\ol{\Q}_p$ which converges to $x_i$ in $\C_p$. Letting $R$ denote the ring of global sections of $Y_{\C_p}\times_{\ol{\DD}} (\G_a^\dagger\times\ol{\DD}\times\ol{\DD})$, first note that, for $A=R$, all the quotients occurring above are non-derived and static. To continue, observe that the quotient $\ol{R}/t$ is equipped with a map from $\C_p$ by the above and we furthermore fix a lift $x_i'$ of $x_i$ from $\ol{R}/t$ to $R$ for each $i\in I$ and a lift $(d'_{in})_n$ of the nullsequence $(x_{in}-x_i)_n$ from $\ol{R}/t$ to $R$. Now note that the composition $\ol{\Q}_p\rightarrow\C_p\rightarrow \ol{R}/t$ lifts uniquely to a map $\ol{\Q}_p\rightarrow R/t$ since $\ol{\Q}_p$ is $\dagger$-formally smooth over $\Q_p$; in particular, we can make sense of the $x_{in}$ in $R/t$. Then 
\begin{equation*}
x_{in}-x_i'-d'_{in}\in\Nil^\dagger(R)/t
\end{equation*}
since the left-hand side vanishes in $\ol{R}/t$ by construction and we fix lifts $\delta'_{in}$ of these elements to $\Nil^\dagger(R)$. Finally, we set $x_{in}'\coloneqq x_i'+d'_{in}+\delta'_{in}$ for all $i\in I$ and $n\geq 0$.

Now let $A$ be any Gelfand ring over $Y_{\C_p}\times_{\ol{\DD}} (\G_a^\dagger\times\ol{\DD}\times\ol{\DD})$. The data of a lift of the given $A$-point of $(\C_p^{\Cone})^\dR$ to an $A$-point of $\C_p^{\Cone}$ can now be described as follows: First, the composition $\ol{\Q}_p\rightarrow\C_p\rightarrow \ol{A}/^\mathbb{L} t$ again lifts uniquely to a map $\ol{\Q}_p\rightarrow A/^\mathbb{L} t$ as above. Once we have fixed an extension of this map to a map
\begin{equation}
\label{eq:perf-qpxilift}
\ol{\Q}_p(x_i: i\in I)\rightarrow A/^\mathbb{L} t
\end{equation}
which lifts the given map $\ol{\Q}_p(x_i: i\in I)\rightarrow\C_p\rightarrow\ol{A}/^\mathbb{L} t$, we obtain a unique lift $\C_p\rightarrow A/^\mathbb{L} t$ of the given map $\C_p\rightarrow \ol{A}/^\mathbb{L} t$ extending the map (\ref{eq:perf-qpxilift}) since $\C_p$ is algebraic and hence $\dagger$-formally étale over $\ol{\Q}_p(x_i: i\in I)$. However, the choice of a map as in (\ref{eq:perf-qpxilift}) amounts to the following data:
\begin{enumerate}[label=(\arabic*)]
\item lifts $\widetilde{x}_i$ of the $x_i$ from $\ol{A}/^\mathbb{L} t$ to $A/^\mathbb{L} t$;

\item lifts $\widetilde{d}_{in}\xrightarrow{n\rightarrow\infty} 0$ of the nullsequences $x_{in}-x_i\xrightarrow{n\rightarrow\infty} 0$ from $\ol{A}/^\mathbb{L} t$ to $A/^\mathbb{L} t$ such that $\widetilde{d}_{in}=x_{in}-\widetilde{x}_i$ for all $i\in I$ and $n\geq 0$.
\end{enumerate}
Recalling that we have already fixed $x_i'$ and nullsequences $(d'_{in})_n$, this amounts to choosing $y_i=\widetilde{x}_i-x_i'\in \Nil^\dagger(A)/^\mathbb{L} t$ and nullsequences $(e_{in})_n=(\widetilde{d}_{in})_n-(d'_{in})_n$ in $\Nil^\dagger(A)/^\mathbb{L} t$ satisfying the relation $e_{in}+y_i=\delta_{in}'$. In other words, the fibre we seek is isomorphic to the pullback of the diagram
\begin{equation*}
\begin{tikzcd}
& *\ar[d] \\
(\G_a^\dagger/^\mathbb{L} t\G_a^\dagger\times \operatorname{Null}(\G_a^\dagger)/^\mathbb{L} t\operatorname{Null}(\G_a^\dagger))^I\ar[r] & (\G_a^\dagger/^\mathbb{L} t\G_a^\dagger)^{\mathbb{N}\times I}\nospacepunct{\;,}
\end{tikzcd}
\end{equation*}
where the bottom map is the sum map and the vertical map on the right is given by the $\delta_{in}'$. Here, $\operatorname{Null}(\G_a^\dagger)$ denotes the Gelfand stack sending $A$ to the anima of nullsequences in $\Nil^\dagger(A)$.

Now first note that $\G_a^\dagger/^\mathbb{L} t\G_a^\dagger$ admits a cover by $\G_a^\dagger$ and $\operatorname{Null}(\G_a^\dagger)/^\mathbb{L} t\operatorname{Null}(\G_a^\dagger)$ admits a cover by $\operatorname{Null}(\G_a^\dagger)$ and this produces a cover of the stack $(\G_a^\dagger/^\mathbb{L} t\G_a^\dagger\times \operatorname{Null}(\G_a^\dagger)/^\mathbb{L} t\operatorname{Null}(\G_a^\dagger))^I$ whose \v{C}ech nerve is given by 
\begin{equation*}
(\G_a^\dagger\times \operatorname{Null}(\G_a^\dagger))^I\times (\G_a^\dagger\times\operatorname{Null}(\G_a^\dagger))^{I\times\bullet}\;.
\end{equation*}
Since $\G_a^\dagger$ and $\operatorname{Null}(\G_a^\dagger)$ are both affine and represented by flat static $\Q_p$-algebras due to \cref{lem:perf-vspflat}, we are reduced to checking that
\begin{equation*}
(\G_a^\dagger\times\operatorname{Null}(\G_a^\dagger))^I\times_{(\G_a^\dagger/^\mathbb{L} t\G_a^\dagger)^{\mathbb{N}\times I}} *
\end{equation*}
is nicely coverable by \cref{lem:perf-covers} and \cref{lem:perf-bc}. For this, in turn, by the proof of \cref{lem:perf-bc}, it suffices to show that the fibre product 
\begin{equation*}
(\G_a^\dagger\times\operatorname{Null}(\G_a^\dagger))\times_{(\G_a^\dagger/^\mathbb{L} t\G_a^\dagger)^{\mathbb{N}}} *
\end{equation*}
is nicely coverable by affine Gelfand stacks corresponding to static flat $R$-algebras for each $i\in I$. 

To get a handle on this latter pullback, observe that the morphism $\G_a^\dagger\times\operatorname{Null}(\G_a^\dagger)\rightarrow (\G_a^\dagger/^\mathbb{L} t\G_a^\dagger)^{\mathbb{N}}$ naturally lifts to $(\G_a^\dagger)^{\mathbb{N}}$ and the same is true for the maps $*\rightarrow (\G_a^\dagger/^\mathbb{L} t\G_a^\dagger)^{\mathbb{N}}$ since the $\delta_{in}'$ already live in $\Nil^\dagger(R)$. Thus, the pullback in question is isomorphic to
\begin{equation*}
\begin{split}
(\G_a^\dagger\times\operatorname{Null}(\G_a^\dagger))&\times_{(\G_a^\dagger)^{\mathbb{N}}} ((\G_a^\dagger)^{\mathbb{N}}\times_{(\G_a^\dagger/^\mathbb{L} t\G_a^\dagger)^{\mathbb{N}}} (\G_a^\dagger)^{\mathbb{N}})\times_{(\G_a^\dagger)^{\mathbb{N}}} * \\
&\cong (\G_a^\dagger\times\operatorname{Null}(\G_a^\dagger))\times_{(\G_a^\dagger)^{\mathbb{N}}} ((\G_a^\dagger)^{\mathbb{N}}\times(\G_a^\dagger)^{\mathbb{N}})\times_{(\G_a^\dagger)^{\mathbb{N}}} * \\
&\cong (\G_a^\dagger\times\operatorname{Null}(\G_a^\dagger))\times_{(\G_a^\dagger)^{\mathbb{N}}} (\G_a^\dagger)^{\mathbb{N}}\;,
\end{split}
\end{equation*}
where the map $(\G_a^\dagger)^{\mathbb{N}}\rightarrow (\G_a^\dagger)^{\mathbb{N}}$ defining the last pullback is given by $z_n\mapsto \delta_{in}'+tz_n$ for $z_0, z_1, \dots$ denoting the coordinates on $(\G_a^\dagger)^{\mathbb{N}}$. However, this last pullback is actually affine relative to the base $Y_{\C_p}\times_{\ol{\DD}} (\G_a^\dagger\times\ol{\DD}\times\ol{\DD})$: It is represented by the base change of
\begin{equation}
\label{eq:perf-cpcoverflatalg}
\colim_N \Q_p\{y_i, (e_{in})_n^\wedge, z_0, z_1, \dots, z_N, \delta_{i0}', \delta_{i1}', \dots, \delta_{iN}'\}^\dagger/(e_{in}+y_i=\delta_{in}'+tz_n: n\leq N)
\end{equation}
along the map
\begin{equation*}
\GSpec \colim_N \Q_p\{\delta_{i0}', \delta_{i1}', \dots, \delta_{iN}'\}^\dagger\rightarrow Y_{\C_p}\times_{\ol{\DD}} (\G_a^\dagger\times\ol{\DD}\times\ol{\DD})\;.
\end{equation*}
Thus, we are done once we can show that (\ref{eq:perf-cpcoverflatalg}) is flat over $\colim_N \Q_p\{\delta_{i0}', \dots, \delta_{iN}'\}^\dagger$ for each $i\in I$.

For this, it suffices to show that 
\begin{equation*}
\Q_p\{y_i, (e_{in})_n^\wedge, z_0, z_1, \dots, z_N, \delta_{i0}', \delta_{i1}', \dots, \delta_{iN}'\}^\dagger/(e_{in}+y_i=\delta_{in}'+tz_n: n\leq N)
\end{equation*}
is flat over $\Q_p\{\delta_{i0}', \dots, \delta_{iN}'\}^\dagger$ for each $N\geq 0$ as flatness is stable under base change and filtered colimits. However, now observe that this algebra is isomorphic to 
\begin{equation*}
\Q_p\{y_i, (e_{in})_{n>N}^\wedge, z_0, z_1, \dots, z_N, \delta_{i0}', \delta_{i1}', \dots, \delta_{iN}'\}^\dagger
\end{equation*}
via the map sending $e_{in}$ to $\delta_{in}'+tz_n-y_i$ for $n\leq N$ and now the claim just follows from the fact that 
\begin{equation*}
\Q_p\{y_i, (e_{in})_{n>N}^\wedge, z_0, z_1, \dots, z_N\}^\dagger
\end{equation*}
is a flat $\Q_p$-algebra by \cref{lem:perf-vspflat}. This finally finishes the proof.
\end{proof}

\begin{cor}
\label{cor:perf-rigidgeomcover}
Let $X$ be a Gelfand stack over $\C_p$ equipped with a map $X\rightarrow \ol{\DD}_{\C_p}^n$ which is the composition of a finite étale map with a closed localisation. Then $(X^\N)_{[r, s]}$ is nicely coverable for any $1<r<s<\infty$. 
\end{cor}
\begin{proof}
This follows by imitating the proof of \cref{prop:perf-zpncover}, but with \cref{prop:perf-cpncover} replacing \cref{prop:perf-zpncover}.
\end{proof}

Throughout, we work in the setting of Gelfand stacks over $\Q_p$ of \cite{dRFF}. Moreover, the term ``arc-stack'' always refers to arc-stacks in the sense of \cite[Def.\ 4.1.5]{dRFF} and would perhaps usually be called ``untilted arc-stack''. As explained in loc.\ cit., there is a forgetful functor
\begin{equation*}
(-)^\diamond: \GelfStk_{\Q_p}\rightarrow \arcStk_{\Q_p}
\end{equation*}
which admits both adjoints. The right-adjoint is the \emph{analytic de Rham stack}
\begin{equation*}
(-)^\dR: \arcStk_{\Q_p}\rightarrow\GelfStk_{\Q_p}\;,
\end{equation*}
while the left-adjoint is obtained by Kan extending $\Spd\ol{A}\mapsto \GSpec\ol{A}$ for perfectoid $\Q_p$-algebras $\ol{A}$. We will often implicitly view an arc-stack over $\Q_p$ as a Gelfand stack over $\Q_p$ via this left-adjoint.

\subsection{Rational analytic prismatisation}

\begin{defi}
\begin{enumerate}[label=(\roman*)]
\item A $\Q_{p, \solid}$-algebra $A$ is called \emph{nilperfectoid} if its $\dagger$-reduction $\ol{A}\coloneqq A/\Nil^\dagger(A)$ is perfectoid. 
\item We denote the category of nilperfectoid $\Q_p$-algebras by $\NilPerfd$.
\end{enumerate}
\end{defi}

\begin{ex}
Let $\ol{A}$ be a perfectoid $\Q_p$-algebra. Then $A=\ol{A}[T]/T^2$ and $A=\ol{A}\{T\}^\dagger$ are both nilperfectoid $\Q_p$-algebras with $\dagger$-reduction $\ol{A}$.
\end{ex}

\begin{rem}
By \cite[Prop.\ 4.6.4]{dRFF}, nilperfectoid $\Q_p$-algebras form a basis of the descendable topology on Gelfand $\Q_p$-algebras, hence stacks on $\NilPerfd$ are equivalent to Gelfand stacks. In particular, in the following, we will pass freely between the two notions.
\end{rem}

Recall that, for any perfectoid $\Q_p$-algebra $\ol{A}$, there is an associated curve $Y_{\ol{A}}$ given by
\begin{equation*}
Y_{\ol{A}}=\{0<|p|<1\}\subseteq \Spd W(\ol{A}^{\flat\circ})[\tfrac{1}{[\pi^\flat]}]\;,
\end{equation*}
where $\pi\in \ol{A}$ is a perfectoid pseudouniformiser, i.e.\ admits $p$-power roots in $\ol{A}$, and we will in the following view this as the Gelfand stack
\begin{equation*}
Y_{\ol{A}}=\colim_k \GSpec W(\ol{A}^{\flat\circ})[\tfrac{1}{[\pi^\flat]}]\langle\tfrac{p^k}{[\pi^\flat]}, \tfrac{[\pi^\flat]^k}{p}\rangle\;.
\end{equation*}
Observe that the Witt vector Frobenius $F$ induces a Frobenius morphism $\phi: Y_{\ol{A}^\flat}\rightarrow Y_{\ol{A}^\flat}$, which is actually an isomorphism. Moreover, recall that there is a canonical map
\begin{equation*}
\iota: \GSpec \ol{A}\rightarrow Y_{\ol{A}}
\end{equation*}
coming from Fontaine's map
\begin{equation*}
\begin{split}
\theta: W(\ol{A}^{\flat\circ})&\rightarrow \ol{A}^\circ \\
\sum_{n\geq 0} [a_n]p^n&\mapsto \sum_{n\geq 0} [a_n^\sharp]p^n
\end{split}
\end{equation*}
and $\iota$ realises $\GSpec \ol{A}$ as a divisor inside $Y_{\ol{A}}$.

Furthermore assuming that $\ol{A}$ is a $\Q_p(p^{1/p^\infty})$-algebra, we have a ``radius map''
\begin{equation*}
\kappa: Y_{\ol{A}}\rightarrow (0, \infty)
\end{equation*}
given by
\begin{equation*}
x\mapsto \frac{\log |[p^\flat]|(x)|}{\log |p(x)|}\;,
\end{equation*}
i.e.\ the preimage of a closed interval $[r, s]\subseteq (0, \infty)$ is given by the locus 
\begin{equation*}
Y_{\ol{A}, [r, s]}\coloneqq \{|[p^\flat]|^{1/r}\leq |p|\leq |[p^\flat]|^{1/s}\}\subseteq Y_{\ol{A}}\;,
\end{equation*}
where $p^\flat=(p, p^{1/p}, p^{1/p^2}, \dots)$. Note that $\kappa$ intertwines the Frobenius with multiplication by $p$, i.e.\ $\kappa\phi=p\kappa$, and that the overconvergent neighbourhood of $\iota$ is precisely given by $Y_{\ol{A}, [1, 1]}$. By descent, we can extend these definitions to all perfectoid $\Q_p$-algebras.

\begin{defi}
Let $A$ be nilperfectoid. For any $r\geq 0$, we define the Gelfand ring $\O(Y_{A, [p^{-r}, p^r]})$ as the pullback
\begin{equation}
\label{eq:defis-pullbackya}
\begin{tikzcd}
\O(Y_{A, [p^{-r}, p^r]})\ar[r]\ar[d] & \O(Y_{\ol{A}, [p^{-r}, p^r]})\ar[d, "\prod_n \theta\circ F^{-n}"] \\
\prod_{n=0}^r A\ar[r] & \prod_{n=0}^r \ol{A}\nospacepunct{\;.}
\end{tikzcd}
\end{equation}
The \emph{relative $Y$-curve} of $A$ is defined as
\begin{equation*}
Y_A\coloneqq \colim_r \GSpec\O(Y_{A, [p^{-r}, p^r]})\;.
\end{equation*}
\end{defi}

\begin{rem}
\begin{enumerate}[label=(\alph*)]
\item Note that $Y_A$ is a derived Berkovich space with the same underlying topological space as $Y_{\ol{A}}$ since $\O(Y_{A, [p^{-r}, p^r]})$ is a $\dagger$-nilpotent thickening of $\O(Y_{\ol{A}, [p^{-r}, p^r]})$ for any $r\geq 0$.

\item There is a natural Frobenius morphism $\phi: Y_A\rightarrow Y_A$ induced by the Frobenius morphism of $Y_{\ol{A}}$ and the left shift
\begin{equation*}
\prod_{n=0}^r A\rightarrow \prod_{n=0}^{r-1} A\;.
\end{equation*}
Evidently, $\phi$ is an isomorphism away from the image of $\iota: \AnSpec A\rightarrow Y_A$.

\item The radius map $\kappa: Y_{\ol{A}}\rightarrow (0, \infty)$ extends to a map $\kappa: Y_A\rightarrow (0, \infty)$.
\end{enumerate}

\item By definition, we obtain a commutative diagram
\begin{equation*}
\begin{tikzcd}
\GSpec\ol{A}\ar[r]\ar[d, "\iota", swap] & \GSpec A\ar[d, "\iota"] \\
Y_{\ol{A}}\ar[r] & Y_A\nospacepunct{\;,}
\end{tikzcd}
\end{equation*}
but the map $\iota: \GSpec A\rightarrow Y_A$ will in general \emph{not} be a divisor anymore.
\end{rem}

\begin{prop}
\label{prop:defis-pushoutberk}
For any nilperfectoid $A$, the diagram
\begin{equation*}
\begin{tikzcd}[column sep=large]
\bigsqcup_{n\geq 0} \GSpec\ol{A}\ar[r, "\sqcup_{n\geq 0} \phi^n\circ\iota"]\ar[d] & Y_{\ol{A}}\ar[d] \\
\bigsqcup_{n\geq 0} \GSpec A\ar[r, "\sqcup_{n\geq 0} \phi^n\circ\iota"] & Y_A
\end{tikzcd}
\end{equation*}
is a pushout in the category of derived Berkovich spaces.
\end{prop}
\begin{proof}
Note that $Y_A$ and $Y_{\ol{A}}$ have the same underlying topological space and the same is true for $\GSpec A$ and $\GSpec \ol{A}$; in the particular, the diagram is a pushout on the level of topological spaces. Now we check the universal property, i.e.\ that the natural map
\begin{equation*}
\Hom(Y_A, X)\rightarrow \Hom(Y_{\ol{A}}, X)\times_{\prod_{n\geq 0} X(\ol{A})} \prod_{n\geq 0} X(A)
\end{equation*}
is an equivalence for any derived Berkovich space $X$. As the corresponding map on the level of topological spaces is an equivalence because the diagram in question induces a pushout of topological spaces, it suffices to check that the induced maps on the fibres of
\begin{equation*}
\Hom(Y_A, X)\rightarrow\Hom(|Y_A|, |X|)
\end{equation*}
are equivalences. 

Thus, fix a map $|Y_A|\rightarrow |X|$ of topological spaces. Taking an affine rational cover of $X$ whose pullback along this map can be refined by a strict rational cover of $|Y_A|$, we can reduce to the case where $X=\GSpec B$ is affine. In that case, however, the claim follows from the fact that each $\O(Y_{A, [p^{-r}, p^r]})$ was defined as a pullback of the corresponding rings of global functions on $Y_{\ol{A}, [p^{-r}, p^r]}$ and $\bigsqcup_{n=0}^r \GSpec A$ over the ring of global functions on $\bigsqcup_{n=0}^r \GSpec\ol{A}$.
\end{proof}

\begin{lem}
\label{lem:defis-frobpushout}
For any nilperfectoid $A$, the diagram
\begin{equation*}
\begin{tikzcd}
\GSpec\ol{A}\ar[r]\ar[d] & \GSpec A\ar[d, "\iota"] \\
Y_A\ar[r, "\phi"] & Y_A
\end{tikzcd}
\end{equation*}
is a pushout in the category of derived Berkovich spaces, where the map $\GSpec\ol{A}\rightarrow Y_A$ is given by the postcomposition of $\phi^{-1}\circ\iota: \GSpec\ol{A}\rightarrow Y_{\ol{A}}$ with $Y_{\ol{A}}\rightarrow Y_A$.
\end{lem}
\begin{proof}
Due to Frobenius inducing an isomorphism $\O(Y_{\ol{A}, [p^{-r}, p^r]})\xrightarrow{\cong} \O(Y_{\ol{A}, [p^{-r-1}, p^{r-1}]})$ for all $r\geq 1$, the definition of $Y_A$ yields corresponding pullback diagrams
\begin{equation*}
\begin{tikzcd}
\O(Y_{A, [p^{-r}, p^r]})\ar[r, "\phi"]\ar[d, "\iota"] & \O(Y_{A, [p^{-r-1}, p^{r-1}]})\ar[d] \\
A\ar[r] & \ol{A}\nospacepunct{\;.}
\end{tikzcd}
\end{equation*}
Now one deduces the desired pushout by the same argument as in the proof of \cref{prop:defis-pushoutberk}.
\end{proof}

Now recall the notion of degree $1$ Cartier divisors on $Y_{\ol{A}}$ from \cite[Def.\ II.1.19]{FarguesScholze}: These are defined to be exactly those Cartier divisors arising as maps
\begin{equation*}
\GSpec \ol{A}^{\flat\sharp}\rightarrow Y_{\ol{A}}
\end{equation*}
for an untilt $A^{\flat\sharp}$ of $\ol{A}^\flat$. Moreover, they are always cut out by an element of the form $p+a[\pi^\flat]$ for some $a\in W(\ol{A}^{\flat\circ})$ and some perfectoid pseudouniformiser $\pi\in\ol{A}$; any such element is called \emph{distinguished}.

\begin{defi}
Let $A$ be nilperfectoid. A Cartier divisor $D\subseteq Y_A$ is called \emph{of degree $1$} if its pullback $D\times_{Y_A} Y_{\ol{A}}$ is a degree $1$ Cartier divisor in $Y_{\ol{A}}$.
\end{defi}

Here, by a Cartier divisor $D$ on $Y_A$, we mean a pullback diagram
\begin{equation*}
\begin{tikzcd}
D\ar[r]\ar[d] & Y_A\ar[d] \\
*/\G_m\ar[r] & \A^1/\G_m\nospacepunct{\;,}
\end{tikzcd}
\end{equation*}
a notion that would perhaps usually be called a ``generalised'' Cartier divisor.

\begin{rem}
Note that any degree $1$ Cartier divisor $D\subseteq Y_A$ is affine. Indeed, any degree $1$ Cartier divisor in $Y_{\ol{A}}$ is fully contained in $Y_{\ol{A}, [r, s]}$ for some $[r, s]\subseteq (0, \infty)$ and hence $D$ is already a Cartier divisor in $Y_{A, [r, s]}$, which is affine.
\end{rem}

\begin{ex}
\label{ex:defis-drdivisor}
The composition of the map $\phi^{-1}\circ\iota: \GSpec\ol{A}\rightarrow Y_{\ol{A}}$ with $Y_{\ol{A}}\rightarrow Y_A$ is a degree $1$ Cartier divisor. This is because, on an open neighbourhood of the image of this composition in $Y_A$, the curve $Y_A$ is isomorphic to $Y_{\ol{A}}$.
\end{ex}

\begin{lem}
Any map $D\rightarrow D'$ between degree $1$ Cartier divisors on $Y_A$ is an isomorphism.
\end{lem}
\begin{proof}
Any such map induces a map between the ideal sheaves $\O(-D')\rightarrow\O(-D)$, which becomes an isomorphism after pulling back to $Y_{\ol{A}}$ by the corresponding statement for degree $1$ Cartier divisors on $Y_{\ol{A}}$. As $Y_{\ol{A}}\rightarrow Y_A$ is a $\dagger$-nilpotent thickening of derived Berkovich spaces, this implies that $\O(-D')\rightarrow\O(-D)$ is an isomorphism itself by \cite[Prop.\ 2.6.16.(2)]{dRStack}, so we are done. (Note that we have implicitly used that Gelfand rings are Fredholm by \cite[Prop.\ 3.3.5]{dRFF} to identify line bundles on affine Gelfand stacks with invertible modules over their rings of global functions.)
\end{proof}

\begin{defi}
The (rational) \emph{analytic prismatisation} of $\Q_p$ is the stack $\Q_p^\prism$ on $\NilPerfd$ defined by
\begin{equation*}
\Q_p^\prism(A)=\{\text{degree $1$ Cartier divisors $D\rightarrow Y_A$}\}\;.
\end{equation*}
For an arbitrary Gelfand stack $X$ over $\Q_p$, we define its (rational) analytic prismatisation as the stack $X^\prism$ over $\Q_p^\prism$ on $\NilPerfd$ given by
\begin{equation*}
X^\prism(\GSpec A\xrightarrow{D} \Q_p^\prism)=\{\text{maps $D\rightarrow X$}\}\;.
\end{equation*}
\end{defi}

\begin{rem}
The Frobenius on $Y_A$ induces a Frobenius $\phi: \Q_p^\prism\rightarrow\Z_p^\prism$ and hence also a Frobenius $\phi: X^\prism\rightarrow X^\prism$ for any $X$.
\end{rem}

\begin{ex}
The degree $1$ divisor $\phi^{-1}\circ\iota: \GSpec\ol{A}\rightarrow Y_A$ from \cref{ex:defis-drdivisor} yields a map $i_\dR: \GSpec\Q_p\rightarrow\Q_p^\prism$ called the \emph{de Rham map}.
\end{ex}

Clearly, pulling back along $\GSpec A\rightarrow Y_A$ induces a map
\begin{equation*}
\mu: \Q_p^\prism\rightarrow \A^1/\G_m\;.
\end{equation*}
However, we can actually do slightly better.

\begin{defi}
The \emph{norm-$1$-torus} $\G_m(1)$ over $\Q_p$ is defined as the locus $\{|T|=1\}$ in $\A^1$, i.e.\
\begin{equation*}
\G_m(1)=\GSpec \colim_{\epsilon>0} \Q_p\langle \pi^\epsilon T, \pi^\epsilon T^{-1}\rangle\;.
\end{equation*}
Moreover, we write $\DD$ for the open unit disk over $\Q_p$.
\end{defi}

\begin{lem}
\label{lem:defis-mu}
The map $\mu: \Q_p^\prism\rightarrow \A^1/\G_m$ can be refined to a map $\mu: \Z_p^\prism\rightarrow \DD/\G_m(1)$.
\end{lem}
\begin{proof}
As being a unit on $Y_A$ may be checked after pullback to $Y_{\ol{A}}$, again by \cite[Prop.\ 2.6.16.(2)]{dRStack}, and the pullback of $D$ to $Y_{\ol{A}}$ is cut out by a distinguished element, we may choose a global section of $\O(D)$ cutting out $D$ which pulls back to a distinguished element on $Y_{\ol{A}}$. We claim that, after pullback along $\iota: \GSpec A\rightarrow Y_A$, all such global sections differ only by units of norm $1$ and have themselves norm less than $1$. This will then define the desired map $\Q_p^\prism\rightarrow \DD/\G_m(1)$.

To prove our claim, it will be enough to show:
\begin{enumerate}[label=(\alph*)]
\item Any distinguished element $p+a[\pi^\flat]$ in $W(\ol{A}^{\flat\circ})$ maps to a topologically nilpotent element under Fontaine's map $\theta: W(\ol{A}^{\flat\circ})\rightarrow \ol{A}^\circ$.
\item If two distinguished elements of $W(\ol{A}^{\flat\circ})$ differ by a global unit $u$ on $Y_{\ol{A}}$, the image of $u$ under Fontaine's map has norm $1$.
\end{enumerate}
For (a), note that the kernel of $\theta$ is generated by a distinguished element $p+b[\varpi^\flat]$ as well and the difference $(p+a[\pi^\flat])-(p+b[\varpi^\flat])=p[\pi^\flat]-b[\varpi^\flat]$ is already topologically nilpotent in $W(\ol{A}^{\flat\circ})$, whence the claim. For (b), we observe that two distinguished elements $p+a[\varpi^\flat]$ and $p+b[\varpi^\flat]$ in $W(\ol{A}^{\flat\circ})$ differing by a global unit $u$ on $Y_{\ol{A}}$ will be nonvanishing in a rational neighbourhood of $V([\pi^\flat])$ in $\Spa W(\ol{A}^{\flat\circ})$ as they cut out a degree $1$ divisor of $Y_{\ol{A}}$. Thus, we can even assume that $u$ is a unit in $W(\ol{A}^{\flat\circ})$ and this will of course imply that $\theta(u)\in (\ol{A}^{\flat\circ})^\times$, i.e.\ $\theta(u)$ has norm $1$.
\end{proof}

\begin{defi}
The \emph{Hodge--Tate divisor} $\Q_p^\HT$ inside $\Q_p^\prism$ is defined by the pullback diagram
\begin{equation*}
\begin{tikzcd}
\Q_p^\HT\ar[r]\ar[d] & \Q_p^\prism\ar[d] \\
*/\G_m(1)\ar[r] & \DD/\G_m(1)
\end{tikzcd}
\end{equation*}
and we let $\cal{I}\in\Pic(\Q_p^\prism)$ denote its ideal sheaf.
\end{defi}

\begin{prop}
The map $\Q_p^\HT\rightarrow \GSpec\Q_p$ is a gerbe banded by $\G_m^\dagger$. In fact, we have
\begin{equation*}
\Q_p^\HT\cong \GSpec\Q_p(\zeta_{p^\infty})/\Z_p^{\times, \la}\;.
\end{equation*}
\end{prop}
\begin{proof}
After base changing along $\Q_p\rightarrow \Q_p(\zeta_{p^\infty})$, the divisor $\iota: \GSpec\ol{A}\subseteq Y_{\ol{A}}$ can be cut out by $\xi=1+[\epsilon^{1/p}]+\dots+[\epsilon^{(p-1)/p}]\in W(\ol{A}^{\flat\circ})$, where $\epsilon=(1, \zeta_p, \zeta_{p^2}, \dots)\in\ol{A}^{\flat\circ}$ as usual. Gluing this global section on $Y_{\ol{A}}$ together with the zero section on $\GSpec A$ then yields a global section of $Y_A$ whose corresponding divisor naturally contains $\GSpec A$; i.e., we obtain a $\Q_p(\zeta_{p^\infty})$-valued point $\eta$ of $\Q_p^\HT$. As $\xi$ is a nonzerodivisor, its automorphism group is trivial while the automorphism group of $0\in A$ is of course given by $\G_m$, and this yields $\Aut(\eta)\cong \G_m^\dagger$, as claimed. Noting that $\eta$ covers $\Q_p^\HT$ shows that $\Q_p^\HT$ is indeed a $\G_m^\dagger$-gerbe.

For the second assertion, one checks that the map $\eta$ descends to a map 
\begin{equation*}
\GSpec\Q_p(\zeta_{p^\infty})/\Z_p^{\times, \la}\rightarrow\Q_p^\HT\;.
\end{equation*}
As this is in fact a morphism of $\G_m^\dagger$-gerbes, it must be an isomorphism and we are done.
\end{proof}

\begin{rem}
The map $\eta$ does \emph{not} further descend to a map 
\begin{equation*}
\GSpec\Q_p\cong \GSpec\Q_p(\zeta_{p^\infty})/\Z_p^{\times, \sm}\rightarrow\Q_p^\HT\;,
\end{equation*}
i.e.\ it does not yield a splitting of the gerbe. In fact, one cannot split the gerbe $\Q_p^\HT$ over any finite extension of $\Q_p$.
\end{rem}

Note that, for any $A$-point $D\subseteq Y_A$ of $\Q_p^\HT$, there is a map $\GSpec A\rightarrow D$ by definition and hence we obtain a natural map $X^\HT\rightarrow X$ for any $X$.

\begin{prop}
If $X$ is a rigid smooth derived Berkovich space over $\Q_p$, then $X^\HT\rightarrow X$ is a gerbe banded by $\cal{T}_X^\dagger\{1\}$. Here, we set $\cal{T}_X\coloneqq \GSpec_X \Sym^\bullet_X \mathbb{L}_X$.
\end{prop}
\begin{proof}
Recall that the divisor $D$ corresponding to the map $\eta: \GSpec\Q_p(\zeta_{p^\infty})\rightarrow \Q_p^\HT$ from the proof of the previous proposition must be affine and its ring of global functions sits in a pullback diagram
\begin{equation*}
\begin{tikzcd}
\O(D)\ar[r]\ar[d] & \ol{A}\ar[d] \\
A/^\mathbb{L} 0\ar[r] & \ol{A}/^\mathbb{L} 0\nospacepunct{\;.}
\end{tikzcd}
\end{equation*}
This implies $D\cong \GSpec(A\oplus \Nil^\dagger(A)[1])$ and now the claim follows by derived deformation theory since $X$ is $\dagger$-formally smooth. Note that we have to keep track of the Breuil--Kisin twist $(-)\{1\}$ in general, but this is trivialised in the above case due to the base change to $\Q_p(\zeta_{p^\infty})$.
\end{proof}

\begin{defi}
The \emph{Breuil--Kisin line bundle} $\O\{1\}$ on $\Q_p^\prism$ is defined by
\begin{equation*}
\O\{1\}\coloneqq \bigotimes_{n\geq 0} (\phi^n)^*\cal{I}\;.
\end{equation*}
We also use $\O\{1\}$ to denote the pullback of this line bundle to any $X^\prism$.
\end{defi}

Note that the infinite tensor product occurring in the definition above makes sense: Indeed, the divisor $\Q_p^\HT\subseteq\Q_p^\prism$ is supported over $\{1\}$ under the radius map $\Q_p^\prism\rightarrow (0, \infty)$ defined below and hence its inverse images under Frobenius are supported over $1/p, 1/p^2, \dots$. In particular, over each $[r, s]\subseteq (0, \infty)$, only finitely many of the line bundles $(\phi^n)^*\cal{I}$ are nontrivial and hence the infinite tensor product above is actually a finite tensor product locally on $\Q_p^\prism$.

\subsection{Relation with the Hyodo--Kato stack and the perfect prismatisation}

Note that, for any Gelfand stack $X$, there is a natural map
\begin{equation*}
X^\prism\rightarrow (X^\prism)^\dR
\end{equation*}
explicitly given by sending a degree $1$ Cartier divisor $D\subseteq Y_A$ equipped with a map $D\rightarrow X$ to the pullback $\ol{D}\coloneqq D\times_{Y_A} Y_{\ol{A}}\subseteq Y_{\ol{A}}$ together with the composite map $\ol{D}\rightarrow D\rightarrow X$.

\begin{prop}
\label{prop:defis-xprismdr}
Let $X$ be any Gelfand stack over $\Q_p$. Then there is a natural map
\begin{equation*}
X^\diamond\times\Spd\Q_p\rightarrow X^\prism
\end{equation*}
which induces an isomorphism on de Rham stacks.
\end{prop}
\begin{proof}
To construct the desired map, we just have to explain how any $\ol{A}$-point of $X^\diamond\times\Spd\Q_p$ gives rise to an $\ol{A}$-point of $X^\prism$. In fact, we will show that these two data are actually equivalent, which then also produces the desired isomorphism on de Rham stacks. To see this, note that giving an $\ol{A}$-point of $(X^\prism)^\dR$ is equivalent to giving an untilt $\ol{A}^{\flat\sharp}$ of $\ol{A}^\flat$ together with a map $\GSpec \ol{A}^{\flat\sharp}\rightarrow X$. However, this is exactly the same as an $\ol{A}^\flat$-point of $X^\diamond\times\Spd\Q_p$ lying above $\Spd \ol{A}\rightarrow\Spd\Q_p$, i.e.\ an $\ol{A}$-point of $X^\diamond\times\Spd\Q_p$, so we are done.
\end{proof}

\begin{lem}
The map $Y_{\ol{A}}\rightarrow Y_A$ induces an isomorphism $Y_{\ol{A}, (0, 1)}\cong Y_{A, (0, 1)}$, where we recall that the radius map is taken with respect to $p$.
\end{lem}
\begin{proof}
This just follows from the fact that the divisor $\iota: \GSpec\ol{A}\rightarrow Y_{\ol{A}}$ sits at $1\in (0, \infty)$ and hence its Frobenius translates sit at $p, p^2, p^3, \dots$. By definition, these are the only loci where the map $Y_{\ol{A}}\rightarrow Y_A$ is not an isomorphism.
\end{proof}

Note that, for any Gelfand stack $X$ over $\Q_p$, the stack $X^\prism$ is equipped with a radius map as well: Indeed, we have maps
\begin{equation*}
X^\prism\rightarrow \Q_p^\prism\rightarrow (\Q_p^\prism)^\dR\cong (\Spd\Q_p\times\Spd\Q_p)^\dR\rightarrow (0, \infty)\;,
\end{equation*}
where the last map is the radius map with respect to the pseudouniformiser $p$, as always.

\begin{prop}
\label{prop:defis-prismffdr}
Let $X$ be any Gelfand stack over $\Q_p$. The natural map 
\begin{equation*}
X^\prism\rightarrow (X^\prism)^\dR\cong (X^\diamond \times\Spd\Q_p)^\dR
\end{equation*}
is an isomorphism over the locus $(1, \infty)$.
\end{prop}
\begin{proof}
We are immediately reduced to the case $X=\GSpec \Q_p$, where the claim follows from the previous lemma by tracing through the identifications.
\end{proof}

\comment{
We want to explicitly point out that the preceding proposition in particular gives a complete description of $X^\prism$ in case $X$ is over $\F_p$. Indeed, in this case, the diamond $X^\diamond\times\Spd\Q_p$ is supported over $\{\infty\}\subseteq (0, \infty]$ and hence we obtain

\begin{cor}
\label{cor:defis-prismfp}
Let $X$ be a Gelfand stack over $\cal{N}_{\F_p}$. Then $X^\prism\cong (X^\diamond\times\Spd\Q_p)^\dR$.
\end{cor}
}

Note that the restriction of $\O\{1\}$ to the locus $(1, \infty)$ becomes trivial since all Frobenius pullbacks of $\Q_p^\HT$ have support in $(0, 1]$. Let us now also say something about the locus over $(0, p)$.

\begin{defi}
For any nilperfectoid $A$, we define the \emph{modified Fargues--Fontaine curve} $\FF_A$ over $A$ as the pushout
\begin{equation*}
\begin{tikzcd}
\GSpec \ol{A}\ar[r]\ar[d] & \FF_{\ol{A}}\ar[d] \\
\GSpec A\ar[r] & \FF_A
\end{tikzcd}
\end{equation*}
in the category of derived Berkovich spaces, where $\FF_{\ol{A}}$ denotes the usual Fargues--Fontaine curve of $\ol{A}$ seen as a Gelfand stack.
\end{defi}

\begin{rem}
By \cref{prop:defis-pushoutberk}, we could equivalently define the rings of global functions on affine pieces of $\FF_A$ by a pullback of rings and then pass to Gelfand stacks.
\end{rem}

As in the case of $Y_A$, we define degree $1$ Cartier divisors on $\FF_A$ to be exactly those whose pullback to $\FF_{\ol{A}}$ has degree $1$. 

\begin{defi}
The \emph{perfect analytic prismatisation} of $\Q_p$ is the stack $\Q_p^{\Div^1}=\Div^1$ on $\NilPerfd$ defined by
\begin{equation*}
\Div^1(A)=\{\text{degree $1$ Cartier divisors $D\rightarrow \FF_A$}\}\;.
\end{equation*}
For an arbitrary Gelfand stack $X$ over $\Q_p$, we define its perfect analytic prismatisation as the stack $X^{\Div^1}$ over $\Div^1$ on $\NilPerfd$ given by
\begin{equation*}
X^{\Div^1}(\GSpec A\xrightarrow{D} \Div^1)=\{\text{maps $D\rightarrow X$}\}\;.
\end{equation*}
\end{defi}

\begin{lem}
The map $\phi: Y_A\rightarrow Y_A$ induces an isomorphism $Y_{A, (1, \infty)}\cong Y_{A, (p, \infty)}$.
\end{lem}
\begin{proof}
By \cref{prop:defis-pushoutberk}, there are pushout diagrams
\begin{equation*}
\begin{tikzcd}[column sep=large]
\bigsqcup_{n\geq 0} \GSpec \ol{A}\ar[r, "\sqcup_{n\geq 0} \phi^n\circ\iota"]\ar[d] & Y_{\ol{A}, (1, \infty)}\ar[d] \\
\bigsqcup_{n\geq 0} \GSpec A\ar[r, "\sqcup_{n\geq 0} \phi^n\circ\iota"] & Y_{A, (1, \infty)}\nospacepunct{\;,}
\end{tikzcd}
\hspace{0.5cm}
\begin{tikzcd}[column sep=large]
\bigsqcup_{n\geq 1} \GSpec \ol{A}\ar[r, "\sqcup_{n\geq 1} \phi^n\circ\iota"]\ar[d] & Y_{\ol{A}, (p, \infty)}\ar[d] \\
\bigsqcup_{n\geq 1} \GSpec A\ar[r, "\sqcup_{n\geq 1} \phi^n\circ\iota"] & Y_{A, (p, \infty)}
\end{tikzcd}
\end{equation*}
in the category of derived Berkovich spaces and the map $\phi$ is exactly given by gluing together $\bigsqcup_{n\geq 0} \GSpec A\xrightarrow{\cong} \bigsqcup_{n\geq 1} \GSpec A$ and $\phi: Y_{\ol{A}}\xrightarrow{\cong} Y_{\ol{A}}$, which implies the claim.
\end{proof}

\begin{prop}
\label{prop:defis-xprismxdiv1}
For any Gelfand stack $X$, the Frobenius $\phi: X^\prism\rightarrow X^\prism$ induces an isomorphism $X^ \prism_{(0, 1)}\cong X^\prism_{(0, p)}$. Moreover, we have
\begin{equation*}
X^\prism_{(0, p)}/\phi\cong X^{\Div^1}\;,
\end{equation*}
where the left-hand side denotes the quotient by the equivalence relation induced by the ``partial Frobenius'' $\phi: X^\prism_{(0, 1)}\rightarrow X^\prism_{(0, p)}$. 
\end{prop}
\begin{proof}
Once again, we can reduce to $X=\GSpec\Q_p$, in which case the claim follows from the previous lemma by tracing through the identifications.
\end{proof}

\subsection{Nygaardifications of smooth derived Berkovich spaces}

\subsubsection{The Nygaardification of $\Q_p$}

\begin{lem}
Over $(1, \infty)\subseteq (0, \infty)$, the map 
\begin{equation*}
[p]_q: \lim_{q\mapsto q^p} (1+\overcirc{\DD})\setminus\{1\}\rightarrow \overcirc{\DD}
\end{equation*}
is $\dagger$-formally étale.
\end{lem}
\begin{proof}
It suffices to check the claim over $(1, p^m)\subseteq (0, \infty)$ for each $m\geq 1$. Moreover, as the $p$-th power map on $1+\overcirc{\DD}$ is $\dagger$-formally étale, it suffices to check that the map
\begin{equation*}
\begin{split}
\{p^{-1/(p-1)}<|x-1|<p^{-1/(p^m(p-1))}\}\subseteq 1+\overcirc{\DD}&\rightarrow \overcirc{\DD} \\
x&\mapsto 1+x^{p^{m+1}}+x^{2p^{m+1}}+\dots+x^{(p-1)p^{m+1}}
\end{split}
\end{equation*}
is $\dagger$-formally étale, where we have written $x\coloneqq q^{1/p^{m+1}}$ for the coordinate on the $(m+2)$-nd copy of $1+\overcirc{\DD}$, and this amounts to checking that the derivative
\begin{equation*}
p^{m+1}x^{p^{m+1}-1}+2p^{m+1}x^{2p^{m+1}-1}+\dots+(p-1)p^{m+1}x^{(p-1)p^{m+1}-1}
\end{equation*}
is nonvanishing. As $p^{m+1}$ and $x$ are units, this amounts to checking that
\begin{equation}
\label{eq:pres-derivativepq}
(1+t)^{p^{m+1}}+2(1+t)^{2p^{m+1}}+\dots+(p-1)(1+t)^{(p-1)p^{m+1}}=0
\end{equation}
has no solution for $p^{-1/(p-1)}<|t|<p^{-1/(p^m(p-1))}$, where we have written $t\coloneqq x-1$.

To see this, first observe that expanding the left-hand side of (\ref{eq:pres-derivativepq}) as a polynomial in $t$ yields a constant term $p(p-1)/2$, which has norm at least $1/p$. We are going to argue that all other terms have norm less than $1/p$, which clearly suffices to establish the claim. Indeed, this is clear for all terms of degree at least $p^{m+1}$ since $|t|^{p^{m+1}}<p^{-p/(p-1)}<p^{-1}$. Moreover, the term of degree $0<d<p^{m+1}$ is given by
\begin{equation*}
t^d\sum_{k=1}^{p-1} k\binom{kp^{m+1}}{d}
\end{equation*}
and observing that each $\binom{kp^{m+1}}{d}$ is divisible by $p$ due to $0<d<p^{m+1}$, we conclude that this has norm at most $|t|^dp^{-1}<p^{-1}$, which is what we wanted to show.
\end{proof}

}

\newpage

\bibliographystyle{alpha}
\bibliography{References}

\newcommand{\etalchar}[1]{$^{#1}$}
\begin{thebibliography}{ABLB{\etalchar{+}}25}

\bibitem[ABLB{\etalchar{+}}25]{dRFF}
Johannes Anschütz, Guido Bosco, Arthur-César Le~Bras, Juan~Esteban
  Rodríguez~Camargo, and Peter Scholze.
\newblock Analytic de {R}ham stacks of {F}argues--{F}ontaine curves, 2025.

\bibitem[ALBRCS]{AnPrism}
Johannes Anschütz, Arthur-César Le~Bras, Juan~Esteban Rodríguez~Camargo, and
  Peter Scholze.
\newblock Analytic prismatisation over $\mathbb{Q}_p$.
\newblock To appear.

\bibitem[Aok26]{AokiThesis}
Ko~Aoki.
\newblock {B}erkovich $2$-motives and normed ring stacks, 2026.

\bibitem[Bei13]{Beilinson}
Alexander Beilinson.
\newblock On the crystalline period map.
\newblock {\em Camb. J. Math.}, 1(1):1--51, 2013.

\bibitem[Ber08]{BergerBPaires}
Laurent Berger.
\newblock Construction de {$(\phi,\Gamma)$}-modules: repr\'esentations
  {$p$}-adiques et {$B$}-paires.
\newblock {\em Algebra Number Theory}, 2(1):91--120, 2008.

\bibitem[Ber16]{BergerMultivariable}
Laurent Berger.
\newblock Multivariable {$(\varphi,\Gamma)$}-modules and locally analytic
  vectors.
\newblock {\em Duke Math. J.}, 165(18):3567--3595, 2016.

\bibitem[BGV24]{BindaGallauerVezzani}
Federico Binda, Martin Gallauer, and Alberto Vezzani.
\newblock Motivic monodromy and $p$-adic cohomology theories, 2024.

\bibitem[Bha22]{FGauges}
Bhargav Bhatt.
\newblock Prismatic {$F$}-gauges, 2022.

\bibitem[BK90]{BlochKato}
Spencer Bloch and Kazuya Kato.
\newblock {$L$}-functions and {T}amagawa numbers of motives.
\newblock In {\em The {G}rothendieck {F}estschrift, {V}ol.\ {I}}, volume~86 of
  {\em Progr. Math.}, pages 333--400. Birkh\"auser Boston, Boston, MA, 1990.

\bibitem[BL22a]{APC}
Bhargav Bhatt and Jacob Lurie.
\newblock Absolute prismatic cohomology, 2022.

\bibitem[BL22b]{PFS}
Bhargav Bhatt and Jacob Lurie.
\newblock The prismatization of $p$-adic formal schemes, 2022.

\bibitem[BMS18]{IntegralpAdicHT}
Bhargav Bhatt, Matthew Morrow, and Peter Scholze.
\newblock Integral {$p$}-adic {H}odge theory.
\newblock {\em Publ. Math. Inst. Hautes \'Etudes Sci.}, 128:219--397, 2018.

\bibitem[BMS19]{THHAndPAdicHT}
Bhargav Bhatt, Matthew Morrow, and Peter Scholze.
\newblock Topological {H}ochschild homology and integral {$p$}-adic {H}odge
  theory.
\newblock {\em Publ. Math. Inst. Hautes \'Etudes Sci.}, 129:199--310, 2019.

\bibitem[Bos23a]{BoscoDrinfeldSpace}
Guido Bosco.
\newblock On the $p$-adic pro-\'etale cohomology of {D}rinfeld symmetric
  spaces, 2023.

\bibitem[Bos23b]{BoscopAdic}
Guido Bosco.
\newblock Rational $p$-adic {H}odge theory for rigid-analytic varieties, 2023.

\bibitem[CN20]{padicComparisons}
Pierre Colmez and Wies{\l}awa Nizio\l.
\newblock On $p$-adic comparison theorems for rigid analytic varieties, {I}.
\newblock {\em Münster J.\ Math.}, 13, 2020.

\bibitem[CN25]{BasicComparison}
Pierre Colmez and Wies{\l}awa Nizio\l.
\newblock On the cohomology of {$p$}-adic analytic spaces, {I}: {T}he basic
  comparison theorem.
\newblock {\em J. Algebraic Geom.}, 34(1):1--108, 2025.

\bibitem[Dri22]{Prismatization}
Vladimir Drinfeld.
\newblock Prismatization, 2022.

\bibitem[Efi25]{KTheory}
Alexander Efimov.
\newblock {$K$}-theory and localizing invariants of large categories, 2025.

\bibitem[EK99]{CrystallineExt}
Matthew Emerton and Mark Kisin.
\newblock Extensions of crystalline representations, 1999.
\newblock \url{https://www.math.uchicago.edu/~emerton/pdffiles/ext.pdf}.

\bibitem[Fal89]{Faltings}
Gerd Faltings.
\newblock Crystalline cohomology and {$p$}-adic {G}alois-representations.
\newblock In {\em Algebraic analysis, geometry, and number theory ({B}altimore,
  {MD}, 1988)}, pages 25--80. Johns Hopkins Univ. Press, Baltimore, MD, 1989.

\bibitem[FF18]{FarguesFontaine}
Laurent Fargues and Jean-Marc Fontaine.
\newblock Courbes et fibrés en théorie de {H}odge $p$-adique.
\newblock {\em Astérisque}, 406, 2018.

\bibitem[FM87]{FontaineMessing}
Jean-Marc Fontaine and William Messing.
\newblock {$p$}-adic periods and {$p$}-adic \'etale cohomology.
\newblock In {\em Current trends in arithmetical algebraic geometry ({A}rcata,
  {C}alif., 1985)}, volume~67 of {\em Contemp. Math.}, pages 179--207. Amer.
  Math. Soc., Providence, RI, 1987.

\bibitem[Fon82]{Fontaine}
Jean-Marc Fontaine.
\newblock Sur certains types de repr\'esentations {$p$}-adiques du groupe de
  {G}alois d'un corps local;\ construction d'un anneau de {B}arsotti-{T}ate.
\newblock {\em Ann. of Math. (2)}, 115(3):529--577, 1982.

\bibitem[FS24]{FarguesScholze}
Laurent Fargues and Peter Scholze.
\newblock Geometrization of the local {L}anglands correspondence, 2024.

\bibitem[GK00]{GrosseKlonne}
Elmar Gro{\ss}e-Kl\"onne.
\newblock Rigid analytic spaces with overconvergent structure sheaf.
\newblock {\em J. Reine Angew. Math.}, 519:73--95, 2000.

\bibitem[GM25]{GardnerMadapusi}
Zachary Gardner and Keerthi Madapusi.
\newblock An algebraicity conjecture of {D}rinfeld and the moduli of
  $p$-divisible groups, 2025.

\bibitem[Guo24]{GuoBdR}
Haoyang Guo.
\newblock Crystalline cohomology of rigid analytic spaces, 2024.

\bibitem[Hau25]{ArithEtaleComp}
Maximilian Hauck.
\newblock An arithmetic \'etale-crystalline comparison with coefficients in
  crystalline local systems, 2025.

\bibitem[HK94]{HyodoKato}
Osamu Hyodo and Kazuya Kato.
\newblock Semi-stable reduction and crystalline cohomology with logarithmic
  poles.
\newblock Number 223, pages 221--268. 1994.
\newblock P\'eriodes $p$-adiques (Bures-sur-Yvette, 1988).

\bibitem[HM24]{HeyerMann}
Claudius Heyer and Lucas Mann.
\newblock 6-functor formalisms and smooth representations, 2024.

\bibitem[Hyo91]{Hyodo}
Osamu Hyodo.
\newblock On the de {R}ham-{W}itt complex attached to a semi-stable family.
\newblock {\em Compositio Math.}, 78(3):241--260, 1991.

\bibitem[Ked10]{KedlayaDiffEqs}
Kiran Kedlaya.
\newblock {\em {$p$}-adic differential equations}, volume 125 of {\em Cambridge
  Studies in Advanced Mathematics}.
\newblock Cambridge University Press, Cambridge, 2010.

\bibitem[Lur17]{HA}
Jacob Lurie.
\newblock Higher algebra, 2017.
\newblock \url{https://www.math.ias.edu/~lurie/papers/HA.pdf}.

\bibitem[Man22]{MannThesis}
Lucas Mann.
\newblock A $p$-adic $6$-functor formalism in rigid-analytic geometry, 2022.

\bibitem[Mik26]{Mikami}
Yutaro Mikami.
\newblock Finiteness and duality of cohomology of $(\varphi,\gamma)$-modules
  and the 6-functor formalism of locally analytic representations, 2026.

\bibitem[Niz08]{Niziol}
Wies{\l}awa Nizio{\l}.
\newblock Semistable conjecture via {$K$}-theory.
\newblock {\em Duke Math. J.}, 141(1):151--178, 2008.

\bibitem[NN16]{NekovarNiziol}
Jan Nekov{\'a}{\v r} and Wies{\l}awa Nizio{\l}.
\newblock Syntomic cohomology and {$p$}-adic regulators for varieties over
  {$p$}-adic fields.
\newblock {\em Algebra Number Theory}, 10(8):1695--1790, 2016.
\newblock With appendices by Laurent Berger and Fr\'ed\'eric D\'eglise.

\bibitem[RC24]{dRStack}
Juan~Esteban Rodríguez~Camargo.
\newblock The analytic de {R}ham stack in rigid geometry, 2024.

\bibitem[RC25]{GeometricSenTheory}
Juan~Esteban Rodríguez~Camargo.
\newblock Geometric {S}en theory over rigid analytic spaces, 2025.

\bibitem[RJRC25]{SolidLocAnReps}
Joaquín Rodrigues~Jacinto and Juan~Esteban Rodríguez~Camargo.
\newblock Solid locally analytic representations, 2025.

\bibitem[Sch12]{Perfectoid}
Peter Scholze.
\newblock Perfectoid spaces.
\newblock {\em Publ. Math. Inst. Hautes \'Etudes Sci.}, 116:245--313, 2012.

\bibitem[Sch13a]{PAdicHodgeTheory}
Peter Scholze.
\newblock {$p$}-adic {H}odge theory for rigid-analytic varieties.
\newblock {\em Forum Math. Pi}, 1:e1, 77, 2013.

\bibitem[Sch13b]{ScholzeSurvey}
Peter Scholze.
\newblock Perfectoid spaces: a survey.
\newblock In {\em Current developments in mathematics 2012}, pages 193--227.
  Int. Press, Somerville, MA, 2013.

\bibitem[Sch16]{padicHodgeTheoryErratum}
Peter Scholze.
\newblock {$p$}-adic {H}odge theory for rigid-analytic varieties---corrigendum.
\newblock {\em Forum Math. Pi}, 4:e6, 4, 2016.

\bibitem[Sch24]{RealLLC}
Peter Scholze.
\newblock Geometrization of the real local {L}anglands correspondence, 2024.

\bibitem[Sch26a]{BerkovichMotives}
Peter Scholze.
\newblock {B}erkovich motives, 2026.

\bibitem[Sch26b]{SixFunctors}
Peter Scholze.
\newblock Six-functor formalisms, 2026.

\bibitem[Shi22]{Shimizu}
Koji Shimizu.
\newblock A {$p$}-adic monodromy theorem for de {R}ham local systems.
\newblock {\em Compos. Math.}, 158(12):2157--2205, 2022.

\bibitem[Sim97]{Simpson}
Carlos Simpson.
\newblock The {H}odge filtration on nonabelian cohomology.
\newblock In {\em Algebraic geometry---{S}anta {C}ruz 1995}, volume 62, Part 2
  of {\em Proc. Sympos. Pure Math.}, pages 217--281. Amer. Math. Soc.,
  Providence, RI, 1997.

\bibitem[ST97]{Simpson2}
Carlos Simpson and Constantin Teleman.
\newblock De {R}ham's theorem for $\infty$-stacks, 1997.
\newblock \url{https://math.berkeley.edu/~teleman/math/simpson.pdf}.

\bibitem[{Sta}25]{Stacks}
The {Stacks project authors}.
\newblock The stacks project, 2025.
\newblock \url{https://stacks.math.columbia.edu}.

\bibitem[Tsu99]{Tsuji}
Takeshi Tsuji.
\newblock {$p$}-adic \'etale cohomology and crystalline cohomology in the
  semi-stable reduction case.
\newblock {\em Invent. Math.}, 137(2):233--411, 1999.

\bibitem[Wie25]{Wiersig1}
Finn Wiersig.
\newblock Solution and de {R}ham functors for {$D$}-cap-modules, 2025.

\bibitem[Zav26]{Zavyalov}
Bogdan Zavyalov.
\newblock Poincar\'e duality in abstract 6-functor formalisms, 2026.

\end{thebibliography}
\end{document}